\documentclass[a4paper,10pt]{article}
\usepackage[utf8]{inputenc}

\usepackage{amssymb,amsmath,mathrsfs}
\usepackage{amsthm}

\numberwithin{equation}{section}
\usepackage{graphicx}
 
\theoremstyle{plain}
 
\newtheorem{theorem}{Theorem}[section] 
\newtheorem{thm}[theorem]{Theorem} 
\newtheorem{Mainthm}{Theorem}
\newtheorem{prop}[theorem]{Proposition}
\newtheorem{cor}[theorem]{Corollary}

\newtheorem{lem}[theorem]{Lemma}

\theoremstyle{definition}
\newtheorem{defi}[theorem]{Definition}

\theoremstyle{remark}
\newtheorem{rem}[theorem]{Remark}
 
\usepackage[margin=2.5cm]{geometry}
 
\usepackage{tikz}

\usetikzlibrary{positioning, arrows.meta, calc}

\usetikzlibrary{matrix}

 \usepackage[colorlinks,
            linkcolor=blue,
            anchorcolor=red,
            citecolor=green
            ]{hyperref}
            
\begin{document}

\title{The focusing Calogero--Moser derivative nonlinear Schr\"odinger equation
in the Hardy--Zhidkov space: explicit flow formula and dark
multi-solitons}

\author{Ruoci Sun\footnote{School of Mathematical Sciences, Laboratory of Mathematics and
Complex Systems, MOE, Beijing Normal University, 100875 Beijing,
China. Email: ruoci.sun16@bnu.edu.cn}}

\maketitle

\bigskip

\begin{abstract} 
This paper studies the focusing Calogero--Moser derivative nonlinear
Schr\"odinger (CMdNLS) equation  under the nonzero boundary condition
\(|u(t,x)|\to 1\) as \(|x|\to+\infty\). In sharp contrast to the
classical paradigm for the
cubic nonlinear Schr\"odinger (NLS) equation,
 namely that a focusing
nonlinearity supports only bright solitons on a zero background, we
show that the focusing CMdNLS equation admits a rich family of
\emph{dark multi-soliton} solutions on a nonzero background, whose
existence is driven by
 a spectral gap \([-1,0]\) in the point spectrum of the chiral Lax
operator,
  together with the analytic constraint of the
Hardy space.  
 
 The main results are threefold. First, we establish an explicit flow
formula for general solutions in the Hardy--Zhidkov space 
$ \mathcal Z^2_+$,
  via a \emph{semigroup method} based on new
commutator estimates between the Lax--Beurling shift semigroup and
Toeplitz operators. This approach circumvents the fundamental
difficulty that the unbounded infinitesimal generator does not commute
with the time derivative, and
it requires no weighted hypothesis on either the solution or the domain
of validity of the commutator formulas.  
  Second, the algebraic definition of a dark multi-soliton potential is
shown to be equivalent to two spectral conditions on the chiral Lax
operator. Combined with the conservation of a generating functional,
this spectral characterization yields the invariance of the dark
\(N\)-soliton set \(\mathcal U_N\) under the CMdNLS flow.
  Third, the action-angle mapping from \(\mathcal U_N\) to an open subset
of \((\mathbb R/2\pi\mathbb Z)\times\mathbb R^N\times\mathbb R^N\) is
constructed and shown to be bijective, which establishes the complete
integrability of the focusing CMdNLS equation on \(\mathcal U_N\).

Combining these results, we obtain an explicit inverse spectral
formula, from which it follows that the poles of dark multi-soliton
solutions satisfy a complexified Calogero--Moser system. The long time
asymptotics and uniform Sobolev estimates for all pure dark
multi-soliton solutions are also derived. 
\end{abstract}

\bigskip

\bigskip


\bigskip


\newpage
\tableofcontents

\newpage

\section{Introduction}

\noindent The focusing  Calogero--Moser derivative nonlinear Schr\"odinger equation (CMdNLS) reads as 
\begin{equation}\label{CMdNLSfocu} 
i \partial_t u + \partial_x^2 u =    u (i - \mathbb{H}) \partial_x   (|u|^2), \qquad  (t,x) \in \mathbb{R} \times \mathbb{R},
\end{equation}where $\mathbb{H} : L^2(\mathbb{R}; \mathbb{C}) \to L^2(\mathbb{R}; \mathbb{C})$ denotes the Hilbert transform on the real line: 
\begin{equation}\label{HilbertTrIntro}
(\mathbb{H} f)(x):= \mathrm{P.V.} \int_{\mathbb{R}}\frac{f(y)}{\pi(x-y)} \mathrm{d}y, \qquad \forall f \in \mathscr{S} .
\end{equation}Here   $\mathscr{S} = \mathscr{S}(\mathbb{R}; \mathbb{C})$ is the Schwartz space. 
By Abanov--Bettelheim--Wiegmann  \cite{AbanovBW2009}, equation \eqref{CMdNLSfocu} can be regarded as the formal thermodynamic limit of the   Calogero--Moser (CM) model for $N$ particles:
\begin{equation}\label{HcmRmodel}
\mathsf{H}^N_{\mathrm{CM}}=\frac{1}{2} \sum_{j=1}^N p_j^2 +  \sum_{1 \leq j \ne k \leq N} \frac{1}{2(x_j-x_k)^2}, \qquad \forall N \in \mathbb{Z} \cap [2, +\infty).
\end{equation} 
In the pioneering work  \cite{GerLenz22CPAM}, 
  P. G\'erard and E. Lenzmann introduced the Lax pair structure for  
\eqref{CMdNLSfocu} and established
 its local wellposedness in high-regularity   
 \emph{Hardy--Sobolev spaces} $H^s_+$,   defined by 
   \begin{equation}\label{HarSobHsintro}
 \small
 H^s_+ : = \Big\lbrace f \in H^s :   \operatorname{supp} (\widehat{f}) \subset \mathbb{R}_+   \Big\rbrace \quad \text{with} \quad   H^s : = \Big\lbrace f \in \mathscr{S}' :  \| f \|_{H^s }^2 : = \int_{\mathbb{R}} |\widehat{f}(\xi) |^2 (1+ |\xi|^2)^s \tfrac{\mathrm{d}\xi}{2 \pi}  <+\infty  \Big\rbrace,   
 \end{equation}where $\mathbb{R}_+ := [0, +\infty)$, 
 for any $ s > \frac{1}{2}$. 
Here $\mathscr{S}' = \mathscr{S}'  (\mathbb{R} ; \mathbb{C})$ denotes the space of all tempered distributions on $\mathbb{R}$ and $\widehat{f} \in \mathscr{S}'$  is the Fourier transform of $f \in \mathscr{S}'$. Moreover, they have completely classified all traveling solitary waves and all (bright) multi-soliton solutions for  \eqref{CMdNLSfocu} posed on $H^1_+$ in Theorem 1.2 and Proposition 6.2 of \cite{GerLenz22CPAM}.  For any $N \geq 2$, a   bright   $N$-soliton solution for \eqref{CMdNLSfocu} takes the form
\begin{equation}\label{BriNsoliCMdNLS}
\mathfrak{u}^N (t,x) = \sum_{j=1}^N \frac{a_j(t)}{x -z_j(t)} , \qquad \mathfrak{u}^N (t) \in H^1_+, \qquad z_j(t) \in \mathbb{C}_- := \{\zeta \in \mathbb{C}  : \operatorname{Im} \zeta <0 \},
\end{equation}such that $\operatorname{Im} z_j(t) <0$ and $(z_j)_{1 \leq j \leq N}$ solves the complexified version of the Calogero--Moser system \eqref{HcmRmodel}.
  Even though  K. Kim, T. Kim  and S. Kwon have constructed finite-time blow-up
solutions to~\eqref{CMdNLSfocu} in~\cite{KimKimKwonBlowup2024},
 Theorem 6.1 and Theorem 1.3 of \cite{GerLenz22CPAM} show that 
 every
bright multi-soliton solution~\eqref{BriNsoliCMdNLS} exists globally in time
and its Sobolev norm enjoys polynomial growth. Specifically, for any $s > 0$, 
there exists a constant $\mathfrak{T}'>1$, which depends on $s$ and the initial datum $\mathfrak{u}^N(0)$, such that the following formula holds:
\begin{equation}\label{GroSobPolBr}
0 <  \inf_{|t| > \mathfrak{T}'} \tfrac{\|\mathfrak{u}^N (t) \|_{H^s }}{|t|^{2s}}   \leq \sup_{|t| > \mathfrak{T}'} \tfrac{\|\mathfrak{u}^N (t) \|_{H^s }}{|t|^{2s}}  <+ \infty .
\end{equation}
    All the aforementioned results on the focusing CMdNLS
equation~\eqref{CMdNLSfocu} including the local wellposedness theory, the
classification of traveling solitary waves and bright multi-solitons, and the construction of finite-time
blow-up solutions,  are developed in the framework of Sobolev
spaces $H^s$ or Hardy--Sobolev spaces $H^s_+$ with $s>\frac{1}{2}$.
  Since every function in $H^s = H^s (\mathbb{R} ; \mathbb{C})$  
  tends to zero at infinity   
   ($s > \frac{1}{2}$), 
these works are exclusively concerned with the \emph{zero boundary}
regime for equation \eqref{CMdNLSfocu}: 
\begin{equation}\label{zeroBounIntr}
u(t,x)\to 0, \quad \text{as} \quad |x|\to + \infty .
\end{equation}
  However,  the \emph{nonzero boundary} case for the focusing CMdNLS
equation~\eqref{CMdNLSfocu}:
\begin{equation}\label{NonzeBoIntr}
|u(t,x)|\to 1 , \quad \text{as} \quad |x|\to + \infty ,
\end{equation}
has remained largely unexplored.  In particular, the existence of
\emph{dark solitary waves} and \emph{dark multi-soliton solutions}, which are 
 the
natural analogues of the bright solitary waves and bright multi-soliton solutions~\eqref{BriNsoliCMdNLS} in
the nonzero boundary setting \eqref{NonzeBoIntr},
has been an open problem.
   The purpose
of the present paper is to establish an explicit formula for general
solutions of~\eqref{CMdNLSfocu} with nonzero boundary conditions, to
construct and classify its dark multi-soliton solutions, and to
analyze their long time  behavior. \\
 
\noindent A convenient functional framework for studying both zero and
nonzero boundary conditions simultaneously is provided by the
\emph{Zhidkov spaces}, which was originally introduced by
  P.\,E.~Zhidkov~\cite{Zhidkov87CauPro, Zhidkov01KdVNLS}
(under the notation $X^n$) in the context of the nonlinear
Schr\"odinger equation.   
  For any  $n\in\mathbb{N}_+:=\mathbb{Z}\cap[1,+\infty)$, the Zhidkov space
$\mathcal{Z}^n = \mathcal{Z}^n (\mathbb{R} ; \mathbb{C})$ and the Zhidkov norm $\| . \|_{\mathcal{Z}^n}$ are   defined as follows,
\begin{equation}\label{calZnIntro}
  \mathcal{Z}^n := \{\mathtt{u}\in L^\infty(\mathbb{R};\mathbb{C}) :
  \partial_x\mathtt{u}\in H^{n-1}\},  \qquad  \|\mathtt{u}\|_{\mathcal{Z}^n} : =    \|\mathtt{u}\|_{L^{\infty}(\mathbb{R}; \mathbb{C})} + \|\partial_x \mathtt{u}\|_{H^{n-1} } .
\end{equation}
Its natural chiral counterpart, the \emph{Hardy--Zhidkov space}
$\mathcal{Z}^n_+$, consists of those $\mathcal{Z}^n$-functions whose
Fourier transform is supported on the positive half-line:
\begin{equation}\label{HarZhidnIntro}
  \mathcal{Z}^n_+ := \{\mathtt{u}\in\mathcal{Z}^n :
  \operatorname{supp}(\widehat{\mathtt{u}})\subset[0,+\infty)\},
  \qquad \forall  n\in\mathbb{N}_+.
\end{equation}
Then $\mathcal{Z}^n_+$ is a closed subspace of the Banach space $ \bigl(\mathcal{Z}^n , \| . \|_{\mathcal{Z}^n} \bigr) $. 
  Using the Kato iterative scheme as implemented by
Chen~\cite{Xchenccm25ExpFor} for the defocusing CMdNLS  equation, we obtain the
following local wellposedness result for~\eqref{CMdNLSfocu} on the
Zhidkov spaces and Hardy--Zhidkov spaces.

\begin{prop}\label{LWPpropCMdNLSZn}
Fix $n \in \{2,3\}$ and $R>0$. There exists $\mathsf{T} = \mathsf{T}(R) > 0$ such that if $u_0 \in \mathcal{Z}^n$ with $\|u_0\|_{\mathcal{Z}^n} \le R$, there exists a unique solution $u \in \mathcal{C}\bigl( (-\mathsf{T}, \mathsf{T}); \mathcal{Z}^n \bigr)$ to the focusing CMdNLS equation  \eqref{CMdNLSfocu} with initial datum $u(0) = u_0 \in  \mathcal{Z}^n$. Moreover, for any $\mathtt{T} \in (0, \mathsf{T})$, the flow map
\begin{equation}\label{Flowmapfor}
\mathbf{F}_{\mathtt{T}} : u_0 \in \overline{\mathsf{B}}_{\mathcal{Z}^n}(0, R) \longmapsto  u \in \mathcal{C}\bigl( [-\mathtt{T}, \mathtt{T}]; \mathcal{Z}^n \bigr),  
\end{equation}
is continuous. In addition, if the initial datum satisfies $u_0 \in \mathcal{Z}^n_+$ and $\|u_0\|_{\mathcal{Z}^n} \le R$, then the corresponding solution satisfies $u(t) \in \mathcal{Z}^n_+$ for all $t \in (-\mathsf{T}, \mathsf{T})$.
\end{prop}

\noindent Proposition~\ref{LWPpropCMdNLSZn} shows that the flow of
\eqref{CMdNLSfocu} is well defined on the full Zhidkov space
$\mathcal{Z}^n$ ($n=2,3$) and preserves the Hardy--Zhidkov subspace
$\mathcal{Z}^n_+$.  In this paper we shall be exclusively interested in
solutions that evolve in $\mathcal{Z}^2_+$; we call them
\emph{chiral solutions} of~\eqref{CMdNLSfocu}, because they are
composed solely of positive frequencies.  All the principal results
stated below --- the explicit formulas, the spectral characterization
of dark multi-solitons, the action-angle variables, 
 the inverse spectral theory, and the long time
behavior --- are established for chiral solutions.
 
 \subsection{Main Results}
\noindent Before stating the explicit formula for \eqref{CMdNLSfocu}, we need to introduce some key  concepts from Hardy space theory. 
  The upper half-plane is denoted by $\mathbb{C}_+ =  \{ \zeta \in \mathbb{C} : \operatorname{Im} \zeta >0  \}$. 
 Fix $p \in [1, +\infty]$. Set 
\begin{equation}\label{LppmintrDefFor}
L^p_+ := \{f \in L^p(\mathbb{R}; \mathbb{C}) : \operatorname{supp} ( \widehat{f} ) \subset [0, +\infty) \} 
 .
\end{equation} 
  For any $f \in L^p (\mathbb{R} ; \mathbb{C})$, its Poisson integral $\underline{f}: z \in \mathbb{C}_+  \mapsto \underline{f}(z) \in \mathbb{C}$
      is  given by  
\begin{equation}\label{IntroPoiIntDef}
\underline{f}(z) = \mathbf{P}_z (f)  := \frac{\operatorname{Im} z}{\pi} \int_{\mathbb{R}} \frac{f(\mathfrak{s})}{|\mathfrak{s} - z|^2} \, \mathrm{d}\mathfrak{s}   , \qquad \forall z   \in \mathbb{C}_+.
\end{equation}
  Let $\Pi  : = \frac{1}{2} \bigl( \mathrm{id}_{L^2} + i \mathbb{H}\bigr) : L^2(\mathbb{R} ; \mathbb{C}) \to L^2(\mathbb{R} ; \mathbb{C})$ denote the Szeg\H{o} projector. It cancels all negative Fourier modes and preserves all nonnegative Fourier modes. Specifically, for any $f \in L^2 (\mathbb{R} ; \mathbb{C})$,  we have
\begin{equation}\label{PiSzeProIntro}
 \widehat{\Pi f} (\xi) = \widehat{f} (\xi), 
 \qquad 
  \widehat{\Pi f} (\eta) = 0, \qquad \forall \xi  \geq  0 > \eta.
\end{equation}
  By the classical Paley--Wiener theorem, the Poisson integral operator 
$f\mapsto\underline{f}$ provides an isometric isomorphism between the filtered space 
$L^2_+ = \Pi(L^2(\mathbb{R};\mathbb{C}))$ and the Hardy space
$\mathfrak{H}^2(\mathbb{C}_+)$ of holomorphic functions on the upper half-plane 
$\mathbb{C}_+$ with square-integrable boundary values.
   For any $\mathtt{u} \in L^{\infty} (\mathbb{R} ; \mathbb{C})$, 
   the Toeplitz operator $\mathbf{T}_{\mathtt{u}} : L^2_+ \to L^2_+$ is defined as follows:      
\begin{equation}\label{ToepOpIntrodef}
 \mathbf{T}_{ \mathtt{u} } (f):=\Pi ( \mathtt{u} f)  \in L^2_+  , \qquad \forall f \in L^2_+.
\end{equation}
  The next operator plays a central role in the explicit formula:
it is the infinitesimal generator of the left shift semigroup $\bigl( \mathbf{T}_{\exp (-i \eta x )}\bigr)_{\eta \geq 0}$ on
$L^2_+$, and its resolvent will be used to encode the time evolution
of the solution.
 This generator is denoted by $\bigl(- i \mathbf{G} , \mathrm{Dom}(\mathbf{G}) \bigr) $. Its domain of definition is given by 
 \begin{equation}\label{DombfGdefIntro}
\mathrm{Dom}(\mathbf{G}):= \{\varphi \in L^2_+: \widehat{\varphi}|_{\mathbb{R}_+^*} \in H^1(\mathbb{R}_+^*)\} \qquad \text{ with } \qquad  \mathbb{R}_+^* = (0,  +\infty). 
\end{equation}
  The operator $\bigl( \mathbf{G} , \mathrm{Dom}(\mathbf{G}) \bigr) $ is defined via the Fourier transform: for any $\varphi \in  \mathrm{Dom}(\mathbf{G}) $, we have 
\begin{equation}\label{bfGvphidefIntro}
  \widehat{\mathbf{G} \varphi} (\xi) : = i  \partial_{\xi} \widehat{\varphi}(\xi), \qquad \widehat{\mathbf{G} \varphi} (\eta )  =0,  \qquad \forall \xi >0 > \eta    . 
\end{equation}
  By the Sobolev embedding theorem $H^1 (\mathbb{R}_+^*  ) \hookrightarrow \mathcal{C} \cap L^{\infty} (\mathbb{R}_+^* )$, for every $\varphi \in  \mathrm{Dom}(\mathbf{G})$, 
its Fourier transform $\widehat{\varphi}$ admits a right limit at $0$. This limit is denoted by  $\mathbf{I}_+ (\varphi)$. In other words, we have
\begin{equation}\label{bfIplusIntro}
\mathbf{I}_+ (\varphi) : = \widehat{\varphi}( 0^+)  = \lim_{\xi  \to 0^+} \widehat{\varphi}( \xi ) , \qquad  \forall \varphi \in  \mathrm{Dom}(\mathbf{G}) .
\end{equation}
  The Schr\"odinger group $\bigl( e^{it \partial_x^2 } \bigr)_{t \in \mathbb{R}}$ is a strongly continuous group of unitary operators on $L^2_+$. For any $t \in \mathbb{R}$, the conjugate of  $\bigl( \mathbf{G} , \mathrm{Dom}(\mathbf{G}) \bigr) $ by   $e^{it \partial_x^2}$ is denoted by $\bigl( \mathbb{G}_t,  \mathrm{Dom}(\mathbb{G}_t) \bigr)$. To be precise, we define
 \begin{equation} \label{bbGtDefintro}
 \mathrm{Dom}(\mathbb{G}_t) : =    e^{-it \partial_x^2} \bigl(  \mathrm{Dom}  (\mathbf{G}) \bigr), \qquad
 \mathbb{G}_t   := e^{- i t \partial_x^2} \mathbf{G} e^{ i t \partial_x^2}    :  \mathrm{Dom}(\mathbb{G}_t) \to L^2_+ .
 \end{equation}
 The Hille--Yosida theorem implies that $\bigl( i \mathbb{G}_t,  \mathrm{Dom}(\mathbb{G}_t) \bigr)$ is a maximal accretive operator on $L^2_+$. 
  Fix $z \in \mathbb{C}_+$. For every $\mathtt{u} \in L^{\infty} (\mathbb{R} ; \mathbb{C})$, its \textit{difference quotient} function $\mathscr{Q}_z(\mathtt{u} )  \in L^2(\mathbb{R}; \mathbb{C}) $ is defined by
\begin{equation}\label{QzdiffQuoDefIntro}
\mathscr{Q}_z( \mathtt{u} ) (x)  :=   (x -z )^{-1} \bigl( \mathtt{u} (x)  - \underline{\mathtt{u}}(z) \bigr) , \qquad \forall x \in \mathbb{R} .
\end{equation} 
 The first main result of this paper is stated as follows.  
\begin{Mainthm}[Explicit Formula]\label{ExpForMainThm}
Let $u\in\mathcal{C}\bigl((T^-,T^+);\mathcal{Z}^2_+\bigr) $ be a solution to the focusing CMdNLS equation \eqref{CMdNLSfocu} with initial datum $u(0)=u_0\in\mathcal{Z}^2_+$, for some   $T^-<0<T^+$. For any $t \in (T^- , T^+)$ and $z \in \mathbb{C}_+$, the Poisson integral $\underline{u}(t, z) = \underline{u(t)}(z) = \mathbf{P}_z \bigl( u(t) \bigr) \in \mathbb{C}$ can be expressed as follows:
\begin{equation}\label{ExpForCMdNLSfoc}
\underline{u}(t, z) 
= \mathbf{P}_z  e^{it \partial_x^2} (u_0)
 + \tfrac{t}{\pi i}\mathbf{I}_+ \bigl( \mathbb{G}_t - 2 t \mathbf{T}_{u_0} \mathbf{T}_{\overline{u}_0} -z \bigr)^{-1}  \mathbf{T}_{u_0} \mathbf{T}_{\overline{u}_0} e^{ - it \partial_x^2} \mathscr{Q}_z  e^{it \partial_x^2} (u_0) .
\end{equation}
 \end{Mainthm} 
  
 \begin{rem}\label{WelDefExpForRem}
  For every $\mathtt{u}\in L^\infty(\mathbb{R};\mathbb{C})$, the operator
  $\mathbf{T}_{\mathtt{u}}\mathbf{T}_{\overline{\mathtt{u}}}$ is a bounded positive operator  on $L^2_+$.  By the Kato--Rellich theorem for accretive
  perturbations, the operator
  $\mathbb{G}_t - 2t\mathbf{T}_{\mathtt{u}}\mathbf{T}_{\overline{\mathtt{u}}} - z
  : \mathrm{Dom}(\mathbb{G}_t) \to L^2_+$
  is bijective for every $z\in\mathbb{C}_+$ and $t\in\mathbb{R}$.
  Since $\mathrm{Dom}(\mathbf{G})$ is defined by the condition
  $\widehat{\varphi}|_{\mathbb{R}_+^*} \in H^1(\mathbb{R}_+^*)$,
  the Sobolev embedding $H^1(\mathbb{R}_+^*) \hookrightarrow \mathcal{C}\cap L^\infty(\mathbb{R}_+^*)$
  guarantees that $\mathbf{I}_+(\varphi) = \widehat{\varphi}(0^+)$ exists for
  every $\varphi\in\mathrm{Dom}(\mathbf{G})$.  The same holds on
  $\mathrm{Dom}(\mathbb{G}_t) = e^{-it\partial_x^2}(\mathrm{Dom}(\mathbf{G}))$
  because multiplication by the continuous function $e^{-it\xi^2}$ preserves
  the right‑hand limit at the origin.
  The Schr\"odinger group $(e^{it\partial_x^2})_{t\in\mathbb{R}}$ acts as a
  strongly continuous group of bounded operators on the Hardy--Zhidkov space $\mathcal{Z}^2_+$, and
  the difference quotient operator $\mathscr{Q}_z$ maps $\mathcal{Z}^2_+$
  boundedly into the Hardy--Sobolev space $H^2_+$.  Consequently, the right‑hand side
  of the explicit flow formula~\eqref{ExpForCMdNLSfoc} is well defined.
\end{rem}
 
 \begin{rem}\label{PoiIntCvRemL2}
  Every function $f\in L^2_+$ can be recovered as the $L^2$-limit of its
  Poisson integral:
  \begin{equation}\label{LimPoitoOrfL2}
    \lim_{y\to 0^+} \|\underline{f}(\cdot + yi) - f\|_{L^2(\mathbb{R};\mathbb{C})} = 0,
    \qquad \forall f\in L^2_+.
  \end{equation}
    For any general solution $u\in\mathcal{C}\bigl((T^-,T^+);\mathcal{Z}^2_+\bigr) $ to   \eqref{CMdNLSfocu}, it satisfies that $u(t) - e^{it\partial_x^2}(u_0) \in H^2_+ \subset L^2_+$ for all
  $t\in(T^-,T^+)$.
   Since the Poisson integral operator $\mathbf{P}_z  :  L^p (\mathbb{R} ; \mathbb{C})  \to \mathbb{C}$ is linear, for any $  z \in \mathbb{C}_+$ and $p \in [1, +\infty]$,  
  formula~\eqref{ExpForCMdNLSfoc} together with the $L^2$-convergence formula 
  \eqref{LimPoitoOrfL2} uniquely determines $u(t) \in \mathcal{Z}^2_+$ itself, for all $ t\in(T^-,T^+) $. So  \eqref{ExpForCMdNLSfoc} is called the \emph{explicit  flow  formula} for general solutions to the focusing  CMdNLS equation \eqref{CMdNLSfocu} posed on  $\mathcal{Z}^2_+$.
\end{rem}
 
\noindent The following corollary specializes the explicit flow formula \eqref{ExpForCMdNLSfoc} to solutions with a constant background $a \in \mathbb{C}$, which is the natural setting for both bright and  dark multi-solitons.

\begin{cor}\label{CorExpforaaddH}
  Let $u\in\mathcal{C}\bigl((T^-,T^+);\mathcal{Z}^2_+\bigr)$ be a solution to \eqref{CMdNLSfocu}
  with initial datum $u(0)=u_0\in\mathcal{Z}^2_+$.  Assume that there exists $a\in\mathbb{C}$
  such that $u_0-a\in L^2_+$.  Then for every $t\in(T^-,T^+)$ and $z\in\mathbb{C}_+$,
  the Poisson integral $\underline{u}(t,z)=\mathbf{P}_z\bigl(u(t)\bigr)$ is given by
  \begin{equation}\label{expForaaddH2}
    \underline{u}(t,z)
    = a + \tfrac{1}{2\pi i}\,
      \mathbf{I}_+\bigl(\mathbb{G}_t-2t\,\mathbf{T}_{u_0}\mathbf{T}_{\overline{u}_0}-z\bigr)^{-1}
      (u_0-a)\, .
  \end{equation}
\end{cor}
 
\begin{rem}\label{remH2sol}
  When $a=0$, the explicit formula for general $H^2_+$-solutions was established
  by G\'erard--Lenzmann~\cite{GerLenz22CPAM}, Sun~\cite{SUNLMP2024} and
  Killip--Laurens--Vi\c{s}an~\cite{KillLV2311}.   
  In Theorem 4.1 of  \cite{KillLV2311},  
    the authors prove formula \eqref{expForaaddH2} under the assumptions  $u_0 \in H^{\infty}_+ := \bigcap_{s \geq 1} H^s_+$  and  $\langle x \rangle u_0 \in L^2 (\mathbb{R} ; \mathbb{C})$,
    where $\langle x \rangle:=\sqrt{1+|x|^2}$.
      Then Theorem~1.7 of~\cite{KillLV2311}
  further extends the formula to $L^2_+$-solutions of equation~\eqref{CMdNLSfocu} and removes the weighted condition $\langle x \rangle u_0 \in L^2 (\mathbb{R} ; \mathbb{C})$ via the density and continuity argument.
\end{rem}
   
   \noindent Closely related to the weighted $L^2$ condition     
   $\langle x \rangle u_0 \in L^2 (\mathbb{R} ; \mathbb{C})$   
   appearing in Proposition 2.3 and Theorem 4.1 of 
\cite{KillLV2311} is the space $\mathrm{Dom}(\mathbf{X})  \subsetneqq \mathrm{Dom}(\mathbf{G})$, defined as
\begin{equation}\label{DomXIntro}
  \mathrm{Dom}(\mathbf{X}) :=
  \{\varphi \in \mathrm{Dom}(\mathbf{G}) : \widehat{\varphi}(0^+)=0\}
  = \Big\lbrace \varphi \in L^2_+ :
    \int_{\mathbb{R}} |\varphi(x)|^2 (1+|x|^2)\,\mathrm{d}x < +\infty \Big\rbrace
 .
\end{equation}
    
\begin{prop}\label{InvSpPropIntro}
For any $a \in \mathbb{C}$, the affine subspaces $a + \bigl( \mathrm{Dom}(\mathbf{G}) \cap H^3_+ \bigr)$ and $a + \bigl( \mathrm{Dom}(\mathbf{X}) \cap H^3_+ \bigr)$ are both invariant under the $\mathcal{Z}^3_+$-flow of the focusing CMdNLS equation \eqref{CMdNLSfocu}.
\end{prop}

\noindent The proofs of Theorem~\ref{ExpForMainThm} and
Proposition~\ref{InvSpPropIntro} both rely crucially on the Lax pair
structure of~\eqref{CMdNLSfocu} on $\mathcal{Z}^2_+$.  For any
$\mathbf{u} \in L^\infty(\mathbb{R};\mathbb{C})$, the \emph{chiral Lax
operator} $\mathbf{L}_{\mathbf{u}} : H^1_+\to L^2_+$ is defined by
\begin{equation}\label{chiralLaxOpIntr}
  \mathbf{L}_{\mathbf{u}} := -i\partial_x - \mathbf{T}_{\mathbf{u}} \mathbf{T}_{\overline{{\mathbf{u}}}},
  \qquad \mathrm{Dom}(\mathbf{L}_{\mathbf{u}} ) = H^1_+.
\end{equation}
The Kato--Rellich theorem for symmetric perturbations shows that
$\mathbf{L}_{\mathbf{u}} $ is an unbounded self-adjoint operator on $L^2_+$.    
  Next, we introduce dark multi-soliton potentials for   equation \eqref{CMdNLSfocu}. 
   
\begin{defi}\label{DefiDarNsoli}
Let $N \in \mathbb{N}_+$ be a positive integer. A function $\mathbf{u}: \mathbb{R} \to \mathbb{C}$ is called a \emph{dark $N$-soliton} of the focusing CMdNLS equation \eqref{CMdNLSfocu} if it can be expressed in the form
\begin{equation}\label{DarkNsolu}
  \mathbf{u}(x) = \tfrac{P(x)}{Q(x)}, \qquad \forall x \in \mathbb{R},
\end{equation}where $Q \in \mathbb{C}[X]$ is a monic polynomial of degree $N$ with all its roots in the open lower half-plane $\mathbb{C}_- = \{z \in \mathbb{C}: \operatorname{Im}(z) < 0\}$, and $P \in \mathbb{C}[X]$ is another polynomial such that the following identity holds,
\begin{equation}\label{PQiden+}
    |P(x)|^2 - |Q(x)|^2 = i\bigl( \overline{Q(x)}  \partial_x Q(x) - Q(x)   \partial_x \overline{Q(x)} \bigr), \qquad \forall  x \in \mathbb{R}.
\end{equation}Let $\mathcal{U}_N$ denote the set of all dark $N$-solitons of the focusing CMdNLS equation \eqref{CMdNLSfocu}.
\end{defi}

\begin{rem}\label{RmkDarNsolsumCstH1}
  For any integer $N\geq 1$, every dark $N$-soliton $\mathbf{u} \in \mathcal{U}_N$ can be
  written as $\mathbf{u} = e^{i\alpha} + v$, where $\alpha\in\mathbb{T}:=\mathbb{R}/2\pi\mathbb{Z}$
  and $v\in H^1_+$.  The constant $e^{i\alpha}$ coincides with the limit of
  $\mathbf{u} (x)$ as $|x| \to + \infty$ (see formulas~\eqref{ueiapluvDarSoli}
  and~\eqref{uinftyDarSoli} in Remark~\ref{Rkuinftyv} for the precise
  algebraic description of $v$).
\end{rem}

\noindent The second main result of this paper gives two equivalent
descriptions for a function of the form $\mathbf{u} = e^{i\alpha} + v$ with $\alpha \in \mathbb{T}$ and 
$v\in H^1_+$ to be a dark multi-soliton potential of equation~\eqref{CMdNLSfocu}.

\begin{Mainthm}[Spectral Characterization]\label{MainThmSpectChar}
Let $N\geq 1$ be a positive integer. Assume that $\mathbf{u} = e^{i \alpha } + v$ for some 
 $\alpha \in \mathbb{T} $ and $v \in H^1_+$. Let $\mathscr{H}_{\mathrm{pp}}(\mathbf{L}_{\mathbf{u}}) \subset H^1_+$ denote the pure point spectral subspace of the Lax operator $\mathbf{L}_{\mathbf{u}} = \mathbf{L}_{\mathbf{u}}^*$ defined by \eqref{chiralLaxOpIntr}. Then the following three statements are equivalent: \\
\noindent $(\mathbf{a}).$ The function $\mathbf{u}$ is a dark $N$-soliton potential of   equation \eqref{CMdNLSfocu}, i.e. $\mathbf{u} \in \mathcal{U}_N$. \\
\noindent $(\mathbf{b}).$ We have $\dim_{\mathbb{C}}\mathscr{H}_{\mathrm{pp}}(\mathbf{L}_{\mathbf{u}}) =N$, $\mathscr{H}_{\mathrm{pp}}(\mathbf{L}_{\mathbf{u}}) \subset \mathrm{Dom}(\mathbf{G})$,  $v = \mathbf{u} - e^{i \alpha} \in \mathscr{H}_{\mathrm{pp}}(\mathbf{L}_{\mathbf{u}})$ and 
\begin{equation}\label{HppLuinvbyG}
\mathbf{G} \bigl( \mathscr{H}_{\mathrm{pp}}(\mathbf{L}_{\mathbf{u}}) \bigr) \subset \mathscr{H}_{\mathrm{pp}}(\mathbf{L}_{\mathbf{u}}).
\end{equation}$(\mathbf{c}).$ The number $-1$ is not an eigenvalue of $\mathbf{L}_{\mathbf{u}}$. Moreover,  we have $\dim_{\mathbb{C}}\mathscr{H}_{\mathrm{pp}}(\mathbf{L}_{\mathbf{u}}) =N$ and  
\begin{equation}\label{vinHppandG-1Hpp}
v = \mathbf{u} - e^{i \alpha} \in \mathscr{H}_{\mathrm{pp}}(\mathbf{L}_{\mathbf{u}}) \cap  \mathrm{Dom}(\mathbf{G}) \cap  \mathbf{G}^{-1} \bigl(\mathscr{H}_{\mathrm{pp}}(\mathbf{L}_{\mathbf{u}})  \bigr).
\end{equation}
\end{Mainthm}

\noindent
For any dark multi-soliton  $\mathbf{u}$, the   subspace $\mathscr{H}_{\mathrm{pp}}(\mathbf{L}_{\mathbf{u}}) $
is finite-dimensional and invariant under the generator
$\mathbf{G}$ by~\eqref{HppLuinvbyG}.  This observation yields the
uniqueness of the polynomial pair $(P,Q)$ in
Definition~\ref{DefiDarNsoli} and identifies the denominator $Q$ with
the characteristic polynomial of
$\mathbf{G}|_{\mathscr{H}_{\mathrm{pp}}(\mathbf{L}_{\mathbf{u}})}$.

\begin{cor}\label{CorQuPu}
  For any \(N\in\mathbb{N}_+\) and \(\mathbf{u}\in\mathcal{U}_N\), there
  exists a unique pair
  \((P_{\mathbf{u}},Q_{\mathbf{u}})\in\mathbb{C}[X]\times\mathbb{C}[X]\)
  such that \(\mathbf{u}=P_{\mathbf{u}}/Q_{\mathbf{u}}\),
  \(\deg P_{\mathbf{u}}=\deg Q_{\mathbf{u}}=N\), \(Q_{\mathbf{u}}\) is
  monic with roots in \(\mathbb{C}_-\), and
  \((P_{\mathbf{u}},Q_{\mathbf{u}})\) satisfies~\eqref{PQiden+}.
  Moreover, \(Q_{\mathbf{u}}\) is precisely the characteristic
  polynomial of the restriction of \(\mathbf{G}\) to
  \(\mathscr{H}_{\mathrm{pp}}(\mathbf{L}_{\mathbf{u}})\), that is,
  \begin{equation}\label{CharPolyQuGHpp1for}
    Q_{\mathbf{u}}(x)
    :=
    \det\bigl(
      x-\mathbf{G}|_{\mathscr{H}_{\mathrm{pp}}(\mathbf{L}_{\mathbf{u}})}
    \bigr),
    \qquad
    \forall x\in\mathbb{R}.
  \end{equation}
\end{cor}

\begin{rem}\label{BriDarkSpecLuComp}
  There is a fundamental contrast between the spectrum of the Lax operator
  for bright multi-soliton potentials in G\'erard--Lenzmann \cite{GerLenz22CPAM} and for the
  dark multi-soliton potentials introduced in Definition~\ref{DefiDarNsoli} of this paper.
  If $\mathbf{u}\in\mathcal{U}_N$ is a dark $N$-soliton potential,
  Lemma~\ref{LemSpLutwoseg} shows that every eigenvalue of
  $\mathbf{L}_{\mathbf{u}}$ lies in $\mathbb{R}\setminus[-1,0]
  =(-\infty,-1)\cup(0,+\infty)$.  In particular, $0$ is never an eigenvalue
  for a dark multi-soliton.  This is in stark contrast with the case of bright
  multi-soliton potentials, where $0$ is always an eigenvalue of the Lax
  operator (see~\cite[Proposition~5.3]{GerLenz22CPAM}), and for every
  $\lambda\in\mathbb{R}$ there exists a bright multi-soliton whose Lax
  operator admits $\lambda$ as an eigenvalue (see~\cite[Proposition~6.2 and
  Subsection~6.3]{GerLenz22CPAM}). 
  We emphasize that the spectral gap $[-1,0]$ appearing here is intrinsically
  tied to the nonzero boundary condition \eqref{NonzeBoIntr}, i.e. $|u(x)|\to 1$ as $|x|\to\infty$.
  The methods developed in the present paper are expected to extend to more
  general nonzero boundary conditions   
  $|u(x)|\to r$ as $|x|\to+\infty$, with $r >0$. 
\end{rem}
     
\noindent Theorem~\ref{MainThmSpectChar} translates the purely algebraic
definition of a dark $N$-soliton potential, i.e. the rational function
$\mathbf{u} = \frac{P}{Q}$ satisfying~\eqref{PQiden+},  into a spectral condition on the
Lax operator $\bigl( \mathbf{L}_{\mathbf{u}}, \mathrm{Dom}(\mathbf{L}_{\mathbf{u}}) \bigr)$
and the generator $\bigl( \mathbf{G}, \mathrm{Dom}(\mathbf{G}) \bigr)$.
 Combining Theorem  $\ref{MainThmSpectChar}$ with the invariant subspace results of Proposition $\ref{InvSpPropIntro}$, we can deduce the invariance of the dark $N$-soliton subset $\mathcal{U}_N$ under the $\mathcal{Z}^2_+$-flow of equation \eqref{CMdNLSfocu}.

\begin{prop}\label{PropInvcalUN}
 Let $u\in\mathcal{C}\bigl((T^-,T^+);\mathcal{Z}^2_+\bigr) $ be a solution to the focusing CMdNLS equation \eqref{CMdNLSfocu}. Assume that  $u(0) =u_0 \in \mathcal{U}_N$, for some $N \in \mathbb{N}_+$, then $u(t)$ belongs to $\mathcal{U}_N$ for all $t \in (T^-, T^+)$.  
 Moreover, every pure dark multi-soliton solution to \eqref{CMdNLSfocu} exists globally in time.  
\end{prop}

\noindent Although the global well‑posedness of the focusing CMdNLS
equation~\eqref{CMdNLSfocu} with the nonzero boundary
condition~\eqref{NonzeBoIntr} remains an open problem, the pure dark
multi‑soliton solutions never blow up in finite time.
The third main result of this paper shows that, when restricted to each
$\mathcal{U}_N$, equation~\eqref{CMdNLSfocu} admits generalized action–angle
variables that linearize its flow.  The set of actions is given by 
\begin{equation}\label{OmegaNIntro}
  \Omega_N:=\bigl\{(\lambda_1,\lambda_2,\cdots,\lambda_N)\in
  \bigl(\mathbb{R}\setminus[-1,0]\bigr)^{1\times N}:
  \lambda_1<\lambda_2<\cdots<\lambda_N\bigr\},
  \qquad N\in\mathbb{N}_+.
\end{equation}
The third main result is stated as follows.

\begin{Mainthm}[Complete Integrability]\label{ActionAngleMainThm}
  For any $N\in\mathbb{N}_+$, the set $\mathcal{U}_N$ admits a global
  parametrization that linearizes the flow of the focusing CMdNLS
  equation~\eqref{CMdNLSfocu}.  More precisely, there exists a bijection
  \begin{equation}\label{PhiNintro}
    \boldsymbol{\Phi}_N: \mathbf{u}\in\mathcal{U}_N
    \longmapsto \boldsymbol{\Phi}_N(\mathbf{u})
    =\bigl(\boldsymbol{\alpha}(\mathbf{u});\,
            \vec{\boldsymbol{\lambda}}(\mathbf{u});\,
            \vec{\boldsymbol{\gamma}}(\mathbf{u})\bigr)
    \in \mathbb{T}\times\Omega_N\times\mathbb{R}^{1\times N},
  \end{equation}
  whose components satisfy the following properties under the evolution
  of~\eqref{CMdNLSfocu}:
  \begin{enumerate}
    \item \textbf{Conserved quantities.}
    The total phase $\boldsymbol{\alpha}:\mathcal{U}_N\to\mathbb{T}$ and the
    spectral data
    $\vec{\boldsymbol{\lambda}}=(\boldsymbol{\lambda}_1,\cdots,
      \boldsymbol{\lambda}_N):\mathcal{U}_N\to\Omega_N$
    are constants of motion.
    \item \textbf{Linear angle evolution.}
    The coordinates
    $\vec{\boldsymbol{\gamma}}=(\boldsymbol{\gamma}_1,\cdots,
      \boldsymbol{\gamma}_N):\mathcal{U}_N\to\mathbb{R}^{1\times N}$
    undergo uniform linear motion.
  \end{enumerate}
  More explicitly, if $u:t\in\mathbb{R}\mapsto u(t)\in\mathcal{U}_N$ solves
  \eqref{CMdNLSfocu} with initial datum $u_0\in\mathcal{U}_N$, then 
  \begin{equation}\label{TimEvoactangl}
    \boldsymbol{\alpha}(u(t))=\boldsymbol{\alpha}(u_0),\qquad
    \boldsymbol{\lambda}_j(u(t))=\boldsymbol{\lambda}_j(u_0),\qquad
    \boldsymbol{\gamma}_j(u(t))=\boldsymbol{\gamma}_j(u_0)+2t\boldsymbol{\lambda}_j(u_0),  \qquad \forall t \in \mathbb{R}, 
  \end{equation}for every
  $j\in\{1,2 , \cdots,N\}$.
     Moreover, for any $\omega\in\mathbb{T}$ and $u_0\in\mathcal{U}_N$, we have 
  \begin{equation}\label{ActAnginvbyPhRot}
    \boldsymbol{\alpha}(e^{i\omega}u_0)=\boldsymbol{\alpha}(u_0)+\omega\in\mathbb{T},
    \qquad
    \vec{\boldsymbol{\lambda}}(e^{i\omega}u_0)=\vec{\boldsymbol{\lambda}}(u_0) \in \Omega_N,
    \qquad
    \vec{\boldsymbol{\gamma}}(e^{i\omega}u_0)=\vec{\boldsymbol{\gamma}}(u_0) \in  \mathbb{R}^{1\times N}. 
  \end{equation}
  Thus the family
  $(\boldsymbol{\lambda}_j,\boldsymbol{\gamma}_j)_{1\leq  j\leq N}$ constitutes
  a set of generalized action–angle variables for the reduced dynamics
  of~\eqref{CMdNLSfocu} on the quotient set   
  $\mathcal{U}_N / \mathbb{T}$.
\end{Mainthm}

\begin{rem}\label{CirGractUNrem}
  For any integer $N\ge 1$, the circle group $(\mathbb{T},+)$ acts on
  $\mathcal{U}_N$ by phase rotation:
  \begin{equation}\label{bbTactcalUNdef}
    (\omega, \mathbf{u}) \in \mathbb{T}\times\mathcal{U}_N
    \;\longmapsto\; \omega\star\mathbf{u} := e^{i\omega}\,\mathbf{u}
    \in \mathcal{U}_N.
  \end{equation}
  The quotient set $\mathcal{U}_N/\mathbb{T}$ consists of all the orbits
  $[\mathbf{u}] = \{e^{i\omega}\,\mathbf{u} : \omega\in\mathbb{T}\}$, for any
  $\mathbf{u}\in\mathcal{U}_N$.
\end{rem}

\noindent By Remark~\ref{RmkDarNsolsumCstH1}, the phase parameter can be
chosen as the limit of the dark multi-soliton potential at infinity:
$e^{i \boldsymbol{\alpha}(\mathbf{u}) } = \mathbf{u}(\infty) = \lim_{|x|\to+\infty}\mathbf{u}(x) \in \mathbb{C}$.
The action variables $(\boldsymbol{\lambda}_j(\mathbf{u}))_{1\leq j\leq N}$ are
exactly the eigenvalues of the Lax operator $\mathbf{L}_{\mathbf{u}}$, while the
generalized angle variables   $(\boldsymbol{\gamma}_j(\mathbf{u}))_{1\leq j\leq N}$ depend on the generator
$\mathbf{G}$ and the eigenvectors of $\mathbf{L}_{\mathbf{u}}$.  For this
reason, the action–angle mapping $\boldsymbol{\Phi}_N$ in~\eqref{PhiNintro} is
referred to as the \emph{direct spectral transform} for dark $N$-soliton
potentials.  Its inverse, the \emph{inverse spectral transform}
$\boldsymbol{\Phi}_N^{-1} : \mathbb{T}\times\Omega_N\times\mathbb{R}^{1\times N} \to  \mathcal{U}_N$,  
 expresses every dark multi‑soliton potential in
terms of its action–angle variables and the phase parameter.  To write $\boldsymbol{\Phi}_N^{-1}$
explicitly, we introduce the inverse spectral matrix and the spectral
coordinate vectors of a dark $N$-soliton $\mathbf{u}\in\mathcal{U}_N$.

\begin{defi}\label{bfMVfrakXYdefintro}
  For any integer $N\geq 1$, let $\mathbf{u}$ be a dark $N$-soliton potential, i.e.  $\mathbf{u}\in\mathcal{U}_N$. Its 
  \emph{inverse spectral matrix}
  $\mathbf{M}(\mathbf{u}) = (\mathbf{M}_{kj}(\mathbf{u}))_{1\le k,j\le N}
  \in\mathbb{C}^{N\times N}$ is defined by
  \begin{equation}\label{InvSpMatIntro}
    \mathbf{M}_{kj}(\mathbf{u})
    = \tfrac{i}{\boldsymbol{\lambda}_k(\mathbf{u})-\boldsymbol{\lambda}_j(\mathbf{u})}
      \sqrt{\tfrac{\boldsymbol{\lambda}_j(\mathbf{u})(1+\boldsymbol{\lambda}_k(\mathbf{u}))}
                 {\boldsymbol{\lambda}_k(\mathbf{u})(1+\boldsymbol{\lambda}_j(\mathbf{u}))}},
    \qquad
    \mathbf{M}_{jj}(\mathbf{u})
    = \boldsymbol{\gamma}_j(\mathbf{u})
      - \tfrac{i}{2\boldsymbol{\lambda}_j(\mathbf{u})(1+\boldsymbol{\lambda}_j(\mathbf{u}))},
  \end{equation}
  with  $k \ne j \in \{1,2, \cdots, N\}$.
  Its  \emph{spectral coordinate vectors}
  $\mathfrak{X}(\mathbf{u}) = (\mathbf{x}_1 (\mathbf{u}), \mathbf{x}_2 (\mathbf{u}), \cdots,\mathbf{x}_N(\mathbf{u}))^{\mathrm{T}}
  \in\mathbb{R}^N$ and
  $\mathfrak{Y}(\mathbf{u}) = (\mathbf{y}_1 (\mathbf{u}), \mathbf{y}_2 (\mathbf{u}) ,\cdots,\mathbf{y}_N(\mathbf{u}))^{\mathrm{T}}
  \in\mathbb{R}^N$ are given by
  \begin{equation}\label{introxjuyjufor}
    \mathbf{x}_j(\mathbf{u}) = \sqrt{\tfrac{1+\boldsymbol{\lambda}_j(\mathbf{u})}{\boldsymbol{\lambda}_j(\mathbf{u})}},
    \qquad
    \mathbf{y}_j(\mathbf{u}) = \tfrac{\mathrm{sign}(\boldsymbol{\lambda}_j(\mathbf{u}))}
                                 {\sqrt{\boldsymbol{\lambda}_j(\mathbf{u})(1+\boldsymbol{\lambda}_j(\mathbf{u}))}} , \qquad \forall j \in \{1,2, \cdots, N\}.
  \end{equation}
  Finally, its \emph{velocity matrix}
  $\mathbf{V}(\mathbf{u})\in   \mathrm{GL}_N (\mathbb{R})$ is the diagonal matrix
  \begin{equation}\label{VeloMatIntro}
    \mathbf{V}(\mathbf{u})
    = \operatorname{diag}\bigl(\boldsymbol{\lambda}_1(\mathbf{u}),\,
                              \boldsymbol{\lambda}_2(\mathbf{u}),\,
                              \cdots,\,
                              \boldsymbol{\lambda}_N(\mathbf{u})\bigr).
  \end{equation}
\end{defi}

\noindent With the inverse spectral matrix and the spectral coordinate
vectors at hand, every pure dark multi‑soliton solution
to the focusing CMdNLS equation~\eqref{CMdNLSfocu} can be expressed explicitly in terms of the time
variable and the action–angle variables of its initial datum.

\begin{prop}\label{InvSpecFortxintro}
  For any $N\in\mathbb{N}_+$, let $u:t\in\mathbb{R}\mapsto u(t)\in\mathcal{U}_N$
 solve~\eqref{CMdNLSfocu} with   $u(0) = u_0\in\mathcal{U}_N$.  Then  $\sigma\bigl(\mathbf{M}(u_0) + 2t\mathbf{V}(u_0)\bigr) \subset \mathbb{C}_-$,
 for any $t\in\mathbb{R}$.  
 Moreover, the following \emph{inverse spectral formula} holds:
  \begin{equation}\label{invSPEFortxIntro}
   u(t,x) = e^{i\boldsymbol{\alpha}(u_0)}
         \Bigl( 1 + i\,
           \mathfrak{Y}(u_0)^{\mathrm{T}}
           (\mathbf{M}(u_0) + 2t\mathbf{V}(u_0) - x)^{-1}
           \mathfrak{X}(u_0)
         \Bigr),
    \qquad \forall (t,x)\in\mathbb{R}\times\mathbb{R}.
  \end{equation}
  \end{prop}

\begin{rem}\label{RemInSpForEXPFlFor}
  Formula~\eqref{invSPEFortxIntro} can be regarded as the special case of the
  explicit formulas~\eqref{ExpForCMdNLSfoc} and~\eqref{expForaaddH2} for the
  CMdNLS equation~\eqref{CMdNLSfocu} when restricted to the dark $N$-soliton
  subset $\mathcal{U}_N$.
\end{rem}

  \noindent
The inverse spectral formula~\eqref{invSPEFortxIntro} reduces the
description of a pure dark $N$-soliton solution of   \eqref{CMdNLSfocu} to two  vectors $ \mathfrak{X}(u_0),  \mathfrak{Y}(u_0) \in \mathbb{R}^N$ and  the resolvent of the
 time-dependent matrix
$\mathbf{M}(u_0)+2t\mathbf{V}(u_0) \in \mathrm{GL}_N(\mathbb{C}) $. 
Since the velocity matrix
$\mathbf{V}(u_0) \in \mathrm{GL}_N(\mathbb{R})$ is diagonal, the large time analysis of this
family falls within the scope of Kato's perturbation theory. 
  For large $|t|$, expanding each  eigenvalue  of the  perturbed matrix $2\mathbf{V} + t^{-1} \mathbf{M}(u_0) \in \mathrm{GL}_N(\mathbb{C}) $ in powers of $t^{-1}$ via the Riesz projector leads to a  Rayleigh--Schr\"{o}dinger series. This yields a complete description of the long time  asymptotics of these eigenvalues, 
  and consequently of the solution itself through a partial
fraction (pole ansatz) expansion.

\begin{prop}\label{propLoTiBeCMsysIntro}
 For any $N \in \mathbb{N}_+$ and $u_0\in\mathcal{U}_N$, there exists a constant
  $\widetilde{\mathcal{T}} = \widetilde{\mathcal{T}}(N, u_0) > 1$ such that
  for all $t \in \mathbb{R}\setminus[-\widetilde{\mathcal{T}},\widetilde{\mathcal{T}}] $,    the matrix
  $\mathbf{M}(u_0)+2t\mathbf{V}(u_0) \in \mathrm{GL}_N(\mathbb{C}) $ possesses $N$ distinct eigenvalues
  $\mathbf{z}_1(t), \mathbf{z}_2(t),\cdots,\mathbf{z}_N(t)\in\mathbb{C}_-$,
    such that  $\mathbf{z}_j : \mathbb{R}\setminus[-\widetilde{\mathcal{T}},\widetilde{\mathcal{T}}] \to \mathbb{C}$ is a real analytic function and   satisfies
  \begin{equation}\label{bfzjtLonTimIntro}
 \sup_{|t|>\widetilde{\mathcal{T}}}
    \bigl|t\bigl(\mathbf{z}_j(t)-2\boldsymbol{\lambda}_j(u_0)t
      -\mathbf{M}_{jj}(u_0)\bigr)\bigr| < +\infty , \qquad \forall j \in \{1,2, \cdots, N\}.
  \end{equation}
  Let $u:t\in\mathbb{R}\mapsto u(t)\in\mathcal{U}_N$ be  the  solution
  to equation \eqref{CMdNLSfocu} with initial datum $u(0)=u_0 \in \mathcal{U}_N$.  For any $ t \in \mathbb{R}\setminus[-\widetilde{\mathcal{T}},\widetilde{\mathcal{T}}]$,   the solution $u$ admits the following partial fraction expansion
  \begin{equation}\label{ParFracExputIntro}
    u(t,x)=e^{i\boldsymbol{\alpha}(u_0)}
           \Bigl(1+\sum_{j=1}^{N}\frac{\mathbf{a}_j(t)}{x-\mathbf{z}_j(t)}\Bigr),
    \qquad\forall\,x\in\mathbb{R},
  \end{equation}
  with coefficients $\mathbf{a}_1, \mathbf{a}_2,\cdots,\mathbf{a}_N\in
  \mathcal{C}^2\bigl(\mathbb{R}\setminus[-\widetilde{\mathcal{T}},\widetilde{\mathcal{T}}];
  \mathbb{C}\bigr)$.  Moreover, the poles
  $\mathbf{z}_1, \mathbf{z}_2 ,\cdots,\mathbf{z}_N  $ satisfy  
  \begin{equation}\label{IntroCaloMoserSys}
    \partial_{t}^{2}\mathbf{z}_j(t)
    = \sum_{k\in\mathbb{I}_{j}^{N}}
      \frac{8}{(\mathbf{z}_j(t)-\mathbf{z}_k(t))^{3}} 
    \qquad  \text{with} \qquad
    \mathbb{I}_{j}^{N}:=\{1,2, \cdots,N\}\setminus\{j\},
  \end{equation}
  for every $j\in\{1,2, \cdots,N\}$ and all
  $t\in\mathbb{R}\setminus[-\widetilde{\mathcal{T}},\widetilde{\mathcal{T}}]$.
\end{prop}
   
\noindent
The pole dynamics of pure dark multi-soliton solutions
of~\eqref{CMdNLSfocu} is governed by the complexified
Calogero--Moser system~\eqref{IntroCaloMoserSys}.  The same is true for
pure bright multi-soliton solutions, as shown by G\'erard--Lenzmann~\cite[formula~(6.4)]{GerLenz22CPAM}.  Thus, from the
viewpoint of pole dynamics, the zero boundary
regime~\eqref{zeroBounIntr} and the nonzero boundary
regime~\eqref{NonzeBoIntr} share a common underlying integrable
structure.
   However, the long time asymptotics of the poles in the two cases are
fundamentally different.  For a pure  dark $N$-soliton solution \eqref{invSPEFortxIntro}, every
pole $\mathbf{z}_j(t)$ lies in the lower half-plane $\mathbb{C}_-$ and 
its imaginary part converges to a strictly negative limit determined by $\boldsymbol{\lambda}_j(u_0)\in\mathbb{R}\setminus[-1,0]$.
   To be precise, for large $|t|$, we have    
\begin{equation}\label{LTAsympImbfzj}
  \operatorname{Im}\mathbf{z}_j(t)
  = - \tfrac{1}{2\boldsymbol{\lambda}_j(u_0)(1+\boldsymbol{\lambda}_j(u_0))}
    + \mathcal{O}(t^{-1}), \qquad \forall j\in\{1, 2, \cdots,N\} .
\end{equation}
  In contrast, for a pure bright $N$-soliton solution
\eqref{BriNsoliCMdNLS},  Lemma~6.2 in Subsection 6.5 of G\'erard--Lenzmann
\cite{GerLenz22CPAM} shows that the situation is completely 
different. When
$N\geq 2$,  let $\lambda_1=0$ and
$\lambda_2<\cdots<\lambda_N$ be the distinct eigenvalues of the Lax
operator associated to the initial bright $N$-soliton profile $\mathfrak{u}^N(0)  \in H^1_+$.  
  Then only the first pole $z_1(t)\in\mathbb{C}_-$ has an
imaginary part tending to a strictly negative number
$-\rho=-\rho \bigl(\mathfrak{u}^N(0) \bigr)<0$, while all the other poles $z_k(t) \in \mathbb{C}_- $,
$2\leq k\leq N$, approach the real axis at a $t^{-2}$ rate:
\begin{equation}\label{BriSolPoLTasym}
  \operatorname{Im}z_1(t)
  = -\rho + \mathcal{O}(t^{-1}),
  \qquad
  \operatorname{Im}z_k(t)
  = -\tfrac{\rho}{4\lambda_k^{4}t^{2}} + \mathcal{O}(t^{-3}), \qquad \forall  k\in\mathbb{I}_{1}^{N} =\{2, \cdots, N\} .
\end{equation} 
  This striking contrast---the same Calogero--Moser  dynamics governing the poles,   yet
completely different long time asymptotic behavior---is one of the most
intriguing features distinguishing the bright and dark multi-soliton
regimes.  For pure dark multi-soliton solutions, the fact that all poles remain uniformly
away from the real axis is the key mechanism behind their uniform
Sobolev bounds at large time.  This is summarized in the following
corollary.

\begin{cor}\label{corLongtimebehave}
  For any integer $N\geq 1$, let
  $u:t\in\mathbb{R}\mapsto u(t)\in\mathcal{U}_N$ be a pure dark
  $N$-soliton solution of the CMdNLS equation~\eqref{CMdNLSfocu} with
  initial datum $u(0)=u_0\in\mathcal{U}_N$.  Then for every $s>0$,
  there exists a constant $\mathfrak{T}=\mathfrak{T}(s,u_0)>1$ such
  that the following estimate holds:
  \begin{equation}\label{LongTimSoboEsti}
    0 < \inf_{|t|>\mathfrak{T}} \|u(t)-u_0(\infty)\|_{H^s_+}
      \leq \sup_{|t|>\mathfrak{T}} \|u(t)-u_0(\infty)\|_{H^s_+}
      < +\infty .
  \end{equation}
\end{cor}

\begin{rem}\label{remDarkBrightSobDiff}
  Corollary~\ref{corLongtimebehave} reveals another significant
  difference between pure dark multi-soliton solutions and pure bright
  multi-soliton solutions of~\eqref{CMdNLSfocu}.
  On the one hand, by estimate~\eqref{LongTimSoboEsti}, every pure dark
  $N$-soliton solution~\eqref{invSPEFortxIntro} can be decomposed as
  the sum of its nonzero boundary value at infinity and a perturbation
  which remains uniformly bounded in time with respect to the Sobolev
  norm $H^s$, for any $s>0$.    
  On the other hand, for every pure bright multi-soliton
  solution~\eqref{BriNsoliCMdNLS} and every $s>0$, the 
  $H^s$-norm grows polynomially at large time, as stated
  in~\eqref{GroSobPolBr}.    
This is yet another manifestation of the contrast between the bright
and dark multi-soliton regimes.   
\end{rem}

\noindent 
The next corollary is another consequence of the inverse spectral formula \eqref{invSPEFortxIntro}.

\begin{cor}\label{DarkSolRescor}
For any $N \in \mathbb{N}_+$, let $u : t \in \mathbb{R} \mapsto u (t) \in \mathcal{U}_N$ be a  pure dark $N$-soliton solution to the focusing CMdNLS equation \eqref{CMdNLSfocu} with initial datum $u(0) = u_0 \in \mathcal{U}_N$.  
  For any $j \in \{1,2, \cdots, N \}$, the total phase $\boldsymbol{\alpha}(u_0) \in \mathbb{T}$ and the action variable $\boldsymbol{\lambda}_j( u_0 )  \in \mathbb{R} \setminus [-1, 0]$ are defined in \eqref{PhiNintro} and  \eqref{TimEvoactangl}. The $jj$-entry of the inverse spectral matrix $\mathbf{M} (u_0) \in \mathbb{C}^{N \times N}$, denoted by   $\mathbf{M}_{jj}(u_0) \in \mathbb{C}_-$, is given by \eqref{InvSpMatIntro}.
    The spectral coordinate vectors $\mathfrak{X}(u_0) , \mathfrak{Y}(u_0)  \in \mathbb{R}^N$ are defined in  \eqref{introxjuyjufor}. We define the diagonal matrix  
      $\mathfrak{K}^{u_0} (t,x) = \operatorname{diag}\bigl(\mathfrak{K}^{u_0}_{1} (t,x)  ,\,
                              \mathfrak{K}^{u_0}_{2} (t,x) ,\,
                              \cdots,\,
                     \mathfrak{K}^{u_0}_{N} (t,x) \bigr) \in  \mathbb{C}^{N \times N}$ by
\begin{equation}\label{bfDeltajdef}
\mathfrak{K}^{u_0}_j  (t,x) : = \mathbf{M}_{jj}(u_0) + 2 t \boldsymbol{\lambda}_j( u_0 ) - x \in \mathbb{C}_- , \qquad \forall j \in \{1,2, \cdots, N \} .
\end{equation} 
 The reference function $u_{\infty} \in \mathcal{C}^{\infty} \bigl( \mathbb{R} \times 
 \mathbb{R} ; \mathbb{C} \bigr) $  is defined as follows,  
  \begin{equation}\label{uinftyDef}
  u_{\infty}  (t,x) 
   : = e^{i\boldsymbol{\alpha}(u_0)}
         \Bigl( 1 + i\,
           \mathfrak{Y}(u_0)^{\mathrm{T}}
           \bigl(  \mathfrak{K}^{u_0} (t,x) \bigr)^{-1} 
           \mathfrak{X}(u_0)
         \Bigr), \qquad \forall (t,x)  \in \mathbb{R} \times 
 \mathbb{R}.
  \end{equation} 
   Then for any $R >0 $, we have $\lim_{|t| \to + \infty} \big\| u(t) -  u_{\infty}  (t) \big\|_{L^2  ( [-R, R] ; \mathbb{C}  )} = 0 $.      
   Moreover, for every $x\in\mathbb{R}$, the denominator
$u_\infty(t,x)-e^{i\boldsymbol{\alpha}(u_0)}$ does not vanish for all
sufficiently large $|t|$, and we have
\begin{equation}\label{DSolResL2PWisCv}
  \lim_{|t|\to+\infty}
  \tfrac{u(t,x)-e^{i\boldsymbol{\alpha}(u_0)}}
       {u_\infty(t,x)-e^{i\boldsymbol{\alpha}(u_0)}}
  = 1 .
\end{equation}  
\end{cor}

\noindent Thus, at large times, every pure dark \(N\)-soliton solution can be
approximated locally by a superposition of \(N\) one-pole terms
$\sum_{j=1}^N
\frac{i e^{i\boldsymbol{\alpha}(u_0)} \boldsymbol{\lambda}_j(u_0)^{-1}}
     {\mathbf{M}_{jj}(u_0)+2t\boldsymbol{\lambda}_j(u_0)-x}$,  
each of which describes the profile of a single dark solitary wave
after subtracting the common background
\(e^{i\boldsymbol{\alpha}(u_0)}\).  
This is a first step toward a dark soliton resolution phenomenon in
the nonzero boundary regime \eqref{NonzeBoIntr} for equation \eqref{CMdNLSfocu}.

\subsection{Strategy of the Proof}
  
\noindent
The proof is organized around three main theorems. 
Theorem~\ref{ExpForMainThm} (Explicit Formula) is proved via the chiral
Lax pair structure and the commutator estimates between the
Lax--Beurling shift semigroup and Toeplitz operators. 
Theorem~\ref{MainThmSpectChar} (Spectral Characterization) builds the
bridge between the algebraic definition of dark multi-soliton
potentials and the spectral theory of the chiral Lax operator
\(\mathbf{L}_{\mathbf{u}}\). As part of its proof, we establish the
spectral decomposition of \(\mathbf{L}_{\mathbf{u}}\) in
Proposition~\ref{propSpecDecomLu}, in the case where \(\mathbf{u}\) is a
dark multi-soliton potential. 
Theorem~\ref{ActionAngleMainThm} (Complete Integrability) combines the
two: it uses the spectral characterization to construct the
action-angle mapping \(\boldsymbol{\Phi}_N\) and prove its bijectivity,
while the time evolution of the angle variables is deduced from the
linearization identity~\eqref{GWtUtisWtUtGadd2tL}   established in the proof of
Theorem~\ref{ExpForMainThm}. 
The invariance of the dark \(N\)-soliton subset \(\mathcal U_N\)
(Proposition~\ref{PropInvcalUN}) follows from
Theorem~\ref{MainThmSpectChar} and the conservation of the generating
functional \(\mathcal H_\nu^a\) defined in \eqref{calHGeneFctdefFor}.
The inverse spectral formula \eqref{invSPEFortxIntro} in
Proposition~\ref{InvSpecFortxintro} is the special case of the explicit
flow formula \eqref{ExpForCMdNLSfoc} in Theorem~\ref{ExpForMainThm} for
dark multi-soliton solutions of equation~\eqref{CMdNLSfocu}. 
The remaining results on time evolution, long time Sobolev estimates,
and dark soliton resolution (Proposition~\ref{propLoTiBeCMsysIntro},
Corollary~\ref{corLongtimebehave}, Corollary~\ref{DarkSolRescor}) then
follow from Theorem~\ref{ActionAngleMainThm},
Proposition~\ref{InvSpecFortxintro}, 
  and the matrix perturbation theory developed in subsection~\ref{SubseMatPert6.3}, which provides the long time asymptotics of the eigenvalues of $\mathbf{M}(\mathbf{u}) + 2t \mathbf{V}(\mathbf{u}) \in \mathbb{C}^{N \times N}$, $\forall \mathbf{u} \in \mathcal{U}_N$. \\

\noindent\textbf{Theorem~\ref{ExpForMainThm} (Explicit Formula).}
The proof relies on the Lax pair structure of~\eqref{CMdNLSfocu}.
Let \(\mathbf{L}_{\mathbf{u}}\) be the chiral Lax operator defined in
\eqref{LuSect2}, 
and let \(\mathbf{B}_{\mathbf{u}}\) be the operator defined in
\eqref{BudefForFoc-}.
 By Lemma~\ref{LaxPairBuLu}, the pair
\((\mathbf{L}_{\mathbf{u}}, \mathbf{B}_{\mathbf{u}})\) is a Lax pair
for~\eqref{CMdNLSfocu}: for any solution
\(u \in \mathcal{C}((T^-,T^+);\mathcal{Z}^3_+) \cap
\mathcal{C}^1((T^-,T^+);\mathcal{Z}^1_+)\), we have  
\begin{equation}\label{LaxPairIntro}
\partial_t \mathbf{L}_{u(t)}
=
\bigl[\mathbf{B}_{u(t)}, \mathbf{L}_{u(t)}\bigr],
\qquad
T^- < t < T^+.
\end{equation}
  Since \(i\mathbf{B}_{u(t)}\) is a bounded and self-adjoint operator on the Hardy space 
\(L^2_+\), it generates a unitary propagator \(\mathbf{W} : t \mapsto \mathbf{W}(t)\) in    \eqref{UniPropagbfWSect4}, uniquely
determined by
\begin{equation}\label{bfWODEintro}
\mathbf{W}(0) = \mathrm{id}_{L^2_+},
\qquad
\partial_t \mathbf{W}(t) = \mathbf{B}_{u(t)} \mathbf{W}(t), \qquad
T^- < t < T^+.
\end{equation}
  The main difficulty in deriving the explicit flow formula
\eqref{ExpForCMdNLSfoc} lies in the appearance of the unbounded
operator \(\bigl(\mathbf{G}, \mathrm{Dom}(\mathbf{G}) \bigr)\) defined in
\eqref{DombfGdefIntro} and \eqref{bfGvphidefIntro}. One must therefore
analyze how \(\mathbf{G}\) interacts with the unitary propagator
\(\mathbf{W}(t)e^{-it\mathbf{L}_{u_0}^2}\),  
  where $u_0 = u(0) \in \mathcal{Z}^3_+$ denotes the initial datum.
   In the earlier works of G\'erard~\cite{BOExp2022},
Chen~\cite{Xchenccm25ExpFor}, and others, the linearization identity 
\begin{equation}\label{LinIdbfGWintro}
\mathbf{G}\,\mathbf{W}(t)e^{-it\mathbf{L}_{u_0}^2}
=
\mathbf{W}(t)e^{-it\mathbf{L}_{u_0}^2}
\bigl(\mathbf{G}+2t\mathbf{L}_{u_0}\bigr)
\end{equation}
was derived by differentiating the conjugate
\(\mathbf{W}(t)^{-1}\mathbf{G}\,\mathbf{W}(t)\) with respect to \(t\). 
  For this computation to be justified, one must verify that
\(\mathrm{Dom}(\mathbf{G})\) is preserved by \(\mathbf{W}(t)\) and that
the time derivative \(\partial_t\) commutes with the unbounded operator
\(\mathbf{G}\). Neither of these properties is established in those
works. Moreover, in~\cite{BOExp2022, Xchenccm25ExpFor} the natural
domains of the operators appearing on the two sides of
\eqref{LinIdbfGWintro} are not identified, and whether they coincide is
left open. \\

\noindent In the present paper, we introduce an alternative approach that avoids
these difficulties entirely by working with the contraction semigroup
\(\bigl(\mathbf{S}(\eta)^*\bigr)_{\eta\ge0}\) given by
\eqref{LaxBeurlingSG}, whose infinitesimal generator is \(-i\mathbf{G}\). 
  To be precise, we have
\begin{equation}\label{bfSstarintro}
\mathbf{S}(\eta)^* = \mathbf{T}_{e^{-i \eta x}} = e^{- i \eta \mathbf{G}} , \qquad \forall \eta \geq 0. 
\end{equation} 
We introduce the commutator limit formula
\eqref{commForSetasttildBuCv} and the estimate
\eqref{SetasttildBuUniEstL2} between \(\mathbf{S}(\eta)^*\) and the
operator \(\widetilde{\mathbf{B}}_{u(t)}\) defined in
\eqref{tildeBudefFoc-}, for any \(t \in (T^-,T^+)\) and \(\eta>0\).
  Then we obtain the limit
\begin{equation}\label{StrtSemigrLimLinid}
\lim_{\eta\to0^+}
\Bigl\|
\tfrac{i}{\eta}\Bigl(
e^{it\mathbf{L}_{u_0}^2}\mathbf{W}(t)^{-1}\mathbf{S}(\eta)^*
\mathbf{W}(t)e^{-it\mathbf{L}_{u_0}^2}-\mathbf{S}(\eta)^*
\Bigr)\varphi
-2t\mathbf{L}_{u_0}(\varphi)
\Bigr\|_{L^2_+}=0,
\qquad \forall \varphi\in H^2_+,
\end{equation}
in Lemma~\ref{LemConjSstcvto2tLuphi}. 
  Since \(\mathbf{S}(\eta)^*\) is a bounded operator on \(L^2_+\)
independent of \(t\), it can be interchanged with the difference
quotient used to define \(\partial_t\). 
  Consequently, no
domain issue arises at this stage.   
  The passage from \eqref{StrtSemigrLimLinid} back to the generator
\(\mathbf{G}\) is the second key step. Since
\(\bigl(\mathbf{S}(\eta)^*\bigr)_{\eta\ge0}\) is a contraction
semigroup, the Hille--Yosida theorem and the strong convergence
characterization \eqref{SemiGrCVdefDomsG} of \(\mathrm{Dom}(\mathbf{G})\)
allow us to identify the limit in \eqref{StrtSemigrLimLinid} with the
action of \(\mathbf{G}\) on the appropriate domain. This yields the
invariance of \(\mathrm{Dom}(\mathbf{G})\cap H^2_+\) under the unitary
propagator \(\mathbf{W}(t)e^{-it\mathbf{L}_{u_0}^2}\), that is,
\begin{equation}\label{IntroDGH2inv}
\mathbf{W}(t)e^{-it\mathbf{L}_{u_0}^2}
\bigl(\mathrm{Dom}(\mathbf{G})\cap H^2_+\bigr)
=
\mathrm{Dom}(\mathbf{G})\cap H^2_+
=
e^{it\mathbf{L}_{u_0}^2}\mathbf{W}(t)^{-1}
\bigl(\mathrm{Dom}(\mathbf{G})\cap H^2_+\bigr).
\end{equation}
Then the linearization identity \eqref{LinIdbfGWintro} holds for all
\(\varphi\in\mathrm{Dom}(\mathbf{G})\cap H^2_+\), now with its natural
domain of validity made explicit. The proof of \eqref{IntroDGH2inv} and
\eqref{LinIdbfGWintro} is given by
Theorem~\ref{invSpDomGHtwoThm}. We refer to
Remark~\ref{remcompariGerChen} for the detailed comparison between the
semigroup method in this paper and the method
in~\cite{BOExp2022, Xchenccm25ExpFor}. \\

\noindent
With Theorem~\ref{invSpDomGHtwoThm} at hand, the derivation of the
explicit formula \eqref{ExpForCMdNLSfoc} proceeds by a standard
resolvent argument: applying \(\mathbf{G}\) to the vector
\((\mathbf{G}-z)^{-1}(u(t)\chi_\epsilon)\), with \(\chi_\epsilon\) as in
\eqref{chiepFor}, and using the linearization identity
\eqref{LinIdbfGWintro} together with the inversion formula
\eqref{InvForNewtonLeibniz}, we can express the Poisson integral
\(\underline{u}(t,z)\) as the limit, as \(\epsilon\to0^+\), of
\(\mathbf{I}_+\) applied to the resolvent
\(\bigl(e^{-it\partial_x^2}\mathbf{G}e^{it\partial_x^2}
-2t\mathbf{T}_{u_0}\mathbf{T}_{\overline{u}_0}-z\bigr)^{-1}\)
acting on the product \(u_0\chi_\epsilon \in H^3_+\).
  Having established the explicit formula \eqref{ExpForCMdNLSfoc} for all
\(\mathcal{Z}^3_+\)-solutions of equation \eqref{CMdNLSfocu}, we use a
density-continuity argument to extend the validity of
\eqref{ExpForCMdNLSfoc} to all \(\mathcal{Z}^2_+\)-solutions. The proof
of Theorem~\ref{ExpForMainThm} and Corollary~\ref{CorExpforaaddH} is
completed in subsection~\ref{DerExpForSubs4.3}.
  In contrast to the approach of Killip--Laurens--Vi\c{s}an
\cite{KillLV2311}, the derivation of the explicit formulas
\eqref{ExpForCMdNLSfoc} and \eqref{expForaaddH2} in this paper does not
require any weighted hypothesis of the form
\(\langle x\rangle f\in L^2(\mathbb{R};\mathbb{C})\) or
\(\langle x\rangle f\in H^2(\mathbb{R};\mathbb{C})\), and no Gr\"onwall
estimate is needed to propagate such a condition in time. The explicit
formulas are obtained directly by using the semigroup method developed
in subsections~\ref{subsUniEquiv4.2} and~\ref{DerExpForSubs4.3}. We
refer to Remark~\ref{ComparwKLVrem} for a detailed comparison between
the semigroup method of this paper and the method in~\cite{KillLV2311}.\\

\noindent The invariant-space properties in Proposition~\ref{InvSpPropIntro}
follow directly from identity~\eqref{IntroDGH2inv} and the dominated
convergence argument in Lemma~\ref{invSpDomGCLbfI}. \\
  
\noindent\textbf{Theorem~\ref{MainThmSpectChar} (Spectral Characterization).}
  The proof proceeds in three steps: 
   $(\mathbf{a}) \Rightarrow (\mathbf{c})$, $(\mathbf{b}) \Rightarrow (\mathbf{a})$, and $(\mathbf{c}) \Rightarrow (\mathbf{b})$.
  First, we prove Proposition~\ref{propSpecDecomLu} on the spectral
decomposition of the self-adjoint operator
\(\mathbf{L}_{\mathbf{u}}\) for any \(\mathbf{u}\in\mathcal{U}_N\). This
proposition is the key ingredient for  
\((\mathbf{a})\Rightarrow(\mathbf{c})\).
 By
Remark~\ref{RmkDarNsolsumCstH1}, every
\(\mathbf{u}\in\mathcal{U}_N\) admits a decomposition
\(\mathbf{u}=\mathbf{u}(\infty)+\mathbf{v}\) with
\(\mathbf{u}(\infty)\in\mathbb{S}^1 = \{\zeta \in \mathbb{C} : |\zeta| =1\}\) and \(\mathbf{v}\in H^1_+\). 
  The main difficulty in the proof of Proposition~\ref{propSpecDecomLu} is to exclude the possibility that $-1$ is an eigenvalue of $\mathbf{L}_{\mathbf{u}}$. 
    The auxiliary identity
\eqref{idenhL2nvhinprod}, which is the key to proving simplicity,
requires the eigenvector \(h\) to belong to
\(\mathrm{Dom}(\mathbf{I}_+)\) as defined in \eqref{DomI+}. 
To be precise, for any $\alpha \in \mathbb{T} = \mathbb{R} \slash 2\pi \mathbb{Z}$,  $\lambda \in \mathbb{R}$ and $\mathtt{u} \in e^{i\alpha} + H^1_+$,  if $h \in \mathrm{Ker} \bigl(\mathbf{L}_{\mathtt{u}} - \lambda \bigr) \cap \mathrm{Dom} (\mathbf{I}_+)$, 
 then the following identity holds: 
 \begin{equation}\label{IntroidenhL2nvhi}
2 \pi \| h \|_{L^2_+}^2 = e^{- i\alpha}   \mathbf{I}_+ (h)  \langle \mathtt{u} - e^{i\alpha} , h \rangle_{L^2_+} + |\langle \mathtt{u} - e^{i\alpha} , h \rangle_{L^2_+}|^2.
 \end{equation} 
  For any 
\(\lambda\in\sigma_{\mathrm{pp}}(\mathbf{L}_{\mathtt{u}})\setminus\{-1\}\),
Lemma~\ref{EigLuhperpvLem} shows that every \(\lambda\)-eigenvector
automatically lies in \(\mathrm{Dom}(\mathbf{I}_+)\), and   identity \eqref{IntroidenhL2nvhi}
then implies that the eigenspace  $\mathrm{Ker}\bigl( \mathbf{L}_{\mathtt{u}} - \lambda \bigr)$ is one-dimensional by Corollary~\ref{CorSimpofEigLu}. For
\(\lambda=-1\), formula \eqref{KerLuadd1DI0} 
 in the same lemma yields that
\(\mathrm{Ker}(\mathbf{L}_{\mathtt{u}}+1)\cap\mathrm{Dom}(\mathbf{I}_+)=\{0\}\).
  Consequently,  for any \(\mathbf{u}\in\mathcal{U}_N\),   once we prove that every eigenvector of
\(\mathbf{L}_{\mathbf{u}}\) belongs to \(\mathrm{Dom}(\mathbf{I}_+)\), then 
  \(-1 \notin \sigma_{\mathrm{pp}} (\mathbf{L}_{\mathbf{u}})\)   and all eigenvalues of $\mathbf{L}_{\mathbf{u}}$ are simple. This
is achieved by   identifying the spectral subspaces
\(\mathscr{H}_{\mathrm{ac}}(\mathbf{L}_{\mathbf{u}})\),
\(\mathscr{H}_{\mathrm{sc}}(\mathbf{L}_{\mathbf{u}})\), and
\(\mathscr{H}_{\mathrm{pp}}(\mathbf{L}_{\mathbf{u}})\) using the Hardy
space tools developed in subsections~\ref{SubsecHarSp2.3}
and~\ref{SubsShifSemiGr} together with the functional calculus. The
inclusion \(\mathscr{H}_{\mathrm{pp}}(\mathbf{L}_{\mathbf{u}})\subset
H^{\infty}_+\cap\mathrm{Dom}(\mathbf{G}) \subset H^{\infty}_+ \cap\mathrm{Dom}(\mathbf{I}_+)\), obtained
from~\eqref{SpecDeHacHscHpp}, \eqref{DomrGsubsetDomsG} and \eqref{DXsubsetDG}, then shows
that for \(\mathbf{u}\in\mathcal{U}_N\) the operator
\(\mathbf{L}_{\mathbf{u}}\) has only finitely many eigenvalues and every
eigenvalue is simple. This completes the proof of Proposition~\ref{propSpecDecomLu}. 
By  Proposition~\ref{KthetapropOrDecom}, the implications \((\mathbf{a})\Rightarrow(\mathbf{b})\) and
\((\mathbf{a})\Rightarrow(\mathbf{c})\) in Theorem~\ref{MainThmSpectChar}  are both valid.  \\
 
\noindent The implication \((\mathbf{b})\Rightarrow(\mathbf{a})\) relies on the
Lax--Beurling theorem characterizing all closed translation-invariant
subspaces of the Hardy space \(L^2_+\)
  (see Lax~\cite{LaxTransinva1959} for the original Lax--Beurling
theorem, and Lemma~3.3 of~\cite{SUNCMP2021} together with
Remark~\ref{RemLaxTrInvSp} for the version used here)
  and on
Lemma~\ref{lemfhinHarSobfL2+}, which characterizes the Hardy space
\(L^2_+\). The proof of Corollary~\ref{CorQuPu} is also completed in this
step. \\

\noindent For the implication \((\mathbf{c})\Rightarrow(\mathbf{b})\), one must
show that \(\mathbf{u}-\mathbf{u}(\infty)\in H^1_+\) is a cyclic vector
of \(\mathbf{L}_{\mathbf{u}}|_{\mathscr{H}_{\mathrm{pp}}(\mathbf{L}_{\mathbf{u}})}\).
This is where the auxiliary identity \eqref{IntroidenhL2nvhi} is used,
and therefore the assertion \((\mathbf{c})\) includes the assumption
that \(-1\) is not an eigenvalue of \(\mathbf{L}_{\mathbf{u}}\). The
implication \((\mathbf{b})\Rightarrow(\mathbf{a})\) does not rely on
\eqref{IntroidenhL2nvhi}, which is why no such assumption is needed in
\((\mathbf{b})\). The proof of Theorem~\ref{MainThmSpectChar} is
completed in subsection~\ref{SpecCharSubsect}.  

\begin{rem}\label{RemBOWuiden}
\noindent The idea of introducing the auxiliary identity \eqref{IntroidenhL2nvhi}
to establish the simplicity of eigenvalues of the chiral Lax operator
\(\mathbf{L}_{\mathtt{u}}\) is inspired by Wu~\cite{WuBO2016}. In
Lemma~2.5 on p.~1353 of~\cite{WuBO2016}, Y.~Wu introduced, under the
additional condition \(\lambda<0\), an analogous auxiliary identity for
the Benjamin--Ono (BO) equation
\begin{equation}\label{BOeq}
\partial_t q = \mathbb{H} \partial_x^2 q - \partial_x (q^2),
\qquad
(t,x) \in \mathbb{R}\times \mathbb{R},
\end{equation}
and deduced from it that every negative eigenvalue of the associated
Lax operator
\(L^{\mathrm{BO}}_q := -i\partial_x - \mathbf{T}_q : H^1_+ \to L^2_+\)
is simple. 
 In Proposition~2.4 of~\cite{SUNCMP2021}, R.~Sun
  removed the condition \(\lambda<0\) from Lemma~2.5
of~\cite{WuBO2016} 
  and showed that the auxiliary identity for the BO equation \eqref{BOeq} holds for
all eigenvalues
\(\lambda\) of \(L^{\mathrm{BO}}_q\), thereby proving that every
eigenvalue of \(L^{\mathrm{BO}}_q\) is simple. We emphasize that, in the
BO setting \eqref{BOeq}, no condition analogous to \(-1\notin
\sigma_{\mathrm{pp}}(\mathbf{L}_{\mathtt{u}})\) is required. In the
focusing CMdNLS case \eqref{CMdNLSfocu} under the nonzero boundary
condition~\eqref{NonzeBoIntr}, however, the auxiliary identity
\eqref{IntroidenhL2nvhi} still requires the assumption \(-1\notin
\sigma_{\mathrm{pp}}(\mathbf{L}_{\mathtt{u}})\).  
\end{rem}

\noindent Let $u \in a + \mathcal{C} \bigl( (T^-, T^+); H^3_+ \bigr) \cap \mathcal{C}^1 \bigl( (T^-, T^+); H^1_+ \bigr)  $ be a solution to equation \eqref{CMdNLSfocu} for some  $a \in \mathbb{C}$. 
  By formula \eqref{CauTrSpeMeaLuisH} and
Lemma~\ref{LemFirsIntGenFct}, the Cauchy transform of the spectral
measure \(\mu^{\mathbf{L}_{u(t)}}_{u(t) - a}\) is conserved along the flow
of~\eqref{CMdNLSfocu}. Together with
Theorem~\ref{MainThmSpectChar} and the linearization identity
\eqref{LinIdbfGWintro}, this yields the invariance of the dark
multi-soliton subset \(\mathcal U_N\) stated in the first assertion of
Proposition~\ref{PropInvcalUN}. The second assertion of the same
proposition, namely the global existence of pure dark multi-soliton
solutions, is deferred to subsection~\ref{SusectDMulSol6.1}, after   formula \eqref{ilapu2uPidmodu2Foc} has been established. \\
 
\noindent\textbf{Theorem~\ref{ActionAngleMainThm} (Complete Integrability).}  The proof is divided into three steps. First, we construct the action-angle mapping   $\boldsymbol{\Phi}_N\colon\mathcal U_N\to
\mathbb{T}\times \Omega_N\times\mathbb{R}^{1\times N}$  in Definition~\ref{InvSPmADef} of  subsection~\ref{AcAngMapSubs5.2}. 
   According to Theorem~\ref{MainThmSpectChar} and Proposition~\ref{propSpecDecomLu}, for any $\mathbf{u} \in \mathcal{U}_N$, the Lax operator $\mathbf{L}_{\mathbf{u}}$ has exactly $N$ distinct eigenvalues and each eigenvalue is simple. The action variables $ 
            \boldsymbol{\lambda}_1 (\mathbf{u}) < \boldsymbol{\lambda}_2 (\mathbf{u}) < \cdots < \boldsymbol{\lambda}_N (\mathbf{u})$ are defined as the eigenvalues of $\mathbf{L}_{\mathbf{u}}$. Since $\dim_{\mathbb{C}} \mathscr{H}_{\mathrm{pp}}(\mathbf{L}_{\mathbf{u}}) = N$            
            by Theorem~\ref{MainThmSpectChar},
             it admits an orthonormal eigenbasis
\(\eth^{\mathbf{u}} := \{\phi_1^{\mathbf{u}}, \phi_2^{\mathbf{u}},\cdots,\phi_N^{\mathbf{u}}\}\)
chosen so that 
\begin{equation}\label{Introphiju}
\phi_j^{\mathbf{u}}\in \mathrm{Ker}(\boldsymbol{\lambda}_j  (\mathbf{u}) -\mathbf{L}_{\mathbf{u}}), 
\quad 
 \|\phi_j^{\mathbf{u}}\|_{L^2_+} = 1
\quad\text{and}\quad
 \mathbf{p}_j( \mathbf{u} ):=  \langle \phi_j^{\mathbf{u}}, \mathbf{u} - \mathbf{u}(\infty) \rangle_{L^2_+} > 0, 
\end{equation}for any $j \in \{ 1, 2, \cdots, N \}$. Then the inverse spectral matrix $\mathbf{M}(\mathbf{u}) \in \mathbb{C}^{N \times N} $ in 
 Definition~\ref{bfMVfrakXYdefintro} is defined to be the matrix of the restriction $\mathbf{G}|_{\mathscr{H}_{\mathrm{pp}}(\mathbf{L}_{\mathbf{u}}) } : \mathscr{H}_{\mathrm{pp}}(\mathbf{L}_{\mathbf{u}})  \to \mathscr{H}_{\mathrm{pp}}(\mathbf{L}_{\mathbf{u}})$ with respect to $\eth^{\mathbf{u}}$.
  The total phase $\boldsymbol{\alpha}(\mathbf{u}) \in \mathbb{T}$ and the generalized angle variables $\vec{\boldsymbol{\gamma}}(\mathbf{u}) =(\boldsymbol{\gamma}_1 (\mathbf{u}), \boldsymbol{\gamma}_2 (\mathbf{u}) , \cdots,
      \boldsymbol{\gamma}_N (\mathbf{u})) \in \mathbb{R}^{1\times N}$ are given by
\begin{equation}\label{bdlamaldefIntr}
e^{i\boldsymbol{\alpha}(\mathbf{u})}  =  \lim_{|x|\to+\infty} \mathbf{u}(x)     = \mathbf{u} (\infty) ,  \qquad \boldsymbol{\gamma}_j( \mathbf{u} ):=   \operatorname{Re}\langle \mathbf{G}\phi_j^{\mathbf{u}}, \phi_j^{\mathbf{u}} \rangle_{L^2_+},  
\end{equation}for any $j \in \{1,2, \cdots, N\}$.  
 Lemma~\ref{MentlamgamFor}
enables us 
to express every entry of $\mathbf{M}(\mathbf{u}) \in \mathbb{C}^{N \times N}$ in terms of the action-angle variables $\{\boldsymbol{\lambda}_j ( \mathbf{u} ) ; \boldsymbol{\gamma}_j( \mathbf{u} )\}_{1 \leq j \leq N}$.   
By Lemma~\ref{LemSpLutwoseg}, the parameter
\(\mathbf{p}_j(\mathbf{u})>0\) in \eqref{Introphiju} can also be
expressed in terms of the action variable
\(\boldsymbol{\lambda}_j(\mathbf{u}) \in \mathbb{R}\). This yields
\(\sigma_{\mathrm{pp}}(\mathbf{L}_{\mathbf{u}}) \cap [-1,0] =
\emptyset\) for every \(\mathbf{u}\in\mathcal{U}_N\).
    This motivates the choice of \(\Omega_N\) in \eqref{OmegaNIntro} as
the set of actions. 
      In Proposition~\ref{PropInvSpecFormula}, we use the inversion formula  \eqref{InvForNewtonLeibniz}   to deduce the   inverse spectral formula   \eqref{invSPEFormain} for $\mathbf{u} \in \mathcal{U}_N$, which implies the injectivity of  $\boldsymbol{\Phi}_N\colon\mathcal U_N\to
\mathbb{T}\times \Omega_N\times\mathbb{R}^{1\times N}$. 
The phase rotation formula \eqref{ActAnginvbyPhRot} follows directly by Lemma~\ref{LemPhaseRotInv}.  \\
 
\noindent Second, we establish the surjectivity of the action-angle mapping $\boldsymbol{\Phi}_N\colon\mathcal U_N\to
\mathbb{T}\times \Omega_N\times\mathbb{R}^{1\times N}$ in subsection~\ref{InSpTrsubsect5.3}. 
  Fix $(\vec{\lambda} , \vec{\gamma}  ) \in \Omega_N \times \mathbb{R}^{1 \times N}$. Let $ M=M(\vec{\lambda}, \vec{\gamma}) \in \mathbb{C}^{N \times N}$ be defined  by \eqref{defMlamgamFoc}. The vectors  
$X=X_{\vec{\lambda}}  \in \mathbb{R}^N$ and $ Y=Y_{\vec{\lambda}} \in \mathbb{R}^N$ are given by \eqref{defXYFoclambda}.  
   Lemma~\ref{SpofMLemfoc} implies that $\sigma (M) \subset \mathbb{C}_-$. For any $\alpha \in \mathbb{T}$, 
   let the potential $\mathsf{u} = \mathsf{u}_{\alpha, \vec{\lambda}, \vec{\gamma}} : \mathbb{R} \to \mathbb{C}$ be defined by   \eqref{ulamgamalpdef}, i.e.
   \begin{equation}\label{ttuallamgamintro}
    \mathsf{u} (x)  
    = e^{i\alpha} \Bigl(1 + i Y^{\mathrm{T}}   (M - x)^{-1} X   \Bigr)
    , \qquad \forall x \in \mathbb{R} .
   \end{equation}
   In Lemma~\ref{ulagaFoclinUNsurj}, we prove that  the potential $\mathsf{u} = \mathsf{u}_{\alpha, \vec{\lambda}, \vec{\gamma}}$ given by  \eqref{ttuallamgamintro}  belongs to $\mathcal{U}_N$. It follows directly from \eqref{ttuallamgamintro} that $\boldsymbol{\alpha}(\mathsf{u}) = \alpha \in \mathbb{T}$.
     Lemma~\ref{fjEigFuncofLuLem} shows that $\vec{\boldsymbol{\lambda}}( \mathsf{u} ) = \vec{\lambda} \in \Omega_N$.     
  We have $ \vec{\boldsymbol{\gamma}}( \mathsf{u} ) = \vec{\gamma} \in \mathbb{R}^{1 \times N}$ by Lemma~\ref{bfGfjfjLem}. So $\boldsymbol{\Phi}_N\colon\mathcal U_N\to
\mathbb{T}\times \Omega_N\times\mathbb{R}^{1\times N}$ is surjective. 
   Moreover, the eigenbasis $\eth^{\mathsf{u}}$ given by \eqref{Introphiju} can be explicitly expressed in terms of the parameters $\bigl(\alpha, \vec{\lambda} , \vec{\gamma} \bigr) \in \mathbb{T} \times \Omega_N \times \mathbb{R}^{1 \times N}$ in \eqref{FocexpExpphij}. \\ 
 
\noindent Finally, the time evolution formula \eqref{TimEvoactangl} follows from
the linearization identity \eqref{LinIdbfGWintro}, which yields the
evolution of the eigenfunctions in Lemma~\ref{LemEigVecTimevo}.  
 Plugging \eqref{TimEvoactangl} into \eqref{invSPEFormain}, we obtain 
the inverse spectral formula \eqref{invSPEFortxIntro}. 
The proofs of Theorem~\ref{ActionAngleMainThm},
Proposition~\ref{PropInvcalUN}, and
Proposition~\ref{InvSpecFortxintro} are completed in
subsection~\ref{SusectDMulSol6.1}.
  The remaining results on time
evolution, long time Sobolev estimates, and dark soliton resolution
(Proposition~\ref{propLoTiBeCMsysIntro}, Corollary~\ref{corLongtimebehave},
Corollary~\ref{DarkSolRescor}) are proved in
subsections~\ref{6.4subsResPody} and~\ref{subsLonTiBe6.5}, using the matrix perturbation theory
developed in subsection~\ref{SubseMatPert6.3}.

\subsection{Other Related Works}

\noindent We first compare the focusing and
defocusing CMdNLS equations under the nonzero boundary condition \eqref{NonzeBoIntr}, and
then turn to other research directions: blow-up and soliton resolution, long time behavior and scattering theory, Hamiltonian structure, the
space-periodic setting, and the direct scattering transform.\\

\noindent
The inverse spectral theory developed in
Section~\ref{DarMulSoliPotSect} and the long time analysis of
Section~\ref{SectTimeEvo6} admit counterparts for the defocusing
CMdNLS equation
\begin{equation}\label{defocuCMdNLS}
i\partial_t u+\partial_x^2u
=
u(\mathbb{H}-i)\partial_x(|u|^2),
\qquad
(t,x)\in\mathbb R\times\mathbb R,
\end{equation}
under the same nonzero boundary condition~\eqref{NonzeBoIntr}. In the
companion paper Sun~\cite{SUNCMdNLSdefo}, a spectral characterization
theorem for the defocusing dark \(N\)-soliton set, denoted by
\(\mathcal U_N^+\), is established, and this set is shown to be
invariant under the flow of equation~\eqref{defocuCMdNLS}. The
action-angle variables of equation~\eqref{defocuCMdNLS} restricted to
each \(\mathcal U_N^+\) are then constructed, together with the
corresponding inverse spectral formula. 
Kato's perturbation theory
further yields the long time Sobolev estimates for pure dark multi-soliton
solutions of~\eqref{defocuCMdNLS}. 
  Moreover, the pole dynamics of every such
solution is shown to be governed by the same complexified
Calogero--Moser system~\eqref{IntroCaloMoserSys}.\\

\noindent However, the focusing and defocusing cases exhibit a fundamental
difference, both in the location of  point spectrum of the chiral
Lax operator and in the pointwise size of   dark multi-soliton
potentials. To distinguish the two, we write
\(\mathcal U_N^-:=\mathcal U_N\) for the dark \(N\)-soliton set of the
focusing CMdNLS equation~\eqref{CMdNLSfocu} defined
in~\ref{DefiDarNsoli}.

\begin{enumerate}
\item For \(\mathbf{u}\in\mathcal U_N^-\) a dark multi-soliton
potential of the focusing equation~\eqref{CMdNLSfocu}, 
  eigenvalues of the chiral Lax operator \(\mathbf{L}_{\mathbf{u}}^-:=\mathbf{L}_{\mathbf{u}}
=-i\partial_x
-\mathbf{T}_{\mathbf{u}}\mathbf{T}_{\overline{\mathbf{u}}}
:H^1_+\to L^2_+\) can be chosen arbitrarily outside the compact
interval \([-1,0]\), i.e.
\(\sigma_{\mathrm{pp}}(\mathbf{L}_{\mathbf{u}}^-)\subset
\mathbb R\setminus[-1,0]\). In particular, embedded eigenvalues may
occur. Moreover, the potential \(\mathbf{u}\in\mathcal U_N^-\)
satisfies
\begin{equation}\label{amplibfubig1}
|\mathbf{u}(x)|>1,
\qquad
\forall x\in\mathbb R .
\end{equation}

\item For \(\mathfrak{u}\in\mathcal U_N^+\) a dark multi-soliton
potential of the defocusing equation~\eqref{defocuCMdNLS}, all
eigenvalues of the chiral Lax operator
\(\mathbf{L}_{\mathfrak{u}}^+:=-i\partial_x
+\mathbf{T}_{\mathfrak{u}}\mathbf{T}_{\overline{\mathfrak{u}}}
:H^1_+\to L^2_+\) are confined to the bounded open interval
\((0,1)\), i.e.
\(\sigma_{\mathrm{pp}}(\mathbf{L}_{\mathfrak{u}}^+)\subset(0,1)\).
Weyl's essential spectrum theorem gives
\(\sigma_{\mathrm{ess}}(\mathbf{L}_{\mathfrak{u}}^+)
=\sigma_{\mathrm{ess}}(D+1)=[1,+\infty)\), so every eigenvalue
of \(\mathbf{L}_{\mathfrak{u}}^+\) is discrete. Moreover, the potential
\(\mathfrak{u}\in\mathcal U_N^+\) satisfies
\begin{equation}\label{amplifrakusmal1}
|\mathfrak{u}(x)|<1,
\qquad
\forall x\in\mathbb R .
\end{equation}
\end{enumerate}
 
\noindent
Physically, equation~\eqref{defocuCMdNLS} arises as the deep-water
limit of the intermediate nonlinear Schr\"odinger equation introduced
by Pelinovsky~\cite{Pelinovsky1995} and
Pelinovsky--Grimshaw~\cite{Pel--Gri1995}, which describes the
propagation of nonlinear wave packets in a two-layer fluid. An
explicit formula for general solutions of~\eqref{defocuCMdNLS} under
the nonzero boundary condition~\eqref{NonzeBoIntr} was established in
Theorem~1.2 of Chen~\cite{Xchenccm25ExpFor}; see
Remark~\ref{remcompariGerChen} for a comparison with the semigroup
method used in this paper. \\
 
\noindent
It is a classical fact that the cubic focusing nonlinear Schr\"odinger
equation
\begin{equation}\label{focNLSonR}
i\partial_t u+\partial_x^2u=-|u|^2u,
\qquad
(t,x)\in\mathbb R\times\mathbb R,
\end{equation}
does not admit dark solitons, nor dark multi-solitons, as exact
solutions on the real line; see Zakharov--Shabat~\cite{ZS1972} and
Subsection~9.3 on p.~697 of
Dubrovin--Grava--Klein--Moro~\cite{DubroGraKlMo15}. The ground state
of~\eqref{focNLSonR} can be computed explicitly by solving an ordinary
differential equation.
  In contrast to~\eqref{focNLSonR}, the focusing CMdNLS
equation~\eqref{CMdNLSfocu} admits dark solitons and dark
multi-solitons. Two structural features are responsible for this
contrast. First, the Szeg\H{o} projector
\(\Pi:L^2(\mathbb R;\mathbb C)\to L^2(\mathbb R;\mathbb C)\) given
by~\eqref{PiSzeProIntro} restricts the dynamics to the Hardy space
\(L^2_+\), whose natural singular functions are Cauchy kernels. The
ground state of~\eqref{CMdNLSfocu} under the nonzero boundary
condition~\eqref{NonzeBoIntr} is of this form. Since \(\Pi\) is a
nonlocal operator, the equation satisfied by the ground state is a
pseudo-differential equation rather than an ordinary differential
equation. Second, the point spectrum of the chiral Lax operator
\(\mathbf{L}_{\mathbf{u}}\) is confined to the complement of the
compact interval \([-1,0]\), and this spectral gap guarantees that the
quantities \(\lambda_j(1+\lambda_j)\) and \((1+\lambda_j)/\lambda_j\)
appearing in the inverse spectral
formula~\eqref{ttuallamgamintro} remain strictly positive.
Consequently, the inverse spectral formula is well defined and dark
multi-soliton solutions of~\eqref{CMdNLSfocu} can be constructed from action-angle variables. \\
 
\noindent
We recall some other results on the qualitative properties of equation~\eqref{CMdNLSfocu} under the
zero boundary condition~\eqref{zeroBounIntr}. 
   In Kim--Kim--Kwon~\cite{KimKimKwonBlowup2024}, finite-time blow-up
solutions of the focusing CMdNLS equation~\eqref{CMdNLSfocu} were
constructed by means of modulation analysis. The quantized blow-up
dynamics of equation~\eqref{CMdNLSfocu} was studied in
Jeong--Kim~\cite{JeongKimT2024}. In~\cite{SoliReso2024KimKwon},
T.~Kim and S.~Kwon   established the soliton resolution for
equation~\eqref{CMdNLSfocu}. Single-bubble blow-up solutions
for~\eqref{CMdNLSfocu} were classified in
Jeong--Kim--Kim--Kwon~\cite{Jeong2026KimKimKwon}.\\
 
\noindent
In~\cite{HogKowal24Tur}, J.~Hogan and M.~Kowalski studied the growth
of Sobolev norms and the restricted orbital stability of the bright
soliton for equation~\eqref{CMdNLSfocu}. Using the explicit
formula~\eqref{expForaaddH2}, R.~Badreddine investigated the
zero-dispersion limit of both \eqref{CMdNLSfocu}
and~\eqref{defocuCMdNLS} in~\cite{Badred24ZerDisp}. The nonlinear
\(L^2\)-scattering of equations~\eqref{CMdNLSfocu}
and~\eqref{defocuCMdNLS} was established in the works of
Chen~\cite{ScaDefChen25}, Hadama~\cite{Hadama26INLS}, and
Killip--Laurens--Murphy--Vi\c{s}an~\cite{KLMV2609DisDecay}.
In~\cite{KMV2604Hamil}, R.~Killip, K.~Marsden, and M.~Vi\c{s}an
analyzed the symplectic form introduced on p.~4011 of
G\'erard--Lenzmann~\cite{GerLenz22CPAM} and established the
Hamiltonian formulation for  
equations~\eqref{CMdNLSfocu} and~\eqref{defocuCMdNLS} in each of the
following two settings: the zero boundary
regime~\eqref{zeroBounIntr}, and the space-periodic regime defined by
\begin{equation}\label{SpacPerioSett}
u(t,x+2\pi)=u(t,x), \qquad \forall t\in(T^-,T^+), \qquad \forall x \in \mathbb{R} .
\end{equation}
   Equations~\eqref{CMdNLSfocu} and~\eqref{defocuCMdNLS} under the
space-periodic condition~\eqref{SpacPerioSett} are called the
Calogero--Sutherland derivative nonlinear Schr\"odinger (CSdNLS)
equations on the \(1\)-dimensional torus
\(\mathbb{T}=\mathbb{R}/2\pi\mathbb{Z}\):
\begin{equation}\label{torusCSdNLSpm}
\partial_t u = i\partial_x^2 u + 2\kappa u
\Pi_{\geq0}^{\mathbb{T}}\partial_x(|u|^2),
\qquad
(t,x)\in\mathbb{R}\times\mathbb{T},
\end{equation}
where \(\kappa=+1\) corresponds to the focusing case and
\(\kappa=-1\) to the defocusing case, and
\(\Pi_{\geq0}^{\mathbb{T}}:L^2(\mathbb{T};\mathbb{C})\to
L^2(\mathbb{T};\mathbb{C})\) denotes the Szeg\H{o} projector on the
torus:
\begin{equation}\label{SzegoProTorus}
\Pi_{\geq 0}^{\mathbb{T}}\Bigl(\sum_{n\in\mathbb{Z}}f_n e^{inx}\Bigr)
:=
\sum_{n\geq 0}f_n e^{inx},
\qquad
\forall f=\sum_{n \in \mathbb{Z}}f_n  e^{inx}\in L^2(\mathbb{T};\mathbb{C}).
\end{equation}
  In~\cite{Badreddine2023}, R.~Badreddine constructed an explicit
formula for general solutions of~\eqref{torusCSdNLSpm} and extended the
flow map of~\eqref{torusCSdNLSpm} to
\(L^2_+(\mathbb{T};\mathbb{C}):=
\Pi_{\geq0}^{\mathbb{T}}\bigl(L^2(\mathbb{T};\mathbb{C})\bigr)\).
Traveling waves were completely classified by
Badreddine~\cite{BadreddineTr}.
  Badreddine~\cite{Badred26BirkMap}  established Birkhoff coordinates for
the defocusing CSdNLS equation~\eqref{torusCSdNLSpm} with \(\kappa = -1\)
and thereby proved its complete integrability.  
  For \(\kappa=1\), X.~Chen and E.~Lenzmann constructed finite-time
blow-up solutions for the focusing CSdNLS
equation~\eqref{torusCSdNLSpm} in~\cite{ChenLenz2026}. \\
  
\noindent Returning to the focusing CMdNLS equation~\eqref{CMdNLSfocu} on the real
line, multiphase solutions were constructed by
Matsuno~\cite{MatSAM2023} using the Hirota bilinear method. 
  In~\cite{ChenPelNonli2025}, J.~Chen and D.~E.~Pelinovsky investigated
traveling periodic waves and breathers for both \eqref{CMdNLSfocu}
and~\eqref{defocuCMdNLS} in each of the two settings: the zero boundary
regime~\eqref{zeroBounIntr} and the nonzero boundary
regime~\eqref{NonzeBoIntr}.     
  Finally, the Jost solutions and the direct scattering transform of the chiral
Lax operator \(\mathbf{L}_u\) given by~\eqref{chiralLaxOpIntr} were
investigated in the independent works of
Frank--Read~\cite{FranRead25CCM} and
Sun--Wang--Zhao~\cite{SWZ26DirScaCMdNLSfoc}.

\subsection{Organization of the Paper}
\noindent
\textbf{Section~\ref{PrelimSect2}} develops the functional-analytic
framework used throughout the paper:  
Hardy spaces on the upper half-plane, Sobolev spaces, Zhidkov spaces,
as well as  the corresponding Hardy--Sobolev spaces \(H^s_+\) and
Hardy--Zhidkov spaces \(\mathcal Z^n_+\).
 It also introduces Toeplitz
operators, the shift semigroup
\((\mathbf S(\eta)^*)_{\eta\ge0}\) and its infinitesimal generator
\((-i\mathbf G,\mathrm{Dom}(\mathbf G))\), Poisson integrals, and the
difference quotient operator \(\mathscr Q_z\). The characterization of
\(\mathrm{Dom}(\mathbf G)\) via this semigroup and the construction of
the unitary propagator underpin the proofs of the explicit formula,
the spectral analysis, and the inverse spectral transform in the
subsequent sections. \\

\noindent
\textbf{Section~\ref{SectLaxPair3}} introduces the Lax pair structure
for the focusing CMdNLS equation \eqref{CMdNLSfocu}. 
The two chiral formulations \((\mathbf{L}_u,\mathbf{B}_u)\) and
\((\mathbf{L}_u,\widetilde{\mathbf{B}}_u)\) on \(H^s_+\) play a
central role in the theory of dark multi-solitons below.  
  A family of
commutator estimates between the shift semigroup
\((\mathbf{S}(\eta)^*)_{\eta\ge 0}\) and the Toeplitz operators, as
well as between \((\mathbf{S}(\eta)^*)_{\eta\ge 0}\) and the operator
\(\widetilde{\mathbf{B}}_u\), is established. These estimates provide
the rigorous foundation for the linearization identity
\eqref{LinIdbfGWintro} and the invariant subspace theorem
(Theorem~\ref{invSpDomGHtwoThm}). The essential spectrum of
\(\mathbf{L}_u\) and the simplicity of its eigenvalues are also
determined. \\

\noindent  \textbf{Section~\ref{LWPExpForSect4}} is dedicated to proving the
explicit formulas in Theorem~\ref{ExpForMainThm} and
Corollary~\ref{CorExpforaaddH},
 after establishing the local well-posedness of
equation~\eqref{CMdNLSfocu} in the Zhidkov space \(\mathcal Z^n\) and
its Hardy--Zhidkov counterpart \(\mathcal Z^n_+\), for \(n\in\{2,3\}\).
 The key ingredient is the
unitary equivalence argument in subsection~\ref{subsUniEquiv4.2} and
the invariant subspace theorem (Theorem~\ref{invSpDomGHtwoThm}). Both
are proved using only the commutator estimates for
\((\mathbf{S}(\eta)^*)_{\eta\ge 0}\) in
subsection~\ref{SsSemiCommEsti3.4}, without differentiating the
unbounded generator \((-i\mathbf G,\mathrm{Dom}(\mathbf G))\).
 The linearization step, namely the limit formula
\eqref{LimforconjSeta-SetaCV} (see also formula~\eqref{StrtSemigrLimLinid}), avoids the domain issues that arise in
the formal differentiation of \(\mathbf{G}\).
  We compare the semigroup method of
this paper with the earlier works
\cite{BOExp2022, Xchenccm25ExpFor, KillLV2311} in
Remarks~\ref{remcompariGerChen} and~\ref{ComparwKLVrem}. In
subsection~\ref{InvSpCoLsubs4.4}, we establish   conservation laws and
invariant subspaces for equation~\eqref{CMdNLSfocu}, thereby completing
the proof of Proposition~\ref{InvSpPropIntro}.   \\
 
\noindent
\textbf{Section~\ref{DarMulSoliPotSect}} 
 studies the dark \(N\)-soliton set \(\mathcal U_N\).  
 For any \(\mathbf{u}\in\mathcal U_N\), we establish the direct and
inverse spectral theory of the chiral Lax operator
\(\mathbf{L}_{\mathbf{u}}\).
    We first apply the spectral analysis from 
subsection~\ref{SubSspanal3.5} together with the functional calculus
to determine the spectral subspaces of the self-adjoint operator
\(\mathbf{L}_{\mathbf{u}}\). Using Proposition~\ref{propSpecDecomLu},
we complete the proof of Theorem~\ref{MainThmSpectChar} in
subsection~\ref{SpecCharSubsect} with the Hardy space tools from
subsections~\ref{SubsecHarSp2.3} and~\ref{SubsShifSemiGr}. Combining
Proposition~\ref{propSpecDecomLu} and
Theorem~\ref{MainThmSpectChar}, we construct the action-angle mapping
\(\boldsymbol{\Phi}_N\colon\mathcal U_N\to
\mathbb T\times\Omega_N\times\mathbb R^{1\times N}\) and prove its
injectivity in subsection~\ref{AcAngMapSubs5.2}.
Subsection~\ref{InSpTrsubsect5.3}  proves the surjectivity of \(\boldsymbol{\Phi}_N\) by establishing the inverse spectral
transform. 
  Section~\ref{DarMulSoliPotSect} does not rely on the results of
Section~\ref{LWPExpForSect4}. \\
  
\noindent
\textbf{Section~\ref{SectTimeEvo6}} studies the dynamics and long time
behavior of dark multi-soliton solutions of
equation~\eqref{CMdNLSfocu}. In
subsection~\ref{SusectDMulSol6.1}, we complete the proofs of
Theorem~\ref{ActionAngleMainThm},
Proposition~\ref{PropInvcalUN}, and
Proposition~\ref{InvSpecFortxintro}.
Subsection~\ref{ssDar2SolSol6.2} treats the case of dark
$2$-soliton solutions.
 In subsection~\ref{SubseMatPert6.3}, we analyze
the long time asymptotics of   eigenvalues of the inverse spectral
matrix
\(\mathbf{M}(u(t))=\mathbf{M}(u_0)+2t\mathbf{V}(u_0)
\in\mathbb{C}^{N\times N}\)
by means of Kato's classical perturbation theory and the
Rayleigh--Schr\"odinger series. 
In
subsection~\ref{6.4subsResPody}, we show that the pole dynamics of all
dark \(N\)-soliton solutions of equation~\eqref{CMdNLSfocu} is governed by the
complexified Calogero--Moser system \eqref{IntroCaloMoserSys}, thereby completing the proof of
Proposition~\ref{propLoTiBeCMsysIntro}. Finally,
subsection~\ref{subsLonTiBe6.5} proves
Corollary~\ref{corLongtimebehave} and
Corollary~\ref{DarkSolRescor}. 
  
 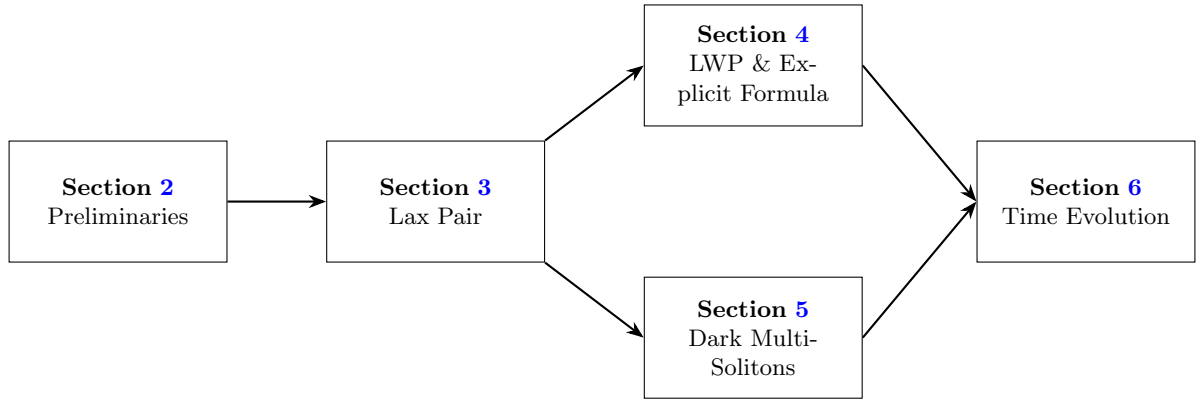
\begin{figure}[htbp]
\centering
\begin{tikzpicture}[
  box/.style={
    draw,
    rectangle,
    text width=2.6cm,
    minimum height=1.6cm,
    align=center,
    inner sep=4pt,
    font=\small
  },
  arrow/.style={->, >=Stealth, thick}
]
  \node[box] (s2) at (0, 0) {\textbf{Section~\ref{PrelimSect2}}\\Preliminaries};

  \node[box] (s3) at (4.2, 0) {\textbf{Section~\ref{SectLaxPair3}}\\Lax Pair};

  \node[box] (s4) at (8.4, 1.8) {\textbf{Section~\ref{LWPExpForSect4}}\\LWP \& Explicit Formula};
  \node[box] (s5) at (8.4, -1.8) {\textbf{Section~\ref{DarMulSoliPotSect}}\\Dark Multi-Solitons};

  \node[box] (s6) at (12.8, 0) {\textbf{Section~\ref{SectTimeEvo6}}\\Time Evolution};

  \draw[arrow] (s2) -- (s3);
  \draw[arrow] (s3.north east) -- (s4.west);
  \draw[arrow] (s3.south east) -- (s5.west);
  \draw[arrow] (s4.east) -- (s6.west);
  \draw[arrow] (s5.east) -- (s6.west);
\end{tikzpicture}
\caption{Organization of the paper.}
\label{fig:organization}
\end{figure}

\subsection*{Acknowledgments}

The author acknowledges financial support from Beijing Normal University, the Ministry of Education of China, and the China Postdoctoral Science Foundation (Fellowship Nos. 2025M773094 and 2026T190814).

\section{Preliminaries}\label{PrelimSect2}
\noindent In this section, we introduce some notations, terminologies and partial results used in this paper. Let $\mathbb{N}  = \mathbb{Z} \cap   [0, +\infty)$ and $\mathbb{N}_+ : = \mathbb{N} \setminus \{0\}$. For any  $N \in \mathbb{N}_+$, we define
\begin{equation}\label{inteintervIjN}
\mathbb{I}_j^N := \{1,2, \cdots, N\}\setminus \{j\}, \quad \forall j \in \{1,2, \cdots, N\}.
\end{equation}
 For any $A \subset \mathbb{R}$, $\mathbf{1}_A$ denotes the indicator function of $A$ (i.e., $\mathbf{1}_A(x) = 1$ for $x \in A$ and $\mathbf{1}_A(x) = 0$ otherwise).  
  Let \(\mathbf{1}\) denote the constant function equal to \(1\).  
 Let $\mathbb{S}^1 := \{ z \in \mathbb{C} : |z| =1\}$ denote the unit circle. The $1$-dimensional torus is denoted by $\mathbb{T} : = \mathbb{R}/ 2 \pi \mathbb{Z}$.  
  Upper and lower half-planes are denoted by \(\mathbb{C}_+\) and
\(\mathbb{C}_-\), respectively;   that is,
\begin{equation}\label{UpperHalfPlane}
\mathbb{C}_+ = \{z \in \mathbb{C} : \operatorname{Im} z > 0\}, \qquad
\mathbb{C}_- = \{z \in \mathbb{C} : \operatorname{Im} z < 0\}.
\end{equation}
Let $\operatorname{Hol}(\mathbb{C}_+)$ denote the set of all holomorphic functions from $\mathbb{C}_+$ to $\mathbb{C}$.  
  The space of all polynomials with complex coefficients is denoted by
\(\mathbb{C}[X]\), and, for \(N\in\mathbb{N}_+\),
\(\mathbb{C}_{\leq N-1}[X]\) denotes the subspace of polynomials of
degree at most \(N-1\).
    For any $a=(a_1, a_2, \cdots, a_N)^{\mathrm{T}} , b=(b_1, b_2, \cdots, b_N)^{\mathrm{T}} \in \mathbb{C}^N$, their inner product in $\mathbb{C}^N$ is given by
\begin{equation}\label{C^Ndotprod}
\langle a, b \rangle_{\mathbb{C}^N} := a^{\mathrm{T}} \overline{b} = b^*a = \overline{\langle b, a \rangle}_{\mathbb{C}^N} \in \mathbb{C}.
\end{equation}
We adopt the convention that the inner product is linear in the first
variable and conjugate-linear in the second.

\begin{defi}\label{matLpFuncdef}
Let \(M,N\) be positive integers. For
\(\mathsf{Z}\in\mathbb C^{M\times N}\), let \(\mathsf{Z}_{kj}\)
denote its $kj$-entry, \(1\le k\le M\), \(1\le j\le N\). The
matrix norm \(\|\mathsf{Z}\|_{\mathbb C^{M\times N}}\) is defined by
\begin{equation}\label{matnormdef}
\|\mathsf{Z}\|_{\mathbb C^{M\times N}}^2
:=
\operatorname{tr}\bigl(\mathsf{Z}\mathsf{Z}^*\bigr)
=
\sum_{k=1}^M\sum_{j=1}^N |\mathsf{Z}_{kj}|^2.
\end{equation}
Given a nonempty Lebesgue measurable subset
\(\mathfrak J\subset\mathbb R\), let
\(F:\mathfrak J\to\mathbb C^{M\times N}\) be a matrix-valued measurable function. For any \(p\in[1,+\infty)\), we set
\begin{equation}\label{LpnormFmatfunc}
\|F\|_{L^p(\mathfrak J\,;\,\mathbb C^{M\times N})}
:=
\Bigl(\int_{\mathfrak J}
\|F(x)\|_{\mathbb C^{M\times N}}^p\,\mathrm dx\Bigr)^{1/p},
\qquad 
\|F\|_{L^\infty(\mathfrak J\,;\,\mathbb C^{M\times N})}
:=
\operatorname*{ess\,sup}_{x\in\mathfrak J}
\|F(x)\|_{\mathbb C^{M\times N}}.
\end{equation}
 For any $p \in [1, +\infty]$,  the space \(L^p(\mathfrak J\,;\,\mathbb C^{M\times N})\) consists of
all functions \(F:\mathfrak J\to\mathbb C^{M\times N}\) such that
\(\|F\|_{L^p(\mathfrak J\,;\,\mathbb C^{M\times N})}<+\infty\). As
usual, functions that agree almost everywhere are identified.
\end{defi}

\begin{rem}\label{RemEquiBocNorm}
With the norm defined in \eqref{LpnormFmatfunc},
the usual H\"older inequality for matrix-valued functions holds:  if \(\frac1r=\frac1p+\frac1q\),
\(U\in L^p(\mathfrak J\,;\,\mathbb C^{M\times N})\) and
\(V\in L^q(\mathfrak J\,;\,\mathbb C^{N\times P})\), then
\(UV\in L^r(\mathfrak J\,;\,\mathbb C^{M\times P})\) and
\begin{equation}\label{HolderMat}
\|UV\|_{L^r(\mathfrak J\,;\,\mathbb C^{M\times P})}
\leq
\|U\|_{L^p(\mathfrak J\,;\,\mathbb C^{M\times N})}
\|V\|_{L^q(\mathfrak J\,;\,\mathbb C^{N\times P})}.
\end{equation}
\end{rem}

\noindent Let $\bigl( \mathfrak{N}_1, \mathfrak{d}_1\bigr)$ and $\bigl( \mathfrak{N}_2, \mathfrak{d}_2\bigr)$ be two metric spaces. We denote by $\mathcal{C} \bigl( \mathfrak{N}_1 ; \mathfrak{N}_2\bigr)$ the set of all continuous functions from $\bigl( \mathfrak{N}_1, \mathfrak{d}_1\bigr)$ to $\bigl( \mathfrak{N}_2, \mathfrak{d}_2\bigr)$.
  Let $\mathscr{S}= \mathscr{S}(\mathbb{R}; \mathbb{C})$ be the Schwartz space. Its dual space is denoted by $\mathscr{S}'$, which consists of all tempered distributions on $\mathbb{R}$. The Fourier transform is denoted by $\mathscr{F} : u \in \mathscr{S}' \mapsto \mathscr{F} (u) = \widehat{u} \in \mathscr{S}'$, where
\begin{equation}\label{FouTrDefS2}
\langle  \widehat{u} , \phi \rangle_{\mathscr{S}', \mathscr{S}} = \langle    u ,  \widehat{\phi} \rangle_{\mathscr{S}', \mathscr{S}} 
\quad \mathrm{with} \quad   \widehat{\phi}(\xi) = \int_{\mathbb{R}} \phi (x) e^{- i x \xi} \mathrm{d}x , \quad \forall \xi \in \mathbb{R}, \quad \forall  \phi  \in  \mathscr{S}.
\end{equation}For any $y \in \mathbb{R}$, the translation operator $\tau_y : u \in \mathscr{S}' \mapsto \tau_y  u \in \mathscr{S}'$ is defined as follows: 
\begin{equation}\label{TransOptauDef}
\langle \tau_y  u , \phi \rangle_{\mathscr{S}', \mathscr{S}} = \langle    u , \tau_{-y}\phi \rangle_{\mathscr{S}', \mathscr{S}}  \quad \mathrm{with} \quad 
\tau_{-y}  \phi (x) : = \phi(x + y ), \quad \forall x \in \mathbb{R}, \quad \forall \phi \in \mathscr{S}.
\end{equation}Set $\mathtt{R} := \frac{1}{2\pi} \mathscr{F}^2$. Then $ \mathtt{R} : u \in \mathscr{S}' \mapsto  \mathtt{R}(u) = \check{u} \in \mathscr{S}'$ denotes the reflection about the vertical axis: 
\begin{equation}\label{RefyAxichecDef}
\langle  \check{u} , \phi \rangle_{\mathscr{S}', \mathscr{S}} = \langle    u , \check{\phi} \rangle_{\mathscr{S}', \mathscr{S}}  \quad \mathrm{with} \quad \check{\phi}(x) = \phi(-x), \quad \forall x \in \mathbb{R}, \quad \forall \phi \in \mathscr{S}.
\end{equation}Let $\mathfrak{U}$ be an arbitrary open nonempty subset of $\mathbb{R}$.  Then $\mathcal{C}_c (\mathfrak{U}; \mathbb{C})$ denotes the space of all continuous functions   $\mathfrak{U} \to \mathbb{C}$ with compact support in $\mathfrak{U}$. For any $k \in \mathbb{N}$, we define 
\begin{equation}\label{CcscrDDefForS2}
\mathcal{C}^k_c (\mathfrak{U}; \mathbb{C}) : = \mathcal{C}^k (\mathfrak{U}; \mathbb{C})\cap    \mathcal{C}_c (\mathfrak{U}; \mathbb{C})
\quad \mathrm{and} \quad 
\mathscr{D}(\mathfrak{U}; \mathbb{C})=\mathcal{C}^{\infty}_c  (\mathfrak{U}; \mathbb{C}) : = \bigcap_{n \in \mathbb{N}} \mathcal{C}^n_c  (\mathfrak{U}; \mathbb{C}) .
\end{equation}For any $p \in [1, +\infty]$, the filtered $L^p$-spaces $L^p_+ \subset L^p (\mathbb{R}; \mathbb{C})$ and $L^p_- \subset L^p (\mathbb{R}; \mathbb{C})$ are given by
\begin{equation}\label{LppmDefForS2}
L^p_+ := \{f \in L^p(\mathbb{R}; \mathbb{C}) : \operatorname{supp} \widehat{f} \subset [0, +\infty) \}, \quad
L^p_- := \{f \in L^p(\mathbb{R}; \mathbb{C}) : \operatorname{supp} \widehat{f} \subset (-\infty, 0] \}.
\end{equation}Since every Cauchy sequence in $L^2 (\mathbb{R} ; \mathbb{C})$ has a subsequence converging pointwise almost everywhere, the Plancherel identity implies that $L^2_+$ and $L^2_-$ are both closed subspaces of the Hilbert space $L^2(\mathbb{R}; \mathbb{C})$. The $L^2$-inner product between $f \in L^2(\mathbb{R}; \mathbb{C})$ and $g \in L^2(\mathbb{R}; \mathbb{C})$ is given by
\begin{equation}\label{inProdL2DefS2}
\langle f, g \rangle_{L^2} = \langle f, g \rangle_{L^2(\mathbb{R}; \mathbb{C})} = \int_{\mathbb{R}} f(x) \overline{g(x)} \, \mathrm{d}x = \frac{1}{2\pi} \int_{\mathbb{R}} \widehat{f}(\xi) \overline{\widehat{g}(\xi)} \, \mathrm{d}\xi \in \mathbb{C}.
\end{equation}
 
\begin{lem}\label{Z1+idealinL2+}
For any $u \in L^{\infty}_+$ and $ \phi  \in L^2_+$, we have $u  \phi  \in L^2_+$.
\end{lem}
\begin{proof} 
Fix $u \in L^{\infty}_+ \subset \mathscr{S}'$. For any $f \in \mathscr{S}_+ : = \mathscr{S} \cap   L^2_+ $, the Fourier transforms $\mathfrak{v}=\widehat{u} \in \mathscr{S}'$ and $g:=\widehat{f} \in \mathscr{S}$ satisfy that $\operatorname{supp} \mathfrak{v}   \cup \operatorname{supp} g \subset \mathbb{R}_+ := [0, +\infty)$ and $2 \pi \widehat{uf} =\mathfrak{v} \ast g    \in \mathcal{C}^{\infty} \cap   L^2(\mathbb{R}; \mathbb{C})$, where
\begin{equation}\label{v*gconvfor}
 \mathfrak{v} \ast g (x)= \mathfrak{v} (\tau_x \check{g}) =\langle \mathfrak{v}, \tau_x \check{g} \rangle_{\mathscr{S}', \mathscr{S}} \in \mathbb{C}, \quad \forall x \in \mathbb{R}.
\end{equation}Fix $x<0$. If $y > x$, then $\tau_x \check{g}(y)= g(x-y)=0$.  So $\operatorname{supp}(\tau_x \check{g}) \subset  (-\infty, x] \subset  (-\infty, 0)$ and
\begin{equation*}
\operatorname{supp} \mathfrak{v} \cap   \operatorname{supp}(\tau_x \check{g})=\emptyset \quad \Longrightarrow \quad  2 \pi \widehat{uf}(x) = \mathfrak{v}(\tau_x \check{g}) = 0, \quad \forall x < 0,
\end{equation*}by formula \eqref{v*gconvfor}. Then $\operatorname{supp}(\widehat{uf}) \subset \mathbb{R}_+$ and we have $uf \in L^2_+$.  It follows that
\begin{equation}\label{Linfty+S+inL2+}
L^{\infty}_+ \cdot   \mathscr{S}_+  \subset L^2_+ \quad \text{with} \quad \mathscr{S}_+ = \mathscr{S} \cap   L^2_+ .
\end{equation}Fix  $\phi \in L^2_+$.  
Since $\mathscr{S}_+ $ is dense in the Hilbert space  $\left(L^2_+, \|.\|_{L^2} \right)$, there exists a sequence $(f_m)_{m\in \mathbb{N}}$ in $\mathscr{S}_+$ such that $\lim_{m \to +\infty} \|f_m - \phi\|_{L^2_+} = 0$. Then  formula \eqref{Linfty+S+inL2+}  implies that $uf_m \in L^2_+$ and $ (\widehat{uf_m}  )_{m \geq 0}$ admits a subsequence that converges pointwise  almost everywhere to $\widehat{u \phi}$. Since $\widehat{uf_m} = 0$ on $  (- \infty, 0)$, we have $\operatorname{supp}(\widehat{u \phi}) \subset   [0, +\infty) $. So $u \phi \in L^2_+$.
\end{proof}

\begin{cor}\label{CorLinfty-L2-}
For any $v \in L^{\infty}_-$ and $\varphi \in L^2_-$, we have $v \varphi \in L^2_-$.
\end{cor}
\begin{proof}
This is a direct consequence of  Lemma $\ref{Z1+idealinL2+}$ and the following identity:
\begin{equation}\label{suppcheu=-suppu}
\mathscr{F}(\overline{u}) = \overline{\mathtt{R}\widehat{u}} \in \mathscr{S}'  \quad \text{and} \quad
\operatorname{supp} u =    \operatorname{supp} \overline{u} =  - \operatorname{supp} \check{u} , \quad  \forall u \in\mathscr{S}' .
\end{equation}
\end{proof}
 
\begin{lem}\label{LinftplconvSchwLem}
For any $u \in L^{\infty}_+$ and $\phi \in \mathscr{S}$, we have $u * \phi \in \mathcal{C}^{\infty} (\mathbb{R} ; \mathbb{C}) \cap L^{\infty}_+$.
\end{lem}

\begin{proof}
Fix $u \in L^{\infty}_+$ and $\phi \in \mathscr{S}$,  Young's convolution inequality and Theorem 7.19 of \cite{RudinFA} imply that
\begin{equation}\label{hatofConv}
u * \phi \in \mathcal{C}^{\infty} \cap L^{\infty}  (\mathbb{R} ; \mathbb{C})  \quad \text{and} \quad \widehat{u * \phi} = \widehat{\phi} \widehat{u} \in \mathscr{S}'.
\end{equation}Then $ \operatorname{supp} \bigl(\widehat{\phi * u} \bigr)= \operatorname{supp} \bigl(\widehat{\phi} \widehat{u} \bigr) \subset \operatorname{supp} (\widehat{u}) \subset [0, +\infty)$ by \eqref{hatofConv} and Definition 6.22 of  \cite{RudinFA}. 
\end{proof}
 
\noindent If $\mathscr{X}$ and $\mathscr{Y}$ are two Banach spaces over $\mathbb{C}$, let $\mathcal{B}_{\mathbb{C}}( \mathscr{X}; \mathscr{Y})$ denote the set of all bounded linear mappings from $\mathscr{X}$ to $\mathscr{Y}$. We use the abbreviation $\mathcal{B}_{\mathbb{C}} (\mathscr{X}) := \mathcal{B}_{\mathbb{C}}( \mathscr{X}; \mathscr{X})$ and $\mathscr{X}^*: = \mathcal{B}_{\mathbb{C}}(\mathscr{X} ; \mathbb{C})$. The operator norm for an arbitrary bounded linear mapping $\mathscr{A} \in \mathcal{B}_{\mathbb{C}}( \mathscr{X}; \mathscr{Y})$ is denoted by
\begin{equation}\label{opNormAinBXY}
\| \mathscr{A} \|_{\mathcal{B}_{\mathbb{C}}( \mathscr{X}; \mathscr{Y})} := \sup\{ \|\mathscr{A} x\|_{\mathscr{Y}} : \|x\|_{\mathscr{X}} \leq 1 \}.
\end{equation}
  For any $a \in \mathscr{X}$ and $R > 0$, the open ball and closed ball centered at $a$ of radius $R$ are denoted by
\begin{equation}\label{OpBallBXaR}
 \mathsf{B}_{\mathscr{X}} (a, R) : = \{x \in \mathscr{X} : \|x -a\|_{\mathscr{X}} < R\}, \qquad \overline{\mathsf{B}}_{\mathscr{X}} (a, R) :  = \{x \in \mathscr{X} : \|x - a \|_{\mathscr{X}} \leq R\} .
\end{equation}
  The Hilbert transform $\mathbb{H}$ given by \eqref{HilbertTrIntro} is a skew-adjoint isometric isomorphism of the Hilbert space $L^2(\mathbb{R}; \mathbb{C})$. The Szeg\H{o} projectors $\Pi_{\geq 0}, \Pi_{< 0} \in \mathcal{B}_{\mathbb{C}}\bigl( L^2(\mathbb{R}; \mathbb{C}) \bigr)$ are given by
\begin{equation}\label{SzegoProjS2}
\Pi= \Pi_{\geq 0} = \tfrac{1}{2} \bigl(  \mathrm{id}_{L^2(\mathbb{R}; \mathbb{C})} + i \mathbb{H} \bigr) = \mathbf{1}_{[0, +\infty)} (D), \qquad \Pi_{< 0} = \tfrac{1}{2} \bigl( \mathrm{id}_{L^2(\mathbb{R}; \mathbb{C})} - i \mathbb{H} \bigr) = \mathbf{1}_{( - \infty, 0)} (D).
\end{equation}Then the focusing CMdNLS equation \eqref{CMdNLSfocu} can be expressed as follows: 
\begin{equation}\label{CMdNLSwithPi} 
\partial_t u = i \partial_x^2 u +  2 u \Pi_{\geq 0}  \partial_x   (|u|^2), \qquad  (t,x) \in \mathbb{R} \times \mathbb{R}.
\end{equation}Moreover, 
\eqref{LppmDefForS2} and \eqref{SzegoProjS2} imply that $L^2_+ = \Pi_{\geq 0} \bigl( L^2 (\mathbb{R} ; \mathbb{C}) \bigr)$ and $L^2_- = \Pi_{< 0} \bigl( L^2 (\mathbb{R} ; \mathbb{C}) \bigr)$.  We also have 
\begin{equation}\label{L2+perpL2-Har}
L^2(\mathbb{R}; \mathbb{C}) = L^2_+ \bigoplus L^2_-, \qquad L^2_+ \perp L^2_- ,
\end{equation}
by the Plancherel identity.
  Using formulas  \eqref{FouTrDefS2}, \eqref{LppmDefForS2} and \eqref{SzegoProjS2},  we obtain that
 \begin{equation}\label{fpi+fPI-fsum}
 f = \Pi_{\geq 0}f + \Pi_{<0}f = \Pi_{\geq 0}f + \overline{\Pi_{\geq 0}(\overline{f})} \Longleftrightarrow \Pi_{<0}f =  \overline{\Pi_{\geq 0}(\overline{f})} \in L^2_-,
 \end{equation}for any $f \in L^2(\mathbb{R}; \mathbb{C})$.
 \begin{lem}\label{L2+capL2r=0Lem}
 We have $L^2(\mathbb{R}; \mathbb{R}) \cap  L^2_+ = L^2(\mathbb{R}; \mathbb{R}) \cap  L^2_- = \{0\}$.
 \end{lem}
 \begin{proof}
 If $f = \Pi_{\geq 0} f \in L^2(\mathbb{R}; \mathbb{R}) \cap  L^2_+$, then \eqref{fpi+fPI-fsum} implies that  $f= 2 \operatorname{Re} \Pi_{\geq 0}f = 2 \operatorname{Re}f = 2f =0$. If $g = \Pi_{<0} g \in L^2(\mathbb{R}; \mathbb{R}) \cap  L^2_-$, then $\overline{g} = \Pi_{\geq 0}(\overline{g}) \in L^2(\mathbb{R}; \mathbb{R}) \cap    L^2_+ =\{0\}$.
 \end{proof}

\subsection{Sobolev Spaces with Weights and Zhidkov Spaces}
\noindent For any $x \in \mathbb{R}$, let $\langle x \rangle:=\sqrt{1+|x|^2}$ denote the Japanese bracket. Given a Lebesgue measurable subset $\mathfrak{J} \subset \mathbb{R}$ and parameters $s \in \mathbb{R}$, $p \in [1,+\infty)$, we define the weighted measure $\mu_s$ to be
\begin{equation}\label{defmus}
\mathrm{d}\mu_s(x):= \langle x \rangle^{2s} \mathrm{d} x  =(1+|x|^2)^s\mathrm{d}x.
\end{equation}Then the weighted $L^p$-space $L^p_s (\mathfrak{J})$ is  defined as follows,
\begin{equation}\label{Lpmussp}
\begin{split}
  L^p_s (\mathfrak{J})  
:=   \left\lbrace f : \mathfrak{J} \to \mathbb{C} \; \mathrm{measurable} \; \Big| \; \|f\|_{L^p_s (\mathfrak{J})}^p : =\int_{\mathfrak{J}}|f |^p\mathrm{d}\mu_s   = \int_{\mathfrak{J}}|f(x)|^p \langle x \rangle^{2s} \mathrm{d} x < +\infty\right\rbrace / \sim,
\end{split}
\end{equation}where $f \sim g$ means that $f=g$ almost everywhere on $\mathfrak{J}$. When $s=0$, we have $L^p_0(\mathfrak{J}) = L^p(\mathfrak{J}; \mathbb{C})$; that is, the weighted space reduces to the usual unweighted Lebesgue space.
 
\begin{lem}\label{TranslationCV}
If $f\in L^p_s(\mathbb{R})$ for some $s \in \mathbb{R}$ and $p \in [1,+\infty)$, then   we have 
\begin{equation}\label{PeetreCorUni}
\lim_{y \to 0}\|\tau_y f - f\|_{L^p_s(\mathbb{R})}=0
\quad \mathrm{and} \quad
\sup_{ |\theta|\leq 1}\|\tau_{\theta} f \|^p_{L^p_s(\mathbb{R})} \leq 8^{|s|} \|  f \|^p_{L^p_s(\mathbb{R})}.
\end{equation}
\end{lem}
\begin{proof}
By Peetre's inequality in Lemma 10.16 of Alazard \cite{AlazardBook2023}, for any $\xi, y, s \in \mathbb{R}$,  we have
\begin{equation}\label{PeetreIne}
\langle y +  \xi   \rangle^s \leq 2^{|s|}  \langle y  \rangle^{|s|}  \langle   \xi \rangle^s.
\end{equation}For any $f \in L^p_s(\mathbb{R})$, we use formulas \eqref{TransOptauDef},  \eqref{Lpmussp} and \eqref{PeetreIne} to obtain that
\begin{equation}\label{InetauyfLps<}
\|\tau_y f \|^p_{L^p_s(\mathbb{R})} \leq \int_{\mathbb{R}} |f(\xi)|^p (1+|y+\xi|^2)^s \mathrm{d}\xi \leq (4+4y^2)^{|s|} \|  f \|^p_{L^p_s(\mathbb{R})}, \quad \forall y \in \mathbb{R}.
\end{equation}So inequality \eqref{PeetreCorUni} is obtained. For any $\mathbf{b} \in \mathcal{C}_c(\mathbb{R}; \mathbb{C})$, there exists $\mathfrak{r}>0$ such that $\operatorname{supp}(\mathbf{b}) \subset [-\mathfrak{r}, \mathfrak{r}]$. Then   $\bigcup_{|t| \leq 1}\operatorname{supp}(\tau_t \mathbf{b}) \subset [-1 -\mathfrak{r}, 1+ \mathfrak{r}]$. We define $\varphi(t, x) : = |\mathbf{b}(x-t) -\mathbf{b}(x)
|^p (1+x^2)^s $, for any $  t \in [-1,1]$ and $x \in \mathbb{R}$. Then  $\lim_{t \to 0} \varphi(t,x) =0$,  $ \forall x \in \mathbb{R}$. Moreover, we have
\begin{equation*}
\sup_{|x| \leq 1+ \mathfrak{r}}\sup_{|t|\leq 1}|\varphi(t, x)| \leq (2+2\mathfrak{r}+ \mathfrak{r}^2)^{|s|}  2^p \|\mathbf{b}\|_{L^{\infty}(\mathbb{R}; \mathbb{C})}^p . 
\end{equation*}Using Lebesgue's dominated convergence theorem, we obtain that
\begin{equation}\label{tautb-bto0inLpsR}
\|\tau_t \mathbf{b} - \mathbf{b}\|_{L^p_s(\mathbb{R})}^p = \int_{-1- \mathfrak{r}}^{1+\mathfrak{r}}|\mathbf{b}(x-t) -\mathbf{b}(x)
|^p (1+x^2)^s \mathrm{d}x = \int_{-1- \mathfrak{r}}^{1+\mathfrak{r}} \varphi(t, x) \mathrm{d}x\to 0, \quad \mathrm{as} \quad t \to 0.
\end{equation}Fix $f\in L^p_s(\mathbb{R})$. For any $\epsilon >0$, there exists a function $b \in \mathcal{C}_c (\mathbb{R}; \mathbb{C})$ such that $\|f-b\|_{L^p_s(\mathbb{R})} \leq \frac{\epsilon}{3(1+8^{|s|})}$.  Formula \eqref{tautb-bto0inLpsR} yields that there exists $ \delta \in (0, 1)$ such that $\sup_{|y| 
\leq \delta }\|\tau_y b - b\|_{L^p_s(\mathbb{R})} \leq \frac{\epsilon}{3}$. Then
\begin{equation*}
\sup_{|y| 
\leq \delta } \|\tau_y f - f\|_{L^p_s(\mathbb{R})} \leq 2 \sup_{|\theta| \leq 1} \|\tau_{\theta} (f - b)\|_{L^p_s(\mathbb{R})}  + \sup_{|y| 
\leq \delta } \|\tau_y b - b\|_{L^p_s(\mathbb{R})} \leq \epsilon. 
\end{equation*}Then $\lim_{y \to 0} \|\tau_y f - f\|_{L^p_s(\mathbb{R})}= 0$ and the proof is completed.
\end{proof}

\noindent Let $\mathfrak{J}:=(a,b) \subset \mathbb{R}$ be an arbitrary open interval (possibly unbounded). For any $s \geq 0$ and $p \in [1, +\infty)$, the weighted Sobolev space $W^{1,p}_s(\mathfrak{J})$ is defined as follows
\begin{equation}\label{WeightedSobolevSp}
\begin{split}
W^{1,p}_s(\mathfrak{J}) := &  \Big\lbrace u \in L^p_s (\mathfrak{J}): \exists f \in L^p_s (\mathfrak{J}) \quad \text{such that} \quad  \int_{\mathfrak{J}}u \partial_x \varphi = - \int_{\mathfrak{J}} f \varphi , \quad \forall \varphi \in \mathcal{C}_c^1(\mathfrak{J}; \mathbb{C}) \Big\rbrace \\
= & \{u \in W^{1,p} (\mathfrak{J}; \mathbb{C}):  u, \; \partial_x u \in L^p_s (\mathfrak{J})  \}. 
\end{split}
\end{equation}For any $u \in W^{1,p}_s(\mathfrak{J}) $, its  $W^{1,p}_s(\mathfrak{J})$-norm   is defined as follows:  
\begin{equation}\label{Wonepnormofudef}
\| u \|_{W^{1,p}_s(\mathfrak{J})}^p : = \| u \|_{L^p_s(\mathfrak{J})}^p + \| \partial_x u \|_{L^p_s(\mathfrak{J})}^p = \int_{ \mathfrak{J} } \bigl( |u(x)|^p + | \partial_x  u(x)|^p \bigr) (1 + |x|^2)^s \mathrm{d}x .
\end{equation}
 
\begin{rem}\label{RmkWonep0nottr}
When the subscript $s=0$, note that the weighted Sobolev space $W^{1, p}_0(\mathfrak{J}) = W^{1, p} (\mathfrak{J}; \mathbb{C})$ denotes the space with weight index zero in this paper, and should not be confused with the closure of $\mathcal{C}_c^{\infty}(\mathfrak{J}; \mathbb{C})$ in $W^{1, p} (\mathfrak{J}; \mathbb{C})$. To be precise, for any $a,b \in \mathbb{R}$ such that $a<b$, we have
\begin{equation*}
W^{1, p} ((a,b); \mathbb{C}) = W^{1, p}_0((a,b)) \neq \bigl\{ f \in W^{1, p}((a,b); \mathbb{C}) : f(a) = f(b) = 0 \bigr\}.
\end{equation*}
\end{rem}

\noindent When $p=2$, we use the abbreviation $H^1_s (\mathfrak{J}):=W^{1,2}_s(\mathfrak{J}) = \{ u \in H^1(\mathfrak{J}; \mathbb{C}) : u ,   \partial_x u \in L^2_s (\mathfrak{J})\}$. We set
\begin{equation}\label{H1sWeiSobnormDef}
 \|u\|_{H^1_s (\mathfrak{J})}^2 : = \|u\|_{L^2_s (\mathfrak{J})}^2 + \|\partial_x u\|_{L^2_s (\mathfrak{J})}^2, \quad \forall  u \in H^1_s (\mathfrak{J}), \quad \forall s \geq 0 .
\end{equation}For any $p \in [1,+\infty)$ and $s \geq 0$, Theorem 8.7 of Brezis \cite{BrezisFA} yields that $\mathcal{C}_c^{\infty}(\mathbb{R})$ is dense in $W^{1,p}_s(\mathbb{R})$. Set
\begin{equation}\label{RstarpmDefForS2}
\mathbb{R}_+^*:= (0, +\infty), \quad \mathbb{R}_+:= [0, +\infty), \quad \mathbb{R}_-^*:= (-\infty, 0), \quad \mathbb{R}_-:= (-\infty, 0].
\end{equation}By Theorem 8.6 of \cite{BrezisFA}, every function in $W^{1,p}_s(\mathbb{R}_+^*)$ can be extended to a function in $W^{1,p}_s(\mathbb{R})$.
\begin{lem}\label{ExtOpLem}
For any $p \in [1, +\infty)$ and $s \geq 0$, there exists a bounded operator $\mathbb{P} : W^{1,p}_s(\mathbb{R}_+^*) \to W^{1,p}_s(\mathbb{R})$ such that $\forall u \in  W^{1,p}_s(\mathbb{R}_+^*)$, we have the following statements:\\
\noindent $(\mathbf{a}).$ $\mathbb{P} u (x) = u(x)$, for any $x > 0$;\\
\noindent $(\mathbf{b}).$ $\partial_x\left(\mathbb{P} u\right) (x) = \partial_x u(x)$, for any $x > 0$ a.e;\\
\noindent $(\mathbf{c}).$ $\|\mathbb{P} u\|_{L^p_s(\mathbb{R})} \leq 2^{\frac{1}{p}} \| u\|_{L^p_s(\mathbb{R}_+^*)}$;\\
\noindent $(\mathbf{d}).$ $\|\mathbb{P} u\|_{W^{1,p}_s(\mathbb{R})} \leq 2^{\frac{1}{p}} \| u\|_{W^{1,p}_s(\mathbb{R}_+^*)}$.
\end{lem}
\noindent The following lemma provides a characterization of the weighted Sobolev space  $W^{1,p}_s(\mathbb{R})$ in terms of the convergence of difference quotients to the weak derivative.
\begin{lem}\label{CharSobDiffQuoCV}
Assume that $f \in L^p_s (\mathbb{R})$ for some $p \in (1, +\infty)$ and $s \geq 0$. Then the following three assertions are equivalent:\\
\noindent $(\mathbf{a}).$ The function $f$ belongs to the weighted Sobolev space $ W^{1,p}_s(\mathbb{R})$.\\
\noindent $(\mathbf{b}).$ There exists a function $g \in L^p_s (\mathbb{R})$ such that $\lim_{\eta \to 0}\|\frac{1}{\eta} (\tau_{-\eta} f - f ) - g \|_{L^p_s (\mathbb{R})} = 0$.\\
 $(\mathbf{c}).$ There exist  a function $g \in L^p_s (\mathbb{R})$ and a sequence $(\eta_n )_{n \geq 0}$ with $\lim_{n \to +\infty} \eta_n = 0$ such that
 \begin{equation}\label{CVquodifffetang}
 \lim_{n \to +\infty}\|\tfrac{1}{\eta_n} (\tau_{-\eta_n} f - f ) - g \|_{L^p_s (\mathbb{R})} = 0.
 \end{equation}Moreover, when these conditions hold, the limit $g$ is equal to the weak derivative $\partial_x f \in  L^p_s(\mathbb{R})$.
\end{lem}
\begin{proof}
$(\mathbf{a}) \Rightarrow (\mathbf{b}).$ Assume that $f \in W^{1,p}_s(\mathbb{R})$ for some $p \in [1, +\infty)$ and $s \geq 0$.  Theorem 8.2 of \cite{BrezisFA} implies that $f$ is continuous and $\frac{1}{\eta}(f(x+\eta) - f(x)) - \partial_x f(x)= \int_0^1  (f' (x+\theta \eta) - f' (x)) \mathrm{d}\theta$, $\forall x \in \mathbb{R}$, $\forall \eta \in \mathbb{R} \setminus \{0\}$. Then H\"older's inequality implies that
$\|\frac{1}{\eta} (\tau_{-\eta} f - f ) - \partial_x f \|_{L^p_s (\mathbb{R})}^p \leq \int_0^1 \|\tau_{-\theta \eta} (\partial_x f) - \partial_x f\|_{L^p_s (\mathbb{R})}^p \mathrm{d}\theta \to 0$, by Lemma $\ref{TranslationCV}$, Peetre's inequality \eqref{PeetreIne} and the dominated convergence theorem. $(\mathbf{b}) \Rightarrow (\mathbf{c})$ is direct. \\

\noindent $(\mathbf{c}) \Rightarrow (\mathbf{a}).$ Fix 
$p \in (1, +\infty)$ and $s \geq 0$. 
Set $q:= \frac{p}{p-1}\in (1, +\infty)$. Assume that $g \in L^p_s(\mathbb{R})$,  $(\eta_n)_{n\in \mathbb{N}}$ is a   sequence that tends to $0$ and  $\lim_{n \to +\infty}\|\tfrac{1}{\eta_n} (\tau_{-\eta_n} f - f ) - g \|_{L^p_s (\mathbb{R})} = 0$.
  For any $\varphi \in \mathcal{C}_c^1 (\mathbb{R} ; \mathbb{C}) \subset W^{1, q } (\mathbb{R} ; \mathbb{C})$, the previous step yields that
$\lim_{n \to +\infty}\|\tfrac{1}{\eta_n} (\tau_{\eta_n} \varphi - \varphi ) + \partial_x \varphi \|_{L^q (\mathbb{R})} = 0$. Then
\begin{equation}
\small
\Big|\int_{\mathbb{R}}  f \partial_x \varphi + g \varphi \Big| =   \Big| \int_{\mathbb{R}} f \left(\tfrac{1}{\eta_n} (\tau_{\eta_n} \varphi - \varphi ) + \partial_x \varphi \right) - \left( \tfrac{1}{\eta_n} (\tau_{-\eta_n} f - f ) - g \right) \varphi \Big| \to 0,
\end{equation}as $n \to +\infty$, by H\"older's inequality. So $ f \in W^{1, p}(\mathbb{R}) \cap  L^p_s(\mathbb{R})$ and $\partial_x f = g \in L^p_s(\mathbb{R})$.
\end{proof}

\noindent The following result is a direct consequence of  Lemma $\ref{ExtOpLem}$ and Lemma $\ref{CharSobDiffQuoCV}$.
\begin{cor}\label{QuoDiffCVR+*Cor}
Assume that $f \in W^{1,p}_s(\mathbb{R}_+^*)$ for some $p \in (1, +\infty)$ and $s \geq 0$, then  
\begin{equation}
\lim_{\eta \to 0^+}\|\tfrac{1}{\eta} (\tau_{-\eta} f - f ) - \partial_x f \|_{L^p_s (\mathbb{R}_+^*)} = 0.
\end{equation} 
\end{cor}
\begin{proof}
If $f \in W^{1,p}_s(\mathbb{R}_+^*)$, it admits an extension $F := \mathbb{P}f \in W^{1,p}_s(\mathbb{R})$, where $\mathbb{P}: W^{1,p}_s(\mathbb{R}_+^*) \to  W^{1,p}_s(\mathbb{R})$ denotes the extension operator constructed in Lemma $\ref{ExtOpLem}$. Then $F(x)=f(x)$ and $\partial_x F (x)= \partial_x f(x)$, $\forall x >0$, a.e. For any $\eta >0$, we have $x+\eta >0$, $\forall x >0$. Consequently, we obtain that
\begin{equation*}
\|\tfrac{1}{\eta}(\tau_{-\eta}f-f) - \partial_x f\|_{L^p_s(\mathbb{R}_+^*)}   \leq \|\tfrac{1}{\eta}(\tau_{-\eta} F  - F ) -\partial_x F \|_{L^p_s(\mathbb{R})}  \to 0, 
\end{equation*}as $\eta \to 0^+$, by using Lemma $\ref{CharSobDiffQuoCV}$.  
\end{proof}

\noindent For any $s \in \mathbb{R}$, the  Sobolev space $H^s (\mathbb{R}; \mathbb{C})$ is defined by   $H^s(\mathbb{R}; \mathbb{C}) := \{ v \in \mathscr{S}' : \widehat{v} \in L^2_s(\mathbb{R})\}$. We define 
  \begin{equation}\label{vHsNormDefFor}
\| v \|_{H^s  }^2 = \| v \|_{H^s (\mathbb{R}; \mathbb{C}) }^2 : = \frac{1}{2 \pi}\|\widehat{v}\|_{L^2_s (\mathbb{R})}^2  =  \frac{1}{2 \pi}\int_{ \mathbb{R}} |\widehat{v}(\xi)|^2 \left(1+ |\xi|^2 \right)^s \mathrm{d}\xi , \quad \forall v \in H^s (\mathbb{R} ; \mathbb{C}).
 \end{equation}Then $ H^s(\mathbb{R}; \mathbb{C})$ is a Hilbert space, equipped with the following  inner product
 \begin{equation}\label{HsInProdDefFor}
\langle u, v \rangle_{ H^s  } = \langle u, v \rangle_{ H^s (\mathbb{R}; \mathbb{C})} : = \frac{1}{2 \pi}\int_{ \mathbb{R}} \widehat{u}(\xi) \overline{\widehat{v}(\xi)}  \left(1+ |\xi|^2 \right)^s \mathrm{d}\xi, \quad \forall u,v \in H^s (\mathbb{R}, \mathbb{C}).
 \end{equation}Recall that $L^2(\mathbb{R}; \mathbb{C}) = H^0 (\mathbb{R} ; \mathbb{C})$ and $W^{1,2} (\mathbb{R}; \mathbb{C}) = H^1 (\mathbb{R}; \mathbb{C})$.   Formula  \eqref{vHsNormDefFor} yields that   
 \begin{equation}\label{varphiHs+normHsHs-1}
 \| v \|_{H^s (\mathbb{R}; \mathbb{C}) }^2  =  \| v \|_{H^{s-1} (\mathbb{R};\mathbb{C})}^2  +  \| \partial_x  v \|_{H^{s-1} (\mathbb{R};\mathbb{C})}^2 , \quad \forall v \in H^s (\mathbb{R}; \mathbb{C}).
 \end{equation}For any $N \in \mathbb{N}_+$ and $u \in H^N (\mathbb{R} ; \mathbb{C})$, we use H\"older's inequality and \eqref{vHsNormDefFor} to obtain that
 \begin{equation}\label{ContrL2dNnorHNnorm}
 \|u \|_{L^2(\mathbb{R}; \mathbb{C})} + \|\partial_x^N u \|_{L^2(\mathbb{R}; \mathbb{C})}  \leq  \sqrt{2} \|u \|_{H^N(\mathbb{R}; \mathbb{C})}  \leq 2 ^{\frac{N}{2}} \left(  \|u \|_{L^2(\mathbb{R}; \mathbb{C})} + \|\partial_x^N u \|_{L^2(\mathbb{R}; \mathbb{C})}  \right).
\end{equation}For any $s >\frac{1}{2}$, 
 there exists $\mathfrak{C}_s' >0$ such that if  $  u, v \in H^s(\mathbb{R}; \mathbb{C})$, then we have
\begin{equation}\label{HNAlgebra}
\|uv\|_{H^s (\mathbb{R}; \mathbb{C})} \leq \mathfrak{C}_s'\|u \|_{H^s (\mathbb{R}; \mathbb{C})}\| v\|_{H^s (\mathbb{R}; \mathbb{C})}.
\end{equation}
  For any $n \in \mathbb{N}_+ =[1, +\infty) \cap  \mathbb{Z}$,  the Zhidkov space $\mathcal{Z}^n=\mathcal{Z}^n (\mathbb{R}; \mathbb{C})$ consists of all essentially bounded functions  $u \in L^{\infty}(\mathbb{R}; \mathbb{C})$ whose weak derivative $\partial_x u$ belongs to the Sobolev space $H^{n-1}(\mathbb{R}; \mathbb{C})$ i.e. 
\begin{equation}\label{ZhiSpdefSect2}
\mathcal{Z}^n := \{u \in L^{\infty} (\mathbb{R}; \mathbb{C}): \partial_x u \in  H^{n-1}(\mathbb{R}; \mathbb{C})\} \subset H^{n}_{\mathrm{loc}}(\mathbb{R}; \mathbb{C}).
\end{equation}Then the Sobolev embedding $H^1 \hookrightarrow L^{\infty}$ yields that $H^n(\mathbb{R}; \mathbb{C})\subset \mathcal{Z}^n$, $\forall n \geq 1$. We define that
\begin{equation}\label{ZhidNormZn}
\|u\|_{\mathcal{Z}^n} :=  \|u\|_{L^{\infty}(\mathbb{R}; \mathbb{C})} + \|\partial_x u\|_{H^{n-1}(\mathbb{R}; \mathbb{C})},
\end{equation}if $u : \mathbb{R} \to \mathbb{C}$ is measurable. Then a measurable function $u \in \mathcal{Z}^n$ if and only if $\|u\|_{\mathcal{Z}^n} < +\infty$.
\begin{lem}\label{ZhidSpisBanach}
For any integer $n \geq 1$, the Zhidkov space $\left(\mathcal{Z}^n, \|.\|_{\mathcal{Z}^n} \right)$ is a Banach space. 
Moreover,  
\begin{equation}\label{ZNcapL2=HN}
\mathcal{Z}^n \cap L^2 (\mathbb{R} ; \mathbb{C}) = H^n (\mathbb{R} ; \mathbb{C}).
\end{equation}
\end{lem}
\begin{proof}
Assume that $(u_m)_{m \geq 0}$ is a Cauchy sequence in $\left(\mathcal{Z}^n, \|.\|_{\mathcal{Z}^n} \right)$. By the completeness of $L^{\infty}(\mathbb{R}; \mathbb{C})$ and $H^{n-1}(\mathbb{R}; \mathbb{C})$, there exist $u \in L^{\infty}(\mathbb{R}; \mathbb{C})$ and $g \in H^{n-1}(\mathbb{R}; \mathbb{C})$ such that
\begin{equation}\label{untouduntog}
\lim_{m \to +\infty}\|u_m - u\|_{L^{\infty}(\mathbb{R}; \mathbb{C})} = \lim_{m \to +\infty}\|\partial_x u_m - g\|_{H^{n-1}(\mathbb{R}; \mathbb{C})} = 0.
\end{equation}For any $\varphi \in \mathcal{C}_c^1(\mathbb{R}; \mathbb{C})$, we have $\int_{\mathbb{R}}g \varphi =\lim_{m \to +\infty}\int_{\mathbb{R}}\varphi \partial_x u_m  =-\lim_{m \to +\infty}\int_{\mathbb{R}}u_m  \partial_x\varphi  = -\int_{\mathbb{R}}u  \partial_x\varphi$, by   using H\"older's inequality and \eqref{untouduntog}. So $g = \partial_x u$ and $\lim_{m \to +\infty}\|u_m - u\|_{\mathcal{Z}^n}=0$.   
Formula \eqref{ZNcapL2=HN} follows directly by  
\eqref{vHsNormDefFor},  \eqref{ContrL2dNnorHNnorm}, 
\eqref{ZhidNormZn} and the Sobolev embedding $H^n (\mathbb{R} ; \mathbb{C}) \hookrightarrow L^{\infty} (\mathbb{R} ; \mathbb{C})$.
\end{proof} 
\noindent For any $k \in \mathbb{N}$ and $\rho \in (0,1]$, the H\"older space $\mathcal{C}^{k, \rho}=\mathcal{C}^{k, \rho}(\mathbb{R}; \mathbb{C})$ consists of all $k$-times continuously differentiable functions $u \in \mathcal{C}^k(\mathbb{R}; \mathbb{C})$ such that the following H\"older norm
\begin{equation}\label{Holdernormfor}
\|u\|_{\mathcal{C}^{k, \rho}} : = \max_{|\alpha| \leq k} \|\partial_x^{\alpha} u\|_{L^{\infty}(\mathbb{R}; \mathbb{C})} + \sup_{x\ne y}\tfrac{|\partial_x^k u(x) - \partial_x^k u(y) |}{|x-y|^{\rho}}
\end{equation}is finite. Then we have the following Zhidkov embedding lemma.
\begin{lem}
For any  integer $n \geq 1$, there exists $\mathscr{C}_n>0$ such that for any $u \in \mathcal{Z}^n$, we have
\begin{equation}\label{ZhidEmbIne}
\|u\|_{\mathcal{C}^{n-1, \frac{1}{2}}} \leq \mathscr{C}_n \|u\|_{\mathcal{Z}^n} .
\end{equation}In other words, the injection $\mathcal{Z}^n \hookrightarrow \mathcal{C}^{n-1, \frac{1}{2}}$ is continuous.
\end{lem}
 
\noindent It is proved
in Proposition 2.2 of Gallo \cite{GalloSchGrZhidSp} and
 the Sobolev embedding theorem $H^N (\mathbb{R}; \mathbb{C}) \hookrightarrow \mathcal{C}^{N-1, \frac{1}{2}} $, $\forall N \in \mathbb{N}_+$, see  Theorem 1.66 of \cite{BaChDaAnalysis} for details.
 We recall Sobolev embedding  $H^1(\mathbb{R}; \mathbb{C}) \hookrightarrow L^{\infty} $  and set 
\begin{equation}\label{InffH1overfLinfty}
\mathfrak{m}^{-1}:= \inf_{f \in H^1 \setminus \{0\}} \tfrac{\|f\|_{H^1(\mathbb{R}; \mathbb{C})}}{\|f\|_{L^{\infty} (\mathbb{R}; \mathbb{C})}} \in \mathbb{R}_+^* = (0, +\infty) .
\end{equation}Set $\varphi(x):=\frac{1}{x + i}$. Then $ \|\varphi\|_{H^1(\mathbb{R}; \mathbb{C})}  = \sqrt{\frac{3 \pi}{2}} \|\varphi\|_{L^{\infty}(\mathbb{R}; \mathbb{C})}$. Since $\int_{\mathbb{R}} \frac{1}{1+x^2}  \mathrm{d}x = \pi$,  we have
\begin{equation}\label{cstmsmallsq2-1}
\sqrt{\tfrac{2}{3 \pi}} \leq  \mathfrak{m} \leq \tfrac{1}{\sqrt{2}},
\end{equation}by the  Cauchy--Schwarz inequality. 
For any $N \in \mathbb{N}_+$. we have $H^N (\mathbb{R}; \mathbb{C}) \hookrightarrow \mathcal{Z}^N$ and 
\begin{equation}\label{HNhookarrZN}
\|f\|_{\mathcal{Z}^N} \leq (1 + \mathfrak{m})\|f\|_{H^N(\mathbb{R}; \mathbb{C})},  \qquad \forall f \in H^N(\mathbb{R} ; \mathbb{C}),
\end{equation}by \eqref{vHsNormDefFor},  \eqref{ZhidNormZn} and \eqref{InffH1overfLinfty}.
The next lemma is another consequence of Sobolev embedding theorem.
\begin{lem}\label{ZnHninHnLem}
For any  integer $n \geq 1$, there exists a constant $C_n>0$ such that
\begin{equation}\label{ufHn<uZnfHn}
\|uf\|_{H^n(\mathbb{R}; \mathbb{C})} \leq C_n \|u\|_{\mathcal{Z}^n} \|f\|_{H^n(\mathbb{R}; \mathbb{C})}, \quad \forall u \in \mathcal{Z}^n, \quad \forall  f \in H^n(\mathbb{R}; \mathbb{C}).
\end{equation}Moreover, we have $\partial_x^m (uf) = \sum_{k=0}^m \frac{m!}{k! (m-k)!} \partial_x^k u \partial_x^{m-k}f \in H^{n-m}(\mathbb{R}; \mathbb{C})$, if $0 \leq m \leq n$.
\end{lem}
\begin{proof}
The proof is based on mathematical induction on $n \in \mathbb{N}_+$.\\

\noindent \textbf{Step \romannumeral 1}. Fix $n=1$. If $u \in \mathcal{Z}^1 \subset H^1_{\mathrm{loc}}(\mathbb{R}; \mathbb{C}) $ and $f \in H^1 \subset  L^{\infty}(\mathbb{R}; \mathbb{C})$, Corollary 8.10 of \cite{BrezisFA} yields that 
\begin{equation}\label{LeibForn=1}
\partial_x (uf) = f \partial_x u + u \partial_x f \in L^2(\mathbb{R}; \mathbb{C}), \quad \|\partial_x (uf)\|_{L^2(\mathbb{R}; \mathbb{C})} \leq \|f\|_{L^{\infty}}\|\partial_x u\|_{L^2} +  \|u\|_{L^{\infty}}\|\partial_x f\|_{L^2}.
\end{equation}Since $\|uf\|_{L^2} \leq \|u\|_{L^{\infty}}\|f\|_{L^2}$,  we use formulas \eqref{ZhidNormZn}, \eqref{InffH1overfLinfty} and \eqref{LeibForn=1} to obtain that
\begin{equation}\label{Z1H1subsetH1}
\| uf \|_{H^1(\mathbb{R}; \mathbb{C})} \leq \| u \|_{L^{\infty}}(\| f \|_{L^2} + \| \partial_x f \|_{L^2}) + \| f \|_{L^{\infty}} \|\partial_x u \|_{L^2} \leq  \sqrt{2} \| u \|_{\mathcal{Z}^1}\| f \|_{H^1(\mathbb{R}; \mathbb{C})}.
\end{equation}\textbf{Step \romannumeral 2}. Assume that Lemma $\ref{ZnHninHnLem}$ holds for $n =N \geq 1$. If $u \in \mathcal{Z}^{N+1}$ and $f \in H^{N+1}(\mathbb{R};\mathbb{C})$, then
\begin{equation}\label{IndAssumpn=NZNHN}
\|uf\|_{H^N } \leq C_N \|u\|_{\mathcal{Z}^N} \|f\|_{H^N }, \quad \partial_x^m (uf) = \sum_{k=0}^m \frac{m!}{k! (m-k)!} \partial_x^k u \partial_x^{m-k}f \in H^{N-m} (\mathbb{R}; \mathbb{C}),
\end{equation}for any $m \in  [0,  N] \cap  \mathbb{Z}$. Since $\partial_x u, \partial_x f \in H^N(\mathbb{R}; \mathbb{C})$, we use \eqref{HNAlgebra}, \eqref{LeibForn=1} and \eqref{IndAssumpn=NZNHN}    to obtain that
\begin{equation}\label{HNfordx(uf)}
\|\partial_x (uf)\|_{H^N (\mathbb{R}; \mathbb{C})} \leq \|f\partial_x  u\|_{H^N}+ \|u \partial_x f\|_{H^N} \leq \mathfrak{C}_N'  \|f \|_{H^N} \| \partial_x  u\|_{H^N} + C_N \|u \|_{\mathcal{Z}^N} \| \partial_x  f\|_{H^N}.
\end{equation}So we have $\| uf \|_{H^{N+1} (\mathbb{R}; \mathbb{C})} \leq \| uf \|_{H^N} +  \|\partial_x (uf)\|_{H^N} \leq (\sqrt{2} C_N +\mathfrak{C}_N') \|u \|_{\mathcal{Z}^{N+1}} \| f\|_{H^{N+1}}$ by \eqref{ZhidNormZn} and \eqref{HNfordx(uf)}. Set $\binom{N}{k} = \frac{N!}{k!(N-k)!}$. Then $\binom{N}{k}  + \binom{N}{k-1} =\binom{N+1}{k} $, if $1 \leq k \leq N$. By \eqref{LeibForn=1} and \eqref{IndAssumpn=NZNHN}, we have $\partial_x^{N+1} (uf) = u \partial_x^{N+1} f + \sum_{k=1}^N \binom{N}{k}  \partial_x^k u \partial_x^{N+1-k} f + \sum_{l=0}^N \binom{N}{l}  \partial_x^{l+1} u \partial_x^{N-l} f =\sum_{k=0}^{N+1} \binom{N+1}{k}  \partial_x^k u \partial_x^{N+1-k} f$.
\end{proof}

\begin{lem}\label{ZnBanachAlgeLem}
For any positive integer $n \in \mathbb{N}_+$, there exists $\mathfrak{C}_n>0$ such that $\forall u, v \in \mathcal{Z}^n$, we have
\begin{equation}\label{uvZn<uZnvZn}
\|u v\|_{\mathcal{Z}^n} \leq \mathfrak{C}_n \|u \|_{\mathcal{Z}^n}\| v \|_{\mathcal{Z}^n}.
\end{equation}Consequently, the Zhidkov space $\left(\mathcal{Z}^n, \mathfrak{C}_n \|.\|_{\mathcal{Z}^n} \right)$ is a commutative  Banach algebra with identity $\mathbf{1} \in  \mathcal{Z}^n$, where $\mathbf{1}(x)=1$, $\forall x \in \mathbb{R}$. The Sobolev space $H^n(\mathbb{R}; \mathbb{C})$ is an ideal of the ring $\left(\mathcal{Z}^n,  +, \times \right)$.
\end{lem}
\begin{proof}
The proof is divided into two steps. \\

\noindent \textbf{Step \romannumeral 1}. Fix $n=1$ and $u,v \in \mathcal{Z}^1 \subset H^1_{\mathrm{loc}} \cap  L^{\infty}(\mathbb{R}; \mathbb{C})$. Then we have $\|uv\|_{L^{\infty}} \leq \|u\|_{L^{\infty}}   \|v\|_{L^{\infty}}$ and
\begin{equation}\label{d(uv)foruvZ1}
\partial_x (uv) = v \partial_x u + u \partial_x v \in L^2(\mathbb{R}; \mathbb{C}), \quad \|\partial_x (uv)\|_{L^2(\mathbb{R}; \mathbb{C})} \leq \|v\|_{L^{\infty}}\|\partial_x u\|_{L^2} +  \|u\|_{L^{\infty}}\|\partial_x v\|_{L^2},
\end{equation}by using Corollary 8.10 in \cite{BrezisFA}. Then formula \eqref{ZhidNormZn} yields that
\begin{equation}\label{Z1Banachalgebra}
 \| u v\|_{\mathcal{Z}^1} \leq  \|u\|_{L^{\infty}}   \|v\|_{L^{\infty}}+ \|v\|_{L^{\infty}}\|\partial_x u\|_{L^2} +  \|u\|_{L^{\infty}}\|\partial_x v\|_{L^2} \leq  \| u \|_{\mathcal{Z}^1} \|  v\|_{\mathcal{Z}^1}.
\end{equation}\textbf{Step \romannumeral 2}. Fix $n :=m+1 \in [2, +\infty) \cap  \mathbb{Z}$ and $u,v \in \mathcal{Z}^n \subset H^n_{\mathrm{loc}} \cap  \mathcal{Z}^1$. Then  \eqref{ufHn<uZnfHn} and \eqref{Z1Banachalgebra} imply that
\begin{equation}\label{dxuvHn-1uvinZn}
\|\partial_x (uv)\|_{H^{n-1}(\mathbb{R}; \mathbb{C})} \leq \|v\partial_x u\|_{H^m} + \|u \partial_x v\|_{H^m} \leq C_m (\| v\|_{\mathcal{Z}^m}\| \partial_x u\|_{H^m} + \| u\|_{\mathcal{Z}^m}\| \partial_x v\|_{H^m} ).
\end{equation}Using  formulas \eqref{ZhidNormZn} and \eqref{dxuvHn-1uvinZn}, we obtain that 
\begin{equation*}
\begin{split}
  \|uv\|_{\mathcal{Z}^n} \leq &  (1+C_{n-1})  \| v\|_{\mathcal{Z}^{n-1}} \left(\| u\|_{L^{\infty}} + \| \partial_x u\|_{H^{n-1}} \right) + C_{n-1}\| u\|_{\mathcal{Z}^{n-1}}\| \partial_x v\|_{H^{n-1}} \\ \leq & (1+ 2 C_{n-1})  \|u \|_{\mathcal{Z}^n} \| v\|_{\mathcal{Z}^n}.
\end{split}
\end{equation*}We choose $\mathfrak{C}_1 = 1$ and $\mathfrak{C}_n := 2  C_{n-1} +1 $, $\forall n \geq 2$. Finally, inequality \eqref{ufHn<uZnfHn} shows that Sobolev space $H^n(\mathbb{R}; \mathbb{C)}$ is an ideal  in the ring $\left( \mathcal{Z}^n, +, \times \right)$ when equipped with pointwise multiplication. That is
\begin{equation}\label{HnZnsubsetHn}
H^n(\mathbb{R}; \mathbb{C)} \cdot \mathcal{Z}^n \subset H^n(\mathbb{R}; \mathbb{C)} , \quad \forall n \in \mathbb{N}_+.
\end{equation}Then the proof is completed.
\end{proof}
\noindent The smooth Sobolev space $H^{\infty}(\mathbb{R}; \mathbb{C})$ and smooth Zhidkov space $\mathcal{Z}^{\infty}$ are given by
\begin{equation}\label{SmoothSobZhidSp}
H^{\infty}(\mathbb{R}; \mathbb{C}):= \bigcap_{n \in \mathbb{N}_+} H^n  (\mathbb{R}; \mathbb{C}), 
\qquad
\mathcal{Z}^{\infty}:= \bigcap_{n \in \mathbb{N}_+} \mathcal{Z}^n.  
\end{equation}
 
\begin{lem}\label{ZinftydenseinZN}
For any $n \in \mathbb{N}_+$, the subspace $\mathcal{Z}^{\infty}$ is dense in the Zhidkov space $\bigl( \mathcal{Z}^n , \| . \|_{\mathcal{Z}^n} \big)$.
\end{lem}   

\begin{proof}
Since $\mathcal{C}_c^{\infty}(\mathbb{R};\mathbb{C}) * \mathcal{Z}^1 \subset \mathscr{S} * \mathcal{Z}^1 \subset\mathcal{Z}^{\infty}$ by Young's convolution inequality, the result follows by the mollification argument used in Proposition~2.1 of \cite{GalloSchGrZhidSp}.
\end{proof}

\noindent For any $t \in \mathbb{R}$, we define $\mathtt{E}_t(\xi) : = \exp (-it \xi^2) \in \mathbb{C}$, $\forall \xi \in \mathbb{R}$. Then $\mathtt{E}_t \in \mathcal{C}^{\infty}(\mathbb{R}; \mathbb{C})$ satisfies that
\begin{equation}\label{SchClinvbyEt}
\mathtt{E}_t \cdot \mathscr{S} \subset \mathscr{S} \quad \mathrm{and} \quad \mathtt{E}_t \cdot \mathscr{S}' \subset \mathscr{S}'.
\end{equation}Recall that $D  = -i \partial_x$. We define $e^{it \partial_x^2} : = \mathtt{E}_t (D) : \mathscr{S}' \to \mathscr{S}'$ as  a Fourier multiplier. Specifically,
\begin{equation}\label{expitdeltaDefFor}
\mathscr{F}  \bigl( e^{it \partial_x^2} u  \bigr) : = \mathtt{E}_t \widehat{u} \in \mathscr{S}', \quad \forall u \in \mathscr{S}'.
\end{equation}Then $\bigl( e^{it \partial_x^2} \bigr)_{t \in \mathbb{R}}$ is a group of operators on $\mathscr{S}'$ and each $ e^{it \partial_x^2}$ commutes with other Fourier multipliers. Specifically, for any $t,r \in \mathbb{R}$, the following identities hold in $\mathscr{S}'$:
\begin{equation}\label{GrouPropexpitdelta}
 e^{i(t +r) \partial_x^2}( u) =
\bigl( e^{it \partial_x^2}e^{ir \partial_x^2} \bigr)u  \in \mathscr{S}' \quad \mathrm{and} \quad     e^{i t \partial_x^2}  \partial_x^m u  = \partial_x^m  e^{i t \partial_x^2} (u ) \in \mathscr{S}' ,  \quad  \forall u \in \mathscr{S}'.
\end{equation}
Since multiplying $\mathtt{E}_t$ leaves the Dirac mass $\delta_0 \in \mathscr{S}'$ invariant, we have
\begin{equation}\label{SchGreitLap1=1}
\mathtt{E}_t \delta_0 = \delta_0 \in \mathscr{S}' \quad 
\Longleftrightarrow \quad 
e^{i t \partial_x^2 } \mathbf{1} = \mathbf{1} \in \mathscr{S}'.
\end{equation}
For any $s \in \mathbb{R}$, 
$\bigl( e^{it \partial_x^2} \bigr)_{t \in \mathbb{R}}$ is a strongly continuous unitary group on the Sobolev space $H^s(\mathbb{R}; \mathbb{C})$. Then  
\begin{equation}\label{SchrGrUnitHs}
\big\|e^{it \partial_x^2}  (\phi) \big\|_{H^s(\mathbb{R} ; \mathbb{C})} =  \big\|\phi\big\|_{H^s(\mathbb{R} ; \mathbb{C})}, \quad \forall \phi \in H^s(\mathbb{R} ; \mathbb{C}), \quad \forall t \in \mathbb{R} .
\end{equation}

\begin{lem}
For any $u \in \mathcal{Z}^1$ and $t \in \mathbb{R}$, we have $e^{it \partial_x^2} u -u \in H^1 (\mathbb{R}; \mathbb{C})$ and
\begin{equation}\label{estiL2EtDuubyu}
\| e^{it \partial_x^2} u - u  \|_{L^2 (\mathbb{R} ; \mathbb{C})} \leq \bigl( 2 + \sqrt{2} \bigr) \sqrt{|t|} \|\partial_x u\|_{L^2(\mathbb{R} ; \mathbb{C})} .
\end{equation}
\end{lem}
\begin{proof}
If $t = 0$, the inequality is trivial. We now assume $t \neq 0$.  
Fix a function $\rho \in \mathcal{C}_c^{\infty}(\mathbb{R}; \mathbb{R})$ such that $0 \leq \rho \leq 1$, $\operatorname{supp} \rho \subset [-2, 2]$ and $\rho(x) = 1$,  $\forall  x \in [-1,1]$. Recall that $\mathtt{E}_t(\xi) = \exp(- i t \xi^2)$. We adopt the idea in 
Lemma 3 of G\'erard \cite{GerGP2006} and 
 Lemma 2.3 of \cite{Xchenccm25ExpFor} to decompose $\xi^{-1} \bigl( \mathtt{E}_t(\xi) -1 \bigr)$ into two parts:  set
\begin{equation}\label{bskappatxi2parts}
\boldsymbol{\kappa} (t, \xi) : = - i t \xi \rho (t \xi^2) \int_0^1 \exp( - i \mathfrak{s} t \xi^2 )\mathrm{d}\mathfrak{s} + \xi^{-1} \bigl( 1 - \rho (t \xi^2) \bigr) \bigl( \mathtt{E}_t(\xi) -1  \bigr)   = \xi^{-1} \bigl( \mathtt{E}_t(\xi) -1 \bigr),
\end{equation}for any $t \in \mathbb{R}$ and $\xi \in \mathbb{R} \setminus \{0\}$. 
Set $\boldsymbol{\kappa} (t, 0) := 0$. 
There are $3$ cases to discuss. \\ 

\noindent \textbf{Case 1.} If $0 < |\xi | \leq \sqrt{|t|}^{-1}$, then $\rho( t \xi^2) = 1$. So formula \eqref{bskappatxi2parts} implies that 
\begin{equation}\label{Case1xi<|t|-1}
|\boldsymbol{\kappa} (t, \xi)| \leq \big| - i t \xi \rho (t \xi^2) \int_0^1 \exp( - i \mathfrak{s} t \xi^2 )\mathrm{d}\mathfrak{s} \big| \leq |t| |\xi| \leq \sqrt{|t|}.
\end{equation}\noindent \textbf{Case 2.}  If $  |\xi | >  \sqrt{ \frac{2}{|t|}}$, then $\rho( t \xi^2) = 0$. So formula \eqref{bskappatxi2parts} implies that 
\begin{equation}\label{Case2xi>2|t|-1}
|\boldsymbol{\kappa} (t, \xi)| \leq \big| \xi^{-1} \bigl( 1 - \rho (t \xi^2) \bigr) \bigl( \exp(-it \xi^2) -1  \bigr) \big| \leq 2 |\xi|^{-1} \leq  \sqrt{2|t|}.
\end{equation}\noindent \textbf{Case 3.}  If $\sqrt{|t|}^{-1} < |\xi | \leq \sqrt{ \frac{2}{|t|}} $, then $0 \leq \rho( t \xi^2)  \leq 1$. So formula \eqref{bskappatxi2parts} implies that
\begin{equation}\label{Case3t-1<xi<2|t|-1}
|\boldsymbol{\kappa} (t, \xi)| \leq \big|  t \xi \rho (t \xi^2)   \big| 
+ 2 \big| \xi^{-1} \bigl( 1 - \rho (t \xi^2) \bigr)   \big| \leq |t||\xi| + 2 |\xi|^{-1}\leq  \bigl( 2 + \sqrt{2} \bigr) \sqrt{|t|}.
\end{equation}From \eqref{Case1xi<|t|-1} and \eqref{Case2xi>2|t|-1}, it follows that inequality \eqref{Case3t-1<xi<2|t|-1} holds for all $(t,  \xi ) \in \mathbb{R} \times \mathbb{R}$. So we have
\begin{equation}\label{Estibskappa123}
\exp(-it \xi^2) - 1 = \xi \boldsymbol{\kappa} (t, \xi) \quad \mathrm{and} \quad |\boldsymbol{\kappa} (t, \xi)| \leq \bigl( 2 + \sqrt{2} \bigr) \sqrt{|t|}, \quad \forall (t, \xi) \in \mathbb{R} \times \mathbb{R}.
\end{equation}Finally, inequality \eqref{estiL2EtDuubyu} follows directly from \eqref{vHsNormDefFor} and \eqref{Estibskappa123}.
\end{proof}
\begin{lem}\label{StrCSchGrZhidlem}
For any $N \in \mathbb{N}_+$, the family $\bigl( e^{it \partial_x^2} \bigr)_{t \in \mathbb{R}}$ is a strongly continuous group of bounded operators on the Zhidkov space $\mathcal{Z}^N$. If $t \in \mathbb{R}$, then $e^{it \partial_x^2} - \mathrm{id}_{\mathcal{Z}^N} \in \mathcal{B}_{\mathbb{C}} \bigl( \mathcal{Z}^N; H^N(\mathbb{R} ; \mathbb{C}) \bigr)$. Moreover,  we have
\begin{equation}\label{expitLap-iduHNBd}
\big\|e^{it \partial_x^2} (u) - u \big\|_{H^N(\mathbb{R} ; \mathbb{C}) } \leq 2^{\frac{N+1}{2}} \bigl( 1 +2 \sqrt{|t|} \bigr)\| u  \|_{\mathcal{Z}^N}\quad \mathrm{and} \quad 
\lim_{t \to 0}\|e^{it \partial_x^2} (u) - u \|_{H^N(\mathbb{R} ; \mathbb{C}) } =0,
 \end{equation}for any $u \in \mathcal{Z}^N$. Recall that $\mathfrak{m}>0$ is given by  \eqref{InffH1overfLinfty}. Then the following estimate holds: 
\begin{equation}\label{BdZNexpitlapwrtt}
 \big\|e^{it \partial_x^2}  \big\|_{\mathcal{B}_{\mathbb{C}} ( \mathcal{Z}^N )} \leq 1 + 2^{\frac{N+1}{2}} (1 + \mathfrak{m}) \bigl( 1 +2 \sqrt{|t|} \bigr) , \quad \forall t \in \mathbb{R}. 
\end{equation}
\end{lem}
\begin{proof}
Fix  $u \in \mathcal{Z}^N$ for some $N \in \mathbb{N}_+$, then  $\partial_x^N u \in L^2 (\mathbb{R} ; \mathbb{C})$. Using formula \eqref{GrouPropexpitdelta}  and  the strong continuity of the unitary group $\bigl( e^{it \partial_x^2} \bigr)_{t \in \mathbb{R}}$ on $L^2(\mathbb{R}; \mathbb{C})$, we obtain that
\begin{equation}\label{StrContexpitLapL2}
 \sup_{ t \in \mathbb{R}}\big\| \partial_x^N   \bigl( e^{i  t  \partial_x^2} (u) - u \bigr)  \big\|_{L^2(\mathbb{R} ; \mathbb{C}) } \leq 2 
\big\| \partial_x^N u \big\|_{L^2(\mathbb{R}; \mathbb{C})} \quad \mathrm{and} \quad \lim_{t \to 0} \big\|     e^{i t \partial_x^2} ( \partial_x^N  u) - \partial_x^N  u  \big\|_{L^2(\mathbb{R} ; \mathbb{C}) }  = 0.
\end{equation}Fix $t \in \mathbb{R}$.  Using formulas \eqref{ContrL2dNnorHNnorm},  \eqref{GrouPropexpitdelta},  \eqref{estiL2EtDuubyu} and \eqref{StrContexpitLapL2}, we obtain that
\begin{equation}\label{EstiHNexpitLapu-u}
\begin{split}
\big\| e^{i t \partial_x^2} (  u) -   u  \big\|_{H^N(\mathbb{R} ; \mathbb{C}) }  
\leq & 2^{\frac{N-1}{2}} \bigl( \bigl( 2 + \sqrt{2} \bigr) \sqrt{|t|} \|\partial_x u\|_{L^2(\mathbb{R} ; \mathbb{C})}
+\big\| e^{i t \partial_x^2} ( \partial_x^N  u) - \partial_x^N  u  \big\|_{L^2(\mathbb{R} ; \mathbb{C}) } \bigr) \\
\leq & 2^{\frac{N+1}{2}} \bigl(1 + 2 \sqrt{|t|} \bigr) \|\partial_x^{N-1} u \|_{H^{N-1}(\mathbb{R}; \mathbb{C})}.
\end{split}
\end{equation}Then formula \eqref{expitLap-iduHNBd}  follows directly from \eqref{StrContexpitLapL2} and \eqref{EstiHNexpitLapu-u}. Recall that $H^N (\mathbb{R}; \mathbb{C}) \hookrightarrow \mathcal{Z}^N$ by \eqref{HNhookarrZN}.  So formula \eqref{BdZNexpitlapwrtt} is obtained by  \eqref{HNhookarrZN} and \eqref{EstiHNexpitLapu-u}. The strong continuity of  $\bigl( e^{it \partial_x^2} \bigr)_{t \in \mathbb{R}}$   on the Zhidkov space $\mathcal{Z}^N$ is a direct consequence of formulas \eqref{HNhookarrZN},  \eqref{expitLap-iduHNBd} and \eqref{BdZNexpitlapwrtt}. 
\end{proof}
  
\begin{lem}\label{GenSchGrlemZN}
For any $N \in \mathbb{N}_+$, the infinitesimal generator of the Schr\"odinger group $(e^{it\partial_x^2})_{t\in\mathbb{R}}$ on the Zhidkov space $\mathcal{Z}^N$ is $\bigl( i\partial_x^2, \mathrm{Dom}_{\mathcal{Z}^N}(i\partial_x^2) \bigr)$ with $\mathrm{Dom}_{\mathcal{Z}^N}(i\partial_x^2) = \mathcal{Z}^{N+2}$. Specifically, we have 
\begin{equation}\label{IniGenSchGrZNZN+2}
\mathcal{Z}^{N+2} = \bigl\{ u \in \mathcal{Z}^N : \exists\, w \in \mathcal{Z}^N \text{ such that } \lim_{t\to 0} \big\| \tfrac{1}{t}(e^{it\partial_x^2}u - u) - w \big\|_{\mathcal{Z}^N} = 0 \bigr\}.
\end{equation}
  Moreover, for any $\phi \in \mathcal{Z}^{N+1}$, the convergence holds in the  $H^{N-1}$-topology:
\begin{equation}\label{CVHNeitLapphiphi}
\lim_{t \to 0} \big\| \tfrac{1}{t} \bigl( e^{i t \partial_x^2} \phi - \phi \bigr) - i \partial_x^2 \phi \big\|_{H^{N-1} (\mathbb{R};\mathbb{C})} = 0.
\end{equation}
\end{lem} 
 
\begin{proof}
Assume that $u , w \in \mathcal{Z}^N$ such that 
$\lim_{t \to 0} \big\| \tfrac{1}{t} \bigl( e^{i  t \partial_x^2} (u) - u \bigr) - w \big\|_{\mathcal{Z}^N} = 0 $. Then \eqref{ZhidNormZn} yields that
\begin{equation}\label{eiLapdxuHN-1}
\big\| \tfrac{1}{t} \bigl( e^{i  t \partial_x^2} ( \partial_x  u) - \partial_x  u \bigr) - \partial_x  w  \big\|_{H^{N-1} (\mathbb{R} ; \mathbb{C})}  
\leq  \big\| \tfrac{1}{t} \bigl( e^{i  t \partial_x^2} (   u) -   u \bigr) -   w  \big\|_{\mathcal{Z}^N} \to 0, \quad \text{as}  \quad t \to 0.
\end{equation}By the Stone theorem (Theorem \uppercase\expandafter{\romannumeral8}.7 of \cite{ReedSimon1FA}), the infinitesimal generator of the strongly continuous unitary group $\bigl( e^{i t \partial_x^2 }\bigr)_{t \in \mathbb{R}}$ on the Sobolev space $H^{N-1} (\mathbb{R} ; \mathbb{C})$ is $\bigl( i \partial_x^2, \mathrm{Dom}_{H^{N-1}  } (i \partial_x^2) \bigr)$ with
\begin{equation}\label{DomiLapHN-1HN+1}
\begin{split}
 & \mathrm{Dom}_{H^{N-1}  } (i \partial_x^2) =     H^{N+1} (\mathbb{R} ; \mathbb{C}) \\
= &  \lbrace f \in H^{N-1} (\mathbb{R} ; \mathbb{C}) : \exists   g \in  H^{N-1} (\mathbb{R} ; \mathbb{C})  \quad  \text{such that} \quad 
\lim_{t \to 0} \big\| \tfrac{1}{t} \bigl( e^{i  t \partial_x^2} (f) - f \bigr) - g \big\|_{H^{N-1}} = 0  \rbrace.
\end{split}
\end{equation}We have $\partial_x u \in H^{N+1}(\mathbb{R} ; \mathbb{C})$ by  formulas 
\eqref{eiLapdxuHN-1} and
\eqref{DomiLapHN-1HN+1}. So $u \in \mathcal{Z}^{N+2}$ by  \eqref{ZhiSpdefSect2}. \\

\noindent Conversely, fix $\phi \in \mathcal{Z}^{N+1}$ and set $h : = \partial_x \phi \in H^N  ( \mathbb{R} ; \mathbb{C}) $. If $t \in \mathbb{R}$, then $e^{i t \partial_x^2} \phi - \phi \in L^2  $ by 
\eqref{estiL2EtDuubyu}.
For any $\xi \in \mathbb{R}$, we have
$e^{-it \xi^2} -1 = -it \xi^2 \int_0^1 e^{- i \mathfrak{s} t \xi^2 }\mathrm{d} \mathfrak{s}$. Since $\partial_x h \in H^{N-1} \subset L^2  ( \mathbb{R} ; \mathbb{C})$,  we have
\begin{equation}\label{eiLapphiLtwocv}
2\pi \big\|  \tfrac{1}{t} \bigl( e^{i  t \partial_x^2} ( \phi ) - \phi \bigr) -  i \partial_x^2 \phi   \big\|_{L^2 (\mathbb{R} ; \mathbb{C})}^2 = \int_{\mathbb{R}} \Big| \xi \widehat{h}(\xi) \Bigl(1 - \int_0^1 e^{- i \mathfrak{s} t \xi^2 }\mathrm{d} \mathfrak{s} \Bigr) \Big|^2 \mathrm{d} \xi \to 0,
\end{equation}as $t \to 0$,  by the dominated convergence theorem.
Then  \eqref{CVHNeitLapphiphi}  holds for $N=1$ and all $\phi \in \mathcal{Z}^2$.      
  Assume that $N \geq 2$.   Using the Stone theorem and formula \eqref{DomiLapHN-1HN+1}, we obtain that
\begin{equation}\label{CveiLapphiHN-1}
\big\|\partial_x \bigl(  \tfrac{1}{t} \bigl( e^{i  t \partial_x^2} ( \phi ) - \phi \bigr) -  i \partial_x^2 \phi  \bigr) \big\|_{H^{N-2} (\mathbb{R} ; \mathbb{C})}  
\leq  \big\| \tfrac{1}{t} \bigl( e^{i  t \partial_x^2} ( h ) -   h \bigr) -  i \partial_x^2 h  \big\|_{H^{N-2}} \to 0, \quad \text{as}  \quad t \to 0.
\end{equation}So 
\eqref{CVHNeitLapphiphi} follows directly from    
\eqref{ContrL2dNnorHNnorm}, 
\eqref{CveiLapphiHN-1} and   \eqref{eiLapphiLtwocv}. The proof is completed by \eqref{HNhookarrZN}.
\end{proof}  

\begin{rem}\label{GalloThmSchGr}
We refer to Theorems 3.1 and 3.2 of   \cite{GalloSchGrZhidSp}
for alternative proofs of Lemmas $\ref{StrCSchGrZhidlem}$ and $\ref{GenSchGrlemZN}$.
\end{rem}

\noindent In order to describe the nonzero boundary condition at infinity for solutions of the focusing CMdNLS equation \eqref{CMdNLSfocu}, we introduce the following subset of the Zhidkov space $\mathcal{Z}^n$:
\begin{equation}\label{calMnmetrSp}
\mathcal{M}^n := \{ u \in \mathcal{Z}^n : |u|^2 - 1 \in L^2(\mathbb{R};\mathbb{C}) \}, \qquad \forall n \in \mathbb{N}_+.
\end{equation}Fix $n \in \mathbb{N}_+$. We recall the classical metric for the set  $ \mathcal{M}^n$, which  is defined as follows: 
\begin{equation}\label{sfdcalMnmetricDef}
\mathsf{d}_{\mathcal{M}^n} ( u , v ) : = \| u -v \|_{\mathcal{Z}^n} + \| |u|^2 - |v|^2 \|_{L^2 (\mathbb{R} ; \mathbb{C})} .
\end{equation}We recall also some basic properties about the metric space $\bigl( \mathcal{M}^n,  \mathsf{d}_{\mathcal{M}^n}  \bigr)$.
\begin{lem}
For any $n \in \mathbb{N}_+$, $\bigl( \mathcal{M}^n,  \mathsf{d}_{\mathcal{M}^n}  \bigr)$ is a complete metric space. Moreover, we have
\begin{equation}\label{calMnbyHndef2}
\mathcal{M}^n = \{ u \in \mathcal{Z}^n : |u|^2 -1 \in H^n (\mathbb{R} ; \mathbb{C})\}.
\end{equation}The inclusions $\iota_n : u \in \mathcal{M}^{n +1} \mapsto u \in \mathcal{M}^n$ and $\mathfrak{i}_n : u \in \mathcal{M}^n  \mapsto u \in \mathcal{Z}^n$ are both Lipschitz continuous.
\end{lem}
\begin{proof}
Formula \eqref{calMnbyHndef2} follows from \eqref{uvZn<uZnvZn} and
\eqref{ZNcapL2=HN}.
\end{proof}
\begin{rem}
 When $u \in \mathcal{M}^n$ for some $n \in \mathbb{N}_+$, the condition $|u|^2 - 1 \in H^n \subset H^1(\mathbb{R} ; \mathbb{C})$ guarantees that $|u(x)| \to 1$ as $|x| \to \infty$, which is the natural setting for dark solitons.
\end{rem}

\begin{lem}\label{rhoHrhosqiden}
For any $\rho \in L^3 (\mathbb{R}; \mathbb{R})$, we have $\int_{\mathbb{R}} \rho^3
= 3 \int_{\mathbb{R}}\rho |\mathbb{H}\rho|^2
= 3  \int_{\mathbb{R}}\rho |\Pi_{\geq 0} \rho|^2 \in \mathbb{R}$. 
\end{lem}
\begin{proof}
The case $\rho \in H^1 (\mathbb{R}; \mathbb{R})$ is proved in Appendix C (Page 53) in \cite{GerLenz22CPAM}.
Since   $H^1(\mathbb{R};\mathbb{R})$ is dense in $L^3(\mathbb{R};\mathbb{R})$ and both the Hilbert transform $\mathbb{H}$ and the Szeg\H{o} projector $\Pi_{\geq 0}$  are bounded linear operators on $L^3(\mathbb{R};\mathbb{C})$, the identity extends to all $\rho\in L^3(\mathbb{R};\mathbb{R})$ by a standard continuity argument.
\end{proof}

\begin{lem}
For any $v \in   L^{\infty}(\mathbb{R}; \mathbb{C})$ such that $|v|^2 - 1 \in L^2 (\mathbb{R} ; \mathbb{C})$, 
we have
\begin{equation}\label{vhvsq-1L2norm}
\|v \mathbb{H} (|v|^2 -1)\|_{L^2 (\mathbb{R} ; \mathbb{C})}^2 = \frac{1}{3} \int_{\mathbb{R}} (|v|^2 -1)^3 + \| |v|^2 -1 \|_{L^2 (\mathbb{R} ; \mathbb{C})}^2 \geq \tfrac{2}{3}  \| |v|^2 -1 \|_{L^2 (\mathbb{R} ; \mathbb{C})}^2.
\end{equation}For any $u \in \mathcal{M}^1$, the following identity holds:
\begin{equation}\label{Dsqusq-1inprod}
2 \operatorname{Re}  \langle \partial_x u, u \mathbb{H} (|u|^2 -1) \rangle_{L^2 (\mathbb{R} ; \mathbb{C})} = - \langle |u|^2 -1 ,   |D|  (|u|^2 -1)   \rangle_{L^2 (\mathbb{R} ; \mathbb{C})} = - \||D|^{\frac{1}{2}} (|u|^2 -1) \|_{L^2 (\mathbb{R}; \mathbb{C})}^2.
\end{equation}
\end{lem}
\begin{proof}
Both identities appear in the work of Chen \cite{Xchenccm25ExpFor}; see the proof of Lemma~3.5 in \cite{Xchenccm25ExpFor} for their detailed derivation. 
For completeness, we note that \eqref{vhvsq-1L2norm} follows by applying Lemma~\ref{rhoHrhosqiden} to $\rho = |v|^2 - 1$ together with the $L^2$-unitary property of the Hilbert transform $\mathbb{H}$.
\end{proof} 

\begin{lem}\label{MnaddHninMnlem}
For any $n \in \mathbb{N}_+$, we have $\mathcal{M}^n + H^n (\mathbb{R} ; \mathbb{C}) \subset \mathcal{M}^n$. The group of operators $\bigl( e^{ it \partial_x^2 } \bigr)_{t \in \mathbb{R}}$ is strongly continuous on the complete  metric space $\bigl( \mathcal{M}^n , \mathsf{d}_{\mathcal{M}^n} \bigr)$. Specifically, for any $u_0 \in \mathcal{M}^n$,  the mapping 
 $u_{\mathrm{L}} : t \in \mathbb{R} \mapsto u_{\mathrm{L}} ( t ) =  e^{ it \partial_x^2 }(u_0) \in \bigl( \mathcal{M}^n , \mathsf{d}_{\mathcal{M}^n} \bigr)$ is continuous. 
\end{lem}  
\begin{proof}
The case $n=2$ is proved in Lemmas~2.2 and 2.4 of \cite{Xchenccm25ExpFor}. 
 The  general case follows the same argument and formulas \eqref{HnZnsubsetHn}, 
\eqref{estiL2EtDuubyu} and 
\eqref{expitLap-iduHNBd}.
\end{proof}

\noindent For any $u \in L^{\infty} (\mathbb{R} ; \mathbb{C})$, the auxiliary function $\mathsf{c}_u :  \mathbb{R} \to \mathbb{C}$ is defined as  
\begin{equation}\label{sfcuFuncDef}
\mathsf{c}_u (x) : =  \exp \bigl( -\frac{i}{2} \int_0^x (|u(y)|^2 -1) \mathrm{d} y \bigr), \qquad \forall x \in \mathbb{R}.
\end{equation}Then $\mathsf{c}_u \in W^{1, \infty} (\mathbb{R} ; \mathbb{C})$ and   $|\mathsf{c}_u (x) | = 1$,  $\forall  x \in \mathbb{R}$. We have the following lemma.
\begin{lem}\label{lemsfcuinZn+1}
For any $n \in \mathbb{N}_+$ and $u \in \mathcal{M}^n$, we have $\mathsf{c}_u \in \mathcal{Z}^{n+1}$. 
\end{lem}
\begin{proof}
If $u \in \mathcal{M}^1$, then $|u|^2 -1 \in  H^1  (\mathbb{R} ; \mathbb{R})$ by 
\eqref{calMnbyHndef2}. Formula \eqref{sfcuFuncDef} yields that
\begin{equation}\label{dxsfcufor}
\partial_x \mathsf{c}_u = -  \tfrac{i}{2} (| u |^2 - 1) \mathsf{c}_u \in L^2 \cap L^{\infty} (\mathbb{R} ; \mathbb{C}) \Longrightarrow \mathsf{c}_u \in \mathcal{Z}^1 \cap W^{1, \infty} (\mathbb{R} ; \mathbb{C}).
\end{equation}Using \eqref{dxsfcufor} and \eqref{HnZnsubsetHn}, we obtain that $\partial_x \mathsf{c}_u \in H^1 \cdot \mathcal{Z}^1 
\subset H^1 (\mathbb{R} ; \mathbb{C})$. Then $\mathsf{c}_u \in \mathcal{Z}^2$.  
 This proves the case $n=1$. Assuming the statement holds for some $n \ge 1$, the case $n+1$ follows by the same reasoning, using \eqref{dxsfcufor}, \eqref{calMnbyHndef2}, and \eqref{HnZnsubsetHn} at order $n+1$. The induction is complete.
\end{proof}

\noindent Fix $u \in L^{\infty}(\mathbb{R} ; \mathbb{C})$. Using the   function $\mathsf{c}_u$ introduced in \eqref{sfcuFuncDef}, we define the gauge transform by 
\begin{equation}\label{GauTrdef}
\mathcal{G}u(x): = u(x) \mathsf{c}_u(x) = u(x) \exp \bigl( -\frac{i}{2} \int_0^x (|u(y)|^2 -1) \mathrm{d} y \bigr), \qquad \forall x \in \mathbb{R}.
\end{equation}The mapping $\mathcal{G} : L^{\infty}(\mathbb{R} ; \mathbb{C}) \to L^{\infty}(\mathbb{R} ; \mathbb{C})$ is called the \textit{gauge transform}. 
 
\begin{lem}\label{BijLemGauTr}
The gauge transform $\mathcal{G} : L^{\infty}(\mathbb{R} ; \mathbb{C}) \to L^{\infty}(\mathbb{R} ; \mathbb{C})$  is  bijective and we have
\begin{equation}\label{invcalGuIsoLinfty}
\mathcal{G}^{-1} (u)  = u \overline{\mathsf{c}}_u  , \quad  
|\mathcal{G}  (u)   | =  | u   | =
  |\mathcal{G}^{-1} (u)   | \in L^{\infty} (\mathbb{R} ; \mathbb{C}), \quad \forall u \in L^{\infty} (\mathbb{R} ; \mathbb{C}).
\end{equation}Consequently, for any $n \in \mathbb{N}_+$, the restriction $\mathcal{G} |_{\mathcal{M}^n} : \mathcal{M}^n \to \mathcal{M}^n$ is bijective.
\end{lem}
\begin{proof}
Formula \eqref{invcalGuIsoLinfty} follows from \eqref{GauTrdef}. 
The last statement is a direct consequence of  formulas \eqref{calMnmetrSp},  \eqref{invcalGuIsoLinfty} and Lemma $\ref{lemsfcuinZn+1}$.\end{proof}

\subsection{Poisson Integrals and Difference Quotients}\label{SubsPoiInDiQuo}

\noindent The Poisson kernel on the upper half-plane $\mathbb{C}_+$ is given by
\begin{equation}\label{PoiKerS2}
\mathfrak{P}_y(x ) =\tfrac{ y }{ \pi ( x^2 + y^2 ) } ,  \qquad \forall x \in \mathbb{R}, \qquad \forall y >0 .
\end{equation}
For any $p \in [1, +\infty]$ and $f \in L^p(\mathbb{R}; \mathbb{C})$, its Poisson integral is defined by 
\begin{equation}\label{PoiIntDefSect2}
\underline{f}(z) = \mathbf{P}_z (f)  := \frac{\operatorname{Im} z}{\pi} \int_{\mathbb{R}} \frac{f(\mathfrak{t})}{|\mathfrak{t} - z|^2} \, \mathrm{d}\mathfrak{t} = \int_{\mathbb{R}} \mathfrak{P}_y(x - \mathfrak{t}) f(\mathfrak{t}) \, \mathrm{d}\mathfrak{t}, \qquad \forall z = x + y i \in \mathbb{C}_+.
\end{equation}Fix $z \in \mathbb{C}_+$. 
Using H\"older's inequality, we obtain that
\begin{equation}\label{PoiIntBddLinFunc}
 |\mathbf{P}_z (f) | \leq (\operatorname{Im} z)^{-\frac{1}{p}} \|\mathfrak{P}_1\|_{L^{\frac{p}{p-1}}(\mathbb{R}; \mathbb{R})} \|f\|_{L^p(\mathbb{R}; \mathbb{C})}, \qquad \forall f \in L^p(\mathbb{R}; \mathbb{C}),
\end{equation}where $\mathfrak{P}_1(x) = \frac{1}{\pi (1 + x^2)}$.
Thus $\mathbf{P}_z : f \in L^p(\mathbb{R} ; \mathbb{C}) \mapsto \underline{f}(z) \in \mathbb{C}$ is a bounded linear functional, and  
\begin{equation}\label{PzL2BdPz1is1}
\mathbf{P}_z(\mathbf{1}) = 1, \qquad
\|\mathbf{P}_z\|_{(L^{\infty}(\mathbb{R} ; \mathbb{C}))^*} = 1, \qquad
\|\mathbf{P}_z\|_{(L^2(\mathbb{R} ; \mathbb{C}))^*} \leq \bigl(2\pi \operatorname{Im} z\bigr)^{-\frac{1}{2}},
\end{equation}
where $\mathbf{1} \in L^{\infty}(\mathbb{R}; \mathbb{C})$ denotes the constant function $\mathbf{1}(x) = 1$ for all $x \in \mathbb{R}$. 
The Szeg\H{o} projector $\Pi_{\geq 0} : L^2 (\mathbb{R} ; \mathbb{C}) \to L^2_+$ is given by \eqref{SzegoProjS2}. 
  For any $\epsilon > 0$, we introduce the auxiliary function
\begin{equation}\label{chiepFor}
\chi_{\epsilon}(x) := \tfrac{2\pi}{\epsilon} \Pi_{\geq 0} (\mathfrak{P}_{\epsilon^{-1}})(x) = (1 - i\epsilon x)^{-1}, \quad \forall x \in \mathbb{R}.
\end{equation}
Clearly $\chi_{\epsilon} \in H^{\infty}(\mathbb{R}; \mathbb{C}) \cap L^2_+$ for every $\epsilon > 0$, and $\chi_{\epsilon}$ converges pointwise to the constant function $\mathbf{1}$ as $\epsilon \to 0^+$. 
The following lemma collects its basic approximation properties. 
\begin{lem}\label{CVchiep-1ghu}
Fix $p \in [1, +\infty)$ and $N \in \mathbb{N}_+$. For any $f \in L^p(\mathbb{R}; \mathbb{C})$ and $h \in H^N(\mathbb{R}; \mathbb{C})$, we have
\begin{equation}\label{chiephCVtohinHpFor}
\lim_{\epsilon \to 0^+} \| \partial_x^N \chi_{\epsilon} \|_{L^2(\mathbb{R} ; \mathbb{C})} = 0, 
\quad 
\lim_{\epsilon \to 0^+} \|f \chi_{\epsilon} - f\|_{L^p(\mathbb{R} ; \mathbb{C})} = 0, \quad
\lim_{\epsilon \to 0^+} \|h \chi_{\epsilon} - h\|_{H^N(\mathbb{R} ; \mathbb{C})} = 0 
.
\end{equation}
If $u \in L^{\infty}(\mathbb{R}; \mathbb{C})$, then $u \chi_{\epsilon} \to u$ pointwise as $\epsilon \to 0^+$. Moreover, for any $z \in \mathbb{C}_+$, we have
\begin{equation}\label{CvPoiIntuchieptou}
\lim_{\epsilon \to 0^+} \mathbf{P}_z(u \chi_{\epsilon}) = \mathbf{P}_z(u) = \underline{u}(z) \in \mathbb{C}.
\end{equation}
\end{lem}
\begin{rem}\label{LinNorchiep-1}
If $u \in L^{\infty}(\mathbb{R}; \mathbb{C})$, the convergence $u\chi_{\epsilon} \to u$ is uniform on every compact subset of $\mathbb{R}$ but not uniform on $\mathbb{R}$ in general. For example, if $u = \mathbf{1}$, then $\sup_{x \in \mathbb{R}} |\chi_{\epsilon}(x) - 1| = 1$, for every $\epsilon > 0$.
\end{rem} 
\begin{proof}[Proof of Lemma $\ref{CVchiep-1ghu}$]
For any $N \in \mathbb{N}_+$, $\epsilon >0$ and $x \in \mathbb{R}$, a direct computation shows that
\begin{equation}\label{dxNchiepL2norm}
1 - \chi_{\epsilon}(x) = \tfrac{-i\epsilon x}{1 - i \epsilon x}, \qquad
\partial_x^N \chi_{\epsilon} (x) = \tfrac{(i\epsilon)^N N !}{(1  - i \epsilon x)^{N+1}}, \qquad \|\partial_x^N \chi_{\epsilon} \|_{L^2( \mathbb{R} ; \mathbb{C})}^2 = \tfrac{\pi (2N)!}{4^N} \epsilon^{2N -1} \to 0, 
\end{equation}as $\epsilon \to 0^+$. For any function $g : \mathbb{R} \to \mathbb{C}$, the first identity in \eqref{dxNchiepL2norm} implies that $\| \chi_{\epsilon} -1 \|_{L^{\infty}(\mathbb{R}; \mathbb{C})} =1$,
\begin{equation}\label{EstCVgchiep-g}
\lim_{\epsilon \to 0^+} |g (x) \chi_{\epsilon} (x)  - g(x) |= 0 \quad \mathrm{and} \quad \sup_{\epsilon >0} |g (x) \chi_{\epsilon} (x)  - g(x) | \leq |g(x)|, \quad \forall x \in \mathbb{R}.
\end{equation}Fix $z \in \mathbb{C}_+$, we set $\boldsymbol{\rho}_z (x) := |x-z|^{-2} \geq 0$, $\forall x \in \mathbb{R}$. For any $u\in L^{\infty}(\mathbb{R}; \mathbb{R})$, formula  \eqref{PoiIntDefSect2} yields that
\begin{equation}\label{PoiIntdefbyrhoz}
\underline{u}(z ) = \mathbf{P}_z (u) = \frac{\operatorname{Im} z}{\pi} \int_{\mathbb{R}} u(x) \boldsymbol{\rho}_z (x) \mathrm{d}x \quad \mathrm{with} \quad \boldsymbol{\rho}_z \in L^1 \cap  L^{\infty} \cap   H^{\infty}(\mathbb{R}; \mathbb{R}).
\end{equation}Then formulas \eqref{CvPoiIntuchieptou} and the second identity in \eqref{chiephCVtohinHpFor} are direct consequences of \eqref{EstCVgchiep-g}, \eqref{PoiIntdefbyrhoz} and the dominated convergence theorem.
 Fix $h \in H^N(\mathbb{R}; \mathbb{C})$. For any $\epsilon >0$ and $k \in \{1,2, \cdots, N\}$,  we have $h \chi_{\epsilon}   \in H^N(\mathbb{R}; \mathbb{C})$ and $\partial_x^{N-k} h \in H^1 \hookrightarrow L^{\infty} (\mathbb{R}; \mathbb{C})$. Using  the 
Leibniz formula \eqref{IndAssumpn=NZNHN}, we obtain that
\begin{equation}\label{LeiFordxNhchiep-h}
 \partial_x^N \bigl(h \chi_{\epsilon} - h \bigr) = \bigl( \chi_{\epsilon} - 1\bigr) \partial_x^N h
+ \sum_{k=1}^N \frac{N!}{k! (N-k)!} \partial_x^k \chi_{\epsilon}  \partial_x^{N-k} h \to 0  \quad \text{in} \quad L^2(\mathbb{R}; \mathbb{C}), 
\end{equation}as $\epsilon \to 0^+$, by the first two identities of \eqref{chiephCVtohinHpFor}.
Finally, we have $\lim_{\epsilon \to 0^+} \|h \chi_{\epsilon} - h\|_{H^N(\mathbb{R} ; \mathbb{C})} = 0 $ by formulas  \eqref{ContrL2dNnorHNnorm}, \eqref{LeiFordxNhchiep-h} and the second identity of \eqref{chiephCVtohinHpFor}.
\end{proof}

\begin{cor}\label{sfbmDefinCor}
For any $m \in \mathbb{N}$, we set $\mathsf{b}_m: = 2^{\frac{m - 1}{2}} \bigl( \sqrt{\frac{\pi}{2}} + \frac{\sqrt{\pi (2 m + 2)!}}{2^{m+1}} \bigr) >0$. Then we have
\begin{equation}\label{Cordxchiepto0}
\lim_{\epsilon \to 0^+}    \| \partial_x \chi_{\epsilon}   \|_{H^m (\mathbb{R}; \mathbb{C})} =0  ,\quad
\sup_{ 0 < \epsilon \leq 1} \big\| \partial_x \chi_{\epsilon} \big\|_{H^m (\mathbb{R}; \mathbb{C})} \leq  \mathsf{b}_m , \quad \sup_{ 0 < \epsilon \leq 1} \big\| \chi_{\epsilon} \big\|_{\mathcal{Z}^{m+1}} \leq 1 + \mathsf{b}_m.
\end{equation}
\end{cor}

\begin{proof}
This is a direct consequence of formulas \eqref{ContrL2dNnorHNnorm}, \eqref{ZhidNormZn}      and 
  \eqref{dxNchiepL2norm}.
\end{proof}
  
\begin{lem}\label{SchrGronchieplem}
For any $t \in \mathbb{R}$, we have $\lim_{\epsilon \to 0^+} \|e^{it \partial_x^2} \chi_{\epsilon} - \chi_{\epsilon}\|_{L^2_+} = 0$.
\end{lem} 

\begin{proof}
The Schr\"odinger group $(e^{it \partial_x^2})_{t \in \mathbb{R}}$ is a strongly continuous unitary group on $L^2 (\mathbb{R} ; \mathbb{C})$. Its infinitesimal generator is given by $ \bigl( i \partial_x^2, \mathrm{Dom}_{L^2 (\mathbb{R} ; \mathbb{C})} (i \partial_x^2 )\bigr)$ with  $\mathrm{Dom}_{L^2 (\mathbb{R} ; \mathbb{C})} (i \partial_x^2  ) = H^2(\mathbb{R} ; \mathbb{C})$. Since $\chi_{\epsilon} \in L^2_+ \cap H^{\infty} (\mathbb{R} ; \mathbb{C})$, $\forall \epsilon >0$, the Newton--Leibniz formula for Bochner integral implies that
\begin{equation*}
e^{i t \partial_x^2} \chi_{\epsilon}  - \chi_{\epsilon}  = i \int_0^t e^{i \mathfrak{s} \partial_x^2} ( \partial_x^2 \chi_{\epsilon})  \mathrm{d}\mathfrak{s} \in L^2_+ , \qquad \forall t \in \mathbb{R}.
\end{equation*}The proof is completed by \eqref{dxNchiepL2norm}, \eqref{SchrGrUnitHs}  and Theorem 3.29 of  Rudin \cite{RudinFA}.
\end{proof}

\noindent For any $\epsilon>0$, $t \in \mathbb{R}$ and $u \in \mathcal{Z}^1$, we define three auxiliary functions:
\begin{equation}\label{fghtepuDefFor}
\mathrm{f}_{\epsilon}^u : = u \chi_{\epsilon} - u \in \mathcal{Z}^1, \qquad \mathrm{g}_{t, \epsilon}^u: = e^{i t \partial_x^2} \bigl( \mathrm{f}_{\epsilon}^u \bigr) - \mathrm{f}_{\epsilon}^u      \in H^1 (\mathbb{R}; \mathbb{C}),  \qquad \mathrm{h}_{t,  \epsilon}^u : = e^{i t \partial_x^2} \bigl( \mathrm{f}_{\epsilon}^u \bigr) \in \mathcal{Z}^1.
\end{equation}These functions satisfy the following lemmas. 

\begin{lem}\label{Lemfepu'to0AN}
Let $N \in \mathbb{N}_+$ and $\mathrm{A}_N:= (1 + \mathfrak{m}) \sum_{k =0 }^N \frac{N!}{k! (N-k)!} \frac{\sqrt{\pi (2k)!}}{2^k} >0 $, where $\mathfrak{m}$ is given by \eqref{InffH1overfLinfty}.
If $u \in \mathcal{Z}^N$, then $\mathrm{f}_{\epsilon}^u \in \mathcal{Z}^N$ and $\mathrm{g}_{t, \epsilon}^u  \in H^N (\mathbb{R}; \mathbb{C})$, for all $t \in \mathbb{R}$ and $\epsilon > 0$. Moreover, we have
\begin{equation}\label{gteputo0andBd}
\lim_{\epsilon \to 0^+} \big\| \partial_x \mathrm{f}_{\epsilon}^u \big\|_{H^{N-1} } = 0, \quad 
 \lim_{\epsilon \to 0^+} \big\|   \mathrm{g}_{t, \epsilon}^u \big\|_{H^N } = 0, \quad \sup_{0 < \epsilon \leq 1} \big\|   \mathrm{g}_{t, \epsilon}^u \big\|_{H^N } \leq 2^{\frac{N+1}{2}} \bigl( \mathrm{A}_N + \sqrt{2 \pi |t|}\bigr) \|u\|_{\mathcal{Z}^N}.
\end{equation}
\end{lem}
\begin{proof}
Fix $N \in \mathbb{N}_+$ and $u \in \mathcal{Z}^N \hookrightarrow \mathcal{Z}^1$, We have $\partial_x \mathrm{f}_{\epsilon}^u = \bigl( \chi_{\epsilon } -1 \bigr) \partial_x u + u \partial_x \chi_{\epsilon} \in H^{N-1} (\mathbb{R}; \mathbb{C})$ and
\begin{equation}\label{dxfeuto0inHN-1Pr}
\big\|\partial_x \mathrm{f}_{\epsilon}^u \big\|_{H^{N-1} (\mathbb{R} ; \mathbb{C})} \leq \|(\chi_{\epsilon} -1 ) \partial_x u\|_{H^{N-1} (\mathbb{R} ; \mathbb{C})} + (1+ C_{N-1} \mathbf{1}_{N \geq 2}) \|u\|_{\mathcal{Z}^N} \|\partial_x \chi_{\epsilon} \|_{H^{N-1}(\mathbb{R}; \mathbb{C})} \to 0,
\end{equation}as $\epsilon \to 0^+$, by   \eqref{ufHn<uZnfHn},
\eqref{chiephCVtohinHpFor}  and 
   \eqref{Cordxchiepto0}. Using \eqref{estiL2EtDuubyu}, \eqref{fghtepuDefFor} and \eqref{dxfeuto0inHN-1Pr},  we obtain that
\begin{equation}\label{gtepuL2normto0}
\big\|\mathrm{g}_{t, \epsilon}^u\big\|_{L^2(\mathbb{R}; \mathbb{C})} = \big\|  e^{i t \partial_x^2} \bigl( \mathrm{f}_{\epsilon}^u \bigr) - \mathrm{f}_{\epsilon}^u  \big\|_{L^2(\mathbb{R}; \mathbb{C})} \leq (2 + \sqrt{2}) \sqrt{|t|}  \big\| \partial_x \mathrm{f}_{ \epsilon}^u\big\|_{L^2(\mathbb{R}; \mathbb{C})} \to 0,  \qquad \forall t \in \mathbb{R}, 
\end{equation}as $\epsilon \to 0^+$. Since $\bigl(e^{i t \partial_x^2}\bigr)_{t \in \mathbb{R}}$ is a strongly continuous unitary group on $L^2 (\mathbb{R}; \mathbb{C})$, we have
\begin{equation}\label{dxNgsmal2dxNf}
\big\| \partial_x^N  \mathrm{g}_{t, \epsilon}^u\big\|_{L^2(\mathbb{R}; \mathbb{C})} \leq 2 \big\| \partial_x^N  \mathrm{f}_{ \epsilon}^u\big\|_{L^2(\mathbb{R}; \mathbb{C})} \quad \mathrm{and} \quad
 \big\| \partial_x^N  \mathrm{g}_{t, \epsilon}^u\big\|_{L^2(\mathbb{R}; \mathbb{C})} = \big\|   e^{i t \partial_x^2} \bigl(\partial_x^N  \mathrm{f}_{\epsilon}^u \bigr) - \partial_x^N  \mathrm{f}_{\epsilon}^u \big\|_{L^2(\mathbb{R}; \mathbb{C})} \to 0,
\end{equation}as $\epsilon \to 0^+$. So    $\lim_{\epsilon \to 0^+} \big\|   \mathrm{g}_{t, \epsilon}^u \big\|_{H^N(\mathbb{R}; \mathbb{C})} = 0$ by
\eqref{ContrL2dNnorHNnorm},  \eqref{gtepuL2normto0} and  \eqref{dxNgsmal2dxNf}. Formula   \eqref{dxNchiepL2norm} implies that $\| \chi_{\epsilon} -1 \|_{L^{\infty}(\mathbb{R} ; \mathbb{C})} =1$  and  $\|\partial_x \chi_{\epsilon} \|_{L^2}^2 = \frac{\pi}{2} \epsilon$. By \eqref{ZhidNormZn} and
\eqref{fghtepuDefFor}, it follows that
\begin{equation}\label{NormL2fepu'}
 \big\| \partial_x  \mathrm{f}_{ \epsilon}^u \big\|_{L^2 (\mathbb{R}; \mathbb{C})} \leq \|(\chi_{\epsilon} -1 ) \partial_x u\|_{L^2(\mathbb{R}; \mathbb{C})} + \|u \|_{L^{\infty} (\mathbb{R}; \mathbb{C})} \| \partial_x \chi_{\epsilon} \|_{L^2_+} \leq \sqrt{ \tfrac{\pi}{2} } \|u\|_{\mathcal{Z}^1}, \quad \forall \epsilon \in (0, 1].
\end{equation}According to \eqref{ZhidNormZn}, \eqref{InffH1overfLinfty}, 
\eqref{dxNchiepL2norm}, 
\eqref{fghtepuDefFor} and the  Leibniz formula \eqref{IndAssumpn=NZNHN},  for any $\epsilon \in (0, 1]$, we have
\begin{equation}\label{EstidxNfepuL2norm}
\begin{split}
\big\| \partial_x^N   \mathrm{f}_{ \epsilon}^u \big\|_{L^2  } 
\leq & \| u \|_{L^{\infty}} \|\partial_x^N \chi_{\epsilon} \|_{L^2_+} + \sum_{n=1}^N \tfrac{N!}{k! (N-k)!} \| \partial_x^{N-k} u \|_{L^{\infty}} \|\partial_x^k \chi_{\epsilon} \|_{L^2_+} + \|\partial_x^N u\|_{L^2 }\| \chi_{\epsilon} -1  \|_{L^{\infty}} \\
\leq & \bigl( \| u \|_{L^{\infty} (\mathbb{R}; \mathbb{C})} +  \| \partial_x u \|_{H^{N-1} (\mathbb{R}; \mathbb{C})} \bigr) (1 + \mathfrak{m}) \sum_{k =0 }^N \tfrac{N!}{k! (N-k)!} \tfrac{\sqrt{\pi (2k)!}}{2^k} = \mathrm{A}_N \| u \|_{\mathcal{Z}^N} .
\end{split}
\end{equation}Plugging \eqref{NormL2fepu'} into \eqref{gtepuL2normto0} and \eqref{EstidxNfepuL2norm} into \eqref{dxNgsmal2dxNf}, we obtain that
\begin{equation}\label{supep<1gtepu}
\sup_{0 < \epsilon \leq 1} \big\|\mathrm{g}_{t, \epsilon}^u\big\|_{L^2(\mathbb{R}; \mathbb{C})} \leq (1 + \sqrt{2}) \sqrt{\pi |t|} \|u\|_{\mathcal{Z}^1}, \qquad
\sup_{0 < \epsilon \leq 1} \big\| \partial_x^N \mathrm{g}_{t, \epsilon}^u\big\|_{L^2(\mathbb{R}; \mathbb{C})} \leq 2 \mathrm{A}_N \|u\|_{\mathcal{Z}^N}.
\end{equation}Then the second estimate in formula \eqref{gteputo0andBd}
follows directly from \eqref{ContrL2dNnorHNnorm} and \eqref{supep<1gtepu}.
\end{proof}

\begin{lem}
For any $t \in \mathbb{R}$, $\epsilon >0$ and $u \in \mathcal{Z}^1$, let $\mathrm{h}_{t,  \epsilon}^u   \in \mathcal{Z}^1$ be defined in \eqref{fghtepuDefFor}. Then we have 
\begin{equation}\label{htepuphito0inL2}
\lim_{\epsilon \to 0^+} \mathbf{P}_z   e^{i t \partial_x^2} (u  \chi_{\epsilon} - u)
= \lim_{\epsilon \to 0^+} \mathbf{P}_z \bigl( \mathrm{h}_{t, \epsilon}^u \bigr) = 0,
 \quad 
\lim_{\epsilon \to 0^+} \big\|\phi \mathrm{h}_{t, \epsilon}^u\big\|_{L^2(\mathbb{R}; \mathbb{C})} =0, \quad \forall \phi \in L^2(\mathbb{R}; \mathbb{C}) .
\end{equation}
\end{lem}
 
\begin{proof}
Fix $t \in \mathbb{R}$ and $u \in \mathcal{Z}^1$.   Then \eqref{fghtepuDefFor} implies that $\mathrm{h}_{t, \epsilon}^u = \mathrm{f}_{\epsilon}^u + \mathrm{g}_{t, \epsilon}^u$, $\forall \epsilon >0$. Let $\mathfrak{m} > 0$ be defined in \eqref{InffH1overfLinfty}. For any $x \in \mathbb{R}$, using \eqref{dxNchiepL2norm}  and the second identity in \eqref{gteputo0andBd}, we obtain that
\begin{equation}\label{PtCVhteputo0}
\big|\mathrm{h}_{t, \epsilon}^u (x) \big| \leq \|u\|_{L^{\infty} (\mathbb{R}; \mathbb{C})} |\chi_{\epsilon}(x) -1 | +   \big\| \mathrm{g}_{t, \epsilon}^u\big\|_{L^{\infty} (\mathbb{R}; \mathbb{C})} \leq \|u\|_{\mathcal{Z}^1} \tfrac{\epsilon |x|}{\sqrt{1 + \epsilon^2 x^2}} + \mathfrak{m} \big\| \mathrm{g}_{t, \epsilon}^u\big\|_{H^1  }  \to 0,
\end{equation}as $\epsilon \to 0^+$.  By formula \eqref{PtCVhteputo0} and the third estimate  in \eqref{gteputo0andBd}, it follows that 
\begin{equation}\label{EsUnhutep<1Linft}
\sup_{ 0 < \epsilon \leq 1} \big\|\mathrm{h}_{t, \epsilon}^u \big\|_{L^{\infty}(\mathbb{R}; \mathbb{C})} \leq \bigl( 1 + 2 \mathfrak{m} \bigl(\mathrm{A}_1 + \sqrt{2 \pi |t|} \bigr)\bigr) \|u\|_{\mathcal{Z}^1} ,
\end{equation}where $\mathrm{A}_1>0$ is defined in Lemma $\ref{Lemfepu'to0AN}$. According to formulas \eqref{PoiIntdefbyrhoz},  \eqref{PtCVhteputo0} and \eqref{EsUnhutep<1Linft}, formula \eqref{htepuphito0inL2} is a direct consequence of the dominated convergence theorem.
\end{proof}
 
\begin{cor}\label{htepuinZNCor}
If $N \in \mathbb{N}_+$ and $u \in \mathcal{Z}^N$, then the function $\mathrm{h}_{t,\epsilon}^u$ given by   \eqref{fghtepuDefFor} belongs to the Zhidkov space $\mathcal{Z}^N$ and satisfies  
\begin{equation}\label{varphihteputo0inHN}
\lim_{\epsilon \to 0^+} \big\| \varphi   \mathrm{h}_{t,\epsilon}^u \big\|_{H^N(\mathbb{R}; \mathbb{C})} = 0, \qquad \forall \varphi \in H^N(\mathbb{R}; \mathbb{C}).
\end{equation}
\end{cor}
\begin{proof}
Fix $u \in \mathcal{Z}^N$ for some $N \in \mathbb{N}_+$. Using \eqref{fghtepuDefFor} and the first two formulas in  \eqref{gteputo0andBd}, we obtain that
\begin{equation}\label{dxhutepto0inHN-1}
\big\| \partial_x \mathrm{h}_{t, \epsilon}^u \big\|_{H^{N-1} (\mathbb{R}; \mathbb{C})}  = \big\|\partial_x \mathrm{f}_{ \epsilon}^u  + \partial_x \mathrm{g}_{t, \epsilon}^u \big\|_{H^{N-1} (\mathbb{R}; \mathbb{C})  } \leq \big\|\partial_x \mathrm{f}_{ \epsilon}^u    \big\|_{H^{N-1} (\mathbb{R}; \mathbb{C}) }  + \big\|   \mathrm{g}_{t, \epsilon}^u \big\|_{H^N   (\mathbb{R}; \mathbb{C})} \to 0,
\end{equation}as $\epsilon \to 0^+$,   $\forall t \in \mathbb{R}$. Fix $\varphi \in H^N   (\mathbb{R}; \mathbb{C})$. For any $n \in \{1,2, \cdots, N\}$, we have $\partial_x^{N-n} \varphi \in H^n \hookrightarrow L^{\infty}(\mathbb{R}; \mathbb{C})$ by \eqref{InffH1overfLinfty}.
 The second formula in \eqref{htepuphito0inL2} implies that $\lim_{\epsilon \to 0^+} \big\|  \mathrm{h}_{t, \epsilon}^u   \partial_x^N   \varphi \big\|_{L^2 (\mathbb{R} ; \mathbb{C})}  =0$.
 So we have
\begin{equation}\label{dxNphihutepto0L2}
\big\| \partial_x^N \bigl( \varphi \mathrm{h}_{t, \epsilon}^u \bigr) \big\|_{L^2  (\mathbb{R}; \mathbb{C})}  \leq \big\|  \mathrm{h}_{t, \epsilon}^u   \partial_x^N   \varphi \big\|_{L^2 }  + \sum_{n =1 }^N \tfrac{N !}{n! (N - n)!} \big\| \partial_x^n \mathrm{h}_{t, \epsilon}^u    \big\|_{L^2 } \big\|    \partial_x^{N - n}   \varphi \big\|_{L^{\infty}} \to 0,
\end{equation}as $\epsilon \to 0^+$, by \eqref{dxhutepto0inHN-1} and   \eqref{IndAssumpn=NZNHN}.
Then    \eqref{varphihteputo0inHN} follows directly from 
\eqref{ContrL2dNnorHNnorm},  \eqref{dxhutepto0inHN-1}
and
 \eqref{dxNphihutepto0L2}.
\end{proof}

\noindent Fix $z = a + bi \in \mathbb{C}_+$ with $b = \operatorname{Im} z > 0$. Let $\chi_{b^{-1}} \in L^2_+ \cap H^{\infty}(\mathbb{R}; \mathbb{C})$ be given by \eqref{chiepFor}. We set 
\begin{equation}\label{bdvarphizHinf-}
\boldsymbol{\varphi}_z(x) := \overline{(bi)^{-1} \tau_a \chi_{b^{-1}}(x)} = (x - z)^{-1} \in \mathbb{C}, \qquad \forall x \in \mathbb{R}.
\end{equation}
Then $\boldsymbol{\varphi}_z \in L^2_- \cap H^{\infty}(\mathbb{R}; \mathbb{C}) \subset \bigcap_{q \in [2, +\infty]} L^q_-$. 
Fix $p \in [2, +\infty]$. By inequality \eqref{PoiIntBddLinFunc}, the Poisson integral operator $\mathbf{P}_z : f \mapsto \mathbf{P}_z(f) = \underline{f}(z)$ given by \eqref{PoiIntDefSect2} is a bounded linear functional from $L^p(\mathbb{R}; \mathbb{C})$ to $\mathbb{C}$. For any $f \in L^p(\mathbb{R}; \mathbb{C})$, its \textit{difference quotient} function $\mathscr{Q}_z(f)$ is defined by
\begin{equation}\label{diffQuoDefForQzS2}
\mathscr{Q}_z(f) := f \boldsymbol{\varphi}_z - \mathbf{P}_z(f) \boldsymbol{\varphi}_z \in L^2(\mathbb{R}; \mathbb{C}).
\end{equation}

\begin{lem}
Fix $p \in [2, +\infty]$ and $z \in \mathbb{C}_+$. The difference quotient operator $\mathscr{Q}_z : f \mapsto \mathscr{Q}_z(f)$ given by \eqref{diffQuoDefForQzS2} is a bounded linear operator from $L^p(\mathbb{R}; \mathbb{C})$ to $L^2(\mathbb{R}; \mathbb{C})$. Precisely, we have
\begin{equation}\label{QzBdLpL2}
\big\| \mathscr{Q}_z \big\|_{\mathcal{B}_{\mathbb{C}} (L^p(\mathbb{R}; \mathbb{C}); L^2(\mathbb{R}; \mathbb{C}))} \leq (\operatorname{Im} z)^{-\frac{1}{2} - \frac{1}{p}} \Bigl( \bigl( \mathtt{J}_{\frac{p}{p-2}} \bigr)^{\frac{1}{2} - \frac{1}{p}} + \tfrac{1}{\sqrt{\pi}} \bigl( \mathtt{J}_{\frac{p}{p-1}} \bigr)^{\frac{p-1}{p}} \Bigr),
\end{equation}
where $\mathtt{J}_{\mathfrak{k}} := \int_{\mathbb{R}} \frac{1}{(1+x^2)^{\mathfrak{k}}} \, \mathrm{d}x > 0$ for any $\mathfrak{k} \ge 1$. In particular, the following estimates hold:
\begin{equation}\label{BdQzLinftyL2}
\big\| \mathscr{Q}_z \big\|_{\mathcal{B}_{\mathbb{C}} (L^{\infty}(\mathbb{R}; \mathbb{C}); L^2(\mathbb{R}; \mathbb{C}))} \leq 2 \sqrt{\tfrac{\pi}{\operatorname{Im} z}}, \qquad
\big\| \mathscr{Q}_z \big\|_{\mathcal{B}_{\mathbb{C}} (L^2(\mathbb{R}; \mathbb{C}))} \leq \Bigl(1 + \tfrac{1}{\sqrt{2}}\Bigr) (\operatorname{Im} z)^{-1}.
\end{equation}
\end{lem}
 
\begin{proof}
Fix $p \in [2, +\infty]$ and $z   \in \mathbb{C}_+$. Set $q : =\frac{ 2 p}{p -2 } \in [2 , +\infty]$. By formulas \eqref{PoiKerS2} and \eqref{bdvarphizHinf-}, we have
\begin{equation}\label{bdvarphnormLq}
\|\boldsymbol{\varphi}_z\|_{L^2_-}^2 =  \tfrac{\pi}{\operatorname{Im} z}, \qquad 
\|\boldsymbol{\varphi}_z\|_{L^q (\mathbb{R}; \mathbb{C})}  = (\operatorname{Im} z)^{-\frac{1}{2} - \frac{1}{p}}   \bigl( \mathtt{J}_{\frac{p}{p-2}} \bigr)^{\frac{1}{2} - \frac{1}{p}} , \qquad 
\|\mathfrak{P}_1\|_{L^{\frac{p}{p-1}} (\mathbb{R}; \mathbb{R})} = 
 \tfrac{1}{ \pi}  \bigl( \mathtt{J}_{\frac{p}{p-1}} \bigr)^{\frac{p-1}{p}}   .
\end{equation}Using H\"older's inequality, we obtain formula \eqref{QzBdLpL2}   by
\eqref{diffQuoDefForQzS2}, 
\eqref{PoiIntBddLinFunc}
and \eqref{bdvarphnormLq}. Since $\mathtt{J}_1 =2 \mathtt{J}_2 = \pi$, formula \eqref{BdQzLinftyL2} is a direct consequence of \eqref{QzBdLpL2}.
\end{proof}
   
\begin{lem}\label{QzuchieptoQzuL2lem}
For any $z \in \mathbb{C}_+$ and $u \in L^{\infty}(\mathbb{R} ; \mathbb{C})$, we have $\lim_{\epsilon \to 0^+} \big\| \mathscr{Q}_z \bigl( u \chi_{\epsilon} - u \bigr) \big\|_{L^2(\mathbb{R}; \mathbb{C})} =0$.
\end{lem}
\begin{proof}
Fix $z \in \mathbb{C}_+$ and $u \in L^{\infty}(\mathbb{R} ; \mathbb{C})$. Since $\boldsymbol{\varphi}_z \in L^2_-$ by formulas \eqref{bdvarphizHinf-} and  \eqref{QzuchieptoQzuL2lem}, 
 we obtain that
\begin{equation*}
\big\| \mathscr{Q}_z \bigl( u \chi_{\epsilon} - u \bigr) \big\|_{L^2(\mathbb{R} ; \mathbb{C})} \leq \| u\|_{L^{\infty}(\mathbb{R}; \mathbb{C})} \|\boldsymbol{\varphi}_z (\chi_{\epsilon} -1 )\|_{L^2 (\mathbb{R}; \mathbb{C}) } + \big\| \boldsymbol{\varphi}_z \big\|_{L^2_-} \big| \mathbf{P}_z ( u \chi_{\epsilon} ) -  \mathbf{P}_z ( u  )\big| \to 0,
\end{equation*}as $\epsilon \to 0^+$, by \eqref{diffQuoDefForQzS2}, 
  \eqref{chiephCVtohinHpFor}
and 
\eqref{CvPoiIntuchieptou}.
\end{proof}
 
\begin{lem}\label{LemKerscrQz}
For any $z \in \mathbb{C}_+$ and $p \in [2, +\infty)$, the restriction $\mathscr{Q}_z|_{L^p (\mathbb{R} ; \mathbb{C})} : L^p (\mathbb{R} ; \mathbb{C}) \to L^2 (\mathbb{R} ; \mathbb{C})$ is injective. The kernel of the restriction $\mathscr{Q}_z|_{L^{\infty} (\mathbb{R} ; \mathbb{C})} : L^{\infty} (\mathbb{R} ; \mathbb{C}) \to L^2 (\mathbb{R} ; \mathbb{C})$ is given by
\begin{equation}\label{KerscrQzinfty}
\mathrm{Ker} \bigl( \mathscr{Q}_z|_{L^{\infty} (\mathbb{R} ; \mathbb{C})} \bigr) = \mathbb{C} \mathbf{1} .
\end{equation} 
\end{lem}

\begin{proof}
Fix $z \in \mathbb{C}_+$ and $p \in [2 , +\infty]$. If $f \in \mathrm{Ker} \bigl( \mathscr{Q}_z|_{L^p (\mathbb{R} ; \mathbb{C})} \bigr) $,   formulas \eqref{diffQuoDefForQzS2}   and \eqref{bdvarphizHinf-} imply that
\begin{equation*}
 f(x) (x - z)^{-1} = f(x) \boldsymbol{\varphi}_z(x)  = \mathbf{P}_z(f) \boldsymbol{\varphi}_z(x) = \underline{f}(z) (x - z)^{-1}, \qquad  \forall x \in \mathbb{R}.
\end{equation*}So $f =  \underline{f}(z)
 = \underline{f}(z) \mathbf{1}$ is a constant function. Conversely, we have $\mathscr{Q}_z (\mathbf{1}) = 0$ by \eqref{PzL2BdPz1is1} and \eqref{diffQuoDefForQzS2}.
\end{proof}

\noindent The following lemma shows that $\mathscr{Q}_z$ remains bounded when acting on the   Zhidkov spaces $\mathcal{Z}^N$.

\begin{lem}\label{LemHNCVQzuchiep}
For any $N \in \mathbb{N}_+$ and $z \in \mathbb{C}_+$, the difference quotient operator $\mathscr{Q}_z$ given by \eqref{diffQuoDefForQzS2} is a bounded linear operator from $\mathcal{Z}^N$ to $H^N(\mathbb{R}; \mathbb{C})$. Precisely, we have 
\begin{equation}\label{NormbdQzZNHN}
\big\| \mathscr{Q}_z \big\|_{\mathcal{B}_{\mathbb{C}} (\mathcal{Z}^N; H^N(\mathbb{R}; \mathbb{C}))} \leq (1 + C_N) \big\| \boldsymbol{\varphi}_z \big\|_{H^N_-} \quad \text{with} \quad  \boldsymbol{\varphi}_z (x) = (x - z)^{-1},
\end{equation}
where $C_N > 0$ denotes the constant from Lemma $\ref{ZnHninHnLem}$. Moreover, if $u \in \mathcal{Z}^N$, then
\begin{equation}\label{HNCVQzuchiep-u}
\lim_{\epsilon \to 0^+} \big\| \mathscr{Q}_z \bigl( u \chi_{\epsilon} - u \bigr) \big\|_{H^N(\mathbb{R}; \mathbb{C})} = 0.
\end{equation}
\end{lem}
\begin{proof}
Fix $z \in \mathbb{C}_+$ and $u \in \mathcal{Z}^N$ for some $N \in \mathbb{N}_+$. Using  \eqref{diffQuoDefForQzS2}, \eqref{ufHn<uZnfHn} and
\eqref{PzL2BdPz1is1},  we obtain that
\begin{equation}
\big\| \mathscr{Q}_z (u) \big\|_{H^N (\mathbb{R} ; \mathbb{C})} 
  \leq C_N \big\| u \big\|_{\mathcal{Z}^N} \big\| \boldsymbol{\varphi}_z \big\|_{H^N_-} +  \big\| u \big\|_{L^{\infty} (\mathbb{R} ; \mathbb{C}) } \big\| \boldsymbol{\varphi}_z \big\|_{H^N_-} 
\leq (1 + C_N) \big\| \boldsymbol{\varphi}_z \big\|_{H^N_-} \big\|  u \big\|_{\mathcal{Z}^N}.
\end{equation}Formula \eqref{HNCVQzuchiep-u} follows directly from 
\eqref{diffQuoDefForQzS2}, 
 \eqref{ufHn<uZnfHn}, \eqref{CvPoiIntuchieptou} and the third formula of \eqref{chiephCVtohinHpFor}. 
\end{proof}
 
\begin{lem}
For any   $z \in \mathbb{C}_+$ and $N \in \mathbb{N}_+$, we have 
\begin{equation}\label{QzeitLapuchiepLimHN}
 \lim_{\epsilon \to 0^+} \big\|\mathscr{Q}_z e^{ i t \partial_x^2} \bigl( u \chi_{\epsilon} \bigr) - \mathscr{Q}_z e^{ i t \partial_x^2}  ( u    )   \big\|_{H^N( \mathbb{R} ; \mathbb{C})} = 0, \quad \forall u \in \mathcal{Z}^N, \quad \forall t \in \mathbb{R} .
\end{equation}
\end{lem}

\begin{proof}
Fix $t \in \mathbb{R}$ and $u \in \mathcal{Z}^N$ for some $N \in \mathbb{N}_+$. For any $\epsilon >0$, let $\mathrm{h}_{t,  \epsilon}^u \in \mathcal{Z}^1$ be defined in
 \eqref{fghtepuDefFor}. Then 
\begin{equation}\label{heptuexpitlauchi-u}
\mathrm{h}_{t,  \epsilon}^u = e^{i t \partial_x^2} \bigl( u \chi_{\epsilon} - u \bigr)  \in \mathcal{Z}^N, 
\end{equation}by Corollary \ref{htepuinZNCor}. Fix $z \in \mathbb{C}_+$. 
 Using   formulas  \eqref{diffQuoDefForQzS2},
 \eqref{htepuphito0inL2} and   \eqref{varphihteputo0inHN}, we obtain that  
 \begin{equation}\label{Qzheptuto0inHN}
 \big\|\mathscr{Q}_z \bigl( \mathrm{h}_{t,  \epsilon}^u \bigr) \big\|_{H^N(\mathbb{R} ; \mathbb{C})} 
 = \big\|  \mathrm{h}_{t,  \epsilon}^u  \boldsymbol{\varphi}_z -  \mathbf{P}_z \bigl( \mathrm{h}_{t,  \epsilon}^u  \bigr) \boldsymbol{\varphi}_z\big\|_{H^N }
 \leq  \big\|  \mathrm{h}_{t,  \epsilon}^u  \boldsymbol{\varphi}_z \big\|_{H^N(\mathbb{R} ; \mathbb{C})} + \big| \mathbf{P}_z \bigl( \mathrm{h}_{t,  \epsilon}^u  \bigr) \big| \big\|   \boldsymbol{\varphi}_z \big\|_{H^N_-}
\to 0,
 \end{equation}as $\epsilon \to 0^+$.  
 Thus formula \eqref{QzeitLapuchiepLimHN} follows directly from \eqref{heptuexpitlauchi-u} and \eqref{Qzheptuto0inHN}.
\end{proof}

\noindent We close this subsection by introducing  a 
nonlinear  operator on $\mathcal{Z}^1$. 
Fix $\epsilon > 0$,   then $\chi_{\epsilon} \in L^2_+ \cap H^{\infty}(\mathbb{R} ; \mathbb{C})$ by \eqref{chiepFor}. For any $u \in \mathcal{Z}^1$, we have $|u|^2 \in \mathcal{Z}^1$ by inequality  \eqref{uvZn<uZnvZn}. Then we define 
\begin{equation}\label{DefbfrepuFor}
\mathbf{r}_{\epsilon} (u) : = - 2 i \partial_x u \partial_x \chi_{\epsilon} - i u \partial_x^2 \chi_{\epsilon} - 2 u \Pi_{\geq 0} \bigl( |u|^2  \partial_x \chi_{\epsilon} \bigr)  -2u \Pi_{\geq 0} \bigl( \chi_{\epsilon} \Pi_{ <0} \partial_x \bigl( |u|^2 \bigr)  \bigr) \in  L^2(\mathbb{R}; \mathbb{C}) .
\end{equation}

\begin{lem}\label{repuBdandto0LemL2}
For any $u \in \mathcal{Z}^1$, we have
\begin{equation}\label{bfrepyto0inL2}
\lim_{\epsilon \to 0^+} \|\mathbf{r}_{\epsilon} (u)\|_{L^2( \mathbb{R} ; \mathbb{C})} =0 \quad \text{and} \quad \sup_{0 < \epsilon \leq 1} \|\mathbf{r}_{\epsilon} (u)\|_{L^2( \mathbb{R} ; \mathbb{C})} \leq 4 \| u \|_{\mathcal{Z}^1} \bigl( 1 + \| u \|_{L^{\infty}(\mathbb{R} ; \mathbb{C})}^2 \bigr).
\end{equation}
\end{lem}
\begin{proof}
Fix $u \in \mathcal{Z}^1$, we set $\mathtt{w}: = \Pi_{<0} \partial_x \bigl( |u|^2 \bigr) \in L^2_-$.  By formulas \eqref{L2+perpL2-Har} and
\eqref{chiephCVtohinHpFor}, we have
\begin{equation}\label{ttwEstiL2}
\big\| \mathtt{w} \big\|_{ L^2_- } \leq 2 \|u \|_{L^{\infty} (\mathbb{R} ; \mathbb{C})} \|\partial_x u\|_{L^2 (\mathbb{R} ; \mathbb{C})}  \quad \mathrm{and} \quad \big\| \Pi_{\geq 0} \bigl( \chi_{\epsilon}  \mathtt{w} \bigr)\big\|_{L^2_+} \leq \big\|   \chi_{\epsilon}  \mathtt{w}  - \mathtt{w}  \big\|_{L^2 (\mathbb{R} ; \mathbb{C})} \to 0,
\end{equation}as $\epsilon \to 0^+$.  For any $\epsilon >0$, formula  \eqref{dxNchiepL2norm} yields that $\| \partial_x \chi_{\epsilon}  \|_{L^{\infty}  } =\epsilon$. By \eqref{DefbfrepuFor}, we have
\begin{equation}\label{estibfrL2}
\|\mathbf{r}_{\epsilon}(u)\|_{L^2} \leq  2 \|\partial_x u\|_{L^2} \| \partial_x \chi_{\epsilon} \|_{L^{\infty}} +  \| u \|_{L^{\infty}}  \|\partial_x^2 \chi_{\epsilon}  \|_{L^2} + 2  \| u \|_{L^{\infty}}^3   \| \partial_x \chi_{\epsilon} \|_{L^2_+} + 2   \| u \|_{L^{\infty}}    \| \Pi_{\geq 0} \bigl(  \chi_{\epsilon}  \mathtt{w} \bigr)\|_{L^2_+} .
\end{equation}Then formula \eqref{bfrepyto0inL2} follows directly by \eqref{estibfrL2},
\eqref{ttwEstiL2},
 \eqref{chiephCVtohinHpFor}  and \eqref{dxNchiepL2norm}.
\end{proof}

\noindent For higher regularity, the arguments used in the proof of Lemma~\ref{repuBdandto0LemL2} carry over directly, now using inequality \eqref{ufHn<uZnfHn} to control the products in higher-order Sobolev norms.
   
\begin{lem}\label{repuHighRegLem}
Let $N \geq 2$ be an integer. Then $\mathbf{r}_{\epsilon}(\mathcal{Z}^N) \subset H^{N-1}(\mathbb{R}; \mathbb{C})$ for any $\epsilon > 0$. 
Moreover, there exists a constant $\mathfrak{A}_N^* > 0$ such that for all $u \in \mathcal{Z}^N$, we have
\begin{equation}\label{bfrepyUniBdinHN-1}
\sup_{0 < \epsilon \le 1} \|\mathbf{r}_{\epsilon}(u)\|_{H^{N-1}(\mathbb{R}; \mathbb{C})} \le \mathfrak{A}_N^* \|u\|_{\mathcal{Z}^N} \bigl( 1 + \|u\|_{\mathcal{Z}^{N-1}}^2 \bigr) \quad \mathrm{and} \quad
\lim_{\epsilon \to 0^+} \|\mathbf{r}_{\epsilon}(u)\|_{H^{N-1}(\mathbb{R}; \mathbb{C})} = 0.
\end{equation}
\end{lem}

\noindent Fix $\epsilon > 0$ and $N \in \mathbb{N}_+$. According to formulas  \eqref{ufHn<uZnfHn}, \eqref{uvZn<uZnvZn} and \eqref{HNAlgebra}, 
the operator $\mathbf{r}_{\epsilon}$ given by \eqref{DefbfrepuFor} is locally Lipschitz continuous  
from $\mathcal{Z}^N$ to $H^{N-1}(\mathbb{R}; \mathbb{C})$. Therefore, we have  
\begin{equation}\label{bfrepContonu}
\mathbf{r}_{\epsilon} \in \mathcal{C} \bigl( \mathcal{Z}^N ; H^{N-1} (\mathbb{R} ; \mathbb{C}) \bigr).
\end{equation}

\subsection{Hardy Spaces}\label{SubsecHarSp2.3}

\noindent 
According to the Paley--Wiener theorem, every $L^2$-function whose Fourier transform is supported on $\mathbb{R}_+ = [0, +\infty)$   admits a holomorphic extension  to the upper half-plane $\mathbb{C}_+$ via its Poisson integral  in \eqref{PoiIntDefSect2}.
 Using the Plancherel identity \eqref{inProdL2DefS2}, we can express the Poisson integral of any function $f \in L^2_+$ by the following holomorphic Fourier transform:
\begin{equation}\label{PoiIntHolFTRL2+}
\underline{f}(z) = \mathbf{P}_z (f)  
= \frac{1}{2 \pi i} \int_{\mathbb{R}} \frac{f(\mathfrak{t})}{\mathfrak{t} - z } \mathrm{d}\mathfrak{t}
 = \frac{1}{2 \pi} \int_0^{+ \infty} e^{iz \xi} \widehat{f} (\xi) \mathrm{d}\xi, \qquad \forall z \in \mathbb{C}_+.
\end{equation}Then the mapping $f \in L^2_+ \mapsto \underline{f} \in \operatorname{Hol} (\mathbb{C}_+)$ is a Hilbert isomorphism between $L^2_+$ and the Hardy space
\begin{equation}\label{DefHardySpFor}
\mathfrak{H}^2 \bigl( \mathbb{C}_+ \bigr):=\{ \mathbf{h} \in \operatorname{Hol} (\mathbb{C}_+) \; : \; \sup_{y >0} \int_{\mathbb{R}} |\mathbf{h}(x+ yi )|^2 \mathrm{d}x <+\infty\}.
\end{equation}So $L^2_+ = \{f \in L^2(\mathbb{R}; \mathbb{C}) \; : \; \widehat{f}(\xi)=0, \; \forall \xi <0  \} $ is also called the Hardy space. For any $f \in L^2_+$, we have
\begin{equation}\label{fundertofinL2}
\sup_{y >0} \big\|\underline{f}( \cdot + y i) \big\|_{L^2(\mathbb{R}; \mathbb{C})} \leq \|f\|_{L^2_+} \quad \text{and} \quad \lim_{y \to 0^+} \big\|\underline{f}( \cdot + y i) - f\big\|_{L^2(\mathbb{R}; \mathbb{C})} =0.
\end{equation}
  Now, we introduce three important operators on the Hardy space $L^2_+$. First, we define
\begin{equation}\label{DomI+}
\small
\mathrm{Dom}(\mathbf{I}_+):= \{\varphi \in L^2_+: \exists a \in \mathbb{C} \quad \mathrm{such} \quad \mathrm{that} \quad \forall \epsilon>0, \quad \exists \delta>0  \quad \mathrm{such} \quad \mathrm{that} \quad   \sup_{0 <t < \delta }|\widehat{\varphi}(t) - a| < \epsilon\}.
\end{equation}For any $\varphi \in \mathrm{Dom}(\mathbf{I}_+)+L^1(\mathbb{R}; \mathbb{C})$, let $\mathbf{I}_+ (\varphi)$ denote the right-hand limit of $\widehat{\varphi}$ at $0$, i.e.
\begin{equation}\label{DefI+}
\mathbf{I}_+ (\varphi):= \widehat{\varphi}(0^+) = \lim_{\xi \to 0^+} \widehat{\varphi}(\xi).
\end{equation}
\begin{lem}\label{LemL2+capL10mode}
For any $f \in L^2_+ \cap L^1(\mathbb{R}; \mathbb{C})$, we have $   f \in \mathrm{Dom}(\mathbf{I}_+)$ and $\mathbf{I}_+(f) = \int
_{\mathbb{R}}f = 0$.
\end{lem}
\begin{proof}
If $f \in L^2_+ \cap  L^1(\mathbb{R}; \mathbb{C})$, its Fourier transform $\widehat{f}$ is continuous and $\operatorname{supp} (\widehat{f}) \subset \mathbb{R}_+ = [0, +\infty)$. Then $ \mathbf{I}_+(f) = \widehat{f}(0^+) = \int_{\mathbb{R}}f = \widehat{f}(0^-) = 0$.
\end{proof}

\begin{defi}\label{DombfGdefSect2}
Set $\mathrm{Dom}(\mathbf{G}):= \{\varphi \in L^2_+: \widehat{\varphi}|_{\mathbb{R}_+^*} \in H^1(\mathbb{R}_+^*)\}$. For any  $\varphi \in \mathrm{Dom}(\mathbf{G})$,  we define
\begin{equation}\label{DefGphiSect2}
\widehat{\mathbf{G} \varphi} (\xi) : = i  \partial_{\xi} \widehat{\varphi}(\xi), \quad \forall \xi >0 .
\end{equation}
\end{defi}
\noindent It is important to distinguish the operator $\mathbf{G}: \mathrm{Dom}(\mathbf{G}) \to L^2_+$ in \eqref{DefGphiSect2} from the operator of multiplying the monomial $x$, which we denote by $\mathbf{X}: f\in L^2_+ \mapsto \mathbf{X}  f \in \mathscr{S}'$, where  $\mathbf{X}f$ is given  by
\begin{equation}\label{defXf-}
\mathbf{X}f(x):=xf(x) \in \mathbb{C}, \quad \forall x\in \mathbb{R}, \quad \forall f \in L^2_+.
\end{equation}To ensure that $\mathbf{X}f$ remains in $L^2_+$, the domain of definition of $\mathbf{X}$ must be  chosen as follows:
\begin{equation}\label{DomX}
\mathrm{Dom}(\mathbf{X}):=\{\varphi \in L^2_+: \widehat{\varphi}  \in H^1(\mathbb{R} ; \mathbb{C})\} = \{\varphi \in \mathrm{Dom}(\mathbf{G}): \widehat{\varphi}(0^+)=0\}.
\end{equation}The Sobolev embedding and extension theorems show that $H^1(\mathbb{R}^*_+; \mathbb{C}) \hookrightarrow H^1(\mathbb{R}; \mathbb{C}) \hookrightarrow 
\mathcal{C}_0(\mathbb{R}; \mathbb{C})$ and 
\begin{equation}\label{DXsubsetDG}
L^2_+ \cap  L^2_1(\mathbb{R}) =\mathrm{Dom}(\mathbf{X}) \subset \mathrm{Dom}(\mathbf{G}) \subset \mathrm{Dom}(\mathbf{I}_+) \subset L^2_+ .
\end{equation}The next lemma shows the difference between $\mathbf{X}:\mathrm{Dom}(\mathbf{X}) \to L^2_+$ and $\mathbf{G}:\mathrm{Dom}(\mathbf{G}) \to L^2_+$.

\begin{lem}\label{LemGphi=Xphi-const}
For any $\varphi \in \mathrm{Dom}(\mathbf{G})$, we have $\mathbf{G} \varphi  = \mathbf{X} \varphi  - \frac{i}{2 \pi}  \widehat{\varphi}(0^+) \in L^2_+$.
\end{lem}
\begin{proof}
If $h \in \mathscr{S}$ and  $\varphi \in \mathrm{Dom}(\mathbf{G})$, then $\partial_{\xi} \widehat{\varphi} \in H^1(\mathbb{R}_+^*)$ and $\lim_{\xi \to + \infty} (|\widehat{\varphi}(\xi)| + |h(\xi)|)=0$. We have $\langle \widehat{\mathbf{X}\varphi}, h\rangle_{\mathscr{S}', \mathscr{S}} =
 - i\int_{\mathbb{R}} \widehat{\varphi} \partial_{\xi} h   = i\widehat{\varphi}(0^+) h(0) + i\int_{0}^{+\infty}h \partial_{\xi} \widehat{\varphi}$ by   Plancherel's identity and integration by parts. Thus $\widehat{\mathbf{X}\varphi} = i \widehat{\varphi}(0^+) \delta_0 + i \mathbf{1}_{\mathbb{R}_+^*} \partial_{\xi} \widehat{\varphi} \in \mathscr{S}'$. We conclude by the invertibility of the Fourier transform.
\end{proof}

\noindent Let $\sigma_{\mathrm{p}}(\mathbf{G})$ denote the set of all eigenvalues of the operator $\mathbf{G} : \mathrm{Dom}(\mathbf{G}) \to L^2_+$.

\begin{lem}\label{LemSpGC_-}
We have $\sigma_{\mathrm{p}}(\mathbf{G}) = \mathbb{C}_-= \{a+bi \in \mathbb{C} : b<0 \} $. Moreover, for any $\zeta \in \mathbb{C}_-$, we have
\begin{equation}\label{KerG-zetadim1}
\mathrm{Ker}\left( \mathbf{G} - \zeta \right) = \mathbb{C}\boldsymbol{\varphi}_{\zeta}, \quad \mathrm{where} \quad \boldsymbol{\varphi}_{\zeta}(x):= (x - \zeta)^{-1}, \quad \forall x \in \mathbb{R}.
\end{equation}
\end{lem}

\begin{proof}
Assume that $\mathbf{G}f = z f$ for some $f \in \mathrm{Dom}(\mathbf{G})\setminus \{0\}$ and $z \in \mathbb{C}$, set $\mathbf{e}_z(\xi):= \exp(iz \xi)$, $\forall \xi \in \mathbb{R}$, then we have $\partial_{\xi} \left(\mathbf{e}_z \widehat{f} \right) = 0$ on $\mathbb{R}_+^*=(0, +\infty)$. So there exists $c \in  \mathbb{C}\setminus \{0\}$ such that $\widehat{f}(\eta)=c\exp(-iz \eta)$, $\forall \eta >0$. Since $\widehat{f}|_{\mathbb{R}_+^*} \in H^1(\mathbb{R}_+^*)$, we have $\operatorname{Im} z <0$ and $f = \frac{ic}{2\pi}\boldsymbol{\varphi}_z \in \mathrm{Dom}(\mathbf{G})$. Conversely, if $\zeta \in \mathbb{C}_-$, we have $\widehat{\boldsymbol{\varphi}}_{\zeta} (\xi)=-2\pi i \exp(-i \zeta \xi)\mathbf{1}_{\mathbb{R}_+^*}(\xi)$. So \eqref{DefGphiSect2} yields that $\widehat{\mathbf{G}\boldsymbol{\varphi}_{\zeta}} = \zeta \widehat{\boldsymbol{\varphi}}_{\zeta} \in H^1(\mathbb{R}_+^*)$.
\end{proof}

\begin{lem}\label{IGXUnboundLem}
The linear mappings $ \mathbf{I}_+: \mathrm{Dom}(\mathbf{I}_+) \to \mathbb{C}$, $\mathbf{G} : \mathrm{Dom}(\mathbf{G}) \to L^2_+$ and $\mathbf{X} : \mathrm{Dom}(\mathbf{X}) \to L^2_+$ are all discontinuous and densely defined on the Hardy space $L^2_+$.
\end{lem}
\begin{proof}
Set $f_n(x):=(n-ix)^{-1}$ and $g_n(x):= \sqrt{n}(n-ix)^{-2}$, $\forall x\in \mathbb{R}$, $\forall n\in \mathbb{N}_+$.  Lemma \ref{LemSpGC_-} yields that $f_n \in \mathrm{Ker}(\mathbf{G}+in) \subset \mathrm{Dom}(\mathbf{G})\subset \mathrm{Dom}(\mathbf{I}_+)$, $g_n \in \mathrm{Dom}(\mathbf{X})$ and $\widehat{f}_n(\xi) = 2 \pi \exp(-n \xi) \mathbf{1}_{\mathbb{R}_+^*}(\xi)$, $\forall \xi \in \mathbb{R}$. So
\begin{equation}\label{normGfntoinfty}
\|\mathbf{X}g_n\|_{L^2_+} = n \|g_n\|_{L^2_+} = \sqrt{\tfrac{\pi}{2}},  \quad \mathbf{I}_+(f_n)=\widehat{f}_n(0^+) = 2\pi, \quad \|\mathbf{G}f_n\|_{L^2_+} = n \|f_n\|_{L^2_+} = \sqrt{n \pi}\to +\infty,
\end{equation}when $n \to +\infty$. So $\mathbf{I}_+$, $\mathbf{G}$ and $\mathbf{X}$ are unbounded operators on $L^2_+$. Since $\mathcal{C}_{c}^{\infty}(\mathbb{R}_+^*; \mathbb{C})$ is dense in $L^2(\mathbb{R}_+^*; \mathbb{C})$, $\mathrm{Dom}(\mathbf{X})$ is dense in the Hardy space $L^2_+$.
\end{proof}
 
\noindent Another connection between $\mathbf{G} : \mathrm{Dom}(\mathbf{G}) \to L^2_+$ and $\mathbf{X} : \mathrm{Dom}(\mathbf{X}) \to L^2_+$ is given by the following lemma.
\begin{lem}\label{XadjisGLem}
We have $\mathbf{G}^* = \mathbf{X}$ and $\mathbf{X}^* = \mathbf{G}$. The operator $(\mathbf{X}, \mathrm{Dom}(\mathbf{X}))$ is   symmetric on $L^2_+$.
\end{lem}
\begin{proof}
The proof is divided into two steps. \textbf{Step \romannumeral 1}. We show that $\mathbf{G}^* = \mathbf{X}$. \\

\noindent Fix $\psi \in \mathrm{Dom}(\mathbf{G}^*)$, there exists a unique function $g \in L^2_+$ such that
\begin{equation}\label{intpsiphiPla}
 \frac{i}{2\pi}\int_0^{+\infty}   \overline{\widehat{\psi}} \partial_{\xi}\widehat{\varphi}=\langle \mathbf{G} \varphi, \psi\rangle_{L^2_+} = \langle  \varphi , g\rangle_{L^2_+}= \frac{1}{2 \pi} \int_0^{+\infty} \widehat{\varphi} \overline{\widehat{g}}, \quad \forall \varphi \in \mathrm{Dom}(\mathbf{G}).
\end{equation}According to  the definition of  $H^1(\mathbb{R}_+^*)$ in \eqref{WeightedSobolevSp}, we have $\widehat{\psi}|_{\mathbb{R}_+^*} \in H^1(\mathbb{R}_+^*)$. So $\psi \in \mathrm{Dom}(\mathbf{G})$ and $\lim_{\xi \to +\infty}\widehat{\psi} (\xi) =0$. Then Corollary 8.10 of \cite{BrezisFA} and \eqref{intpsiphiPla} imply that
\begin{equation}\label{IpsiintbParts}
2 \pi \langle \mathbf{G} \varphi, \psi\rangle_{L^2_+} = -i\overline{\widehat{\psi}(0^+)}\widehat{\varphi}(0^+) -i\int_0^{+\infty}   \widehat{\varphi}\overline{\partial_{\xi}\widehat{\psi}}  = 2 \pi \langle  \varphi, \mathbf{G}\psi\rangle_{L^2_+} -i \overline{\mathbf{I}_+(\psi)}\mathbf{I}_+(\varphi), 
\end{equation}for any $\varphi\in \mathrm{Dom}(\mathbf{G})$. The discontinuity of $\mathbf{I}_+ : \varphi \in \mathrm{Dom}(\mathbf{G}) \subset \mathrm{Dom}(\mathbf{I}_+) \mapsto \mathbf{I}_+(\varphi) \in \mathbb{C}$ is given by  Lemma $\ref{IGXUnboundLem}$. Since the function $\varphi \in (\mathrm{Dom}(\mathbf{G}), \|.\|_{L^2_+}) \mapsto \langle \mathbf{G} \varphi, \psi\rangle_{L^2_+}\in \mathbb{C}$ is continuous, we have $\widehat{\psi}(0^+)=\mathbf{I}_+(\psi) = 0$. So $\mathrm{Dom}(\mathbf{G}^*) \subset \mathrm{Dom}(\mathbf{X})$. Using \eqref{DomX} and Lemma $\ref{LemGphi=Xphi-const}$, we obtain that
\begin{equation}\label{SymofXinL2+}
\mathbf{G}^*(\psi) = \mathbf{G} (\psi) = \mathbf{X}\psi \in L^2_+, \quad \forall \psi \in \mathrm{Dom}(\mathbf{G}^*) \subset \mathrm{Dom}(\mathbf{X}) \subset \mathrm{Dom}(\mathbf{G}).
\end{equation}Conversely, if $h \in \mathrm{Dom}(\mathbf{X})$, then $\widehat{h} \in H^1(\mathbb{R}; \mathbb{C})$ and $\widehat{h}(0^+) = 0=\lim_{\xi \to +\infty}\widehat{h}(\xi)$. So we have
\begin{equation}\label{XhGphi-}
2 \pi \langle \mathbf{G} \varphi, h\rangle_{L^2_+} =i\int_0^{+\infty}   \overline{\widehat{h}} \partial_{\xi}\widehat{\varphi} =
 -i\int_0^{+\infty}   \widehat{\varphi}\overline{\partial_{\xi}\widehat{h}}  = 2 \pi \langle  \varphi, \mathbf{X}h\rangle_{L^2_+} \in \mathbb{C}, \quad \forall\varphi \in \mathrm{Dom}(\mathbf{G}),
\end{equation}by Lemma $\ref{LemGphi=Xphi-const}$. Thus, $h \in \mathrm{Dom}(\mathbf{G}^*)$ and $\mathrm{Dom}(\mathbf{G}^*) = \mathrm{Dom}(\mathbf{X})$. So formula  \eqref{SymofXinL2+} yields that $\mathbf{G}^* = \mathbf{X}$.\\

\noindent \textbf{Step \romannumeral 2}. We claim that $\left(\mathbf{G}, \mathrm{Dom}(\mathbf{G})\right)$ is a closed operator on $L^2_+$.
\begin{proof}[Proof of the closedness of graph of $\mathbf{G}$]Assume that $\left(\varphi_n, \mathbf{G}\varphi_n \right)_{n \geq 0}$ is a convergent sequence in $L^2_+ \times L^2_+$, then $(\widehat{\varphi}_n|_{\mathbb{R}_+^*})_{n \geq 0}$ is a Cauchy sequence in the Sobolev space $H^1(\mathbb{R}_+^*)$. Let $f:= \lim_{n \to +\infty}\widehat{\varphi}_n|_{\mathbb{R}_+^*} $ denote its limit in $H^1(\mathbb{R}_+^*)$. We define $\varphi \in L^2_+$ by using Fourier transform: $\widehat{\varphi}(\xi)=f(\xi)$, $\forall \xi >0$. So  
\begin{equation}\label{phin-phiconv}
\varphi \in \mathrm{Dom}(\mathbf{G})  \quad \mathrm{and} \quad \|\varphi_n - \varphi\|_{L^2_+}^2+\|\mathbf{G}\varphi_n - \mathbf{G}\varphi\|_{L^2_+}^2 = \tfrac{1}{2\pi} \|\widehat{\varphi}_n - f\|_{H^1(\mathbb{R}_+^*)}^2 \to 0,
\end{equation}as $n \to +\infty$. Then $\Gamma (\mathbf{G}):=\{(\phi, \mathbf{G}\phi)\in L^2_+ \times L^2_+ : \forall \phi \in \mathrm{Dom}(\mathbf{G})\}$ is closed in $L^2_+ \times L^2_+$.
\end{proof}
\noindent According to Lemma $\ref{IGXUnboundLem}$, the operator $(\mathbf{G}, \mathrm{Dom}(\mathbf{G}) )$ is both densely defined and closed in $L^2_+$. So we have $\mathbf{X}^* = \mathbf{G}^{**}= \mathbf{G}$ by Theorem 13.12 of  \cite{RudinFA}.
\end{proof}
 
\noindent Another useful characterization of $\mathrm{Dom}(\mathbf{G})$  is provided by the following lemma.

\begin{lem}\label{LemDomGXphi+ainL2}
We have $\mathrm{Dom}(\mathbf{G}) = \{\varphi \in L^2_+ : \exists a \in \mathbb{C} \text{ such that } \mathbf{X}\varphi + a \in L^2(\mathbb{R};\mathbb{C})\}$.
\end{lem}
\begin{proof}
If $\varphi \in \mathrm{Dom}(\mathbf{G})$, Lemma~\ref{LemGphi=Xphi-const} gives
$\mathbf{X}\varphi + (2\pi i)^{-1} \widehat{\varphi}(0^+) = \mathbf{G}\varphi \in L^2_+$.
Thus $\varphi$ belongs to the set on the right with $a = (2\pi i)^{-1} \widehat{\varphi}(0^+) \in \mathbb{C}$.  
   Conversely, assume that $\phi \in L^2_+$ and there exists $\mathsf{a} \in \mathbb{C}$ such that $ h := \mathbf{X} \phi +\mathsf{a}  \in L^2 (\mathbb{R} ; \mathbb{C})$.  
 Fix an arbitrary  test function  $\mathtt{v} \in \mathcal{C}_c^{\infty} (\mathbb{R} ; \mathbb{C})$ such that  $\operatorname{supp} (\mathtt{v}) \subset \mathbb{R}_-^* = ( -\infty, 0)$. Since $\operatorname{supp}(\widehat{\phi}) \subset \mathbb{R}_+ = [0,+\infty)$, using the Fourier transform \eqref{FouTrDefS2}, we have
\begin{equation*}
\int_{\mathbb{R}} \widehat{h} \mathtt{v} 
= \langle \widehat{h} , \mathtt{v} \rangle_{\mathscr{S}' ,  \mathscr{S}} = - i \langle \widehat{\phi} , \partial_{\xi} \mathtt{v} \rangle_{\mathscr{S}' , \mathscr{S}} + 2 \pi \mathsf{a} \mathtt{v}(0) = 0 .
\end{equation*}
  So $\mathbf{X} \phi +\mathsf{a}  = h \in L^2_+$. For any $f \in \mathrm{Dom}( \mathbf{X} )$, formula  \eqref{DXsubsetDG} yields that $ f \in  L^2_+ \cap L^2_1 (\mathbb{R}) \subset L^2_+ \cap L^1 (\mathbb{R};\mathbb{C})$. Then $\int_{\mathbb{R}}f = \mathbf{I}_+ (f ) =0$ by Lemma  $\ref{LemL2+capL10mode}$. By  \eqref{defXf-}, we have 
\begin{equation*}
\langle \mathbf{X}f , \phi \rangle_{L^2_+} = \int_{\mathbb{R}} f \bigl( \overline{\mathbf{X}\phi  + \mathsf{a} - \mathsf{a}} \bigr) = \int_{\mathbb{R}} f \overline{h}  - \overline{\mathsf{a}} \int_{\mathbb{R}} f  = \langle  f , h \rangle_{L^2_+} , \quad \forall f \in \mathrm{Dom}(\mathbf{X}) .
\end{equation*}
Therefore, Lemma $\ref{XadjisGLem}$ implies $\phi \in \mathrm{Dom}(\mathbf{X}^*) = \mathrm{Dom}(\mathbf{G})$.
\end{proof}
   
\noindent The next lemma shows that the family $(\chi_\epsilon)_{\epsilon>0}$ given by \eqref{chiepFor} can   be used to compute the boundary value $\mathbf{I}_+(\varphi) = \widehat{\varphi}(0^+)$ for functions in the domain of $\mathbf{G}$.

\begin{lem}\label{phichiepCVtoI0Lem}
For any $ \varphi \in \mathrm{Dom}(\mathbf{G}) \subset \mathrm{Dom}(\mathbf{I}_+)$,   we have  $\lim_{\epsilon \to 0^+} \langle \varphi, \chi_{\epsilon} \rangle_{L^2_+} = \widehat{\varphi}(0^+) = \mathbf{I}_+ (\varphi)$. 
\end{lem}
\begin{proof}
See formula $(2.3)$ of \cite{GerPush2023}.
\end{proof}

\noindent For any $s \geq 0$, the \textit{Hardy--Sobolev spaces} $H^s_+$ and $H^s_-$ are defined by
\begin{equation}\label{HardySobHspmS2}
H^s_+ := L^2_+ \cap H^s(\mathbb{R}; \mathbb{C}) = \Pi_{\geq 0}\bigl( H^s(\mathbb{R}; \mathbb{C}) \bigr), \qquad
H^s_- := L^2_- \cap H^s(\mathbb{R}; \mathbb{C}) = \Pi_{< 0}\bigl( H^s(\mathbb{R}; \mathbb{C}) \bigr).
\end{equation}Both $H^s_+$ and $H^s_-$ are closed subspaces in the Hilbert space $\bigl( H^s(\mathbb{R}; \mathbb{C}), \| .\|_{H^s }\bigr)$. If $t > s \geq 0$, then  we have $H^t_+ \subset H^s_+$ and $H^t_- \subset H^s_-$. The smooth Hardy--Sobolev spaces are denoted by
\begin{equation}\label{HinftypmDef}
H^{\infty}_+ : = \bigcap_{s \geq 0} H^s_+, \qquad H^{\infty}_- : = \bigcap_{s \geq 0} H^s_-.
\end{equation}
The following criterion allows us to verify that a function belongs to the Hardy space $L^2_+$ by testing against smooth positive-frequency multipliers.

\begin{lem}\label{lemfhinHarSobfL2+}
 Let $f \in L^2(\mathbb{R}; \mathbb{C})$ be a function such that its Fourier transform $\widehat{f}$ is continuous on $\mathbb{R}^*_- = (-\infty, 0)$.   
 If $f h \in L^2_+$ for every $h \in H^\infty_+$, then $f \in L^2_+$.
\end{lem}

\begin{proof}
Fix  $f \in L^2 (\mathbb{R}; \mathbb{C})$ such that $\widehat{f} \in \mathcal{C}( \mathbb{R}_-^* ; \mathbb{C})$. Assume that $\exists \xi_0 <0$ such that $\widehat{f}(\xi_0) \ne 0$. We define
\begin{equation}\label{alpphideffor1}
\alpha= \alpha(f, \xi_0) := \tfrac{\overline{\widehat{f}(\xi_0)}}{|\widehat{f}(\xi_0)|} \in \mathbb{S}^1 \quad \mathrm{and} \quad \varphi (\xi):= \operatorname{Re}(\alpha \widehat{f}(\xi))\in \mathbb{R}, \quad \forall \xi \in \mathbb{R}.
\end{equation}Then $\varphi(\xi_0)= \alpha \widehat{f}(\xi_0) = |\widehat{f}(\xi_0)|>0$  and $\varphi \in L^2(\mathbb{R}; \mathbb{R}) \cap  \mathcal{C}(\mathbb{R}^*_-; \mathbb{R})$. So there exists $\delta \in (0, |\xi_0|)$ such that 
\begin{equation}\label{minphi>0.5fhatxi0}
\min\{\varphi(\zeta):  \forall \zeta \in [\xi_0 - \delta, \xi_0 + \delta]\} \geq \tfrac{1}{2}|\widehat{f}(\xi_0)|.
\end{equation}We define a function $g = g_{\delta} \in H^{\infty}_+$ via its Fourier transform: $\widehat{g}(\xi) = \widehat{g}_{\delta} (\xi) : = \mathbf{1}_{[0, \frac{1}{2}\delta]}(\xi)$, $\forall \xi \in \mathbb{R}$. For any $\xi \in [\xi_0 -\frac{\delta}{2}, \xi_0 + \frac{\delta}{2}]$ and $\eta \in [0, \frac{\delta}{2}]$, we have $|\xi-\eta -\xi_0| \leq \delta$. Then formula  \eqref{minphi>0.5fhatxi0} yields that
\begin{equation}
2 \pi \operatorname{Re}\left(\alpha \widehat{fg}(\xi) \right) = \int_{\mathbb{R}} \varphi(\xi - \eta) \widehat{g}(\eta) \mathrm{d} \eta = \int_0^{\frac{\delta}{2}} \varphi(\xi - \eta) \mathrm{d}\eta \geq \tfrac{\delta}{4}|\widehat{f}(\xi_0)|.
\end{equation}We thus have $\min_{|\xi - \xi_0| \leq \frac{\delta}{2}} \operatorname{Re}(\alpha \widehat{fg}(\xi)) \geq \frac{\delta}{8 \pi}|\widehat{f}(\xi_0)|$ with $\xi_0 <0$ and $\alpha \in \mathbb{S}^1$. Hence, formula  \eqref{LppmDefForS2}  implies  that $fg \notin L^2_+$  with $g = g_{\delta} \in H^{\infty}_+$.  This completes the proof by contradiction.
\end{proof}

\noindent The next lemma provides a quantitative control on the low‑frequency
behaviour of the nonlinear term appearing in the focusing CMdNLS equation \eqref{CMdNLSfocu}.  This
estimate will be essential in the construction of the conservation law
involving the functional $\mathbf{I}_+$ in Subsection~\ref{InvSpCoLsubs4.4}.

\begin{lem}\label{LemNLtermCVEsti}
For any $a \in \mathbb{C}$ and $f \in \mathrm{Dom} (\mathbf{G}) \cap H^1_+$, we define $u:= a + f \in \mathcal{Z}^1 \cap L^{\infty}_+$ and
\begin{equation}\label{varpiudeffor}
\boldsymbol{\varpi}_u : = u \Pi_{\geq 0} \partial_x \bigl( |u|^2 \bigr) \in L^2_+ .
\end{equation}Let $\widetilde{\mathfrak{m}} >0$ denote the Sobolev constant of  $H^1(\mathbb{R}_+^* ; \mathbb{C}) \hookrightarrow L^{\infty}(\mathbb{R}_+^* ; \mathbb{C})$  in Theorem 8.8 of \cite{BrezisFA}, i.e. 
\begin{equation}\label{tilfrakmdefSoCo}
(\widetilde{\mathfrak{m}})^{-1} : = \inf \Big\lbrace \tfrac{\| v \|_{H^1(\mathbb{R}_+^* ; \mathbb{C})}}{\| v \|_{L^{\infty}(\mathbb{R}_+^* ; \mathbb{C}) }} : \forall v \in H^1(\mathbb{R}_+^* ; \mathbb{C}) \Big\rbrace \in \mathbb{R}_+^* = (0 , + \infty).
\end{equation}
 Then we have $\lim_{\xi \to 0^+} \widehat{\boldsymbol{\varpi}_u}  (\xi) = 0$.  Moreover, 
 the following estimate holds: 
\begin{equation}\label{Esthatvrpiuxi}
\sup_{0 < \xi \leq 1} \big| \widehat{\boldsymbol{\varpi}_u} (\xi) \big| \leq \sqrt{2 \pi} \widetilde{\mathfrak{m}} |a|^2 \bigl( \|f\|_{L^2_+} + \| \mathbf{G} f \|_{L^2_+} \bigr) + 2 \bigl( |a|+ \|u\|_{L^{\infty} (\mathbb{R}; \mathbb{C})}\bigr) \|f\|_{L^2_+} \|\partial_x f\|_{L^2_+} .
\end{equation}
\end{lem}

\begin{proof}
Fix $a \in \mathbb{C}$ and $f \in \mathrm{Dom} (\mathbf{G}) \cap H^1_+$, we set $u: = a + f \in a + H^1_+ \subset \mathcal{Z}^1 \cap L^{\infty}_+$. Then we have
\begin{equation}\label{dxmodusquaf}
\partial_x \bigl( |u|^2 \bigr)
= 2 \operatorname{Re} \bigl( \overline{u} \partial_x u  \bigr)
  = a \partial_x \overline{f} + \overline{a} \partial_x f + 2 \operatorname{Re} \bigl( \overline{f} \partial_x f  \bigr) \in L^2 (\mathbb{R} ; \mathbb{R}).
\end{equation}The H\"older's inequality and   \eqref{HNAlgebra} imply that 
$|f|^2 \in H^1 \cap W^{1,1} (\mathbb{R} ; \mathbb{R})$. Moreover, we have
\begin{equation}\label{fmodsqW11norm}
 \big\|  \widehat{|f|^2}   \big\|_{L^{\infty}  }
\leq   \| f \|_{L^2_+}^2 ,
\quad
\big\| \partial_x \bigl( |f|^2 \bigr) \big\|_{L^1}  \leq 2 \| f \|_{L^2_+} \| \partial_x f \|_{L^2_+} , \quad 
\big\| \partial_x \bigl( |u|^2 \bigr) \big\|_{L^2 } \leq 2  \| u \|_{L^{\infty}} \| \partial_x u \|_{L^2_+}.
\end{equation}
  We define $\mathcal{N}_1 : = \partial_x f \in L^2_+$, $\mathcal{N}_2 : = \Pi_{\geq 0} \partial_x \bigl( |f|^2 \bigr) \in L^2_+$ and $\mathcal{N}_3 : = f \Pi_{\geq 0} \partial_x \bigl( | u |^2 \bigr) \in L^2_+ \cap L^1 (\mathbb{R}; \mathbb{C})$, then
\begin{equation}\label{calN123def}
\boldsymbol{\varpi}_u  = u \Pi_{\geq 0} \partial_x \bigl( |u|^2 \bigr)
 = |a|^2 \partial_x f + a \Pi_{\geq 0} \partial_x \bigl( |f|^2 \bigr) + f \Pi_{\geq 0} \partial_x \bigl( | u |^2 \bigr) 
 = |a|^2 \mathcal{N}_1 + a \mathcal{N}_2 + \mathcal{N}_3 \in L^2_+,
\end{equation}by \eqref{dxmodusquaf}. 
 Since $f \in \mathrm{Dom} (\mathbf{G})$, we have $\widehat{f}|_{\mathbb{R}_+^*} \in H^1 (\mathbb{R}_+^*; \mathbb{C})$ by Definition $\ref{DombfGdefSect2}$. By formula  \eqref{DXsubsetDG}, the limit $\widehat{f} (0^+) = \lim_{\xi \to 0^+} \widehat{f} (\xi)$ exists. Using  
 \eqref{tilfrakmdefSoCo},   
 \eqref{H1sWeiSobnormDef} and  \eqref{DefGphiSect2}, we obtain that 
 \begin{equation}\label{calN1limesti}
 \lim_{\xi \to 0^+} \widehat{\mathcal{N}}_1 (\xi) =  \lim_{\xi \to 0^+} \bigl( i \xi \widehat{f}(\xi) \bigr) =0, 
 \quad  
  \sup_{0 < \xi \leq 1} \big| \widehat{\mathcal{N}}_1 (\xi) \big| \leq  \sup_{ \xi > 0} \big| \widehat{f} (\xi) \big|  \leq \sqrt{2 \pi} \widetilde{\mathfrak{m}} \bigl( \|f\|_{L^2_+} + \|\mathbf{G} f\|_{L^2_+}  \bigr) .
 \end{equation}For any $\xi > 0$, we have $\big| \widehat{\mathcal{N}}_2 (\xi) \big| \leq \min \{ \big\|  \widehat{|f|^2}   \big\|_{L^{\infty} (\mathbb{R} ; \mathbb{C})} |\xi| \; , \; \big\| \partial_x \bigl( |f|^2 \bigr) \big\|_{L^1(\mathbb{R} ; \mathbb{R})}\}$. 
 It follows that
 \begin{equation}\label{limcalN2Esti}
 \lim_{\xi \to 0^+} \widehat{\mathcal{N}}_2 (\xi) =  0
  \quad \text{and} \quad \sup_{\xi > 0} \big| \widehat{\mathcal{N}}_2 (\xi) \big| \leq 2 \| f \|_{L^2_+} \| \partial_x f \|_{L^2_+}, 
 \end{equation}by the first two formulas in \eqref{fmodsqW11norm}. Since $\mathcal{N}_3  = f \Pi_{\geq 0} \partial_x \bigl( | u |^2 \bigr) \in L^2_+ \cap L^1 (\mathbb{R}; \mathbb{C})$, we have
\begin{equation}\label{EstiLimcalN3}
 \lim_{\xi \to 0^+} \widehat{\mathcal{N}}_3 (\xi) = \mathbf{I}_+ (\mathcal{N}_3 ) = 0,
  \qquad  \sup_{\xi > 0} \big| \widehat{\mathcal{N}}_3 (\xi) \big| \leq 2 \| f \|_{L^2_+} \| \partial_x u \|_{L^2_+} \| u \|_{L^{\infty} (\mathbb{R} ; \mathbb{C})},
\end{equation}by  Lemma  $\ref{LemL2+capL10mode}$ and the third estimate in \eqref{fmodsqW11norm}.   
  Plugging  
 \eqref{calN1limesti}, \eqref{limcalN2Esti} and
\eqref{EstiLimcalN3} into formula \eqref{calN123def}, we obtain estimate  \eqref{Esthatvrpiuxi}.
  This completes the proof.
\end{proof}

 \begin{lem}\label{DNSelfadjinHsLem}
 Assume that $s \geq 0$ and $N \in \mathbb{N}_+$, then $ D^N : \mathrm{Dom}^s(D^N):=H^{s+N}_+ \to H^s_+$ is an unbounded positive self-adjoint operator on the Hardy--Sobolev space  $\bigl( H^s_+, \langle \cdot , \cdot \rangle_{H^s_+} \bigr)$.
 \end{lem}
 \begin{proof}
Fix $s\geq 0$. For any $f \in H^{s+N}_+$, we have $2 \pi \langle D^N f, f \rangle_{H^s_+} = \int_0^{+\infty} |\widehat{f}(\xi)|^2 \xi^N (1+\xi^2)^s \mathrm{d}\xi \geq 0$. Then $D^N \geq 0$ and it is automatically a symmetric operator, i.e.  $D^N \subset  (D^N)^*$, $\forall N \in \mathbb{N}_+$.  
 Assume that $g, \psi \in H^s_+$ such that for any $\varphi \in \mathrm{Dom}^s(D^N):=H^{s+N}_+$, we have
\begin{equation*}
\int_0^{+\infty} \widehat{\varphi}(\xi) \overline{\widehat{\psi}(\xi)} (1+\xi^2)^s \mathrm{d}\xi = 2 \pi \langle \varphi , \psi \rangle_{H^s_+} =  2 \pi \langle D^N \varphi , g \rangle_{H^s_+} =\int_0^{+\infty} \widehat{\varphi}(\xi) \overline{\widehat{g}(\xi)} \xi^N (1+\xi^2)^s \mathrm{d}\xi.
\end{equation*}Then $\psi = D^N g \in H^s_+$ and $ \|g\|_{H^{s+N}_+}^2 \leq 2^{N-1}  (\|g\|_{H^s_+}^2 + \| \psi \|_{H^s_+}^2 )$. So $D^N  = (D^N)^* \geq 0$. We set 
\begin{equation*}
h_n(x) : = (n x +i )^{-1}  \quad \text{and} \quad
\varphi_n (x) := \bigl( \langle D \rangle^{-s} h_n \bigr) (x) , \quad \forall x \in \mathbb{R}, \quad \forall n \in \mathbb{N}_+ , 
\end{equation*}where $\langle D \rangle^{-s}$ denotes the Fourier multiplier with symbol $\langle \xi \rangle^{-\frac{s}{2}}=(1+\xi^2)^{-\frac{s}{2}}$.
Then $h_n, \varphi_n \in H^{\infty}_+$. Set $\mathtt{J}_{N+1} = \int_{\mathbb{R}} (1 +x^2 )^{-N -1}  \mathrm{d} x >0$. A direct computation shows that 
\begin{equation*}
\| \varphi_n \|_{H^s_+}^2 
 = \| h_n \|_{L^2_+}^2  = \tfrac{\pi}{n} \quad \text{and} \quad 
 \| D^N \varphi_n \|_{H^s_+}^2 
 = \| \partial_x^N h_n \|_{L^2_+}^2  = \mathtt{J}_{N+1} (N! )^2 n^{2N-1} ,  \quad \forall n \in \mathbb{N}_+.
\end{equation*}
Since $\lim_{n \to +\infty} \|\varphi_n\|_{H^s_+} = 0$ and $\lim_{n \to +\infty} \|D^N \varphi_n\|_{H^s_+} = +\infty$, the operator $D^N$ is unbounded on $H^s_+$.
\end{proof}
  
\noindent Now we introduce a uniform estimate for time-dependent Cauchy kernels whose poles lie in a compact subset of the lower half-plane $\mathbb{C}_-$. Such bounds will be essential when studying the long time asymptotics of dark multi-solitons in Section $\ref{SectTimeEvo6}$. 
 
\begin{lem}\label{BdvphiLongtimeLem}
Let $\mathtt{T}>1$ and $\mathtt{K}>1$ be given. Assume that a function $\mathfrak{z} : \mathbb{R} \to \mathbb{C}_-$  satisfies
\begin{equation}\label{ImztboundKKinv}
\tfrac{1}{\mathtt{K}} \leq | \operatorname{Im}\mathfrak{z} (t) | = - \operatorname{Im}\mathfrak{z}  (t) \leq \mathtt{K}, \quad \forall t \in \mathbb{R}\setminus [-\mathtt{T}, \mathtt{T}].
\end{equation}Set $\boldsymbol{\varphi}(t,x):=\boldsymbol{\varphi}_{\mathfrak{z}(t)} (x)  = (x - \mathfrak{z}(t))^{-1}$, $\forall (t,x) \in \mathbb{R} \times \mathbb{R}$. Then for any $s \geq 0$,  $\exists \mathrm{B}(s, \mathtt{K})>1$ such that
\begin{equation}\label{Longtivaphitbound}
\mathrm{B}(s, \mathtt{K})^{-1} \leq \|\boldsymbol{\varphi}(t)\|_{H^s_+} \leq \mathrm{B}(s, \mathtt{K}), \quad \forall t \in \mathbb{R}\setminus [-\mathtt{T}, \mathtt{T}].
\end{equation}
\end{lem}
\begin{proof}
For any $t \in \mathbb{R}$, since $\operatorname{Im}\mathfrak{z}(t) <0$, we have $\boldsymbol{\varphi}(t) \in H^{\infty}_+$. Denote by $\widehat{\boldsymbol{\varphi}}(t, \xi)$ the Fourier transform of $\boldsymbol{\varphi}(t) $ with respect to the spatial variable $x$, then
\begin{equation}\label{FouTrx-ztinv}
\widehat{\boldsymbol{\varphi}}(t, \xi) = -2 \pi ie^{-i \mathfrak{z}(t) \xi} \mathbf{1}_{\xi \geq 0}, \quad \forall \xi \in \mathbb{R}.
\end{equation}For any $s>0$, the constants $\mathtt{C}(s, \mathtt{K} )>0$ and $\widetilde{\mathtt{C}}(s, \mathtt{K}) >0$
are defined as follows:
\begin{equation}\label{CsktildCskfor}
\mathtt{C}(s, \mathtt{K})^2:=2 \pi \int_0^{\infty} (1+\xi^2)^s e^{-2 \mathtt{K}^{-1}\xi}\mathrm{d}\xi, \quad \left(\widetilde{\mathtt{C}}(s, \mathtt{K})\right)^2:=2 \pi \int_0^{\infty} (1+\xi^2)^s e^{-2 \mathtt{K} \xi}\mathrm{d}\xi.
\end{equation}For any $ t \in \mathbb{R}\setminus [-\mathtt{T}, \mathtt{T}]$, using formulas \eqref{ImztboundKKinv} and \eqref{CsktildCskfor}, we obtain that
\begin{equation}\label{tilCsmaHsnorphismalC}
\left(\widetilde{\mathtt{C}}(s, \mathtt{K})\right)^2 \leq \| \boldsymbol{\varphi} (t)\|_{H^s_+}^2 = 2 \pi \int_0^{+\infty}(1+\xi^2)^s \exp(2 \xi \operatorname{Im}\mathfrak{z}(t)) \mathrm{d}\xi \leq \mathtt{C}(s, \mathtt{K})^2.
\end{equation}Set $\mathrm{B}(s, \mathtt{K}) : = \max\{\mathtt{C}(s, \mathtt{K}), \left(\widetilde{\mathtt{C}}(s, \mathtt{K})\right)^{-1}\}+1$. Then formula \eqref{Longtivaphitbound}  follows from \eqref{tilCsmaHsnorphismalC}.
\end{proof}

\noindent For any positive integer $n \geq 1$, the \textit{Hardy--Zhidkov space} $\mathcal{Z}_+^n$ is defined as the subspace of functions in the Zhidkov space $\mathcal{Z}^n$ whose Fourier transforms are supported on  $\mathbb{R}_+ =[0, +\infty)$, i.e. 
\begin{equation}\label{HarZhidSpZnDefFor}
\mathcal{Z}_+^n  :=\{u \in \mathcal{Z}^n : \operatorname{supp} \widehat{u} \subset \mathbb{R}_+\} = \{u \in L^{\infty}(\mathbb{R}; \mathbb{C}) : \partial_x u \in H^{n-1}_+ \} = \mathcal{Z}^n \cap L^{\infty}_+.
\end{equation}Let $\mathcal{Z}_-^n$ denote the subspace of  $\mathcal{Z}^n$ consisting of functions whose Fourier transforms are supported on the negative real line $\mathbb{R}_-=(-\infty, 0]$, i.e.  $\mathcal{Z}_-^n :=\{u \in \mathcal{Z}^n : \operatorname{supp} \widehat{u} \subset \mathbb{R}_-\}$. For any $f \in \mathcal{Z}^n$, we have 
\begin{equation}\label{ZnnegisbarZnposi}
f \in \mathcal{Z}^n_- \Longleftrightarrow \overline{f}  \in \mathcal{Z}^n_+ \Longleftrightarrow \partial_x f \in H^{n-1}_- ,
\end{equation}where $\overline{f}(x) =\overline{f(x)} $ denotes the complex conjugate of $f(x)$. The Sobolev embedding  \eqref{InffH1overfLinfty}  implies that
\begin{equation}\label{Hn+subsetZn+}
H^n_+    \subset \mathcal{Z}^n \cap   L^{\infty}_+ = \mathcal{Z}^n_+ \subset H^n_{\mathrm{loc}}(\mathbb{R}; \mathbb{C})   \cap  \mathcal{Z}^1 , \quad \forall n \in \mathbb{N}_+.
\end{equation}If a tempered distribution $\Lambda \in \mathscr{S}'$ has Fourier transform supported at $\{0\}$, then $\Lambda$ is a polynomial. So 
\begin{equation}\label{Z1poscapZ1neg}
\mathcal{Z}^n_+ \cap  \mathcal{Z}^n_- = \{ u \in \mathcal{Z}^n : \operatorname{supp} \widehat{u}  \subset \{0\}\} = \mathbb{C}\mathbf{1}, \quad \forall n \in \mathbb{N}_+.
\end{equation} 
\begin{lem}\label{HarZhidisBanachLem}
For any integer $n \geq 1$, the Hardy--Zhidkov spaces $\mathcal{Z}^n_+$ and $\mathcal{Z}^n_-$ are closed subspaces of the Zhidkov space $\left(\mathcal{Z}^n, \|.\|_{\mathcal{Z}^n} \right)$. Moreover, we have  
\begin{equation}\label{ZnHardBanAlgebra}
 L^{\infty}_+ \cdot L^2_+ \subset L^2_+ , \quad \mathcal{Z}^n_+ \cdot H^n_+ \subset H^n_+  \quad \mathrm{and} \quad \mathcal{Z}^n_+ \cdot \mathcal{Z}^n_+  \subset \mathcal{Z}^n_+ .
\end{equation}Consequently, the Hardy--Zhidkov space $\left(\mathcal{Z}^n_+, \mathfrak{C}_n \|.\|_{\mathcal{Z}^n} \right)$ is a commutative Banach algebra with identity $\mathbf{1} \in  \mathcal{Z}^n_+ \cap  \mathcal{Z}^n_-$. The Hardy--Sobolev space $H^n_+$   is an ideal of the ring $\left(\mathcal{Z}^n_+, +, \times \right)$.  
\end{lem}   

\begin{proof}Fix $n \in \mathbb{N}_+$. Assume that $\phi \in \mathcal{Z}^n$ and $\phi_m \in \mathcal{Z}^n_+$ such that $\lim_{m \to +\infty}\|\phi_m - \phi\|_{\mathcal{Z}^n} =0$. By formula \eqref{ZhidNormZn}, the sequence  $(\widehat{\partial_x \phi_m})_{m \in \mathbb{N}}$ admits a subsequence that converges to $\widehat{\partial_x \phi}$ almost everywhere in $\mathbb{R}$. Since $\partial_x \phi_m \in H^{n-1}_+$, we have $\operatorname{supp} (\widehat{\partial_x \phi}) \subset \mathbb{R}_+ = [0, +\infty)$ and  $\partial_x\phi \in H^{n-1}_+$. So $\left(\mathcal{Z}^n_+, \|.\|_{\mathcal{Z}^n} \right)$ is a Banach space.   
 By Lemma $\ref{Z1+idealinL2+}$,  we have $\mathcal{Z}^n_+ \cdot H^n_+ \subset L^{\infty}_+ \cdot L^2_+ \subset L^2_+$. So formula  \eqref{ufHn<uZnfHn} implies that 
\begin{equation}\label{ZnHardyHninHnfor}
\mathcal{Z}^n_+ \cdot H^n_+ \subset L^2_+ \cap  \left(\mathcal{Z}^n  \cdot H^n  \right)\subset L^2_+ \cap  H^n (\mathbb{R} ; \mathbb{C}) =H^n_+.
\end{equation}For any $u, v \in \mathcal{Z}^n_+$, Lemma $\ref{ZnBanachAlgeLem}$ yields that $uv \in \mathcal{Z}^n$. Using   \eqref{d(uv)foruvZ1} and \eqref{ZnHardyHninHnfor}, we obtain that
\begin{equation*}
\partial_x (uv) = v \partial_x u + u \partial_x v  \in \mathcal{Z}^n_+ \cdot H^n_+  + \mathcal{Z}^n_+ \cdot H^n_+ \subset H^n_+.
\end{equation*}Then formula \eqref{HarZhidSpZnDefFor} yields that $\operatorname{supp}\left( \widehat{uv} \right) \subset \mathbb{R}_+$ and $uv \in \mathcal{Z}^n_+$. So $\mathcal{Z}^n_+ \cdot \mathcal{Z}^n_+ \subset \mathcal{Z}^n_+$.
\end{proof} 
\begin{cor}\label{Zn-corBanaalgebra}
For any integer $n\geq 1$, the Hardy--Zhidkov space $\left(\mathcal{Z}^n_-, \mathfrak{C}_n \|.\|_{\mathcal{Z}^n} \right)$ is   a commutative Banach algebra with identity $\mathbf{1}$.  The Hardy--Sobolev space $H^n_-$   is an ideal of the ring $\left(\mathcal{Z}^n_-, +, \times \right)$, i.e.
\begin{equation}\label{BanAlgebraZn-}
  L^{\infty}_- \cdot L^2_- \subset L^2_- , \quad\mathcal{Z}^n_- \cdot H^n_- \subset H^n_-  \quad \mathrm{and} \quad \mathcal{Z}^n_- \cdot \mathcal{Z}^n_- \subset \mathcal{Z}^n_- .
\end{equation}
\end{cor}
\begin{proof}
The proof is given by formula \eqref{ZnnegisbarZnposi} and Lemma $\ref{HarZhidisBanachLem}$.
\end{proof}

\noindent The Hardy--Zhidkov space $\mathcal{Z}^n_+$ inherits the properties of the Schr\"odinger group $(e^{it\partial_x^2})_{t\in\mathbb{R}}$ established for the full Zhidkov space $\mathcal{Z}^n$ in Lemma~\ref{StrCSchGrZhidlem} and Lemma~\ref{GenSchGrlemZN}.

\begin{lem}\label{lemSchGrHardZn}
For any $n \in \mathbb{N}_+$ and $\mathfrak{t} \in \mathbb{R}$, we have $e^{i \mathfrak{t} \partial_x^2} - \mathrm{id}_{\mathcal{Z}^n_+} \in \mathcal{B}_{\mathbb{C}} (\mathcal{Z}^n_+ ; H^n_+)$. The Schr\"odinger group $(e^{it\partial_x^2})_{t\in\mathbb{R}}$ is a strongly continuous group of bounded operators on the Hardy--Zhidkov space $\bigl(\mathcal{Z}^n_+, \|\cdot\|_{\mathcal{Z}^n_+} \bigr)$. Its infinitesimal generator is given by $\bigl( i \partial_x^2 , \mathrm{Dom}_{\mathcal{Z}^n_+} (i \partial_x^2) \bigr)$ with $\mathrm{Dom}_{\mathcal{Z}^n_+} (i \partial_x^2) = \mathcal{Z}^{n+2}_+$.
\end{lem}

\begin{proof}
For any $f \in \mathcal{Z}^n_+$ and $\mathfrak{t} \in \mathbb{R}$, we have $\partial_x e^{i \mathfrak{t} \partial_x^2}(f) = e^{i \mathfrak{t} \partial_x^2} (\partial_x f) \in H^{n-1}_+$. Then Lemma~\ref{StrCSchGrZhidlem} yields that $e^{i \mathfrak{t} \partial_x^2}(f) \in \mathcal{Z}^n_+$ and $e^{i \mathfrak{t} \partial_x^2}(f) - f \in H^n_+$. The proof is completed by Lemma~\ref{GenSchGrlemZN}.
\end{proof}

\noindent For any $n \in \mathbb{N}_+$, we have $\mathcal{Z}^{n+1} \subset \mathcal{Z}^n$. The smooth Hardy--Zhidkov spaces $\mathcal{Z}^{\infty}_+$, $\mathcal{Z}^{\infty}_-$ are denoted by 
\begin{equation}\label{SmoothZhidkHa}
\mathcal{Z}^{\infty}_+:= \bigcap_{n \in \mathbb{N}_+} \mathcal{Z}^n_+ = \mathcal{Z}^{\infty} \cap L^{\infty}_+, \qquad
\mathcal{Z}^{\infty}_-:= \bigcap_{n \in \mathbb{N}_+} \mathcal{Z}^n_- = \mathcal{Z}^{\infty} \cap L^{\infty}_-.
\end{equation}
The next lemma is a direct consequence of Lemma~\ref{ZinftydenseinZN} and Lemma~\ref{LinftplconvSchwLem}.
\begin{lem}\label{HarZhiSmoodense}
For any $n \in \mathbb{N}_+$, the smooth Hardy--Zhidkov space $\mathcal{Z}^{\infty}_+$ is dense in $\bigl( \mathcal{Z}^n_+, \|\cdot\|_{\mathcal{Z}^n_+} \bigr)$.
\end{lem}

\subsection{Toeplitz Operators and Shift Semigroups}\label{SubsShifSemiGr}

\noindent For any $u \in L^2(\mathbb{R};\mathbb{C}) + L^{\infty}(\mathbb{R};\mathbb{C})$,  the Toeplitz operator of symbol $u$ is denoted by $\mathbf{T}_u : H^1_+ \to  L^2_+$, i.e. 
\begin{equation}\label{ToepdefForSec2}
 \mathbf{T}_u (h):=\Pi (uh) = \Pi_{\geq 0} (uh) \in L^2_+, \quad \forall h \in H^1_+.
\end{equation}If $\mathfrak{a} \in L^{\infty}(\mathbb{R}; \mathbb{C})$, then the Toeplitz operator    $\mathbf{T}_{\mathfrak{a}}$ in \eqref{ToepdefForSec2} extends to a bounded $\mathbb{C}$-linear operator on the Hardy space  $L^2_+$. 
 Its $L^2_+$-adjoint is given by $\left( \mathbf{T}_{\mathfrak{a}}\right)^* =  \mathbf{T}_{\overline{\mathfrak{a}}} \in \mathcal{B}_{\mathbb{C}}(L^2_+)$. Specifically, we have
\begin{equation}\label{AdjofToepOpinL2+}
 \langle \mathbf{T}_{\mathfrak{a}}(f), g \rangle_{L^2_+} = \int_{\mathbb{R}} \mathfrak{a} f \overline{g} = \langle  f , \mathbf{T}_{\overline{\mathfrak{a}}}(g)  \rangle_{L^2_+}, \quad \forall f, g \in L^2_+.
\end{equation}Assume that $\mathfrak{a} \in L^{\infty}(\mathbb{R}; \mathbb{C})$, $\mathfrak{b} \in H^s (\mathbb{R}; \mathbb{C})$ and  $u \in \mathcal{Z}^n$ for some $s> \frac{1}{2}$ and $n \in \mathbb{N}_+$,  we have
\begin{equation}\label{BounToeplitzBHs}
\|\mathbf{T}_{\mathfrak{a}}\|_{\mathcal{B}_{\mathbb{C}}(L^2_+)} \leq  \|\mathfrak{a}\|_{L^{\infty}(\mathbb{R}; \mathbb{C})}, \quad \|\mathbf{T}_{\mathfrak{b}}\|_{\mathcal{B}_{\mathbb{C}}(H^s_+)} \leq \mathfrak{C}_s' \|\mathfrak{b}\|_{H^s(\mathbb{R}; \mathbb{C})} \quad \mathrm{and} \quad \|\mathbf{T}_{u}\|_{\mathcal{B}_{\mathbb{C}}(H^n_+)} \leq C_n \|u\|_{\mathcal{Z}^n}.
\end{equation}by using  \eqref{HNAlgebra}, 
\eqref{ufHn<uZnfHn} and \eqref{ToepdefForSec2}. 
For any $\eta \geq 0$, $x\in \mathbb{R}$, $z \in \mathbb{C}$ and $f \in L^2_+$, we define that  
\begin{equation}\label{LaxBeurlingSG}
\mathbf{e}_{z}(x): = e^{i z x}\in \mathbb{C}, \quad \mathbf{S}(\eta)f := \mathbf{T}_{\mathbf{e}_{\eta}}(f)= \mathbf{e}_{\eta}f\in L^2_+,  \quad \mathbf{S}(\eta)^*f := \mathbf{T}_{\mathbf{e}_{-\eta}}(f)=\Pi  \left(\mathbf{e}_{-\eta}f \right)\in L^2_+,  
\end{equation}where $\mathbf{S}(\eta)^* \in \mathcal{B}_{\mathbb{C}}(L^2_+)$ denotes the $L^2_+$-adjoint of $\mathbf{S}(\eta) \in \mathcal{B}_{\mathbb{C}}(L^2_+)$. Then the following lemma holds.

\begin{lem}\label{ContSemiGrLem}
For any $s \geq 0$, both $(\mathbf{S}(\eta))_{\eta \geq 0}$ and $(\mathbf{S}(\eta)^*)_{\eta \geq 0}$ are strongly continuous semigroups on the Hardy--Sobolev space $H^s_+$. Moreover, $(\mathbf{S}(\eta))_{\eta \geq 0}$ is  a contraction semigroup on $L^2_+$ and $(\mathbf{S}(\eta)^*)_{\eta \geq 0}$ is a contraction semigroup on $H^s_+$, $\forall s \geq 0$.  
\end{lem}

\begin{proof}
For any $\xi_1, \xi_2 \geq 0$ and $f \in L^2_+$, we have $\mathbf{S}(\xi_1) \mathbf{S}(\xi_2) f = \mathbf{e}_{\xi_1+\xi_2}f = \mathbf{S}(\xi_1+ \xi_2) f \in L^2_+$. So  
\begin{equation}\label{SemigrIdent}
\mathbf{S}(\xi_1)^*\mathbf{S}(\xi_2)^* = \left(\mathbf{S}(\xi_2)\mathbf{S}(\xi_1) \right)^* =  \mathbf{S}(\xi_1+\xi_2)^* \in \mathcal{B}_{\mathbb{C}}(L^2_+) 
\end{equation}and $\mathbf{S}(0)=\mathrm{id}_{L^2_+}$. Fix $s \geq 0$ and $\phi \in H^s_+$. Peetre's inequality \eqref{PeetreIne} and formula \eqref{Lpmussp} imply that $\forall \eta > 0$,
\begin{equation}\label{SetainBHs}
\small
\|\mathbf{S}(\eta)\phi\|_{H^s_+}^2 = \int_{0}^{+\infty} |\widehat{\phi}(y)|^2 (1 + |y +\eta|^2)^s \frac{\mathrm{d}y}{2 \pi }   \leq \frac{(4 + 4\eta^2)^s}{2 \pi } \int_0^{+\infty} |\widehat{\phi}(y)|^2(1+y^2)^s \mathrm{d}y =(4 + 4\eta^2)^s \|\phi\|_{H^s_+}^2.
\end{equation}Assume that $\eta, \eta_n \in \mathbb{R}_+^*= (0, +\infty)$ such that $\lim_{n \to +\infty}\eta_n = \eta$. Set $t_n:= \eta_n -\eta$.  Then we have
\begin{equation}\label{Seta-SetanphiHs}
2 \pi \|\mathbf{S}(\eta_n)\phi - \mathbf{S}(\eta)\phi\|_{H^s_+}^2 \leq \int_{\mathbb{R}} |\widehat{\phi}(y -t_n)-\widehat{\phi}(y)|^2 \left(1+|y + \eta|^2\right)^s\mathrm{d}y \leq (4+4\eta^2)^s\|\tau_{t_n}\widehat{\phi}-\widehat{\phi}\|_{L^2_s(\mathbb{R})}^2,
\end{equation}by using  \eqref{Lpmussp} and \eqref{PeetreIne}. Lemma $\ref{TranslationCV}$ yields that $\lim_{n \to +\infty}\|\mathbf{S}(\eta_n)\phi - \mathbf{S}(\eta)\phi\|_{H^s_+} = 0$. So 
$\left(\mathbf{S}(\eta)\right)_{\eta \geq 0}$ is a strongly continuous semigroup on $H^s_+$, $\forall s\geq 0$. \\

\noindent For any $\xi, \eta \geq 0$, we have $1+|\xi|^2 \leq 1+|\xi+\eta|^2$. If $\varphi \in H^s_+$ for some $s\geq 0$, then $\forall \eta \geq 0$, we have
\begin{equation}\label{SetastcontraSemiGr}
2\pi \|\mathbf{S}(\eta)^*\varphi\|_{H^s_+}^2 \leq \int_0^{+\infty} |\widehat{\varphi}(\xi+\eta)|^2(1+|\xi+\eta|^2)^s \mathrm{d}\xi = \int_{\eta}^{+\infty} |\widehat{\varphi}(y)|^2(1+|y|^2)^s \mathrm{d}y \leq 2\pi \| \varphi\|_{H^s_+}^2.
\end{equation}So $\sup_{s \geq 0} \sup_{\eta \geq 0}\|\mathbf{S}(\eta)^*\|_{\mathcal{B}_{\mathbb{C}}(H^s_+)} = 1$  and $\sup_{\eta \geq 0}\|\mathbf{S}(\eta) \|_{\mathcal{B}_{\mathbb{C}}(L^2_+)} = 1$. We still assume that  $\eta, \eta_n \in \mathbb{R}_+^*$ such that $\lim_{n \to +\infty}\eta_n = \eta$, and $t_n:= \eta_n -\eta$.   Then 
\begin{equation}\label{StrContSeta*}
\small
2 \pi \|\mathbf{S}(\eta_n)^*\varphi - \mathbf{S}(\eta)^*\varphi\|_{H^s_+}^2 \leq \int_0^{+\infty}(1+|\xi +\eta|^2)^s |\widehat{\varphi}(\xi+\eta +t_n)-\widehat{\varphi}(\xi+\eta)|^2 \mathrm{d}\xi \leq \|\tau_{-t_n}\widehat{\varphi}-\widehat{\varphi}\|_{L^2_s(\mathbb{R})}^2.
\end{equation}The strong continuity of $\left(\mathbf{S}(\eta)^*\right)_{\eta \geq 0}$ is given by \eqref{StrContSeta*} and Lemma $\ref{TranslationCV}$.  The proof is completed.
\end{proof}
\begin{rem}
Fix $s>0$ and $\eta>0$. Assume that $\phi \in H^s_+$ such that $\widehat{\phi}|_{\mathbb{R}_+^*} \in \mathcal{C}(\mathbb{R}_+^*, \mathbb{C})$ and $\widehat{\phi}(\xi)\ne 0$ for any $\xi >0$. Since $1 + y^2 < 1+|y +\eta|^2$ holds for all $y \geq 0$, formula \eqref{SetainBHs} yields that  $\|\mathbf{S}(\eta)\phi\|_{H^s_+} > \|\phi\|_{H^s_+}$. So $(\mathbf{S}(\eta))_{\eta \geq 0}$ is not a contraction semigroup on the Hardy--Sobolev space $H^s_+$, for all $s>0$.
\end{rem}
 
\noindent The semigroups $(\mathbf{S}(\eta))_{\eta \geq 0}$ and $(\mathbf{S}(\eta)^*)_{\eta \geq 0}$ given by formula  \eqref{LaxBeurlingSG} originate from Lax's study of translation invariant spaces \cite{LaxTransinva1959} and are the natural counterpart on the line of Beurling's characterization of invariant subspaces of the Hardy space on the unit disc. In view of this, we call them the \emph{Lax--Beurling right and left shift semigroups}. \\

\noindent The next lemma identifies the infinitesimal generator of the right shift semigroup $(\mathbf{S}(\eta))_{\eta \geq 0}$ on $L^2_+$.

\begin{lem}\label{X=InfiniGeL2SetaLem}
A function $\varphi \in L^2_+$ belongs to $\mathrm{Dom}(\mathbf{X})$ if and only if the limit 
$
\lim_{\eta \to 0^+} \frac{1}{\eta}\bigl( \mathbf{S}(\eta) \varphi - \varphi \bigr)
$
exists in $L^2_+$. In this case, the limit equals $i\mathbf{X}\varphi \in L^2_+$. Consequently, the operator $(i\mathbf{X}, \mathrm{Dom}(\mathbf{X}))$ defined in \eqref{defXf-} and \eqref{DomX} is the infinitesimal generator of $(\mathbf{S}(\eta))_{\eta \ge 0}$ on $L^2_+$.
\end{lem}
\begin{proof}
For any $\varphi, g \in L^2_+$ and $\eta >0$, we have $\widehat{\mathbf{S}(\eta)\varphi}(\xi)=\widehat{\mathbf{e}_{\eta} \varphi}(\xi)=\widehat{\varphi}(\xi-\eta)$, $\forall \xi\in \mathbb{R}$, and
\begin{equation}\label{PlanSetaphi-phiid}
2 \pi \|\tfrac{1}{\eta}\left(\mathbf{S}(\eta)\varphi - \varphi \right) - g\|_{L^2_+}^2 = \int_{\mathbb{R}}|\tfrac{1}{\eta}(\widehat{\varphi}(\xi - \eta) - \widehat{\varphi}(\xi))-\widehat{g}(\xi)|^2\mathrm{d}\xi = \|\tfrac{1}{-\eta}(\tau_{-(-\eta)} \widehat{\varphi} - \widehat{\varphi})+ \widehat{g}\|_{L^2(\mathbb{R}; \mathbb{C})}^2.
\end{equation}$(\mathbf{ \romannumeral 1})$. If $\varphi \in \mathrm{Dom}(\mathbf{X})$, then $\widehat{\varphi} \in H^1(\mathbb{R}; \mathbb{C}) = W^{1,2}_0(\mathbb{R})$ and $i \widehat{\mathbf{X}\varphi}=-\partial_{\xi} \widehat{\varphi}\in L^2(\mathbb{R}; \mathbb{C})$. By Lemma $\ref{CharSobDiffQuoCV}$ and   \eqref{PlanSetaphi-phiid}, we have $\sqrt{2 \pi}\|\tfrac{1}{\eta}\left(\mathbf{S}(\eta)\varphi - \varphi \right) - i\mathbf{X}\varphi\|_{L^2_+}= \|\tfrac{1}{-\eta}(\tau_{-(-\eta)} \widehat{\varphi} - \widehat{\varphi}) - \partial_{\xi} \widehat{\varphi}\|_{L^2(\mathbb{R}; \mathbb{C})} \to 0$, as $\eta \to 0^+$.\\
\noindent $(\mathbf{ \romannumeral 2})$. Assume that $\varphi, g \in L^2_+$ and  $\lim_{n \to +\infty} \|\tfrac{1}{\eta_n}\left(\mathbf{S}(\eta_n)\varphi - \varphi \right) - g\|_{L^2_+}=0$, where $\eta_n := n^{-1} \to  0^+$, as $n \to +\infty$. 
By \eqref{PlanSetaphi-phiid} and  Lemma $\ref{CharSobDiffQuoCV}$, we have $\widehat{\varphi}\in  W^{1,2}_0(\mathbb{R}) = H^1(\mathbb{R}; \mathbb{C})$. So $\varphi \in \mathrm{Dom}(\mathbf{X})$.
\end{proof}
 
\noindent For any $s\geq 0$, let $H^1_s(\mathbb{R}_+^*) =W^{1,2}_s(\mathbb{R}_+^*)$  be defined in \eqref{WeightedSobolevSp} with $\mathbb{R}_+^* = (0 , +\infty)$. We have
\begin{equation}\label{H1snormWeiSobo}
\|f\|_{H^1_s(\mathbb{R}_+^*)}^2 = \|f\|^2_{L^2_s(\mathbb{R}_+^*)}+ \|\partial_{\xi} f\|^2_{L^2_s(\mathbb{R}_+^*)} =\int_0^{+\infty} \left(|f(\xi)|^2+ |\partial_{\xi} f(\xi)|^2 \right) (1+\xi^2)^s \mathrm{d}\xi \in [0, +\infty],
\end{equation}by formula \eqref{H1sWeiSobnormDef}.
 A measurable function $f$ on $\mathbb{R}_+^*$ belongs to $H^1_s(\mathbb{R}_+^*)$ if and only if $\|f\|_{H^1_s(\mathbb{R}_+^*)} < +\infty$. 
  Let $\mathrm{Dom}^s(\mathbf{G})$ be the set of all functions $\varphi \in H^s_+ \cap \mathrm{Dom}(\mathbf{G})$ for which $\mathbf{G}\varphi \in H^s_+$. Specifically, we have   
\begin{equation}\label{DomsGDefSect2}
\mathrm{Dom}^s(\mathbf{G}) : =  H^s_+ \cap  \mathrm{Dom}(\mathbf{G}) \cap  \mathbf{G}^{-1}\left( H^s_+ \right) = \{\varphi \in L^2_+: \widehat{\varphi}|_{\mathbb{R}_+^*} \in H^1_s(\mathbb{R}_+^*)\} ,
\end{equation}by Definition $\ref{DombfGdefSect2}$, formulas \eqref{DefGphiSect2} and  \eqref{H1snormWeiSobo}.
 For any $r> s  \geq 0$, we have
  \begin{equation}\label{DomrGsubsetDomsG}
\mathrm{Dom}^{\infty}(\mathbf{G}) := \bigcap_{\alpha \geq 0} \mathrm{Dom}^{\alpha}(\mathbf{G}) \subset \mathrm{Dom}^r (\mathbf{G}) \subset \mathrm{Dom}^s (\mathbf{G}) \subset \mathrm{Dom}^0 (\mathbf{G}) =\mathrm{Dom} (\mathbf{G}).
\end{equation}

\begin{prop}\label{GisInfiniGeneHsProp}
For any $s\geq 0$, the operator $\left(-i \mathbf{G}, \mathrm{Dom}^s(\mathbf{G}) \right)$ is the infinitesimal generator of the contraction semigroup $(\mathbf{S}(\eta)^*)_{\eta \geq 0}$ on the Hardy--Sobolev space $H^s_+$. Specifically, we have
\begin{equation}\label{SemiGrCVdefDomsG}
 \mathrm{Dom}^s(\mathbf{G}) =\{\phi \in H^s_+: \exists v\in H^s_+ \quad \mathrm{such} \quad \mathrm{that} \quad \lim_{\eta \to 0^+}\|\tfrac{1}{\eta}(\mathbf{S}(\eta)^* \phi - \phi)-v\|_{H^s_+}=0\}.
\end{equation}Furthermore, if $\phi \in \mathrm{Dom}^s(\mathbf{G})$, then $\tfrac{1}{\eta}(\mathbf{S}(\eta)^* \phi - \phi) \to -i \mathbf{G} \phi$ in $H^s_+$, as $\eta \to 0^+$.
\end{prop}
\begin{rem}
The exponential expression for the left shift semigroup is given by
\begin{equation}\label{ForSeta*=e-ietaG}
\mathbf{S}(\eta)^*= \mathbf{T}_{\mathbf{e}_{-\eta}}= \exp(-i\eta \mathbf{G}) \in \bigcap_{s \geq 0}\mathcal{B}_{\mathbb{C}}(H^s_+), \quad \forall \eta \geq 0.
\end{equation}
\end{rem}

\begin{proof}[Proof of Proposition $\ref{GisInfiniGeneHsProp}$]
The proof proceeds in two steps.\\

\noindent \textbf{Step \romannumeral 1}. First, we assume that $s=0$. If $\varphi \in \mathrm{Dom}(\mathbf{G})$, then  $\widehat{\varphi}|_{\mathbb{R}_+^*} \in H^1(\mathbb{R}_+^*)$. For any $\xi, \eta >0$,  we have $\widehat{\mathbf{G}  \varphi }(\xi) = i \partial_{\xi} \widehat{\varphi}(\xi)$  and $\widehat{\mathbf{S}(\eta)^*\varphi}(\xi) = \widehat{\varphi}(\xi +\eta)$. By Plancherel identity and Corollary $\ref{QuoDiffCVR+*Cor}$, we have
\begin{equation}\label{Set*phi-phi+iGphi}
\sqrt{2\pi} \|\tfrac{1}{\eta}(\mathbf{S}(\eta)^* \varphi - \varphi) + i \mathbf{G}\varphi\|_{L^2_+} = \|\tfrac{1}{\eta}(\tau_{-\eta}\widehat{\varphi} - \widehat{\varphi}) - \partial_{\xi} \widehat{\varphi}  \|_{L^2 (\mathbb{R}_+^*; \mathbb{C})}  \to 0, \quad \mathrm{when} \quad  \eta \to 0^+.
\end{equation}Conversely, assume that $\lim_{\eta \to 0^+}\|\tfrac{1}{\eta}(\mathbf{S}(\eta)^* f - f) - h\|_{L^2_+}=0$ for some $f, h \in L^2_+$. Then we have
\begin{equation}\label{Pleta-1(S*f-f)-h}
\small
  \|\tfrac{1}{\eta}(\tau_{-\eta} \widehat{f} - \widehat{f}) - \widehat{h}\|_{L^2(\mathbb{R}_+^*; \mathbb{C})}^2 = \int_0^{+\infty}\Big|\tfrac{1}{\eta}\left(\widehat{f}(\xi + \eta) - \widehat{f}(\xi)\right) - \widehat{h}(\xi)\Big|^2 \mathrm{d}\xi = 2\pi \|\tfrac{1}{\eta}(\mathbf{S}(\eta)^* f - f) - h\|_{L^2_+}^2 \to 0,
\end{equation}when $\eta \to 0^+$, by using Plancherel identity. If $\psi \in \mathcal{C}_c^1(\mathbb{R}; \mathbb{C})$ such that $\operatorname{supp}\psi \subset \mathbb{R}_+^*=(0, +\infty) $, we have
\begin{equation}\label{cvdiffquopsi}
\lim_{\eta \to 0^+}\| \tfrac{1}{\eta}(\tau_{\eta} \psi  - \psi) + \partial_{\xi} \psi\|_{L^2(\mathbb{R}; \mathbb{C})} = 0,
\end{equation}by Lemma $\ref{CharSobDiffQuoCV}$. Since $\psi|_{\mathbb{R}_-}=0$ with $\mathbb{R}_-=(-\infty, 0]$, we use \eqref{Pleta-1(S*f-f)-h} and  \eqref{cvdiffquopsi} to obtain that
\begin{equation}\label{IntbyPartfhpsi}
-\int_{\mathbb{R}_+^*}\widehat{f} \partial_{\xi} \psi = \lim_{\eta \to 0^+} \tfrac{1}{\eta}\int_{\mathbb{R}}\widehat{f} (\tau_{\eta} \psi - \psi)  =\lim_{\eta \to 0^+} \tfrac{1}{\eta} \int_0^{+\infty}(\tau_{-\eta}\widehat{f} - \widehat{f})  \psi =\int_{\mathbb{R}_+^*} \widehat{h}\psi.
\end{equation}Therefore, formula \eqref{WeightedSobolevSp} implies that $\widehat{f}|_{\mathbb{R}_+^*} \in W^{1,2}_0(\mathbb{R}_+^*) = H^1(\mathbb{R}_+^*)$, and hence $f \in \mathrm{Dom}(\mathbf{G})$. We have proved Proposition $\ref{GisInfiniGeneHsProp}$ in the case $s=0$.\\

\noindent \textbf{Step \romannumeral 2}. Fix $s>0$. If $\varphi \in \mathrm{Dom}^s(\mathbf{G})$, then \eqref{DomsGDefSect2} yields that $\widehat{\varphi}|_{\mathbb{R}_+^*} \in H^1_s(\mathbb{R}_+^*)=W^{1,2}_s(\mathbb{R}_+^*)$. According to Corollary $\ref{QuoDiffCVR+*Cor}$ and Plancherel identity, we use formula  \eqref{Lpmussp}  to obtain that
\begin{equation}\label{Seta*p-p+iGphiHsEst}
\begin{split}
2\pi \|\tfrac{1}{\eta}(\mathbf{S}(\eta)^* \varphi - \varphi) + i \mathbf{G}\varphi\|_{H^s_+}^2 = & \int_0^{+\infty}|\tfrac{1}{\eta}(\widehat{\varphi}(\xi+\eta)-\widehat{\varphi}(\xi)) - \partial_{\xi} \widehat{\varphi}(\xi)|^2 (1+|\xi|^2)^s\mathrm{d}\xi \\
= & \|\tfrac{1}{\eta}(\tau_{-\eta}\widehat{\varphi} - \widehat{\varphi}) - \partial_{\xi} \widehat{\varphi}  \|_{L^2_s (\mathbb{R}_+^*)}^2 \to 0 , \quad \mathrm{as} \quad \eta \to 0^+.
\end{split}
\end{equation}Conversely, assume that there exist $f, h \in H^s_+$ such that $\lim_{\eta \to 0^+}\|\tfrac{1}{\eta}(\mathbf{S}(\eta)^* f - f) - h\|_{H^s_+}=0$. Then 
\begin{equation}\label{L2norm<=Hsnorm1}
 \|\tfrac{1}{\eta}(\mathbf{S}(\eta)^* f - f) - h\|_{L^2_+} \leq \|\tfrac{1}{\eta}(\mathbf{S}(\eta)^* f - f) - h\|_{H^s_+} \to 0, \quad \mathrm{as} \quad \eta \to 0^+.
\end{equation}According to  
\textbf{Step \romannumeral 1}, Proposition \ref{GisInfiniGeneHsProp}  holds in the case $s=0$. Then   \eqref{SemiGrCVdefDomsG} yields that
\begin{equation}\label{Forh=-iGf}
f \in \mathrm{Dom}^0(\mathbf{G}) = \mathrm{Dom}(\mathbf{G})  \quad \mathrm{and} \quad h = \lim_{\eta \to 0^+} \tfrac{1}{\eta} (\mathbf{S}(\eta)^* f -f) = -i \mathbf{G}f \in L^2_+.
\end{equation}Since $ h \in H^s_+$, we have $f \in  H^s_+ \cap  \mathrm{Dom}(\mathbf{G}) \cap  \mathbf{G}^{-1}\left( H^s_+ \right) =\mathrm{Dom}^s(\mathbf{G}) $  by using \eqref{DomsGDefSect2} and \eqref{Forh=-iGf}. 
\end{proof} 

\begin{rem}\label{DifferenceX&GRem}
Assume that $f \in L^2_+$   such that $\widehat{f}$ has a discontinuity at $0$. Then we have $\widehat{f} \notin H^1(\mathbb{R}; \mathbb{C})$. By Lemma $\ref{CharSobDiffQuoCV}$, for any sequence $(\eta_n)_{n \geq 1}$ with $\lim_{n \to +\infty}\eta_n = 0$, the difference quotients $\left(\eta_n^{-1}(\tau_{-\eta_n}\widehat{f} -\widehat{f})\right)_{n \geq 1}$ fail to converge in $L^2(\mathbb{R}; \mathbb{C})$. The key observation is that for test functions $\psi \in \mathcal{C}_c^1 (\mathbb{R}; \mathbb{C})$ whose support intersects $\mathbb{R}_-^*=(-\infty,0)$,  the integral $\int_{\mathbb{R}}(\tau_{-\eta}\widehat{f} - \widehat{f})  \psi$ may differ from $\int_0^{+\infty}(\tau_{-\eta}\widehat{f} - \widehat{f})  \psi$ due to the discontinuity of $\widehat{f}$ at $0$, when $\eta>0$. Consequently, even if \eqref{Pleta-1(S*f-f)-h} holds for some $h \in L^2_+$, there exists $\phi \in \mathcal{C}_c^1 (\mathbb{R}; \mathbb{C})$ with $\phi(\xi_0) \ne 0$ for some $\xi_0 < 0$ such that $-\int_{\mathbb{R}}\widehat{f} \partial_{\xi} \phi \ne \int_{\mathbb{R}} \widehat{h}\phi $. Specifically, we have
\begin{equation}\label{f^notinH1RFor}
\small
-\int_{\mathbb{R}}\widehat{f} \partial_{\xi} \phi 
= \lim_{\eta \to 0^+} \tfrac{1}{\eta}\int_{\mathbb{R}}\widehat{f} (\tau_{\eta} \phi - \phi)  
= \lim_{\eta \to 0^+} \tfrac{1}{\eta} \int_{\mathbb{R}}(\tau_{-\eta}\widehat{f} - \widehat{f})  \phi 
\ne   \lim_{\eta \to 0^+} \tfrac{1}{\eta} \int_0^{+\infty}(\tau_{-\eta}\widehat{f} - \widehat{f})  \phi =\int_{\mathbb{R}} \widehat{h}\phi .
\end{equation}In contrast, formula \eqref{IntbyPartfhpsi} holds for any $\psi \in \mathcal{C}_c^1 (\mathbb{R}; \mathbb{C})$ with $\operatorname{supp}(\psi) \subset \mathbb{R}_+^*$, under  assumption \eqref{Pleta-1(S*f-f)-h}. This demonstrates the essential difference in the proof strategies between Lemma $\ref{X=InfiniGeL2SetaLem}$ and Proposition $\ref{GisInfiniGeneHsProp}$.
\end{rem}

\noindent For any $s \geq 0$, the operator $\left(i\mathbf{G}, \mathrm{Dom}^s(\mathbf{G}) \right)$ is maximal accretive on $H^s_+$. In particular, it is closed and densely defined. If $z \in \mathbb{C}_+=\{a + bi \in \mathbb{C}:  b >0\}$, then $\mathbf{G} - z: \mathrm{Dom}^s(\mathbf{G}) \to H^s_+$ is bijective, and 
\begin{equation}\label{EstG-z-1BHs}
\sup_{s \geq 0}\|(\mathbf{G}-z)^{-1}\|_{\mathcal{B}_{\mathbb{C}}(H^s_+)} \leq \tfrac{1}{\operatorname{Im} z}.
\end{equation}Its accretiveness can also be proven directly by the next lemma.
\begin{lem}\label{ImGffinprodLem}
Assume that $\varphi \in \mathrm{Dom}^s(\mathbf{G})$ for some $s \geq 0$, then we have
\begin{equation}\label{fo<Gphi,phi>Hs}
2 \pi \operatorname{Re} \langle i \mathbf{G}\varphi , \varphi \rangle_{H^s_+} = \frac{1}{2} |\widehat{\varphi}(0^+)|^2 + s \int_0^{+\infty} |\widehat{\varphi}(\xi)|^2 (1+\xi^2)^{s-1} \xi \mathrm{d} \xi.
\end{equation}In particular, for any $f\in \mathrm{Dom}(\mathbf{G})$, we have $\operatorname{Im} \langle  \mathbf{G} f , f \rangle_{L^2_+} = - \frac{1}{4 \pi}|\widehat{f}(0^+)|^2 \leq 0$.
\end{lem}
\begin{proof}
Fix $s \geq 0$ and $\varphi \in \mathrm{Dom}^s(\mathbf{G})$, set $h(\xi)=h_{\varphi, s}(\xi):=\widehat{\varphi}(\xi)\langle \xi \rangle^s$, $\forall \xi \in \mathbb{R}$. Then $h \in H^1(\mathbb{R}_+^*)$,
\begin{equation}\label{fohL2norm}
\|h\|^2_{L^2(\mathbb{R}_+^*)} = 2 \pi \|\varphi\|^2_{H^s_+}, \quad \mathrm{and} \quad \|\partial_{\xi} h\|^2_{L^2(\mathbb{R}_+^*)}\leq 4 \pi \left(\|\mathbf{G} \varphi\|_{H^s_+}^2 + s^2 \|\varphi\|_{H^s_+}^2 \right).
\end{equation}Using the Plancherel identity, formulas \eqref{DefGphiSect2} and  \eqref{fohL2norm}, we obtain that 
\begin{equation*}
2 \pi \operatorname{Re} \langle i \mathbf{G}\varphi , \varphi \rangle_{H^s_+} =  \frac{1}{2} \int_0^{+\infty} |\widehat{\varphi}(\xi)|^2 2s(1+ \xi^2)^{s-1}\xi - \partial_{\xi}(|h|^2)(\xi)\mathrm{d}\xi.
\end{equation*}We conclude by using $\widehat{\varphi}(0^+)=h(0^+)$ and $\lim_{\xi \to +\infty} h(\xi) =0$.
\end{proof}

\noindent For $s \ge 1$, the space $\mathrm{Dom}^s(\mathbf{G})$ admits an alternative description.
\begin{lem}\label{DsphiinDomGLem}
For any $s \geq 1$, we have $\mathrm{Dom}^s(\mathbf{G}) = \{\phi \in H^s_+ \cap  \mathrm{Dom}(\mathbf{G}): D^s \phi  \in \mathrm{Dom}(\mathbf{G}) \}$ and 
\begin{equation}\label{DDsGsubsetDs-1Gfo}
\partial_x \mathrm{Dom}^s(\mathbf{G}) = D \left( \mathrm{Dom}^s(\mathbf{G}) \right) \subset \mathrm{Dom}^{s-1}(\mathbf{G}).
\end{equation}Besides, if $h \in \mathrm{Dom}^1(\mathbf{G})$,  then $[\mathbf{G}, D] h =\mathbf{G}  D h - D\mathbf{G} h = i h \in \mathrm{Dom}^1(\mathbf{G})$.
\end{lem}
\begin{proof}
Fix $s \geq 1$, set $\mathfrak{p}_s(\xi):=\xi^s$, $\forall \xi >0$. For any $\varphi \in H^s_+ \cap  \mathrm{Dom} (\mathbf{G})$, we set $g:=D^s \varphi \in L^2_+$. Then 
\begin{equation}\label{dxihatgfor}
\partial_{\xi}\widehat{g} =\partial_{\xi} \left(\mathfrak{p}_s \widehat{\varphi} \right)  = s \mathfrak{p}_{s-1} \widehat{\varphi} + \mathfrak{p}_s \partial_{\xi} \widehat{\varphi}  \quad \mathrm{on} \quad \mathbb{R}_+^* = (0, +\infty).
\end{equation}The proof of Lemma $\ref{DsphiinDomGLem}$ is divided into three steps.\\

\noindent \textbf{Step \romannumeral 1}. Assume that $g = D^s \varphi \in \mathrm{Dom}(\mathbf{G})$, then $\widehat{g}|_{\mathbb{R}_+^*} \in H^1(\mathbb{R}_+^*)$. Plancherel identity yields that
\begin{equation*}
\|\mathfrak{p}_{s-1} \widehat{\varphi}\|_{L^2(\mathbb{R}_+^*)}^2 \leq \int_0^{+\infty} |\widehat{\varphi}(\xi)|^2 \xi^{2s-2} \mathrm{d}\xi \leq 2 \pi \|\varphi\|_{H^{s-1}_+}^2.
\end{equation*}By  \eqref{dxihatgfor}, we have $\mathfrak{p}_s \partial_{\xi} \widehat{\varphi} = \partial_{\xi}\widehat{g} -s \mathfrak{p}_{s-1} \widehat{\varphi} \in L^2(\mathbb{R}_+^*)$. Since $(1+|\xi|^2)^s \leq 2^{s-1} (1+|\xi|^{2s})$, we have 
\begin{equation*}
2 \pi \|\mathbf{G} \varphi\|_{H^s_+}^2 \leq 2^{s-1} \int_0^{+\infty}(1+|\xi|^{2s})|\partial_{\xi} \widehat{\varphi}(\xi)|^2 \mathrm{d} \xi \leq 2^s \pi \|\mathbf{G}\varphi\|_{L^2}^2 + 2^{s-1} \|\mathfrak{p}_s \partial_{\xi} \widehat{\varphi}\|^2_{L^2(\mathbb{R}_+^*)} <+\infty.
\end{equation*}So $\varphi \in H^s_+ \cap  \mathrm{Dom}(\mathbf{G}) \cap  \mathbf{G}^{-1}(H^s_+) = \mathrm{Dom}^s (\mathbf{G})$. \\

\noindent \textbf{Step \romannumeral 2}. Conversely, assume that $\varphi \in \mathrm{Dom}^s (\mathbf{G})$, then $g:= D^s \varphi \in L^2_+$ and \eqref{dxihatgfor} holds. So we have
\begin{equation}
\|\partial_{\xi} \widehat{g}\|_{L^2(\mathbb{R}_+^*)}^2 \leq 2 \int_0^{+\infty}  |s \xi^{s-1} \widehat{\varphi}(\xi)|^2 +|\xi^s \partial_{\xi} \widehat{\varphi}(\xi)|^2 \mathrm{d}\xi \leq 4s^2 \pi \| \varphi\|_{H^{s-1}_+}^2 + 4 \pi \|\mathbf{G} \varphi\|_{H^s_+}^2 < +\infty.
\end{equation}Then $\widehat{g}|_{\mathbb{R}_+^*} \in H^1(\mathbb{R}_+^*)$ and $D^s \varphi = g \in \mathrm{Dom}(\mathbf{G})$. As a consequence, we have  
\begin{equation}\label{DomsGCharbyDsinDomG}
\mathrm{Dom}^s (\mathbf{G}) = H^s_+ \cap  \mathrm{Dom}(\mathbf{G}) \cap  (D^s)^{-1}\left(\mathrm{Dom}(\mathbf{G}) \right).
\end{equation}\textbf{Step \romannumeral 3}. If $h \in \mathrm{Dom}^1 (\mathbf{G})$, we choose $s=1$ in \textbf{Step \romannumeral 2}, then $D h = -i \partial_x h \in \mathrm{Dom} (\mathbf{G})$. Formula \eqref{dxihatgfor} yields that $\widehat{\mathbf{G}Dh}(\xi) = i \partial_{\xi}(\mathfrak{p}_1 \widehat{h})(\xi) = i \widehat{h}(\xi) + \xi\widehat{\mathbf{G}h} (\xi)$, $\xi >0$. So $\mathbf{G}Dh = ih + D \mathbf{G} h \in L^2_+$. We have
\begin{equation}\label{dxG]commfor}
\partial_x \mathrm{Dom}^1 (\mathbf{G}) \subset \mathrm{Dom}(\mathbf{G})\quad \mathrm{and} \quad \left[\partial_x, \mathbf{G} \right] = \mathrm{id}_{L^2_+} : \mathrm{Dom}^1 (\mathbf{G}) \to \mathrm{Dom}^1 (\mathbf{G}).
\end{equation}Fix $s>1$. If $\varphi \in \mathrm{Dom}^s (\mathbf{G})$, then we have $f:= \partial_x \varphi \in \mathrm{Dom}(\mathbf{G}) \cap  H^{s-1}_+$ by using  \eqref{DomrGsubsetDomsG} and \eqref{dxG]commfor}. \textbf{Step \romannumeral 2} yields that $D^{s-1}f = D^s \varphi \in \mathrm{Dom}(\mathbf{G})$. Using   \eqref{DomsGCharbyDsinDomG}, we obtain that $\partial_x \varphi = f \in   \mathrm{Dom}^{s-1} (\mathbf{G})$. 
\end{proof} 

\begin{rem}\label{RkDGI=0DX}
Assume that $\phi \in \mathrm{Dom}(\mathbf{G}) \cap  H^{\alpha}_+$ for some $\alpha >0$, then $\widehat{D^{\alpha} \phi}(\xi) = \xi^{\alpha} \widehat{\phi}(\xi ) \to 0 \widehat{\phi}(0^+) = 0$, as $\xi \to 0^+$. For any $s \geq 1$, Lemma $\ref{DsphiinDomGLem}$ yields that $\partial_x \mathrm{Dom}^s(\mathbf{G})   \subset \mathrm{Dom}^{s-1}(\mathbf{G}) \cap  \mathrm{Dom}(\mathbf{X}) $ and
\begin{equation}\label{DomsGcharDsphiinDomX}
\mathrm{Dom}^s(\mathbf{G}) = \{\phi \in H^s_+ \cap  \mathrm{Dom}(\mathbf{G}): D^s \phi  \in \mathrm{Dom}(\mathbf{X}) \} .
\end{equation} 
\end{rem}

\noindent For any $f \in L^2_+$, its Poisson integral can be represented in terms of the generator $(-i\mathbf{G}, \mathrm{Dom}(\mathbf{G}))$ and the operator $(\mathbf{I}_+, \mathrm{Dom}(\mathbf{I}_+))$ via the following inversion formula of Newton--Leibniz type:
\begin{equation}\label{InvForNewtonLeibniz}
\mathbf{P}_z (f) = \underline{f}(z) = \tfrac{1}{2\pi i}  \mathbf{I}_+ \bigl( (\mathbf{G}-z)^{-1} f \bigr), \qquad \forall z \in \mathbb{C}_+.
\end{equation}Several proofs of the identity  \eqref{InvForNewtonLeibniz} exist in the literature.
\begin{itemize}
    \item It first appeared in the work of G\'erard~\cite[ p.~599]{BOExp2022}, proved via a formal Laplace transform argument.
    \item A proof using the identification $\mathbf{X} = \mathbf{G}^*$ was subsequently given by G\'erard--Pushnitski~\cite[formula~(2.4)]{GerPush2023}.
\end{itemize}

\noindent Both of the proofs above rely on the    approximation argument 
 (Lemma $\ref{phichiepCVtoI0Lem}$ and formula \eqref{bfIeitLapisbfI})  employing the auxiliary functions $(\chi_{\epsilon})_{\epsilon >0}$ defined in~\eqref{chiepFor}, followed by passing to the limit $\epsilon \to 0^+$.  
A more direct proof, which bypasses this approximation, is presented in Lemma~2.3 of~\cite{SUNMSzegoR}. 
It uses only the fundamental theorem of calculus, justifying the name \textit{Newton--Leibniz type formula}. 
For the sake of completeness, we recall its proof, following the idea of Lemma~2.3 in~\cite{SUNMSzegoR}.

\begin{proof}[Proof of formula \eqref{InvForNewtonLeibniz}]
Fix $z \in \mathbb{C}_+$ and $f \in L^2_+$, and set $\varphi := (\mathbf{G} - z)^{-1} f \in \mathrm{Dom}(\mathbf{G})$. Then we have
\begin{equation}\label{idxibfezhatphi=ezhatf}
i\partial_{\xi} \bigl( \mathbf{e}_z \widehat{\varphi} \bigr) = \mathbf{e}_z \bigl( i \partial_{\xi} \widehat{\varphi} - z \widehat{\varphi} \bigr) = \mathbf{e}_z \widehat{f} \quad \text{on} \quad \mathbb{R}_+^* = (0, +\infty), \quad \text{with} \quad \mathbf{e}_z (\xi) = e^{i z \xi} \in \mathbb{C},
\end{equation}
by \eqref{DefGphiSect2} and \eqref{LaxBeurlingSG}. Since $\operatorname{Im} z > 0$ and $\widehat{\varphi}|_{\mathbb{R}_+^*} \in H^1(\mathbb{R}_+^*)$ by Definition $\ref{DombfGdefSect2}$, we have $\lim_{\xi \to +\infty} \widehat{\varphi}(\xi) = \lim_{\xi \to +\infty} \mathbf{e}_z(\xi) = 0$. Integrating \eqref{idxibfezhatphi=ezhatf} with respect to $\xi$ on $\mathbb{R}_+^*$, we obtain $2 \pi i  \underline{f}(z) = i  \int_0^{+\infty} \mathbf{e}_z \widehat{f}  =   \widehat{\varphi}(0^+) =   \mathbf{I}_+ \bigl( ( \mathbf{G} - z )^{-1} f \bigr)$, 
by formulas \eqref{PoiIntHolFTRL2+} and \eqref{DefI+}.
\end{proof}

\noindent The following lemma identifies the difference quotient operator $\mathscr{Q}_z$ given by \eqref{diffQuoDefForQzS2} with the resolvent of the generator $\mathbf{G}$ when acting on the Hardy space $L^2_+$, for any $z \in \mathbb{C}_+$.

\begin{lem}\label{scQzinL2+LemG-z}
For any $z \in \mathbb{C}_+$ and $f \in L^2_+$, we have $\mathscr{Q}_z(f) = (\mathbf{G} - z)^{-1} f \in \mathrm{Dom}(\mathbf{G})$.  
Consequently, for every $s \ge 0$, the restriction $\mathscr{Q}_z|_{H^s_+}$ coincides with the resolvent $(\mathbf{G} - z)^{-1}: H^s_+ \to \mathrm{Dom}^s(\mathbf{G})$, which is therefore a bijection. Moreover,   the following  uniform estimate holds:
\begin{equation}\label{EstiuniQzHs+}
\sup_{s \geq 0}\|\mathscr{Q}_z\|_{\mathcal{B}_{\mathbb{C}}(H^s_+)} \le (\operatorname{Im} z)^{-1}, \qquad \forall z \in \mathbb{C}_+. 
\end{equation}
\end{lem}

\begin{proof}
Fix $z \in \mathbb{C}_+$ and $f \in L^2_+$. Set $\phi := (\mathbf{G} - z)^{-1} f \in \mathrm{Dom}(\mathbf{G})$. For almost every $x \in \mathbb{R}$, we have
\begin{equation}\label{fx-funz=x-zphix}
f(x) = (\mathbf{G}\phi)(x) - z\phi(x) = x\phi(x) - \tfrac{i}{2\pi}\widehat{\phi}(0^+) - z\phi(x) = (x - z)\phi(x) + \underline{f}(z),
\end{equation}
where we used Lemma~\ref{LemGphi=Xphi-const} and formula~\eqref{InvForNewtonLeibniz}.
By \eqref{bdvarphizHinf-}, \eqref{diffQuoDefForQzS2} and \eqref{fx-funz=x-zphix}, it follows that
\begin{equation}\label{qzfisvarphiinDomG}
\mathscr{Q}_z(f) = (\mathbf{G} - z)^{-1} f \in \mathrm{Dom}(\mathbf{G}), \quad \forall f \in L^2_+, \quad \forall z \in \mathbb{C}_+.
\end{equation}
The estimate~\eqref{EstiuniQzHs+} then follows directly from \eqref{qzfisvarphiinDomG} and \eqref{EstG-z-1BHs}.
\end{proof}

\noindent The next lemma shows that the restriction of the difference quotient operator $\mathscr{Q}_z$ to $L^2_+$ coincides with the Toeplitz operator whose symbol is the Cauchy kernel $\boldsymbol{\varphi}_z \in H^{\infty}_-$ given by \eqref{bdvarphizHinf-}, for any $z \in \mathbb{C}_+$. Consequently, such Toeplitz operators always map into the domain of $\mathbf{G}$.   
   
\begin{lem}
For any $z \in \mathbb{C}_+$, let $\boldsymbol{\varphi}_z(x) = (x - z)^{-1}$ be given by \eqref{bdvarphizHinf-}. Then we have
\begin{equation}\label{TbdvarphiisQzfor}
\mathbf{T}_{\boldsymbol{\varphi}_z}(f) = \mathscr{Q}_z(f) = (\mathbf{G} - z)^{-1} f \in \mathrm{Dom}(\mathbf{G}), \qquad \forall f \in L^2_+.
\end{equation}
\end{lem}
\begin{proof}
Fix $z \in \mathbb{C}_+$ and set $\zeta := \overline{z} \in \mathbb{C}_-$. Define $\boldsymbol{\varphi}_{\zeta}(x) := (x - \zeta)^{-1} = \overline{\boldsymbol{\varphi}_z(x)}$ for all $x \in \mathbb{R}$, so that $\boldsymbol{\varphi}_{\zeta} \in H^{\infty}_+$. 
 The Hankel operator of symbol  $\boldsymbol{\varphi}_{\zeta}  $ is denoted by $\mathbf{H}_{\boldsymbol{\varphi}_{\zeta}} : f \in L^2_+ \mapsto \mathbf{H}_{\boldsymbol{\varphi}_{\zeta}} (f) : = \Pi_{\geq 0} (\boldsymbol{\varphi}_{\zeta} \overline{f})  \in L^2_+$. The   Kronecker theorem (Theorem 2.1 of   Pocovnicu \cite{pocoAPDE2011}) shows that
\begin{equation}\label{KronThmHankOprk}
\operatorname{Ran} (\mathbf{H}_{\boldsymbol{\varphi}_{\zeta}}) = \mathbb{C}\boldsymbol{\varphi}_{\zeta}.
\end{equation}Fix $f \in L^2_+$. By \eqref{KronThmHankOprk}, there exists a unique constant  $\mathtt{c} = \mathtt{c}(f ,\zeta) \in \mathbb{C}$ such that  
\begin{equation}\label{Piphizefbariscphi}
\Pi_{\geq 0} (\boldsymbol{\varphi}_{\zeta} \overline{f})  = \mathbf{H}_{\boldsymbol{\varphi}_{\zeta}} (f)  = \mathtt{c} \boldsymbol{\varphi}_{\zeta} \in  H^{\infty} _+. 
\end{equation}Since $f \boldsymbol{\varphi}_{\zeta} \in L^2_+ \cap L^1 (\mathbb{R} ; \mathbb{C})$ by Lemma $\ref{Z1+idealinL2+}$,   Lemma   $\ref{LemL2+capL10mode}$  implies that $\int_{\mathbb{R}} f \boldsymbol{\varphi}_{\zeta}  = 0$. Moreover, we have
\begin{equation}\label{Poiisintfvarphiz}
|\boldsymbol{\varphi}_{\zeta}|^2 = \tfrac{i}{2 \operatorname{Im} z} (\boldsymbol{\varphi}_{\zeta} - \overline{\boldsymbol{\varphi}_{\zeta}}) \in H^{\infty} (\mathbb{R} ; \mathbb{R}) 
\quad  \text{and} \quad 
2 \pi i \underline{f}(z) = \int_{\mathbb{R}} f  \boldsymbol{\varphi}_z  = \int_{\mathbb{R}}f \overline{\boldsymbol{\varphi}_{\zeta}} \in \mathbb{C}  
,
\end{equation}by \eqref{PoiIntHolFTRL2+}. 
  Using formulas \eqref{Piphizefbariscphi},  \eqref{Poiisintfvarphiz} and \eqref{bdvarphnormLq}, we obtain that
\begin{equation}\label{ttcbarisfz}
\mathtt{c}\| \boldsymbol{\varphi}_{\zeta} \|_{L^2_+}^2 = \int_{\mathbb{R}} \overline{f} |\boldsymbol{\varphi}_{\zeta}|^2 = \frac{i}{2 \operatorname{Im} z} \left(  \int_{\mathbb{R}} \overline{f}  \boldsymbol{\varphi}_{\zeta} - \int_{\mathbb{R}} \overline{f \boldsymbol{\varphi}_{\zeta}} \right) = \frac{\pi}{\operatorname{Im} z} \overline{\underline{f}(z)}
=  \overline{\underline{f}(z)} \|  \boldsymbol{\varphi}_z \|_{L^2_-}^2 
.
\end{equation}So  $\overline{\mathtt{c}} = \underline{f}(z)$. By formulas \eqref{fpi+fPI-fsum} and   \eqref{Piphizefbariscphi}, it follows that
\begin{equation}\label{Pinegphizffzpz}
\Pi_{ < 0} (\boldsymbol{\varphi}_z f ) = \overline{\Pi_{\geq 0} (\overline{ f \boldsymbol{\varphi}_z  })}
= \overline{\Pi_{\geq 0} (\boldsymbol{\varphi}_{\zeta} \overline{f})} = \overline{\mathtt{c} \boldsymbol{\varphi}_{\zeta}} = \underline{f}(z) \boldsymbol{\varphi}_z \in H^{\infty}_-.
\end{equation}Then $
 \mathbf{T}_{\boldsymbol{\varphi}_z}(f)  = \boldsymbol{\varphi}_z f - \Pi_{<0} (\boldsymbol{\varphi}_z f) =
(f - \underline{f}(z))\boldsymbol{\varphi}_z = \mathscr{Q}_z (f) =  (\mathbf{G} - z)^{-1}(f) \in \mathrm{Dom}(\mathbf{G})$ follows directly from formulas \eqref{ToepdefForSec2},  \eqref{fpi+fPI-fsum}, \eqref{Pinegphizffzpz} and Lemma $\ref{scQzinL2+LemG-z}$.
\end{proof}
  
\noindent A remarkable feature of the difference quotient operator $\mathscr{Q}_z$ is that, although its multiplier $\boldsymbol{\varphi}_z$ given by \eqref{bdvarphizHinf-} belongs to the \emph{negative} frequency space $H^\infty_-$, the operator itself preserves \emph{positive} frequency regularity when acting on $L^{\infty}_+$ or on the Hardy--Zhidkov spaces $\mathcal{Z}^N_+$ ($N \in \mathbb{N}_+$).
 
\begin{lem}\label{QzLinftyinHardyLem}
For any $z \in \mathbb{C}_+$ and $N \in \mathbb{N}_+$, we have $\mathscr{Q}_z(L^{\infty}_+) \subset L^2_+$ and $\mathscr{Q}_z(\mathcal{Z}^N_+) \subset H^N_+$.
\end{lem}
 
\begin{proof}
Fix $z \in \mathbb{C}_+$ and $u \in L^{\infty}_+$. For every $\epsilon>0$, we have $u \chi_{\epsilon} \in L^2_+$ by Lemma \ref{Z1+idealinL2+} and    \eqref{chiepFor}. Lemma $\ref{scQzinL2+LemG-z}$ implies that $\mathscr{Q}_z (u\chi_{\epsilon}) = (\mathbf{G} - z )^{-1} (u\chi_{\epsilon})  \in \mathrm{Dom}(\mathbf{G}) $. Since   $L^2_+$ is a closed subspace of $L^2(\mathbb{R} ; \mathbb{C})$, we have $\mathscr{Q}_z (u) = \lim_{\epsilon \to 0^+}\mathscr{Q}_z (u \chi_{\epsilon}) \in L^2_+$ by  
Lemma $\ref{QzuchieptoQzuL2lem}$. So  $\mathscr{Q}_z(L^{\infty}_+) \subset L^2_+$. Fix $N \in  \mathbb{N}_+$. Using   \eqref{HarZhidSpZnDefFor}, \eqref{HardySobHspmS2} and  Lemma  $\ref{LemHNCVQzuchiep}$, we obtain that
 $\mathscr{Q}_z (\mathcal{Z}^N_+) \subset  \mathscr{Q}_z (L^{\infty}_+) \cap \mathscr{Q}_z (\mathcal{Z}^N )    \subset L^2_+ \cap  H^N (\mathbb{R} ; \mathbb{C}) = H^N_+$.
\end{proof}

\noindent Inspired by Lax \cite{LaxTransinva1959} and Lemma 3.3  of \cite{SUNCMP2021},
 we introduce a distinguished family of finite-dimensional translation-invariant subspaces in the Hardy space $L^2_+$.
 
\begin{prop}\label{KthetapropOrDecom}
If $Q \in \mathbb{C}[X]$ is a monic polynomial of degree $N$ such that $Q^{-1}\{0\} \subset \mathbb{C}_-$, we define an inner function $\theta_Q:= \frac{\overline{Q}}{Q} \in 1+ H^{\infty}_+$. Then we have the following orthogonal decomposition of $L^2_+$,
\begin{equation}\label{OrthDecomInnKthe}
L^2_+ = \theta_Q L^2_+ \oplus \tfrac{\mathbb{C}_{\leq N-1}[X]}{Q} \quad \mathrm{and} \quad \theta_Q L^2_+ \perp \tfrac{\mathbb{C}_{\leq N-1}[X]}{Q}  \quad \mathrm{in} \quad L^2_+.
\end{equation}Moreover, the subspace $\theta_Q L^2_+$ is invariant under the action of the shift semigroup $(\mathbf{S}(\eta))_{\eta \geq 0}$, while its orthogonal complement $\frac{\mathbb{C}_{\leq N-1}[X]}{Q}$ is invariant under the adjoint semigroup $(\mathbf{S}(\eta)^*)_{\eta \geq 0}$. That is,
\begin{equation}\label{ShiftInvProperty}
\mathbf{S}(\eta) (\theta_Q L^2_+) \subset \theta_Q L^2_+, \quad \mathbf{S}(\eta)^* \bigl(\tfrac{\mathbb{C}_{\leq N-1}[X]}{Q} \bigr) \subset \tfrac{\mathbb{C}_{\leq N-1}[X]}{Q} \subset \mathrm{Dom}(\mathbf{G}) \cap  H^{\infty}_+, \quad \forall \eta \geq 0.
\end{equation}As a consequence, we have 
\begin{equation}\label{KthetainDominftyG}
\mathbf{G}\bigl(\tfrac{\mathbb{C}_{\leq N-1}[X]}{Q} \bigr)  = \tfrac{\mathbb{C}_{\leq N-1}[X]}{Q}
\quad \text{and} \quad 
\theta_Q - 1 \in \tfrac{\mathbb{C}_{\leq N-1}[X]}{Q} \subset \mathrm{Dom}^{\infty}(\mathbf{G}).
\end{equation} 
\end{prop}
\begin{proof}
According to Proposition 1.2 of \cite{SUNCMP2021}, there is a bijection between the set of monic polynomials $Q \in \mathbb{C}[X]$ of degree $N$ with roots in $\mathbb{C}_-$ and the $N$-soliton manifold of the Benjamin--Ono equation. Then the orthogonal decomposition formula  \eqref{OrthDecomInnKthe} follows from Proposition 4.3 of \cite{SUNCMP2021}. The invariant subspace property \eqref{ShiftInvProperty} is a direct consequence of formula  \eqref{SemiGrCVdefDomsG} and the following facts
\begin{equation}
\dim_{\mathbb{C}} \bigl(\tfrac{\mathbb{C}_{\leq N-1}[X]}{Q} \bigr)   = N \quad \mathrm{and} \quad  \mathbf{S}(\eta) \theta_Q f = \theta_Q \mathbf{e}_{\eta} f = \theta_Q\mathbf{S}(\eta) f \in \theta_Q L^2_+, \quad \forall f \in L^2_+, \quad \forall \eta \geq 0.
\end{equation}Moreover, Proposition $\ref{GisInfiniGeneHsProp}$ and  formula  \eqref{ShiftInvProperty} imply that for any $\varphi \in \tfrac{\mathbb{C}_{\leq N-1}[X]}{Q}$, we have
\begin{equation}\label{GinvKtheta}
-i \mathbf{G}\varphi = \lim_{\eta \to 0^+} \tfrac{1}{\eta}(\mathbf{S}(\eta)^* \varphi - \varphi) \in \tfrac{\mathbb{C}_{\leq N-1}[X]}{Q} \quad \Longrightarrow \quad  \mathbf{G}\bigl(\tfrac{\mathbb{C}_{\leq N-1}[X]}{Q} \bigr)  \subset \tfrac{\mathbb{C}_{\leq N-1}[X]}{Q}.
\end{equation}The bijectivity of $\mathbf{G}: \tfrac{\mathbb{C}_{\leq N-1}[X]}{Q} \to \tfrac{\mathbb{C}_{\leq N-1}[X]}{Q}$ is established by Lemma $\ref{LemSpGC_-}$. Using formulas \eqref{DomsGDefSect2} and \eqref{DomrGsubsetDomsG}, we obtain that $\tfrac{\mathbb{C}_{\leq N-1}[X]}{Q} \subset \mathrm{Dom}(\mathbf{G}) \cap  H^{\infty}_+ \cap  \mathbf{G}^{-1} (H^{\infty}_+) = \mathrm{Dom}^{\infty}(\mathbf{G})$.
\end{proof}

\begin{rem}\label{RemLaxTrInvSp}
By Lax \cite{LaxTransinva1959} and Lemma 3.3  of \cite{SUNCMP2021}, if $\{0\} \ne \mathsf{M} \subset \mathrm{Dom}(\mathbf{G})$ is a   finite-dimensional subspace such that $\mathbf{G} (\mathsf{M}) \subset \mathsf{M}$, then $\mathsf{M}$ must be of the form $\mathsf{M} = \frac{\mathbb{C}_{\leq N-1}[X]}{Q}$, for some monic polynomial $Q \in \mathbb{C}[X]$ of degree $N = \dim_{\mathbb{C}} \mathsf{M} \in \mathbb{N}_+$ with zeros inside the lower half-plane $\mathbb{C}_-$. Specifically, we have
\begin{equation}\label{QCharPolyGresttoM}
Q (x) = \det \bigl( x - \mathbf{G}|_{\mathsf{M}} \bigr), \qquad x \in \mathbb{R}.
\end{equation}  
\end{rem}
  
\begin{lem}\label{LemqtoC[X]ovqbij}
Assume that $Q, q \in \mathbb{C}[X]$ are two polynomials of degree $N\geq 1$ such that 
\begin{equation}\label{Condforq<->Q}
Q^{-1}\{0\}  \cup q^{-1}\{0\} \subset \mathbb{C}_- \quad \mathrm{and} \quad \tfrac{\mathbb{C}_{\leq N-1}[X]}{Q} = \tfrac{\mathbb{C}_{\leq N-1}[X]}{q},
\end{equation}then $Q$ and $q$ are associates in $\mathbb{C}[X]$, i.e. $Q=cq$ with $c=\frac{Q(0)}{q(0)} \in \mathbb{C}^*=\mathbb{C}\setminus \{0\}$.
\end{lem}
\begin{proof}
Since $\frac{1}{Q} \in \tfrac{\mathbb{C}_{\leq N-1}[X]}{Q} = \tfrac{\mathbb{C}_{\leq N-1}[X]}{q}$ by formula \eqref{Condforq<->Q}, there exists $\mathfrak{b} \in \mathbb{C}[X]$ such that $\deg \mathfrak{b} \leq N-1$ and $q = Q \mathfrak{b} \in Q \mathbb{C}[X]$. By symmetry, swapping $Q$ and $q$ in the above argument, we obtain  $q \mid Q$.
\end{proof}

\noindent Fix $s \geq 0$. We now complete our discussion of the shift semigroups by identifying the precise subspace of $H^s_+$ on which the right shift semigroup $(\mathbf{S}(\eta))_{\eta \geq 0}$ given by \eqref{LaxBeurlingSG} is strongly differentiable. This subspace, denoted by $\mathrm{Dom}^s(\mathbf{X})$, is the natural higher-order analogue of $\mathrm{Dom}(\mathbf{X})$ defined in \eqref{DomX}. \\ 

\noindent Recall that $H^1_s(\mathbb{R}) = W^{1,2}_s(\mathbb{R})$ denotes the global weighted Sobolev space in \eqref{WeightedSobolevSp}. Set
\begin{equation}\label{WeiGlobSoboH1s(R)}
\|f\|_{H^1_s(\mathbb{R})}^2:= \|f\|^2_{L^2_s(\mathbb{R})}+ \|\partial_{\xi} f\|^2_{L^2_s(\mathbb{R})} =\int_{\mathbb{R}} \left(|f(\xi)|^2+ |\partial_{\xi} f(\xi)|^2 \right) (1+\xi^2)^s \mathrm{d}\xi \in [0, +\infty],
\end{equation}for any measurable function $f: \mathbb{R} \to \mathbb{C}$.  Then $f\in H^1_s(\mathbb{R})$ if and only if $\|f\|_{H^1_s(\mathbb{R})} < +\infty$. We define \begin{equation}\label{DomsXSect2defFor}
\mathrm{Dom}^s(\mathbf{X}):=    \{\varphi \in L^2_+: \widehat{\varphi} \in H^1_s(\mathbb{R})\}  
\end{equation}to be the space of all functions in $L^2_+$ whose Fourier transforms belong to   $H^1_s(\mathbb{R})$. If $r >s \geq 0$, then
\begin{equation}\label{DomrX<DomsXfor}
\mathrm{Dom}^r (\mathbf{X}) \subset \mathrm{Dom}^s (\mathbf{X}) \subset \mathrm{Dom}^0 (\mathbf{X}) =\mathrm{Dom} (\mathbf{X}).
\end{equation}The next lemma shows the relation between $\mathrm{Dom}^s(\mathbf{X})$ and $\mathrm{Dom}^s(\mathbf{G})$, which is defined by \eqref{DomsGDefSect2}. 

\begin{lem}
For any $s \geq 0$, the space $\mathrm{Dom}^s(\mathbf{X})$ admits the following equivalent characterizations:
\begin{equation}\label{DomsXCharbyDomsG}
\small
\mathrm{Dom}^s(\mathbf{X})    = \mathrm{Dom}^{s}(\mathbf{G}) \cap  \mathrm{Dom}(\mathbf{X}) =\{\varphi \in \mathrm{Dom}^s (\mathbf{G}): \widehat{\varphi}(0^+)=0\}   =  H^s_+ \cap  \mathrm{Dom}(\mathbf{X}) \cap  \mathbf{X}^{-1}\left( H^s_+ \right).
\end{equation}
\end{lem}
\begin{proof}
The case $s=0$ is given by \eqref{DomX} and \eqref{DXsubsetDG}. Fix $s>0$. The proof is divided into $3$ steps. \\

\noindent \textbf{Step \romannumeral 1}. If $f \in H^1_s(\mathbb{R})$, then \eqref{WeightedSobolevSp} yields that $f|_{\mathbb{R}_+^*} \in H^1_s(\mathbb{R}_+^*)$. By \eqref{DomsGDefSect2}, \eqref{DomsXSect2defFor}, \eqref{DomrX<DomsXfor}, we have
\begin{equation*}
\small
\mathrm{Dom}^s(\mathbf{X})    \subset \mathrm{Dom}^{s}(\mathbf{G}) \cap  \mathrm{Dom}(\mathbf{X}) \subset  \{\varphi \in \mathrm{Dom}^s (\mathbf{G}): \widehat{\varphi}(0^+)=0\} .
\end{equation*}\textbf{Step \romannumeral 2}. Assume that $\varphi \in \mathrm{Dom}^s (\mathbf{G})  \subset \mathrm{Dom}  (\mathbf{G}) \cap  H^s_+$ such that $\widehat{\varphi}(0^+)=0$,  then we have $\varphi \in \mathrm{Dom}(\mathbf{X})$ by  \eqref{DomX}. Lemma $\ref{LemGphi=Xphi-const}$ yields that $\mathbf{X} \varphi = \mathbf{G}\varphi \in H^s_+$. So $\varphi \in  H^s_+ \cap  \mathrm{Dom}(\mathbf{X}) \cap  \mathbf{X}^{-1}\left( H^s_+ \right)$. \\

\noindent \textbf{Step \romannumeral 3}. For any $\phi \in  H^s_+ \cap  \mathrm{Dom}(\mathbf{X}) \cap  \mathbf{X}^{-1}\left( H^s_+ \right)$, we have $\mathbf{X}\phi   \in H^s_+$. Since $\widehat{\mathbf{X}\phi} = i \partial_{\xi} \widehat{\phi} \in H^s_+$, then
\begin{equation*}
\small
\|\widehat{\phi}\|_{H^1_s(\mathbb{R})} = \int_{\mathbb{R}}\left(|\widehat{\phi}(\xi)|^2+|\widehat{\mathbf{X}\phi}(\xi)|^2 \right) (1+\xi^2)^s \mathrm{d}\xi = 2 \pi \left(\|\phi\|_{H^s_+}^2 + \|\mathbf{X}\phi\|_{H^s_+}^2 \right) <+\infty,
\end{equation*}by \eqref{WeiGlobSoboH1s(R)} and Plancherel identity. So $\phi \in \mathrm{Dom}^s(\mathbf{X})$ and formula \eqref{DomsXCharbyDomsG} is obtained.
\end{proof}

\noindent Then Lemma $\ref{X=InfiniGeL2SetaLem}$ can be generalized to the following proposition.

\begin{prop}\label{XDomsXinfiGeneforSeProp}
For any $s \geq 0$,  the operator $\left(i\mathbf{X}, \mathrm{Dom}^s(\mathbf{X})\right)$, defined in \eqref{defXf-} and \eqref{DomsXSect2defFor}, is the infinitesimal generator of the right shift semigroup $(\mathbf{S}(\eta))_{\eta \geq 0}$ on the Hardy--Sobolev space $H^s_+$, i.e.
\begin{equation}
\label{DomsXSemiGrCVHs}
 \mathrm{Dom}^s(\mathbf{X}) =\{\phi \in H^s_+: \exists v\in H^s_+ \quad \mathrm{such} \quad \mathrm{that} \quad \lim_{\eta \to 0^+}\|\tfrac{1}{\eta}(\mathbf{S}(\eta)  \phi - \phi)-v\|_{H^s_+}=0\}.
\end{equation}Moreover, if $\phi \in \mathrm{Dom}^s(\mathbf{X})$, then $\tfrac{1}{\eta}(\mathbf{S}(\eta)^* \phi - \phi) \to  i \mathbf{X} \phi$ in $H^s_+$, as $\eta \to 0^+$.
\end{prop}
\begin{proof}
Fix $s >0$. If $\varphi \in \mathrm{Dom}^s(\mathbf{X})$, then \eqref{DomsXSect2defFor} yields that $\widehat{\varphi} \in H^1_s(\mathbb{R}) = W^{1,2}_s(\mathbb{R})$. We have $\widehat{X \varphi} = i \partial_{\xi} \widehat{\varphi}$ and $\widehat{\mathbf{S}(\eta) \varphi} = \tau_{\eta} \widehat{\varphi}$, $\forall \eta >0$. Using Plancherel identity and Lemma $\ref{CharSobDiffQuoCV}$, we obtain that
\begin{equation*}
\begin{split}
2 \pi \|\tfrac{1}{\eta}\left(\mathbf{S}(\eta)\varphi - \varphi \right) - i\mathbf{X} \varphi \|_{H^s_+}^2 = & \int_{\mathbb{R}}|\tfrac{1}{\eta}(\widehat{\varphi}(\xi - \eta)  - \widehat{\varphi}(\xi) ) + \partial_{\xi} \widehat{\varphi}(\xi)|^2 (1+ \xi^2)^s \mathrm{d}\xi \\
= & \|\tfrac{1}{ \eta}(\tau_{ \eta } \widehat{\varphi} - \widehat{\varphi}) + \partial_{\xi} \widehat{\varphi} \|_{L^2_s(\mathbb{R})}^2 \to 0, \quad \mathrm{as} \quad \eta \to 0^+.
\end{split}
\end{equation*}Conversely, assume that $\phi, v \in H^s_+$ satisfy $\lim_{\eta \to 0^+}\|\tfrac{1}{\eta}(\mathbf{S}(\eta)  \phi - \phi)-v\|_{H^s_+}=0$. Then  
\begin{equation*}
\|\tfrac{1}{\eta}(\mathbf{S}(\eta)  \phi - \phi)-v\|_{L^2_+}   \leq   \|\tfrac{1}{\eta}(\mathbf{S}(\eta)  \phi - \phi)-v\|_{H^s_+} \to 0, 
\end{equation*}as $\eta \to 0^+$.  Lemma $\ref{X=InfiniGeL2SetaLem}$ implies that $\phi \in \mathrm{Dom}(\mathbf{X})$. Since $ \tfrac{1}{\eta}(\mathbf{S}(\eta)  \phi - \phi)$ converges to $i \mathbf{X}\phi$  in $L^2_+$, as $\eta \to 0^+$, we have $\mathbf{X}\phi = -i v \in  H^s_+$. So $\phi \in H^s_+ \cap  \mathrm{Dom}(\mathbf{X}) \cap  \mathbf{X}^{-1}(H^s_+) = \mathrm{Dom}^s(\mathbf{X})$ by \eqref{DomsXCharbyDomsG}.
\end{proof}

\begin{lem}
For any $s \geq 1$, we have $\partial_x \mathrm{Dom}^s (\mathbf{X}) \subset \partial_x \mathrm{Dom}^s (\mathbf{G}) \subset   \mathrm{Dom}^{s-1} (\mathbf{X}) $. Moreover,
\begin{equation}\label{forDomsXCharbyDsinDX}
\mathrm{Dom}^s(\mathbf{X}) = \{\phi \in  H^s_+ \cap  \mathrm{Dom}(\mathbf{X}) :  D^s \phi \in \mathrm{Dom}(\mathbf{X}) \}.
\end{equation}
\end{lem}
\begin{proof}
Fix $s \geq 1$. According to Remark $\ref{RkDGI=0DX}$ and formula \eqref{DomsXCharbyDomsG}, we have  
\begin{equation*}
\partial_x \mathrm{Dom}^s (\mathbf{X}) \subset \partial_x \mathrm{Dom}^s (\mathbf{G}) \subset  \mathrm{Dom}^{s-1} (\mathbf{G}) \cap  \mathrm{Dom}(\mathbf{X}) =  \mathrm{Dom}^{s-1} (\mathbf{X}).
\end{equation*}Then formula \eqref{forDomsXCharbyDsinDX} can be deduced directly by formulas \eqref{DomsGcharDsphiinDomX} and \eqref{DomsXCharbyDomsG}.
\end{proof}

\noindent We close this subsection with an algebraic identity that
reformulates the nonlinearity of~\eqref{CMdNLSfocu} in the language of
Toeplitz operators, when the potential $u$ is a constant plus a Hardy--Sobolev
perturbation.

\begin{lem}
For any $a \in \mathbb{C}$ and $f \in   H^1_+$, we define $u:= a + f \in \mathcal{Z}^1_+$. Let $\mathbf{T}_u \in \mathcal{B}_{\mathbb{C}}(L^2_+)$ denote the Toeplitz operator given by \eqref{ToepdefForSec2}. Then the following identity holds:  
\begin{equation}\label{uPidxumodsqToep}
u \Pi_{\geq 0} \partial_x \bigl( |u|^2 \bigr) = \mathbf{T}_u \partial_x \mathbf{T}_{\overline{u}} (f) \in L^2_+ .
\end{equation}
\end{lem}

\begin{proof}
Since $a \in \mathbb{C}$,  $f \in H^1_+$ and $u =a +f \in \mathcal{Z}^1_+$, we have $\partial_x \overline{u} = \partial_x \overline{f} \in L^2_-$. Then it follows that 
\begin{equation}\label{gPidxumodsqfor}
g : = \Pi_{\geq 0} \partial_x \bigl( |u|^2 \bigr) = a \Pi_{\geq 0} \partial_x \overline{u} + \partial_x \Pi_{\geq 0} (\overline{u} f) = \partial_x \mathbf{T}_{\overline{u}} (f ) \in L^2_+, \quad \mathbf{T}_u (g) = u g \in L^2_+,  
\end{equation}by \eqref{ZnHardBanAlgebra} and  \eqref{ToepdefForSec2}. 
Then  formula \eqref{uPidxumodsqToep} 
follows directly from \eqref{gPidxumodsqfor}. 
\end{proof}

  \subsection{Unitary Propagator Lemma}
 
\noindent The Hardy--Sobolev spaces $\{H^s_+\}_{s\ge 0}$ built in the previous subsections are the natural arena for the operators appearing in the Lax pair formulation of chiral solutions to the focusing CMdNLS equation 
\eqref{CMdNLSfocu}.  
To propagate such operators along the nonlinear flow, we need a unitary propagator that, besides preserving the $L^2_+$ inner product, also leaves every $H^s_+$ invariant up to a prescribed order.  \\

\noindent Before introducing the unitary propagator lemma, we recall the following classical result on the existence and uniqueness of solutions to linear ordinary differential equations in a Banach space.

\begin{lem}[Cauchy]\label{CauchyThmODE}
Let $\mathcal{E}$ be an arbitrary Banach space over $\mathbb{C}$, $\mathcal{I}$ an open interval of $\mathbb{R}$, and $A \in \mathcal{C}(\mathcal{I}; \mathcal{B}(\mathcal{E}))$. If $(t_0, x_0) \in \mathcal{I} \times \mathcal{E}$, then there exists a unique function $x \in \mathcal{C}^1(\mathcal{I}; \mathcal{E})$ satisfying 
\begin{equation}\label{ODEgeneral}
x(t_0) = x_0, \quad \partial_t x(t) = A(t) \bigl( x(t) \bigr), \quad \forall t \in \mathcal{I}.
\end{equation} 
\end{lem}
\begin{proof}
See Theorem 1.1.1 of Chemin  \cite{CheminNote2016}.
\end{proof}
 
\begin{defi}
For any $s \geq 0$, let $\mathcal{E}^s:=\mathcal{B}_{\mathbb{C}}(H^s_+)$ denote the Banach algebra  of all bounded $\mathbb{C}$-linear operators on the Hardy--Sobolev space $H^s_+$.
\end{defi}
 
\noindent The following lemma provides the unitary propagator: starting from a   continuous time-dependent family of bounded skew-adjoint operators on $L^2_+$, it produces a norm-differentiable family of unitary operators that act simultaneously on all $H^p_+$ with $p \in \{0,1,2, \cdots, m\}$.

\begin{lem}\label{UniPropaLemW}
Assume that $\mathbb{B}\in \bigcap_{p=0}^m \mathcal{C} \left( (T^-, T^+); \mathcal{B}_{\mathbb{C}}(H^p_+) \right)$ for some   $m \in \mathbb{N}_+$ and $T^- <0 < T^+$. If $i\mathbb{B}(t) $ is a bounded self-adjoint   operator on $L^2_+$, $\forall t \in (T^-, T^+)$, then there exists a unique mapping
\begin{equation}\label{Wdeffor}
\mathbb{W} =  \mathbb{W}_{\mathbb{B}}  \in \bigcap_{p=0}^m \mathcal{C}^1 \left( (T^-, T^+); \mathcal{B}_{\mathbb{C}}(H^p_+) \right)
\end{equation}satisfying the following:\\
\noindent $(\mathbf{a}).$ For any $f\in L^2_+$, we have $\mathbb{W}(0)f=f$.\\
\noindent $(\mathbf{b}).$  For any $t \in (T^-, T^+)$, let $\mathbb{W}(t)^*\in \mathcal{B}_{\mathbb{C}}(L^2_+)$ denote the $L^2_+$-adjoint of $\mathbb{W}(t) \in \mathcal{B}_{\mathbb{C}}(L^2_+)$, then 
\begin{equation}\label{ForuniofW}
\mathbb{W}(t)^* \mathbb{W}(t) = \mathbb{W}(t) \mathbb{W}(t)^* = \mathrm{id}_{L^2_+}, \quad T^- < t < T^+.
\end{equation}$(\mathbf{c}).$ Set $\mathbb{W}^* :   t \in (T^-, T^+) \mapsto \mathbb{W}^*(t):= \mathbb{W}(t)^* \in \mathcal{B}_{\mathbb{C}}(L^2_+)$. Then $\mathbb{W}^* \in \bigcap_{p=0}^m \mathcal{C}^1 \left( (T^-, T^+); \mathcal{B}_{\mathbb{C}}(H^p_+) \right)$.\\
\noindent $(\mathbf{d}).$ For any $t \in (T^-, T^+)$,   the following limits hold:
\begin{equation}\label{Wt+eta-Wt-BWCVfor}
\begin{cases}
\lim_{r \to 0} \sum_{p=0}^m\|r^{-1}(\mathbb{W}(t+r)-\mathbb{W}(t)) - \mathbb{B}(t) \mathbb{W}(t)\|_{\mathcal{B}_{\mathbb{C}}(H^p_+)} = 0,\\
\lim_{r \to 0} \sum_{p=0}^m\|r^{-1}(\mathbb{W}(t+r)^{-1}-\mathbb{W}(t)^{-1}) + \mathbb{W}(t)^{-1}\mathbb{B}(t) \|_{\mathcal{B}_{\mathbb{C}}(H^p_+)} = 0.
\end{cases}
\end{equation}
\end{lem}

\begin{rem}
Assertions $(\mathbf{a})$ and $(\mathbf{d})$ in Lemma $\ref{UniPropaLemW}$ show that  $\mathbb{W}, \mathbb{W}^* \in \bigcap_{p=0}^m \mathcal{C}^1 \left( (T^-, T^+); \mathcal{B}_{\mathbb{C}}(H^p_+) \right)$  satisfy that $\mathbb{W}(0)=\mathrm{id}_{L^2_+}=\mathbb{W}(0)^*$ and for any $ t \in (T^-, T^+)$, we have 
\begin{equation}\label{ODEforW}
\partial_t \mathbb{W}(t) = \mathbb{B}(t)\mathbb{W}(t) \in \mathcal{B}_{\mathbb{C}}(H^p_+), \qquad \partial_t (\mathbb{W}^*)(t)  = -\mathbb{W}(t)^{-1}\mathbb{B}(t)  \in \mathcal{B}_{\mathbb{C}}(H^p_+),
\end{equation}for any $ p \in \{0,1,2, \cdots, m\}$. Formula \eqref{ODEforW} is to be understood in the sense of   \eqref{Wt+eta-Wt-BWCVfor}.
\end{rem}
\begin{proof}[Proof of Lemma $\ref{UniPropaLemW}$]
Fix $p \in [0, m] \cap  \mathbb{Z}$. For any $t \in (T^-, T^+)$ and $\mathscr{W} \in \mathcal{E}^p = \mathcal{B}_{\mathbb{C}}(H^p_+)$, we set
\begin{equation}\label{Ajdeffor}
\mathcal{A}_1(t)(\mathscr{W}):= \mathbb{B}(t) \mathscr{W} \in \mathcal{E}^p, \qquad \mathcal{A}_2(t)(\mathscr{W}):= -\mathscr{W}\mathbb{B}(t)  \in \mathcal{E}^p.
\end{equation}Then the mapping $\mathcal{A}_j : t \in (T^-, T^+) \mapsto\mathcal{A}_j(t) \in \mathcal{E}^p$ is continuous with respect to the uniform operator topology. Lemma $\ref{CauchyThmODE}$ implies that there exists a unique solution to the following ODE
\begin{equation}\label{ODEWp}
\mathscr{W}_p(0)= \mathrm{id}_{H^p_+}, \quad \partial_t \mathscr{W}_p(t) = \mathcal{A}_1(t)\left(\mathscr{W}_p(t)\right), \quad T^- < t <T^+.
\end{equation}\textbf{Step \romannumeral 1}. We claim that $\mathscr{W}_0(t)|_{H^p_+}=\mathscr{W}_p(t)$,  for any $p \in [1,m]\bigcap\mathbb{Z}$ and  $t \in (T^-, T^+)$. 
\begin{proof}[Proof of the claim]Fix $p\in [1,m]\bigcap\mathbb{Z}$ and $\phi \in H^p_+$. Set $\phi_p(t):= \mathscr{W}_p(t) \phi \in H^p_+$ and $\phi_0(t):=\mathscr{W}_0(t) \phi \in L^2_+$, for any $t \in (T^-, T^+)$. Then $\phi_p \in  \mathcal{C}^1\left( (T^-, T^+);  H^p_+ \right)$ and $\phi_0 \in  \mathcal{C}^1\left( (T^-, T^+);  L^2_+ \right)$ both solve the ODE 
\begin{equation}\label{Odevarphifor}
\boldsymbol{\varphi} (0)=\phi, \quad \partial_t \boldsymbol{\varphi} (t) = \mathbb{B}(t) \left(\boldsymbol{\varphi} (t)\right)\in L^2_+,  \quad T^- < t <T^+.
\end{equation}By Lemma $\ref{CauchyThmODE}$, equation \eqref{Odevarphifor} admits only one solution. So $\mathscr{W}_p(t) \phi=\mathscr{W}_0(t) \phi \in H^p_+$. 
\end{proof}
\noindent Set $\mathbb{W}(t):= \mathscr{W}_0(t) \in \bigcap_{p=0}^m \mathcal{B}_{\mathbb{C}}(H^p_+)$,   $\forall t \in (T^-, T^+)$. Then formulas \eqref{Ajdeffor} and \eqref{ODEWp} yield that 
\begin{equation}\label{DerForW'-BW}
 \lim_{r \to 0}\sum_{p=0}^m\|r^{-1}(\mathbb{W}(t+r)-\mathbb{W}(t)) - \mathbb{B}(t) \mathbb{W}(t)\|_{\mathcal{B}_{\mathbb{C}}(H^p_+)} = 0.
\end{equation}\textbf{Step \romannumeral 2}. We aim to show that $\mathbb{W}(t) \in \mathcal{E}^0 = \mathcal{B}_{\mathbb{C}}(L^2_+)$ is $L^2_+$-unitary, $\forall t \in (T^-, T^+)$.
\begin{proof}[Proof of formula \eqref{ForuniofW}]
By using the self-adjointness of $i\mathbb{B}(t) \in \mathcal{E}^0$ and formulas \eqref{Ajdeffor}, \eqref{ODEWp}, we obtain that $\partial_t(\mathbb{W}^*\mathbb{W})(t) = \mathbb{W}(t)^*\left( \mathbb{B}(t)^*+ \mathbb{B}(t) \right)\mathbb{W} (t)=0 \in \mathcal{E}^0$ and the $\mathcal{C}^1$-mapping  $\mathbb{W}\mathbb{W}^* : t \in (T^-, T^+) \mapsto \mathbb{W}(t)\mathbb{W}(t)^* \in \mathcal{E}^0$ solves the following ODE: for any $t \in (T^-, T^+)$,
\begin{equation}
\mathscr{U}(0)=\mathrm{id}_{L^2_+}, \quad \partial_t \mathscr{U} (t) = [\mathbb{B}(t), \mathscr{U}(t)] =\left( \mathcal{A}_1(t) + \mathcal{A}_2(t)\right) (\mathscr{U} (t) ) \in \mathcal{E}^0.
\end{equation}Combining the continuity of $\mathcal{A}_1  + \mathcal{A}_2 : t \in (T^-, T^+) \mapsto   \mathcal{A}_1(t) + \mathcal{A}_2(t)  \in \mathcal{E}^0$ and the uniqueness result in Lemma $\ref{CauchyThmODE}$, we obtain that $\mathbb{W}(t)\mathbb{W}(t)^* = \mathrm{id}_{L^2_+}$, $\forall t \in (T^-, T^+)$. So $\mathbb{W}(t)^*=\mathbb{W}(t)^{-1} \in \mathcal{B}_{\mathbb{C}}(L^2_+)$.
\end{proof}
\noindent \textbf{Step \romannumeral 3}. Show that   $\mathbb{W}(t): H^p_+ \to H^p_+$ is invertible, $\forall p \in [1,m] \cap  \mathbb{Z}$. It suffices to prove that
 \begin{equation}\label{surjForWtHp}
\mathbb{W}(t) (H^p_+) =  H^p_+ = \mathbb{W}(t)^{-1} (H^p_+), \quad \forall t \in \mathcal{I}:=(T^-, T^+).
\end{equation}
\begin{proof}[Proof of formula \eqref{surjForWtHp}]
The proof uses a uniqueness argument via flow reversal, a standard technique in ODE theory. Fix $\mathbf{t}_1 \in \mathcal{I}= (T^-, T^+)$. Set $\mathcal{I}_1:=(\mathbf{t}_1-T^+, \mathbf{t}_1-T^-)$ and $\mathbb{B}_1(\mathfrak{t}):=\mathbb{B}(\mathbf{t}_1 - \mathfrak{t}) \in \mathcal{E}^0$, $\forall \mathfrak{t} \in \mathcal{I}_1$. Since $\mathbb{B}_1\in \bigcap_{p=0}^m \mathcal{C}\left( \mathcal{I}_1; \mathcal{B}_{\mathbb{C}}(H^p_+) \right)$, we apply the previous  \textbf{Step \romannumeral 1, Step \romannumeral 2} and formula \eqref{DerForW'-BW} to deduce that there exists a unique function $\mathscr{M} \in \bigcap_{p=0}^m \mathcal{C}^1\left( \mathcal{I}_1; \mathcal{B}_{\mathbb{C}}(H^p_+) \right)$ such that $\mathscr{M}(0)=\mathrm{id}_{L^2_+}$ and
\begin{equation}\label{DerforM'+B1M}
\lim_{r \to 0}\sum_{p=0}^m\|r^{-1}(\mathscr{M}(\mathfrak{t}+r)-\mathscr{M}(\mathfrak{t})) + \mathbb{B}_1(\mathfrak{t})\mathscr{M}(\mathfrak{t})\|_{\mathcal{B}_{\mathbb{C}}(H^p_+)} = 0, \quad \forall t \in (T^-, T^+).
\end{equation}Moreover $\mathscr{M}(t) \in \mathcal{B}_{\mathbb{C}}(L^2_+)$ is $L^2_+$-unitary. Fix $f \in L^2_+$. Set $\mathbf{w}(t):= \mathbb{W}(t) f \in L^2_+$, $\forall t \in (T^-, T^+)$. Then $\mathbf{w} \in   \mathcal{C}^1 \left( (T^-, T^+); L^2_+ \right)$ solves equation \eqref{Odevarphifor} with initial datum $\mathbf{w}(0)=f$. We define that
\begin{equation}\label{ydefForm}
\mathbf{y}(\mathfrak{t})=\mathbf{w}(\mathbf{t}_1 - \mathfrak{t}) \in L^2_+, \quad \forall \mathfrak{t} \in \mathcal{I}_1:=(\mathbf{t}_1-T^+, \mathbf{t}_1-T^-).
\end{equation}Then $\mathbf{y}\in \mathcal{C}^1 \left( \mathcal{I}_1; L^2_+ \right)$ satisfies that $\mathbf{y}(0)=\mathbf{w}(\mathbf{t}_1)=\mathbb{W}(\mathbf{t}_1)f\in L^2_+$ and it solves the following ODE
\begin{equation}\label{Odefory'=-B1y}
\partial_t \mathbf{y}(\mathfrak{t}) = - \mathbb{B}_1 (\mathfrak{t}) \left(\mathbf{y}(\mathfrak{t})\right) \in L^2_+, \quad \forall \mathfrak{t} \in \mathcal{I}_1.
\end{equation}For any $\mathfrak{t} \in \mathcal{I}_1$, we set $\boldsymbol{\xi}(\mathfrak{t}):= \mathscr{M}(\mathfrak{t})(\mathbf{y}(0)) \in L^2_+$. Then  \eqref{DerforM'+B1M} yields that the function $\boldsymbol{\xi} \in \mathcal{C}^1 \left( \mathcal{I}_1; L^2_+ \right)$ also solves equation \eqref{Odefory'=-B1y} with initial datum $\boldsymbol{\xi} (0) = \mathbf{y}(0)=\mathbb{W}(\mathbf{t}_1)f \in L^2_+$. By using the uniqueness result in Lemma $\ref{CauchyThmODE}$, we have $\boldsymbol{\xi} = \mathbf{y}\in \mathcal{C}^1 \left( \mathcal{I}_1; L^2_+ \right)$ and $f=\mathbf{y}(\mathbf{t}_1)=\boldsymbol{\xi}(\mathbf{t}_1)=\mathscr{M}(\mathbf{t}_1)\left(\mathbf{W}(\mathbf{t}_1)f  \right)\in L^2_+$. Since $f \in L^2_+$ can be arbitrarily chosen, $\mathbb{W}(\mathbf{t}_1)^{-1} = \mathscr{M}(\mathbf{t}_1) \in \bigcap_{p=0}^m \mathcal{B}_{\mathbb{C}}(H^p_+)$.  As a consequence, the operator $\mathbb{W}(\mathbf{t}_1) : H^p_+ \to H^p_+$ is bijective, $\forall p \in [0,m]\cap  \mathbb{Z}$, $\forall \mathbf{t}_1 \in \mathcal{I}=(T^-, T^+)$. 
\end{proof}
\noindent Fix $p \in [0,m]\cap  \mathbb{Z}$.
  Let $\mathscr{G}(\mathcal{E}^p)=\mathscr{G}\left(\mathcal{B}_{\mathbb{C}}(H^p_+)\right)$ denote the set of all bounded invertible $\mathbb{C}$-linear operators on $H^p_+$. Then $\mathscr{G}(\mathcal{E}^p)$ is an open subset of the Banach algebra $\mathcal{E}^p = \mathcal{B}_{\mathbb{C}}(H^p_+)$ and the mapping 
\begin{equation}\label{InvBanaAlg}
\mathbf{Inv}: \mathscr{W} \in \left(\mathscr{G}(\mathcal{E}^p), \|.\|_{\mathcal{E}^p}  \right)\mapsto  \mathbf{Inv}(\mathscr{W}):=\mathscr{W}^{-1} \in \left(\mathscr{G}(\mathcal{E}^p), \|.\|_{\mathcal{E}^p}  \right)
\end{equation}is a homeomorphism. Fix $t \in (T^-, T^+)$. Then there exists $\delta=\delta_{t,p} \in (0, \min(t-T^-, T^+-t))$ such that
\begin{equation}\label{W-1pertnorm<2}
\sup_{|r| \leq \delta}\|\mathbb{W}(t +r)^{-1}\|_{\mathcal{E}^p} \leq 2\|\mathbb{W}(t)^{-1}\|_{\mathcal{E}^p}.
\end{equation}Then formulas \eqref{Ajdeffor}, \eqref{ODEWp} and \eqref{W-1pertnorm<2} imply that for any $r \in (-\delta, \delta)$, we have
\begin{equation*}
\small
\begin{split}
& \|r^{-1}(\mathbb{W}(t+r)^{-1}-\mathbb{W}(t)^{-1}) + \mathbb{W}(t)^{-1}\mathbb{B}(t) \|_{\mathcal{B}_{\mathbb{C}}(H^p_+)} \\
\leq &\|\mathbb{W}(t +r)^{-1}\|_{\mathcal{E}^p}  \|r^{-1}\left(\mathbb{W}(t +r) -\mathbb{W}(t)\right) - \mathbb{B}(t) \mathbb{W}(t)\|_{\mathcal{E}^p} \|\mathbb{W}(t )^{-1}\|_{\mathcal{E}^p} \\
& \qquad \qquad + \|\mathbb{W}(t )^{-1} - \mathbb{W}(t + r )^{-1}\|_{\mathcal{E}^p} \|\mathbb{B}(t ) \|_{\mathcal{E}^p} \\
\to & 0,
\end{split}
\end{equation*}as $r \to 0$. Finally, assertions $(\mathbf{c})$ and $(\mathbf{d})$ of Lemma $\ref{UniPropaLemW}$ are obtained and the proof is completed.
\end{proof}

\section{The Lax Pair Structure}\label{SectLaxPair3}

\noindent This section is devoted to the Lax pair formulation $\partial_t \mathbf{L}_u = [\mathbf{B}_u, \mathbf{L}_u]$ of the focusing CMdNLS equation \eqref{CMdNLSfocu}. 
We show that the sum $-i (\mathbf{G} + 2 t \mathbf{L}_u)$, where $\mathbf{G}$ is given by \eqref{DefGphiSect2} and $\mathbf{L}_u$ denotes the chiral Lax operator, can be interpreted as a perturbation of the infinitesimal generator of the conjugate shift semigroup $(e^{-it\partial_x^2}\mathbf{S}(\eta)^* e^{it\partial_x^2})_{\eta \geq 0}$ on $L^2_+$. 
We then establish several fundamental commutation relations of the semigroup $(\mathbf{S}(\eta)^*)_{\eta \geq 0}$ defined in \eqref{LaxBeurlingSG} with Toeplitz operators, the Lax operator $\mathbf{L}_u$, and the operator $\widetilde{\mathbf{B}}_u := \mathbf{B}_u - i \mathbf{L}_u^2$. 
These identities and estimates are essential for the explicit formula, the invariant subspaces, and the action-angle mapping developed in later sections. 
Finally, we carry out a spectral analysis that determines the essential spectrum and proves the simplicity of  eigenvalues of $\mathbf{L}_u$.

 \subsection{The Non-Chiral Lax Pair Formulation}

\noindent The focusing CMdNLS equation \eqref{CMdNLSfocu} admits a Lax pair formulation whose operators can be studied either on the full Sobolev space $H^1(\mathbb{R};\mathbb{C})$ or on its Hardy subspace $H^1_+ = \Pi_{\geq 0}(H^1(\mathbb{R};\mathbb{C}))$. 
In this subsection we consider the first setting, which we refer to as the \emph{non-chiral} formulation. \\ 

\noindent We begin by introducing the multiplication operator. For any $u \in L^2(\mathbb{R};\mathbb{C}) + L^{\infty}(\mathbb{R};\mathbb{C})$, we define
\begin{equation}\label{MultiOpbfmDef}
\mathbf{m}_u : h \in H^1(\mathbb{R};\mathbb{C}) \mapsto \mathbf{m}_u(h) := uh \in L^2(\mathbb{R};\mathbb{C}).
\end{equation}
When $\mathfrak{a} \in L^{\infty}(\mathbb{R};\mathbb{C})$, the operator $\mathbf{m}_{\mathfrak{a}}$ extends to a bounded $\mathbb{C}$-linear operator on $L^2(\mathbb{R};\mathbb{C})$. Its $L^2(\mathbb{R};\mathbb{C})$-adjoint is given by $(\mathbf{m}_{\mathfrak{a}})^* = \mathbf{m}_{\overline{\mathfrak{a}}} \in \mathcal{B}_{\mathbb{C}}(L^2(\mathbb{R};\mathbb{C}))$. Specifically, for any $f, g \in L^2(\mathbb{R};\mathbb{C})$, we have
\begin{equation}\label{AdjL2ofbfma}
\langle \mathbf{m}_{\mathfrak{a}}(f), g \rangle_{L^2(\mathbb{R};\mathbb{C})} = \int_{\mathbb{R}} \mathfrak{a} f \overline{g} \, \mathrm{d}x = \langle f, \mathbf{m}_{\overline{\mathfrak{a}}}(g) \rangle_{L^2(\mathbb{R};\mathbb{C})} \quad \text{and} \quad \|\mathbf{m}_{\mathfrak{a}}\|_{\mathcal{B}_{\mathbb{C}}(L^2(\mathbb{R};\mathbb{C}))} \leq \|\mathfrak{a}\|_{L^{\infty}(\mathbb{R};\mathbb{C})}.
\end{equation}
The filtered multiplication operator is denoted by $\mathscr{T}_u := \Pi_{\geq 0} \circ \mathbf{m}_u$, where $\Pi_{\geq 0} \in \mathcal{B}_{\mathbb{C}}(L^2(\mathbb{R};\mathbb{C}))$ is the Szeg\H{o} projector given by \eqref{SzegoProjS2}. The Toeplitz operator $\mathbf{T}_u$ given by \eqref{ToepdefForSec2} is the restriction of $\mathscr{T}_u$ to the Hardy--Sobolev space $H^1_+$. In other words, for any $u \in L^2(\mathbb{R};\mathbb{C}) + L^{\infty}(\mathbb{R};\mathbb{C})$, we have
\begin{equation}\label{FilMulOpscrTu}
\mathbf{T}_u = \mathscr{T}_u|_{H^1_+} : H^1_+ \to L^2_+, \quad \text{where} \quad \mathscr{T}_u(h) = \Pi_{\geq 0}(uh) \in L^2_+, \quad \forall h \in H^1(\mathbb{R};\mathbb{C}).
\end{equation}
For any $\mathfrak{a} \in L^{\infty} (\mathbb{R} ; \mathbb{C})$ and $\varphi \in L^2 (\mathbb{R} ; \mathbb{C})$, it follows directly from \eqref{MultiOpbfmDef} and \eqref{FilMulOpscrTu}  that
\begin{equation}\label{forTamphi=Tphima}
\mathscr{T}_{\mathfrak{a}} \mathbf{m}_{\varphi} = \mathscr{T}_{\mathfrak{a} \varphi} = \mathscr{T}_{\varphi} \mathbf{m}_{\mathfrak{a}} : H^1 (\mathbb{R} ; \mathbb{C}) \to L^2_+ \subset L^2 (\mathbb{R} ; \mathbb{C}).
\end{equation}

\begin{lem}
For any $u, v, f, g \in \mathcal{Z}^1$, the following identity holds:
\begin{equation}\label{bfmscrTTtomT+mTFor}
\begin{split}
\mathbf{m}_f \mathscr{T}_{ \partial_x (\overline{u} v)} \mathscr{T}_{\overline{g}}
= & \mathbf{m}_f \bigl( \mathscr{T}_{\overline{u}} \mathbf{m}_{\partial_x v} + \mathscr{T}_{ \partial_x \overline{u}} \mathbf{m}_v  \bigr) \mathscr{T}_{\overline{g}} \\
 = & \mathbf{m}_{f \Pi_{\geq 0} \partial_x (\overline{u} v)} \mathscr{T}_{\overline{g}} + \mathbf{m}_f \mathscr{T}_{\overline{g \Pi_{\geq 0} \partial_x (\overline{v}u)}}
: H^1 (\mathbb{R} ; \mathbb{C}) \to L^2 (\mathbb{R} ; \mathbb{C}).
\end{split}
\end{equation}
\end{lem} 
\begin{proof}
Fix $f, g, u, v \in \mathcal{Z}^1$. Set $\mathtt{h}: = \partial_x (\overline{u} v) = v \partial_x \overline{u} + \overline{u} \partial_x v \in L^2 (\mathbb{R} ; \mathbb{C})$. By \eqref{fpi+fPI-fsum} and \eqref{Z1Banachalgebra}, 
we have
\begin{equation}\label{tthpmdeffor}
\mathtt{h}_+ : = \Pi_{\geq 0} \mathtt{h} = \Pi_{\geq 0}   \partial_x (\overline{u} v) \in L^2_+ \quad \text{and} \quad \mathtt{h}_- : =\Pi_{< 0} \mathtt{h} = \Pi_{<0  }   \partial_x (\overline{u} v) = \overline{\Pi_{\geq 0}   \partial_x (\overline{v} u ) } \in  L^2_-.
\end{equation}For any $q \in H^1 (\mathbb{R} ; \mathbb{C})$, set $\varphi : = \mathscr{T}_g (q) = \Pi_{\geq 0} (gq) \in H^1_+$ by \eqref{Z1H1subsetH1}.  Using  \eqref{tthpmdeffor} and  Lemma $\ref{Z1+idealinL2+}$, we obtain 
\begin{equation}\label{bfmscrTh+sctTgq}
\bigl( \mathbf{m}_f \mathscr{T}_{\mathtt{h}_+} \mathscr{T}_g  \bigr) (q) =  f \Pi_{\geq 0} (\mathtt{h}_+ \varphi) =  \varphi  f \Pi_{\geq 0}   \partial_x (\overline{u} v)  = \mathbf{m}_{f \Pi_{\geq 0}   \partial_x (\overline{u} v)} \mathscr{T}_g (q) \in L^2 (\mathbb{R} ; \mathbb{C}) .
\end{equation}Since $\Pi_{< 0 } (\overline{g} q) \in H^1_-$, Corollary $\ref{CorLinfty-L2-}$ implies that $\mathtt{h}_- \Pi_{< 0 } (\overline{g} q) \in L^2_-$. By \eqref{fpi+fPI-fsum} and  \eqref{tthpmdeffor}, we have 
\begin{equation}\label{scrThnegTovgqto1T}
\mathscr{T}_{\mathtt{h}_-} \mathscr{T}_{\overline{g}} (q) = \Pi_{\geq 0} \bigl( 
\mathtt{h}_- \overline{g} q  - \mathtt{h}_- \Pi_{< 0} (\overline{g} q) \bigr) = \Pi_{\geq 0} \bigl( 
\mathtt{h}_- \overline{g} q \bigr)  = 
\mathscr{T}_{\overline{g \Pi_{\geq 0} \partial_x (\overline{v}u)}}  (q) \in L^2_+ , 
\end{equation}for any   $q \in H^1 (\mathbb{R} ; \mathbb{C})$. 
Using formulas \eqref{FilMulOpscrTu}, \eqref{forTamphi=Tphima}  and  \eqref{tthpmdeffor}, we obtain that
\begin{equation}\label{TubarbfmdxvsymFor}
\mathscr{T}_{\overline{u}} \mathbf{m}_{\partial_x v} + \mathscr{T}_{ \partial_x \overline{u}} \mathbf{m}_v  = 
\mathscr{T}_{ \partial_x ( \overline{u} v)} =  \mathscr{T}_{ \mathtt{h}} = \mathscr{T}_{ \mathtt{h}_+} + \mathscr{T}_{ \mathtt{h}_-} 
: H^1 (\mathbb{R} ; \mathbb{C}) \to L^2 (\mathbb{R} ; \mathbb{C}).
\end{equation}Then formula  \eqref{bfmscrTTtomT+mTFor} 
follows from \eqref{tthpmdeffor},  \eqref{bfmscrTh+sctTgq},  
\eqref{scrThnegTovgqto1T}  and \eqref{TubarbfmdxvsymFor}. 
\end{proof}
 
\noindent According to Sun \cite{SUNLMP2024}, the focusing CMdNLS equation \eqref{CMdNLSfocu} can be regarded as the diagonal reduction $u=v$ of the following intertwined derivative nonlinear Schr\"odinger system (idNLS):
\begin{equation}\label{idNLS}
\begin{cases}
\partial_t u = i \partial_x^2 u + u \Pi_{\geq 0} \partial_x (\overline{v}u) + v \Pi_{\geq 0} \partial_x (|u|^2),\\[4pt]
\partial_t v = i \partial_x^2 v + v \Pi_{\geq 0} \partial_x (\overline{u}v) + u \Pi_{\geq 0} \partial_x (|v|^2),
\end{cases}\qquad (t,x) \in \mathbb{R} \times \mathbb{R}.
\end{equation}
We now introduce the non-chiral Lax pair structure of  system  \eqref{idNLS}. For every $n \in \mathbb{N}_+$, let $\mathcal{Z}^n$ be the Zhidkov space given by \eqref{ZhiSpdefSect2}.
 For any $\psi = (u,v) \in \mathcal{Z}^1 \times \mathcal{Z}^1$, we define two auxiliary operators:
\begin{equation}\label{scrApsiOmegaDef}
\begin{cases}
\mathscr{A}_{\psi} : = \mathbf{m}_u \mathscr{T}_{\overline{v}}  + \mathbf{m}_v    \mathscr{T}_{\overline{u}}  =\mathbf{m}_u \Pi_{\geq 0} \mathbf{m}_{\overline{v}}  + \mathbf{m}_v \Pi_{\geq 0} \mathbf{m}_{\overline{u}} 
 \in \mathcal{B}_{\mathbb{C}}(L^2 (\mathbb{R} ; \mathbb{C})) \cap \mathcal{B}_{\mathbb{C}}(H^1 (\mathbb{R} ; \mathbb{C}) ) ; \\
 \Omega_{\psi}:= \mathbf{m}_u \mathscr{T}_{\partial_x\overline{v}} 
+ \mathbf{m}_v \mathscr{T}_{\partial_x\overline{u}} - \mathbf{m}_{\partial_x u} \mathscr{T}_{ \overline{v}} 
 - \mathbf{m}_{\partial_x v} \mathscr{T}_{ \overline{u}} \; : \;  H^1 (\mathbb{R} ; \mathbb{C}) \to L^2 (\mathbb{R} ; \mathbb{C}).
 \end{cases}
\end{equation}Then $\mathscr{A}_{\psi}$ is a bounded self-adjoint operator on $L^2(\mathbb{R} ; \mathbb{C})$ by \eqref{AdjL2ofbfma}.
 Recall that $D = -i\partial_x  $ is a  self-adjoint operator on $L^2(\mathbb{R}; \mathbb{C})$. 
With these preparations, the Lax pair $(\mathscr{L}_{\psi}, \mathscr{B}_{\psi})$ for   \eqref{idNLS} is defined by
\begin{equation}\label{scrLscrBidNLS}
\mathscr{L}_{\psi} : = D - \tfrac{1}{2} \mathscr{A}_{\psi} : H^1 (\mathbb{R} ; \mathbb{C}) \to L^2 (\mathbb{R} ; \mathbb{C}) 
\quad \text{and}   \quad 
\mathscr{B}_{\psi} : = \tfrac{1}{2} \Omega_{\psi} +  \tfrac{i}{4} \mathscr{A}_{\psi} ^2 :  H^1 (\mathbb{R} ; \mathbb{C}) \to L^2 (\mathbb{R} ; \mathbb{C}).
\end{equation}The next lemma is a  direct consequence   of  the Kato--Rellich theorem and   
 von Neumann's theorem on the self-adjointness of $A^*A$ (Theorem \uppercase\expandafter{\romannumeral 10}. 25  in  Reed--Simon \cite{Reed2SimFouSAop}).
\begin{lem}\label{LemSelAdjscrLpsi}
For any  $\psi = (u, v) \in \mathcal{Z}^1 \times \mathcal{Z}^1$, the operators $\mathscr{L}_{\psi}$ and $\mathscr{L}_{\psi}^2$ are unbounded self-adjoint operators on $L^2(\mathbb{R} ; \mathbb{C})$ with $\mathrm{Dom}(\mathscr{L}_{\psi}) = H^1(\mathbb{R} ; \mathbb{C})$ and  $\mathrm{Dom}(\mathscr{L}_{\psi}^2) = H^2(\mathbb{R} ; \mathbb{C})$.
\end{lem}
\begin{lem}\label{tilscrBpsiLem}
 For any $\psi = (u, v) \in \mathcal{Z}^1 \times \mathcal{Z}^1$,   set  $\widetilde{\mathscr{B}}_{\psi} :=   \mathscr{B}_{\psi} - i \mathscr{L}_{\psi}^2 $. Then we have
\begin{equation}\label{tildscrBpsiDef}
 \widetilde{\mathscr{B}}_{\psi} = i \partial_x^2 + \mathbf{m}_u \partial_x \mathscr{T}_{\overline{v}} + \mathbf{m}_v \partial_x \mathscr{T}_{\overline{u}} : H^1 (\mathbb{R} ; \mathbb{C}) \to L^2(\mathbb{R} ; \mathbb{C}).
\end{equation}Moreover, the operator $i  \widetilde{\mathscr{B}}_{\psi}$ is a self-adjoint operator on $L^2(\mathbb{R} ; \mathbb{C})$ with $\mathrm{Dom} (i  \widetilde{\mathscr{B}}_{\psi}) = H^2 (\mathbb{R} ; \mathbb{C})$.
\end{lem}
\begin{proof}[Proof of Lemma $\ref{tilscrBpsiLem}$]
Fix $\psi = (u,v) \in \mathcal{Z}^1 \times \mathcal{Z}^1$.
Using formula   \eqref{scrLscrBidNLS}, we obtain that
\begin{equation}\label{scrLpsisqFor}
 i \mathscr{L}_{\psi}^2 = - i \partial_x^2 - \tfrac{1}{2} \mathscr{A}_{\psi} \partial_x - \tfrac{1}{2} \partial_x \mathscr{A}_{\psi} + \tfrac{i}{4} \mathscr{A}_{\psi}^2   : H^2 (\mathbb{R}; \mathbb{C}) \to L^2 (\mathbb{R}; \mathbb{C}).
\end{equation}Recall that $\Pi_{\geq 0 } \partial_x = \partial_x \Pi_{\geq 0 } : H^1 (\mathbb{R}; \mathbb{C}) \to L^2 (\mathbb{R}; \mathbb{C})$. From \eqref{MultiOpbfmDef}, \eqref{AdjL2ofbfma} and \eqref{FilMulOpscrTu}, 
it follows that
\begin{equation}\label{imudxTvbarissym}
\langle \mathbf{m}_u \partial_x \mathscr{T}_{\overline{v}} (h_1), h_2 \rangle_{ L^2( \mathbb{R} ; \mathbb{C})} = - \langle h_1 , \mathbf{m}_v \partial_x \mathscr{T}_{\overline{u}} (h_2) \rangle_{ L^2( \mathbb{R} ; \mathbb{C})} , \quad \forall  h_1, h_2 \in H^1 (\mathbb{R} ; \mathbb{C}).
\end{equation}
Since  $\widetilde{\mathscr{B}}_{\psi}  =  \mathscr{B}_{\psi} - i \mathscr{L}_{\psi}^2 $, formula  \eqref{tildscrBpsiDef} follows directly from \eqref{scrApsiOmegaDef},  \eqref{scrLpsisqFor} and the following identity:
\begin{equation}\label{mudxscrTfor}
\mathbf{m}_u \partial_x \mathscr{T}_{\overline{v}}  = \mathbf{m}_u \mathscr{T}_{\partial_x \overline{v}} + \mathbf{m}_u \mathscr{T}_{\overline{v}} \partial_x = \partial_x \mathbf{m}_u \mathscr{T}_{\overline{v}} - \mathbf{m}_{\partial_x u} \mathscr{T}_{\overline{v}} : H^1 (\mathbb{R}; \mathbb{C}) \to L^2 (\mathbb{R}; \mathbb{C}).
\end{equation}By   \eqref{mudxscrTfor}, \eqref{ZhidNormZn}, \eqref{varphiHs+normHsHs-1},  \eqref{InffH1overfLinfty}, \eqref{cstmsmallsq2-1} and the interpolation inequality, we have
\begin{equation}\label{relBoundofmudxTvbar}
  \big\| \mathbf{m}_u \partial_x \mathscr{T}_{\overline{v}} (h) \big\|_{L^2 (\mathbb{R} ; \mathbb{C})} \leq   \| u \|_{L^{\infty}} \| v \|_{\mathcal{Z}^1} \|   h \|_{H^1 } 
\leq  \| u \|_{\mathcal{Z}^1} \| v \|_{\mathcal{Z}^1} \sqrt{\|   h \|_{L^2  (\mathbb{R} ; \mathbb{C})} \|   h \|_{H^2  (\mathbb{R} ; \mathbb{C})}} ,
\end{equation}for any  $h \in H^2 (\mathbb{R}; \mathbb{C})$.  According to  formulas \eqref{imudxTvbarissym} and \eqref{relBoundofmudxTvbar},
   $ i \bigl( \mathbf{m}_u \partial_x \mathscr{T}_{\overline{v}} + \mathbf{m}_v \partial_x \mathscr{T}_{\overline{u}}  \bigr)$ is a symmetric  operator on $L^2 (\mathbb{R} ; \mathbb{C})$ with domain $H^1(\mathbb{R} ; \mathbb{C})$ and it is infinitesimally small with respect to the self-adjoint operator $-\partial_x^2$. Therefore,  the self-adjointness of  the operator $i \widetilde{\mathscr{B}}_{\psi}$ on $L^2 (\mathbb{R} ; \mathbb{C})$ with domain $H^2(\mathbb{R} ; \mathbb{C})$ follows from the Kato--Rellich theorem.
\end{proof}

\noindent For any $\psi = (u,v) \in \mathcal{Z}^2 \times \mathcal{Z}^2$, the operators $  \Omega_{\psi}$ in  \eqref{scrApsiOmegaDef} and $ \mathscr{B}_{\psi}$ in \eqref{scrLscrBidNLS}   both extend to skew-adjoint operators on $L^2 (\mathbb{R} ; \mathbb{C})$ by formulas \eqref{AdjL2ofbfma} and \eqref{FilMulOpscrTu}.  The following proposition shows that if $\psi = (u,v)$ solves   the idNLS system \eqref{idNLS}, then the pair $(\mathscr{L}_{\psi}, \mathscr{B}_{\psi})$ satisfies the Heisenberg equation of motion, thereby constituting a Lax pair for \eqref{idNLS}.
\begin{prop}\label{LaxPairidNLSprop}
For any $\psi = (u,v) \in \mathcal{Z}^2 \times \mathcal{Z}^2$, let $\mathtt{p}_{\psi}$ and $\mathtt{q}_{\psi}$ be defined by
\begin{equation}\label{ttppsittqpsiDef}
\begin{cases}
\mathtt{p}_{\psi} := i \partial_x^2 u + u \Pi_{\geq 0} \partial_x (\overline{v}u) + v \Pi_{\geq 0} \partial_x (|u|^2), \\[4pt]
\mathtt{q}_{\psi} := i \partial_x^2 v + v \Pi_{\geq 0} \partial_x (\overline{u}v) + u \Pi_{\geq 0} \partial_x (|v|^2).
\end{cases}
\end{equation}
Then $\mathscr{B}_{\psi} \in \mathcal{B}_{\mathbb{C}}(H^1(\mathbb{R};\mathbb{C})) \cap \mathcal{B}_{\mathbb{C}}(L^2(\mathbb{R};\mathbb{C}))$ and the following commutation identity holds on $H^1(\mathbb{R};\mathbb{C})$:
\begin{equation}\label{CommscrBpsiscrLpsi}
-2 \bigl[ \mathscr{B}_{\psi}, \mathscr{L}_{\psi} \bigr] = \mathbf{m}_{\mathtt{p}_{\psi}} \mathscr{T}_{\overline{v}} + \mathbf{m}_v \mathscr{T}_{ \overline{\mathtt{p}_{\psi}}} + \mathbf{m}_{\mathtt{q}_{\psi}} \mathscr{T}_{\overline{u}} +   \mathbf{m}_u \mathscr{T}_{ \overline{\mathtt{q}_{\psi}}} : H^1 (\mathbb{R} ; \mathbb{C}) \to L^2 (\mathbb{R} ; \mathbb{C}) .
\end{equation}
Consequently, if $\psi \in \mathcal{C}\bigl( (T^-, T^+) ; \mathcal{Z}^3 \times \mathcal{Z}^3 \bigr)$ is a solution to \eqref{idNLS}, then we have
\begin{equation}\label{HeiEqscrLscrBpsi}
\partial_t \mathscr{L}_{\psi} = [\mathscr{B}_{\psi}, \mathscr{L}_{\psi}]  , \qquad \forall t \in (T^-, T^+).
\end{equation}
\end{prop}
\begin{proof}
We shall adopt the idea in  Theorem 1.2 of \cite{SUNLMP2024}. Fix $\psi = (u , v) \in \mathcal{Z}^1 \times \mathcal{Z}^1$. It follows that $\mathbf{m}_u , \mathscr{T}_{u} \in \mathcal{B}_{\mathbb{C}}(H^1(\mathbb{R};\mathbb{C})) \cap \mathcal{B}_{\mathbb{C}}(L^2(\mathbb{R};\mathbb{C}))$,  $[\mathbf{m}_u , D] = i  \mathbf{m}_{\partial_x u}$ and $[\mathscr{T}_u , D ] = i\mathscr{T}_{\partial_x u}$. Moreover, we have
\begin{equation}\label{muscrTvbar,DCommFor}
[\mathbf{m}_u \mathscr{T}_{\overline{v}} , D] = i \bigl( \mathbf{m}_u \mathscr{T}_{\partial_x \overline{v}} + \mathbf{m}_{\partial_x u} \mathscr{T}_{\overline{v}} \bigr)  :   H^1 (\mathbb{R} ; \mathbb{C}) \to L^2 (\mathbb{R} ; \mathbb{C}) .
\end{equation}Formula \eqref{scrApsiOmegaDef} yields  that $-i[\mathscr{A}, \mathrm{D}]= \mathbf{m}_u \mathscr{T}_{\partial_x \overline{v}}  + \mathbf{m}_v \mathscr{T}_{\partial_x \overline{u}}  +
 \mathbf{m}_{\partial_x u} \mathscr{T}_{\overline{v}} + \mathbf{m}_{\partial_x v} \mathscr{T}_{\overline{u}} : H^1 \to L^2$. We have
\begin{equation}\label{CommscrAsq,D]For}
\begin{split}
& i [\mathscr{A}_{\psi}^2 , D] =    i \mathscr{A}_{\psi} [\mathscr{A}_{\psi}  , D] + i [\mathscr{A}_{\psi}, D] \mathscr{A}_{\psi} :  H^1 (\mathbb{R} ; \mathbb{C}) \to L^2 (\mathbb{R} ; \mathbb{C}) \\
= & - \mathscr{A}_{\psi}  \mathbf{m}_u \mathscr{T}_{\partial_x \overline{v}}  - \mathscr{A}_{\psi}  \mathbf{m}_v \mathscr{T}_{\partial_x \overline{u}}  - \mathscr{A}_{\psi} 
 \mathbf{m}_{\partial_x u} \mathscr{T}_{\overline{v}} - \mathscr{A}_{\psi}  \mathbf{m}_{\partial_x v} \mathscr{T}_{\overline{u}} \\
 & \quad \; \; - \mathbf{m}_u \mathscr{T}_{\partial_x \overline{v}} \mathscr{A}_{\psi} - \mathbf{m}_v \mathscr{T}_{\partial_x \overline{u}} \mathscr{A}_{\psi} -
 \mathbf{m}_{\partial_x u} \mathscr{T}_{\overline{v}} \mathscr{A}_{\psi}
 - \mathbf{m}_{\partial_x v} \mathscr{T}_{\overline{u}}\mathscr{A}_{\psi} .
\end{split} 
\end{equation}by \eqref{muscrTvbar,DCommFor}. 
Recall that $\Omega_{\psi} : H^1 \to L^2 $ is defined in \eqref{scrApsiOmegaDef}. Then $[\Omega_{\psi}  , \mathscr{A}_{\psi}]  : H^1 \to L^2 $ is given by     
\begin{equation}\label{CommOmescrApsifor}
\begin{split}
   [\Omega_{\psi}  , \mathscr{A}_{\psi}]   
=  & - \mathscr{A}_{\psi}  \mathbf{m}_u \mathscr{T}_{\partial_x \overline{v}}  - \mathscr{A}_{\psi}  \mathbf{m}_v \mathscr{T}_{\partial_x \overline{u}}  + \mathscr{A}_{\psi} 
 \mathbf{m}_{\partial_x u} \mathscr{T}_{\overline{v}} + \mathscr{A}_{\psi}  \mathbf{m}_{\partial_x v} \mathscr{T}_{\overline{u}} \\
 & \quad \; \; + \mathbf{m}_u \mathscr{T}_{\partial_x \overline{v}} \mathscr{A}_{\psi} + \mathbf{m}_v \mathscr{T}_{\partial_x \overline{u}} \mathscr{A}_{\psi} -
 \mathbf{m}_{\partial_x u} \mathscr{T}_{\overline{v}} \mathscr{A}_{\psi}
 - \mathbf{m}_{\partial_x v} \mathscr{T}_{\overline{u}}\mathscr{A}_{\psi} .
\end{split} 
\end{equation}Plugging \eqref{bfmscrTTtomT+mTFor} into formulas 
\eqref{CommscrAsq,D]For}
and 
\eqref{CommOmescrApsifor}, we obtain that
\begin{equation}\label{sumiAsqD-comomeA}
\begin{split}
&  i [\mathscr{A}_{\psi}^2 , D]  -  [\Omega_{\psi}  , \mathscr{A}_{\psi}] : H^1 (\mathbb{R} ; \mathbb{C}) \to L^2 (\mathbb{R} ; \mathbb{C})  \\ 
  =&  -2 \bigl( \mathbf{m}_u \mathscr{T}_{\partial_x (\overline{v} u)}  \mathscr{T}_{\overline{v}}  + \mathbf{m}_v \mathscr{T}_{\partial_x (|u|^2)}  \mathscr{T}_{\overline{v}} 
   + \mathbf{m}_v \mathscr{T}_{\partial_x (\overline{u} v)}  \mathscr{T}_{\overline{u}} 
   + \mathbf{m}_u \mathscr{T}_{\partial_x (|v|^2)}  \mathscr{T}_{\overline{u}} 
  \bigr) \\
  = & -2   \mathbf{m}_{u \Pi_{\geq 0} \partial_x (\overline{v}u) + v \Pi_{\geq 0} \partial_x (|u|^2)} \mathscr{T}_{\overline{v}}  -2   \mathbf{m}_v \mathscr{T}_{\overline{u \Pi_{\geq 0} \partial_x (\overline{v}u) + v \Pi_{\geq 0} \partial_x (|u|^2)}}  \\ 
  & - 2  \mathbf{m}_{v \Pi_{\geq 0} \partial_x (\overline{u}v) + u \Pi_{\geq 0} \partial_x (|v|^2)} \mathscr{T}_{\overline{u}} 
  -   2  \mathbf{m}_u  \mathscr{T}_{\overline{v \Pi_{\geq 0} \partial_x (\overline{u}v) + u \Pi_{\geq 0} \partial_x (|v|^2)}} .
\end{split} 
\end{equation}Now we assume that $\psi = (u,v) \in \mathcal{Z}^2 \times \mathcal{Z}^2$, then $\partial_x u, \partial_x v \in H^1 (\mathbb{R}; \mathbb{C})$. 
  Formula \eqref{muscrTvbar,DCommFor} implies that
\begin{equation}\label{CommOmepsi,D}
[\Omega_{\psi}, D] = - \mathbf{m}_{i \partial_x^2 u} \mathscr{T}_{\overline{v}} 
-   \mathbf{m}_v \mathscr{T}_{\overline{i \partial_x^2 u}}   
- \mathbf{m}_{i \partial_x^2 v} \mathscr{T}_{\overline{u}} 
-   \mathbf{m}_u \mathscr{T}_{\overline{i \partial_x^2 v}}   
: H^1 (\mathbb{R} ; \mathbb{C}) \to L^2 (\mathbb{R} ; \mathbb{C})  .
\end{equation}  
So  \eqref{CommscrBpsiscrLpsi} follows from 
\eqref{scrLscrBidNLS}, \eqref{ttppsittqpsiDef},  \eqref{sumiAsqD-comomeA} and  \eqref{CommOmepsi,D}.
If $\psi = (u,v) \in \mathcal{C}\bigl( (T^-, T^+) ; \mathcal{Z}^3 \times \mathcal{Z}^3 \bigr)$   solves the idNLS system  \eqref{idNLS}, formula 
\eqref{HeiEqscrLscrBpsi} is a direct consequence of \eqref{CommscrBpsiscrLpsi}.
\end{proof}

\noindent Specializing the idNLS Lax pair \eqref{scrLscrBidNLS} to the case 
$u=v$  yields the non-chiral Lax pair for the focusing CMdNLS equation \eqref{CMdNLSfocu}. For any $u \in \mathcal{Z}^1$, we define
\begin{equation}\label{fracLuBuLaxPfoc}
\begin{cases}
\mathfrak{L}_u := \mathscr{L}_{\psi} \big|_{\psi = (u  ,u)} = D - \mathbf{m}_u \mathscr{T}_{\overline{u}} =  D - \mathbf{m}_u \Pi_{\geq 0}  \mathbf{m}_{\overline{u}}     \; \;  \;\; : H^1 (\mathbb{R} ; \mathbb{C}) \to L^2 (\mathbb{R} ; \mathbb{C}) ; \\
\mathfrak{B}_u : = \mathscr{B}_{\psi} \big|_{\psi = (u , u)}  = \mathbf{m}_u \mathscr{T}_{\partial_x \overline{u}}  - \mathbf{m}_{\partial_x u } \mathscr{T}_{  \overline{u}} + i (\mathbf{m}_u \mathscr{T}_{\overline{u}})^2  : H^1 (\mathbb{R} ; \mathbb{C}) \to L^2 (\mathbb{R} ; \mathbb{C}) .
\end{cases}
\end{equation}We set $\widetilde{\mathfrak{B}}_u 
 : =  \widetilde{\mathscr{B}}_{\psi} \big|_{\psi = (u , u)}$. Then Lemma 
$\ref{tilscrBpsiLem}$ and formula 
 \eqref{tildscrBpsiDef} 
 imply that 
\begin{equation}\label{tilBuNCDef}
\widetilde{\mathfrak{B}}_u 
 = \mathfrak{B}_u - i \mathfrak{L}_u^2  
= i \partial_x^2 + 2 \mathbf{m}_u \partial_x \mathscr{T}_{\overline{u}} : H^2 (\mathbb{R} ; \mathbb{C}) \to L^2 (\mathbb{R} ; \mathbb{C}).
\end{equation}
By    
Lemma $\ref{LemSelAdjscrLpsi}$,      $\mathfrak{L}_u^2 $ is a self-adjoint operator on $L^2 (\mathbb{R} ; \mathbb{C})$ with $\mathrm{Dom} (\mathfrak{L}_u^2) = H^2 (\mathbb{R} ; \mathbb{C})$. We have
\begin{equation}\label{fracLusqFor}
\mathfrak{L}_u^2 = - \partial_x^2 + i \partial_x \mathbf{m}_u \mathscr{T}_{\overline{u}} +   i \mathbf{m}_u \mathscr{T}_{\overline{u}} \partial_x  + (\mathbf{m}_u \mathscr{T}_{\overline{u}} )^2 : H^2 (\mathbb{R} ; \mathbb{C}) \to L^2 (\mathbb{R} ; \mathbb{C}), 
\end{equation}by   \eqref{scrLpsisqFor}. 
If $w \in \mathcal{Z}^2$, then  $\mathfrak{B}_w \in \mathcal{B}_{\mathbb{C}}(H^1 (\mathbb{R} ; \mathbb{C})) \cap \mathcal{B}_{\mathbb{C}}(L^2  (\mathbb{R} ; \mathbb{C}))$.  
  Using  \eqref{ttppsittqpsiDef} and 
 \eqref{CommscrBpsiscrLpsi}, we obtain  
\begin{equation}\label{CommfrakBufrLufor}
- [\mathfrak{B}_w , \mathfrak{L}_w] =
   \mathbf{m}_{ i \partial_x^2 w + 2 w \Pi_{\geq 0} \partial_x (|w|^2)} \mathscr{T}_{\overline{w}} + \mathbf{m}_w \mathscr{T}_{\overline{i \partial_x^2 w + 2 w \Pi_{\geq 0} \partial_x (|w|^2)}}: H^1 (\mathbb{R} ; \mathbb{C}) \to L^2 (\mathbb{R} ; \mathbb{C}) .
\end{equation}

\begin{lem}
For any $\epsilon > 0$ and $ u \in \mathcal{Z}^2$, let $\mathbf{r}_{\epsilon}(u)$ be defined in \eqref{DefbfrepuFor}. Then we have
\begin{equation}\label{bfrepuexprebytilBu}
\mathbf{r}_{\epsilon} (u) = \chi_{\epsilon} \bigl( i \partial_x^2 u + 2 u \Pi_{\geq 0} \partial_x \bigl( |u|^2 \bigr) \bigr) - \widetilde{\mathfrak{B}}_u \bigl(  u \chi_{\epsilon} \bigr) \in H^1 (\mathbb{R} ; \mathbb{C}) ,
\end{equation}where $\chi_{\epsilon} \in H^{\infty}_+$ is given by formula \eqref{chiepFor} and $u \chi_{\epsilon} \in H^2 (\mathbb{R} ; \mathbb{C})$. 
\end{lem}
\begin{proof}
Fix $u \in \mathcal{Z}^2$,  then  $ \partial_x (|u|^2) \in H^1 (\mathbb{R} ; \mathbb{R})$ by \eqref{uvZn<uZnvZn} and \eqref{ZhidNormZn}. For any $\epsilon > 0$,  Lemma $\ref{Z1+idealinL2+}$ and   \eqref{HNAlgebra} imply that  $\chi_{\epsilon} \Pi_{\geq 0} \partial_x (|u|^2) \in H^1_+$. Lemma  $\ref{repuHighRegLem}$  implies that $\mathbf{r}_{\epsilon} (u) \in H^1 (\mathbb{R} ; \mathbb{C})$. 
By  \eqref{fpi+fPI-fsum},  we have
\begin{equation}\label{uchiPidxumsqfor}
u \chi_{\epsilon} \Pi_{\geq 0} \partial_x (|u|^2) - u \Pi_{\geq 0} \bigl( \chi_{\epsilon}  \partial_x (|u|^2)  \bigr) = - u \Pi_{\geq 0} \bigl( \chi_{\epsilon} \Pi_{<0}  \partial_x (|u|^2)  \bigr) \in H^1 (\mathbb{R} ; \mathbb{C}).
\end{equation}
Inserting \eqref{uchiPidxumsqfor} into the expansion of $\widetilde{\mathfrak{B}}_u(u\chi_{\epsilon})$, we use \eqref{tilBuNCDef} and   \eqref{DefbfrepuFor}  to obtain \eqref{bfrepuexprebytilBu}.
\end{proof}

\noindent Fix  $\epsilon >0$ and $u \in \mathcal{Z}^2$. 
We define three auxiliary functions as follows:
\begin{equation}\label{bfjepsukDef}
\small 
\mathbf{j}^{(0)}_{\epsilon}(u) : = 2 \operatorname{Re} \langle u \chi_{\epsilon} , \mathbf{r}_{\epsilon} (u)  \rangle_{L^2  }, \quad \mathbf{j}^{(1)}_{\epsilon}(u) : = 2 \operatorname{Re} \langle \mathfrak{L}_u ( u \chi_{\epsilon} ) ,\mathbf{r}_{\epsilon} (u)  \rangle_{L^2  }, \quad  
\mathbf{j}^{(2)}_{\epsilon}(u) : = 2 \operatorname{Re} \langle \mathfrak{L}_u^2  ( u \chi_{\epsilon} ) , \mathbf{r}_{\epsilon} (u)  \rangle_{L^2  }.
\end{equation}
 
\begin{lem}
There exists a universal constant 
$\mathfrak{K}_0 > 0$ such that the following estimate holds:
\begin{equation}\label{bfjepsuUniEsti012}
\sup_{0 < \epsilon \leq 1} \bigl( \big|\mathbf{j}^{(2)}_{\epsilon}(u) + \mathbf{j}^{(1)}_{\epsilon}(u) \big|  +  \big| \mathbf{j}^{(1)}_{\epsilon}(u) + \mathbf{j}^{(0)}_{\epsilon}(u)  \big| \bigr) \leq \mathfrak{K}_0 \bigl( 1 + \|u\|_{\mathcal{Z}^2}^2 \bigr)^3 \bigl( 1 + \||u|^2 -1\|_{L^2 (\mathbb{R} ; \mathbb{C})}  \bigr) ,
\end{equation}for any $ u \in \mathcal{M}^2$,  where $\mathcal{M}^2$ is given by formula  \eqref{calMnmetrSp}. Moreover,   we have 
\begin{equation}\label{bfjepuLim012}
\lim_{\epsilon \to 0^+} \big| \mathbf{j}^{(2)}_{\epsilon}(u) + \mathbf{j}^{(1)}_{\epsilon}(u)  \big|  = \lim_{\epsilon \to 0^+} \big| \mathbf{j}^{(1)}_{\epsilon}(u) + \mathbf{j}^{(0)}_{\epsilon}(u)  \big| = 0 , \qquad \forall u \in \mathcal{M}^2 .
\end{equation}
\end{lem}
\begin{proof}
For any $\epsilon >0$, formula \eqref{dxNchiepL2norm} implies that $\|\chi_{\epsilon}\|_{L^{\infty}  } =1$ and $\|\partial_x \chi_{\epsilon} \|_{L^2_+}^2 = \frac{\pi}{2} \epsilon$. Fix $u \in \mathcal{M}^2$.  Since $\mathfrak{L}_u (u \chi_{\epsilon}) + u \chi_{\epsilon} = D (u \chi_{\epsilon} )  - u \Pi_{\geq 0} ((| u |^2 -1 )\chi_{\epsilon}) \in H^1(\mathbb{R} ; \mathbb{C})$ by 
\eqref{fracLuBuLaxPfoc}, we have
\begin{equation}\label{frLuuchi+uchiL2norm}
\sup_{0 < \epsilon \leq 1} \|
\partial_x (u \chi_{\epsilon} )   \|_{L^2} \leq \sqrt{\tfrac{\pi}{2}} \| u \|_{\mathcal{Z}^1} ,  \quad   \sup_{0 < \epsilon \leq 1}  \| \mathfrak{L}_u (u \chi_{\epsilon}) + u \chi_{\epsilon} \|_{L^2} \leq \sqrt{\tfrac{\pi}{2}} \| u \|_{\mathcal{Z}^1} (1 + \| |
u|^2 -1\|_{L^2}).  
\end{equation}Using formulas \eqref{fracLuBuLaxPfoc}, 
\eqref{bfrepyto0inL2} and \eqref{bfrepyUniBdinHN-1}, we obtain that
\begin{equation}\label{frLubfrepL2norm}
\lim_{\epsilon \to 0^+ } \|\mathfrak{L}_u (\mathbf{r}_{\epsilon} (u) )\|_{L^2  } = 0 ,\quad 
\sup_{0 < \epsilon \leq 1} \|\mathfrak{L}_u (\mathbf{r}_{\epsilon} (u) )\|_{L^2  } \leq \sqrt{2} \mathfrak{A}_2^* \| u \|_{\mathcal{Z}^2} \bigl( 1 + \| u \|_{\mathcal{Z}^1}^2 \bigr) \bigl( 1 + \| u \|_{L^{\infty}  }^2 \bigr) .
\end{equation}By formula  \eqref{bfjepsukDef},  for any $\epsilon  \in (0, 1]$ and $u \in \mathcal{M}^2$, we have
\begin{equation}\label{bfjepu1+0,2+1}
\small
\mathbf{j}_{\epsilon}^{(1)} (u) + \mathbf{j}_{\epsilon}^{(0)} (u) = 2 \operatorname{Re} \langle \mathfrak{L}_u (u \chi_{\epsilon}) + u \chi_{\epsilon}   , \mathbf{r}_{\epsilon} (u) \rangle_{L^2} , \quad \mathbf{j}_{\epsilon}^{(2)} (u) + \mathbf{j}_{\epsilon}^{(1)} (u) = 2 \operatorname{Re} \langle \mathfrak{L}_u (u \chi_{\epsilon}) + u \chi_{\epsilon}   , \mathfrak{L}_u (\mathbf{r}_{\epsilon} (u) )\rangle_{L^2}.
\end{equation}Set $\mathfrak{K}_0 : = 4 \sqrt{2 \pi} + 2 \sqrt{\pi}   \mathfrak{A}_2^* > 0 $.  Then the estimate \eqref{bfjepsuUniEsti012}  follows directly by 
\eqref{frLuuchi+uchiL2norm},  \eqref{frLubfrepL2norm},  \eqref{bfjepu1+0,2+1} and 
the Cauchy--Schwarz inequality. Formula \eqref{bfjepuLim012} is a direct consequence of 
 \eqref{bfrepyto0inL2}, 
 \eqref{frLubfrepL2norm} and    \eqref{bfjepu1+0,2+1}.
\end{proof}
 
\noindent Fix  $\epsilon >0$ and $u \in \mathcal{Z}^3$. We have $u \chi_{\epsilon} \in H^3 (\mathbb{R} ; \mathbb{C})$ and $\mathbf{r}_{\epsilon} (u) \in H^2 (\mathbb{R} ; \mathbb{C})$ by Lemma $\ref{repuHighRegLem}$. We set
\begin{equation}\label{defbfjepsu4}
\mathbf{j}^{(4)}_{\epsilon}(u) : = 2 \operatorname{Re} \langle \mathfrak{L}_u^2  ( u \chi_{\epsilon} ) , \mathfrak{L}_u^2  ( \mathbf{r}_{\epsilon} (u))  \rangle_{L^2  }.
\end{equation}
 \begin{lem}
There exists a universal constant 
$\mathfrak{K}_1 > 0$ such that the following estimate holds:
\begin{equation}\label{UniEsti24bfjepsu}
\sup_{0 < \epsilon \leq 1}   \big|\mathbf{j}^{(4)}_{\epsilon}(u) -  \mathbf{j}^{(2)}_{\epsilon}(u) \big|    \leq \mathfrak{K}_1 \bigl( 1 + \|u\|_{\mathcal{Z}^3}^2 \bigr)^5 \bigl( 1 + \||u|^2 -1\|_{H^1 (\mathbb{R} ; \mathbb{C})}  \bigr) ,
\end{equation}for any $ u \in \mathcal{M}^3$, where $ \mathcal{M}^3  = \mathcal{M}^2 \cap \mathcal{Z}^3$ is given by \eqref{calMnmetrSp}. Moreover,   we have 
\begin{equation}\label{Lim42bfjepu}
\lim_{\epsilon \to 0^+} \big| \mathbf{j}^{(4)}_{\epsilon}(u) -  \mathbf{j}^{(2)}_{\epsilon}(u)  \big|  = 0 , \qquad \forall u \in \mathcal{M}^3 .
\end{equation}
\end{lem}
\begin{proof}
If $\epsilon > 0$, formula \eqref{dxNchiepL2norm} yields that $ \| \partial
_x^2 \chi_{\epsilon} \|_{L^2}^2 = \frac{3 \pi}{2} \epsilon^3$, $\| \partial
_x  \chi_{\epsilon} \|_{L^{\infty}} = \epsilon $ and $\|\chi_{\epsilon} \|_{\mathcal{Z}^1} = 1 + \sqrt{\frac{\pi}{2} \epsilon}$. For any $u \in \mathcal{M}^3$, we have  $|u|^2 -1 \in H^3(\mathbb{R} ; \mathbb{C})$ by \eqref{calMnbyHndef2}.  
  Using  \eqref{dxNchiepL2norm} and 
\eqref{ufHn<uZnfHn}, we obtain that
\begin{equation}\label{Lapluchiep}
\small 
\sup_{ 0 < \epsilon \leq 1} \| \partial_x^2 (u \chi_{\epsilon})\|_{L^2} \leq 2 \sqrt{2} \| u \|_{\mathcal{Z}^2}, \quad
\sup_{ 0 < \epsilon \leq 1} \|D \bigl( u \Pi_{\geq 0} \bigl( ( |u |^2 -1 ) \chi_{\epsilon} \bigr) \bigr) \|_{L^2} \leq (1 + \sqrt{\tfrac{\pi}{2}}) C_1^2 \|u\|_{\mathcal{Z}^1} \| |u|^2  - 1 \|_{H^1}.
\end{equation}By    \eqref{fracLuBuLaxPfoc} and \eqref{fracLusqFor}, 
 the term 
$\mathfrak{L}_u^2 ( u \chi_{\epsilon} ) - u \chi_{\epsilon}  \in H^1 (\mathbb{R} ; \mathbb{C})$ can be expressed as follows:
\begin{equation}\label{frLusquchi-uchi}
\begin{split}
\mathfrak{L}_u^2 ( u \chi_{\epsilon} ) - u \chi_{\epsilon}  
= & D^2 (u \chi_{\epsilon} ) - D \bigl( u \Pi_{\geq 0} \bigl( ( |u |^2 -1 ) \chi_{\epsilon} \bigr) \bigr) -  D  (u \chi_{\epsilon} )  - u \Pi_{\geq 0} \bigl( \overline{u} D (u \chi_{\epsilon} ) \bigr) \\ 
& \qquad \; \; \; \quad + u \Pi_{\geq 0} \bigl( |u|^2 \Pi_{\geq 0} \bigl( ( |u |^2 -1 ) \chi_{\epsilon} \bigr) \bigr) +  u \Pi_{\geq 0} \bigl( ( |u |^2 -1 ) \chi_{\epsilon} \bigr).
\end{split}
\end{equation}Set $\widetilde{\mathsf{K}} : = 2 \sqrt{2} + (1 + \sqrt{\frac{\pi}{2}}) C_1^2 + \sqrt{\frac{\pi}{2}} >0$. 
Using 
\eqref{Lapluchiep}, \eqref{frLusquchi-uchi} and
\eqref{frLuuchi+uchiL2norm}, we obtain that \begin{equation}\label{Lusquchiep-uchiepL2}
\small
\begin{split}
  \sup_{0 < \epsilon \leq 1} \| \mathfrak{L}_u^2 ( u \chi_{\epsilon} ) - u \chi_{\epsilon} \|_{L^2}  
 \leq &  2 \sqrt{2} \| u \|_{\mathcal{Z}^2} + (1 + \sqrt{\tfrac{\pi}{2}}) C_1^2 \|u\|_{\mathcal{Z}^1} \| |u|^2  - 1 \|_{H^1} + \sqrt{\tfrac{\pi}{2}} \| u \|_{\mathcal{Z}^1} \\
 & \qquad    +  \sqrt{\tfrac{\pi}{2}} \| u \|_{\mathcal{Z}^1} \|u \|_{L^{\infty}}^2   + \| u \|_{L^{\infty}}^3 \||u|^2 - 1 \|_{L^2} + \| u \|_{L^{\infty}} \||u|^2 - 1 \|_{L^2}  \\
 \leq &  \widetilde{\mathsf{K}} \| u \|_{\mathcal{Z}^2} \bigl( 1 + \| u \|_{L^{\infty}}^2 \bigr) \bigl( 1 + \| |u|^2 - 1  \|_{H^1} \bigr)   .
\end{split}
\end{equation}
Since $u \in   \mathcal{M}^3 \subset \mathcal{Z}^3$, Lemma  $\ref{repuHighRegLem}$ yields that $\mathbf{r}_{\epsilon} (u) \in H^2 (\mathbb{R} ; \mathbb{C})$. 
From   \eqref{ufHn<uZnfHn} and    \eqref{fracLusqFor}, it follows that 
\begin{equation}\label{fracLusqbfrepuL2smal}
\|\mathfrak{L}_u^2 (\mathbf{r}_{\epsilon} (u) )\|_{L^2 (\mathbb{R} ; \mathbb{C})} \leq \| \mathbf{r}_{\epsilon} (u) \|_{H^2 (\mathbb{R} ; \mathbb{C})} \bigl( 1 + C_1^2 \| u \|_{\mathcal{Z}^1}^2 + \| u \|_{L^{\infty} (\mathbb{R} ; \mathbb{C})}^2  + \| u \|_{L^{\infty} (\mathbb{R} ; \mathbb{C})}^4 \bigr), \quad \forall \epsilon >0. 
\end{equation}By formulas  \eqref{ZhidNormZn}, 
\eqref{fracLusqbfrepuL2smal}  and   \eqref{bfrepyUniBdinHN-1}, 
for any $u \in \mathcal{M}^3$, 
we have
\begin{equation}\label{CVfrakLusqreputo0L2}
\sup_{ 0 < \epsilon \leq 1} \|\mathfrak{L}_u^2 (\mathbf{r}_{\epsilon} (u) )\|_{L^2 (\mathbb{R} ; \mathbb{C})} \leq \mathfrak{A}_3^* (1 + C_1^2 ) \| u \|_{\mathcal{Z}^3} \bigl( 1 + \| u \|_{\mathcal{Z}^2}^2  \bigr)^3, \quad  
\lim_{\epsilon \to  0^+ } \|\mathfrak{L}_u^2 (\mathbf{r}_{\epsilon} (u) )\|_{L^2 (\mathbb{R} ; \mathbb{C})} = 0 .
\end{equation}Fix $u \in \mathcal{M}^3$ and $\epsilon \in (0 , 1]$.
 Formulas  \eqref{defbfjepsu4}  and  \eqref{bfjepsukDef}  imply that 
\begin{equation}\label{forbfjepu4-2}
\mathbf{j}_{\epsilon}^{(4)} (u)  -  \mathbf{j}_{\epsilon}^{(2)} (u) = 2 \operatorname{Re} \langle \mathfrak{L}_u^2  (u \chi_{\epsilon}) -  u \chi_{\epsilon}   , \mathfrak{L}_u^2  (\mathbf{r}_{\epsilon} (u) ) \rangle_{L^2} .
\end{equation}Set $\mathfrak{K}_1 : = 2 \mathfrak{A}_3^* (1 + C_1^2) \widetilde{\mathsf{K}} =    2 \mathfrak{A}_3^* (1 + C_1^2) \bigl( 2 \sqrt{2} + (1 + \sqrt{\frac{\pi}{2}}) C_1^2 + \sqrt{\frac{\pi}{2}}  \bigr) > 0$. Consequently, formulas
\eqref{UniEsti24bfjepsu} and
\eqref{Lim42bfjepu} follow directly from 
\eqref{Lusquchiep-uchiepL2}, 
\eqref{CVfrakLusqreputo0L2}, 
 \eqref{forbfjepu4-2} and the Cauchy--Schwarz inequality.
\end{proof}

\begin{lem}\label{sfwu1defLem}
 For any $ u \in \mathcal{M}^2$, we set $ \mathsf{w}^{(1)}_u : = D u - u \Pi_{\geq 0} (|u|^2 -1) \in H^1  (\mathbb{R} ; \mathbb{C})$. Then we have
\begin{equation}\label{frLu+1uchieLim}
 \lim_{\epsilon \to 0^+}  \| \mathfrak{L}_u (u \chi_{\epsilon}) + u \chi_{\epsilon} - \mathsf{w}^{(1)}_u \|_{H^1 (\mathbb{R} ; \mathbb{C}) }= 0 .
\end{equation}
\end{lem}
\begin{proof}
Fix $u \in \mathcal{M}^2$, then $u \in \mathcal{Z}^2 $ and $| u |^2 - 1 \in H^2 (\mathbb{R} ; \mathbb{C})$. From  \eqref{chiephCVtohinHpFor} and 
\eqref{gteputo0andBd},  it follows that
\begin{equation}\label{cvdxuchi-uunsq-1}
\lim_{\epsilon \to 0^+} \| \partial_x (u \chi_{\epsilon} - u )\|_{H^1 (\mathbb{R}; \mathbb{C})} = \lim_{\epsilon \to 0^+} \| \partial_x \mathrm{f}_{\epsilon}^u \|_{H^1 (\mathbb{R}; \mathbb{C})} = 0 ,\quad \lim_{\epsilon \to 0^+} \| (|u|^2 -1) (\chi_{\epsilon} - 1 )\|_{H^2 (\mathbb{R}; \mathbb{C})} =  0 .
\end{equation}Since $ \mathfrak{L}_u (u \chi_{\epsilon}) + u \chi_{\epsilon} - \mathsf{w}^{(1)}_u  =  D (u \chi_{\epsilon} - u ) + u \Pi_{\geq 0} \bigl((|u|^2 -1) (\chi_{\epsilon} - 1 ) \bigr) \in H^1 (\mathbb{R} ; \mathbb{C})$ by   \eqref{fracLuBuLaxPfoc},  the proof follows directly from formulas  \eqref{cvdxuchi-uunsq-1} and \eqref{ufHn<uZnfHn}.
\end{proof}
  
\begin{lem}\label{LemfrLusq-1uchielim}
 For any $ u \in \mathcal{M}^2$, set $\mathsf{w}^{(2)}_u : = \mathfrak{L}_u (\mathsf{w}^{(1)}_u) - \mathsf{w}^{(1)}_u \in L^2 (\mathbb{R} ; \mathbb{C})$.   
 Equivalently, we have  
\begin{equation}\label{forsfwu2def}
\small
 \mathsf{w}^{(2)}_u  
 =   D^2 u - D \bigl( u \Pi_{\geq 0} (|u|^2 -1) \bigr)  - Du - u \Pi_{\geq 0} (\overline{u} D u) + u \Pi_{\geq 0} \bigl( |u|^2 \Pi_{\geq 0} (|u|^2 -1)  \bigr) + u \Pi_{\geq 0} (|u|^2 -1).
\end{equation}Then  $\lim_{\epsilon \to 0^+} \|\mathfrak{L}_u^2 (u \chi_{\epsilon}) - u \chi_{\epsilon} - \mathsf{w}^{(2)}_u \|_{L^2 (\mathbb{R} ; \mathbb{C})} = 0$.
\end{lem}
\begin{proof}
This is a direct consequence of
 \eqref{fracLuBuLaxPfoc},  \eqref{cvdxuchi-uunsq-1},  \eqref{ufHn<uZnfHn},   \eqref{frLusquchi-uchi},  \eqref{forsfwu2def} and Lemma $\ref{sfwu1defLem}$.
\end{proof}

\subsection{The Chiral Lax Pair Formulation}

\noindent In the previous subsection, we studied the non-chiral Lax pair formulation, where the operators act on the full Sobolev space $H^1(\mathbb{R};\mathbb{C})$. From this subsection to the end of Section $\ref{SectLaxPair3}$, we restrict our attention to the case where all operators act on the Hardy--Sobolev space $H^s_+$, for some $s \geq 0$. We refer to this as the \emph{chiral formulation}. \\

\noindent For any $u \in L^{\infty}(\mathbb{R}; \mathbb{C})$, formula \eqref{AdjofToepOpinL2+} implies that $\mathbf{T}_u \mathbf{T}_{\overline{u}}$ is a bounded positive operator on the Hardy space $L^2_+$.  
According to Lemma $\ref{DNSelfadjinHsLem}$,  $D = - i\partial_x $ is an unbounded positive self-adjoint operator on $L^2_+$ with $\mathrm{Dom} (D) = \mathrm{Dom}^0 (D) = H^1_+$.
 The chiral Lax operator $\mathbf{L}_u : H^1_+ \to L^2_+$ is defined by
\begin{equation}\label{LuSect2}
\mathbf{L}_u := D - \mathbf{T}_u \mathbf{T}_{\overline{u}}, \qquad \mathrm{Dom}(\mathbf{L}_u) = H^1_+.
\end{equation}
\begin{lem}\label{lemLuselfadj}
For any $u \in L^{\infty}(\mathbb{R}; \mathbb{C})$, the chiral Lax operator $\mathbf{L}_u :  H^1_+ \to L^2_+$ is an unbounded self-adjoint operator on the Hardy space $L^2_+$.
\end{lem}

\noindent Assume that $u \in \mathcal{Z}^1$. Then formula \eqref{BounToeplitzBHs} yields that $\mathbf{T}_{\overline{u}} \in \mathcal{B}_{\mathbb{C}}(L^2_+) \cap \mathcal{B}_{\mathbb{C}}(H^1_+)$. We set 
\begin{equation}\label{BudefForFoc-}
\mathbf{B}_u  := \mathbf{T}_u \mathbf{T}_{\partial_x \overline{u}} - \mathbf{T}_{\partial_x u} \mathbf{T}_{\overline{u}} + i \left( \mathbf{T}_u \mathbf{T}_{\overline{u}}\right)^2 : H^1_+ \to L^2_+ .
\end{equation}
Since $u \in \mathcal{Z}^1$, we have 
$\mathbf{L}_u \in \mathcal{B}_{\mathbb{C}}(H^2_+;  H^1_+) \cap \mathcal{B}_{\mathbb{C}}(H^1_+;  L^2_+)$ by \eqref{BounToeplitzBHs}. By formula \eqref{mudxscrTfor}, we have
\begin{equation}\label{bfTudxTubarfor}
\mathbf{T}_u \partial_x \mathbf{T}_{\overline{u}} = \mathbf{T}_u \mathbf{T}_{\partial_x \overline{u}} +  \mathbf{T}_u \mathbf{T}_{\overline{u}} \partial_x  = \partial_x \mathbf{T}_u \mathbf{T}_{\overline{u}} - \mathbf{T}_{\partial_x u}  \mathbf{T}_{\overline{u}}  \in \mathcal{B}_{\mathbb{C}}(H^1_+;  L^2_+) .
\end{equation}For any $u \in \mathcal{Z}^1$,  we  set $\widetilde{\mathbf{B}}_u  := \mathbf{B}_u - i \mathbf{L}_u^2$. From  formulas 
\eqref{LuSect2}, 
\eqref{BudefForFoc-} and 
\eqref{bfTudxTubarfor}, it follows that
\begin{equation}\label{tildeBudefFoc-}
\widetilde{\mathbf{B}}_u    = i\partial_x^2 + 2 \mathbf{T}_u \partial_x \mathbf{T}_{  \overline{u}} \in \mathcal{B}_{\mathbb{C}}(H^2_+; L^2_+)  .
\end{equation}The next lemma is a direct consequence of  the Kato--Rellich theorem and estimate
\eqref{relBoundofmudxTvbar}.

\begin{lem}\label{LusqiLemselfadj>0}
For any $u \in \mathcal{Z}^1$,  both $\mathbf{L}_u^2$ and $ i \widetilde{\mathbf{B}}_u   $ are  unbounded positive self-adjoint operator on $L^2_+$, with  $\mathrm{Dom}(\mathbf{L}_u^2) = \mathrm{Dom}( i \widetilde{\mathbf{B}}_u   ) = H^2_+$. Moreover, the following identity holds:
\begin{equation}\label{Lu2expForTuTubarpartx}
 \mathbf{L}_u^2 = - \partial_x^2 + i \partial_x \mathbf{T}_u \mathbf{T}_{\overline{u}} + i \mathbf{T}_u \mathbf{T}_{\overline{u}} \partial_x + \left(  \mathbf{T}_u \mathbf{T}_{\overline{u}}   \right)^2 : H^2_+ \to L^2_+ .
 \end{equation}
\end{lem}

\noindent Fix $u \in \mathcal{Z}^1$. The non-chiral Lax pairs $\bigl( \mathfrak{L}_u,  \mathfrak{B}_u \bigr)$ and $\bigl( \mathfrak{L}_u,  \widetilde{\mathfrak{B}}_u \bigr)$ are defined by \eqref{fracLuBuLaxPfoc} and \eqref{tilBuNCDef}. Then 
\begin{equation}\label{bfLufrakLuRel}
\mathbf{L}_u = \Pi_{\geq 0} \mathfrak{L}_u |_{H^1_+} \in \mathcal{B}_{\mathbb{C}} \bigl( H^1_+ ; L^2_+ \bigr), \qquad \widetilde{\mathbf{B}}_u = \Pi_{\geq 0} \widetilde{\mathfrak{B}}_u  |_{H^1_+} \in \mathcal{B}_{\mathbb{C}} \bigl( H^2_+ ; L^2_+ \bigr) .
\end{equation}
If $w \in L^{\infty}_+$ and $v \in L^{\infty} ( \mathbb{R} ; \mathbb{C})$, 
Lemma $
\ref{Z1+idealinL2+}$ implies that $\mathbf{T}_w \mathbf{T}_{\overline{v}} = \mathbf{m}_w \mathscr{T}_{\overline{v}} \big|_{L^2_+} \in \mathcal{B}_{\mathbb{C}}( L^2_+ ) $. Then we have
\begin{equation}\label{bfLufrakLuHarSob}
\mathbf{L}_{\mathtt{u}} = \mathfrak{L}_{\mathtt{u}} \big|_{H^1_+} : H^1_+ \to L^2_+,  \quad \mathbf{B}_{\mathtt{u}} = \mathfrak{B}_{\mathtt{u}} \big|_{H^1_+} : H^1_+ \to L^2_+, \quad \widetilde{\mathbf{B}}_{\mathtt{u}} = \widetilde{\mathfrak{B}}_{\mathtt{u}} \big|_{H^2_+} : H^2_+ \to L^2_+, \quad \forall \mathtt{u} \in \mathcal{Z}^1_+ .
\end{equation} 
The following lemma provides a useful estimate for the action of $\widetilde{\mathbf{B}}_u$ on the approximating functions $\chi_\epsilon$, which will be essential in the derivation of explicit formula of \eqref{CMdNLSfocu}.
 
\begin{lem}
Let $\chi_{\epsilon}\in \mathrm{Dom}^{\infty}(\mathbf{G})$ be defined in \eqref{chiepFor}, $\forall \epsilon >0$. For any $u \in \mathcal{Z}^1_+$, we have
\begin{equation}\label{wtilBuchiepto0Uni}
\lim_{\epsilon \to 0^+} \|\widetilde{\mathbf{B}}_u (\chi_{\epsilon})\|_{L^2_+} =0 \quad \mathrm{and} \quad
\sup_{0< \epsilon \leq 1} \|\widetilde{\mathbf{B}}_u (\chi_{\epsilon})\|_{L^2_+} 
\leq \sqrt{2 \pi}  \left( 1+ \|u\|_{Z^1_+}^2 \right) .
\end{equation}
\end{lem}
\begin{proof}
Since $u \in \mathcal{Z}^1_+$, we have $\partial_x u \in L^2_+$ and $\Pi_{\geq 0}(\partial_x \overline{u})  = 0$. Then we have
\begin{equation}\label{Pichiepdxubarto0}
\limsup_{\epsilon \to 0^+} \|\Pi_{\geq 0} \left( \chi_{\epsilon} \partial_x \overline{u} \right)\|_{L^2_+} 
\leq \lim_{\epsilon \to 0^+} \| \chi_{\epsilon} \partial_x \overline{u} - \overline{u}\|_{L^2 (\mathbb{R}; \mathbb{C})} = 0,
\end{equation}by \eqref{chiephCVtohinHpFor}.
Using formula  \eqref{tildeBudefFoc-} and Lemma $\ref{Z1+idealinL2+}$, for any $\epsilon \in (0,1]$, we obtain that
\begin{equation}\label{tilBuchiepEsti}
\begin{split}
& \|\widetilde{\mathbf{B}}_u (\chi_{\epsilon})\|_{L^2_+}  = \| i \partial_x^2 \chi_{\epsilon} + 2 u \Pi_{\geq 0} \left( \chi_{\epsilon} \partial_x \overline{u} \right) + 2 u \Pi_{\geq 0}\left( \overline{u}  \partial_x  \chi_{\epsilon} \right)\|_{L^2_+} \\
\leq & \|\partial_x^2 \chi_{\epsilon} \|_{L^2_+} 
+ 2 \|u\|_{L^{\infty}(\mathbb{R}; \mathbb{C})} \left( \|\Pi_{\geq 0} \left( \chi_{\epsilon} \partial_x \overline{u} \right)\|_{L^2_+} 
+  \|u \|_{L^{\infty}(\mathbb{R}; \mathbb{C})} \| \partial_x \chi_{\epsilon} \|_{L^2_+} \right) .
\end{split}
\end{equation}Formula \eqref{chiephCVtohinHpFor} yields that $\lim_{\epsilon \to 0^+}\| \partial_x^p \chi_{\epsilon} \|_{L^2_+} =0$, for any $p \in \{1,2\}$. By using \eqref{Pichiepdxubarto0} and \eqref{tilBuchiepEsti}, we have $\lim_{\epsilon \to 0^+} \|\widetilde{\mathbf{B}}_u (\chi_{\epsilon})\|_{L^2_+} =0$. Formula \eqref{chiepFor} implies that $\sup_{\epsilon >0}\|\chi_{\epsilon}\|_{L^{\infty}(\mathbb{R}; \mathbb{C})} \leq 1$. So we have
\begin{equation}\label{PichiepdxubarL2estim}
\|\Pi_{\geq 0}\left( \chi_{\epsilon} \partial_x \overline{u} \right)\|_{L^2_+}  \leq \|\chi_{\epsilon}\|_{L^{\infty}(\mathbb{R}; \mathbb{C})}\| \partial_x u\|_{L^2_+} \leq \| \partial_x u\|_{L^2_+}, \quad \forall \epsilon >0.
\end{equation}Then the second formula of \eqref{wtilBuchiepto0Uni}  is a direct consequence of \eqref{ZhidNormZn}, \eqref{chiephCVtohinHpFor},  \eqref{tilBuchiepEsti} and  \eqref{PichiepdxubarL2estim}.
\end{proof}

 \noindent For any $a \in \mathbb{C}$, $u \in a + H^1_+ \subset \mathcal{Z}^1_+$ and $\nu \in \mathbb{C} \setminus \sigma (\mathbf{L}_u)$, the generating functional $\mathcal{H}_{\nu}^a (u)$ is defined by
 \begin{equation}\label{calHGeneFctdefFor}
 \mathcal{H}_{\nu}^a (u) = \langle \bigl( \mathbf{L}_u - 
 \nu \bigr)^{-1} (u - a) , u - a\rangle_{L^2_+} \in \mathbb{C}.
 \end{equation}
 Lemma $\ref{lemLuselfadj}$ shows  that $\mathbf{L}_u$ is a self-adjoint operator on $L^2_+$.    
 Let $\mu^{\mathbf{L}_u}_{ u - a }$ denote the spectral measure of  $\mathbf{L}_u$ associated to  $u-a \in H^1_+ \subset L^2_+$. Then   $\zeta \in \mathbb{C} \setminus \mathbb{R} \mapsto  \mathcal{H}_{ \zeta }^a (u)  \in \mathbb{C}$ is the Cauchy transform of  $\mu^{\mathbf{L}_u}_{ u - a }$, i.e. 
 \begin{equation}\label{CauTrSpeMeaLuisH}
  \mathcal{H}_{\zeta }^a (u) = \int_{\mathbb{R}} (\xi -\zeta )^{-1}   \mathrm{d}\mu^{\mathbf{L}_u}_{ u - a } (\xi) , \qquad \forall \zeta  \in \mathbb{C} \setminus \mathbb{R}.
 \end{equation}
   For any $u \in \mathcal{Z}^2$,  Lemma $\ref{LusqiLemselfadj>0}$ and the Stone theorem imply that $\bigl(\exp (i t \mathbf{L}_u^2 ) \bigr)_{t \in \mathbb{R}}$ is a strongly continuous unitary group on the Hardy space $L^2_+$ such that $e^{i t \mathbf{L}_u^2}  (H^2_+) = H^2_+ = \mathrm{Dom} (\mathbf{L}_u^2)$, for any $t \in \mathbb{R}$. 
 To determine whether  $H^1_+ = \mathrm{Dom} (\mathbf{L}_u)$  is also invariant under this group, we first establish a general lemma showing that if two self-adjoint operators commute (in the sense that their generated unitary groups commute), then each group preserves the domain of the other operator. 
 
\begin{lem}\label{LemABcommuteexpCom}
 Let $\left(\mathtt{A}, \mathrm{Dom}(\mathtt{A}) \right)$ and $\left(\mathtt{B}, \mathrm{Dom}(\mathtt{B}) \right)$ be two self-adjoint operators on some Hilbert space $\mathtt{H}$. Assume that $\mathtt{A}$ and $\mathtt{B}$ commute in the following sense: for any $t,s \in \mathbb{R}$, we have
\begin{equation}\label{expitAcommuteexpisB}
\exp(it \mathtt{A}) \exp(is \mathtt{B}) =  \exp(is \mathtt{B}) \exp(it \mathtt{A}) \in \mathcal{B}_{\mathbb{C}}(\mathtt{H}).
\end{equation}Then  for any $t ,s \in \mathbb{R}$,  we have $\exp(it \mathtt{A}) \left(\mathrm{Dom}(\mathtt{B}) \right) = \mathrm{Dom}(\mathtt{B})$, 
$\exp(is \mathtt{B}) \left(\mathrm{Dom}(\mathtt{A}) \right) = \mathrm{Dom}(\mathtt{A})$, and
\begin{equation}\label{eitABcommuteFor}
\exp(it \mathtt{A}) \mathtt{B} = \mathtt{B} \exp(it \mathtt{A}) : \mathrm{Dom}(\mathtt{B}) \to \mathtt{H}, \quad \exp(is \mathtt{B}) \mathtt{A} = \mathtt{A} \exp(is \mathtt{B}) : \mathrm{Dom}(\mathtt{A}) \to \mathtt{H}.
\end{equation}
\end{lem}
 
\begin{proof}
According to Stone theorem (Theorem \uppercase\expandafter{\romannumeral8}.7 of \cite{ReedSimon1FA}), it follows that
\begin{equation}\label{DomttBCVUniGrStonethm}
\mathrm{Dom}(\mathtt{B}) = \{\varphi \in \mathtt{H}: \exists g \in \mathtt{H}  \quad \text{such that} \quad \lim_{s \to 0} \|\tfrac{1}{s}(\exp(is \mathtt{B}) \varphi - \varphi) - g\|_{\mathtt{H}} = 0\}.
\end{equation}Fix $t \in \mathbb{R}$ and $\varphi \in \mathrm{Dom}(\mathtt{B})$. We set $f: = \exp(it \mathtt{A} ) \varphi \in \mathtt{H}$. Since $ \exp(it \mathtt{A}) \in \mathcal{B}_{\mathbb{C}}(\mathtt{H})$, we have
\begin{equation}
\tfrac{1}{s}\left(\exp(is \mathtt{B}) f -f  \right) = \exp(it \mathtt{A}) \left( \tfrac{\exp(is \mathtt{B})  \varphi - \varphi }{s}  \right) \to \exp(it \mathtt{A}) i \mathtt{B} \varphi  \quad \text{in} \quad \mathtt{H}, \quad \text{as} \quad s \to 0,
\end{equation}by \eqref{expitAcommuteexpisB}  and \eqref{DomttBCVUniGrStonethm}. So $\exp(it \mathtt{A}) \varphi = f \in \mathrm{Dom}(\mathtt{B})$ and $i\mathtt{B}\exp(it \mathtt{A}) \varphi = i  \mathtt{B}f = \exp(it \mathtt{A}) i \mathtt{B} \varphi  \in \mathtt{H}$. For any $t \in \mathbb{R}$, since $\exp(it \mathtt{A}) \left(\mathrm{Dom} ( \mathtt{B}) \right) + \exp(-it \mathtt{A}) \left(\mathrm{Dom} ( \mathtt{B}) \right) \subset \mathrm{Dom} ( \mathtt{B})$, we have  
\begin{equation}\label{expittAttBfor}
\exp(it \mathtt{A}) \left(\mathrm{Dom} ( \mathtt{B}) \right) = \mathrm{Dom} ( \mathtt{B}) \quad \mathrm{and} \quad  \mathtt{B}\exp(it \mathtt{A})   =  \exp(it \mathtt{A})   \mathtt{B}  :  \mathrm{Dom} ( \mathtt{B}) \to \mathtt{H}.
\end{equation}The proof is completed by swapping the roles of $\mathtt{A}$ and $\mathtt{B}$ in formula \eqref{expittAttBfor}.
\end{proof}
  
\noindent Applying the previous lemma to $\bigl(\mathbf{L}_u, \mathrm{Dom}(\mathbf{L}_u) \bigr)$ and $\bigl(\mathbf{L}_u^2, \mathrm{Dom}(\mathbf{L}_u^2) \bigr)$, we obtain the following result, which shows in particular that $H^1_+$ is preserved under the group $\bigl(\exp (i t \mathbf{L}_u^2 ) \bigr)_{t \in \mathbb{R}}$. 

\begin{lem}\label{LucommutLusqlem}
 For any $u \in \mathcal{Z}^2$, the two  self-adjoint operators $\mathbf{L}_u : H^1_+ \to L^2_+$ and $\mathbf{L}_u^2 : H^2_+ \to L^2_+$ commute in the following sense: for any $z , \mu \in \mathbb{C} \setminus \mathbb{R}$, we have
 \begin{equation}\label{resolvLuLusqcommute}
 \left(z - \mathbf{L}_u \right)^{-1}  \left(\mu - \mathbf{L}_u^2 \right)^{-1} =  \left(\mu - \mathbf{L}_u^2 \right)^{-1}  \left(z - \mathbf{L}_u \right)^{-1}  : L^2_+ \to H^3_+.
 \end{equation}As a consequence, for any $t,s \in \mathbb{R}$, the following identity holds:
 \begin{equation}\label{eitLueisLusqcommute}
 \exp(it \mathbf{L}_u) \exp(is \mathbf{L}_u^2) = \exp(is \mathbf{L}_u^2) \exp(it \mathbf{L}_u) \in \mathcal{B}_{\mathbb{C}}(L^2_+).
 \end{equation}Moreover, for any $s \in \mathbb{R}$, $a \in \mathbb{C}$ and $\kappa \in \mathbb{C}\setminus \sigma(\mathbf{L}_u)$, we have $\exp(i s \mathbf{L}_u^2) (H^1_+ ) = H^1_+$ and  
\begin{equation}\label{expisLusqcommutLu-a}
 \exp(i s \mathbf{L}_u^2) \left( \mathbf{L}_u - a  \right) = \left( \mathbf{L}_u - a  \right)\exp(i s \mathbf{L}_u^2),    \quad 
 \exp(i s \mathbf{L}_u^2) \left( \mathbf{L}_u - \kappa  \right)^{-1} = \left( \mathbf{L}_u - \kappa \right)^{-1}\exp(i s \mathbf{L}_u^2).
\end{equation}  
\end{lem}
\begin{proof}
Fix $z, \mu \in \mathbb{C} \setminus \mathbb{R}$ and $f \in L^2_+$. Set $g : = \left( z - \mathbf{L}_u\right)^{-1} f \in H^1_+$ and $h : = \left( \mu - \mathbf{L}_u^2 \right)^{-1} g \in H^2_+$. Since $u \in \mathcal{Z}^2$, we have $\mathbf{T}_u \mathbf{T}_{\overline{u}}(H^2_+) \subset H^2_+$ and $\mathbf{T}_u \mathbf{T}_{\overline{u}}(H^1_+) \subset H^1_+$ by \eqref{HnZnsubsetHn}. So formula \eqref{Lu2expForTuTubarpartx} yields that
\begin{equation}\label{hinH3+fromLaphinH1+}
\partial_x^2 h = g -\mu h + i \partial_x \left( \mathbf{T}_u \mathbf{T}_{\overline{u}} (h) \right) + i \mathbf{T}_u \mathbf{T}_{\overline{u}} (\partial_x h) + \left(\mathbf{T}_u \mathbf{T}_{\overline{u}}   \right)^2 (h) \in H^1_+ \Longrightarrow h \in H^3_+.
\end{equation}Recall that $\mathbf{L}_u (H^3_+) \subset H^2_+ $. We set $\phi: = \left( z - \mathbf{L}_u \right) h \in H^2_+$. Since $g = \left( \mu - \mathbf{L}_u^2 \right)  h  \in H^1_+$, we obtain that
$\left( \mu - \mathbf{L}_u^2 \right) \phi
 =    z \mu h - \mu \mathbf{L}_u (h) - z \mathbf{L}_u^2 (h) + \mathbf{L}_u^3 (h) = z \left( \mu - \mathbf{L}_u^2 \right) h  -  \mathbf{L}_u \left( \mu - \mathbf{L}_u^2 \right)  h  = \left( z - \mathbf{L}_u \right) g  =f \in L^2_+$. So $\left( z - \mathbf{L}_u  \right)^{-1} \left( \mu - \mathbf{L}_u^2 \right)^{-1} f = h = \left( \mu - \mathbf{L}_u^2 \right)^{-1}  \left( z - \mathbf{L}_u  \right)^{-1} f \in H^3_+$. Since formula  \eqref{resolvLuLusqcommute} holds for any $z, \mu \in \mathbb{C} \setminus \mathbb{R}$, 
  Theorem \uppercase\expandafter{\romannumeral8}.13 of \cite{ReedSimon1FA} yields that \eqref{eitLueisLusqcommute} holds for any $t,s \in \mathbb{R}$.  
  Finally, formula  
  \eqref{expisLusqcommutLu-a} follows from \eqref{eitLueisLusqcommute} by Lemma \ref{LemABcommuteexpCom}.
\end{proof}
  
\noindent The following lemma upgrades the invariance of the Hardy--Sobolev spaces $H^p_+$ ($p = 0,1,2,3$) under $\exp(-it\mathbf{L}_u^2)$ to strong continuity on each of these spaces.
 
\begin{lem}\label{expitLusqStrContGroup}
For any $u \in \mathcal{Z}^2$, the family $\left(\exp(- i t \mathbf{L}_u^2) \right)_{t \in \mathbb{R}}$ forms a strongly continuous group of bounded operators on the Hardy--Sobolev space $H^p_+$, for any $p \in \{0,1,2,3\}$. In particular, when $p = 0$, it is a group of unitary operators on $L^2_+$.
\end{lem}
\begin{proof}
Fix $u \in \mathcal{Z}^2$. For any $t \in \mathbb{R}$, we set $\mathbb{U}(t):= \exp(-it \mathbf{L}_u^2)$. 
Since $\mathbf{L}_u^2$ is a self-adjoint operator on $L^2_+$ by Lemma $\ref{LusqiLemselfadj>0}$, Theorem \uppercase\expandafter{\romannumeral8}.7 of \cite{ReedSimon1FA} implies that $\left( \mathbb{U}(t)  \right)_{t \in \mathbb{R}}$ is a strongly continuous group of unitary operators on $L^2_+$. Recall that $\mathcal{E}^s = \mathcal{B}_{\mathbb{C}}(H^s_+)$, $\forall s \geq 0$. For any $\mathbf{q} \in L^2_+$,  we have
\begin{equation}\label{StrContiUtL2}
\lim_{t \to 0} \|\mathbb{U}(t) \mathbf{q} - \mathbf{q}\|_{L^2_+} = 0 \quad 
\mathrm{and} \quad 
\|\mathbb{U}(\mathfrak{t} ) \mathbf{q}\|_{L^2_+ } = \| \mathbf{q} \|_{L^2_+}  , \quad \forall \mathfrak{t} \in \mathbb{R}.
\end{equation}Let $C_1, C_2>0$ be defined in formula 
\eqref{BounToeplitzBHs}. Then \eqref{BounToeplitzBHs} yields that
\begin{equation}\label{BdTuTubar<uZ1Z2esti}
\|\mathbf{T}_u \mathbf{T}_{\overline{u}}\|_{\mathcal{E}^0} \leq  \|u\|_{L^{\infty}(\mathbb{R}; \mathbb{C})}^2,
\quad \|\mathbf{T}_u \mathbf{T}_{\overline{u}}\|_{\mathcal{E}^1} \leq C_1^2 \|u\|_{\mathcal{Z}^1}^2, \quad \|\mathbf{T}_u \mathbf{T}_{\overline{u}}\|_{\mathcal{E}^2} \leq C_2^2 \|u\|_{\mathcal{Z}^2}^2.
\end{equation}For any $g \in H^1_+$ and $ t \in \mathbb{R}$, it follows from formulas  \eqref{BdTuTubar<uZ1Z2esti} and  \eqref{expisLusqcommutLu-a}    that
\begin{equation}\label{LugnormUtLucomm2for}
\| \mathbf{L}_u ( g ) \|_{L^2_+} \leq \|\partial_x g\|_{L^2_+}  + \| u\|_{L^{\infty}(\mathbb{R}; \mathbb{C})}^2 \| g\|_{L^2_+}  ,
\quad
D \left(\mathbb{U}(t) g \right) = \mathbb{U}(t) \mathbf{L}_u (g) + \mathbf{T}_u \mathbf{T}_{\overline{u}} \left( \mathbb{U}(t) g \right) \in L^2_+.
\end{equation}Using formulas \eqref{varphiHs+normHsHs-1},  \eqref{LuSect2}, \eqref{StrContiUtL2} and \eqref{LugnormUtLucomm2for}, we obtain that
\begin{equation}\label{EstiUtgH1norm}
\|\mathbb{U}(t) g\|_{H^1_+} \leq \| g \|_{L^2_+} + \|\partial_x  \left( \mathbb{U}(t) g \right)\|_{L^2_+} \leq \sqrt{2}\left( 1 +2 \|u\|_{L^{\infty}(\mathbb{R}; \mathbb{C})}^2 \right) \|g\|_{H^1_+}, \quad \forall g \in H^1_+, \quad \forall t \in \mathbb{R}.
\end{equation}Fix $g \in H^1_+$, we set $\mathfrak{f}:=\mathbf{L}_u (g) \in L^2_+$. By 
\eqref{StrContiUtL2}, \eqref{BdTuTubar<uZ1Z2esti}  and \eqref{LugnormUtLucomm2for}, we have
\begin{equation}\label{CVforDUtg-gcvto0}
\|D \left( \mathbb{U}(t) g - g \right) \|_{L^2_+} \leq \|   \mathbb{U}(t) \mathfrak{f} - \mathfrak{f}   \|_{L^2_+} + \|u\|_{L^{\infty}(\mathbb{R}; \mathbb{C})}^2 \|  \mathbb{U}(t) g - g   \|_{L^2_+} \to 0,  \quad \text{as} \quad t \to 0.
\end{equation}The strong continuity of $\left( \mathbb{U}(t) \right)_{t \in \mathbb{R}}$ on $H^1_+$ follows from \eqref{StrContiUtL2} and  \eqref{CVforDUtg-gcvto0}. Specifically, we have
\begin{equation}\label{StrCONUtonH1For}
\lim_{t \to 0} \| \mathbb{U}(t) g - g   \|_{H^1_+} = 0, \qquad \forall g \in H^1_+.
\end{equation}For any  $h \in H^2_+$, we have 
$\| \mathbf{L}_u ( h ) \|_{H^1_+} \leq  \|\partial_x h\|_{H^1_+}  +   C_1^2 \| u\|_{\mathcal{Z}^1}^2 \| h\|_{H^1_+}$ by \eqref{BdTuTubar<uZ1Z2esti}.  
It follows that
\begin{equation}\label{estidxUthH1norm}
\|\partial_x \left( \mathbb{U}(t) h \right)\|_{H^1_+} 
\leq  \|  \mathbb{U}(t)  \|_{ \mathcal{E}^1 } \left( \| \partial_x h \|_{H^1_+}  +  2 C_1^2 \| u\|_{\mathcal{Z}^1}^2 \| h\|_{H^1_+}  \right),  \quad \forall t \in \mathbb{R}, 
\end{equation}by \eqref{LugnormUtLucomm2for}.
Formula \eqref{EstiUtgH1norm} yields that $\sup_{t \in \mathbb{R}}\| \mathbb{U}(t)\|_{\mathcal{E}^1} \leq \sqrt{2}\left( 1 + 2 \|u \|_{L^{\infty}(\mathbb{R} ; \mathbb{C})}^2 \right)$. For any $t \in \mathbb{R}$, we have
\begin{equation}\label{EstiUthH2norm}
  \|\mathbb{U}(t) h \|_{H^2_+} \leq \| h \|_{H^1_+} + \|\partial_x  \left( \mathbb{U}(t) h \right)\|_{H^1_+}  
\leq   2 \left( 1 + 2 (1+ C_1)^2 \|u\|_{\mathcal{Z}^1}^2 \right)^2 \|h\|_{H^2_+}, \quad \forall h \in H^2_+,  
\end{equation}by \eqref{estidxUthH1norm}. Fix $h \in H^2_+$, we set $\mathfrak{g}: = \mathbf{L}_u (h) \in H^1_+$. Using \eqref{BdTuTubar<uZ1Z2esti}, 
\eqref{LugnormUtLucomm2for} and  
\eqref{StrCONUtonH1For}, we obtain that
\begin{equation}\label{CVforH1normDUth-h}
\|D \left( \mathbb{U}(t) h - h \right) \|_{H^1_+} \leq \|   \mathbb{U}(t) \mathfrak{g} - \mathfrak{g}   \|_{H^1_+} + C_1^2 \|u\|_{ \mathcal{Z}^1 }^2 \|  \mathbb{U}(t) h - h   \|_{H^1_+} \to 0,  \quad \text{as} \quad t \to 0.
\end{equation}The strong continuity of $\left( \mathbb{U}(t) \right)_{t \in \mathbb{R}}$ on $H^2_+$ is a direct consequence of \eqref{StrCONUtonH1For} and  \eqref{CVforH1normDUth-h}, i.e. 
\begin{equation}\label{StrconvUtonH2}
\lim_{t \to 0} \| \mathbb{U}(t) h - h   \|_{H^2_+} = 0, \qquad \forall h \in H^2_+.
\end{equation}For any $\phi \in H^3_+$, formula \eqref{BdTuTubar<uZ1Z2esti} yields that $\| \mathbf{L}_u (\phi) \|_{H^2_+} \leq \|  \partial_x \phi  \|_{H^2_+} + C_2^2 \| u \|_{\mathcal{Z}^2}^2 \| \phi \|_{H^2_+}$. Then we have 
\begin{equation}\label{estpartialUtphiH2norm}
\|\partial_x \left( \mathbb{U}(t) \phi \right)\|_{H^2_+} 
\leq  \|  \mathbb{U}(t)  \|_{ \mathcal{E}^2 } \left( \| \partial_x \phi \|_{H^2_+}  +  2 C_2^2 \| u\|_{\mathcal{Z}^2}^2 \| \phi \|_{H^2_+}  \right),  \quad \forall t \in \mathbb{R}, 
\end{equation}by \eqref{LugnormUtLucomm2for}.
Recall that 
$\sup_{t \in \mathbb{R}}\| \mathbb{U}(t)\|_{\mathcal{E}^2} \leq  2 \left( 1 + 2 (1+ C_1)^2 \|u\|_{\mathcal{Z}^1}^2 \right)^2 $ by  \eqref{EstiUthH2norm}.  For any $t \in \mathbb{R}$, we have
\begin{equation}\label{estiUtphiH3norm}
  \|\mathbb{U}(t) \phi \|_{H^3_+} \leq \| \phi \|_{H^2_+} + \|\partial_x  \left( \mathbb{U}(t) \phi \right)\|_{H^2_+}  
\leq   2 \sqrt{2} \left( 1 + 2 (1+ C_1 + C_2 )^2 \|u\|_{\mathcal{Z}^2}^3 \right)^2 \| \phi \|_{H^3_+}, 
\end{equation}by \eqref{estpartialUtphiH2norm}.
Fix $\phi \in H^2_+$, we set $\mathfrak{h}: = \mathbf{L}_u (\phi) \in H^2_+$. By formulas \eqref{BdTuTubar<uZ1Z2esti}, 
\eqref{LugnormUtLucomm2for} and  
\eqref{StrconvUtonH2}, we have
\begin{equation}\label{CvforH2normDUtphi-phi}
\|D \left( \mathbb{U}(t) \phi - \phi \right) \|_{H^2_+} \leq \|   \mathbb{U}(t) \mathfrak{h} - \mathfrak{h}   \|_{H^2_+} + C_2^2 \|u\|_{ \mathcal{Z}^2 }^2 \|  \mathbb{U}(t) \phi - \phi \|_{H^2_+} \to 0,  \quad \text{as} \quad t \to 0.
\end{equation}The strong continuity of $\left( \mathbb{U}(t) \right)_{t \in \mathbb{R}}$ on $H^3_+$ is obtained by \eqref{StrconvUtonH2} and \eqref{CvforH2normDUtphi-phi}.
\end{proof}

\begin{defi}\label{AudefTuTbaru}
For any $v \in L^{\infty} (\mathbb{R}; \mathbb{C})$, we define $\mathbf{A}(v)=\mathbf{A}_v := \mathbf{T}_v \mathbf{T}_{\overline{v}} \in \mathcal{B}_{\mathbb{C}}(L^2_+)$. Then the mapping $\mathbf{A} : v \in L^{\infty} (\mathbb{R}; \mathbb{C}) \mapsto \mathbf{A}(v) \in \mathcal{B}_{\mathbb{C}}(L^2_+)$ is real analytic.  Its Fr\'echet differential is given by 
\begin{equation}\label{dAuhFreDiffFor}
\mathrm{d}\mathbf{A}(v) (h) = \mathbf{T}_v \mathbf{T}_{\overline{h}} + \mathbf{T}_h \mathbf{T}_{\overline{v}} \in \mathcal{B}_{\mathbb{C}}(L^2_+) , \qquad  \forall h \in L^{\infty}(\mathbb{R} ; \mathbb{C}).
\end{equation}
\end{defi}

\noindent The next lemma follows directly from formulas \eqref{CommfrakBufrLufor} and \eqref{dAuhFreDiffFor}. It rewrites the commutation relation of the chiral Lax pair in terms of the Fr\'echet differential of the operator-valued map $\mathbf{A}$.

\begin{lem}\label{LaxPairBuLu}
If $u \in \mathcal{Z}^2_+$, let $\mathbf{L}_u , \mathbf{B}_u$ be defined in \eqref{LuSect2} and \eqref{BudefForFoc-}.
Then $\mathbf{B}_u \in \bigcap_{p=0}^1 \mathcal{B}_{\mathbb{C}} (H^p_+)$ and 
\begin{equation}\label{LaxPairBuLuchi}
 - \mathrm{d}\mathbf{A}(u) \left(i\partial_x^2 u + 2 u \Pi_{\geq 0} \partial_x \left( |u|^2 \right) \right) = \left[ \mathbf{B}_u , \mathbf{L}_u \right] \;  : \; H^1_+ \to L^2_+ .
\end{equation}
\end{lem}

 \noindent For later use, we also note the expression of the nonlinear remainder $\mathbf{r}_{\epsilon}(u)$ in the chiral setting. For any $u \in \mathcal{Z}^2_+$ and $\epsilon >0$, let $\mathbf{r}_{\epsilon} (u) \in H^1 (\mathbb{R} ; \mathbb{C})$ be defined in \eqref{DefbfrepuFor} and \eqref{bfrepuexprebytilBu}. Let $\chi_{\epsilon} \in \mathrm{Dom}^{\infty}(\mathbf{G})$ be given by \eqref{chiepFor}. Lemma~\ref{Z1+idealinL2+} yields that $u \chi_{\epsilon} \in H^2_+$. By formulas  \eqref{bfrepuexprebytilBu} and \eqref{bfLufrakLuHarSob},  we have
\begin{equation}\label{bfreputilBuHardZhi}
\mathbf{r}_{\epsilon} (u) = \chi_{\epsilon} \bigl( i \partial_x^2 u + 2 u \Pi_{\geq 0} \partial_x \bigl( |u|^2 \bigr) \bigr) - \widetilde{\mathbf{B}}_u \bigl( u \chi_{\epsilon} \bigr) \in H^1_+ .
\end{equation}
  To apply the unitary propagator lemma (Lemma~\ref{UniPropaLemW}) in the chiral setting, we need the continuity of the operator-valued maps appearing in the Lax pair. The following results establish the required regularity.

\begin{defi}
We define $\mathbf{T}: u \in L^{\infty}(\mathbb{R} ; \mathbb{C}) \mapsto \mathbf{T}(u) := \mathbf{T}_u \in \mathcal{B}_{\mathbb{C}}(L^2_+)$, where $\mathbf{T}_u$ denotes the Toeplitz operator of symbol $u$ given by \eqref{ToepdefForSec2}.
\end{defi}

\noindent Fix an integer $m \geq 1$. For any $p \in \{0,1,2, \cdots, m\}$, the restriction $\mathbf{T}|_{ \mathcal{Z}^m} : u \in \mathcal{Z}^m \mapsto \mathbf{T}_u \in \mathcal{B}_{\mathbb{C}}(H^p_+)$ is $\mathbb{C}$-linear and bounded by the third inequality of \eqref{BounToeplitzBHs}.
The following lemma establishes the continuity of these operator-valued maps.

\begin{lem}\label{LemTcircuPuvcontinu}
Let $a<b$ and $m \ge 1$ be an integer. Assume that $u, v \in \mathcal{C}\left( (a,b) ; \mathcal{Z}^m \right)$. Then $\mathbf{T}\circ u$, $\mathbf{T}\circ v$ and $\mathscr{P}_{u,v} := (\mathbf{T}\circ u)(\mathbf{T}\circ v)$ are continuous mappings from $(a, b)$ to $\mathcal{B}_{\mathbb{C}}(H^p_+)$, for any $p \in \{0,1,2, \cdots, m\}$. Here the product is taken pointwise:
\begin{equation}\label{scrPu,vdefFor}
\mathscr{P}_{u,v}(t) = \left( (\mathbf{T}\circ u)(\mathbf{T}\circ v) \right)(t) = \mathbf{T}_{u(t)} \mathbf{T}_{v(t)} \in \bigcap_{p=0}^m \mathcal{B}_{\mathbb{C}}(H^p_+), \quad \forall t \in (a, b).
\end{equation}
\end{lem}

\begin{lem}\label{quvttoquvtLemContinu}
Assume that $u \in \mathcal{C}\left( (a,b) ; L^{\infty}(\mathbb{R}; \mathbb{C}) \right)$ and $v \in \mathcal{C}\left( (a,b) ; \mathcal{Z}^1 \right)$, for some $a, b \in \mathbb{R}$ with $a<b$. For any $t \in (a,b)$, set $\mathfrak{q}_{u,v}(t) := \mathbf{T}_{u(t)} \partial_x \mathbf{T}_{\overline{v(t)}} : H^1_+ \to L^2_+$. Then the mapping $\mathfrak{q}_{u,v} : t \mapsto \mathfrak{q}_{u,v}(t)$ is continuous from $(a, b)$ to $\mathcal{B}_{\mathbb{C}}(H^1_+ ; L^2_+)$.
\end{lem}
\begin{proof}
Formula  \eqref{BounToeplitzBHs} implies $\mathbf{T}\circ u \in \mathcal{C}\left( (a,b) ; \mathcal{B}_{\mathbb{C}}(L^2_+) \right)$. By Lemma~\ref{LemTcircuPuvcontinu}, $\mathbf{T}\circ \overline{v}$ is continuous from $(a, b)$ to  $\mathcal{B}_{\mathbb{C}}(H^1_+)$. Since $\partial_x \in \mathcal{B}_{\mathbb{C}}(H^1_+ ; L^2_+)$, we have  $
\mathfrak{q}_{u,v} = (\mathbf{T}\circ u) \partial_x (\mathbf{T}\circ \overline{v}) \in \mathcal{C}\left( (a,b); \mathcal{B}_{\mathbb{C}}(H^1_+ ; L^2_+) \right)$.
\end{proof}

\begin{cor}\label{tildButtotilButcontiCor}
Assume that $u \in \mathcal{C}\left( (a,b) ; \mathcal{Z}^1 \right)$ for some $a<b$. By formula  \eqref{tildeBudefFoc-}, we have 
\begin{equation}\label{tilbfBuCont}
\widetilde{\mathbf{B}}_{u(t)} = i\partial_x^2 + 2 \mathbf{T}_{u(t)} \partial_x \mathbf{T}_{\overline{u(t)}} \in \mathcal{B}_{\mathbb{C}}(H^2_+; L^2_+), \qquad \forall  t \in (a, b).
\end{equation} 
  Then $\widetilde{\mathbf{B}}_u : t \mapsto \widetilde{\mathbf{B}}_{u(t)}$ is a continuous mapping from $(a, b)$ to  $\mathcal{B}_{\mathbb{C}}(H^2_+; L^2_+)$.
\end{cor}
\begin{proof}
Since $i\partial_x^2 \in \mathcal{B}_{\mathbb{C}}(H^2_+ ; L^2_+)$ and $\widetilde{\mathbf{B}}_u = i\partial_x^2 + 2 \mathfrak{q}_{u,u}$, the result follows from Lemma~\ref{quvttoquvtLemContinu} and the continuity of the embedding $\mathcal{B}_{\mathbb{C}}(H^1_+; L^2_+) \hookrightarrow \mathcal{B}_{\mathbb{C}}(H^2_+; L^2_+)$.
\end{proof}

 \subsection{Conjugate Shift Semigroup and its Perturbation}
 
 \noindent The Lax--Beurling shift semigroups $(\mathbf{S}(\eta))_{\eta \geq 0}$ and $(\mathbf{S}(\eta)^*)_{\eta \geq 0}$ have been  studied in Subsection~\ref{SubsShifSemiGr}.  
  We now introduce the conjugate shift semigroup $(e^{-it\partial_x^2}\mathbf{S}(\eta)^* e^{it\partial_x^2})_{\eta \geq 0}$, whose infinitesimal generator allows us to interpret $-i(\mathbf{G} + 2t\mathbf{L}_u)$ as a perturbation of that generator. \\
  
\noindent Let $(\mathbf{S}(\eta)^*)_{\eta \geq 0}$ be defined in formula  \eqref{LaxBeurlingSG}. For any  $t \in \mathbb{R}$ and $\eta \geq 0$, we define
\begin{equation}\label{frakStetaDef}
\mathfrak{S}_t(\eta):= e^{-it\partial_x^2}\mathbf{S}(\eta)^* e^{it\partial_x^2} \in \mathcal{B}_{\mathbb{C}}(L^2_+)  .
\end{equation}
For any $s \geq 0$,   the Schr\"odinger group $(e^{it\partial_x^2})_{t \in \mathbb{R}}$ is a strongly continuous unitary group on the Hardy--Sobolev space $H^s_+$. Then the following lemma is a direct consequence of \eqref{SchrGrUnitHs} and  Lemma $\ref{ContSemiGrLem}$.
 
\begin{lem}\label{LemContrSefrakSt}
For any $t \in \mathbb{R}$ and $s \geq 0$, the family $(\mathfrak{S}_t(\eta))_{\eta \geq 0}$  forms a contraction semigroup  on $H^s_+$.
\end{lem}
 
\noindent For any $s \geq 0$, let $\mathrm{Dom}^s (\mathbf{G})$ be defined by \eqref{DomsGDefSect2}.  The operator $\bigl( \mathbf{G} , \mathrm{Dom}(\mathbf{G}) \bigr)$ is given by Definition $\ref{DombfGdefSect2}$ and formula  \eqref{DefGphiSect2}. 
 Then 
 $\mathrm{Dom}(\mathbf{G}) = \mathrm{Dom}^0 (\mathbf{G}) $. 
 For any $t \in \mathbb{R}$, we set 
\begin{equation}\label{DomsbbGtdefFor}
\mathrm{Dom}^s (\mathbb{G}_t): = e^{-it \partial_x^2} \bigl(  \mathrm{Dom}^s (\mathbf{G}) \bigr) = \{ \varphi \in H^s_+ : e^{it \partial_x^2} \varphi \in \mathrm{Dom}  (\mathbf{G}) \quad \text{and} \quad \mathbf{G} e^{it \partial_x^2} \varphi \in H^s_+ \} .
\end{equation}

\begin{defi}\label{DomGtbbGtdef}
For any $t\in \mathbb{R}$, we define  $\mathrm{Dom}(\mathbb{G}_t) : = \mathrm{Dom}^0 (\mathbb{G}_t)  = e^{-it \partial_x^2} \bigl(  \mathrm{Dom}  (\mathbf{G}) \bigr) \subset L^2_+$ and
\begin{equation}\label{ForbbGtDef}
\mathbb{G}_t (\varphi) := e^{- i t \partial_x^2} \mathbf{G} e^{ i t \partial_x^2} (\varphi) \in L^2_+ , \qquad \forall \varphi \in \mathrm{Dom}(\mathbb{G}_t) .
\end{equation}
\end{defi}
  
\begin{prop}\label{PropinGenefrakSt}
For any $t \in \mathbb{R}$ and $s \geq 0$, the operator $\bigl( -i \mathbb{G}_t, \mathrm{Dom}^s (\mathbb{G}_t) \bigr)$ is the infinitesimal generator of the contraction semigroup $(\mathfrak{S}_t(\eta))_{\eta \geq 0}$ on $H^s_+$. Specifically, we have
\begin{equation}\label{CvdefDomsbbGt}
\mathrm{Dom}^s (\mathbb{G}_t)  = \{ \phi \in H^s_+ : \exists \mathfrak{v} \in H^s_+ \quad \text{such that} \quad \lim_{\eta \to 0^+} \| \tfrac{1}{\eta} \bigl( \mathfrak{S}_t(\eta) \phi -\phi \bigr) - \mathfrak{v}\|_{H^s_+} = 0\}.
\end{equation}Moreover, if  $\varphi \in \mathrm{Dom}^s (\mathbb{G}_t)$, then $\tfrac{1}{\eta} \bigl( \mathfrak{S}_t(\eta) \varphi -\varphi \bigr) \to - i \mathbb{G}_t (\varphi) $ in $H^s_+$, as $\eta \to 0^+$. Consequently, the operator  $\bigl( i \mathbb{G}_t, \mathrm{Dom}^s (\mathbb{G}_t) \bigr)$ is maximal accretive on $H^s_+$. In other words, the following assertions hold.
\begin{enumerate}
\item For any $\varphi \in \mathrm{Dom}^s (\mathbb{G}_t)$,  we have $ \operatorname{Re}  \langle i  \mathbb{G}_t (\varphi) , \varphi \rangle_{H^s_+} \geq 0$.
\item For any $z \in \mathbb{C}_+$, the operator $\mathbb{G}_t - z : \mathrm{Dom}^s (\mathbb{G}_t) \to H^s_+$ is bijective and we have
\begin{equation}\label{estiResolbbGt-z}
\sup_{s \geq 0} \big\| (\mathbb{G}_t - z )^{-1} \big\|_{\mathcal{B}_{\mathbb{C}} (H^s_+) } \leq (\operatorname{Im} z)^{-1} .
\end{equation} 
\end{enumerate}
\end{prop}
\begin{proof}
Fix $s \geq 0$.  By Proposition $\ref{GisInfiniGeneHsProp}$, the operator $\left(-i \mathbf{G}, \mathrm{Dom}^s(\mathbf{G}) \right)$ generates the shift semigroup $(\mathbf{S}(\eta)^*)_{\eta \geq 0}$ on $H^s_+$. 
  For any $t \in \mathbb{R}$, the domain characterization \eqref{CvdefDomsbbGt} follows from the corresponding characterization for $\mathbf{G}$ in \eqref{SemiGrCVdefDomsG} together with \eqref{frakStetaDef}, \eqref{DomsbbGtdefFor}, \eqref{ForbbGtDef}, and the unitary property \eqref{SchrGrUnitHs} of $e^{it\partial_x^2}$ on $H^s_+$. 
 So $\bigl( -i \mathbb{G}_t, \mathrm{Dom}^s (\mathbb{G}_t) \bigr)$  generates the contraction semigroup $(\mathfrak{S}_t(\eta))_{\eta \geq 0}$ on $H^s_+$, i.e. 
\begin{equation}\label{frakStinBHsexponen}
\mathfrak{S}_t(\eta) = \exp\bigl( - i \eta \mathbb{G}_t \bigr) \in \bigcap_{ s \geq 0} \mathcal{B}_{\mathbb{C}} (H^s_+), \qquad \forall \eta \geq 0.
\end{equation}Theorem \uppercase\expandafter{\romannumeral 10}. 48 of   \cite{Reed2SimFouSAop} yields that  $\bigl( i \mathbb{G}_t, \mathrm{Dom}^s (\mathbb{G}_t) \bigr)$ is maximal accretive on $H^s_+$. Fix $z \in \mathbb{C}_+$. Then
\begin{equation}\label{bbGtminuzconju}
\mathbb{G}_t - z = e^{-it \partial_x^2} (\mathbf{G} - z)
 e^{it \partial_x^2} : \mathrm{Dom}^s (\mathbb{G}_t)
 = e^{-it \partial_x^2} \bigl( \mathrm{Dom}^s (\mathbf{G}) \bigr) \to H^s_+
\end{equation}is bijective by \eqref{SchrGrUnitHs} and \eqref{ForbbGtDef}.  
  Finally, the resolvent estimate \eqref{estiResolbbGt-z} is  a direct consequence of formulas  \eqref{bbGtminuzconju} and \eqref{EstG-z-1BHs}.
\end{proof}

\noindent The following lemma shows that a bounded skew-adjoint perturbation of $\bigl( i\mathbb{G}_t, \mathrm{Dom} (\mathbb{G}_t) \bigr)$ remains maximal accretive on the Hardy space $L^2_+$, with the same uniform resolvent bound.

\begin{lem}\label{LemPertGtMaxAcc}
For any $t, \mathrm{r} \in \mathbb{R}$ and $u \in L^{\infty}(\mathbb{R};\mathbb{C})$, we define 
\begin{equation}\label{DombbGtpertrTuTub}
\mathrm{Dom} (\mathbb{G}_t + 2 \mathrm{r} \mathbf{T}_u \mathbf{T}_{\overline{u}} ) : =\mathrm{Dom}  (\mathbb{G}_t)
 = e^{-it \partial_x^2} \bigl( \mathrm{Dom}  (\mathbf{G}) \bigr).
\end{equation}
 Then  $\bigl(i\mathbb{G}_t + 2 i \mathrm{r} \mathbf{T}_u \mathbf{T}_{\overline{u}}  , \mathrm{Dom} (\mathbb{G}_t + 2 \mathrm{r} \mathbf{T}_u \mathbf{T}_{\overline{u}} ) \bigr)$ is maximal accretive on $L^2_+$. In particular, for any $z \in \mathbb{C}_+$, the operator $\mathbb{G}_t + 2 \mathrm{r} \mathbf{T}_u \mathbf{T}_{\overline{u}} - z  : e^{-it \partial_x^2} \bigl( \mathrm{Dom}  (\mathbf{G})  \bigr) \to L^2_+$ is bijective. Moreover, we have
 \begin{equation}\label{EstuniresolGtpluTuTub}
\sup_{t, r \in \mathbb{R}}\sup_{u \in L^{\infty} (\mathbb{R} ; \mathbb{C})} \big\|\bigl( \mathbb{G}_t + 2 \mathrm{r} \mathbf{T}_u \mathbf{T}_{\overline{u}} -z \bigr)^{-1} \big\|_{\mathcal{B}_{\mathbb{C}}(L^2_+)} \leq (\operatorname{Im} z)^{-1} .
 \end{equation}
\end{lem} 
 
\begin{proof}
Fix $u \in L^{\infty}(\mathbb{R};\mathbb{C})$, then $\mathbf{A}_u = \mathbf{T}_u \mathbf{T}_{\overline{u}}$ is a bounded self-adjoint operator on $L^2_+$. Fix $t, \mathrm{r} \in \mathbb{R}$. We have
$\operatorname{Re} \langle 2 i \mathrm{r} \mathbf{T}_u \mathbf{T}_{\overline{u}} (f), f \rangle_{L^2_+} = 0$ for all $f \in L^2_+$. The Kato--Rellich theorem for accretive perturbations (Lemma and Theorem~X.50 on p.~244 of \cite{Reed2SimFouSAop}) implies that $i\mathbb{G}_t + 2 i \mathrm{r} \mathbf{T}_u \mathbf{T}_{\overline{u}} : e^{-it \partial_x^2} \bigl( \mathrm{Dom} (\mathbf{G}) \bigr) \to L^2_+$ is maximal accretive.
Then \eqref{EstuniresolGtpluTuTub} follows from Proposition~7.1 of \cite{BrezisFA} and Theorem~X.48 of \cite{Reed2SimFouSAop}.
\end{proof} 
  
\noindent By conjugating the previous lemma with the Schr\"odinger group $(e^{it\partial_x^2})_{ t \in \mathbb{R}}$, we immediately obtain a corresponding statement for the original generator $\bigl( - i  \mathbf{G}, \mathrm{Dom} (\mathbf{G}) \bigr)$.
   
\begin{lem}\label{ttZutdefLem}
 For any $t, \mathrm{r} \in \mathbb{R}$ and $u \in L^{\infty} (\mathbb{R} ; \mathbb{C})$, we define that
\begin{equation}\label{dombfGpluttZut}
\mho_u^{(t)}: =  e^{it \partial_x^2} \mathbf{T}_u \mathbf{T}_{\overline{u}}  e^{- it \partial_x^2} \in \mathcal{B}_{\mathbb{C}}(L^2_+)  
\quad  \text{and}    \quad
\mathrm{Dom} (\mathbf{G}+ 2 \mathrm{r} \mho_u^{(t)} ) : =\mathrm{Dom}  (\mathbf{G}) .
\end{equation}Then $\mho_u^{(t)}$ is self-adjoint on $L^2_+$ and $ \bigl(  i \mathbf{G}+ 2 i \mathrm{r} \mho_u^{(t)} , \mathrm{Dom} (\mathbf{G}+ 2 \mathrm{r} \mho_u^{(t)} ) \bigr)$ is maximal accretive on $L^2_+$. For any $z \in \mathbb{C}_+$, the operator $\mathbf{G}+ 2  \mathrm{r} \mho_u^{(t)}  - z  :   \mathrm{Dom}  (\mathbf{G})  \to L^2_+$ is bijective. Moreover, we have
 \begin{equation}\label{resolGpluZutEstuni}
\sup_{t, r \in \mathbb{R}}\sup_{u \in L^{\infty} (\mathbb{R} ; \mathbb{C})} \big\|\bigl( \mathbf{G}+ 2  \mathrm{r} \mho_u^{(t)}  - z \bigr)^{-1} \big\|_{\mathcal{B}_{\mathbb{C}}(L^2_+)} \leq (\operatorname{Im} z)^{-1} .
 \end{equation}
\end{lem}

\begin{proof}
The proof follows from  
Lemma $\ref{LemPertGtMaxAcc}$,  the unitary property \eqref{SchrGrUnitHs} and the following identity:
\begin{equation}\label{CjRebbGtplusTutub}
\mathbb{G}_t + 2 \mathrm{r} \mathbf{T}_u \mathbf{T}_{\overline{u}} - a = 
e^{ - it\partial_x^2} \bigl( \mathbf{G}+ 2  \mathrm{r} \mho_u^{(t)}  - a \bigr)e^{it\partial_x^2} : e^{-it \partial_x^2} \bigl( \mathrm{Dom}  (\mathbf{G}) \bigr) \to L^2_+,
\end{equation}for any $u \in L^{\infty}(\mathbb{R} ; \mathbb{C})$, $a \in \mathbb{C}$ and $t ,\mathrm{r} \in \mathbb{R}$.
\end{proof}

\noindent The next lemma describes the intersection of the domains $\mathrm{Dom}^s(\mathbf{G})$ under the group 
$( e^{i t \partial_x^2 })_{t \in \mathbb{R}}$ with $s \geq 0$. It shows in particular that a function belongs to this intersection if and only if it gains one extra order of Sobolev regularity.

\begin{lem}\label{lembigcpDomsG}
Fix $s \geq 0$. For any $\mathtt{t}_1, \mathtt{t}_2 \in \mathbb{R}$ such that $\mathtt{t}_1 \neq \mathtt{t}_2$, the following identity holds: 
\begin{equation}\label{bigcapeitLapDomsG}
\begin{split}
\bigcap_{\mathrm{r} \in \mathbb{R}} e^{i \mathrm{r} \partial_x^2 } \bigl( \mathrm{Dom}^s (\mathbf{G}) \bigr) 
= & e^{i \mathtt{t}_1 \partial_x^2 } \bigl( \mathrm{Dom}^s (\mathbf{G}) \bigr) \cap e^{i \mathtt{t}_2 \partial_x^2 } \bigl( \mathrm{Dom}^s (\mathbf{G}) \bigr) 
=  \mathrm{Dom}^s (\mathbf{G}) \cap H^{s+1}_+  \\
= &  H^{s+1}_+   \cap  e^{ i \mathtt{t}_1 \partial_x^2} \bigl(\mathrm{Dom}^s (\mathbf{G})\bigr) 
= e^{ i \mathtt{t}_1 \partial_x^2} \bigl(\mathrm{Dom}^s (\mathbf{G}) \cap H^{s+1}_+ \bigr) .
\end{split}
\end{equation}
Consequently, the restriction $e^{i \mathtt{t}_1 \partial_x^2 }|_{\mathrm{Dom}^s (\mathbf{G}) \cap H^{s+1}_+} : \mathrm{Dom}^s (\mathbf{G}) \cap H^{s+1}_+ \to \mathrm{Dom}^s (\mathbf{G}) \cap H^{s+1}_+$ is bijective.
\end{lem}

\begin{proof}
Fix $s \geq 0$ and $t \in \mathbb{R}$. For any $\varphi \in   e^{-it \partial_x^2} \bigl( \mathrm{Dom}^s (\mathbf{G}) \bigr)$, we define 
\begin{equation}\label{bfhtvphidef}
\mathbf{h}_{t , \varphi} : =e^{it \partial_x^2} \varphi \in \mathrm{Dom}^s (\mathbf{G}).
\end{equation} 
 Formula  \eqref{DomsGDefSect2} implies that 
$\widehat{\mathbf{h}}_{t , \varphi} |_{\mathbb{R}_+^*} \in H^1_s (\mathbb{R}_+^*)$ with $\mathbb{R}_+^* = (0, +\infty)$. Then $\widehat{\varphi}|_{\mathbb{R}_+^*} \in H^1_{\mathrm{loc}} (\mathbb{R}_+^*)$ and 
\begin{equation}\label{dxihtvphibfhtvarp}
e^{- i t \xi^2} \partial_{\xi} \widehat{\varphi}(\xi) 
= 2 i t \xi \widehat{\mathbf{h}}_{t , \varphi} (\xi) + \partial_{\xi} \widehat{\mathbf{h}}_{t , \varphi} (\xi), \qquad \forall \xi >0, 
\end{equation}by  Corollary 8.10 of \cite{BrezisFA}. The rest of the proof is divided into two steps. \\
  
\noindent \textbf{Step \romannumeral 1}. Assume that  $\varphi \in   H^{s+1}_+   \cap  e^{ - i  t \partial_x^2} \bigl(\mathrm{Dom}^s (\mathbf{G})\bigr) $ for some $s \geq 0$ and $t \in \mathbb{R}$.  Formula \eqref{SchrGrUnitHs} yields that
$\mathbf{h}_{t , \varphi}   =e^{it \partial_x^2} \varphi \in \mathrm{Dom}^s (\mathbf{G}) \cap H^{s+1}_+$. 
 Using  \eqref{dxihtvphibfhtvarp},  \eqref{Lpmussp},  \eqref{DefGphiSect2} and \eqref{DomsGDefSect2}, we obtain that
\begin{equation*}
\| \partial_{\xi}  \widehat{\varphi} \|_{L^2_s (\mathbb{R}_+^*)}^2 \leq 16 \pi t^2  \| \mathbf{h}_{t , \varphi} \|_{H^{s+1}_+}^2  + 4 \pi \| \mathbf{G} \mathbf{h}_{t , \varphi} \|_{H^s_+}^2 <+ \infty.
\end{equation*}Then $\widehat{\varphi}|_{\mathbb{R}_+^*} \in H^1_s (\mathbb{R}_+^*)$
by \eqref{H1snormWeiSobo}. Formula \eqref{DomsGDefSect2} implies that $\varphi \in \mathrm{Dom}^s (\mathbf{G})$. Therefore, we have 
\begin{equation}\label{eitLapDsGandHsset}
e^{ - i t \partial_x^2} \bigl( \mathrm{Dom}^s (\mathbf{G}) \cap H^{s+1}_+   \bigr) 
\subset   H^{s+1}_+  \cap  e^{ - i t \partial_x^2} \bigl( \mathrm{Dom}^s  (\mathbf{G}) \bigr)    
 \subset \mathrm{Dom}^s (\mathbf{G}) \cap H^{s+1}_+  , \qquad \forall t \in \mathbb{R}.
\end{equation}
Using the group property  \eqref{GrouPropexpitdelta} and substituting $t$ for $-t$ in formula  \eqref{eitLapDsGandHsset}, we obtain that
\begin{equation}\label{eitDsandHsbij}
e^{   i t \partial_x^2} \bigl( \mathrm{Dom}^s (\mathbf{G}) \cap H^{s+1}_+   \bigr) 
 =   H^{s+1}_+  \cap  e^{  i t \partial_x^2} \bigl( \mathrm{Dom}^s  (\mathbf{G}) \bigr)    
 = \mathrm{Dom}^s (\mathbf{G}) \cap H^{s+1}_+, \qquad \forall t \in \mathbb{R} .
\end{equation}So $e^{i t  \partial_x^2 }|_{\mathrm{Dom}^s (\mathbf{G}) \cap H^{s+1}_+} : \mathrm{Dom}^s (\mathbf{G}) \cap H^{s+1}_+ \to \mathrm{Dom}^s (\mathbf{G}) \cap H^{s+1}_+$ is bijective. Moreover, formula \eqref{eitDsandHsbij} implies that $\mathrm{Dom}^s (\mathbf{G}) \cap H^{s+1}_+   = e^{ i \mathrm{r} \partial_x^2} \bigl( \mathrm{Dom}^s (\mathbf{G}) \cap H^{s+1}_+ \bigr) \subset e^{ i \mathrm{r} \partial_x^2} \bigl( \mathrm{Dom}^s (\mathbf{G})  \bigr)$,  $\forall \mathrm{r} \in  \mathbb{R}$. Hence we have
\begin{equation}\label{bceirLaDsGandHs}
 \mathrm{Dom}^s (\mathbf{G}) \cap H^{s+1}_+  
  \subset \bigcap_{\mathrm{r} \in \mathbb{R}} e^{i \mathrm{r} \partial_x^2 } \bigl( \mathrm{Dom}^s (\mathbf{G}) \bigr) .
\end{equation} 

\noindent \textbf{Step \romannumeral 2}. Fix $\mathtt{t}_1, \mathtt{t}_2 \in \mathbb{R}$ such that $\mathbf{t}: = \mathtt{t}_1 - \mathtt{t}_2 \ne 0$. 
If $\phi \in e^{i \mathtt{t}_1 \partial_x^2 } \bigl( \mathrm{Dom}^s (\mathbf{G}) \bigr) \cap e^{i \mathtt{t}_2 \partial_x^2 } \bigl( \mathrm{Dom}^s (\mathbf{G}) \bigr)$, we set 
\begin{equation}\label{bffiseitLapnehi}
\mathbf{f} : = \mathbf{h}_{ - \mathtt{t}_1 , \phi} = e^{- i \mathtt{t}_1 \partial_x^2} (\phi) \in \mathrm{Dom}^s (\mathbf{G})  \quad \text{and} \quad \mathbf{g} : = \mathbf{h}_{ - \mathtt{t}_2 , \phi}   
= e^{ i \mathbf{t} \partial_x^2} (\mathbf{f}) = \mathbf{h}_{\mathbf{t} , \mathbf{f}} \in \mathrm{Dom}^s (\mathbf{G}). 
\end{equation}Formula \eqref{dxihtvphibfhtvarp} implies that $e^{ - i \mathbf{t} \xi^2} \partial_{\xi} \widehat{\mathbf{f}}(\xi) = 2 i \mathbf{t} \xi  \widehat{\mathbf{g}}(\xi) + \partial_{\xi} \widehat{\mathbf{g}}(\xi)$, $\forall \xi >0$. Then we have
\begin{equation}\label{normbfgHs+1}
\big\| \mathbf{g} \big\|_{H^{s+1}_+}^2 
 - \big\| \mathbf{g} \big\|_{H^s_+}^2  = \tfrac{1}{2 \pi} \int_0^{+\infty} |\widehat{\mathbf{g}} ( \xi ) |^2 \xi^2 (1 + \xi^2)^s \mathrm{d}\xi
\leq   \tfrac{1}{2 \mathbf{t}^2} \bigl( \big\| \mathbf{G}\mathbf{f} \big\|_{H^s_+}^2 + \big\| \mathbf{G} \mathbf{g} \big\|_{H^s_+}^2  \bigr) <+\infty , 
\end{equation}by \eqref{varphiHs+normHsHs-1},   \eqref{H1snormWeiSobo}
 and   \eqref{DomsGDefSect2}. So $\phi = e^{ i \mathtt{t}_2 \partial_x^2 } (\mathbf{g}) \in H^{s+1}_+$ by \eqref{SchrGrUnitHs} and \eqref{normbfgHs+1}. Consequently, we have
 \begin{equation}\label{bcapeitLapDsGsubs}
 \bigcap_{\mathrm{r} \in \mathbb{R}} e^{i \mathrm{r}  \partial_x^2 } \bigl( \mathrm{Dom}^s (\mathbf{G}) \bigr) 
\subset  e^{i \mathtt{t}_1 \partial_x^2 } \bigl( \mathrm{Dom}^s (\mathbf{G}) \bigr) \cap e^{i \mathtt{t}_2 \partial_x^2 } \bigl( \mathrm{Dom}^s (\mathbf{G}) \bigr) 
\subset   H^{s+1}_+   \cap  e^{ i \mathtt{t}_1 \partial_x^2} \bigl(\mathrm{Dom}^s (\mathbf{G})\bigr). 
 \end{equation} 
  Finally, formula \eqref{bigcapeitLapDomsG} follows directly from  \eqref{eitDsandHsbij}, \eqref{bceirLaDsGandHs} and  \eqref{bcapeitLapDsGsubs}.
\end{proof}

\begin{cor}\label{CorDinfGinvbySchGr}
The subspace $\mathrm{Dom}^{\infty}(\mathbf{G})$ given by \eqref{DomrGsubsetDomsG} is invariant under the Schr\"odinger group $(e^{it \partial_x^2})_{t \in \mathbb{R}}$. Moreover, for any $s \geq 0$ and $t \in \mathbb{R}$, we have 
\begin{equation}\label{DinftybbGt}
\bigcap_{\alpha \geq 0} \mathrm{Dom}^{\alpha} (\mathbb{G}_t) = \mathrm{Dom}^{\infty} (\mathbf{G})  \subset \mathrm{Dom}^{s+1} (\mathbf{G}) \subset \mathrm{Dom}^s (\mathbf{G}) \cap H^{s+1}_+  = \bigcap_{ \mathrm{r} \in \mathbb{R}} \mathrm{Dom}^s (\mathbb{G}_{\mathrm{r}}) .
\end{equation}
\end{cor}

\begin{proof}
If $t \in \mathbb{R}$ and $\phi \in \mathrm{Dom}^{\infty} (\mathbf{G})$, we have $e^{ i t \partial_x^2} \phi \in e^{ i t \partial_x^2} \bigl( \mathrm{Dom}^s (\mathbf{G}) \cap H^{s+1}_+  \bigr) = \mathrm{Dom}^s (\mathbf{G}) \cap H^{s+1}_+$, for any $s \geq 0$,  by \eqref{DomrGsubsetDomsG} and \eqref{bigcapeitLapDomsG}.
 Then $e^{ i t \partial_x^2} \phi \in \mathrm{Dom}^{\infty} (\mathbf{G})$, $\forall t \in \mathbb{R}$. Therefore, we have
\begin{equation}\label{DinGinvSchGrfor}
e^{ i t \partial_x^2} \bigl( \mathrm{Dom}^{\infty} (\mathbf{G}) \bigr) = \bigl( \mathrm{Dom}^{\infty} (\mathbf{G}) \bigr), \qquad \forall t \in \mathbb{R}.
\end{equation}
 Formula \eqref{DinftybbGt} is a direct consequence of \eqref{DomsGDefSect2},   \eqref{DomrGsubsetDomsG},  \eqref{DomsbbGtdefFor},  \eqref{bigcapeitLapDomsG} and \eqref{DinGinvSchGrfor}.
\end{proof} 

\noindent The following lemma gives an explicit expression for $\mathbb{G}_t - \mathbf{G}$ when restricted to $\bigcap_{\mathrm{r} \in \mathbb{R}} \mathrm{Dom} (\mathbb{G}_{\mathrm{r}})$. 
\begin{lem}
For any $t \in \mathbb{R}$ and  $\varphi \in \bigcap_{\mathrm{r} \in \mathbb{R}} \mathrm{Dom} (\mathbb{G}_{\mathrm{r}}) =  \mathrm{Dom}(\mathbf{G}) \cap H^1_+ $, the following identity holds:
\begin{equation}\label{bbGtissumtwoGtD}
\mathbb{G}_t (\varphi) = e^{- i t \partial_x^2} \mathbf{G} e^{ i t \partial_x^2} (\varphi) = \mathbf{G} \varphi + 2 t D \varphi \in L^2_+. 
\end{equation}
\end{lem}
 
\begin{proof}
 Formula \eqref{bbGtissumtwoGtD}  follows directly from  \eqref{dxihtvphibfhtvarp}, \eqref{bfhtvphidef},  \eqref{ForbbGtDef} and 
\eqref{DefGphiSect2}.
\end{proof}

\noindent Combining   formula \eqref{bbGtissumtwoGtD}   with the definition of the chiral Lax operator $\mathbf{L}_u : H^1_+ \to L^2_+$  given by \eqref{LuSect2}, we can now express $\mathbf{G} + 2t\mathbf{L}_u$ as a bounded self-adjoint perturbation of $\mathbb{G}_t$.

\begin{lem}\label{lem:sumbfG2tLubbGt}
For any $t \in \mathbb{R}$ and $u \in L^{\infty}(\mathbb{R};\mathbb{C})$,   we define the common domain
\begin{equation}\label{DomsumG2Lu}
\mathrm{Dom} (\mathbf{G} + 2 t \mathbf{L}_u) = \mathrm{Dom} (\mathbf{G} + 2 t D) := \bigcap_{\mathrm{r} \in \mathbb{R}} \mathrm{Dom} (\mathbb{G}_{\mathrm{r}}) = \mathrm{Dom}(\mathbf{G}) \cap H^1_+ .
\end{equation}
Then for any $\varphi \in \mathrm{Dom} (\mathbf{G} + 2 t \mathbf{L}_u)  \subset \mathrm{Dom} (\mathbb{G}_t - 2 t \mathbf{T}_u \mathbf{T}_{\overline{u}} )$, the following identity holds:
\begin{equation}\label{sumbfG2tLubbGt}
\bigl( \mathbf{G} + 2 t \mathbf{L}_u \bigr) \varphi 
= \mathbf{G} \varphi + 2 t D \varphi - 2 t \mathbf{T}_u \mathbf{T}_{\overline{u}} (\varphi) 
= \bigl( \mathbb{G}_t - 2 t \mathbf{T}_u \mathbf{T}_{\overline{u}} \bigr) \varphi \in L^2_+ .
\end{equation}
Consequently, for any $z \in \mathbb{C}_+$, the operator $\mathbf{G} + 2 t \mathbf{L}_u - z : \mathrm{Dom}(\mathbf{G}) \cap H^1_+ \to L^2_+$ is injective.
\end{lem}

\begin{proof}
Formula \eqref{sumbfG2tLubbGt} follows directly from \eqref{bbGtissumtwoGtD} and \eqref{LuSect2}. 
 The injectivity of $\mathbf{G} + 2 t \mathbf{L}_u - z$ is then a consequence of \eqref{sumbfG2tLubbGt} and Lemma~\ref{LemPertGtMaxAcc}.
\end{proof}

 \noindent The next lemma shows that the operator $\bigl( \mathbf{I}_+  , \mathrm{Dom}(\mathbf{I}_+)\bigr)$ given by \eqref{DomI+} and \eqref{DefI+} is invariant under the Schr\"odinger group $\bigl( e^{it \partial_x^2 } \bigr)_{t \in \mathbb{R}}$. For any $\epsilon >0$, let $\chi_{\epsilon} \in \mathrm{Dom}^{\infty}(\mathbf{G})$ be defined in \eqref{chiepFor}. 
  
\begin{lem}\label{LemInvbfISchGr}
For any $t \in \mathbb{R}$, we have $\mathrm{Dom}(\mathbb{G}_t) \subset \mathrm{Dom}(\mathbf{I}_+)$. For any $\varphi \in \mathrm{Dom}(\mathbb{G}_t)$, we have
\begin{equation}\label{bfIeitLapisbfI}
\mathbf{I}_+( \varphi ) = \mathbf{I}_+\bigl( e^{it \partial_x^2}( \varphi ) \bigr) = \lim_{\epsilon \to 0^+} \langle \varphi, \chi_{\epsilon} \rangle_{L^2_+}.
\end{equation}
\end{lem}
 
\begin{proof}
Fix  $t \in \mathbb{R}$ and $\varphi \in \mathrm{Dom}(\mathbb{G}_t)$. Set $g : = e^{ i t \partial_x^2} \varphi \in \mathrm{Dom}(\mathbf{G})$. By Definition $\ref{DombfGdefSect2}$ and Lemma $\ref{ExtOpLem}$, there exists $\mathfrak{h} \in H^1 (\mathbb{R} ; \mathbb{C}) \hookrightarrow \mathcal{C} \cap L^{\infty}  (\mathbb{R} ; \mathbb{C})$ such that $\mathfrak{h} (\xi ) = \widehat{g}(\xi) $, $\forall \xi >0$. Then we have
\begin{equation*}
\widehat{ \varphi }(\xi) = e^{-it \xi^2} \widehat{g}(\xi) = e^{-it \xi^2} \mathfrak{h}(\xi) \to \mathfrak{h}(0) = \widehat{g}(0^+) , 
\end{equation*}as $\xi \to 0^+$.  So $\mathbf{I}_+( \varphi ) = \mathbf{I}_+ (g) = \mathbf{I}_+\bigl( e^{it \partial_x^2}( \varphi ) \bigr)$. Since $e^{ i t \partial_x^2}$ is a unitary operator on $L^2_+$ by formula \eqref{SchrGrUnitHs},   the second identity of \eqref{bfIeitLapisbfI}  is a direct consequence of  Lemma  $\ref{phichiepCVtoI0Lem}$ and  
 Lemma $\ref{SchrGronchieplem}$. 
\end{proof}

\noindent The following lemma extends the linear functional $\mathbf{I}_+ \bigl( \mathbb{G}_t + 2 \mathrm{r} \mathbf{T}_u \mathbf{T}_{\overline{u}} -z \bigr)^{-1}$, which is a priori only defined on the Hardy space $L^2_+$, to the Hardy--Zhidkov space $\mathcal{Z}^1_+$. The extension is expressed through the Poisson integral and the difference quotient operator given by Subsection $\ref{SubsPoiInDiQuo}$.
 
\begin{lem}\label{LembfIinvbbGtper}
For any $t ,\mathrm{r} \in \mathbb{R}$, $z \in \mathbb{C}_+$ and $u \in L^{\infty} (\mathbb{R} ; \mathbb{C})$, the mapping $\mathbf{I}_+ \bigl( \mathbb{G}_t + 2 \mathrm{r} \mathbf{T}_u \mathbf{T}_{\overline{u}} -z \bigr)^{-1}$ is a bounded linear functional on $L^2_+$. Let $\mathbf{P}_z \in \mathcal{B}_{\mathbb{C}} (L^2_+ ; \mathbb{C})$ denote the Poisson integral operator given by 
\eqref{PoiIntDefSect2}. The difference quotient operator $\mathscr{Q}_z \in \mathcal{B}_{\mathbb{C}} (L^2_+)$ is defined in  
\eqref{diffQuoDefForQzS2}. Then we have
\begin{equation}\label{bfIbbGtpertExpsum}
\tfrac{1}{2 \pi i } \mathbf{I}_+ \bigl( \mathbb{G}_t + 2 \mathrm{r} \mathbf{T}_u \mathbf{T}_{\overline{u}} -z \bigr)^{-1}  = \Bigl( \mathbf{P}_z e^{it \partial_x^2} - \tfrac{\mathrm{r}}{\pi i}\mathbf{I}_+ \bigl( \mathbb{G}_t + 2 \mathrm{r} \mathbf{T}_u \mathbf{T}_{\overline{u}} -z \bigr)^{-1}  \mathbf{T}_u \mathbf{T}_{\overline{u}} e^{ - it \partial_x^2} \mathscr{Q}_z   e^{it \partial_x^2} \Bigr) \Big|_{L^2_+} .
\end{equation}
\end{lem}

\begin{proof}
Fix $t , \mathrm{r} \in \mathbb{R}$ and $u \in L^{\infty} (\mathbb{R} ; \mathbb{C})$, let $\mho_u^{(t)} \in \mathcal{B}_{\mathbb{C}}(L^2_+)$ be defined in \eqref{dombfGpluttZut}.  
 Fix $z \in \mathbb{C}_+$. Then
\begin{equation}\label{sumidandmho}
\mathrm{id}_{L^2_+} + 2 \mathrm{r} \mho_u^{(t)} (\mathbf{G} - z)^{-1} =  \bigl( \mathbf{G} + 2 \mathrm{r} \mho_u^{(t)}   - z\bigr) (\mathbf{G} - z)^{-1} : L^2_+ \to L^2_+
\end{equation}is a  bijective continuous linear mapping by estimates 
\eqref{EstG-z-1BHs} and 
 \eqref{resolGpluZutEstuni}. By the open mapping theorem, its inverse is also a   bounded linear operator on $L^2_+$. Moreover, we have
 \begin{equation}\label{idplmhoresolGinv}
\small
  \bigl( \mathrm{id}_{L^2_+} + 2 \mathrm{r} \mho_u^{(t)} (\mathbf{G} - z)^{-1}  \bigr)^{-1}
   = \mathrm{id}_{L^2_+}  - 2 \mathrm{r}  \bigl( \mathrm{id}_{L^2_+} + 2 \mathrm{r} \mho_u^{(t)} (\mathbf{G} - z)^{-1}  \bigr)^{-1} \mho_u^{(t)} (\mathbf{G} - z)^{-1} \in \mathcal{B}_{\mathbb{C}} (L^2_+) .
 \end{equation}Formula   \eqref{InvForNewtonLeibniz} implies that $\mathbf{I}_+ (\mathbf{G} - z )^{-1} = 2 \pi i \mathbf{P}_z \in \mathcal{B}_{\mathbb{C}} (L^2_+ ; \mathbb{C})$. So we have
\begin{equation}\label{bfIGtperttobfPz}
\begin{split}
& \mathbf{I}_+ \bigl( \mathbb{G}_t + 2 \mathrm{r} \mathbf{T}_u \mathbf{T}_{\overline{u}} -z \bigr)^{-1}  
 =    \mathbf{I}_+ \bigl( \mathbf{G} + 2 \mathrm{r}  \mho_u^{(t)} -z \bigr)^{-1} e^{it \partial_x^2}    \\
= & 2 \pi i \mathbf{P}_z \bigl( \mathrm{id}_{L^2_+} + 2 \mathrm{r} \mho_u^{(t)} (\mathbf{G} - z)^{-1}  \bigr)^{-1} e^{it \partial_x^2}  \in \mathcal{B}_{\mathbb{C}}(L^2_+ ; \mathbb{C}) .
\end{split}
\end{equation}by   \eqref{CjRebbGtplusTutub},  \eqref{bfIeitLapisbfI} and   \eqref{sumidandmho}. 
Plugging   \eqref{idplmhoresolGinv} and 
\eqref{dombfGpluttZut}  into  \eqref{bfIGtperttobfPz}, we obtain that
\begin{equation}\label{bfIbbGtpertG-z}
\small
 \mathbf{I}_+ \bigl( \mathbb{G}_t + 2 \mathrm{r} \mathbf{T}_u \mathbf{T}_{\overline{u}} -z \bigr)^{-1}  = (2 \pi i )\mathbf{P}_z e^{it \partial_x^2} - 2\mathrm{r} \mathbf{I}_+ \bigl( \mathbb{G}_t + 2 \mathrm{r} \mathbf{T}_u \mathbf{T}_{\overline{u}} -z \bigr)^{-1}  \mathbf{T}_u \mathbf{T}_{\overline{u}} e^{ - it \partial_x^2} (\mathbf{G} - z)^{-1}   e^{it \partial_x^2}  .
\end{equation}Finally, formula 
 \eqref{bfIbbGtpertExpsum} follows from \eqref{bfIbbGtpertG-z} and
\eqref{qzfisvarphiinDomG}.
\end{proof}

\begin{rem}\label{remReaConvG-ztoQz}
The derivation of the explicit formula for \eqref{CMdNLSfocu} on the Hardy--Zhidkov space $\mathcal{Z}^3_+$ relies on approximating the initial datum $u_0 \in \mathcal{Z}^3_+$ by $u_0\chi_{\epsilon}$, which belongs to $H^3_+ \subset L^2_+$. 
However, since $u_0$ itself need not lie in $L^2_+$, the $L^2_+$-inner product and the resolvent $(\mathbf{G}-z)^{-1}$ cannot be applied to $u_0$ directly. 
The approximation argument, based on \eqref{htepuphito0inL2} and \eqref{QzeitLapuchiepLimHN}, circumvents this difficulty.  
 To take advantage of this approximation, we replace  the resolvent $(\mathbf{G}-z)^{-1}$ by the difference quotient operator $\mathscr{Q}_z$ and  transform formula \eqref{bfIbbGtpertG-z} into \eqref{bfIbbGtpertExpsum}.  
By Lemma $\ref{scQzinL2+LemG-z}$, $\mathscr{Q}_z$ extends the action of $(\mathbf{G}-z)^{-1}$ to the Hardy--Zhidkov space  $\mathcal{Z}^3_+$. 
Consequently, 
  $\lim_{\epsilon \to 0^+}\mathbf{I}_+(\mathbb{G}_t + 2\mathrm{r}\mathbf{T}_{u_0} \mathbf{T}_{\overline{u}_0} - z)^{-1}(u_0\chi_{\epsilon})$   can be expressed  explicitly in terms of the initial datum $u_0$ and the time $t$ from the right-hand side of \eqref{bfIbbGtpertExpsum}.
\end{rem}
  
\begin{lem}
For any $t ,\mathrm{r} \in \mathbb{R}$, $z \in \mathbb{C}_+$ and $u \in  \mathcal{Z}^1_+$, the following formula holds: 
\begin{equation}\label{limbfIbbGtpertInv}
\begin{split}
 & \tfrac{1}{2 \pi i } \lim_{\epsilon \to 0^+} \mathbf{I}_+ \bigl( \mathbb{G}_t + 2 \mathrm{r} \mathbf{T}_u \mathbf{T}_{\overline{u}} -z \bigr)^{-1} (u \chi_{\epsilon} )\\
  = & \mathbf{P}_z e^{it \partial_x^2} (u) 
  - \tfrac{\mathrm{r}}{\pi i}\mathbf{I}_+ \bigl( \mathbb{G}_t + 2 \mathrm{r} \mathbf{T}_u \mathbf{T}_{\overline{u}} -z \bigr)^{-1}  \mathbf{T}_u \mathbf{T}_{\overline{u}} e^{ - it \partial_x^2} \mathscr{Q}_z   e^{it \partial_x^2} (u) .
\end{split}
\end{equation}
\end{lem}

\begin{proof}
Formula \eqref{limbfIbbGtpertInv}  follows directly from \eqref{bfIbbGtpertExpsum},  \eqref{htepuphito0inL2} and \eqref{QzeitLapuchiepLimHN}.
\end{proof}
  
 \noindent The following lemma establishes the continuity of the nonlinear functional appearing in the right-hand side of explicit formula for general solutions to \eqref{CMdNLSfocu}.

\begin{lem}\label{lemContRHSExpFor}
For any $t, \mathrm{r} \in \mathbb{R}$, $z \in \mathbb{C}_+$ and $u \in \mathcal{Z}^1_+$, we define
\begin{equation}\label{defbethtrzufor}
\beth_{t, \mathrm{r}}^z (u) := \mathbf{P}_z e^{it \partial_x^2}(u) 
- \tfrac{\mathrm{r}}{\pi i} \mathbf{I}_+ \bigl( \mathbb{G}_t + 2 \mathrm{r} \mathbf{T}_u \mathbf{T}_{\overline{u}} - z \bigr)^{-1} \mathbf{T}_u \mathbf{T}_{\overline{u}} e^{-it \partial_x^2} \mathscr{Q}_z e^{it \partial_x^2}(u) \in \mathbb{C}.
\end{equation}
Then the mapping $\beth_{t, \mathrm{r}}^z : \mathcal{Z}^1_+ \to \mathbb{C}$ is a continuous nonlinear functional.
\end{lem}
 
\begin{proof} 
By    
Definition $\ref{AudefTuTbaru}$,  the mapping $\mathbf{A} : u \in L^{\infty} (\mathbb{R} ; \mathbb{C}) \mapsto  \mathbf{A}(u) = \mathbf{A}_u =\mathbf{T}_u \mathbf{T}_{\overline{u}} \in \mathcal{B}_{\mathbb{C}} ( L^2_+ )$ is real-analytic.   
 Fix $t, \mathrm{r} \in \mathbb{R}$ and $z \in \mathbb{C}_+$.  
  Using Lemma $\ref{LemPertGtMaxAcc}$,
 we obtain that
\begin{equation}\label{ResolForGtTTuz}
\begin{split}
& \bigl( \mathbb{G}_t + 2 \mathrm{r} \mathbf{T}_u \mathbf{T}_{\overline{u}} - z \bigr)^{-1} - \bigl( \mathbb{G}_t + 2 \mathrm{r} \mathbf{T}_v \mathbf{T}_{\overline{v}} - z \bigr)^{-1} \in \mathcal{B}_{\mathbb{C}} ( L^2_+ )\\
= & 2 \mathrm{r} \bigl( \mathbb{G}_t + 2 \mathrm{r} \mathbf{T}_u \mathbf{T}_{\overline{u}} - z \bigr)^{-1} \bigl(\mathbf{A}_v - \mathbf{A}_u \bigr) \bigl( \mathbb{G}_t + 2 \mathrm{r} \mathbf{T}_v \mathbf{T}_{\overline{v}} - z \bigr)^{-1} , \quad \forall u , v \in L^{\infty} (\mathbb{R} ; \mathbb{C}) .
\end{split}
\end{equation}
 Then  $\mathrm{Y}_{t, \mathrm{r}}^z :  u \in L^{\infty} (\mathbb{R} ; \mathbb{C}) \mapsto  \mathrm{Y}_{t, \mathrm{r}}^z (u) : = \mathbf{I}_+ \bigl( \mathbb{G}_t + 2 \mathrm{r} \mathbf{T}_u \mathbf{T}_{\overline{u}} -z \bigr)^{-1}  \in 
\mathcal{B}_{\mathbb{C}} (L^2_+ ; \mathbb{C})$ is a continuous mapping by \eqref{ResolForGtTTuz} and    \eqref{EstuniresolGtpluTuTub}.
 Formula  \eqref{SchrGrUnitHs} and 
Lemma $\ref{lemSchGrHardZn}$   imply that $e^{it \partial_x^2} \in \mathcal{B}_{\mathbb{C}}(H^1_+) \cap \mathcal{B}_{\mathbb{C}}(\mathcal{Z}^1_+)$, for any $t \in \mathbb{R}$. 
 It follows from  \eqref{NormbdQzZNHN} and 
Lemma 
$\ref{QzLinftyinHardyLem}$ that $\mathscr{Q}_z     \in \mathcal{B}_{\mathbb{C}} (\mathcal{Z}^1_+ ;  H^1_+ )$. 
 So $\mathtt{g}_t^z : = e^{ - it \partial_x^2} \mathscr{Q}_z  e^{it \partial_x^2}  |_{\mathcal{Z}^1_+}$ is a bounded linear operator from $\mathcal{Z}^1_+$ to $H^1_+$. Consequently, the mapping  $\boldsymbol{\Sigma}_{t ,\mathrm{r}}^z :u \in \mathcal{Z}^1_+  \mapsto \boldsymbol{\Sigma}_{t ,\mathrm{r}}^z (u) \in \mathbb{C}$ with
\begin{equation}\label{bdSigtzudeffor}
 \boldsymbol{\Sigma}_{t ,\mathrm{r}}^z (u)  : =
\mathrm{Y}_{t, \mathrm{r}}^z (u) \mathbf{A}_u \bigl( \mathtt{g}_t^z (u) \bigr)
 = \mathbf{I}_+ \bigl( \mathbb{G}_t + 2 \mathrm{r} \mathbf{T}_u \mathbf{T}_{\overline{u}} - z \bigr)^{-1} \mathbf{T}_u \mathbf{T}_{\overline{u}} e^{-it \partial_x^2} \mathscr{Q}_z e^{it \partial_x^2}(u) \in \mathbb{C},
\end{equation}is a continuous nonlinear functional. 
Since $\mathbf{P}_z e^{it \partial_x^2} \in  \mathcal{B}_{\mathbb{C}}  ( \mathcal{Z}^1_+ ; \mathbb{C} ) $ by   Lemma $\ref{lemSchGrHardZn}$ and formula  \eqref{PzL2BdPz1is1}, we have  
$\beth_{t, \mathrm{r}}^z =  \mathbf{P}_z e^{it \partial_x^2} + \boldsymbol{\Sigma}_{t ,\mathrm{r}}^z \in \mathcal{C}  ( \mathcal{Z}^1_+ ; \mathbb{C} ) $.
\end{proof}

\subsection{Commutator Formulas}\label{SsSemiCommEsti3.4}

\noindent This subsection is devoted to the commutation relations between the shift semigroup $(\mathbf{S}(\eta)^*)_{\eta \geq 0}$ and the Toeplitz operators, the chiral Lax operator $\mathbf{L}_u$, and the operator $\widetilde{\mathbf{B}}_u$. These formulas will provide essential tools for establishing the simplicity of eigenvalues, characterizing invariant subspaces, deriving the explicit formula, and constructing action-angle variables. \\

\noindent As a preparation, we first compute the commutators of  $(\mathbf{S}(\eta)^*)_{\eta \geq 0}$ defined in  
\eqref{LaxBeurlingSG}  with the differential operators $D = -i \partial_x$ and $D^2$, which will be building blocks for  more complicated commutators with Toeplitz and Lax operators.
  For any $ p \geq 0$,   Lemma $\ref{ContSemiGrLem}$ shows that $(\mathbf{S}(\eta)^*)_{\eta \geq 0}$  is a contraction semigroup on $H^p_+$.  
  For any $s \geq 1$ and $\xi, \eta >0$, if $g \in H^s_+$, we have $\frac{1}{\eta}\left(\left[\mathbf{S}(\eta)^*, D\right]g \right)^{\wedge} (\xi) = \widehat{g}(\xi +\eta)$.  Then   
\begin{equation}\label{CommSetastDisSst}
\tfrac{1}{\eta} \left[\mathbf{S}(\eta)^*, D\right]g = \mathbf{S}(\eta)^* g \to g \quad \mathrm{in} \quad  H^s_+, \quad \mathrm{as} \quad \eta \to 0^+,  
\end{equation}for any $g \in H^s_+$, by the strong continuity of the semigroup $(\mathbf{S}(\eta)^*)_{\eta \geq 0}$ on $H^s_+$, $\forall s \geq 1$.

  \begin{lem}
  For any $s \geq 1$ and $\eta >0$, we have $\frac{1}{\eta} \left[\mathbf{S}(\eta)^*, D\right] = \mathbf{S}(\eta)^* \in \mathcal{B}_{\mathbb{C}}(H^s_+)$. Moreover, we have
  \begin{equation}\label{BHsnorm[S*,D]-1}
  1 \leq \|\tfrac{1}{\eta} \left[\mathbf{S}(\eta)^*, D\right] - \mathrm{id}_{H^s_+}\|_{\mathcal{B}_{\mathbb{C}}(H^s_+)} \leq 2, 
  \quad \forall \eta >0, \quad \forall s \geq 1.
  \end{equation}
  \end{lem}
\begin{proof}
Assume that $g \in H^s_+$ for some $s \geq 1$. By formulas \eqref{SetastcontraSemiGr} and \eqref{CommSetastDisSst},  it follows that 
\begin{equation*}
\|\tfrac{1}{\eta} \left[\mathbf{S}(\eta)^*, D\right] g - g\|_{H^s_+} =  \| \mathbf{S}(\eta)^* g -g  \|_{H^s_+} \leq 2 \|g \|_{H^s_+}, \quad \forall \eta >0.
\end{equation*}Fix $\eta >0$, there exists $\mathfrak{g} \in H^s_+$ such that
\begin{equation*}
\widehat{\mathfrak{g}}(\xi)=\begin{cases}
1, \quad \mathrm{if} \quad \xi  \in (\eta,  \frac{3}{2}\eta);\\
0, \quad \mathrm{if} \quad \xi  \in  \mathbb{R}\setminus (\eta,  \frac{3}{2}\eta).\\
\end{cases}
\end{equation*}Then  $\operatorname{supp}(\widehat{\mathfrak{g}}) \cap  \operatorname{supp}(\tau_{-\eta}\widehat{\mathfrak{g}}) = \emptyset$. So $\|\tfrac{1}{\eta} \left[\mathbf{S}(\eta)^*, D\right] \mathfrak{g} - \mathfrak{g}\|_{H^s_+}^2 = \|\mathfrak{g}\|_{H^s_+}^2 + \|\mathbf{S}(\eta)^* \mathfrak{g}\|_{H^s_+}^2 \geq \|\mathfrak{g}\|_{H^s_+}^2$.
\end{proof}

\begin{lem}\label{LemcvComm[S*,D2]-2D}
For any $h\in H^2_+$, we have $\lim_{\eta \to 0^+}\|\frac{1}{\eta}\left[\mathbf{S}(\eta)^* , D^2 \right] h - 2 Dh\|_{H^1_+} = 0$. Moreover, we have
\begin{equation}\label{Comm[S*,D2]-2DEsti}
\sup_{\eta >0 }\|\tfrac{1}{\eta}\left[\mathbf{S}(\eta)^* , D^2 \right]   - 2 D \|_{\mathcal{B}_{\mathbb{C}}(H^2_+; H^1_+)} \leq 4 .
\end{equation}
\end{lem}
\begin{proof}
If $h \in H^2_+$, set $g:=Dh \in H^1_+$. Then formula \eqref{CommSetastDisSst} yields that
\begin{equation}\label{comm[S*,D2htog}
\|\tfrac{1}{\eta}\left[\mathbf{S}(\eta)^* , D^2 \right] h - 2 Dh \|_{H^1_+} \leq \| \tfrac{1}{\eta}\left[\mathbf{S}(\eta)^* , D  \right] g -  g \|_{H^1_+} + \|D \left(\tfrac{1}{\eta}\left[\mathbf{S}(\eta)^* , D  \right] h - h  \right)\|_{H^1_+} \to 0,
\end{equation}as $\eta \to 0^+$.  Then  \eqref{Comm[S*,D2]-2DEsti}  is deduced from  \eqref{BHsnorm[S*,D]-1}, \eqref{comm[S*,D2htog} and the inequality $\| D h \|_{H^1_+} \leq \|h\|_{H^2_+}$.
\end{proof}
 
\begin{defi}\label{defibolGammSet*,Dsq-2D}
For any $\eta > 0$, we set $\boldsymbol{\Gamma}^{(\eta)} : = \tfrac{1}{\eta}\left[\mathbf{S}(\eta)^* , D^2 \right]   - 2 D : H^2_+ \to L^2_+$.
\end{defi}
\noindent By Lemma $\ref{LemcvComm[S*,D2]-2D}$ and inequality  \eqref{Comm[S*,D2]-2DEsti}, it follows that $\boldsymbol{\Gamma}^{(\eta)} \in \mathcal{B}_{\mathbb{C}}(H^2_+; H^1_+ )$, for any $\eta >0$,  
\begin{equation}\label{boldGammaetEstimH2H1}
\sup_{\eta >0} \| \boldsymbol{\Gamma}^{(\eta)} \|_{\mathcal{B}_{\mathbb{C}}(H^2_+; H^1_+ )} \leq 4 
\qquad \text{and} \qquad 
\lim_{\eta \to 0^+} \|\boldsymbol{\Gamma}^{(\eta)} (h)\|_{H^1_+} = 0, \qquad \forall h \in H^2_+.
\end{equation}
  With these preliminary commutators at hand, we now turn to the commutation relations between the shift semigroup $(\mathbf{S}(\eta)^*)_{\eta \geq 0}$ and  Toeplitz operators given by \eqref{ToepdefForSec2}.

\begin{lem}
For any $\eta \geq 0$, $u \in L^{\infty}(\mathbb{R}; \mathbb{C})$ and $f \in L^2_+ $,   the following identity holds: 
\begin{equation}\label{comfosetastTu}
\left[\mathbf{S}(\eta)^*, \mathbf{T}_u \right](f) = \Pi_{\geq 0} \left( u \Pi_{<0} \left( \mathbf{e}_{-\eta} f \right) \right)\in L^2_+  \quad \text{with} \quad \mathbf{e}_{-\eta} (x) =e^{-i\eta x}, \quad \forall x \in \mathbb{R}.
\end{equation}Here  $\Pi_{\geq 0} : L^2 (\mathbb{R}; \mathbb{C}) \to L^2_+$ and  $\Pi_{<0}  : L^2(\mathbb{R}; \mathbb{C}) \to L^2_-$ are the Szeg\H{o} projectors given by  \eqref{SzegoProjS2}.
\end{lem}

\begin{proof}
Fix $u \in L^{\infty}(\mathbb{R}; \mathbb{C})$ and $\eta \geq 0$. For any $f, g \in L^2_+$, we have $\Pi_{<0} (uf) = uf - \mathbf{T}_u (f) \in L^2_- \perp L^2_+$, $\mathbf{S}(\eta) g = \mathbf{e}_{\eta} g \in L^2_+$ and $\mathbf{S}(\eta)^* f = \Pi_{\geq 0}( \mathbf{e}_{-\eta} f) \in L^2_+$ by \eqref{LaxBeurlingSG}. Using \eqref{L2+perpL2-Har} and \eqref{AdjofToepOpinL2+}, we obtain that
\begin{equation}
\label{Set*Tuf,gL2+}
\langle \mathbf{S}(\eta)^* \mathbf{T}_u (f), g\rangle_{L^2_+} = \langle  f,\mathbf{e}_{\eta} \overline{u} g\rangle_{L^2(\mathbb{R}; \mathbb{C})}  \quad \mathrm{and} \quad \langle  \mathbf{T}_u \mathbf{S}(\eta)^*(f) , g\rangle_{L^2_+} = \langle f , \mathbf{e}_{\eta}\Pi_{\geq 0}(\overline{u} g)\rangle_{L^2_+}.
\end{equation}Both $\Pi_{\geq 0}$ and $\Pi_{<0}$ are bounded self-adjoint operators on $L^2(\mathbb{R}; \mathbb{C})$. As a consequence, we have
\begin{equation}
\langle \left[ \mathbf{S}(\eta)^*,  \mathbf{T}_u \right] (f), g\rangle_{L^2_+} = \langle  f,\mathbf{e}_{\eta} \Pi_{<0}( \overline{u} g)\rangle_{L^2(\mathbb{R}; \mathbb{C})}   = \langle \Pi_{\geq 0} \left( u \Pi_{<0}(\mathbf{e}_{-\eta} f) \right) ,  g \rangle_{L^2_+}, \quad \forall g \in L^2_+,  
\end{equation}by using \eqref{Set*Tuf,gL2+}.  Then formula \eqref{comfosetastTu} is obtained.
\end{proof}

\begin{cor}\label{SadjTbnegcomuLem}
For any $b \in L^{\infty}_- $ and $\eta \geq 0$, we have $ \mathbf{S}(\eta)^*  \mathbf{T}_b =  \mathbf{T}_b \mathbf{S}(\eta)^* \in \mathcal{B}_{\mathbb{C}}(L^2_+)$.
\end{cor}
\begin{proof}
For any $b \in L^{\infty}_- $ and $f\in L^2_+$, we have $b \Pi_{<0} \left( \mathbf{e}_{-\eta} f \right)  \in L^{\infty}_-  \cdot L^2_- \subset L^2_-$ by Corollary $\ref{CorLinfty-L2-}$ and
 formula \eqref{BanAlgebraZn-}.
 Then   \eqref{comfosetastTu} implies that $\left[\mathbf{S}(\eta)^*, \mathbf{T}_b \right](f) =0$, $\forall \eta \geq 0$.
\end{proof}

\begin{rem}
Assume that $\eta > 0$, $\mathfrak{b} \in L^2 \cap  L^{\infty}(\mathbb{R}; \mathbb{C})$ and $f \in L^2_+$. Then $\Pi_{<0}(\mathbf{e}_{- \eta} f) \in L^2_-$. According to the convolution formula for Fourier transform and formula \eqref{comfosetastTu}, we have
\begin{equation}\label{Fou(Setaast,Tu]f)xi}
\left(\left[\mathbf{S}(\eta)^*, \mathbf{T}_{\mathfrak{b}} \right](f)\right)^{\wedge}(\xi) =   \frac{1}{2\pi} \int_{\xi}^{\xi +\eta} \widehat{\mathfrak{b}}(\zeta) \widehat{f}(\xi + \eta - \zeta) \mathrm{d} \zeta, \quad \forall \xi >0.
\end{equation}Note that   \eqref{Fou(Setaast,Tu]f)xi} is the same as formula (3.4) appearing in Lemma 3.1 of \cite{SUNCMP2021}. While  \eqref{Fou(Setaast,Tu]f)xi} requires the symbol $\mathfrak{b}$ to be in $L^2 \cap  L^{\infty}(\mathbb{R}; \mathbb{C})$, formula \eqref{comfosetastTu} only requires $u \in L^{\infty}(\mathbb{R}; \mathbb{C})$. Since $L^2 \cap  L^\infty \subset L^\infty (\mathbb{R}; \mathbb{C})$, formula \eqref{comfosetastTu} is valid for a broader class of symbols and thus constitutes a generalization of \eqref{Fou(Setaast,Tu]f)xi}.  Lemma 3.1 of \cite{SUNCMP2021} yields that $\mathbf{T}_{\mathfrak{b}} (\mathrm{Dom}(\mathbf{G}))\subset \mathrm{Dom}(\mathbf{G})$ and the following identity holds,
\begin{equation}\label{Comm[G,Tb]For}
2 \pi \left[\mathbf{G}, \mathbf{T}_{\mathfrak{b}}\right] = i \widehat{\varphi}(0^+) \Pi_{\geq 0} \mathfrak{b}, \quad \forall \varphi \in \mathrm{Dom}  (\mathbf{G}).
\end{equation}If $\mathfrak{v} \in H^s(\mathbb{R}; \mathbb{C})$ and $\phi \in \mathrm{Dom}^s(\mathbf{G})$ for some $s > \frac{1}{2}$, we have $\mathbf{G}\mathbf{T}_{\mathfrak{v}} \phi =   \tfrac{i}{2\pi}  \widehat{\phi}(0^+)\Pi_{\geq 0} \mathfrak{v} + \mathbf{T}_{\mathfrak{v}} \mathbf{G}\phi \in H^s_+$, by using \eqref{BounToeplitzBHs}, \eqref{DomsGDefSect2} and \eqref{Comm[G,Tb]For}.  We use \eqref{DomsGDefSect2} again to obtain that $\mathbf{T}_{\mathfrak{v}} \phi \in \mathrm{Dom}^s(\mathbf{G})$. So we have
\begin{equation}\label{TvDsGinDsGs>0.5}
\mathbf{T}_{\mathfrak{v}} (\mathrm{Dom}^s(\mathbf{G}))\subset \mathrm{Dom}^s (\mathbf{G}), \quad \forall \mathfrak{v} \in H^s(\mathbb{R}; \mathbb{C}), \quad \forall s >\tfrac{1}{2}.
\end{equation}
\end{rem}
\noindent Then Lemma 3.1 of \cite{SUNCMP2021} can be generalized to the following lemma.
\begin{lem}
For any $s \geq 0$, $\mathfrak{b} \in H^s \cap  L^{\infty}(\mathbb{R}; \mathbb{C})$ and $f \in \mathrm{Dom}(\mathbf{I}_+)$, we have
\begin{equation}\label{HsCVforSe,Tb]fs>0}
\lim_{\eta \to 0^+}\|\tfrac{2 \pi}{\eta} \left[\mathbf{S}(\eta)^*, \mathbf{T}_{\mathfrak{b}} \right](f) - \widehat{f}(0^+) \Pi_{\geq 0} \mathfrak{b} \|_{H^s_+} =0.
\end{equation}
\end{lem}
\begin{proof}
Fix $s  \geq 0$, $\mathfrak{b} \in H^s \cap  L^{\infty}(\mathbb{R}; \mathbb{C})$ and $f \in \mathrm{Dom}(\mathbf{I}_+) \subset L^2_+$. Formula \eqref{Fou(Setaast,Tu]f)xi} yields that 
\begin{equation}\label{FouChavarInttheta}
\tfrac{1}{\eta}\left(\left[\mathbf{S}(\eta)^*, \mathbf{T}_{\mathfrak{b}} \right](f)\right)^{\wedge}(\xi) =   \frac{1}{2\pi} \int_0^1 \widehat{\mathfrak{b}}(\xi + \theta \eta) \widehat{f}((1-\theta) \eta)\mathrm{d} \theta, \quad \forall \xi >0, \quad \forall \eta >0.
\end{equation}For any $\eta >0$, set $\mathfrak{a}_{\eta} = \mathfrak{a}_{\eta}(f) : = \frac{1}{2\pi} \int_0^1 \widehat{f}((1-\theta) \eta)\mathrm{d} \theta \in \mathbb{C}$. We define $\phi_{\eta} = \phi_{\eta}(\mathfrak{b}, f)\in L^2_+$ as follows,
\begin{equation}\label{phietaFouTrdef}
\widehat{\phi_{\eta}}(\xi) := \frac{1}{2\pi} \int_0^1 \left( \widehat{\mathfrak{b}}(\xi + \theta \eta) - \widehat{\mathfrak{b}}(\xi) \right) \widehat{f}((1-\theta) \eta)\mathrm{d} \theta, \quad \forall \xi>0.
\end{equation}Since $f \in \mathrm{Dom}(\mathbf{I}_+)$,  the limit $ \lim_{y \to 0^+} \widehat{f}(y) =\widehat{f}(0^+)$ exists. Specifically, for any $\epsilon \in (0, 1]$, there exists $\delta = \delta(f, \epsilon) \in (0, 1)$ such that $\sup_{0 < y < \delta}|\widehat{f}(y) - \widehat{f}(0^+)| \leq \epsilon$. If $0 < \eta < \delta$, then we have
\begin{equation}\label{aetatohatf0+limit}
|2 \pi \mathfrak{a}_{\eta} - \widehat{f}(0^+)| \leq \int_0^1 |\widehat{f}((1-\theta)\eta)  - \widehat{f}(0^+) | \mathrm{d}\theta \leq \epsilon \Longrightarrow 2 \pi\lim_{\eta \to 0^+}\mathfrak{a}_{\eta} = \widehat{f}(0^+).
\end{equation}Using Plancherel identity, Cauchy--Schwarz inequality and Lemma $\ref{TranslationCV}$, we obtain that
\begin{equation}\label{phietaEsti}
8 \pi^3 \|\phi_{\eta}\|_{H^s_+}^2 \leq   \int_0^1 |\widehat{f}((1 -\theta)\eta)|^2 \left(\int_0^{+\infty} |\widehat{\mathfrak{b}}(\xi + \theta \eta) - \widehat{\mathfrak{b}}(\xi)|^2 (1+ \xi^2 )^s \mathrm{d}\xi \right) \mathrm{d}\theta = \int_0^1 \mathfrak{p}(\theta, \eta) \mathrm{d}\theta,
\end{equation}with $\mathfrak{p}(\theta, \eta):=  |\widehat{f}((1 -\theta)\eta)|^2 \|\tau_{-\theta \eta} \widehat{\mathfrak{b}} - \widehat{\mathfrak{b}}\|_{L^2_s(\mathbb{R}_+^*)}^2 \geq 0$, for any $\theta, \eta \in (0, 1)$. Since $\mathfrak{b} \in H^s(\mathbb{R}; \mathbb{C})$, we have $\widehat{\mathfrak{b}} \in L^2_s(\mathbb{R})$.  Fix $ \theta \in (0, 1)$,   Lemma $\ref{TranslationCV}$ yields that
\begin{equation}\label{LimP(theta,0+)=0+}
 0 \leq \mathfrak{p} (\theta, \eta) \leq |\widehat{f}((1 -\theta)\eta)|^2 \|\tau_{-\theta \eta} \widehat{\mathfrak{b}} - \widehat{\mathfrak{b}}\|_{L^2_s(\mathbb{R})}^2 \to |\widehat{f}(0^+)|^2 0 = 0, \quad \mathrm{as} \quad \eta \to 0^+.
\end{equation}For any $\eta \in (0, \delta(f, 1))\subset (0,1)$ and $\theta \in (0,1)$, we have $|\widehat{f}((1 - \theta)\eta)| \leq 1+ |\widehat{f}(0^+)|$. Inequality  \eqref{PeetreCorUni} implies that $\sup_{|y| \leq 1}\|\tau_{y}\widehat{\mathfrak{b}}\|_{L^2_s(\mathbb{R})}^2 \leq 2 \pi 8^s \|\mathfrak{b}\|^2_{H^s(\mathbb{R}; \mathbb{C})}$. As a consequence, we have
\begin{equation}\label{unifBounPthetaeta}
\sup_{0 < \theta <1} \sup_{0 < \eta < \delta(f, 1)} \mathfrak{p} (\theta, \eta) \leq 2 \pi \left( 1+ |\widehat{f}(0^+)| \right)^2\left( 1+ (2\sqrt{2})^s\right)^2 \|\mathfrak{b}\|_{H^s(\mathbb{R}; \mathbb{C})}^2 < +\infty.
\end{equation}According to   \eqref{phietaEsti}, \eqref{LimP(theta,0+)=0+} and  \eqref{unifBounPthetaeta}, we use the dominated convergence theorem to obtain that
\begin{equation}\label{phietato0inHs}
8 \pi^3 \|\phi_{\eta}\|_{H^s_+}^2 \leq \int_0^1 \mathfrak{p} (\theta, \eta) \mathrm{d}\theta \to 0, \quad \mathrm{as} \quad \eta \to 0^+.
\end{equation}So the function $\tfrac{2 \pi}{\eta} \left[\mathbf{S}(\eta)^*, \mathbf{T}_{\mathfrak{b}} \right](f) = 2 \pi\left( \mathfrak{a}_{\eta} \Pi_{\geq 0} b + \phi_{\eta}\right)$ converges to $\widehat{f}(0^+) \Pi_{\geq 0} b$ in $ H^s_+ $, as $\eta \to 0^+$, according to formulas \eqref{FouChavarInttheta}, \eqref{phietaFouTrdef}, \eqref{aetatohatf0+limit} and \eqref{phietato0inHs}.
\end{proof}

\begin{lem}
For any $s \geq 0$, assume that $\mathfrak{b} \in H^s  \cap  L^{\infty}(\mathbb{R}; \mathbb{C})$  and $\mathfrak{v}  \in L^2_+ \cap  L^{\infty}(\mathbb{R}; \mathbb{C})$, then
\begin{equation}\label{CvSeta*,TbTvbaronHs}
\lim_{\eta \to 0^+}\|\tfrac{2 \pi}{\eta} \left[\mathbf{S}(\eta)^*, \mathbf{T}_{\mathfrak{b}} \mathbf{T}_{\overline{\mathfrak{v}}} \right](g) - \langle g, \mathfrak{v} \rangle_{L^2_+} \Pi_{\geq 0} \mathfrak{b} \|_{H^s_+} =0, \quad \forall g \in L^2_+.
\end{equation}
\end{lem}
\begin{proof}
Since $\mathfrak{v} \in L^2_+ \cap  L^{\infty}(\mathbb{R}; \mathbb{C})$ and $\mathfrak{b} \in  H^s\cap  L^{\infty}(\mathbb{R} ; \mathbb{C})$ for some $s \geq 0$, Corollary $\ref{SadjTbnegcomuLem}$ yields that
\begin{equation}\label{Set*,TbTvbar]for}
\left[\mathbf{S}(\eta)^*, \mathbf{T}_{\mathfrak{b}} \mathbf{T}_{\overline{\mathfrak{v}}} \right] = \left[\mathbf{S}(\eta)^*, \mathbf{T}_{\mathfrak{b}} \right] \mathbf{T}_{\overline{\mathfrak{v}}} +\mathbf{T}_{\mathfrak{b}}  \left[\mathbf{S}(\eta)^*,   \mathbf{T}_{\overline{\mathfrak{v}}} \right] = \left[\mathbf{S}(\eta)^*, \mathbf{T}_{\mathfrak{b}} \right] \mathbf{T}_{\overline{\mathfrak{v}}} \in \mathcal{B}_{\mathbb{C}}(L^2_+).
\end{equation}Fix $g \in L^2_+$, set $f : = \mathbf{T}_{\overline{\mathfrak{v}}} (g) = \Pi_{\geq 0} (\overline{\mathfrak{v}}g) \in L^2_+$. Then $\widehat{f}\in \mathcal{C}(\mathbb{R}_+^*; \mathbb{C})\cap  L^2(\mathbb{R}; \mathbb{C})$ and we have
\begin{equation}\label{limhhat0+=intgvbar}
\widehat{f}(0^+) = \lim_{\xi \to 0^+}  \widehat{\overline{\mathfrak{v}}g}(\xi) = \int_{\mathbb{R}} g \overline{\mathfrak{v}} = \langle g, \mathfrak{v} \rangle_{L^2_+}.
\end{equation}So $\mathbf{T}_{\overline{\mathfrak{v}}} (g) = f \in \mathrm{Dom}(\mathbf{I}_+)$ and formula \eqref{CvSeta*,TbTvbaronHs} is obtained by inserting \eqref{Set*,TbTvbar]for} and \eqref{limhhat0+=intgvbar} into \eqref{HsCVforSe,Tb]fs>0}.
\end{proof}

\begin{lem}
Assume that $u = e^{i \alpha} + v$ for some $\alpha \in \mathbb{T} = \mathbb{R}/ 2 \pi \mathbb{Z}$ and $v \in H^1_+$, then we have
\begin{equation}\label{CommS*,TuTubaru=1+vfor}
\lim_{\eta \to 0^+} \Bigl\Vert\tfrac{2 \pi}{\eta} \left[\mathbf{S}(\eta)^*, \mathbf{T}_{u} \mathbf{T}_{\overline{u}} \right](\varphi) - \left( e^{-i\alpha} \widehat{\varphi}(0^+) + \langle \varphi, v \rangle_{L^2_+} \right) v  \Bigl\rVert_{H^1_+} =0, \quad \forall \varphi \in \mathrm{Dom}(\mathbf{I}_+).  
\end{equation}Let $\mathbf{L}_u : H^1_+ \to L^2_+$ denote the chiral Lax operator in  \eqref{LuSect2}. For any $h \in \mathrm{Dom}(\mathbf{I}_+) \cap  H^1_+$, we have
\begin{equation}\label{ForcommS*,L]phi}
\lim_{\eta \to 0^+} \Big\| \tfrac{1}{\eta}\left[\mathbf{S}(\eta)^*, \mathbf{L}_u \right] h - h  + \left(e^{-i \alpha} \widehat{h}(0^+) + \langle h ,v \rangle_{L^2_+} \right) \tfrac{v}{2\pi} \Big\|_{H^1_+} = 0.
\end{equation}
\end{lem}
 
\begin{proof}
If $v \in H^1_+$, then $\left[\mathbf{S}(\eta)^*,  \mathrm{id}_{L^2_+} \right]=\left[\mathbf{S}(\eta)^*, \mathbf{T}_{\overline{v}} \right]=0$ by Corollary $\ref{SadjTbnegcomuLem}$. Since $u=e^{i\alpha}+v$, we have 
\begin{equation}\label{TuTubaris1andTv}
\mathbf{T}_{u} \mathbf{T}_{\overline{u}}  = \mathrm{id}_{L^2_+} + e^{-i\alpha} \mathbf{T}_v + e^{ i\alpha} \mathbf{T}_{\overline{v}} +  \mathbf{T}_{v} \mathbf{T}_{\overline{v}} \in \mathcal{B}_{\mathbb{C}}(H^1_+) \cap  \mathcal{B}_{\mathbb{C}}(L^2_+).
\end{equation}According to formulas \eqref{HsCVforSe,Tb]fs>0}, \eqref{CvSeta*,TbTvbaronHs} and \eqref{TuTubaris1andTv}, for any $\varphi \in \mathrm{Dom}(\mathbf{I}_+)$, it follows that
\begin{equation*}
\tfrac{1}{\eta}\left[\mathbf{S}(\eta)^*,  \mathbf{T}_{u} \mathbf{T}_{\overline{u}} \right] \varphi = \tfrac{e^{-i\alpha}}{\eta} \left[\mathbf{S}(\eta)^*, \mathbf{T}_v \right] \varphi + \tfrac{1}{\eta} \left[\mathbf{S}(\eta)^*, \mathbf{T}_v \mathbf{T}_{\overline{v}} \right] \varphi \to \left(e^{-i\alpha} \widehat{\varphi}(0^+) + \langle \varphi, v\rangle_{L^2_+} \right)\tfrac{v}{2\pi} \quad \mathrm{in} \quad H^1_+,
\end{equation*}as $\eta \to 0^+$. So  \eqref{CommS*,TuTubaru=1+vfor}    is obtained. 
 Formula \eqref{ForcommS*,L]phi} follows directly from    \eqref{LuSect2},  \eqref{CommSetastDisSst}  and  \eqref{CommS*,TuTubaru=1+vfor}.
\end{proof}

\noindent For any $s \geq 0$, let $\mathrm{Dom}^s (\mathbf{G})$ be defined in formula \eqref{DomsGDefSect2} and $\mathrm{Dom}(\mathbf{G}) = \mathrm{Dom}^0(\mathbf{G})$.

 \begin{lem}\label{LuD1GinDGLem}
If $u = e^{i \alpha } + v$ for some 
 $\alpha \in \mathbb{T}$ and $v \in H^1_+$, then $\mathbf{L}_u (\mathrm{Dom}^1(\mathbf{G})) \subset \mathrm{Dom}(\mathbf{G})$ and 
 \begin{equation}\label{ForcommGbfLu}
\left[\mathbf{G}, \mathbf{L}_u \right]  \phi = i \phi - \left(e^{-i \alpha} \widehat{\phi}(0^+) + \langle \phi , v \rangle_{L^2_+} \right) \tfrac{iv}{2\pi} \in H^1_+, \qquad \forall \phi \in \mathrm{Dom}^1(\mathbf{G}).
\end{equation} 
\end{lem}
\begin{proof}
Since $u = e^{i \alpha} + v \in \mathcal{Z}^1_+$, formula \eqref{BounToeplitzBHs} yields that $\mathbf{L}_u = D - \mathbf{T}_u \mathbf{T}_{\overline{u}} \in \mathcal{B}_{\mathbb{C}}(H^1_+; L^2_+)$. We have
\begin{equation}\label{LuCVS*phi-phi}
\Big\|\mathbf{L}_u \left(\tfrac{1}{\eta}(\mathbf{S}(\eta)^* \phi - \phi)\right) + i \mathbf{L}_u \mathbf{G} \phi\Big\|_{L^2_+} \leq \|\mathbf{L}_u\|_{\mathcal{B}_{\mathbb{C}}(H^1_+; L^2_+)}  \big\|   \tfrac{1}{\eta}(\mathbf{S}(\eta)^* \phi - \phi) + i  \mathbf{G} \phi\big\|_{H^1_+} \to 0,
\end{equation}as $\eta \to 0^+$, by  \eqref{SemiGrCVdefDomsG}, for any $\phi \in \mathrm{Dom}^1(\mathbf{G})$. Combining   \eqref{ForcommS*,L]phi} and \eqref{LuCVS*phi-phi}, we obtain that 
\begin{equation}\label{Seta-1LuPhi[G,Lu]}
\begin{split}
\tfrac{1}{\eta}\left(\mathbf{S}(\eta)^* - \mathrm{id}_{L^2_+} \right)\mathbf{L}_u (\phi) = & \tfrac{1}{\eta} \left[\mathbf{S}(\eta)^* , \mathbf{L}_u \right] \phi + \mathbf{L}_u \left(\tfrac{1}{\eta}(\mathbf{S}(\eta)^* \phi - \phi)\right)\\
\to & \phi  - \left(e^{-i \alpha} \widehat{\phi}(0^+) + \langle \phi ,v \rangle_{L^2_+} \right) \tfrac{v}{2\pi} - i \mathbf{L}_u \mathbf{G}\phi \quad \mathrm{in} \quad L^2_+,
\end{split}
\end{equation}as $\eta \to 0^+$. So $\mathbf{L}_u (\phi) \in \mathrm{Dom}(\mathbf{G})$ and formula \eqref{ForcommGbfLu} is deduced from  \eqref{SemiGrCVdefDomsG}  and \eqref{Seta-1LuPhi[G,Lu]}.
\end{proof}

\noindent The following estimates provide quantitative bounds for the commutators studied above, which will be essential in the subsequent analysis of the Lax operator and its spectral properties.

\begin{lem}\label{LemSetaTudhL}
Assume that $u \in L^{\infty}(\mathbb{R}; \mathbb{C})$ and $h \in H^1_+ $, then the following inequality holds
\begin{equation}\label{foSeta*Tu]dhL2est}
2 \pi \eta^{-2}  \|\left[\mathbf{S}(\eta)^*, \mathbf{T}_u \right] \partial_x h \|_{L^2_+ }^2 \leq  \|u\|_{L^{\infty}(\mathbb{R}; \mathbb{C})}^2  \int_0^{\eta} |\widehat{h}(\zeta)|^2\mathrm{d}\zeta , \quad \forall \eta >0.
\end{equation}As a consequence, we have $\lim_{\eta \to 0^+} \|\eta^{-1}\left[\mathbf{S}(\eta)^*, \mathbf{T}_u \right] \partial_x h \|_{L^2_+ } =0$ and
\begin{equation}\label{supSeta*Tu]dhL2norm}
\sup_{\eta >0}\|\eta^{-1}\left[\mathbf{S}(\eta)^*, \mathbf{T}_u \right] \partial_x h \|_{L^2_+ } \leq \|u\|_{L^{\infty}(\mathbb{R}; \mathbb{C})} \|h\|_{L^2_+ }.
\end{equation}
\end{lem}
\begin{proof}
Fix $u \in L^{\infty}(\mathbb{R}; \mathbb{C})$ and $h \in H^1_+$. Set $f: = Dh \in L^2_+$ and $\mathfrak{F}_{\eta}: = \eta^{-1} \Pi_{<0} (\mathbf{e}_{-\eta }f) \in L^2_-$, $\forall \eta>0$. Since $\operatorname{supp}\widehat{f} \subset \mathbb{R}_+ = [0, +\infty)$, we have $\operatorname{supp}\widehat{\mathfrak{F}}_{\eta} \subset [-\eta, 0]$ and $\eta \widehat{\mathfrak{F}}_{\eta} (\xi) = \widehat{f}(\xi + \eta)$, $\forall \xi \in [-\eta, 0]$. Then 
\begin{equation}\label{FetaL2norm}
\|\mathfrak{F}_{\eta}\|_{L^2_-}^2 = \frac{1}{2\pi \eta^2} \int_0^{\eta} \zeta^2 |\widehat{h}(\zeta)|^2 \mathrm{d}\zeta \leq \frac{1}{2\pi} \int_0^{\eta}   |\widehat{h}(\zeta)|^2 \mathrm{d}\zeta 
\leq \|h\|_{L^2_+}^2.
\end{equation}Using Plancherel identity, Cauchy--Schwarz inequality and formula \eqref{comfosetastTu}, we obtain that
\begin{equation}\label{uFetaboundL2}
\|\tfrac{1}{\eta}\left[\mathbf{S}(\eta)^*, \mathbf{T}_u \right] \partial_x h \|_{L^2_+ }^2 = \eta^{-2} \|\Pi_{\geq 0}\left(u \Pi_{<0}(\mathbf{e}_{-\eta}f ) \right)\|_{L^2_+ }^2 = \|\Pi_{\geq 0}\left(u \mathfrak{F}_{\eta} \right)\|_{L^2_+ }^2 \leq \|u\|_{L^{\infty}(\mathbb{R}; \mathbb{C})}^2 \|\mathfrak{F}_{\eta} \|_{L^2_-}^2.
\end{equation}Then   \eqref{foSeta*Tu]dhL2est} and \eqref{supSeta*Tu]dhL2norm} are deduced from \eqref{FetaL2norm} and \eqref{uFetaboundL2}. Finally, since $\widehat{h} \in L^2_1(\mathbb{R}_+^*) $, we have
\begin{equation}\label{inthathto0domcv}
2 \pi \|\mathfrak{F}_{\eta}\|_{L^2_-}^2 
\leq    \int_0^{\eta}   |\widehat{h}(\zeta)|^2 \mathrm{d}\zeta \to 0, \quad \mathrm{as} \quad \eta \to 0^+,  
\end{equation}by using the dominated convergence theorem. The proof is completed.
\end{proof}
\begin{rem}
For any  $\eta >0$ and $u \in L^{\infty}(\mathbb{R}; \mathbb{C})$, the operator $\eta^{-1} \left[\mathbf{S}(\eta)^*, \mathbf{T}_u \right] \partial_x : H^1_+ \to L^2_+$ can be extended to a bounded $\mathbb{C}$-linear operator on $L^2_+$,  according to inequality \eqref{supSeta*Tu]dhL2norm}. Specifically, there exists a unique operator $\mathfrak{Z}_u^{\eta}\in  \mathcal{B}_{\mathbb{C}}(L^2_+)$ such that  if $h \in H^1_+$, we have
\begin{equation}\label{boundforZuetadef}
 \mathfrak{Z}_u^{\eta} (h) = \eta^{-1} \left[\mathbf{S}(\eta)^*, \mathbf{T}_u \right]\partial_x h \in L^2_+ \quad  \mathrm{and}  \quad \sup_{\eta >0}\| \mathfrak{Z}_u^{\eta} \|_{\mathcal{B}_{\mathbb{C}}(L^2_+)} \leq \|u\|_{L^{\infty}(\mathbb{R}; \mathbb{C})}.
\end{equation}We define that $\mathfrak{Z}^{\eta} : u \in L^{\infty}(\mathbb{R}; \mathbb{C}) \mapsto \mathfrak{Z}^{\eta}(u): = \mathfrak{Z}^{\eta}_u \in \mathcal{B}_{\mathbb{C}}(L^2_+) $. Then $(\mathfrak{Z}^{\eta})_{\eta >0}$ is a uniformly equicontinuous family of $\mathbb{C}$-linear mappings from $L^{\infty}(\mathbb{R}; \mathbb{C})$ into $\mathcal{B}_{\mathbb{C}}(L^2_+)$ with $\sup_{\eta >0}\| \mathfrak{Z}^{\eta} \|_{\mathcal{B}_{\mathbb{C}}\left( L^{\infty}(\mathbb{R}; \mathbb{C}); \mathcal{B}_{\mathbb{C}}(L^2_+) \right)}  \leq 1$.
\end{rem}

\begin{cor}
Assume that $u \in L^{\infty}(\mathbb{R}; \mathbb{C})$ and $h \in H^1_+$ such that $\partial_x h \in \mathrm{Dom}(\mathbf{G})$. Then 
\begin{equation}\label{GTudh=TuGdh}
\mathbf{T}_u (\partial_x h ) \in \mathrm{Dom}(\mathbf{G})\quad \mathrm{and} \quad \mathbf{G}\mathbf{T}_u  (\partial_x h )  = \mathbf{T}_u \mathbf{G}  (\partial_x h ) \in L^2_+.
\end{equation}
\end{cor}
\begin{proof}
Fix $u \in L^{\infty}(\mathbb{R}; \mathbb{C})$ and $h \in H^1_+ \cap  \partial_x^{-1} \left( \mathrm{Dom}(\mathbf{G}) \right)$. Proposition $\ref{GisInfiniGeneHsProp}$ yields that
\begin{equation}\label{Tu(S-1)cvto-iTuG}
\Big\|\tfrac{1}{\eta}\mathbf{T}_u \left(\mathbf{S}(\eta)^* - \mathrm{id}_{L^2_+} \right) \partial_x h + i \mathbf{T}_u\mathbf{G}h \Big\|_{L^2_+} \leq \|u\|_{L^{\infty}(\mathbb{R}; \mathbb{C})} \Big\|\tfrac{1}{\eta}  \left(\mathbf{S}(\eta)^* - \mathrm{id}_{L^2_+} \right) \partial_x h + i  \mathbf{G}h \Big\|_{L^2_+} \to 0,
\end{equation}as $\eta \to 0^+$. Since $\left[\mathbf{S}(\eta)^* ,\mathbf{T}_u \right] = \left[\mathbf{S}(\eta)^* -\mathrm{id}_{L^2_+},\mathbf{T}_u \right] \in \mathcal{B}_{\mathbb{C}}(L^2_+)$, $\forall \eta >0$, we have
\begin{equation*}
\tfrac{1}{\eta}\left(\mathbf{S}(\eta)^* -\mathrm{id}_{L^2_+}\right) \mathbf{T}_u \partial_x h = \tfrac{1}{\eta} \left[\mathbf{S}(\eta)^* ,\mathbf{T}_u \right]  \partial_x h +  \tfrac{1}{\eta}\mathbf{T}_u \left(\mathbf{S}(\eta)^* - \mathrm{id}_{L^2_+} \right) \partial_x h \to -i \mathbf{T}_u \mathbf{G} h \quad \mathrm{in} \quad L^2_+,
\end{equation*}as $\eta \to 0^+$, by Lemma $\ref{LemSetaTudhL}$ and formula \eqref{Tu(S-1)cvto-iTuG}. Then \eqref{GTudh=TuGdh} is obtained by Proposition $\ref{GisInfiniGeneHsProp}$.
\end{proof}
 
\begin{cor}
If $u \in L^{\infty}(\mathbb{R}; \mathbb{C})$, $v \in \mathcal{Z}^1$, and $\varphi \in H^1_+ $, then we have $\mathbf{T}_{\overline{v}}(\varphi) \in H^1_+ $ and  
\begin{equation}\label{L2SetaTudxTvbarphi}
\lim_{\eta \to 0^+}\|\eta^{-1}\left[\mathbf{S}(\eta)^*, \mathbf{T}_u \right] \partial_x \mathbf{T}_{\overline{v}}(\varphi) \|_{L^2_+ } = 0.
\end{equation}Moreover, we also have the following uniform estimate  
\begin{equation}\label{UnEsSestTudxTvbphie}
\sup_{\eta >0} \|\eta^{-1}\left[\mathbf{S}(\eta)^*, \mathbf{T}_u \right] \partial_x \mathbf{T}_{\overline{v}}(\varphi) \|_{L^2_+} \leq \|u\|_{L^{\infty}(\mathbb{R}; \mathbb{C})}\|v\|_{L^{\infty}(\mathbb{R}; \mathbb{C})} \|\varphi\|_{L^2_+ }.
\end{equation}
\end{cor}
\begin{proof}
Since $v \in \mathcal{Z}^1$ and $\varphi \in H^1_+$,  we have $h:= \mathbf{T}_{\overline{v}}(\varphi) =\Pi_{\geq 0}(\overline{v}\varphi) \in   H^1_+$ by   \eqref{ufHn<uZnfHn}. Then \eqref{L2SetaTudxTvbarphi} is obtained by \eqref{foSeta*Tu]dhL2est} and \eqref{inthathto0domcv}. 
Formula  \eqref{UnEsSestTudxTvbphie} is a direct consequence of  \eqref{BounToeplitzBHs} and \eqref{supSeta*Tu]dhL2norm}.
\end{proof}

 \begin{lem}
For any $u, v \in \mathcal{Z}^1_+$  and $\varphi \in H^1_+ $,   we have   $\mathbf{T}_u D    \mathbf{T}_{\overline{v}} (\varphi) \in L^2_+$ and 
\begin{equation}\label{SetastTuDTvbarCvto0}
\lim_{\eta \to 0^+}\| \left( \eta^{-1}\left[\mathbf{S}(\eta)^*, \mathbf{T}_u D    \mathbf{T}_{\overline{v}} \right]- \mathbf{T}_u \mathbf{T}_{\overline{v}} \right) \varphi \|_{L^2_+ } = 0.
\end{equation}Moreover, we also have the following uniform estimate 
\begin{equation}\label{UNIEstiSetastTuDTvbarL2}
\sup_{\eta >0}\| \left( \eta^{-1}\left[\mathbf{S}(\eta)^*, \mathbf{T}_u D    \mathbf{T}_{\overline{v}} \right]- \mathbf{T}_u \mathbf{T}_{\overline{v}} \right) \varphi \|_{L^2_+ } \leq 3 \|u\|_{L^{\infty}(\mathbb{R}; \mathbb{C})} \|v\|_{L^{\infty}(\mathbb{R}; \mathbb{C})}  \| \varphi \|_{L^2_+} 
\end{equation}
\end{lem}
\begin{proof}
Fix $u,v \in \mathcal{Z}^1_+$, then $\overline{v} \in \mathcal{Z}^1_-$.  Corollary $\ref{SadjTbnegcomuLem}$ implies  that
\begin{equation}\label{COMMForeqsetastTuDTvbar}
\tfrac{1}{\eta}\left[\mathbf{S}(\eta)^*, \mathbf{T}_u D    \mathbf{T}_{\overline{v}} \right] - \mathbf{T}_u \mathbf{T}_{\overline{v}} = \tfrac{1}{\eta}\left[\mathbf{S}(\eta)^*, \mathbf{T}_u  \right] D    \mathbf{T}_{\overline{v}} + \mathbf{T}_u \left( \tfrac{1}{\eta}\left[\mathbf{S}(\eta)^*,  D     \right] - \mathrm{id}_{H^1_+}\right) \mathbf{T}_{\overline{v}} :H^1_+ \to L^2_+.
\end{equation}Fix $\varphi \in H^1_+$, we set $h := \mathbf{T}_{\overline{v}} (\varphi) \in H^1_+$. By formulas  \eqref{BounToeplitzBHs}, 
\eqref{CommSetastDisSst}
 and
\eqref{L2SetaTudxTvbarphi},
it follows that
\begin{equation}\label{Esti[SeTuDTv]spto2}
\begin{split}
& \Big\|\left(\tfrac{1}{\eta}\left[\mathbf{S}(\eta)^*, \mathbf{T}_u D    \mathbf{T}_{\overline{v}} \right] - \mathbf{T}_u \mathbf{T}_{\overline{v}} \right) \varphi \Big\|_{L^2_+} \\
 \leq  & \big\|\tfrac{1}{\eta}\left[\mathbf{S}(\eta)^*, \mathbf{T}_u  \right]  \partial_x  \mathbf{T}_{\overline{v}} (\varphi) \big\|_{L^2_+} + 
\Big\|\mathbf{T}_u \left( \tfrac{1}{\eta}\left[\mathbf{S}(\eta)^*,  D     \right] - \mathrm{id}_{L^2_+}\right) \mathbf{T}_{\overline{v}} (\varphi) \Big\|_{H^1_+} \to 0,
\end{split}
\end{equation}as $\eta \to 0^+$.  Using formulas
\eqref{BounToeplitzBHs} and \eqref{BHsnorm[S*,D]-1}, we obtain that
\begin{equation}\label{ESTITuSetst-1hL2norm}
\begin{split}
\Big\| \mathbf{T}_u \left( \tfrac{1}{\eta}\left[\mathbf{S}(\eta)^*,  D     \right] - \mathrm{id}_{H^1_+}\right) \mathbf{T}_{\overline{v}} (\varphi)  \Big\|_{L^2_+} \leq & 
\|u \|_{L^{\infty}(\mathbb{R}; \mathbb{C}) }  \|\mathbf{S}(\eta)^* (h) - h \|_{L^2_+ } \\
 \leq  & 2 \|u \|_{L^{\infty}(\mathbb{R}; \mathbb{C}) } \|v \|_{L^{\infty}(\mathbb{R}; \mathbb{C}) } \|\varphi \|_{L^2_+ },
\end{split}
\end{equation}for any $\eta >0$. Finally, formula \eqref{UNIEstiSetastTuDTvbarL2} follows from  \eqref{UnEsSestTudxTvbphie},  \eqref{Esti[SeTuDTv]spto2} and
\eqref{ESTITuSetst-1hL2norm}.
\end{proof}
 
\noindent Recall that  $\widetilde{\mathbf{B}}_u   \in \mathcal{B}_{\mathbb{C}}(H^2_+; L^2_+)$ is defined by formula \eqref{tildeBudefFoc-}. 

\begin{lem}\label{SemiGrLemtiBL}
 For any $u \in \mathcal{Z}^1_+$ and $\phi \in H^2_+$, we have $\mathbf{L}_u (\phi) \in H^1_+$, $\widetilde{\mathbf{B}}_u   (\phi) \in L^2_+$ and  
 \begin{equation}\label{commForSetasttildBuCv}
 \lim_{\eta \to 0^+} \Big\|\tfrac{i}{\eta} \left[ \mathbf{S}(\eta)^* , \widetilde{\mathbf{B}}_u  \right] \phi - 2 \mathbf{L}_u (\phi) \Big\|_{L^2_+} = 0.
 \end{equation}Moreover, the following uniform estimate  holds:
 \begin{equation}\label{SetasttildBuUniEstL2}
\sup_{\eta >0} \Big\|\tfrac{i}{\eta} \left[ \mathbf{S}(\eta)^* , \widetilde{\mathbf{B}}_u  \right] \phi - 2 \mathbf{L}_u (\phi) \Big\|_{L^2_+} \leq 4 \| \phi \|_{H^2_+} + 6 \| u \|_{L^{\infty}(\mathbb{R} ; \mathbb{C})}^2 \| \phi \|_{L^2_+} .
 \end{equation}
 \end{lem}

\begin{proof}
Since $u \in \mathcal{Z}^1_+$ and $\phi \in H^2_+$,  it follows from formulas 
\eqref{LuSect2}, \eqref{tildeBudefFoc-}, \eqref{comm[S*,D2htog} and  \eqref{L2SetaTudxTvbarphi} that
\begin{equation}\label{SetstTILDBUuPhisplittwo}
\begin{split}
 & \Big\|\tfrac{i}{\eta} \left[ \mathbf{S}(\eta)^* , \widetilde{\mathbf{B}}_u  \right] \phi - 2 \mathbf{L}_u (\phi) \Big\|_{L^2_+}\\
  \leq &  \Big\|\tfrac{1}{\eta} \left[ \mathbf{S}(\eta)^* , D^2 \right] \phi - 2 D \phi \Big\|_{L^2_+}
 + 
2 \Big\| \left( \tfrac{1}{\eta}\left[\mathbf{S}(\eta)^*, \mathbf{T}_u D    \mathbf{T}_{\overline{u}} \right]- \mathbf{T}_u \mathbf{T}_{\overline{u}} \right) \phi \Big\|_{L^2_+ }  \to 0.
\end{split}
\end{equation}Formula \eqref{SetasttildBuUniEstL2} follows from  \eqref{Comm[S*,D2]-2DEsti}, 
\eqref{UNIEstiSetastTuDTvbarL2} 
and 
 \eqref{SetstTILDBUuPhisplittwo}.
 \end{proof}

\begin{cor}\label{CorcommGBu}
For any $u \in \mathcal{Z}^1_+$, we have $\widetilde{\mathbf{B}}_u \left( \mathrm{Dom}^2(\mathbf{G}) \right) \subset \mathrm{Dom} (\mathbf{G})$ and 
\begin{equation}\label{commbfGtildBuFor}
\left[ \mathbf{G}, \widetilde{\mathbf{B}}_u \right] = 2 \mathbf{L}_u : \mathrm{Dom}^2(\mathbf{G}) \to H^1_+.
\end{equation}
\end{cor}
\begin{proof}
Fix $u \in \mathcal{Z}^1_+$ and $ \phi \in \mathrm{Dom}^2(\mathbf{G})$. By Proposition $\ref{GisInfiniGeneHsProp}$ and  formula  \eqref{tildeBudefFoc-}, we have
\begin{equation}\label{tilBuS-1phicv}
\lim_{\eta \to 0^+} \Big\|\widetilde{\mathbf{B}}_u \left(\tfrac{1}{\eta}(\mathbf{S}(\eta)^* \phi - \phi) + i \mathbf{G}\phi \right) \Big\|_{L^2_+} = 0.
\end{equation}Using formulas  \eqref{commForSetasttildBuCv}
and 
\eqref{tilBuS-1phicv}, we obtain that
\begin{equation}\label{Seta-1tildBuphiLim}
\small
\tfrac{1}{\eta} \left(\mathbf{S}(\eta)^* - \mathrm{id}_{L^2_+} \right) \widetilde{\mathbf{B}}_u (\phi)  = \tfrac{1}{\eta}\left[\mathbf{S}(\eta)^* , \widetilde{\mathbf{B}}_u \right] \phi + \widetilde{\mathbf{B}}_u  \left(  \tfrac{1}{\eta}\left(\mathbf{S}(\eta)^* \phi - \phi \right) \right) \to - 2 i \mathbf{L}_u (\phi) - i \widetilde{\mathbf{B}}_u  \mathbf{G} \phi \quad \mathrm{in} \quad L^2_+,
\end{equation}as $\eta \to 0^+$. Formula \eqref{Seta-1tildBuphiLim} yields that $\widetilde{\mathbf{B}}_u (\phi)  \in \mathrm{Dom}(\mathbf{G})$ and $  \mathbf{G} \widetilde{\mathbf{B}}_u (\phi) =  2   \mathbf{L}_u (\phi) +  \widetilde{\mathbf{B}}_u  \mathbf{G} (\phi) \in L^2_+$.
\end{proof}

\begin{rem}\label{remKLV}
  When the symbol $u$ additionally satisfies $u\in H^{\infty}_+$ and
  $\langle x\rangle u\in L^2(\mathbb{R};\mathbb{C})$, the commutator formula
  \eqref{commbfGtildBuFor} was also obtained by Killip--Laurens--Vi\c{s}an
  \cite[Lemma~4.6 on p.~304]{KillLV2311}, under the notational correspondence
  $\mathbf{G}=X$ and $\widetilde{\mathbf{B}}_u=\mathcal{P}(u)$.  
  In their proof, the commutator $[X,\mathcal{P}(u)]$ is taken to act only
  on functions $f\in L^2_+$ with $\langle x\rangle f\in H^2(\mathbb{R};\mathbb{C})$.
    By expanding the derivatives $\partial_x\langle x\rangle$,
$\partial_x^2 \langle x\rangle$, using the Leibniz rule,  formulas \eqref{DomX}  and  \eqref{forDomsXCharbyDsinDX},
  this space coincides with $\mathrm{Dom}^2(\mathbf{X})  =\{\varphi \in \mathrm{Dom}^2 (\mathbf{G}): \widehat{\varphi}(0^+)=0\}$, and is therefore
  strictly contained in $\mathrm{Dom}^2(\mathbf{G})$, according to  \eqref{DomsXCharbyDomsG}. Specifically,  we have 
\begin{equation}\label{DomtwoXbyjbxf}
\mathrm{Dom}^2 (\mathbf{X}) = \{f  \in L^2_+ : \langle x \rangle f \in H^2 (\mathbb{R} ; \mathbb{C})\} \subsetneqq \mathrm{Dom}^2(\mathbf{G}). 
\end{equation}
  In contrast, Corollary~\ref{CorcommGBu} establishes the same commutator
  identity  on the larger domain $\mathrm{Dom}^2(\mathbf{G})$ under the
  sole assumption $u\in\mathcal{Z}^1_+$.  This shows that neither the
  weighted condition $\langle x\rangle f\in H^2 (\mathbb{R} ; \mathbb{C})$ nor the higher regularity with zero boundary condition 
  $u\in H^{\infty}_+$ is intrinsic to the commutator formula \eqref{commbfGtildBuFor}; both are
  artefacts of the generator-based derivation employed in~\cite{KillLV2311}.
  The semigroup approach developed in Lemma $\ref{SemiGrLemtiBL}$ removes these unnecessary hypotheses
  and simultaneously yields the additional information
  $\widetilde{\mathbf{B}}_u(\mathrm{Dom}^2(\mathbf{G}))\subset\mathrm{Dom}(\mathbf{G})$.
\end{rem}

\begin{rem}\label{AdvTSemiEstGenrem}
The commutator identity \eqref{commbfGtildBuFor} involves three unbounded operators $\mathbf{G}, \widetilde{\mathbf{B}}_u, \mathbf{L}_u$ on $L^2_+$, and its precise formulation requires the restriction to $\mathrm{Dom}^2(\mathbf{G})$ in order to ensure that the commutator is well defined.  
However, when one later studies invariant subspaces and derives the explicit formula, the discontinuity of $\mathbf{G}$ (see Lemma $\ref{IGXUnboundLem}$ and formula \eqref{normGfntoinfty}) introduces a serious obstacle: the limit
\begin{equation}
\lim_{r \to 0} \mathbf{G}\Bigl( \tfrac{1}{r} \bigl( \mathbf{f} (t+r) - \mathbf{f} (t) \bigr) \Bigr)
\end{equation} 
may fail to exist even when $\mathbf{f} : \mathbb{R} \to L^2_+$ is differentiable at $t$. In other words, the time derivative does not commute with $\mathbf{G}$.   The commutator estimates \eqref{commForSetasttildBuCv} and \eqref{SetasttildBuUniEstL2} for the contraction semigroup $(\mathbf{S}(\eta)^*)_{\eta\ge 0}$ completely avoid this difficulty.  
  Since each $\mathbf{S}(\eta)^*$ is a bounded linear operator on $L^2_+$, it commutes with the time derivative.
  We can therefore use the characterization \eqref{SemiGrCVdefDomsG} of $\mathrm{Dom} (\mathbf{G}) = \mathrm{Dom}^0(\mathbf{G})$ and pass to the limit $\eta \to 0^+$ to determine which functions remain in $\mathrm{Dom}(\mathbf{G})$ under the time evolution, without ever needing to verify a priori that the test functions belong to the domain of the unbounded generator $\bigl( -i \mathbf{G} , \mathrm{Dom}(\mathbf{G})\bigr)$.  
Moreover, these semigroup estimates provide uniform bounds that are perfectly suited for the dominated convergence arguments used to prove that $\mathrm{Dom}(\mathbf{G})\cap H^2_+$ is invariant under the unitary propagator generated by the Lax pair $(\mathbf{L}_u,\widetilde{\mathbf{B}}_u)$.  
In this sense, \eqref{commForSetasttildBuCv} and \eqref{SetasttildBuUniEstL2} yield more precise information than the generator commutator formula \eqref{commbfGtildBuFor}.
\end{rem}

\noindent For any $\psi= (u,v) \in L^{\infty}(\mathbb{R}; \mathbb{C}^{1 \times 2})$, we set
 \begin{equation}\label{Linfty1times2norm}
 \|\psi\|^2_{L^{\infty}(\mathbb{R}; \mathbb{C}^{1 \times 2})} =  
 \| u \|^2_{L^{\infty}(\mathbb{R}; \mathbb{C})} +  \| v \|^2_{L^{\infty}(\mathbb{R}; \mathbb{C})}
 \end{equation}For any integer $n \geq 1$, if $\psi = (u,v) \in \mathcal{Z}^n_+ \times \mathcal{Z}^n_+$, we define
\begin{equation}\label{ProdHardyZhidNorm}
\|\psi\|_{\mathcal{Z}^n_+ \times \mathcal{Z}^n_+}^2  : = \|u\|_{\mathcal{Z}^n_+ }^2 + \|v\|_{\mathcal{Z}^n_+ }^2.
\end{equation}
   
\begin{defi}
For any $\eta >0$  and $\psi = ( u, v)  \in \mathcal{Z}^1_+ \times  \mathcal{Z}^1_+$,  we define 
\begin{equation}\label{MetapsidefFor}
\mathfrak{M}^{(\eta)}_{\psi} := \eta^{-1}\left[\mathbf{S}(\eta)^*, \mathbf{T}_u D    \mathbf{T}_{\overline{v}} \right]   - \mathbf{T}_u \mathbf{T}_{\overline{v}} : H^1_+ \to L^2_+ .
\end{equation}The mapping $\mathfrak{M}^{(\eta)}$ is defined by $\mathfrak{M}^{(\eta)}  : \psi  \in \mathcal{Z}^1_+ \times  \mathcal{Z}^1_+ \mapsto \mathfrak{M}^{(\eta)}_{\psi} \in \mathcal{B}_{\mathbb{C}}(H^1_+ ; L^2_+)$.
\end{defi}

\noindent For any $\psi= (u,v) \in \mathcal{Z}^1_+ \times \mathcal{Z}^1_+ $,  estimate \eqref{UNIEstiSetastTuDTvbarL2} implies that 
\begin{equation}\label{EstisupMpsietaH1toL2}
\sup_{ \eta >0 } \|\mathfrak{M}^{(\eta)}_{\psi} \|_{\mathcal{B}_{\mathbb{C}}(H^1_+ ; L^2_+)} \leq 3 \| u \|_{L^{\infty}(\mathbb{R}; \mathbb{C})}  \| v \|_{L^{\infty}(\mathbb{R}; \mathbb{C})}.
\end{equation}For any  $\eta >0$, the operator $\mathfrak{M}^{(\eta)}_{\psi} : H^1_+ \to L^2_+ $ can be extended uniquely to a bounded operator on $L^2_+$.  
  The following results establish the uniform equicontinuity of the commutator operators introduced above. This will be essential for the Bochner integrability arguments in the subsequent analysis of the unitary propagator induced by the Lax pair formulation $\bigl(\mathbf{L}_u , \widetilde{\mathbf{B}}_u \bigr)$.
 
\begin{lem}
If $A>0$ and $\psi,\phi \in \mathcal{Z}^1_+ \times \mathcal{Z}^1_+$ such that $ \max\{
\|\psi\|_{\mathcal{Z}^1_+ \times \mathcal{Z}^1_+} , \|\phi\|_{\mathcal{Z}^1_+ \times \mathcal{Z}^1_+}\} \leq A$ ,  then
\begin{equation}\label{LocLipschofMphieta}
\sup_{\eta>0} \bigl\| \mathfrak{M}^{(\eta)}_{\psi} - \mathfrak{M}^{(\eta)}_{\phi} \bigr\|_{\mathcal{B}_{\mathbb{C}}(H^1_+;L^2_+)} \le 6A \; \|\psi - \phi\|_{\mathcal{Z}^1_+ \times \mathcal{Z}^1_+}.
\end{equation}
\end{lem} 
\begin{proof}
Assume that $\eta >0$, $\psi = (u,v), \; \phi = (\mathfrak{a},\mathfrak{b})  \in \mathcal{Z}^1_+ \times \mathcal{Z}^1_+$. For any $\mathbf{h} \in H^1_+$, we have
\begin{equation}\label{estiMpsi-MphihL2}
\begin{split}
& \|\left(\mathfrak{M}^{(\eta)}_{\psi} - \mathfrak{M}^{(\eta)}_{\phi}  \right) \mathbf{h}\|_{L^2_+} \\ \leq 
& \| \eta^{-1} \left[\mathbf{S}(\eta)^* , \mathbf{T}_{u - \mathfrak{a}} \right] \partial_x \mathbf{T}_{\overline{v}}( \mathbf{h})\|_{L^2_+} +
 \| \eta^{-1} \left[\mathbf{S}(\eta)^* , \mathbf{T}_{\mathfrak{a}} \right] \partial_x \mathbf{T}_{\overline{v - \mathfrak{b}}}( \mathbf{h})\|_{L^2_+} \\
 & \quad + \| \mathbf{T}_{u - \mathfrak{a}} \left( \tfrac{1}{\eta}  \left[\mathbf{S}(\eta)^* , D \right] - \mathrm{id}_{L^2_+} \right) \mathbf{T}_{\overline{v}}( \mathbf{h})\|_{L^2_+} 
 + \| \mathbf{T}_{ \mathfrak{a}} \left( \tfrac{1}{\eta}  \left[\mathbf{S}(\eta)^* , D \right] - \mathrm{id}_{L^2_+} \right) \mathbf{T}_{\overline{v - \mathfrak{b}}}( \mathbf{h})\|_{L^2_+} , 
 \end{split}
\end{equation}by formulas \eqref{COMMForeqsetastTuDTvbar} and \eqref{MetapsidefFor}. Plugging  formulas \eqref{UnEsSestTudxTvbphie} and \eqref{ESTITuSetst-1hL2norm}
into  \eqref{estiMpsi-MphihL2},  we obtain that
\begin{equation}\label{EstiMpsi-MphiBH1L2}
\begin{split}
  \| \mathfrak{M}^{(\eta)}_{\psi} - \mathfrak{M}^{(\eta)}_{\phi}   \|_{\mathcal{B}_{\mathbb{C}}(H^1_+ ; L^2_+)} 
  \leq  & 
  3 \| u - \mathfrak{a}\|_{L^{\infty}(\mathbb{R}; \mathbb{C})} \| v \|_{L^{\infty}(\mathbb{R}; \mathbb{C})} + 3 \|   \mathfrak{a}\|_{L^{\infty}(\mathbb{R}; \mathbb{C})} \| v - \mathfrak{b}\|_{L^{\infty}(\mathbb{R}; \mathbb{C})} \\
  \leq & 3 \| \psi -\phi \|_{L^{\infty}(\mathbb{R}; \mathbb{C}^{1 \times 2})} \left( \| \psi   \|_{L^{\infty}(\mathbb{R}; \mathbb{C}^{1 \times 2})} + \| \phi   \|_{L^{\infty}(\mathbb{R}; \mathbb{C}^{1 \times 2})} \right). 
\end{split}
\end{equation}Then formula \eqref{LocLipschofMphieta} follows from \eqref{ZhidNormZn}, 
\eqref{Linfty1times2norm},  \eqref{ProdHardyZhidNorm} and \eqref{EstiMpsi-MphiBH1L2}.
\end{proof}

\begin{defi}\label{DeffbfKeta}
For any $\eta >0$ and $u \in \mathcal{Z}^1_+$, we define
\begin{equation}\label{bbKetaudefFor}
\mathbf{K}^{(\eta)}_u := \tfrac{i}{\eta} \left[ \mathbf{S}(\eta)^* , \widetilde{\mathbf{B}}_u  \right] - 2 \mathbf{L}_u : H^2_+ \to L^2_+.
\end{equation}The mapping $\mathbf{K}^{(\eta)}$ is defined by $\mathbf{K}^{(\eta)} : u \in \mathcal{Z}^1_+ \mapsto \mathbf{K}^{(\eta)}_u \in \mathcal{B}_{\mathbb{C}}(H^2_+ ; L^2_+)$. 
\end{defi}
\noindent Using formulas   \eqref{boldGammaetEstimH2H1}, 
\eqref{tildeBudefFoc-},  \eqref{MetapsidefFor} 
and Definition $\ref{defibolGammSet*,Dsq-2D}$, we obtain that
\begin{equation}\label{bbKetaurelGammafrM}
\mathbf{K}^{(\eta)}_u = \tfrac{1}{\eta} \left[ \mathbf{S}(\eta)^* , D^2 \right] -2 D -2 \left( \tfrac{1}{\eta} \left[\mathbf{S}(\eta)^*, \mathbf{T}_u D    \mathbf{T}_{\overline{u}} \right]   - \mathbf{T}_u \mathbf{T}_{\overline{u}}\right) =  \boldsymbol{\Gamma}^{(\eta)} - 2 \mathfrak{M}^{(\eta)}_{(u,u)} \in \mathcal{B}_{\mathbb{C}}(H^2_+ ; L^2_+) .
\end{equation}By formulas \eqref{boldGammaetEstimH2H1},  \eqref{EstisupMpsietaH1toL2}
and \eqref{bbKetaurelGammafrM}, it follows that
\begin{equation}\label{UniEstbbKueta>0}
\sup_{\eta >0} \|\mathbf{K}^{(\eta)}_u \|_{\mathcal{B}_{\mathbb{C}}(H^2_+ ; L^2_+)} \leq 4 + 6 \|u\|_{L^{\infty}(\mathbb{R} ; \mathbb{C})}^2 .
\end{equation}

\begin{lem}
For any $A>0$ and $u, v \in \mathcal{Z}^1_+$ such that $\max\{\|u\|_{\mathcal{Z}^1_+}, \|v\|_{\mathcal{Z}^1_+}\} \leq A$, we have
\begin{equation}\label{UniEquiContKeta}
\sup_{\eta > 0 } \|\mathbf{K}^{(\eta)}_u - \mathbf{K}^{(\eta)}_v \|_{\mathcal{B}_{\mathbb{C}}(H^2_+ ; L^2_+)} \leq 24 A \|u - v\|_{\mathcal{Z}^1_+}.
\end{equation}
\end{lem}
\begin{proof}
Assume that $\eta >0$ and $u,v \in \mathcal{Z}^1_+$. Formula \eqref{ProdHardyZhidNorm} yields that $\|(u,u)\|_{\mathcal{Z}^1_+ \times \mathcal{Z}^1_+}=\sqrt{2}\|u\|_{\mathcal{Z}^1_+}$ and $\|(v , v )\|_{\mathcal{Z}^1_+ \times \mathcal{Z}^1_+}=\sqrt{2}\|v\|_{\mathcal{Z}^1_+}$. 
Using formulas \eqref{LocLipschofMphieta}  and \eqref{bbKetaurelGammafrM},  we obtain that $\|\mathbf{K}^{(\eta)}_u - \mathbf{K}^{(\eta)}_v \|_{\mathcal{B}_{\mathbb{C}}(H^2_+ ; L^2_+)} \leq 2  \|\mathfrak{M}^{(\eta)}_{(u,u)} - \mathfrak{M}^{(\eta)}_{(v , v)} \|_{\mathcal{B}_{\mathbb{C}}(H^1_+ ; L^2_+)} \leq 12  \|(u,u)-(v,v)\|_{\mathcal{Z}^1_+ \times \mathcal{Z}^1_+}   \max\{\|(u,u)\|_{\mathcal{Z}^1_+ \times \mathcal{Z}^1_+}, \|(v , v )\|_{\mathcal{Z}^1_+ \times \mathcal{Z}^1_+}\}$. So we obtain \eqref{UniEquiContKeta},    for any $A \geq \max\{\|u\|_{\mathcal{Z}^1_+}, \|v\|_{\mathcal{Z}^1_+}\}$.
\end{proof}

\begin{cor}\label{CorbbKcircuuniEquiCont}
Assume that $a <b$ and $u \in \mathcal{C} \left( (a, b) ; \mathcal{Z}^1_+ \right)$. For any $\eta >0$, we define the mapping 
\begin{equation}\label{CompKetaudefFor}
\mathbf{K}^{(\eta)} \circ u : t \in (a,b ) \mapsto \left( \mathbf{K}^{(\eta)}\circ u  \right) (t) = \mathbf{K}^{(\eta)}_{u(t)} \in \mathcal{B}_{\mathbb{C}}(H^2_+; L^2_+).
\end{equation}Then the family $\left(\mathbf{K}^{(\eta)}\circ u \right)_{\eta >0}$ is uniformly equicontinuous on every compact subset of $(a,b)$. 
\end{cor}
\begin{proof}
Fix $\alpha , \beta \in (a, b)$ such that $\alpha < \beta$. Set $A: = \max\{\|u(t)\|_{\mathcal{Z}^1_+} : \alpha \leq t \leq \beta \}>0$. Then we have
\begin{equation}\label{CompKetaututau24A}
\sup_{\eta >0 } \| \mathbf{K}^{(\eta)} \circ u (t) - \mathbf{K}^{(\eta)} \circ u (\tau )\|_{\mathcal{B}_{\mathbb{C}} (H^2_+ ; L^2_+)} \leq 24 A \| u (t) - u(\tau)  \|_{\mathcal{Z}^1_+}, \quad \forall t, \tau \in [\alpha, \beta],
\end{equation}by   \eqref{UniEquiContKeta} and \eqref{CompKetaudefFor}.
  The proof is completed by using the continuity of $u : (a,b) \to \mathcal{Z}^1_+$ and   \eqref{CompKetaututau24A}.
\end{proof}

\subsection{Spectral analysis of the Chiral Lax Operator}\label{SubSspanal3.5}

\noindent We close this section by describing the essential spectrum and the simplicity of eigenvalues of the chiral Lax operator $\mathbf{L}_u$ when $u = e^{i \alpha} + v$ for some $\alpha \in \mathbb{T} = \mathbb{R}/ 2\pi \mathbb{Z} $ and $v \in H^1_+$. These results will provide essential tools for investigating the dark multi-soliton potentials in Section $\ref{DarMulSoliPotSect}$. \\

\noindent Lemma $\ref{DNSelfadjinHsLem}$ shows that $D = -i \partial_x$ is an  unbounded self-adjoint operator on $L^2_+$ with domain $H^1_+$. For any  $v \in L^2(\mathbb{R} ; \mathbb{C})$, the densely defined Toeplitz operator $\mathbf{T}_v : H^1_+ \to L^2_+$ given by \eqref{ToepdefForSec2} is not necessarily bounded. However,  it is relatively compact with respect to   $D : H^1_+ \to L^2_+$, as shown in the next lemma.
  
\begin{lem}\label{bfTvrelcompwrtDlem}
If  $r \in \mathbb{R}$ and  $v \in L^2(\mathbb{R} ; \mathbb{C})$, then $\mathbf{T}_v \bigl( D + r +i \bigr)^{-1}$ is a Hilbert--Schmidt operator on $L^2_+$.
\end{lem}
\begin{proof} 
For any $r \in \mathbb{R}$, using the first resolvent formula,  we obtain that	
\begin{equation}\label{firstResolFor}
(D +  r + i)^{-1} =  (D +   i)^{-1} \bigl( \mathrm{id}_{L^2_+} - r (D +  r + i)^{-1} \bigr) \in \mathcal{B}_{\mathbb{C}} (L^2_+) .
\end{equation}Recall that $\mathbb{R}_+ = [0 , +\infty)$. Fix $v \in L^2 (\mathbb{R} ; \mathbb{C})$. For any $f \in L^2 (\mathbb{R} ; \mathbb{C})$, we define
\begin{equation}\label{bbKkernelOp}
\mathbb{K}_v (f) (\xi) : = \int_{\mathbb{R}} K_v (\xi, \eta) f(\eta) \mathrm{d} \eta 
\quad \text{with} \quad 
K_v (\xi , \eta) = \mathbf{1}_{\mathbb{R}_+}(\xi) \mathbf{1}_{\mathbb{R}_+}( \eta )  \tfrac{\widehat{v}(\xi -\eta) }{\eta + i}, 
\end{equation}for any $\xi, \eta \in \mathbb{R}$. Let $\mathscr{F} : \mathscr{S}' \to  \mathscr{S}'$ denote the Fourier transform in \eqref{FouTrDefS2}. Then we have
\begin{equation}\label{bbKvisHSonHardy}
\mathbf{T}_v (D + i)^{-1} = \tfrac{1}{2 \pi} \mathscr{F}^{-1} \mathbb{K}_v \mathscr{F} :  L^2_+  \to L^2_+ , \quad \|K_v\|_{L^2 (\mathbb{R} \times 
\mathbb{R}  ; \mathbb{C} )}  \leq \pi \|v\|_{L^2 (\mathbb{R} ; \mathbb{C})} ,
\end{equation}by 
 Proposition 2.3 of \cite{SUNCMP2021}. So $\mathbf{T}_v (D + i)^{-1}$ is a Hilbert--Schmidt operator on $L^2_+$. The proof is completed by formulas  \eqref{bbKvisHSonHardy} and \eqref{firstResolFor}.
\end{proof}

\noindent The following lemma identifies the essential spectrum of the chiral Lax operator.
 
\begin{lem}\label{EssSpebfLulem}
Assume that $u = e^{i \alpha} + v$ for some $\alpha \in \mathbb{T} $ and $v \in H^1_+$.
 Let  $\mathbf{L}_u : H^1_+ \to L^2_+$ be defined in formula \eqref{LuSect2}.  
  Its essential spectrum is given by $\sigma_{\mathrm{ess}} (\mathbf{L}_u) = [ - 1 , + \infty)$. 
\end{lem}

\begin{proof}
Proposition 2.3 of \cite{SUNCMP2021} gives the spectrum of the positive operator $D : H^1_+ \to L^2_+$: 
\begin{equation}\label{sigmaDfor}
\sigma (D) = \sigma_{\mathrm{ess}}(D) = \sigma_{\mathrm{ac}}(D) = [0 , +\infty) .
\end{equation}
The proof follows directly from formula 
\eqref{TuTubaris1andTv}, Lemma $\ref{bfTvrelcompwrtDlem}$ and Weyl's essential spectrum theorem.
\end{proof}

 \noindent After establishing the essential spectrum, we now turn to the study of    discrete spectrum of   $\mathbf{L}_u$. The following auxiliary identity will play a central role in  establishing the simplicity of its eigenvalues.
   
\begin{lem}\label{SimIdhvhlem}
Let $\mathrm{Dom} (\mathbf{I}_+)$ be defined in \eqref{DomI+}. 
  Assume that $u = e^{i \alpha } + v$ for some 
 $\alpha \in \mathbb{T} $ and $v \in H^1_+$. If there exist $\lambda \in \mathbb{R}$ and $h \in \mathrm{Dom} (\mathbf{I}_+) \cap  H^1_+$ such that $\mathbf{L}_u (h) = \lambda h$, then the following identity holds:  
 \begin{equation}\label{idenhL2nvhinprod}
2 \pi \| h \|_{L^2_+}^2 = e^{-i \alpha }\widehat{h}(0^+)   \langle v , h \rangle_{L^2_+} + |\langle v , h \rangle_{L^2_+}|^2.
 \end{equation}
\end{lem}
\begin{proof}
Since $\mathbf{L}_u $ is an unbounded self-adjoint operator on $L^2_+$, and $\mathbf{L}_u (h)=\lambda h \in H^1_+$, we have 
\begin{equation}\label{CommForSstLhinproh}
\langle \left[\mathbf{S}(\eta)^*, \mathbf{L}_u \right] h, h \rangle_{L^2_+} = \langle  \mathbf{S}(\eta)^* (\mathbf{L}_u  - \lambda) h, h \rangle_{L^2_+} - \langle  \mathbf{S}(\eta)^* (h), (\mathbf{L}_u  - \lambda) h \rangle_{L^2_+} = 0.
\end{equation}Since  we assume that $h \in \mathrm{Dom}(\mathbf{I}_+)$ in addition, formulas  \eqref{ForcommS*,L]phi} and \eqref{CommForSstLhinproh}   imply  that
\begin{equation*}
 0 = \lim_{\eta \to 0^+}\langle \tfrac{1}{\eta}\left[\mathbf{S}(\eta)^*, \mathbf{L}_u \right] h, h \rangle_{L^2_+} = \langle h- \left(e^{-i \alpha} \widehat{h}(0^+) + \langle h,v \rangle_{L^2_+} \right) \tfrac{v}{2\pi} , h \rangle_{L^2_+},
\end{equation*}which is equivalent to formula 
\eqref{idenhL2nvhinprod}.
\end{proof}

\noindent For any $\lambda \in \mathbb{R} \setminus \{-1\}$, the next lemma shows that the $\lambda$-eigenspace of $\mathbf{L}_u$ is contained in $\mathrm{Dom} (\mathbf{I}_+)$.

\begin{lem}\label{EigLuhperpvLem}
For any $\alpha \in \mathbb{T}$ and $v \in H^1_+$, we define $u = e^{i \alpha } + v \in \mathcal{Z}^1_+$.  Assume that there exist $\lambda \in \mathbb{R}$ and $h \in   H^1_+$ such that $\mathbf{L}_u (h) = \lambda h \in H^1_+$, then we have 
\begin{equation}\label{Limlapls1mixihath}
\lim_{\xi \to 0^+} (\lambda + 1 -\xi) \widehat{h}(\xi) = - e^{i \alpha} \langle h, v \rangle_{L^2_+}.
\end{equation}In addition, if $\lambda=-1$, we have $\langle h, v \rangle_{L^2_+} = 0$. Consequently, it follows that
\begin{equation}\label{KerLuadd1DI0}
\mathrm{Ker} \bigl( \mathbf{L}_u + 1 \bigr) \cap  \mathrm{Dom} (\mathbf{I}_+) = \{0\} .
\end{equation}
\end{lem}

\begin{proof}
Since $v, h \in H^1_+$, we have $vh, v \Pi(\overline{v}h) \in H^1_+ \cap  L^1(\mathbb{R}; \mathbb{C})$ and $\overline{v}h \in H^1 \cap  L^1(\mathbb{R}; \mathbb{C})$ by \eqref{HnZnsubsetHn} and \eqref{ZnHardBanAlgebra}. Then $\widehat{v h}(0)= \left( v \Pi(\overline{v}h)\right)^{\wedge}(0)=0$ and $\widehat{\overline{v}h}(0)= \langle h, v \rangle_{L^2_+}$. For any $\xi \in \mathbb{R}$, we define 
\begin{equation}\label{betaxidef}
\boldsymbol{\beta}(\xi) := \left(e^{-i \alpha} vh + e^{i\alpha} \Pi(\overline{v}h) + v \Pi(\overline{v}h) \right)^{\wedge}(\xi).
\end{equation}Then $\boldsymbol{\beta} \in \mathcal{C}(\mathbb{R}^*; \mathbb{C}) \cap  L^{\infty}(\mathbb{R}; \mathbb{C})$ and $\mathscr{F}^{-1} \boldsymbol{\beta} \in H^1_+$, where $\mathbb{R}^* = \mathbb{R}\setminus \{0\}$ and $\mathscr{F} : \Lambda \in  \mathscr{S}' \mapsto \mathscr{F}(\Lambda):=\widehat{\Lambda} \in  \mathscr{S}'$ denotes the Fourier transform. Besides, we have    $\mathscr{F}^{-1} \boldsymbol{\beta} \in \mathrm{Dom}(\mathbf{I}_+)$ and
\begin{equation}\label{beta0+=<h,v>}
\boldsymbol{\beta} (0^+) = \lim_{\xi \to 0^+} \boldsymbol{\beta} (\xi)  = e^{i \alpha} \langle h, v \rangle_{L^2_+}.
\end{equation}Since $\mathbf{L}_u (h) = \lambda h \in H^1_+$ for some $\lambda \in \mathbb{R}$, formula \eqref{betaxidef} yields that
\begin{equation}\label{DhbetahRel}
Dh = (\lambda +1)h + e^{-i \alpha} vh + e^{i\alpha}  \Pi(\overline{v}h) + v \Pi(\overline{v}h) = (\lambda +1)h + \mathscr{F}^{-1} \boldsymbol{\beta} \in H^1_+.
\end{equation}Then formula \eqref{Limlapls1mixihath} is deduced from \eqref{beta0+=<h,v>} and \eqref{DhbetahRel}.\\

\noindent Assume that $\lambda = -1$. Formula  \eqref{DhbetahRel} yields that $\boldsymbol{\beta} (\xi)= \xi \widehat{h}(\xi)$. We have $\boldsymbol{\beta}, \widehat{h} \in \mathcal{C}(\mathbb{R}^*;\mathbb{C}) \cap  L^2(\mathbb{R}; \mathbb{C})$. If $\rho:=|  \boldsymbol{\beta}(0^+)|>0$, then there exists $\delta>0$ such that $\inf\{|\boldsymbol{\beta}(\xi)|:0<\xi \leq \delta\} \geq \frac{1}{2}\rho$. Then
\begin{equation*}
\int_{\mathbb{R}}|\widehat{h}|^2 \geq \int_0^{\delta} \frac{\rho^2}{4 \xi^2} \mathrm{d}\xi=+\infty,
\end{equation*}which contradicts the assumption that $h \in H^1_+$. So we have $\langle h , v \rangle_{L^2_+} = e^{-i \alpha}\boldsymbol{\beta} (0^+)=0$.
  Finally, formula \eqref{KerLuadd1DI0} is a direct consequence of identity \eqref{idenhL2nvhinprod}. 
\end{proof}
 
\begin{cor}\label{CorSimpofEigLu}
If $u = e^{i \alpha } + v$ for some 
 $\alpha \in \mathbb{T}$ and $v \in H^1_+$, for any $\lambda \in \sigma_{\mathrm{pp}}(\mathbf{L}_u) \setminus \{-1\}$, we have $\mathrm{Ker}(\lambda - \mathbf{L}_u) \subset \mathrm{Dom}(\mathbf{I}_+) \cap  H^1_+$ and $\dim_{\mathbb{C}}\mathrm{Ker}(\lambda - \mathbf{L}_u)  = 1$.
\end{cor}
\begin{proof}
Assume that $\lambda \in \sigma_{\mathrm{pp}}(\mathbf{L}_u) \setminus \{-1\}$, set $\mathbb{V}_{\lambda}:= \mathrm{Ker}(\lambda - \mathbf{L}_u) \subset  H^1_+$. Then $\dim_{\mathbb{C}} \mathbb{V}_{\lambda} \geq 1$. For any $h \in \mathbb{V}_{\lambda}$, we have $\mathbf{L}_u (h) = \lambda h$ and $\lim_{\xi \to 0^+} \widehat{h}(\xi) = -\frac{e^{i \alpha}}{\lambda+1} \langle h, v \rangle_{L^2_+}$ by using  \eqref{Limlapls1mixihath}. So $h \in \mathrm{Dom}(\mathbf{I}_+)$. We define a linear form $\mathfrak{j}_v : \varphi \in \mathbb{V}_{\lambda} \mapsto \mathfrak{j}_v (\varphi) := \langle \varphi, v \rangle_{L^2_+} \in \mathbb{C}$. If $\phi \in \mathrm{Ker}(\mathfrak{j}_v) \subset \mathbb{V}_{\lambda} \subset \mathrm{Dom}(\mathbf{I}_+) \cap  H^1_+$, then $\langle \phi, v \rangle_{L^2_+} = 0$ and formula  \eqref{idenhL2nvhinprod} yields that $2 \pi \|\phi \|_{L^2_+}^2 = e^{-i \alpha } \widehat{\phi}(0^+) \langle v,  \phi  \rangle_{L^2_+} +|\langle v ,  \phi\rangle_{L^2_+} |^2 =0$. So $\phi =0$ and the linear form $\mathfrak{j}_v: \mathbb{V}_{\lambda} \to \mathbb{C}$ is injective. Since $\dim_{\mathbb{C}} \mathbb{V}_{\lambda} \geq 1 = \dim_{\mathbb{C}} \mathbb{C}$, we have $\dim_{\mathbb{C}}\mathrm{Ker}(\lambda - \mathbf{L}_u) = \dim_{\mathbb{C}} \mathbb{V}_{\lambda}  = 1$, for any $\lambda \in \sigma_{\mathrm{pp}}(\mathbf{L}_u) \setminus \{-1\}$.
\end{proof}

\section{Well‑Posedness and the Explicit Formula}\label{LWPExpForSect4}

\noindent In this section, we first show the local well-posedness of the focusing CMdNLS equation \eqref{CMdNLSfocu} on the Zhidkov space $\mathcal{Z}^n$ (defined in \eqref{ZhiSpdefSect2}) and on the Hardy--Zhidkov space $\mathcal{Z}^n_+$ (defined in \eqref{HarZhidSpZnDefFor}), for $n \in \{2, 3\}$. 
We then derive the explicit formula for its general solutions in the Hardy--Zhidkov space $\mathcal{Z}^2_+$ via the unitary equivalence property induced by the chiral Lax pair formulation $(\mathbf{L}_u, \mathbf{B}_u)$ given by \eqref{LuSect2} and \eqref{BudefForFoc-}. 
Finally, we introduce several invariant subspaces and conservation laws of equation \eqref{CMdNLSfocu}.

\subsection{Local Well‑Posedness}

\noindent The local well-posedness of equation \eqref{CMdNLSfocu} is established by Kato's iterative scheme presented in Subsection 3.1 of \cite{Xchenccm25ExpFor}. The focusing CMdNLS equation \eqref{CMdNLSwithPi} can be written as 
\begin{equation}\label{withReCMdNLS}
\partial_t u = i \partial_x^2 u +  4 u \Pi_{\geq 0} \bigl( \operatorname{Re} (\overline{u} \partial_x u ) \bigr), \qquad  (t,x) \in \mathbb{R} \times \mathbb{R}.
\end{equation}Using the idea of formula (3.2) of \cite{Xchenccm25ExpFor}, we adopt the following Kato's iterative scheme,
\begin{equation}\label{KatoIteSchPDE}
\partial_t u^{k+1} = i \partial_x^2 u^{k+1}  +  4 u^k \Pi_{\geq 0} \bigl( \operatorname{Re} (\overline{u^k} \partial_x u^{k+1}  ) \bigr) , 
\end{equation}with $u^{k+1} (0) = u_0 \in \mathcal{Z}^2$, $\forall k \in \mathbb{N}_+$. 
  The following lemma gives  the well-posedness of equation \eqref{KatoIteSchPDE}. 
   
\begin{lem}\label{kapLemLiPart}
Assume that $u \in \mathcal{C} \bigl( (- \mathsf{T}, \mathsf{T}) ; \mathcal{Z}^2 \bigr)$ for some $\mathsf{T} >0$. For any $\varkappa_0 \in \mathcal{Z}^2$, there exists a unique function $\varkappa \in \mathcal{C} \bigl( (- \mathsf{T}, \mathsf{T}) ; \mathcal{Z}^2 \bigr)$ such that  $\varkappa(0) = \varkappa_0$ and 
\begin{equation}\label{vkappaPDEdef}
\partial_t \varkappa = i \partial_x^2 \varkappa + 4 u \Pi_{\geq 0} \operatorname{Re} \bigl( (\overline{u} \partial_x \varkappa) \bigr)  , \qquad \forall t \in (- \mathsf{T}, \mathsf{T}) .
\end{equation}Moreover, we set $ \mathsf{v} (t) : = \varkappa(t) - e^{it \partial_x^2} \varkappa_0$. Then $\mathsf{v} \in \mathcal{C} \bigl( (- \mathsf{T}, \mathsf{T}); H^2 (\mathbb{R} ; \mathbb{C}) \bigr) \cap \mathcal{C}^1 \bigl( (- \mathsf{T}, \mathsf{T}); L^2 (\mathbb{R} ; \mathbb{C})  \bigr)$ and there exists a universal constant  $\mathrm{C} >0$ such that the following estimate holds:
\begin{equation}\label{UnEstsfvTH2}
\sup_{|t| \leq \mathsf{T}} \|\mathsf{v} (t)\|_{H^2 (\mathbb{R} ; \mathbb{C})}  \leq \mathrm{C} \|\partial_x \varkappa_0\|_{H^1 (\mathbb{R} ; \mathbb{C})}  \Bigl( 1 + \int_{-\mathsf{T}}^{\mathsf{T}} \|u( \mathfrak{s})\|_{\mathcal{Z}^2}^2 \mathrm{d} \mathfrak{s} \Bigr) e^{\mathrm{C} \int_{-\mathsf{T}}^{\mathsf{T}} \|u( \mathfrak{s})\|_{\mathcal{Z}^2}^2 \mathrm{d} \mathfrak{s} }.
\end{equation}In addition, if $u \in \mathcal{C} \bigl( (- \mathsf{T}, \mathsf{T}) ; \mathcal{Z}^2_+ \bigr)$ and $\varkappa_0 \in \mathcal{Z}^2_+$, then the  solution $\varkappa$ given by \eqref{vkappaPDEdef} and $\mathsf{v}$ satisfy that
\begin{equation}\label{sfvinCtHx2}
\varkappa \in \mathcal{C} \bigl( (- \mathsf{T}, \mathsf{T}) ; \mathcal{Z}^2_+  \bigr)
 \qquad \text{and}  \qquad
 \mathsf{v} \in \mathcal{C} \bigl( (- \mathsf{T}, \mathsf{T}); H^2_+ \bigr) \cap \mathcal{C}^1 \bigl( (- \mathsf{T}, \mathsf{T}); L^2_+ \bigr) .
\end{equation}
\end{lem}
 
\noindent  Lemma~\ref{kapLemLiPart} is proved by the same energy method as in Lemma~3.1 of \cite{Xchenccm25ExpFor}. The only difference is the sign of the nonlinear term.
  All estimates in Subsection 3.1 of \cite{Xchenccm25ExpFor}   rely only on the Zhidkov space norm $\|\cdot\|_{\mathcal{Z}^2}$, and the condition $|u|^2-1 \in L^2 (\mathbb{R} ; \mathbb{C})$ appearing in the definition of the space $E$ in \cite{Xchenccm25ExpFor} (i.e. the metric space $\bigl( \mathcal{M}^2 , \mathsf{d}_{\mathcal{M}^2}  \bigr)$ given by \eqref{calMnmetrSp} of this paper) is not used in the local well-posedness argument.
 
\begin{proof}[Proof of Proposition $\ref{LWPpropCMdNLSZn}$]
The proof follows the Kato iterative scheme \eqref{KatoIteSchPDE} used in \cite[Proposition~3.2]{Xchenccm25ExpFor} for the defocusing CMdNLS equation \eqref{defocuCMdNLS}. For any $t \in (- \mathsf{T} , \mathsf{T})$, we have 
\begin{equation}\label{sfvDefFor}
\mathsf{v} (t) =\mathsf{v}_{u, \varkappa} (t) = \int_0^t e^{ i (t-s) \partial_x^2} \bigl(4 u \Pi_{\geq 0}  \operatorname{Re} ( \overline{u} \partial_x \varkappa )\bigr) (\mathfrak{s}) \mathrm{d} \mathfrak{s} \in H^2 (\mathbb{R} ; \mathbb{C}) .
\end{equation} 
 The presence of the Hardy--Zhidkov space is handled by observing that the nonlinearity preserves the positive-frequency condition, thanks to the ideal property of $H^n_+$ in $\mathcal{Z}^n_+$, $\forall 
n \in \{2 , 3 \}$. (See also Lemma $\ref{HarZhidisBanachLem}$). 
We therefore omit the details. 
\end{proof}
 
\noindent  Assume that $n \in \{2, 3\}$ and $T^- < 0 < T^+$.  Let $u \in \mathcal{C}\bigl((T^-,T^+);\mathcal{Z}^n \bigr)$ be a fixed solution to the focusing CMdNLS equation \eqref{CMdNLSfocu} with initial datum $u(0) = u_0 \in \mathcal{Z}^n$. 
  We define
\begin{equation}\label{LinParturmL}
u_{\mathrm{L}}(t) := e^{it \partial_x^2} u_0 \in \mathcal{Z}^n, \qquad \forall t \in (T^- , T^+).
\end{equation}By formulas 
 \eqref{sfvDefFor} and \eqref{LinParturmL}, we have the  Duhamel formula for the solution $u$: 
\begin{equation}\label{Duhamelufor}
u (t) = e^{it \partial_x^2} u_0  +   \int_0^t e^{ i (t-s) \partial_x^2} \bigl(4 u \Pi_{\geq 0}  \operatorname{Re} ( \overline{u} \partial_x  u)\bigr) (\mathfrak{s}) \mathrm{d} \mathfrak{s} \in \mathcal{Z}^n,  \qquad \forall t \in (T^- , T^+).
\end{equation}
  The following lemma  is a direct consequence  of Lemma~\ref{kapLemLiPart}.

\begin{lem}\label{lemLinPartsol} 
We have $u - u_{\mathrm{L}} \in \mathcal{C}\bigl((T^-,T^+) ; H^n (\mathbb{R} ; \mathbb{C}) \bigr) \cap \mathcal{C}^1 \bigl((T^-,T^+) ; H^{n-2}  (\mathbb{R} ; \mathbb{C}) \bigr)$ and 
\begin{equation}\label{soluqudifCVSob}
\lim_{r \to 0} \big\| \tfrac{1}{r}  \bigl(u(t+r) - u(t) \bigr) - \bigl( i \partial_x^2 u + 2 u \Pi_{\geq 0} \partial_x  (|u|^2  ) \bigr) (t)  \big\|_{H^{n-2}(\mathbb{R} ; \mathbb{C})} = 0.
\end{equation}
  In addition, if $u_0 \in \mathcal{Z}^3$, we have $u \in \mathcal{C}\bigl((T^-,T^+);\mathcal{Z}^3 \bigr) \cap  \mathcal{C}^1 \bigl((T^-,T^+);\mathcal{Z}^1 \bigr)$.
\end{lem}

\begin{proof}
Set $\mathsf{F} = 4 u \Pi_{\geq 0} \operatorname{Re} (\overline{u} \partial_x u) \in \mathcal{C} \bigl((T^-,T^+) ; H^{n-1} (\mathbb{R} ; \mathbb{C}) \bigr)$.  
For any $t \in (T^-,T^+)$, we set 
\begin{equation}\label{sfvuuinHn}
\mathsf{v}_{u, u}  (t) = u (t) - u_{\mathrm{L}} (t)  = \int_0^t e^{ i (t-s) \partial_x^2} \bigl(4 u \Pi_{\geq 0}  \operatorname{Re} ( \overline{u} \partial_x  u )\bigr) (\mathfrak{s}) \mathrm{d} \mathfrak{s} \in H^2 \cap \mathcal{Z}^n =  H^n (\mathbb{R} ; \mathbb{C}) .
\end{equation}By 
Lemma $\ref{GenSchGrlemZN}$, formula \eqref{sfvuuinHn} and  the dominated convergence theorem, we have
\begin{equation}\label{vuudifQuocvfor}
\lim_{r \to 0} \|\tfrac{1}{r} \bigl( \mathsf{v}_{u, u}  (t + r ) - \mathsf{v}_{u, u}  (t)\bigr) -  \bigl( i \partial_x^2 \mathsf{v}_{u, u} + \mathsf{F} \bigr) (t)\|_{H^{n-2}(\mathbb{R} ; \mathbb{C})} =0 .
\end{equation}Then   \eqref{soluqudifCVSob} follows directly from \eqref{vuudifQuocvfor}, \eqref{sfvuuinHn} and  \eqref{CVHNeitLapphiphi}.   The proof is completed  by  \eqref{HNhookarrZN}.
\end{proof}

\noindent The next lemma handles the important case where the initial datum is a constant plus a Sobolev perturbation, which is the natural setting for dark multi-solitons.

\begin{lem}\label{lemuminainCtHnx}
If $a \in \mathbb{C}$ and $u_0 \in a + H^n(\mathbb{R} ; \mathbb{C})$ with $n \in \{2,3\}$, then the solution $u$ satisfies that $u - a \in \mathcal{C}\bigl((T^-,T^+) ; H^n (\mathbb{R} ; \mathbb{C}) \bigr) \cap \mathcal{C}^1 \bigl((T^-,T^+) ; H^{n-2}  (\mathbb{R} ; \mathbb{C}) \bigr)$.
\end{lem}

\begin{proof}
Since $u_0 - a \in H^n(\mathbb{R} ; \mathbb{C})$, we have $u_{\mathrm{L}}(t) - a = e^{it \partial_x^2} (u_0 - a) \in H^n(\mathbb{R} ; \mathbb{C})$ by formula \eqref{SchGreitLap1=1}, for any $t \in (T^-,T^+)$. Then $u_{\mathrm{L}} - a \in \mathcal{C}\bigl(\mathbb{R} ; H^n(\mathbb{R} ; \mathbb{C}) \bigr) \cap \mathcal{C}^1 \bigl(\mathbb{R} ; H^{n-2}(\mathbb{R} ; \mathbb{C}) \bigr)$ by the Stone theorem. The proof is completed by using Lemma~\ref{lemLinPartsol}.
\end{proof}

\noindent Then we show that 
  $\bigl(\mathcal{M}^n , \mathsf{d}_{\mathcal{M}^n} \bigr)$ defined by \eqref{calMnmetrSp} and  \eqref{sfdcalMnmetricDef} is invariant under the flow of \eqref{CMdNLSfocu}.

\begin{lem}\label{lemCMdNLSMnLWP}
Fix $n \in \{2,3\}$. If the initial datum $u(0) = u_0 \in \mathcal{M}^n$, then the solution $u$ satisfies that  $u \in \mathcal{C}\bigl( (T^-, T^+) ; \mathcal{M}^n \bigr)$. 
\end{lem}

\begin{proof}
 Since $\mathsf{v}_{u,u} = u - u_{\mathrm{L}}   \in \mathcal{C}\bigl((T^-,T^+) ; H^n (\mathbb{R} ; \mathbb{C}) \bigr) \cap \mathcal{C}^1 \bigl((T^-,T^+) ; H^{n-2}  (\mathbb{R} ; \mathbb{C}) \bigr)$, we have
\begin{equation}\label{uttoufrtinL2}
 \| u (t) - u (\mathfrak{t})\|_{L^2 (\mathbb{R} ; \mathbb{C})} \leq (2 + \sqrt{2}) \sqrt{|t  - \mathfrak{t} |} \|\partial_x u_0\|_{L^2 (\mathbb{R} ; \mathbb{C})} + \| \mathsf{v}_{u,u} (t) - \mathsf{v}_{u,u} (\mathfrak{t})\|_{L^2 (\mathbb{R} ; \mathbb{C})} \to 0,
\end{equation}as $\mathfrak{t}  \to t$, by formulas     \eqref{LinParturmL} and   \eqref{estiL2EtDuubyu},   for any $t\in (T^-, T^+)$. It follows from  \eqref{sfdcalMnmetricDef} and \eqref{uttoufrtinL2} that 
\begin{equation*}
\mathsf{d}_{\mathcal{M}^n} \bigl( u (t) , u (\mathfrak{t})\bigr)  \leq \|u (t) - u (\mathfrak{t})\|_{\mathcal{Z}^n}  + \bigl( \| u (t) \|_{L^2 }  + \| u (\mathfrak{t})\|_{L^2  }  \bigr) \| u (t) - u (\mathfrak{t})\|_{L^2 (\mathbb{R} ; \mathbb{C})}  \to 0, \quad \text{as} \quad \mathfrak{t}  \to t.
\end{equation*} 
  So $u : t \in (T^-, T^+) \mapsto u(t) \in \bigl(\mathcal{M}^n , \mathsf{d}_{\mathcal{M}^n} \bigr)$ is continuous by  \eqref{sfvuuinHn},  \eqref{LinParturmL} and  Lemma $\ref{MnaddHninMnlem}$.
\end{proof}

\noindent For any $u \in \mathcal{M}^2$, let $\mathsf{w}^{(1)}_u \in H^1 (\mathbb{R} ; \mathbb{C})$ be defined in Lemma $\ref{sfwu1defLem}$ and  $\mathsf{w}^{(2)}_u \in L^2 (\mathbb{R} ; \mathbb{C})$ is given by formula \eqref{forsfwu2def}. We define 
$\mathfrak{W}^{(1)}(u) : = \|\mathsf{w}^{(1)}_u\|_{L^2(\mathbb{R} ; \mathbb{C})}$ and $\mathfrak{W}^{(2)}(u) : = \|\mathsf{w}^{(2)}_u\|_{L^2(\mathbb{R} ; \mathbb{C})}$. Specifically, 
\begin{equation}\label{consvLawfrW12def}
\small
\begin{split}
& \mathfrak{W}^{(1)}(u) = \|\mathsf{w}^{(1)}_u\|_{L^2(\mathbb{R} ; \mathbb{C})} =   \|D u - u \Pi_{\geq 0} (|u|^2 -1) \|_{L^2 (\mathbb{R} ; \mathbb{C})}, \\
& \mathfrak{W}^{(2)}(u) = \| \mathfrak{L}_u (\mathsf{w}^{(1)}_u) - \mathsf{w}^{(1)}_u  \|_{L^2 (\mathbb{R} ; \mathbb{C})}  \\ = & \| D^2 u - D \bigl( u \Pi_{\geq 0} (|u|^2 -1) \bigr)  - Du - u \Pi_{\geq 0} (\overline{u} D u) + u \Pi_{\geq 0} \bigl( |u|^2 \Pi_{\geq 0} (|u|^2 -1)  \bigr) + u \Pi_{\geq 0} (|u|^2 -1) \|_{L^2 }.
\end{split}
\end{equation}The next lemma shows that these quantities are conserved along the flow of   equation \eqref{CMdNLSwithPi}.

\begin{lem}\label{LemConservLaws}
Let $u \in \mathcal{C} \bigl( (T^-, T^+) ; \mathcal{M}^3 \bigr)$ be a solution to \eqref{CMdNLSwithPi} with initial datum $u(0) = u_0 \in \mathcal{M}^3$. Then $\mathfrak{W}^{(1)}(u(t)) = \mathfrak{W}^{(1)}(u_0)$ and $\mathfrak{W}^{(2)}(u(t)) = \mathfrak{W}^{(2)}(u_0)$, for any $t \in (T^-, T^+)$.
\end{lem}

\begin{proof}
Since $u \in \mathcal{C} \bigl( (T^-, T^+) ; \mathcal{M}^3 \bigr)$ solves \eqref{CMdNLSwithPi}, by the Lax pair formulation established in \eqref{fracLuBuLaxPfoc} and \eqref{CommfrakBufrLufor},   the operators $\mathfrak{L}_u$ and $\mathfrak{B}_u$ satisfy the evolution equation $\partial_t \mathfrak{L}_u = [\mathfrak{B}_u, \mathfrak{L}_u]$.  
For any $t \in (T^-, T^+)$, we know from Lemmas~\ref{LemSelAdjscrLpsi} and $\ref{tilscrBpsiLem}$ that $\mathfrak{L}_{u(t)}$, $\mathfrak{L}_{u(t)}^2$ and $i\widetilde{\mathfrak{B}}_{u(t)}$ are unbounded self-adjoint operators on $L^2 (\mathbb{R} ; \mathbb{C})$. The operator $\mathfrak{B}_{u(t)}$ is skew-adjoint and bounded on $L^2 (\mathbb{R} ; \mathbb{C})$.
 For any $\epsilon >0$, let $\chi_{\epsilon} \in \mathrm{Dom}^{\infty} (\mathbf{G})$ be defined in \eqref{chiepFor}. Since $u(t) \in \mathcal{M}^3$, we have $ u(t) \chi_{\epsilon} \in H^3 (\mathbb{R} ; \mathbb{C})$ by \eqref{ufHn<uZnfHn}. 
Set
\begin{equation}\label{sfqDefEnertep}
\mathsf{q}_{\epsilon}^{(0)} (t): = \| u (t) \chi_{\epsilon}\|_{L^2}^2, \quad 
 \mathsf{q}_{\epsilon}^{(1)} (t): = \langle \mathfrak{L}_{u(t)}  ( u (t) \chi_{\epsilon} ) , u (t) \chi_{\epsilon} \rangle_{L^2} , \quad
 \mathsf{q}_{\epsilon}^{(2)} (t): = \|  \mathfrak{L}_{u(t)}  ( u (t) \chi_{\epsilon} )  \|_{L^2}^2.
\end{equation}
For any $k \in \{0,1,2, 4\}$ and $t \in (T^-, T^+)$, let  $\mathbf{j}^{(k)}_{\epsilon}(u (t)) \in \mathbb{C}$ be defined in \eqref{bfjepsukDef} and \eqref{defbfjepsu4}.
 Taking derivative for \eqref{sfqDefEnertep} and using 
\eqref{CMdNLSwithPi}, 
\eqref{tilBuNCDef}, 
\eqref{CommfrakBufrLufor}, 
\eqref{bfrepuexprebytilBu} and
\eqref{bfjepsukDef}, we obtain that
\begin{equation}\label{dtsfqepktfor}
\partial_t \mathsf{q}_{\epsilon}^{(0)} (t) = \mathbf{j}^{(0)}_{\epsilon}(u (t)), \quad 
\partial_t \mathsf{q}_{\epsilon}^{(1)} (t) = \mathbf{j}^{(1)}_{\epsilon}(u (t)), \quad 
\partial_t \mathsf{q}_{\epsilon}^{(2)} (t) = \mathbf{j}^{(2)}_{\epsilon}(u (t))  .
\end{equation}Set $\mathsf{q}_{\epsilon}^{(4)} (t) :=  \|  \mathfrak{L}_{u(t)}^2  ( u (t) \chi_{\epsilon} )  \|_{L^2}^2$, $\forall t   \in (T^-, T^+)$. By \eqref{CMdNLSwithPi},  \eqref{bfjepsukDef},  \eqref{CommfrakBufrLufor}, \eqref{bfrepuexprebytilBu} and \eqref{defbfjepsu4},  we have
\begin{equation}\label{fordtsfqep4}
\begin{split}
\partial_t \mathsf{q}_{\epsilon}^{(4)} (t) 
=  &  2 \operatorname{Re} \langle [\mathfrak{B}_{u(t)} , \mathfrak{L}_{u(t)}^2 ] (u(t) \chi_{\epsilon}) , \mathfrak{L}_{u(t)}^2  (u(t) \chi_{\epsilon}) \rangle_{L^2  } \\
& \qquad
+ 2 \operatorname{Re} \langle   \mathfrak{L}_{u(t)} \bigl(\mathbf{r}_{\epsilon}(u(t) )   + \widetilde{\mathfrak{B}}_{u(t)} (u(t) \chi_{\epsilon}) \bigr)  , \mathfrak{L}_{u(t)}^3  (u(t) \chi_{\epsilon}) \rangle_{L^2  }  \\
= & \mathbf{j}^{(4)}_{\epsilon}(u(t)) .
\end{split}
\end{equation}For any $\epsilon >0 $ and $t \in (T^-, T^+)$,  formula  \eqref{sfqDefEnertep} implies that
\begin{equation}\label{EnerAuxisfQ12}
\begin{cases}
\mathsf{Q}_{\epsilon}^{(1)}(t) := \| \mathfrak{L}_{u(t)} (u(t) \chi_{\epsilon} )  + u(t) \chi_{\epsilon} \|_{L^2 (\mathbb{R} ; \mathbb{C})}^2 = \mathsf{q}_{\epsilon}^{(2)} (t) + 2 \mathsf{q}_{\epsilon}^{(1)} (t) + \mathsf{q}_{\epsilon}^{( 0 )} (t) , \\
\mathsf{Q}_{\epsilon}^{(2)}(t) := \| \mathfrak{L}_{u(t)}^2 (u(t) \chi_{\epsilon} )  - u(t) \chi_{\epsilon} \|_{L^2 (\mathbb{R} ; \mathbb{C})}^2 = \mathsf{q}_{\epsilon}^{(4)} (t) - 2 \mathsf{q}_{\epsilon}^{(2)} (t) + \mathsf{q}_{\epsilon}^{( 0 )} (t) .
\end{cases}
\end{equation}Plugging  \eqref{dtsfqepktfor} and \eqref{fordtsfqep4} into  \eqref{EnerAuxisfQ12} with the fundamental theorem of calculus, we have
\begin{equation}\label{NewLeibsfQktep}
\begin{cases}
\mathsf{Q}_{\epsilon}^{(1)}(t) - \mathsf{Q}_{\epsilon}^{(1)}(0) = \int_0^t  \bigl( \mathbf{j}^{(2)}_{\epsilon}(u(\mathfrak{s}))  +  2 \mathbf{j}^{(1)}_{\epsilon}(u(\mathfrak{s})) + \mathbf{j}^{(0)}_{\epsilon}(u(\mathfrak{s}))\bigr) \mathrm{d}\mathfrak{s}, \\
\mathsf{Q}_{\epsilon}^{(2)}(t) - \mathsf{Q}_{\epsilon}^{(2)}(0) = \int_0^t  \bigl( \mathbf{j}^{(4)}_{\epsilon}(u(\mathfrak{s}))  -  2 \mathbf{j}^{(2)}_{\epsilon}(u(\mathfrak{s})) + \mathbf{j}^{(0)}_{\epsilon}(u(\mathfrak{s}))\bigr) \mathrm{d}\mathfrak{s} .
\end{cases}
\end{equation}For any $ \mathfrak{s} \in (T^-, T^+)$, formulas 
\eqref{bfjepuLim012} and 
\eqref{Lim42bfjepu} imply that
\begin{equation}\label{integrCVto0bfjksum}
\lim_{\epsilon \to 0^+} \big| \mathbf{j}^{(2)}_{\epsilon}(u(\mathfrak{s}))  +  2 \mathbf{j}^{(1)}_{\epsilon}(u(\mathfrak{s})) + \mathbf{j}^{(0)}_{\epsilon}(u(\mathfrak{s})) \big| = \lim_{\epsilon \to 0^+} \big|  \mathbf{j}^{(4)}_{\epsilon}(u(\mathfrak{s}))  -  2 \mathbf{j}^{(2)}_{\epsilon}(u(\mathfrak{s})) + \mathbf{j}^{(0)}_{\epsilon}(u(\mathfrak{s})) \big| =0.
\end{equation}Fix $t \in (T^-, T^+)$. By  
\eqref{bfjepsuUniEsti012}, 
\eqref{UniEsti24bfjepsu} and the continuity of $u : (T^-, T^+) \to \bigl( \mathcal{M}^3 , \mathsf{d}_{\mathcal{M}^3}  \bigr)$, we have 
\begin{equation}\label{bfjksumUniEst}
\small 
 \sup_{|\mathfrak{s} |\leq |t|} \sup_{0 < \epsilon \leq 1} \left( \big| \mathbf{j}^{(2)}_{\epsilon}(u(\mathfrak{s}))  +  2 \mathbf{j}^{(1)}_{\epsilon}(u(\mathfrak{s})) + \mathbf{j}^{(0)}_{\epsilon}(u(\mathfrak{s})) \big| + \big|  \mathbf{j}^{(4)}_{\epsilon}(u(\mathfrak{s}))  -  2 \mathbf{j}^{(2)}_{\epsilon}(u(\mathfrak{s})) + \mathbf{j}^{(0)}_{\epsilon}(u(\mathfrak{s})) \big| \right) < +\infty .
\end{equation}Using  \eqref{NewLeibsfQktep}, 
\eqref{integrCVto0bfjksum}, 
\eqref{bfjksumUniEst} and the dominated convergence theorem, we obtain that
\begin{equation}\label{sfQtminusqQzer0}
\lim_{\epsilon \to 0^+} \bigl( \mathsf{Q}_{\epsilon}^{(1)}(t) - \mathsf{Q}_{\epsilon}^{(1)}(0) \bigr) = 0 \quad \text{and} \quad \lim_{\epsilon \to 0^+} \bigl( \mathsf{Q}_{\epsilon}^{(2)}(t) - \mathsf{Q}_{\epsilon}^{(2)}(0) \bigr) = 0 .
\end{equation}From formulas 
\eqref{consvLawfrW12def}, 
  \eqref{EnerAuxisfQ12},   
\eqref{frLu+1uchieLim},  and 
Lemma $\ref{LemfrLusq-1uchielim}$, it follows that
\begin{equation}\label{limsqsfQkept}
\lim_{\epsilon \to 0^+} \sqrt{\mathsf{Q}_{\epsilon}^{(1)}(t) } = \mathfrak{W}^{(1)}(u(t)) \quad \text{and} \quad
\lim_{\epsilon \to 0^+} \sqrt{\mathsf{Q}_{\epsilon}^{(2)}(t) } = \mathfrak{W}^{(2)}(u(t))  , \quad \forall t \in (T^-, T^+).
\end{equation}Consequently, the proof is completed  by  \eqref{limsqsfQkept} and \eqref{sfQtminusqQzer0}.
\end{proof}

\begin{rem}\label{GWPopenRem}
For any $u \in \mathcal{M}^1$, let $\mathsf{g} = \mathsf{g}_u := \mathcal{G}(u) \in \mathcal{M}^1$, where $\mathcal{G} : \mathcal{M}^1 \to \mathcal{M}^1$ is the gauge transform given by \eqref{GauTrdef}. Using \eqref{vhvsq-1L2norm}, \eqref{Dsqusq-1inprod}, \eqref{sfcuFuncDef} and \eqref{dxsfcufor}, we obtain
\begin{equation}\label{DuuPiusqGautr}
\begin{split}
& \|Du - u \Pi_{\geq 0}(|u|^2 -1)\|_{L^2(\mathbb{R};\mathbb{C})}^2 
= \|\partial_x \mathsf{g} + \tfrac{1}{2} \mathsf{g} \mathbb{H}(|\mathsf{g}|^2 -1)\|_{L^2(\mathbb{R};\mathbb{C})}^2 \\
  \geq  &\|\partial_x \mathsf{g}\|_{L^2(\mathbb{R};\mathbb{C})}^2 - \| |D|^{\frac{1}{2}}(|\mathsf{g}|^2 -1)\|_{L^2(\mathbb{R};\mathbb{C})}^2 + \tfrac{1}{6}\| |\mathsf{g}|^2 -1\|_{L^2(\mathbb{R};\mathbb{C})}^2 .
\end{split}
\end{equation}
In the defocusing case \eqref{defocuCMdNLS}, the term $\| |D|^{\frac{1}{2}}(|\mathsf{g}|^2 -1)\|_{L^2 (\mathbb{R} ; \mathbb{C})}^2$ appears with a positive sign, which allows it to be discarded to obtain a uniform estimate for $\| |\mathsf{g}|^2 -1\|_{L^2 (\mathbb{R} ; \mathbb{C})} = \| |u|^2 -1\|_{L^2 (\mathbb{R} ; \mathbb{C})}$. Together with the first conservation law $\mathfrak{W}^{(1)}$, this yields a-priori bounds for both $\| |u|^2 -1\|_{L^2 (\mathbb{R} ; \mathbb{C})}$ and $\|u\|_{\mathcal{Z}^1}$. The second conservation law $\mathfrak{W}^{(2)}$ then provides an a-priori bound for $\|\partial_x^2 u\|_{L^2 (\mathbb{R} ; \mathbb{C})}$, and the global well-posedness of equation \eqref{defocuCMdNLS} follows; see Subsection~3.2 of \cite{Xchenccm25ExpFor} for the detailed argument. 
In the focusing case \eqref{CMdNLSfocu}, however, the sign of $\| |D|^{\frac{1}{2}}(|\mathsf{g}|^2 -1)\|_{L^2 (\mathbb{R} ; \mathbb{C})}^2$ in \eqref{DuuPiusqGautr} is negative, and it cannot be cancelled by the remaining terms. Consequently, the method used for the defocusing equation does not carry over, and the global well-posedness of \eqref{CMdNLSfocu} on $\bigl( \mathcal{M}^2, \mathsf{d}_{\mathcal{M}^2} \bigr)$ remains an open problem.
\end{rem}

\subsection{Unitary Equivalence}\label{subsUniEquiv4.2}

\noindent In this subsection, assume that $T^-<0<T^+$ and let $u\in\mathcal{C}\bigl((T^-,T^+);\mathcal{Z}^3_+\bigr) \cap \mathcal{C}^1 \bigl((T^-,T^+);\mathcal{Z}^1_+\bigr)$ be a fixed solution to the focusing CMdNLS equation \eqref{CMdNLSfocu} with initial datum $u(0)=u_0\in\mathcal{Z}^3_+$.  We define
\begin{equation}\label{bbB=bfButfor}
\mathbb{B}(t):= \mathbf{B}_{u(t)} =  \mathbf{T}_{u(t)} \mathbf{T}_{\partial_x \overline{u(t)}} - \mathbf{T}_{\partial_x u(t)} \mathbf{T}_{\overline{u(t)}} + i \left( \mathbf{T}_{u(t)} \mathbf{T}_{\overline{u(t)}}\right)^2  , \quad \forall t \in (T^-, T^+).
\end{equation}Then
$\mathbb{B} =  \mathscr{P}_{u, \partial_x \overline{u}} - \mathscr{P}_{\partial_x u,  \overline{u}} + i \left( \mathscr{P}_{u,   \overline{u}}\right)^2 $ by   \eqref{BudefForFoc-} and \eqref{scrPu,vdefFor}. Lemma $\ref{LemTcircuPuvcontinu}$ yields that $\mathbb{B} : t \in (T^-, T^+) \mapsto \mathbb{B}(t) = \mathbf{B}_{u(t)} \in \mathcal{B}_{\mathbb{C}}(H^p_+) $ is a continuous mapping, for any $p \in \{0,1,2\}$.
 Since $\mathbb{B}(t) = \mathbf{B}_{u(t)}$ is a skew-adjoint operator on $L^2_+$, $\forall t \in (T^-, T^+)$, 
 Lemma $\ref{UniPropaLemW}$ implies that there exists a unique mapping
\begin{equation}\label{UniPropagbfWSect4}
\mathbf{W} = \mathbb{W}_{\mathbb{B}} \in \bigcap_{p=0}^2 \mathcal{C}^1 \left( (T^-, T^+); \mathcal{B}_{\mathbb{C}}(H^p_+) \right) , 
\end{equation}satisfying the following properties:
\begin{enumerate}
\item For any $t \in (T^-, T^+)$, $\mathbf{W}(t)$ is a unitary operator on $L^2_+$. Specifically,  we have
\begin{equation}\label{bfWtisUnitaryinL2+}
\mathbf{W}(t)^* = \mathbf{W}(t)^{-1} \in \mathcal{B}_{\mathbb{C}}(L^2_+),
\end{equation}where $\mathbf{W}(t)^*$ denotes the $L^2_+$-adjoint of $\mathbf{W}(t) \in \mathcal{B}_{\mathbb{C}}(L^2_+)$. 
Moreover,  $\mathbf{W}(0)=\mathrm{id}_{L^2_+}$.
\item Set $\mathbf{W}^* : t \in (T^-, T^+) \mapsto \mathbf{W}^*(t) : = \mathbf{W}(t)^* = \mathbf{W}(t)^{-1} \in \mathcal{B}_{\mathbb{C}}(L^2_+)$.  Then we have
\begin{equation}\label{bfWstarisCone}
\mathbf{W}^* \in \bigcap_{p=0}^2 \mathcal{C}^1 \left( (T^-, T^+); \mathcal{B}_{\mathbb{C}}(H^p_+) \right) .
\end{equation} 
\item For any $p \in \{0, 1,2\}$ and $t \in (T^-, T^+)$, we have the following limits:
\begin{equation}\label{dtbfWt=ButbfWtlim}
\begin{cases}
\lim_{r \to 0}  \|r^{-1}(\mathbf{W}(t+ r)-\mathbf{W}(t)) - \mathbf{B}_{u(t)} \mathbf{W}(t)\|_{\mathcal{B}_{\mathbb{C}}(H^p_+)} = 0,\\
\lim_{r \to 0}  \|r^{-1}(\mathbf{W}(t+r)^{-1}-\mathbf{W}(t)^{-1}) + \mathbf{W}(t)^{-1}\mathbf{B}_{u(t)} \|_{\mathcal{B}_{\mathbb{C}}(H^p_+)} = 0.
\end{cases}
\end{equation}
\end{enumerate}

\noindent For any $t \in (T^-, T^+)$, Definition~\ref{AudefTuTbaru} and estimate \eqref{BounToeplitzBHs} imply that $\mathbf{A}_{u(t)} = \mathbf{T}_{u(t)} \mathbf{T}_{\overline{u(t)}} \in \bigcap_{p = 0}^3 \mathcal{B}_{\mathbb{C}} (H^p_+)$. 
By Lemma~\ref{lemLuselfadj}, the chiral Lax operator $\mathbf{L}_{u(t)} = D - \mathbf{A}_{u(t)}$, defined in \eqref{LuSect2}, is self-adjoint on $L^2_+$ with domain $H^1_+$. 
The following lemma establishes the unitary equivalence between $\mathbf{L}_{u(t)}$ and $\mathbf{L}_{u_0}$ via the $L^2_+$-unitary propagator $\mathbf{W}(t)$ constructed in \eqref{UniPropagbfWSect4}.
  
\begin{lem}
For any $p \in \{0, 1, 2\}$ and $t \in (T^-, T^+)$, we have $\mathbf{W}(t) \left( H^p_+ \right) = H^p_+ = \mathrm{Dom}(\mathbf{L}_{u(t)}^p)$ and 
\begin{equation}\label{UniEquLunFor}
\mathbf{L}_{u(t)}^p \mathbf{W}(t) = \mathbf{W}(t) \mathbf{L}_{u_0}^p  \;  : \;  H^p_+ \to L^2_+ .
\end{equation}Here we  choose $\mathbf{L}_{u(t)}^0 =\mathrm{id}_{L^2_+}$. As a consequence, for any bounded Borel function $\mathtt{g} : \mathbb{R} \to \mathbb{C}$, we have
\begin{equation}\label{UniEquivf(Lu(t))For}
\mathtt{g} \left( \mathbf{L}_{u(t)} \right) \mathbf{W}(t) = \mathbf{W}(t) \mathtt{g} \left( \mathbf{L}_{u_0} \right) \in \mathcal{B}_{\mathbb{C}} (L^2_+) , \quad \forall t \in (T^-, T^+). 
\end{equation}
  Moreover, for every $\mathrm{xx} \in \{ \mathrm{ac}, \mathrm{sc}, \mathrm{pp} \}$ and $t \in (T^-, T^+)$, the following identities hold:
\begin{equation}\label{SpeSubSpxxevot}
\sigma \bigl( \mathbf{L}_{ u (t) } \bigr) = \sigma \bigl( \mathbf{L}_{ u_0 } \bigr), 
\qquad
\sigma_{\mathrm{xx}} \bigl( \mathbf{L}_{ u (t) } \bigr) = \sigma_{\mathrm{xx}} \bigl( \mathbf{L}_{ u_0 } \bigr), \qquad 
\mathscr{H}_{\mathrm{xx}}\bigl( \mathbf{L}_{ u (t) } \bigr) = \mathbf{W}(t) \Bigl( \mathscr{H}_{\mathrm{xx}}\bigl( \mathbf{L}_{ u_0 } \bigr) \Bigr).
\end{equation}
  For any $\mu \in \mathbb{C} \setminus \sigma \bigl( \mathbf{L}_{ u_0 } \bigr)$, the operator $\mathbf{L}_{u (t)} - \mu : H^1_+ \to L^2_+$ is bijective and we have
\begin{equation}\label{ResLutuniequiv}
\bigl( \mathbf{L}_{u (t)} - \mu  \bigr)^{-1} \mathbf{W}(t) = \mathbf{W}(t) \bigl( \mathbf{L}_{u_0} - \mu  \bigr)^{-1}  : L^2_+ \to H^1_+, \qquad  \forall  t \in (T^-, T^+) .
\end{equation}
\end{lem}

\begin{proof}
Fix $h , \varphi \in H^1_+$. We set $\mathtt{b}(t) = \mathtt{b}_{h, \varphi} (t) := \langle \mathbf{L}_{u(t)} \mathbf{W}(t) h, \mathbf{W}(t) \varphi \rangle_{L^2_+}$, for any $t \in (T^-, T^+)$. Since $u \in \mathcal{C}\left( (T^-, T^+); \mathcal{Z}^3_+ \right)$ solves the CMdNLS equation \eqref{CMdNLSfocu}, the Lax pair structure \eqref{LaxPairBuLuchi} implies that
\begin{equation}\label{LaxPairdAutdtu=[Lu,Bu]}
\mathrm{d}\mathbf{A}(u(t)) \left( \partial_t u (t) \right) = \left[  \mathbf{L}_{u(t)} , \mathbf{B}_{u(t)}  \right] \;  : \; H^1_+ \to L^2_+, \quad \forall t \in (T^-, T^+).
\end{equation}Recall that $\mathbf{B}_{u(t)}^* = - \mathbf{B}_{u(t)}$. Using   \eqref{dtbfWt=ButbfWtlim} and \eqref{LaxPairdAutdtu=[Lu,Bu]}, we obtain that $\mathtt{b} : \mathbb{R}  \to \mathbb{C}$ is differentiable and
\begin{equation*}
\partial_t \mathtt{b}(t) = \langle \left[\mathbf{L}_{u(t)}, \mathbf{B}_{u(t)} \right] \mathbf{W}(t) h , \mathbf{W}(t) \varphi \rangle_{L^2_+} 
- 
\langle \left( \mathrm{d}\mathbf{A}(u(t)) \left( \partial_t u(t) \right) \right) \mathbf{W}(t) h, \mathbf{W}(t) \varphi \rangle_{L^2_+} = 0,
\end{equation*}for any $t \in (T^-, T^+)$. So $\mathbf{L}_{u(t)} =
\mathbf{W}(t) \mathbf{L}_{u_0} \mathbf{W}(t)^* : H^1_+ \to L^2_+$. Since $\mathbf{W}(t)$ is a unitary operator on $L^2_+$ and  $\mathbf{W}(t) (H^p_+) = H^p_+$ is bijective, for any $p \in \{1,2\}$, formula \eqref{UniEquLunFor} holds. Consequently, the two self-adjoint operators $\bigl( \mathbf{L}_{u(t)} , \mathrm{Dom} (\mathbf{L}_{u(t)})   \bigr)$ and $\bigl( \mathbf{L}_{u_0} , \mathrm{Dom} (\mathbf{L}_{u_0})   \bigr)$ are unitarily equivalent. 
 Formulas \eqref{UniEquivf(Lu(t))For},    \eqref{SpeSubSpxxevot} and 
\eqref{ResLutuniequiv} 
 follow directly from   \eqref{UniEquLunFor} and Proposition 2.1 of \cite{SUNCMP2021}.
\end{proof}

\noindent To simplify the  action of the unitary propagator $\mathbf{W}(t) \exp \left(-i t \mathbf{L}_{u_0}^2 \right)$, we introduce the following notation.

\begin{defi}\label{bbUt=exp(-itLu02)Def}
For any $\mathbf{t}\in \mathbb{R}$, we define $\mathbf{U}(\mathbf{t}):= \exp \left(-i \mathbf{t}  \mathbf{L}_{u_0}^2 \right) \in   \mathcal{B}_{\mathbb{C}} (L^2_+)$.
\end{defi}
\noindent Lemma $\ref{expitLusqStrContGroup}$ implies that $\left(\mathbf{U}(\mathbf{t}) \right)_{\mathbf{t} \in \mathbb{R}}$ is a strongly  continuous unitary group on $L^2_+$. Then we have
\begin{equation}\label{bbUt=e-itLu02QuoDiffCV}
\lim_{\mathbf{t} \to 0} \| \tfrac{1}{\mathbf{t}}\left(  \mathbf{U}(\mathbf{t}) \phi - \phi\right) + i \mathbf{L}_{u_0}^2 (\phi)\|_{L^2_+} = 0, \qquad \forall \phi \in H^2_+ = \mathrm{Dom}(\mathbf{L}_{u_0}^2),
\end{equation}for any $\mathbf{t} \in \mathbb{R}$. The next lemma identifies the composition $\mathbf{W}(t)e^{-it\mathbf{L}_{u_0}^2}\varphi$ as the unique solution to a linear evolution equation driven by $\widetilde{\mathbf{B}}_{u(t)}$, which is defined in  \eqref{tildeBudefFoc-}.

\begin{lem}\label{LemWtexp-itL2phiPDE}
For any $\varphi \in H^2_+$, there exists a unique function
\begin{equation}\label{fuvarphidefFor}
\mathbf{f} \in \mathcal{C} \left( (T^-, T^+) ; H^2_+ \right) \cap     \mathcal{C}^1 \left( (T^-, T^+) ; L^2_+ \right),
\end{equation}such that $\mathbf{f}(0) = \varphi \in H^2_+$ and the following linear equation holds:
\begin{equation}\label{PDEdtfistilButf}
\partial_t \mathbf{f}   =     i\partial_x^2 \mathbf{f} + 2  u \partial_x \Pi_{\geq 0} \left( \overline{u} \mathbf{f} \right),  \qquad t \in (T^-,T^+), \qquad x \in \mathbb{R} . 
\end{equation}Moreover, this solution satisfies that
 \begin{equation}\label{ExpForftsWteitLusq}
\mathbf{f}(t) = \mathbf{f}_{u, \varphi}(t) = \mathbf{W}(t) e^{-it\mathbf{L}_{u_0}^2} ( \varphi )  \in H^2_+, \qquad \forall t \in (T^-, T^+) .
\end{equation}
\end{lem}
\begin{rem}
By \eqref{tildeBudefFoc-}, equation \eqref{PDEdtfistilButf} 
is to be understood in the following sense
\begin{equation}\label{limformPDEdtf-tiButf}
\lim_{r \to 0} \|\tfrac{1}{r} \left(\mathbf{f}(t+r) -\mathbf{f}(t)\right) -  \widetilde{\mathbf{B}}_{u(t)} \left( \mathbf{f}(t) \right) \|_{L^2_+} = 0, \quad \forall t \in (T^-, T^+).
\end{equation}
\end{rem}
  
\begin{proof}[Proof of Lemma $\ref{LemWtexp-itL2phiPDE}$]Fix $\varphi \in H^2_+$. The proof is divided into two steps. \\

\noindent \textbf{Step \romannumeral 1}. Let $\mathbf{f} = \mathbf{f}_{u, \varphi}$ be defined in \eqref{ExpForftsWteitLusq}. Since $u_0 \in \mathcal{Z}^3_+$, Lemma $\ref{expitLusqStrContGroup}$ shows that $\left( \mathbf{U}(\mathbf{t}) \right)_{\mathbf{t} \in  \mathbb{R}}$ is a strongly continuous group of bounded operators on $H^2_+$. So we have $\mathbf{f} \in \mathcal{C} \left(  (T^-, T^+); H^2_+\right)$ by  formula \eqref{UniPropagbfWSect4}.
 Using \eqref{tildeBudefFoc-}, \eqref{dtbfWt=ButbfWtlim}, \eqref{UniEquLunFor}
and \eqref{bbUt=e-itLu02QuoDiffCV}, we obtain that
\begin{equation}\label{bffquodiffcvtilB}
\small
\begin{split}
 &\|\tfrac{1}{r}\left(\mathbf{f}(t + r) - \mathbf{f}(t) \right) -\widetilde{\mathbf{B}}_{u(t)} \left( \mathbf{f}(t) \right) \|_{L^2_+} \\
\leq & \|\left( \mathbf{W}(t+r) - \mathbf{W}(t) \right)\tfrac{1}{r} \left(\mathbf{U}(t+r) \varphi - \mathbf{U}(t) \varphi \right) \|_{L^2_+} 
+ \|\tfrac{1}{r} \left( \left(\mathbf{W}(t+r) - \mathbf{W}(t)\right) - \mathbf{B}_{u(t)} \mathbf{W}(t) \right) \mathbf{U}(t) \varphi \|_{L^2_+}\\
& \quad \quad  + \|\mathbf{W}(t) \left( \tfrac{1}{r} \left( \mathbf{U}(t+r) \varphi - \mathbf{U}(t) \varphi \right) + i \mathbf{U}(t) \mathbf{L}_{u_0}^2 (\varphi) \right) \|_{L^2_+}\\
\to &  0, \quad \mathrm{as} \quad r \to 0.
\end{split}
\end{equation}Then $\mathbf{f} : t \in (T^-, T^+) \mapsto \mathbf{f}(t) = \mathbf{W}(t) \mathbf{U}(t) \varphi \in H^2_+$ solves  \eqref{PDEdtfistilButf}.
Lemma $\ref{quvttoquvtLemContinu}$ implies that the mapping $\partial_t \mathbf{f}: t \in (T^-, T^+) \mapsto \partial_t \mathbf{f}(t) = \left( i \partial_x^2   + 2  \mathfrak{q}_{u, u}(t) \right) (\mathbf{f}(t)) \in L^2_+$ is continuous. So $\mathbf{f} \in \mathcal{C}^1 \left( (T^-, T^+); L^2_+ \right)$. \\

\noindent \textbf{Step \romannumeral 2}. Let $f \in \mathcal{C} \left( (T^-, T^+); H^2_+ \right) \cap  \mathcal{C}^1 \left( (T^-, T^+); L^2_+ \right)$ be a solution to equation  \eqref{PDEdtfistilButf} with $f(0) = \varphi$. 
Then the function $f$ satisfies  \eqref{limformPDEdtf-tiButf}. Fix $t \in (T^-, T^+)$. 
By  \eqref{tildeBudefFoc-},   \eqref{dtbfWt=ButbfWtlim},
and \eqref{UniEquLunFor}, 
it follows that
\begin{equation}\label{dtWt-1ft+iLu2Wt-1ftFor}
\small
\begin{split}
 & \|\tfrac{1}{r}\left( \mathbf{W}(t+r)^{-1} f(t+r) 
-\mathbf{W}(t)^{-1} f(t) \right)  + i \mathbf{L}_{u_0}^2 \mathbf{W}(t)^{-1} \left( f(t) \right)\|_{L^2_+} \\
\leq & \|\left(\mathbf{W}(t+r)^{-1} - \mathbf{W}(t)^{-1} \right) \tfrac{1}{r}\left(f(t+r) -f(t) \right)  \|_{L^2_+}
+ \|\mathbf{W}(t)^{-1} \left( \tfrac{1}{r} \left( f(t+r) - f (t) \right)  - \widetilde{\mathbf{B}}_{u(t)} \left( f(t) \right) \right)\|_{L^2_+}\\
&  \qquad \qquad + \| \left(\tfrac{1}{r} \left(\mathbf{W}(t+r)^{-1} - \mathbf{W}(t)^{-1} \right)  + \mathbf{W}(t)^{-1}  \mathbf{B}_{u(t)} \right) f(t) \|_{L^2_+}\\
\to &  0,    
\end{split}
\end{equation}as $r \to 0$.
 Recall that $\mathbf{W}(t)^{-1} f(t) \in H^2_+ $ by   \eqref{bfWstarisCone}.
 Using \eqref{dtbfWt=ButbfWtlim}, \eqref{bbUt=e-itLu02QuoDiffCV} 
and
\eqref{dtWt-1ft+iLu2Wt-1ftFor},
we obtain that 
\begin{equation}
\begin{split}
& \| \tfrac{1}{r} \left( \mathbf{U}(-t-r) \mathbf{W}(t+r)^{-1} f(t+r) - \mathbf{U}(-t) \mathbf{W}(t)^{-1} f(t)  \right)  \|_{L^2_+}\\
\leq & \| \left( \mathbf{U}(-t-r) - \mathbf{U}(-t) \right) \tfrac{1}{r} \left(  \mathbf{W}(t+r)^{-1} f(t+r) -   \mathbf{W}(t)^{-1} f(t)  \right) \|_{L^2_+} \\
& \qquad + \|\mathbf{U}(-t)\left(\tfrac{1}{r}\left( \mathbf{W}(t+r)^{-1} f(t+r) 
-\mathbf{W}(t)^{-1} f(t) \right)  + i \mathbf{L}_{u_0}^2 \mathbf{W}(t)^{-1} \left( f(t) \right) \right)\|_{L^2_+} \\
& \qquad + \|\mathbf{U}(-t) \left( \tfrac{1}{r} \left( \mathbf{U}(-r) - \mathrm{id}_{L^2_+} \right) \mathbf{W}(t)^{-1} f(t)   - i \mathbf{L}_{u_0}^2 \mathbf{W}(t)^{-1} f(t) \right)\|_{L^2_+} \\
\to & 0 ,
\end{split}
\end{equation}as $r \to 0$. 
So the mapping $g :t \in (T^-, T^+) \mapsto g(t) : = \mathbf{U}(-t) \mathbf{W}(t)^{-1}  \left(f(t) \right) \in H^2_+$ is constant. Since $\mathbf{U}(0) = \mathbf{W}(0) = \mathrm{id}_{L^2_+}$, we have $f(t) = \mathbf{W}(t) \mathbf{U}(t) \left( f(0) \right) = \mathbf{W}(t) \exp \left(- i t \mathbf{L}_{u_0}^2  \right) \varphi = \mathbf{f}_{u, \varphi} (t) \in H^2_+$.  
\end{proof}

\noindent For any $t \in (T^-, T^+)$ and $\eta >0$, recall that $\mathbf{S}(\eta)^* = \mathbf{T}_{\mathbf{e}_{-\eta}} \in \mathcal{B}_{\mathbb{C}}(L^2_+)$ is defined in  \eqref{LaxBeurlingSG}.

\begin{lem}\label{LemConjSstcvto2tLuphi}
For any $t\in(T^-,T^+)$ and $\varphi\in H^2_+$, the following formula holds:
\begin{equation}\label{LimforconjSeta-SetaCV}
\lim_{\eta\to 0^+} \Biggl\| \tfrac{i}{\eta}\Bigl( e^{it\mathbf{L}_{u_0}^2} \mathbf{W}(t)^{-1} \mathbf{S}(\eta)^* \mathbf{W}(t) e^{-it\mathbf{L}_{u_0}^2} - \mathbf{S}(\eta)^* \Bigr)\varphi - 2t\mathbf{L}_{u_0}(\varphi) \Biggr\|_{L^2_+} = 0.
\end{equation}
\end{lem}
\begin{proof}
Fix two arbitrary functions $\varphi_1, \varphi_2 \in H^2_+$. For any $t \in (T^-, T^+)$ and $j \in \{1, 2 \}$, we define
\begin{equation}\label{sffjDefFor}
\mathsf{f}_j (t) :=\mathbf{f}_{u, \varphi_j} (t) = \mathbf{W}(t) \exp\left(- i t \mathbf{L}_{u_0}^2 \right) \varphi_j = \mathbf{W}(t)  \mathbf{U}(t) \varphi_j  \in H^2_+.
\end{equation}By Lemma $\ref{LemWtexp-itL2phiPDE}$, the function $\mathsf{f}_j \in  \mathcal{C} \left( (T^-, T^+) ; H^2_+ \right) \cap     \mathcal{C}^1 \left( (T^-, T^+) ; L^2_+ \right)$ satisfies   \eqref{limformPDEdtf-tiButf}. So we have   
\begin{equation}\label{sffj(0)=varphij}
\partial_t \mathsf{f}_j(t) = \widetilde{\mathbf{B}}_{u(t)} \left( \mathsf{f}_j(t) \right),\quad
\mathsf{f}_j(0) = \varphi_j \in H^2_+,   \quad \forall t \in (T^-, T^+), \quad \forall j \in \{1,2\}.
\end{equation}For any $\eta >0$, the function $\boldsymbol{\nu}_{\eta} =
\boldsymbol{\nu}^{u, \varphi_1, \varphi_2}_{\eta} \in \mathcal{C}^1  \left( (T^-, T^+) ; \mathbb{C} \right)$ is defined as follows: 
\begin{equation}\label{bfnuDefFor}
 \boldsymbol{\nu}_{\eta}(t)  : = \langle \mathbf{S}(\eta)^* \left(\mathsf{f}_1 (t) \right),   \mathsf{f}_2 (t) \rangle_{L^2_+} = \langle e^{it\mathbf{L}_{u_0}^2} \mathbf{W}(t)^{-1} \mathbf{S}(\eta)^* \mathbf{W}(t) e^{-it\mathbf{L}_{u_0}^2}   (\varphi_1) , \varphi_2 \rangle_{L^2_+} ,
\end{equation}for any $t \in (T^-, T^+)$. Since $u(t) \in \mathcal{Z}^3_+$, the operator $i \widetilde{\mathbf{B}}_{u(t)}  = \mathbf{L}_{u(t)}^2 + i \mathbf{B}_{u(t)}$ is a self-adjoint operator on $L^2_+$ with domain   $H^2_+$ by Lemma \ref{LusqiLemselfadj>0} and the Kato--Rellich theorem. Then \eqref{sffj(0)=varphij} yields that
\begin{equation}\label{dtbfnuFor}
\partial_t  \boldsymbol{\nu}_{\eta}(t) = \langle \left[ \mathbf{S}(\eta)^*, \widetilde{\mathbf{B}}_{u(t)}  \right] \left( \mathsf{f}_1 (t) \right),  \mathsf{f}_2 (t) \rangle_{L^2_+} = \langle \mathbf{U}(-t) \mathbf{W}(t)^{-1} 
\left[ \mathbf{S}(\eta)^*, \widetilde{\mathbf{B}}_{u(t)}  \right]
 \mathbf{W}(t) \mathbf{U}(t) \varphi_1 , \varphi_2 \rangle_{L^2_+}, 
\end{equation}for any $ t \in (T^-, T^+)$. Recall that $\mathbf{U}(-t) = \exp \left( it \mathbf{L}_{u_0}^2 \right) \in \bigcap_{p=0}^3 \mathcal{B}_{\mathbb{C}}(H^p_+)$ by Lemma $\ref{expitLusqStrContGroup}$. 
Then  we have 
\begin{equation}\label{CForU-tWt-1LutWtUtLu0}
 \langle \mathbf{U}(-t) \mathbf{W}(t)^{-1} 
 \mathbf{L}_{u(t)}
 \mathbf{W}(t) \mathbf{U}(t) \varphi_1 , \varphi_2 \rangle_{L^2_+} = \langle  \mathbf{L}_{u_0} ( \varphi_1) , \varphi_2 \rangle_{L^2_+} , \quad \forall t \in (T^-, T^+),
\end{equation}by using   \eqref{UniEquLunFor} and
\eqref{expisLusqcommutLu-a}.  For any $\eta >0$, let $\mathbf{K}^{(\eta)} \in \mathcal{C} \bigl( \mathcal{Z}^1_+ ; \mathcal{B}_{\mathbb{C}} (H^2_+ ; L^2_+) \bigr)$  be given by formula \eqref{bbKetaudefFor} and Definition $\ref{DeffbfKeta}$. For any $t \in (T^-, T^+)$,  
 we define that 
\begin{equation}\label{getauvarphidefFor}
  g_{\eta}^{u, \varphi_1} (t) : = \Bigl( \mathbf{K}^{(\eta)} \circ u(t) \Bigr) \left(\mathsf{f}_1 (t) \right)   = 
\left( \tfrac{i}{\eta}\left[ \mathbf{S}(\eta)^*, \widetilde{\mathbf{B}}_{u(t)}  \right] - 2 \mathbf{L}_{u(t)} \right)
 \mathbf{W}(t) \mathbf{U}(t) \varphi_1 \in L^2_+.
\end{equation}Since $u \in \mathcal{C}\left((T^-, T^+) ;\mathcal{Z}^3_+ \right)$, Corollary $\ref{CorbbKcircuuniEquiCont}$ yields that $g_{\eta}^{u, \varphi_1}  \in \mathcal{C} \left( (T^-, T^+); L^2_+ \right)$. Using the continuity of $\mathbf{W}^* : t \in (T^-, T^+) \mapsto \mathbf{W}(t)^* =\mathbf{W}(t)^{-1} \in \mathcal{B}_{\mathbb{C}}(L^2_+)$, we obtain that
\begin{equation}\label{hetauphi1DEFfor}
\mathsf{h}_{\eta}^{u, \varphi_1} : t \in (T^-, T^+) \mapsto   \mathsf{h}_{\eta}^{u, \varphi_1} ( t) := \mathbf{W}(t)^{-1} \left( g_{\eta}^{u, \varphi_1} (t) \right) \in L^2_+
\end{equation}is continuous. For any $\eta>0$ and $t \in (T^-, T^+)$, we define $\mathtt{y}_{\eta}^{u, \varphi_1}(t):= \mathbf{U}(-t) \left(\mathsf{h}_{\eta}^{u, \varphi_1} ( t) \right) $. Then 
\begin{equation}\label{ttyetauphi1tDefFor}
 \mathtt{y}_{\eta}^{u, \varphi_1}(t) = 
\mathbf{U}(-t) \mathbf{W}(t)^{-1} 
\left( \tfrac{i}{\eta}\left[ \mathbf{S}(\eta)^*, \widetilde{\mathbf{B}}_{u(t)}  \right] - 2 \mathbf{L}_{u(t)} \right)
 \mathbf{W}(t) \mathbf{U}(t) \varphi_1  \in L^2_+, 
\end{equation}by \eqref{getauvarphidefFor} and \eqref{hetauphi1DEFfor}. Since $\left( \mathbf{U}(\mathbf{t}) \right)_{\mathbf{t} \in \mathbb{R}}$ is a  strongly continuous unitary group  on $L^2_+$ by Lemma $\ref{expitLusqStrContGroup}$, the mapping $\mathtt{y}_{\eta}^{u, \varphi_1} : t \in (T^-, T^+) \mapsto \mathtt{y}_{\eta}^{u, \varphi_1} ( t) \in  L^2_+$ is continuous.
So it is Bochner integrable over any compact subset of $(T^-, T^+)$ by  Theorem $3.27$ of  \cite{RudinFA}. 
 Fix $t \in (T^-, T^+)\setminus \{0\}$ and set
\begin{equation}\label{calItDef}
\mathtt{I}_t:=\begin{cases}
[0, t], \qquad \mathrm{if} \quad t >0; \\
[t, 0], \qquad \mathrm{if} \quad t <0.
\end{cases}
\end{equation}By Definition $3.26$ of \cite{RudinFA}, the $L^2_+$-valued integral  $\mathscr{Z}_{\eta}^{u, \varphi_1} (t) : = \int_{\mathtt{I}_t} \mathtt{y}_{\eta}^{u, \varphi_1}  (\mathfrak{s}) \mathrm{d} \mathfrak{s} \in L^2_+$   is well defined, for any $\eta >0$. Moreover, since $L^2_+$
  is a Hilbert space, we have the following characterization:  
\begin{equation}\label{BochIntZChar}
\langle \mathscr{Z}_{\eta}^{u, \varphi_1} (t) , \phi \rangle_{L^2_+} = \int_{\mathtt{I}_t}
\langle \mathtt{y}_{\eta}^{u, \varphi_1}  (\mathfrak{s}), \phi
\rangle_{L^2_+} \mathrm{d} \mathfrak{s},   \qquad \forall \phi \in L^2_+ .
\end{equation}By   \eqref{EstiUthH2norm} and \eqref{UniPropagbfWSect4}, we have $\mathsf{A}_1(u,t):= \sup_{\mathfrak{s} \in \mathtt{I}_t} \| \mathbf{W}(\mathfrak{s}) \|_{\mathcal{B}_{\mathbb{C}}(H^2_+)} \sup_{\mathbf{t} \in \mathbb{R}} \| \mathbf{U}(\mathbf{t}) \|_{\mathcal{B}_{\mathbb{C}}(H^2_+)} <+\infty$.  Set
\begin{equation}\label{fraknuphisetaDefFor}
\mathfrak{n}_{u, \varphi_1}(\mathfrak{s}, \eta):= \|g_{\eta}^{u, \varphi_1} (\mathfrak{s})\|_{L^2_+} 
=\Big\|\left( \tfrac{i}{\eta}\left[ \mathbf{S}(\eta)^*, \widetilde{\mathbf{B}}_{u(\mathfrak{s} )}  \right] - 2 \mathbf{L}_{u(\mathfrak{s} )} \right)
\mathsf{f}_1(\mathfrak{s}) \Big\|_{L^2_+}, \quad \forall \mathfrak{s} \in \mathtt{I}_t, \quad \forall \eta >0.
\end{equation}Formula \eqref{commForSetasttildBuCv} yields that $\lim_{\eta \to 0^+}  \mathfrak{n}_{u, \varphi_1}(\mathfrak{s}, \eta) = 0$.   Plugging \eqref{sffjDefFor} into \eqref{fraknuphisetaDefFor},   we obtain that
\begin{equation}
\sup_{\mathfrak{s} \in \mathtt{I}_t} \sup_{\eta >0 } \mathfrak{n}_{u, \varphi_1}(\mathfrak{s}, \eta) \leq  \left(4 \mathsf{A}_1(u,t) + 6 \sup_{\mathfrak{s} \in \mathtt{I}_t} \|u(\mathfrak{s})\|_{L^{\infty}(\mathbb{R}; \mathbb{C})}^2  \right) \|\varphi_1 \|_{H^2_+} <+\infty ,
\end{equation}by   \eqref{SetasttildBuUniEstL2}.   
Both $\mathbf{U}(-\mathfrak{s})$ and $\mathbf{W}(\mathfrak{s})^{-1}$ are unitary operators on $L^2_+$, so $\|\mathtt{y}_{\eta}^{u, \varphi_1}  (\mathfrak{s})\|_{L^2_+} = \|g_{\eta}^{u, \varphi_1}  (\mathfrak{s})\|_{L^2_+}$, $\forall \mathfrak{s} \in \mathtt{I}_t$.
By Theorem $3.29$ of   \cite{RudinFA} and the dominated convergence theorem, we have
 \begin{equation}\label{CVDomZuphietatto0}
  \|\mathscr{Z}_{\eta}^{u, \varphi_1} (t)  \|_{L^2_+}
   \leq   \int_{\mathtt{I}_t}  \|\mathtt{y}_{\eta}^{u, \varphi_1}  (\mathfrak{s})\|_{L^2_+} \mathrm{d}\mathfrak{s} 
  =
  \int_{\mathtt{I}_t} \mathfrak{n}_{u, \varphi_1}(\mathfrak{s}, \eta) \mathrm{d}\mathfrak{s} 
   \to 
  0, \qquad  \mathrm{as} \quad \eta \to 0^+.
 \end{equation}Formula  \eqref{bfnuDefFor} implies that
$\boldsymbol{\nu}_{\eta}(0) = \langle
\mathbf{S}(\eta)^*\varphi_1 , \varphi_2
 \rangle_{L^2_+}$.  
Plugging \eqref{ttyetauphi1tDefFor} into \eqref{dtbfnuFor}, we obtain that
\begin{equation}\label{NeLeibbfnuFor}
\begin{split}
& \Bigl\langle  \tfrac{i}{\eta}\Bigl( e^{it\mathbf{L}_{u_0}^2} \mathbf{W}(t)^{-1} \mathbf{S}(\eta)^* \mathbf{W}(t) e^{-it\mathbf{L}_{u_0}^2} - \mathbf{S}(\eta)^* \Bigr)\varphi_1  , \varphi_2  \Bigr\rangle_{L^2_+} = \tfrac{i}{\eta} \left( \boldsymbol{\nu}_{\eta}(t) - \boldsymbol{\nu}_{\eta}(0) \right) \\
= &\tfrac{i}{\eta} \int_0^t \partial_t \boldsymbol{\nu}_{\eta} (\mathfrak{s}) \mathrm{d} \mathfrak{s}=  \mathrm{sign}(t)\int_{\mathtt{I}_t}
\Bigl\langle \mathbf{U}(-\mathfrak{s}) \mathbf{W}(\mathfrak{s})^{-1} 
\tfrac{i}{\eta}\left[ \mathbf{S}(\eta)^*, \widetilde{\mathbf{B}}_{u(\mathfrak{s})}  \right]
 \mathbf{W}(\mathfrak{s}) \mathbf{U}(\mathfrak{s}) \varphi_1 , \varphi_2 \Bigr\rangle_{L^2_+}
 \mathrm{d} \mathfrak{s}\\
 = & \mathrm{sign}(t) \int_{\mathtt{I}_t}
\langle \mathtt{y}_{\eta}^{u, \varphi_1}  (\mathfrak{s}), \varphi_2
\rangle_{L^2_+} \mathrm{d} \mathfrak{s} + 2t \langle \mathbf{L}_{u_0} (\varphi_1), \varphi_2 \rangle_{L^2_+} ,
\end{split}
\end{equation}by \eqref{CForU-tWt-1LutWtUtLu0} and \eqref{calItDef}. Since formula \eqref{NeLeibbfnuFor} holds for every $\varphi_2 \in H^2_+$, it follows that
\begin{equation}\label{scrZidenConjSet-Set}
  \tfrac{i}{\eta}\Bigl( e^{it\mathbf{L}_{u_0}^2} \mathbf{W}(t)^{-1} \mathbf{S}(\eta)^* \mathbf{W}(t) e^{-it\mathbf{L}_{u_0}^2} - \mathbf{S}(\eta)^* \Bigr)\varphi_1  - 2 t \mathbf{L}_{u_0} (\varphi_1)
=
\mathrm{sign}(t)\mathscr{Z}_{\eta}^{u, \varphi_1} (t)  \in L^2_+, 
\end{equation}by \eqref{BochIntZChar}, for any $\eta>0$. 
Then formula
\eqref{LimforconjSeta-SetaCV} is a direct consequence of \eqref{CVDomZuphietatto0} and  \eqref{scrZidenConjSet-Set}.
\end{proof}

\noindent As shown in Proposition  $\ref{GisInfiniGeneHsProp}$, the operator  $\left( -i \mathbf{G} , \mathrm{Dom}(\mathbf{G}) \right)$, given by \eqref{DefGphiSect2}, is the infinitesimal generator of   $\left( \mathbf{S}(\eta)^* \right)_{\eta \geq 0}$ on $L^2_+$.  In the next theorem, we show that $\mathrm{Dom}(\mathbf{G}) \cap  H^2_+$ is invariant under the unitary operator $\mathbf{W}(t)e^{-it \mathbf{L}_{u_0}^2 }\in \mathcal{B}_{\mathbb{C}}(L^2_+)$, for any $t \in (T^-, T^+)$.

\begin{thm}\label{invSpDomGHtwoThm}
For any $t \in (T^-, T^+)$,   the following identity holds:
\begin{equation}\label{DomGcapH2InvWU}
\mathbf{W}(t)e^{-it \mathbf{L}_{u_0}^2}\left( \mathrm{Dom}(\mathbf{G}) \cap  H^2_+ \right) = \mathrm{Dom}(\mathbf{G}) \cap  H^2_+    = e^{ i t \mathbf{L}_{u_0}^2} \mathbf{W}(t)^{-1} \left( \mathrm{Dom}(\mathbf{G}) \cap  H^2_+ \right).
\end{equation}Moreover, for any   $\phi \in \mathrm{Dom}(\mathbf{G}) \cap  H^2_+$, we have
\begin{equation}\label{GWtUtisWtUtGadd2tL}
\mathbf{G} \mathbf{W}(t)e^{-it \mathbf{L}_{u_0}^2} (\phi) = \mathbf{W}(t)e^{-it \mathbf{L}_{u_0}^2} \left( \mathbf{G} + 2 t \mathbf{L}_{u_0}\right) \phi  \in L^2_+.
\end{equation}
\end{thm}

\begin{proof}
Fix $t \in (T^-, T^+)$.  For any $\phi \in \mathrm{Dom}(\mathbf{G}) \cap  H^2_+$, recall that $\mathbf{f}_{u, \phi}(t) = \mathbf{W}(t) e^{-it\mathbf{L}_{u_0}^2} ( \phi )  \in H^2_+$ by formula  \eqref{ExpForftsWteitLusq}.  Using Proposition  $\ref{GisInfiniGeneHsProp}$ and formula  \eqref{LimforconjSeta-SetaCV}, we obtain that
\begin{equation}\label{UtWt-1Setaf-fCVL2}
\small 
\begin{split}
& \Big\| e^{it\mathbf{L}_{u_0}^2}   \mathbf{W}(t) ^{-1} \tfrac{i}{\eta} \left( \mathbf{S}(\eta)^* \mathbf{f}_{u, \phi}(t)-  \mathbf{f}_{u, \phi}(t)\right) -\left( \mathbf{G} + 2 t \mathbf{L}_{u_0}\right) \phi \Big\|_{L^2_+} \\
\leq & \Big\| \tfrac{i}{\eta}\Bigl( e^{it\mathbf{L}_{u_0}^2} \mathbf{W}(t)^{-1} \mathbf{S}(\eta)^* \mathbf{W}(t) e^{-it\mathbf{L}_{u_0}^2} - \mathbf{S}(\eta)^* \Bigr)\phi -2 t \mathbf{L}_{u_0} (\phi) \Big\|_{L^2_+} 
+ \Big\|\tfrac{i}{\eta} \left(\mathbf{S}(\eta)^* \phi - \phi \right)  
- \mathbf{G} \phi \Big\|_{L^2_+} \to 0,
\end{split}
\end{equation}as $\eta \to 0^+$. Recall that $e^{it\mathbf{L}_{u_0}^2}   \mathbf{W}(t) ^{-1}$ is a unitary operator on $L^2_+$, by  Lemma $\ref{expitLusqStrContGroup}$ and  formula \eqref{bfWtisUnitaryinL2+}.    
   Then $ \mathbf{W}(t) e^{-it\mathbf{L}_{u_0}^2} ( \phi ) = \mathbf{f}_{u, \phi}(t)  \in  \mathrm{Dom}^0 (\mathbf{G}) = \mathrm{Dom}(\mathbf{G}) $ and formula \eqref{GWtUtisWtUtGadd2tL} follows from \eqref{UtWt-1Setaf-fCVL2}. So   
 \begin{equation}\label{WtUtDGcapH2subsetDG}
 \mathbf{W}(t)e^{-it \mathbf{L}_{u_0}^2}\left( \mathrm{Dom}(\mathbf{G}) \cap  H^2_+ \right) \subset  \mathrm{Dom}(\mathbf{G}) \cap  H^2_+ .
\end{equation}By Lemma $\ref{expitLusqStrContGroup}$ and 
formula \eqref{UniPropagbfWSect4}, 
  the operator $\mathbf{W}(t) e^{-it\mathbf{L}_{u_0}^2} : H^2_+ \to H^2_+$ is bijective and continuous. 
 For any $f \in \mathrm{Dom}(\mathbf{G}) \cap   H^2_+$, set $\varphi=\varphi_{u, f, t}:=   e^{it \mathbf{L}_{u_0}^2} \mathbf{W}(t)^{-1} (f) \in H^2_+ $.   Then  $f = \mathbf{W}(t) e^{ -it \mathbf{L}_{u_0}^2} ( \varphi ) = \mathbf{f}_{u, \varphi}(t) $ by \eqref{ExpForftsWteitLusq}. 
Since $e^{it\mathbf{L}_{u_0}^2}   \mathbf{W}(t) ^{-1}  $ is a unitary  operator on $L^2_+$, it follows that
\begin{equation}\label{CVSetavPhiInvSp}
\small 
\begin{split}
& \Big\| \tfrac{i}{\eta} \left( \mathbf{S}(\eta)^* \varphi -  \varphi \right) 
 + 2 t \mathbf{L}_{u_0} ( \varphi)  
  - e^{it\mathbf{L}_{u_0}^2}   \mathbf{W}(t) ^{-1}  \mathbf{G}f \Big\|_{L^2_+} \\
\leq & \Big\|2 t \mathbf{L}_{u_0} (\varphi) - \tfrac{i}{\eta}\Bigl( e^{it\mathbf{L}_{u_0}^2} \mathbf{W}(t)^{-1} \mathbf{S}(\eta)^* \mathbf{W}(t) e^{-it\mathbf{L}_{u_0}^2} - \mathbf{S}(\eta)^* \Bigr)\varphi  \Big\|_{L^2_+} 
+ \Big\|\tfrac{i}{\eta} \left(\mathbf{S}(\eta)^* f - f \right)  
- \mathbf{G} f \Big\|_{L^2_+} \to 0,
\end{split}
\end{equation}as $\eta \to 0^+$,  by Proposition 
$\ref{GisInfiniGeneHsProp}$
and formula 
\eqref{LimforconjSeta-SetaCV}. Then 
$e^{it \mathbf{L}_{u_0}^2} \mathbf{W}(t)^{-1} (f) = \varphi \in \mathrm{Dom}(\mathbf{G})$.  So 
\begin{equation}\label{DGcapH2<WtUtDGH2}
  \mathrm{Dom}(\mathbf{G}) \cap  H^2_+  \subset  \mathbf{W}(t)e^{-it \mathbf{L}_{u_0}^2}\left( \mathrm{Dom}(\mathbf{G}) \cap  H^2_+ \right) .
\end{equation}Formula \eqref{DomGcapH2InvWU} is obtained by   \eqref{WtUtDGcapH2subsetDG}
and \eqref{DGcapH2<WtUtDGH2}.
\end{proof}

\begin{rem}\label{WtstDomGHtwoisRem}
Since $\mathbf{W}(t)$ is a unitary operator on $L^2_+$, formula \eqref{DomGcapH2InvWU} is equivalent to
\begin{equation}\label{DomGcapH2InvWUequiv}
 \mathbf{W}(t)^{-1}\bigl( \mathrm{Dom}(\mathbf{G}) \cap H^2_+ \bigr) = e^{-it \mathbf{L}_{u_0}^2}\bigl( \mathrm{Dom}(\mathbf{G}) \cap H^2_+ \bigr) , \qquad \forall t \in (T^-, T^+).
\end{equation}
Although no explicit formula for the unitary propagator $\mathbf{W}(t) \in \mathcal{B}_{\mathbb{C}} (L^2_+)$ is available, the identity \eqref{DomGcapH2InvWUequiv} provides precise information on how the subspace $\mathrm{Dom}(\mathbf{G}) \cap H^2_+$ is rotated by $\mathbf{W}(t)^{-1}$:
this rotation is described explicitly by the action of the known unitary operator $e^{-it \mathbf{L}_{u_0}^2} \in \mathcal{B}_{\mathbb{C}}(L^2_+)$ on the same subspace.
\end{rem}

\begin{rem}\label{RelyShSeGrComForem} 
The proofs of Lemma~\ref{LemConjSstcvto2tLuphi} and Theorem~\ref{invSpDomGHtwoThm} rely exclusively on the commutator estimates 
\eqref{commForSetasttildBuCv} and \eqref{SetasttildBuUniEstL2} for the shift semigroup $(\mathbf{S}(\eta)^*)_{\eta \ge 0}$. 
In particular, the derivation of formulas \eqref{LimforconjSeta-SetaCV}, \eqref{DomGcapH2InvWU}, and \eqref{GWtUtisWtUtGadd2tL} makes no use of the commutator formula \eqref{commbfGtildBuFor} for the unbounded generator $\mathbf{G}$.
\end{rem}

\begin{rem}\label{remcompariGerChen}
Formulas~\eqref{DomGcapH2InvWU} and~\eqref{GWtUtisWtUtGadd2tL} provide an exact linearization of the conjugate of $\mathbf{G}$ by the unitary propagator $\mathbf{W}(t)e^{-it \mathbf{L}_{u_0}^2}$ along the flow of the focusing CMdNLS equation~\eqref{CMdNLSfocu}.
We now compare the approach developed in the present paper with the arguments originally used by G\'erard~\cite{BOExp2022} and subsequently adopted by Chen~\cite{Xchenccm25ExpFor}    
to derive explicit formulas for the Benjamin--Ono equation~\eqref{BOeq} on $\mathbb{R}$, and for the defocusing Calogero--Moser derivative NLS equation~\eqref{defocuCMdNLS} under the nonzero boundary condition \eqref{NonzeBoIntr}.  
  In those works, the authors introduced an analogous linearization identity 
 (see~\cite[the formulas below Lemma 9 on p.~599]{BOExp2022} and~\cite[the formulas below (4.2) on p.~19 $\&$ 20 of the arXiv version]{Xchenccm25ExpFor}):
\begin{equation}\label{GerardIdeUstGU}
U(t)^{*} \mathbf{G} U(t) = 2 t L_{u_0} + \mathrm{e}^{-i t L_{u_0}^{2}} \mathbf{G} \mathrm{e}^{i t L_{u_0}^{2}},
\end{equation} 
where $U(t)$ denotes the unitary propagator playing the same role as $\mathbf{W}(t) \in \mathcal{B}_{\mathbb{C}}(L^2_+)$ in the present paper, and $L_{u_0}$ is the corresponding Lax operator. 
We emphasize, however, that the identity~\eqref{GerardIdeUstGU} involves unbounded operators on both sides, and the works~\cite{BOExp2022, Xchenccm25ExpFor} did not discuss whether the domains of the left- and right-hand sides coincide. 
 The derivation of \eqref{GerardIdeUstGU} in those works proceeds by differentiating the left-hand side with respect to time and relating the resulting expression to the commutator $[\mathbf{G}, B_{u(t)}]$, where $(L_u, B_u)$ denotes the Lax pair of the corresponding equation. 
More precisely, one formally computes
\begin{equation}\label{formalDerivGer}
\begin{split}
\tfrac{\mathrm{d}}{\mathrm{d}t} \bigl( U(t)^* \mathbf{G} U(t) \bigr)
&= \lim_{r \to 0} \tfrac{1}{r} \bigl( U(t + r)^* \mathbf{G} U(t + r) - U(t)^* \mathbf{G} U(t) \bigr) \\
&= U(t)^* [\mathbf{G}, B_{u(t)}] U(t)  
= 2 L_{u_0} + i [U(t)^* \mathbf{G} U(t), L_{u_0}^2],
\end{split}
\end{equation}
and then integrates the resulting differential equation \eqref{formalDerivGer}  to obtain \eqref{GerardIdeUstGU}. 
  However, this differentiation is performed on an operator expression involving the unbounded operator $\mathbf{G}$, without verifying that  $\mathrm{Dom}(\mathbf{G})$ given by Definition $\ref{DombfGdefSect2}$ is preserved along the trajectory, nor that the limit of the difference quotient $\lim_{r \to 0} \tfrac{1}{r} \bigl( U(t + r)^* \mathbf{G} U(t + r) - U(t)^* \mathbf{G} U(t) \bigr)$  converges in the appropriate operator topology. 
More precisely, the formal differentiation process \eqref{formalDerivGer} suffers from the following logical gaps.
 
 \begin{enumerate}
    \item For the expression $U(t)^* \mathbf{G} U(t)\phi \in L^2_+$ to be well-defined, one must ensure that for a sufficiently large set of vectors $\phi \in L^2_+$, the orbit $U(t)\phi$ stays inside $\mathrm{Dom}(\mathbf{G})$ for all $t$. When this condition is violated, the term $\mathbf{G} U(t) \phi$ is not defined, and the entire computation that follows is meaningless. The works \cite{BOExp2022, Xchenccm25ExpFor} do not address this point.

    \item Even if $U(t)^* \mathbf{G} U(t)\phi \in L^2_+$ is well-defined, its time derivative can only be taken when the map $t\mapsto U(t)^* \mathbf{G} U(t)\phi$ is differentiable, and the time derivative $\partial_t$ commutes with the unbounded operator $\mathbf{G}$. 
    Recall that $\lim_{r \to 0} \|\frac{1}{r} (U(t+r) - U(t)) - B_{u(t)} U(t)\|_{\mathcal{B}_{\mathbb{C}}(L^2_+)} = 0$     
    for all $t$ in its maximal interval of   existence.     
    However, as demonstrated in Remark~\ref{AdvTSemiEstGenrem}, the discontinuity of $\mathbf{G}$ implies that there might exist a function $\varphi \in L^2_+$ for which 
\begin{equation*}
    \mathbf{G} \bigl( \tfrac{1}{r} (U(t+r) - U(t)) - B_{u(t)} U(t) \bigr) \varphi
\end{equation*}
    does not converge to $0$ in $L^2_+$. Hence one cannot assert that $\partial_t$ commutes with $\mathbf{G}$ without further justification.
    None of these verifications is carried out in \cite{BOExp2022, Xchenccm25ExpFor}; their computation \eqref{formalDerivGer} is performed as if $\mathbf{G}$ were a bounded operator on $L^2_+$.

    \item The commutator $[U(t)^* \mathbf{G} U(t), L_{u_0}^2]$ that appears in the last equality of \eqref{formalDerivGer} involves two unbounded operators: $U(t)^* \mathbf{G} U(t)$, whose natural domain is $U(t)^{-1}(\mathrm{Dom}(\mathbf{G}))$, and $L_{u_0}^2$, whose domain is $\mathrm{Dom}(L_{u_0}^2) = H^2_+$. 
For such a commutator to be well-defined, one must verify both compatibility conditions 
$\operatorname{Ran}(U(t)^* \mathbf{G} U(t)) \subset \mathrm{Dom}(L_{u_0}^2)$ and 
$\operatorname{Ran}(L_{u_0}^2) \subset U(t)^{-1}(\mathrm{Dom}(\mathbf{G}))$. 
Equivalently, one should ensure that the intersection $U(t)^{-1}(\mathrm{Dom}(\mathbf{G})) \cap \mathrm{Dom}(L_{u_0}^2)$ is sufficiently large. 
The works \cite{BOExp2022, Xchenccm25ExpFor} do not discuss this issue.
 \end{enumerate}

 \noindent Similar issues arise in other works that employ the same approach to
derive explicit (flow) formulas.  For instance, the linearization step in the
study of the half-wave maps equation in G\'erard--Lenzmann~\cite{GerLenz2024HWM}
(see~\cite[formula~(5.3) on p.~28 of the arXiv version 3]{GerLenz2024HWM}) and
in the treatment of the cubic Szeg\H{o} equation in
G\'erard--Pushnitski~\cite{GerPush2023} (see~\cite[the formulas below (2.10)
on p.~10 of the arXiv version]{GerPush2023}) follows the same formal
differentiation of an unbounded operator without addressing the domain issues
discussed above.  Moreover, the same formal derivation of formula \eqref{GerardIdeUstGU},
without any discussion of its domain of validity, is reproduced in the
recent ICM proceedings article G\'erard~\cite{GerICM2026Proc} (see formulas
(3.11)--(3.12) on p.~88).  This confirms that the logical gaps identified
in items (1)--(3) above are not confined to the original
papers~\cite{BOExp2022, Xchenccm25ExpFor}, but reflect a systematic
difficulty inherent in the direct differentiation method that leads
from~\eqref{formalDerivGer} to~\eqref{GerardIdeUstGU}. \\

 \noindent In contrast, the proof of Theorem~\ref{invSpDomGHtwoThm} in this paper circumvents all of the above difficulties by working entirely with the contraction semigroup $(\mathbf{S}(\eta)^*)_{\eta \geq 0}$ instead of the unbounded generator $\mathbf{G}$. 
No derivative of an expression containing $\mathbf{G}$ is ever taken; instead, the computations in Lemma~\ref{LemConjSstcvto2tLuphi} and Theorem~\ref{invSpDomGHtwoThm} are based solely on difference quotients composed with the bounded operators $\mathbf{S}(\eta)^* \in \mathcal{B}_{\mathbb{C}}(L^2_+)$, a procedure that is perfectly legitimate (see, for instance, formula~\eqref{dtbfnuFor}). The unbounded operator $\mathbf{G}$ is identified only at the very end, through the strong convergence provided by the Hille--Yosida theorem and formula~\eqref{SemiGrCVdefDomsG} (see~\eqref{UtWt-1Setaf-fCVL2} and~\eqref{CVSetavPhiInvSp}).
In this way we obtain, in a completely rigorous fashion, the invariant subspace property~\eqref{DomGcapH2InvWU} and the exact commutation relation~\eqref{GWtUtisWtUtGadd2tL}. 
The latter identity is the precise analogue of the linearization identity~\eqref{GerardIdeUstGU} in our context, with the essential additional information that its natural domain of validity is exactly $\mathrm{Dom}(\mathbf{G}) \cap H^2_+$ (under the assumption that $u$ is a $\mathcal{Z}^3_+$-solution to \eqref{CMdNLSfocu}).  
This crucial information about the domain is invisible in the formal manipulation of the earlier works~\cite{BOExp2022, Xchenccm25ExpFor, SUNLMP2024, SUNMSzegoR, GerPush2023, GerLenz2024HWM, GerICM2026Proc}. \\

\noindent
Moreover, the semigroup method developed in this paper is not tied to the specific structure of the CMdNLS or Benjamin--Ono equations. 
It relies on general principles of functional analysis — generators of contraction semigroups, the Hille--Yosida theorem, and the characterization of domains via strong derivatives — together with commutator estimates tailored to the specific equation. 
As a consequence, it can be systematically applied to other integrable equations posed on the Hardy space $L^2_+$ that admit a Lax pair formulation. 
This provides a new, rigorous paradigm for the derivation of explicit formulas and the identification of invariant subspaces in such settings, replacing the formal differentiation of unbounded operators that has been used in the past.
\end{rem}

\noindent Lemma $\ref{LemWtexp-itL2phiPDE}$  characterizes the forward propagator $\mathbf{W}(t)e^{-it\mathbf{L}_{u_0}^2} (\varphi) $ as the unique solution to the linear evolution equation \eqref{PDEdtfistilButf} driven by $\widetilde{\mathbf{B}}_{u(t)}$. We now study the inverse propagator $e^{it\mathbf{L}_{u_0}^2}\mathbf{W}(t)^{-1}$.

\begin{lem}
For any $t \in (T^-, T^+)$ and  $\varphi \in H^2_+$, the following identity holds:
\begin{equation}\label{NewLeibeitL2W-1For}
e^{it \mathbf{L}_{u_0}^2} \mathbf{W}(t)^{-1}  \varphi  - \varphi = - \int_0^t  e^{i \mathfrak{s} \mathbf{L}_{u_0}^2} \mathbf{W}(\mathfrak{s})^{-1}  \widetilde{\mathbf{B}}_{u( \mathfrak{s} )}    ( \varphi) \mathrm{d}\mathfrak{s} \in H^2_+.
\end{equation}Recall that $\chi_{\epsilon} \in \mathrm{Dom}^{\infty}(\mathbf{G})$ is given by  \eqref{chiepFor}, $\forall \epsilon >0$. Then we have
\begin{equation}\label{CVeitL2W-1chiep}
\lim_{\epsilon \to 0^+} \Big\|e^{it \mathbf{L}_{u_0}^2} \mathbf{W}(t)^{-1}  \chi_{\epsilon} - \chi_{\epsilon}  \Big\|_{L^2_+} = 0.
\end{equation}
\end{lem}

\begin{proof}
For any $t \in (T^-, T^+)$, we have $\mathbf{W}(t), \mathbf{W}(t)^{-1} \in \bigcap_{p = 0}^2 \mathcal{B}_{\mathbb{C}}(H^p_+)$ by    formulas  \eqref{UniPropagbfWSect4} and   \eqref{bfWstarisCone}.
Lemma $\ref{expitLusqStrContGroup}$ implies that $\mathbf{U}(-t) = e^{it \mathbf{L}_{u_0}^2}  \in \bigcap_{p = 0}^2 \mathcal{B}_{\mathbb{C}}(H^p_+)$ and $\left( \mathbf{U}( \mathbf{t} ) \right)_{\mathbf{t} \in \mathbb{R}}$ is a strongly continuous  unitary group on $L^2_+$. 
Fix $\varphi \in H^2_+$. Then $\mathbf{W}(t)^{-1} \varphi \in H^2_+$. For any $t \in (T^-, T^+)$, we set 
\begin{equation}\label{bfguphgitDefFor}
\mathbf{g}(t) = \mathbf{g}_{u, \varphi}(t) : = e^{it \mathbf{L}_{u_0}^2} \mathbf{W}(t)^{-1}  \varphi \in H^2_+.
\end{equation}For any $t \in (T^-, T^+)$, formula  \eqref{dtbfWt=ButbfWtlim}   implies that 
\begin{equation}\label{CVForUtrUtovr}
\small
\begin{split}
& \Big\| \tfrac{1}{r} \left(\mathbf{U}( -r) - \mathrm{id}_{L^2_+} \right) \left( \mathbf{W}(t+r)^{-1}  -\mathbf{W}(t)^{-1}  \right) \varphi \Big\|_{L^2_+} \\ 
\leq & 2 \Big\|   \left( \tfrac{1}{r} \left( \mathbf{W}(t+r)^{-1}  -\mathbf{W}(t)^{-1}  \right) + \mathbf{W}(t)^{-1} \mathbf{B}_{u(t)} \right)\varphi \Big\|_{L^2_+} + 
\Big\| \left(\mathbf{U}(-r) -  \mathrm{id}_{L^2_+} \right) \mathbf{W}(t)^{-1} \mathbf{B}_{u(t)} ( \varphi ) \Big\|_{L^2_+} \to 0,
\end{split}
\end{equation}as $r \to 0$.  
Using  formulas 
\eqref{tildeBudefFoc-}, \eqref{dtbfWt=ButbfWtlim}  
and  \eqref{UniEquLunFor}, we obtain that 
\begin{equation}\label{quodiffgdxgFor}
\small
\begin{split}
& \Big\| \tfrac{1}{r} \left( \mathbf{g}(t+r)- \mathbf{g}(t) \right) + 
e^{it \mathbf{L}_{u_0}^2} \mathbf{W}(t)^{-1}  \widetilde{\mathbf{B}}_{u(t)} (  \varphi  )  \Big\|_{L^2_+} \\
\leq & \Big\| \tfrac{1}{r} \left(\mathbf{U}( -r) - \mathrm{id}_{L^2_+} \right) \left( \mathbf{W}(t+r)^{-1}  -\mathbf{W}(t)^{-1}  \right) \varphi \Big\|_{L^2_+}    +  \Big\| \left( \tfrac{1}{r} \left( \mathbf{U}(- r ) - \mathrm{id}_{L^2_+} 
 \right)   - i \mathbf{L}_{u_0}^2 \right) \mathbf{W}(t)^{-1} \varphi \Big\|_{L^2_+}  \\
&  \qquad + \Big \|   \left( \tfrac{1}{r} \left( \mathbf{W}(t+r)^{-1}  -\mathbf{W}(t)^{-1}  \right) + \mathbf{W}(t)^{-1} \mathbf{B}_{u(t)} \right)\varphi \Big\|_{L^2_+} \\
\to  & 0,
\end{split}
\end{equation}as $r \to 0$. So $\mathbf{g}  : t \in (T^-, T^+) \mapsto \mathbf{g}(t) \in L^2_+$ is differentiable.  
Corollary $\ref{tildButtotilButcontiCor}$ implies that the mapping
\begin{equation}\label{dtgExpFore-itL2Wt-1}
\partial_t \mathbf{g} :(T^-, T^+) \mapsto \partial_t \mathbf{g}(t) = - e^{it \mathbf{L}_{u_0}^2} \mathbf{W}(t)^{-1}  \widetilde{\mathbf{B}}_{u(t)} (  \varphi  )  \in    L^2_+ 
\end{equation}is continuous. So it is Bochner integrable over any compact subset of $(T^-, T^+)$ by 
 Theorem $3.27$ of   \cite{RudinFA}. Since $\mathbf{g} \in \mathcal{C}^1 \left( (T^-, T^+) ; L^2_+ \right)$, formula  \eqref{NewLeibeitL2W-1For}
follows from \eqref{dtgExpFore-itL2Wt-1}. Fix $t \in (T^-, T^+) \setminus \{0\}$. We set 
\begin{equation}\label{frakmepsDefFor}
\mathfrak{m}( \epsilon , \mathfrak{s}) : = \|\widetilde{\mathbf{B}}_{u(\mathfrak{s})} (  \chi_{\epsilon}  ) \|_{L^2_+} \geq 0, \quad \forall \epsilon \in (0, 1], \quad \forall \mathfrak{s} \in \mathtt{I}_t,
\end{equation}where  $\mathtt{I}_t$ is given by \eqref{calItDef}.  
By   \eqref{wtilBuchiepto0Uni} and \eqref{frakmepsDefFor}, we have  
$\lim_{\epsilon \to 0^+} \mathfrak{m}( \epsilon , \mathfrak{s}) = 0$,  $ \forall \mathfrak{s} \in \mathtt{I}_t$, and 
\begin{equation}
 \sup_{\mathfrak{s} \in \mathtt{I}_t} \sup_{0 < \epsilon \leq 1} \mathfrak{m}( \epsilon , \mathfrak{s})
 = \sup_{\mathfrak{s} \in \mathtt{I}_t} \sup_{0 < \epsilon \leq 1} \|\widetilde{\mathbf{B}}_{u(\mathfrak{s})} (  \chi_{\epsilon}  ) \|_{L^2_+} \leq \sqrt{2 \pi} \left( 1 +  \sup_{\mathfrak{s} \in \mathtt{I}_t}  \|u(\mathfrak{s}) \|_{\mathcal{Z}^1_+}^2\right) < +\infty.
\end{equation}By  Theorem $3.29$ of \cite{RudinFA} and the dominated convergence theorem, it follows that 
\begin{equation*}
\Big\|e^{it \mathbf{L}_{u_0}^2} \mathbf{W}(t)^{-1}  \chi_{\epsilon} - \chi_{\epsilon}  \Big\|_{L^2_+} \leq \int_{\mathtt{I}_t}  \Big\|  e^{i \mathfrak{s} \mathbf{L}_{u_0}^2} \mathbf{W}( \mathfrak{s} )^{-1}  \widetilde{\mathbf{B}}_{u( \mathfrak{s} )} (  \chi_{\epsilon}   )  \Big\|_{L^2_+} \mathrm{d} \mathfrak{s} = \int_{\mathtt{I}_t} \mathfrak{m}( \epsilon , \mathfrak{s}) \mathrm{d} \mathfrak{s} \to 0,
\end{equation*}as $\epsilon \to 0^+$. The proof is completed.
\end{proof}

\noindent The next two lemmas prepare the last convergence arguments needed for the derivation of the explicit formula \eqref{ExpForCMdNLSfoc}. 
 Lemma~\ref{ttkuqContfct} establishes the continuity of the inverse propagator $e^{it\mathbf{L}_{u_0}^2}\mathbf{W}(t)^{-1}$  along time-dependent curves, while Lemma~\ref{invPropagutchilim} uses this continuity together with the dynamics of   solution 
\begin{equation}\label{recaluisZ3sol}
u\in\mathcal{C}\bigl((T^-,T^+);\mathcal{Z}^3_+\bigr) \cap \mathcal{C}^1 \bigl((T^-,T^+);\mathcal{Z}^1_+\bigr)  \qquad \text{with} \qquad  u(0) = u_0 \in \mathcal{Z}^3_+,
\end{equation}
 to compute the limit of $e^{it \mathbf{L}_{u_0}^2} \mathbf{W}(t)^{-1}(u(t)\chi_{\epsilon})$ as $\epsilon \to 0^+$.

\begin{lem}\label{ttkuqContfct}
For any $\mathrm{q} \in \mathcal{C}  \bigl( (T^-, T^+) ; H^2_+ \bigr)$,  we set
 $\mathtt{k}_{u , \mathrm{q}}(t) : = \exp (it \mathbf{L}_{u_0}^2) \mathbf{W}(t)^{-1} \bigl( \mathrm{q} (t) \bigr) \in H^2_+$, for any $ t \in (T^-, T^+)$. Then $\mathtt{k}_{u , \mathrm{q}}$ is a continuous function from $(T^-, T^+)$ to $H^2_+$.
\end{lem}
\begin{proof}
For any $ t \in (T^-, T^+)$, set $\mathfrak{l}_{u , \mathrm{q}}(t) : =\mathbf{W}(t)^{-1} \bigl( \mathrm{q} (t) \bigr) \in H^2_+ $. Then  formula \eqref{bfWstarisCone}
 implies that 
\begin{equation}\label{frakluqconti}
\mathfrak{l}_{u , \mathrm{q}} \in \mathcal{C}  \bigl( (T^-, T^+) ; H^2_+ \bigr).
\end{equation}Lemma $\ref{expitLusqStrContGroup}$ shows that $\bigl( \mathbf{U}(   \mathbf{t}) \bigr)_{\mathbf{t} \in \mathbb{R}}$ is a strongly continuous group on $H^2_+$, where
$\mathbf{U}( \mathbf{t} ) = \exp (  - i\mathbf{t} \mathbf{L}_{u_0}^2)$. So 
$ \mathtt{k}_{u , \mathrm{q}} (t ) = \mathbf{U}( - t )  \mathfrak{l}_{u , \mathrm{q}}(t) $, for every $t \in (T^-, T^+) $.
 By estimate \eqref{EstiUthH2norm}, we have 
 \begin{equation}\label{EstunibfUtinRHtwo}
 \sup_{\mathbf{t}\in \mathbb{R}}  \big\|\mathbf{U}( \mathbf{t} ) \big\|_{\mathcal{B}_{\mathbb{C}}(H^2_+)}   
\leq   2 \left( 1 + 2 (1+ C_1)^2 \| u_0 \|_{\mathcal{Z}^1}^2 \right)^2
 \end{equation}Using  \eqref{frakluqconti}, \eqref{EstunibfUtinRHtwo} and the strong continuity of $\bigl( \mathbf{U}(   \mathbf{t}) \bigr)_{\mathbf{t} \in \mathbb{R}}$ on $H^2_+$, we obtain that  
 \begin{equation}
\begin{split}
 & \|\mathtt{k}_{u , \mathrm{q}}(t + r) - \mathtt{k}_{u , \mathrm{q}}(t) \|_{H^2_+} \\
 \leq & \sup_{\mathbf{t}\in \mathbb{R}} \big\|\mathbf{U}( \mathbf{t} ) \big\|_{\mathcal{B}_{\mathbb{C}}(H^2_+)}   \left( \| \mathfrak{l}_{u , \mathrm{q}} (t + r) -  \mathfrak{l}_{u , \mathrm{q}} (t) \|_{H^2_+}  + \|  \mathbf{U}( - \mathbf{r} ) \mathfrak{l}_{u , \mathrm{q}} (t)  - \mathfrak{l}_{u , \mathrm{q}} (t) \|_{H^2_+} \right)  
\to  0 ,
\end{split}
\end{equation}as $r \to 0$. So $\mathtt{k}_{u , \mathrm{q}}$ is continuous at $t$, for any $ t \in (T^-, T^+)$. 
\end{proof}

\noindent With the continuity of the inverse propagator $e^{it\mathbf{L}_{u_0}^2}\mathbf{W}(t)^{-1}$ acting on continuous curves established, we can now prove the key convergence result needed for the explicit formula.

\begin{lem}\label{invPropagutchilim}
For any  $\epsilon >0$, let $\chi_{\epsilon} \in \mathrm{Dom}^{\infty}(\mathbf{G})$ be defined in   \eqref{chiepFor}. Then  we have
\begin{equation}\label{eitLusqWstutchiepLim}
\lim_{\epsilon \to 0^+} \big\| \exp (it \mathbf{L}_{u_0}^2) \mathbf{W}(t)^{-1} \bigl( u (t)  \chi_{\epsilon}  \bigr) - u_0 \chi_{\epsilon} \big\|_{L^2_+} = 0 ,  \qquad \forall t \in (T^-, T^+) .
\end{equation} 
\end{lem}

\begin{proof}
For any $\epsilon >0$ and $t \in (T^-, T^+)$, we set 
  $\mathfrak{k}_{u, \epsilon}   (t) : = \mathtt{k}_{u , u\chi_{\epsilon}}(t) = \exp (it \mathbf{L}_{u_0}^2) \mathbf{W}(t)^{-1} \bigl(  u(t) \chi_{\epsilon} \bigr) \in H^2_+$. Recall that $u\in\mathcal{C}\bigl((T^-,T^+);\mathcal{Z}^3_+\bigr) \cap \mathcal{C}^1 \bigl((T^-,T^+);\mathcal{Z}^1_+\bigr)$.  Let $\mathbf{r}_{\epsilon}(u(t)) \in H^1_+$ be defined in \eqref{bfreputilBuHardZhi}.
 Using formulas  \eqref{bfreputilBuHardZhi}, \eqref{CMdNLSwithPi} and  
Lemma  $\ref{repuHighRegLem}$, we obtain that  
\begin{equation}\label{bfreputdtutchitiB}
\mathbf{r}_{\epsilon}(u(t)) = \chi_{\epsilon} \partial_t u(t) - \widetilde{\mathbf{B}}_{u(t)} (u(t) \chi_{\epsilon}) \in H^2_+ .
\end{equation}
 From formulas \eqref{soluqudifCVSob}, 
\eqref{bfguphgitDefFor},  
\eqref{quodiffgdxgFor} and  \eqref{bfreputdtutchitiB}, it follows that  
\begin{equation}\label{varkadifeitLWstrep}
\begin{split}
 & \big\|\tfrac{1}{r} \bigl( \mathfrak{k}_{u, \epsilon}   (t+r) - \mathfrak{k}_{u, \epsilon}   (t) \bigr) - \exp(it \mathbf{L}_{u_0}^2 )\mathbf{W}(t)^{-1} \bigl( \mathbf{r}_{\epsilon}(u(t)) \bigr) \big\|_{L^2_+} \\
 \leq &  \big\| \tfrac{1}{r} \bigl( e^{i(t+r) \mathbf{L}_{u_0}^2} \mathbf{W}(t + r)^{-1}   -  e^{it \mathbf{L}_{u_0}^2} \mathbf{W}(t)^{-1} \bigr)(u(t)\chi_{\epsilon})  
 +  e^{it \mathbf{L}_{u_0}^2} \mathbf{W}(t)^{-1}  \widetilde{\mathbf{B}}_{u(t)} (u(t)\chi_{\epsilon})\big\|_{L^2_+}
  \\
 & \qquad + \big\|  e^{i(t+r) \mathbf{L}_{u_0}^2} \mathbf{W}(t + r)^{-1} \bigl(  \tfrac{\chi_{\epsilon}}{r}\bigl( u(t+r) -u(t)  \bigr)   - \chi_{\epsilon} \partial_t u(t) \bigr) \big\|_{L^2_+}
  \\
 & \qquad + \big\| e^{i(t+r) \mathbf{L}_{u_0}^2} \mathbf{W}(t + r)^{-1} (\chi_{\epsilon} \partial_t u(t)) -  e^{it \mathbf{L}_{u_0}^2} \mathbf{W}(t  )^{-1} (\chi_{\epsilon} \partial_t u(t)) \big\|_{L^2_+}\\
 \to & 0, \quad \text{as} \quad r \to 0.
\end{split}
\end{equation}Then $\mathfrak{k}_{u, \epsilon}$ is differentiable on $(T^-, T^+)$ and $\partial_t \mathfrak{k}_{u, \epsilon} (t) = \exp(it \mathbf{L}_{u_0}^2 )\mathbf{W}(t)^{-1} \bigl( \mathbf{r}_{\epsilon}(u(t)) \bigr) \in H^2_+$, $\forall t \in (T^-, T^+)$. Formula  \eqref{bfrepContonu} and Lemma $\ref{ttkuqContfct}$ imply that $ \mathfrak{k}_{u, \epsilon} \in \mathcal{C}^1 \bigl( (T^-, T^+) ; L^2_+ \bigr)$. Consequently, we have  
\begin{equation}\label{eitLsqWstutchiNeLe}
 e^{it \mathbf{L}_{u_0}^2}\mathbf{W}(t)^{-1} \bigl(  u(t) \chi_{\epsilon} \bigr) -  u_0 \chi_{\epsilon}  
  = \int_0^t \partial_t \mathfrak{k}_{u, \epsilon} (\mathfrak{s})  \mathrm{d}\mathfrak{s}
  = \mathrm{sign}(t) \int_{\mathtt{I}_t}  e^{ i \mathfrak{s} \mathbf{L}_{u_0}^2 } \mathbf{W}(\mathfrak{s})^{-1} \bigl( \mathbf{r}_{\epsilon}(u(\mathfrak{s})) \bigr) \mathrm{d}\mathfrak{s},
\end{equation}where $\mathtt{I}_t$ is given by \eqref{calItDef}.
Since $e^{ i \mathfrak{s} \mathbf{L}_{u_0}^2 } \mathbf{W}(\mathfrak{s})^{-1} $ is a unitary operator on $L^2_+$,  we have
\begin{equation*}
\limsup_{\epsilon \to 0^+} \big\|  e^{it \mathbf{L}_{u_0}^2}\mathbf{W}(t)^{-1} \bigl(  u(t) \chi_{\epsilon} \bigr) -  u_0 \chi_{\epsilon}  \big\|_{L^2_+} 
\leq  \lim_{\epsilon \to 0^+} \int_{\mathtt{I}_t} \big\| \mathbf{r}_{\epsilon}(u(\mathfrak{s}))\big\|_{L^2_+} \mathrm{d}\mathfrak{s} =  0  , \quad \forall t \in (T^-, T^+), 
\end{equation*}by formulas 
\eqref{eitLsqWstutchiNeLe}, 
\eqref{bfrepyto0inL2} and the dominated convergence theorem.
\end{proof}

 \subsection{Derivation of the Explicit Formulas}\label{DerExpForSubs4.3}

\noindent This subsection is devoted to proving Theorem~\ref{ExpForMainThm}
and Corollary~\ref{CorExpforaaddH}.  After the proof, we compare the semigroup
method developed in this paper with the approach used by
Killip--Laurens--Vi\c{s}an~\cite{KillLV2311} to derive the explicit formula for solutions in $H^{\infty}_+ \cap \mathrm{Dom} (\mathbf{X})$.
We begin with the proof of Theorem~\ref{ExpForMainThm}.

\begin{proof}[Proof of Theorem $\ref{ExpForMainThm}$]

The proof is divided into two steps. First, we establish formula~\eqref{ExpForCMdNLSfoc} under the additional assumption that $u$ is a $\mathcal{Z}^3_+$-solution to \eqref{CMdNLSfocu} with $u_0 \in \mathcal{Z}^3_+$. Then, since $\mathcal{Z}^3_+$ is dense in $(\mathcal{Z}^2_+, \| . \|_{\mathcal{Z}^2_+})$ by Lemma~\ref{HarZhiSmoodense} and the right-hand side of~\eqref{ExpForCMdNLSfoc} defines a continuous functional on $\mathcal{Z}^2_+$ by Lemma~\ref{lemContRHSExpFor}, the formula extends to all $\mathcal{Z}^2_+$-solutions by a standard continuity argument. \\

\noindent \textbf{Step \romannumeral 1}.   
Assume that $T^- < 0 <T^+$ and $u\in\mathcal{C}\bigl((T^-,T^+);\mathcal{Z}^3_+\bigr) \cap \mathcal{C}^1 \bigl((T^-,T^+);\mathcal{Z}^1_+\bigr)$ is a solution to \eqref{CMdNLSfocu} with   $u(0)=u_0\in\mathcal{Z}^3_+$. For any $\epsilon >0$, let $\chi_{\epsilon} \in \mathrm{Dom}^{\infty} (\mathbf{G})$ be defined in \eqref{chiepFor}. For any $t \in (T^-, T^+)$, 
the Lax pair $\bigl(\mathbf{L}_{u(t)}, \mathbf{B}_{u(t)} \bigr)$ is given by   \eqref{LuSect2} and
\eqref{BudefForFoc-}. Let $\mathbf{W} (t) \in \bigcap_{p=0}^2  \mathcal{B}_{\mathbb{C}}(H^p_+)$ be the unitary propagator defined   in \eqref{UniPropagbfWSect4}.  For any $\epsilon >0$, we have $u(t) \chi_{\epsilon} \in H^3_+$. Fix $z \in \mathbb{C}_+$.  We set
\begin{equation}\label{bfqeptzdef}
\mathbf{q}^{(\epsilon)}_{t, z} : =  (\mathbf{G} - z )^{-1} \bigl( u(t) \chi_{\epsilon} \bigr) \in \mathrm{Dom}^3 (\mathbf{G})  , \qquad \boldsymbol{\phi}^{(\epsilon)}_{t, z} : = e^{ i t \mathbf{L}_{u_0}^2 } \mathbf{W}(t)^{-1} \bigl( \mathbf{q}^{(\epsilon)}_{t, z} \bigr) \in  L^2_+.
\end{equation}Formula  \eqref{DomsGDefSect2} implies that $\mathrm{Dom}^3 (\mathbf{G}) \subset  \mathrm{Dom}(\mathbf{G}) \cap H^2_+$.
By  \eqref{bfqeptzdef} and  \eqref{DomGcapH2InvWU}, we  have
\begin{equation}\label{bdsphiinDGcapH2}
\boldsymbol{\phi}^{(\epsilon)}_{t, z} = e^{ i t \mathbf{L}_{u_0}^2 } \mathbf{W}(t)^{-1} (\mathbf{G} - z )^{-1} \bigl( u(t) \chi_{\epsilon} \bigr) \in e^{ i t \mathbf{L}_{u_0}^2 } \mathbf{W}(t)^{-1} \bigl( \mathrm{Dom}^3 (\mathbf{G})  \bigr)  \subset  \mathrm{Dom}(\mathbf{G}) \cap H^2_+.
\end{equation}Then \eqref{sumbfG2tLubbGt} yields that 
$\bigl(\mathbf{G}  + 2 t \mathbf{L}_{u_0}  \bigr)   \boldsymbol{\phi}^{(\epsilon)}_{t, z} = \bigl(\mathbb{G}_t  -  2 t \mathbf{T}_{u_0} \mathbf{T}_{\overline{u}_0}  \bigr)   \boldsymbol{\phi}^{(\epsilon)}_{t, z} \in L^2_+$. Consequently, 
we obtain  
\begin{equation}\label{utchiepisWeitLGp}
\begin{split}
u(t) \chi_{\epsilon} 
= & (\mathbf{G} - z ) \mathbf{W}(t) e^{ - i t \mathbf{L}_{u_0}^2 }   \bigl( \boldsymbol{\phi}^{(\epsilon)}_{t, z} \bigr)
= \mathbf{W}(t) e^{ - i t \mathbf{L}_{u_0}^2 }   \bigl(\mathbf{G}  + 2 t \mathbf{L}_{u_0} - z \bigr)   \boldsymbol{\phi}^{(\epsilon)}_{t, z} \in H^3_+\\
= &  \mathbf{W}(t) e^{ - i t \mathbf{L}_{u_0}^2 }   \bigl(\mathbb{G}_t  -  2 t \mathbf{T}_{u_0} \mathbf{T}_{\overline{u}_0} - z \bigr)   \boldsymbol{\phi}^{(\epsilon)}_{t, z} ,
\end{split}
\end{equation}by formulas 
  \eqref{bdsphiinDGcapH2} and \eqref{GWtUtisWtUtGadd2tL}. Since $\operatorname{Im} z >0$, 
Lemma $\ref{LemPertGtMaxAcc}$ implies that  $\mathbb{G}_t  -  2 t \mathbf{T}_{u_0} \mathbf{T}_{\overline{u}_0} - z $ is a bijective mapping from $ e^{-it \partial_x^2} \bigl( \mathrm{Dom}  (\mathbf{G}) \bigr)$ to  $L^2_+$. Then it follows from \eqref{utchiepisWeitLGp} and \eqref{bdsphiinDGcapH2}  that 
\begin{equation}\label{sfRueptdefFor}
\begin{split}
\boldsymbol{\phi}^{(\epsilon)}_{t, z} 
= & \bigl(\mathbb{G}_t  -  2 t \mathbf{T}_{u_0} \mathbf{T}_{\overline{u}_0} - z \bigr)^{-1}  e^{  i t \mathbf{L}_{u_0}^2 }   \mathbf{W}(t)^{-1} \bigl(u(t) \chi_{\epsilon}  \bigr) \in \mathrm{Dom}(\mathbf{G}) \cap H^2_+ \\
  = & \bigl(\mathbb{G}_t  -  2 t \mathbf{T}_{u_0} \mathbf{T}_{\overline{u}_0} - z \bigr)^{-1}    \bigl( \mathsf{R}_{u, \epsilon} (t)  +  u_0 \chi_{\epsilon}  \bigr),
\end{split}
\end{equation}with $\mathsf{R}_{u, \epsilon} (t)  : = e^{  i t \mathbf{L}_{u_0}^2 }   \mathbf{W}(t)^{-1} \bigl(u(t) \chi_{\epsilon}  \bigr)  -  u_0 \chi_{\epsilon}  \in H^2_+$.
Formulas 
\eqref{eitLusqWstutchiepLim} and \eqref{CVeitL2W-1chiep} imply that
\begin{equation}\label{ForsfRurto0lim}
\lim_{r \to 0^+} \|\mathsf{R}_{u, r} (t) \|_{L^2_+} = 0, \qquad
\lim_{r \to 0^+} \langle \mathbf{q}^{(\epsilon)}_{t, z} \; , \;  \mathbf{W}(t) e^{ - i t \mathbf{L}_{u_0}^2 }   (\chi_r) - \chi_r \rangle_{L^2_+} =0 , \qquad \forall \epsilon > 0 .
\end{equation}
By Lemma $\ref{expitLusqStrContGroup}$ and formula  \eqref{bfWtisUnitaryinL2+},  both $e^{ i t \mathbf{L}_{u_0}^2 }$ and $\mathbf{W}(t)$ are unitary operators on $L^2_+$. 
According to  Lemma  $\ref{phichiepCVtoI0Lem}$, \eqref{InvForNewtonLeibniz}, 
\eqref{bfqeptzdef} and the second formula of \eqref{ForsfRurto0lim}, for any $\epsilon >0$,  we have 
\begin{equation}\label{PoiutchiepbfIphi}
\small 
2 \pi i \mathbf{P}_z \bigl( u(t) \chi_{\epsilon} \bigr)
= \mathbf{I}_+  \bigl( \mathbf{q}^{(\epsilon)}_{t, z} \bigr)
 = \lim_{r \to 0^+} \langle \mathbf{q}^{(\epsilon)}_{t, z} , \chi_r \rangle_{L^2_+}
 = \lim_{r \to 0^+} \langle e^{ i t \mathbf{L}_{u_0}^2 } \mathbf{W}(t)^{-1} \bigl( \mathbf{q}^{(\epsilon)}_{t, z} \bigr) , \chi_r \rangle_{L^2_+} = \mathbf{I}_+ \bigl( \boldsymbol{\phi}^{(\epsilon)}_{t, z} \bigr) .
\end{equation}
 Lemma $\ref{LembfIinvbbGtper}$
yields that 
$\mathbf{I}_+ \bigl(\mathbb{G}_t  -  2 t \mathbf{T}_{u_0} \mathbf{T}_{\overline{u}_0} - z \bigr)^{-1} \in \mathcal{B}_{\mathbb{C}} (L^2_+ ; \mathbb{C})$.
  By the first formula of \eqref{ForsfRurto0lim}, we have 
\begin{equation}\label{cvbfIGrpertsfR}
\lim_{\epsilon \to 0^+}\mathbf{I}_+ \bigl(\mathbb{G}_t  -  2 t \mathbf{T}_{u_0} \mathbf{T}_{\overline{u}_0} - z \bigr)^{-1}   \bigl( \mathsf{R}_{u, \epsilon} (t)    \bigr) =0 .
\end{equation}
Using 
\eqref{CvPoiIntuchieptou} and plugging  
\eqref{sfRueptdefFor} and 
\eqref{cvbfIGrpertsfR} into  formula \eqref{PoiutchiepbfIphi}, we obtain that
\begin{equation}\label{utzislimbfIGtper}
\begin{split}
\underline{u}(t,z ) = & \mathbf{P}_z \bigl( u(t) \bigr) 
= \lim_{\epsilon \to 0^+} \mathbf{P}_z \bigl( u(t) \chi_{\epsilon} \bigr) = \tfrac{1}{2 \pi i} \lim_{\epsilon \to 0^+}  \mathbf{I}_+ \bigl( \boldsymbol{\phi}^{(\epsilon)}_{t, z} \bigr) \\
  = & \tfrac{1}{2 \pi i} \lim_{\epsilon \to 0^+}  \mathbf{I}_+  \bigl(\mathbb{G}_t  -  2 t \mathbf{T}_{u_0} \mathbf{T}_{\overline{u}_0} - z \bigr)^{-1}    \bigl(  u_0 \chi_{\epsilon}  \bigr) .
\end{split}
\end{equation} 
  Finally, formula~\eqref{ExpForCMdNLSfoc} follows from~\eqref{utzislimbfIGtper} and~\eqref{limbfIbbGtpertInv}. Thus~\eqref{ExpForCMdNLSfoc} holds for all $\mathcal{Z}^3_+$-solutions to~\eqref{CMdNLSfocu}.\\

\noindent \textbf{Step \romannumeral 2}. Let $u\in\mathcal{C}\bigl((T^-,T^+);\mathcal{Z}^2_+\bigr)$ be a fixed maximal solution to   \eqref{CMdNLSfocu} with initial datum $u(0)=u_0\in\mathcal{Z}^2_+$. Here
$T^- = T_2^-(\|u_0\|_{\mathcal{Z}^2_+}) < 0 < T^+ = T_2^+(\|u_0\|_{\mathcal{Z}^2_+})$, and 
  $(T^-,T^+) = \mathcal{I}^2_{\max}(\|u_0\|_{\mathcal{Z}^2_+})$ denotes the maximal lifespan of  $u$ with respect to the $\mathcal{Z}^2_+$-norm, whose existence is guaranteed by the local well-posedness theory established in Proposition~\ref{LWPpropCMdNLSZn} together with the usual blow-up alternative. 
  Let $(\rho_m)_{m \in \mathbb{N}_+}$ be the sequence of mollifiers as defined on p.~108 of Brezis \cite{BrezisFA} (see the formulas above Proposition~4.21). For every $m \in \mathbb{N}_+$, we have $\rho_m \in \mathcal{C}_c^{\infty} (\mathbb{R} ; \mathbb{C})$, $0 \leq \rho_m \leq 1$, $\operatorname{supp} (\rho_m) \subset [- \frac{1}{m}  , \frac{1}{m} ]$ and $\int_{\mathbb{R}} \rho_m = 1$. Set
  \begin{equation}\label{defum0conv}
  u_0^m := \rho_m * u_0 \in \mathcal{Z}^{\infty}_+,
  \end{equation} 
where the inclusion $\rho_m * u_0 \in \mathcal{Z}^{\infty}_+$ follows from \eqref{HarZhidSpZnDefFor}, \eqref{SmoothZhidkHa}, 
Lemma $\ref{LinftplconvSchwLem}$, 
Theorem 6.30 of Rudin \cite{RudinFA} and 
 Young's convolution inequality. Moreover, we have the following 
uniform bound
\begin{equation}\label{molliun0Zhinor}
\sup_{m \in \mathbb{N}_+} \|u_0^m\|_{\mathcal{Z}^2_+} \leq \|u_0\|_{\mathcal{Z}^2_+},
\end{equation}by formula  \eqref{ZhidNormZn}. 
 Proposition 4.21 and Theorem 4.22 of Brezis \cite{BrezisFA} imply that
 \begin{equation}\label{umzeroCvtou0}
\lim_{m \to +\infty} \|u_0^m - u_0\|_{\mathcal{Z}^2_+} = 0.
\end{equation}
For each $m \in \mathbb{N}_+$, since $u_0^m \in \mathcal{Z}^{\infty}_+ \subset \mathcal{Z}^3_+$,
  Proposition~\ref{LWPpropCMdNLSZn} and Lemma~\ref{lemLinPartsol} yield  
  a unique maximal solution
\begin{equation}\label{umdefSol}
u^m \in \mathcal{C}\bigl(\mathcal{I}^3_{\max}(\|u_0^m\|_{\mathcal{Z}^3_+}); \mathcal{Z}^3_+\bigr) \cap \mathcal{C}^1\bigl(\mathcal{I}^3_{\max}(\|u_0^m\|_{\mathcal{Z}^3_+}); \mathcal{Z}^1_+\bigr)
\end{equation} 
to \eqref{CMdNLSfocu} with initial datum $u^m(0) = u_0^m$. Here $\mathcal{I}^3_{\max}(\|u_0^m\|_{\mathcal{Z}^3_+}) = (T_3^-(\|u_0^m\|_{\mathcal{Z}^3_+}), T_3^+(\|u_0^m\|_{\mathcal{Z}^3_+}))$ denotes its maximal lifespan with respect to the $\mathcal{Z}^3_+$-norm. Since $\mathcal{Z}^3_+ \subset \mathcal{Z}^2_+$, the function $u^m$ is also a $\mathcal{Z}^2_+$-solution; let $\mathcal{I}^2_{\max}(\|u_0^m\|_{\mathcal{Z}^2_+})$ denote its maximal lifespan with respect to the $\mathcal{Z}^2_+$-norm.  
 By the persistence of regularity for dispersive equations (a classical argument based on the blow-up criterion and the fact that the nonlinear term $4 u \Pi_{\geq 0} \bigl( \operatorname{Re} (\overline{u} \partial_x u ) \bigr)$ is a continuous polynomial in the solution and its derivatives up to order one), the $\mathcal{Z}^3_+$-maximal lifespan coincides with the $\mathcal{Z}^2_+$-maximal lifespan:
\begin{equation}\label{persofRegum}
\mathcal{I}^3_{\max}(\|u_0^m\|_{\mathcal{Z}^3_+}) = \mathcal{I}^2_{\max}(\|u_0^m\|_{\mathcal{Z}^2_+}).
\end{equation}
Indeed, since the nonlinear term in \eqref{CMdNLSfocu} involves at most one derivative of $u$, a standard energy estimate shows that the $\mathcal{Z}^3_+$-norm of $u^m$ satisfies a Gr\"onwall inequality whose coefficients depend only on $\|u^m (t) \|_{  \mathcal{Z}^2_+}$. Consequently, as long as the $\mathcal{Z}^2_+$-norm remains finite, the $\mathcal{Z}^3_+$-norm cannot blow up. \\

\noindent Then estimate \eqref{molliun0Zhinor} together with the monotonicity of the maximal lifespan implies that
\begin{equation}\label{TneTposisubsetI3m}
(T^-, T^+) = \mathcal{I}^2_{\max}(\|u_0\|_{\mathcal{Z}^2_+}) \subset \mathcal{I}^2_{\max}(\|u_0^m\|_{\mathcal{Z}^2_+}) = \mathcal{I}^3_{\max}(\|u_0^m\|_{\mathcal{Z}^3_+}).
\end{equation}
Thus, for  every $m \in \mathbb{N}_+$, the solution $u^m$ is well-defined and belongs to $\mathcal{Z}^3_+$ in $(T^-, T^+)$.  
  Fix $z \in \mathbb{C}_+$ and $t \in (T^-, T^+)$. Since $u^m$ is a $\mathcal{Z}^3_+$-solution to \eqref{CMdNLSfocu} with initial datum $u_0^m \in \mathcal{Z}^3_+$, the result of \textbf{Step \romannumeral 1} applies directly, leading to  the explicit formula
\begin{equation}\label{PzumtisbethZ3}
\mathbf{P}_z(u^m(t)) = \underline{u^m}(t,z) = \beth_{t, -t}^z(u_0^m),
\end{equation}
where the nonlinear functional $\beth_{t, -t}^z \in \mathcal{C}(\mathcal{Z}^1_+; \mathbb{C})$ is defined in   \eqref{defbethtrzufor}. 
 By formula \eqref{umzeroCvtou0} and Lemma~\ref{lemContRHSExpFor}, the right-hand side converges:
\begin{equation}\label{RHSCV}
\lim_{m \to +\infty} \beth_{t, -t}^z(u_0^m) = \beth_{t, -t}^z(u_0).
\end{equation}
  The flow map $\mathbf{F}_{\mathtt{T}}$ given by  \eqref{Flowmapfor} is continuous by  Proposition~\ref{LWPpropCMdNLSZn}, so  $
\lim_{m \to +\infty}   \|u^m(t) - u(t)\|_{\mathcal{Z}^2_+} = 0 $ by \eqref{umzeroCvtou0}.
  Recall that the Poisson integral operator $\mathbf{P}_z \in (L^{\infty}(\mathbb{R};\mathbb{C}))^*$ by \eqref{PzL2BdPz1is1}. Then we obtain
\begin{equation}\label{umtCVtoutPoiZ2}
\lim_{m \to +\infty} \mathbf{P}_z(u^m(t)) = \mathbf{P}_z(u(t)) = \underline{u}(t,z).
\end{equation}
  Combining \eqref{PzumtisbethZ3}, \eqref{RHSCV}, and \eqref{umtCVtoutPoiZ2}, we conclude that
$\underline{u}(t,z) = \beth_{t, -t}^z(u_0)$, 
which is exactly the explicit formula \eqref{ExpForCMdNLSfoc}. This completes the proof of Theorem~\ref{ExpForMainThm}.
\end{proof}

\noindent Theorem~\ref{ExpForMainThm} immediately yields a simpler
expression for solutions with a constant background, which is the
setting needed for dark multi-solitons.  

\begin{proof}[Proof of Corollary $\ref{CorExpforaaddH}$]
Since the initial datum satisfies that $u_0 -a \in L^2_+ \cap \mathcal{Z}^2_+$, we have $v_0 : =u_0 -a  \in H^2_+$ by \eqref{ZNcapL2=HN}. 
Lemma $\ref{lemuminainCtHnx}$ implies that $v:=u-a \in \bigl((T^-,T^+); H^2_+ \bigr) \cap \mathcal{C}^1 \bigl((T^-,T^+); L^2_+ \bigr)$.  
  Then we have
\begin{equation}\label{bfPzeitLapuav}
\mathbf{P}_z e^{i t \partial_x^2} (u_0) = a + \mathbf{P}_z e^{i t \partial_x^2} (v_0) \in \mathbb{C}, \qquad \mathscr{Q}_z e^{i t \partial_x^2} (u_0) = \mathscr{Q}_z e^{i t \partial_x^2} (v_0)  \in \mathrm{Dom}^2 (\mathbf{G}) ,  
\end{equation}for any
$z \in \mathbb{C}_+$ and $t \in ( T^-, T^+ )$,  by  formula \eqref{SchGreitLap1=1}  and  Lemma $\ref{LemKerscrQz}$. Then  formula \eqref{expForaaddH2}  is a direct consequence of
\eqref{ExpForCMdNLSfoc},   \eqref{bfPzeitLapuav} and  \eqref{bfIbbGtpertExpsum}.
\end{proof}

\begin{rem}\label{ComparwKLVrem}
  We compare the derivation of the explicit formulas \eqref{ExpForCMdNLSfoc}
  and \eqref{expForaaddH2} in the present paper with the one given by
  Killip--Laurens--Vi\c{s}an~\cite{KillLV2311}.
  The proof strategy in~\cite{KillLV2311} proceeds in two stages. First,
  Theorem~4.1 of~\cite{KillLV2311} establishes the explicit formula
  \begin{equation}\label{smooEFbyKLV}
    q(t,z) = \tfrac{1}{2\pi i} I_+\big( (X + 2t\mathcal{L}_{q^0} - z)^{-1} q^0 \bigr),
    \qquad \forall z \in \mathbb{C}_+,
  \end{equation}
  under the additional hypotheses that the initial datum $q(0)=q^0$ belongs to
  the smooth Hardy--Sobolev space $H^\infty_+$ and satisfies the weighted condition
  $\langle x\rangle q^0\in L^2(\mathbb{R};\mathbb{C})$.  Subsequently, a density
  and continuity argument removes the weighted hypothesis, extending \eqref{smooEFbyKLV}
  to all $L^2_+$-solutions (Theorem~1.7 of~\cite{KillLV2311}). \\

\noindent In the proof of Theorem~4.1 (see pp.~309--310 of~\cite{KillLV2311}),
the computation between formulas (4.46) and (4.47) applies the commutator
identity $[X,\mathcal{P}(q(t))] = 2\mathcal{L}_{q(t)}$ (Lemma~4.6
of~\cite{KillLV2311}; see also Remark~\ref{remKLV} of the present paper)
to the vector $U(t;0)\psi_0$, where $U(t;0)$ denotes the unitary propagator
corresponding to $\mathbf{W}(t) e^{-i t \mathbf{L}_{u_0}^2}$ in our notation (see  
\eqref{ExpForftsWteitLusq} and  \eqref{limformPDEdtf-tiButf} for verification), 
in order to obtain that
\begin{equation}\label{dtpsiKLV}
  \tfrac{d}{dt} \psi(t)
  = \mathcal{P}(q(t)) \psi(t)
    + \bigl[ X, \mathcal{P}(q(t))\bigr] U(t;0)\psi_0
    - 2 \mathcal{L}_{q(t)} U(t;0)\psi_0
  = \mathcal{P}(q(t)) \psi(t).
\end{equation}
However, Lemma~4.6 of~\cite{KillLV2311} requires the commutator $[X,\mathcal{P}(q(t))]$  to act on
functions $f\in L^2_+$ satisfying
$\langle x\rangle f\in H^2(\mathbb{R};\mathbb{C})$, i.e., functions in
$\mathrm{Dom}^2(\mathbf{X})$ by formula~\eqref{DomtwoXbyjbxf} of the present paper.
For \eqref{dtpsiKLV} to be justified, one must therefore verify the
stronger hypothesis
\begin{equation}\label{xUtpsiinH2klv}
  \langle x\rangle U(t;0)\psi_0\in H^2(\mathbb{R};\mathbb{C}), \qquad  \forall t \in (T^-, T^+) \setminus \{0\}.
\end{equation}
     Proposition~2.3 of~\cite{KillLV2311} guarantees, via Gr\"onwall's  inequality,
  that $U(t;0)\psi_0\in H^\infty_+$ and
  $\langle x\rangle U(t;0)\psi_0\in L^2(\mathbb{R};\mathbb{C})$,
  but does not verify the stronger hypothesis \eqref{xUtpsiinH2klv}.
  Consequently, the validity of \eqref{dtpsiKLV} relies on a regularity
  property of the propagated vector that is not established
  in~\cite{KillLV2311}. \\

\noindent In contrast, the proof of Theorem~\ref{ExpForMainThm} in the present paper
  does not require the weighted condition $u_0\in\mathrm{Dom}(\mathbf{X})$ to
  be imposed and later removed via density and continuity arguments.
  For any $a \in \mathbb{C}$, the invariance of $a + \bigl( H^3_+\cap\mathrm{Dom}(\mathbf{X}) \bigr)$ under the
  $\mathcal{Z}^3_+$-flow is proved in   
  Corollary $\ref{CorDomXinvsp}$   
  of the next subsection without recourse to Gr\"onwall's
  inequality.  The proof of Theorem~\ref{ExpForMainThm} employs only  the semigroup
  commutator estimates \eqref{commForSetasttildBuCv}, \eqref{SetasttildBuUniEstL2}
  together with the invariant subspace theorem (Theorem~\ref{invSpDomGHtwoThm}),
  operating entirely with the bounded operators $(\mathbf{S}(\eta)^*)_{\eta\ge 0}$
  and deferring the appearance   of the unbounded generator $\mathbf{G}$ to the
  final step, where it is extracted via the strong convergence characterization
  \eqref{SemiGrCVdefDomsG} of $\mathrm{Dom}(\mathbf{G})$.    
    The term $(\mathbf{G}-z)^{-1}(u(t)\chi_\epsilon)$ appearing
  in formula~\eqref{bfqeptzdef} lies automatically in
  $\mathrm{Dom}^3(\mathbf{G})\subset\mathrm{Dom}(\mathbf{G})\cap H^2_+$
  by~\eqref{DomsGDefSect2}.    
  Thus the present approach does not require
  any weighted hypotheses of the form
  $\langle x\rangle f\in L^2(\mathbb{R};\mathbb{C})$ or
  $\langle x\rangle f\in H^2(\mathbb{R};\mathbb{C})$, nor does it rely on
  Gr\"onwall estimates to propagate such conditions in time.
  As a result, the explicit formula~\eqref{ExpForCMdNLSfoc} is obtained
  directly for all $\mathcal{Z}^2_+$-solutions, and its specialization to the
  constant background case (Corollary~\ref{CorExpforaaddH}) holds under the
  sole hypothesis $u_0-a\in L^2_+$.
  \end{rem}

   \subsection{Invariant Subspaces and Conservation Laws}\label{InvSpCoLsubs4.4}
    
\noindent Throughout this subsection, we fix  
$T^{-}<0<T^{+}$, and let
$u\in\mathcal{C}\bigl((T^{-},T^{+});\mathcal{Z}^{3}_+\bigr)
   \cap \mathcal{C}^{1}\bigl((T^{-},T^{+});\mathcal{Z}^{1}_+\bigr)$
be a solution to equation~\eqref{CMdNLSfocu} with initial datum \(u(0)=u_{0}\in\mathcal{Z}^{3}_+\).
We assume that there exists a constant $a \in \mathbb{C}$ such that $\mathtt{f}_0 : = u_0 - a \in H^3_+$. 
 Lemma~\ref{lemuminainCtHnx} then guarantees that
\begin{equation}\label{ttfisu-a}
  \mathtt{f} : = u - a  \in \mathcal{C}\bigl((T^{-},T^{+});H^{3}_+\bigr)
          \cap \mathcal{C}^{1}\bigl((T^{-},T^{+});H^{1}_+\bigr).
\end{equation}   
  Let $\mathbf{W} \in \bigcap_{p=0}^{2} \mathcal{C}^{1}
\bigl((T^{-},T^{+}); \mathcal{B}_{\mathbb{C}}(H^{p}_+)\bigr)$ be the unitary propagator
constructed in~\eqref{UniPropagbfWSect4} of Subsection~\ref{subsUniEquiv4.2}.  
 The chiral Lax operator $\mathbf{L}_{u_0}: H^1_+ \to L^2_+$ is given by \eqref{LuSect2}. 
 The next lemma expresses   $ u (t) $ directly
in terms of the initial datum via the unitary propagator $\mathbf{W}(t)$ and the
Lax operator $\mathbf{L}_{u_0}$.

\begin{lem}
The solution $u \in a +  \mathcal{C}\bigl( (T^- ,T^+) ; H^3_+ \bigr)
          \cap \mathcal{C}^{1}\bigl( (T^- ,T^+) ; H^1_+\bigr)$ to \eqref{CMdNLSfocu}  satisfies that
 \begin{equation}\label{timevoutbfWu0}
 u(t) - a = \mathbf{W}(t) e^{-it\mathbf{L}_{u_0}^2}  (u (0) -a)  \in H^3_+  ,   
 \qquad   \forall t \in (T^-, T^+).
 \end{equation}
\end{lem}

\begin{proof}
Fix $ t  \in (T^- ,T^+) $. Let $\widetilde{\mathbf{B}}_{u(t)}\in \mathcal{B}_{\mathbb{C}}(H^2_+; L^2_+) $ be defined by  
\eqref{tildeBudefFoc-}.  
 Since $u$ solves  equation  \eqref{CMdNLSfocu} such that  $\mathtt{f}(t) = u(t) -a \in H^3_+$,  
formulas  
\eqref{uPidxumodsqToep} and \eqref{soluqudifCVSob}  imply that
\begin{equation*} 
\small
  \big\| \tfrac{1}{r}  \bigl(\mathtt{f}(t+r) - \mathtt{f}(t) \bigr) - \widetilde{\mathbf{B}}_{u(t)} \bigl( \mathtt{f}(t) \bigr)  \big\|_{H^1_+ } 
  =
  \big\| \tfrac{1}{r}  \bigl(u(t+r) - u(t) \bigr) - \bigl( i \partial_x^2 u + 2 u \Pi_{\geq 0} \partial_x  (|u|^2  ) \bigr) (t)  \big\|_{H^1_+ }  \to  0,
\end{equation*}as $r \to 0$. 
By  Lemma $\ref{LemWtexp-itL2phiPDE}$, formula \eqref{timevoutbfWu0} follows directly from   
  \eqref{ExpForftsWteitLusq}. 
\end{proof}

\noindent The next lemma shows that the generating functional $ \mathcal{H}_{\nu}^a  $ given by 
\eqref{calHGeneFctdefFor} is a conservation law of equation  \eqref{CMdNLSfocu}, for any $\nu \in \mathbb{C} \setminus \sigma \bigl(  \mathbf{L}_{u_0}  \bigr)$, 
under the same hypotheses as above: $u(0) = u_0  \in  a +  H^3_+$. 

 \begin{lem}\label{LemFirsIntGenFct}
 For every $\nu \in \mathbb{C} \setminus \sigma \bigl(  \mathbf{L}_{u_0}  \bigr)$ and $ t  \in (T^- ,T^+)$, we have $\mathcal{H}_{\nu}^a \bigl( u(t) \bigr)  = \mathcal{H}_{\nu}^a \bigl( u_0 \bigr)$. 
 \end{lem}
 
\begin{proof}
Fix $t \in (T^- ,T^+) $ and $\nu \in \mathbb{C} \setminus \sigma \bigl(  \mathbf{L}_{u_0}  \bigr) = \mathbb{C} \setminus \sigma \bigl(  \mathbf{L}_{u (t)}  \bigr)$. By 
\eqref{timevoutbfWu0},  \eqref{ResLutuniequiv} and \eqref{expisLusqcommutLu-a}, we have 
\begin{equation}\label{resLunuutabfWe}
\bigl( \mathbf{L}_{ u( t)} - 
 \nu \bigr)^{-1} (u ( t) - a) =   \mathbf{W}(t) e^{-it\mathbf{L}_{u_0}^2} \bigl( \mathbf{L}_{ u_0 } - 
 \nu \bigr)^{-1}  (u_0 - a).
\end{equation}
 Plugging \eqref{resLunuutabfWe} and \eqref{timevoutbfWu0} into formula \eqref{calHGeneFctdefFor}, we complete the proof.
\end{proof}

\noindent The conservation of $\mathcal{H}_{\nu}^a$ established above relies
solely on the unitary equivalence formula \eqref{UniEquLunFor} 
 between $\mathbf{L}_{u(t)}$ and $\mathbf{L}_{u_0}$  
  and does not require
any information about the domain of $\mathbf{G}$.  To obtain the conservation
of the functional $\mathbf{I}_+ (\cdot - a)$, which involves the boundary value
$\widehat{\mathtt{f}}(0^+)$ of the Fourier transform, we need the stronger hypothesis
that the initial perturbation $\mathtt{f}_0 : = u_0 - a $  belongs to $\mathrm{Dom}(\mathbf{G})$.
The next lemma shows that this property is propagated by the flow and that
$\mathbf{I}_+(u(t)-a)$ is independent of $t$.

\begin{lem}\label{invSpDomGCLbfI}
Assume that the initial datum $u(0) = u_0$ satisfies that $u_0 - a \in \mathrm{Dom} (\mathbf{G}) \cap H^3_+$, then 
\begin{equation}\label{bfIConserLaw}
u (t)  - a \in \mathrm{Dom} (\mathbf{G}) \cap H^3_+ 
\quad \text{and}  \quad
 \mathbf{I}_+ \bigl( u(t) - a  \bigr) = \mathbf{I}_+ \bigl( u_0 - a  \bigr) , \quad  \forall t \in (T^-, T^+).
\end{equation}
\end{lem}

\begin{proof}
 Fix $t \in (T^- ,T^+)$. 
 Since $ \mathtt{f}_0 = u_0 - a \in \mathrm{Dom} (\mathbf{G}) \cap H^3_+$,
 it follows that  
 \begin{equation}\label{ttftinDomGandH3}
u (t)  - a  = \mathtt{f} (t) = \mathbf{W}(t) e^{-it\mathbf{L}_{u_0}^2}  ( \mathtt{f}_0 )  \in H^3_+ \cap  \mathbf{W}(t)e^{-it \mathbf{L}_{u_0}^2}\left( \mathrm{Dom}(\mathbf{G}) \cap  H^2_+ \right) = \mathrm{Dom}(\mathbf{G}) \cap  H^3_+ .    
 \end{equation}
by  formulas 
 \eqref{ttfisu-a}, \eqref{timevoutbfWu0} and  
\eqref{DomGcapH2InvWU}.
 Since $\mathbf{W}(t) e^{-it\mathbf{L}_{u_0}^2} $ is a unitary operator on $L^2_+$,       
  we have 
\begin{equation}\label{EstibfGttfs}
\big\| \mathbf{G} \bigl( \mathtt{f} ( \mathfrak{s}) \bigr) \big\|_{L^2_+} 
= \big\|  \mathbf{W}( \mathfrak{s} ) e^{-i \mathfrak{s} \mathbf{L}_{u_0}^2} \bigl( \mathbf{G}   ( \mathtt{f}_0 ) + 2 \mathfrak{s} \mathbf{L}_{u_0} (  \mathtt{f}_0 ) \bigr) \big\|_{L^2_+} 
\leq \|\mathbf{G}   ( \mathtt{f}_0 ) \|_{L^2_+} + 2 |\mathfrak{s}| \|\mathbf{L}_{u_0} (  \mathtt{f}_0 )\|_{L^2_+}, 
\end{equation}by  
\eqref{GWtUtisWtUtGadd2tL}.
  The time interval $\mathtt{I}_t$ is given by \eqref{calItDef}.
  For any $\mathfrak{s} \in (T^- ,T^+)$, we define   
  \begin{equation}\label{boldupsideffor}
  \boldsymbol{\upsilon}^t( \mathfrak{s} ) : = e^{i ( t -  \mathfrak{s})  \partial_x^2} \bigl(  \boldsymbol{\varpi}_{u (\mathfrak{s})} \bigr) \in H^2_+
   \quad \text{with} \quad 
 \boldsymbol{\varpi}_{u (\mathfrak{s})} =  \Bigl( u \Pi_{\geq 0} \partial_x \bigl( |u|^2 \bigr) \Bigr) (\mathfrak{s}) \in H^2_+, 
  \end{equation}
  by formula     \eqref{varpiudeffor}. 
  Since $\boldsymbol{\upsilon}^t \in \mathcal{C} \bigl( (T^-, T^+) ; H^2_+ \bigr)$, it is Bochner integrable over the compact interval $\mathtt{I}_t$  by Theorem 3.27 of \cite{RudinFA}. 
  The Duhamel formula \eqref{Duhamelufor} yields that
\begin{equation}\label{DuhaForttfumisa}
\mathtt{f} (t) = e^{ i t \partial_x^2} \mathtt{f}_0 + 2 \mathrm{sign}(t)   \int_{\mathtt{I}_t} \boldsymbol{\upsilon}^t( \mathfrak{s} ) \mathrm{d} \mathfrak{s} \in H^3_+ \cap \mathrm{Dom} (\mathbf{G}) .
\end{equation}
Let $\widehat{\boldsymbol{\upsilon}^t} (\mathfrak{s} , \xi) $ denote the Fourier transform of $\boldsymbol{\upsilon}^t  (\mathfrak{s}) $ with respect to the spatial variable $x$. Then  
\begin{equation}\label{CVhatbdupsxito0}
\lim_{\xi \to 0} \big| \widehat{\boldsymbol{\upsilon}^t} (\mathfrak{s} , \xi)  \big| = \lim_{\xi \to 0} \big| e^{-i (t - \mathfrak{s}) \xi^2}\widehat{\boldsymbol{\varpi}_{u (\mathfrak{s})}} ( \xi)  \big|   = 0, \qquad \forall \mathfrak{s} \in (T^-, T^+),
\end{equation}by Lemma $\ref{LemNLtermCVEsti}$. Moreover,   formulas  \eqref{Esthatvrpiuxi} and \eqref{EstibfGttfs} imply  that 
\begin{equation}\label{hatuptsxiUniEsti}
\begin{split}
  \sup_{ \mathfrak{s} \in \mathtt{I}_t } \sup_{ 0 <  \xi \leq 1  } \big| \widehat{\boldsymbol{\upsilon}^t} (\mathfrak{s} , \xi)  \big|  
\leq  &  \sqrt{2 \pi} \widetilde{\mathfrak{m}} |a|^2 \bigl( \big\| \mathtt{f}_0  \big\|_{L^2_+} + \big\| \mathbf{G}(\mathtt{f}_0) \big\|_{L^2_+} + 2 |t| \big\| \mathbf{L}_{u_0} (\mathtt{f}_0) \big\|_{L^2_+}  \bigr) \\
& \qquad + 2 \bigl( |a| + \sup_{ \mathfrak{s} \in \mathtt{I}_t } \| u (\mathfrak{s}) \|_{L^{\infty} (\mathbb{R}; \mathbb{C})}\bigr) \big\| \mathtt{f}_0  \big\|_{L^2_+}   \sup_{ \mathfrak{s} \in \mathtt{I}_t } \big\|\partial_x \mathtt{f} (\mathfrak{s} ) \big\|_{L^2_+} .
\end{split}
\end{equation}
  Set $\widehat{\mathtt{f}} (t, \xi) := \widehat{\mathtt{f} (t)} (\xi)$. Using formulas    \eqref{DuhaForttfumisa},     \eqref{CVhatbdupsxito0} and  \eqref{hatuptsxiUniEsti},  we obtain that 
   \begin{equation}\label{CVdomttftxi}
   \big| \widehat{\mathtt{f}} (t, \xi) - e^{ - i t \xi^2} \widehat{\mathtt{f}}_0 (\xi)  \big| \leq 2 \int_{\mathtt{I}_t} \big| \widehat{\boldsymbol{\upsilon}^t} (\mathfrak{s} , \xi)  \big|   \mathrm{d} \mathfrak{s}  \to 0,  \quad \text{as} \quad \xi \to 0^+ ,
   \end{equation}by  the dominated convergence theorem.
      Consequently, the second formula of    \eqref{bfIConserLaw}
   follows directly from formulas    \eqref{ttfisu-a} and  \eqref{CVdomttftxi}.
\end{proof}

\noindent The next corollary is a direct consequence of Lemma $\ref{invSpDomGCLbfI}$  and formula \eqref{DomX}.

\begin{cor}\label{CorDomXinvsp}
If the initial datum $u(0) = u_0$ satisfies that $u_0 - a \in \mathrm{Dom} (\mathbf{X}) \cap H^3_+$, then 
\begin{equation}\label{DombfXinvarSp}
u (t)  - a \in \mathrm{Dom} (\mathbf{X}) \cap H^3_+ 
 ,  \qquad  \forall t \in (T^-, T^+).
\end{equation}
\end{cor}
  
 \begin{rem}\label{RemInvSpWeiKLV}
  When $a = 0$, the invariance of $\mathrm{Dom}(\mathbf{X}) \cap H^{\infty}_+$ under
  the $H^{\infty}_+$-flow of~\eqref{CMdNLSfocu} was also established in
  Proposition~2.3 on p. 295 of  Killip--Laurens--Vi\c{s}an~\cite{KillLV2311}.  The proof in~\cite{KillLV2311} relies
  on a sequence of bounded approximating weights combined with Gr\"onwall's
  inequality.  In contrast, Corollary~\ref{CorDomXinvsp} follows directly from
  the invariance of $a + \bigl(\mathrm{Dom}(\mathbf{G})\cap H^3_+\bigr)$ and the
  conservation of $\mathbf{I}_+ (\cdot - a)$, both established in
  Lemma~\ref{invSpDomGCLbfI}.  The invariance of
  $a + \bigl(\mathrm{Dom}(\mathbf{G})\cap H^3_+\bigr)$ is a direct consequence of
  Theorem~\ref{invSpDomGHtwoThm}, which is obtained via the semigroup
  commutator estimates \eqref{commForSetasttildBuCv} and
  \eqref{SetasttildBuUniEstL2} together with the strong convergence
  characterization \eqref{SemiGrCVdefDomsG} of $\mathrm{Dom}(\mathbf{G})$.
  This method does not rely on Gr\"onwall's estimate, nor does it require a
  cut-off and monotone convergence argument.  No differentiation of weighted
  norms is required.  Moreover, the result is stronger: it holds even for any
  nonzero constant background $a\in\mathbb{C}\setminus\{0\}$.
\end{rem}
 
\noindent Proposition $\ref{InvSpPropIntro}$ follows  directly from Lemma $\ref{invSpDomGCLbfI}$ and  Corollary $\ref{CorDomXinvsp}$.

\section{Dark Multi‑Soliton Potentials}\label{DarMulSoliPotSect}
 
\noindent This section is devoted to the study of dark multi-soliton potentials
of the focusing CMdNLS equation~\eqref{CMdNLSfocu}. 
 Let $N \geq 1$ be an arbitrary positive integer. 
 Recall from
Definition~\ref{DefiDarNsoli} that a function $u$ belongs to
$\mathcal{U}_N$ if and only if it admits a rational representation
$u=P/Q$, where $Q\in\mathbb{C}[X]$ is monic of degree $N$ with
$Q^{-1}\{0\}\subset\mathbb{C}_-$,  
and the following algebraic identity
holds (see formula~\eqref{PQiden+}):
 \begin{equation}\label{PQidenSec5}
  |P(x)|^{2}-|Q(x)|^{2}
  =
  i\bigl(
    \overline{Q(x)}\partial_{x}Q(x)
    -
    Q(x)\partial_{x}\overline{Q(x)}
  \bigr),
  \qquad
  \forall x\in\mathbb{R}.
\end{equation}

\begin{rem}\label{Rkuinftyv}
 Assume that $Q\in\mathbb{C}[X]$ is  
  a monic polynomial of degree $N$ such that
  $Q^{-1}\{0\}\subset\mathbb{C}_-$.  Then
  $\overline{Q}\partial_{x}Q-Q\partial_{x}\overline{Q}\in\mathbb{C}[X]$
  is a polynomial of degree $2N-2$.  If
  $u=P/Q\in\mathcal{U}_N$ is a dark $N$-soliton
  of equation~\eqref{CMdNLSfocu}, then \eqref{PQidenSec5} yields the existence of
  $\alpha\in\mathbb{T} = \mathbb{R}/2\pi\mathbb{Z}$ and
  $a_{0},a_{1},\cdots,a_{N-1}\in\mathbb{C}$ such that
  \begin{equation}\label{ueiapluvDarSoli}
    P(x)=e^{i\alpha}x^{N}+\sum_{j=0}^{N-1}a_{j}x^{j}
    \qquad\text{and}\qquad
    u(x)=e^{i\alpha}+v(x),
    \qquad
    \forall x\in\mathbb{R},
  \end{equation}
  where
  $v\in\frac{\mathbb{C}_{N-1}[X]}{Q}
     \subset\mathrm{Dom}^{\infty}(\mathbf{G})
     \subset H^{\infty}_{+}\cap\mathrm{Dom}(\mathbf{G})$
  by formula~\eqref{KthetainDominftyG}.
  We define
  \begin{equation}\label{uinftyDarSoli}
    u(\infty)
    :=
    \lim_{|x|\to+\infty}u(x)
    =
    e^{i\alpha}
    \in
    \mathbb{S}^{1}.
  \end{equation}
  Every dark $N$-soliton $u\in\mathcal{U}_N$ admits nonzero
  limits as $x\to-\infty$ and $x\to+\infty$, and these limits coincide.
  Consequently, every dark multi-soliton belongs to the smooth
  Hardy--Zhidkov space $\mathcal{Z}^{\infty}_{+}$, that is,
  \begin{equation}\label{NdarSoliinZinfty}
    \bigcup_{N\in\mathbb{N}_{+}}\mathcal{U}_{N}
    \subset
    \bigcup_{\alpha\in\mathbb{T}}
      \bigl(
        e^{i\alpha}+\mathrm{Dom}^{\infty}(\mathbf{G})
      \bigr)
    \subset
    \mathcal{Z}^{\infty}_{+}.
  \end{equation}
\end{rem}
 
\noindent The purpose of this section is twofold: first, to prove the spectral characterization theorem (Theorem~\ref{MainThmSpectChar}); second, to construct the action-angle mapping $\boldsymbol{\Phi}_N\colon\mathcal U_N\to
\mathbb{T}\times \Omega_N\times\mathbb{R}^{1\times N}$  in~\eqref{PhiNintro} and establish its bijectivity.  All objects in this section are time-independent; the proof of the time evolution identity \eqref{TimEvoactangl} is deferred to Section~\ref{SectTimeEvo6}.
  For notational simplicity, in this section and the next, we denote the
Szeg\H{o} projector $ \Pi_{\geq 0} \in \mathcal{B}_{\mathbb{C}} \bigl( L^2 (\mathbb{R} ; \mathbb{C}) \bigr) $ given by  \eqref{SzegoProjS2}
simply by $\Pi$.

\subsection{Spectral Characterization}\label{SpecCharSubsect}

\noindent For any $N \in \mathbb{N}_+$ and $u \in \mathcal{U}_N$, Lemma $\ref{lemLuselfadj}$ shows  that  the Lax operator $\mathbf{L}_u$ is a self-adjoint operator on $L^2_+$ with domain $H^1_+$. Then $\mathbf{L}_u$   admits the following spectral decomposition
\begin{equation}\label{SpDecbfLu}
L^2_+ = \mathscr{H}_{\mathrm{ac}}(\mathbf{L}_u)  \bigoplus  \mathscr{H}_{\mathrm{sc}}(\mathbf{L}_u)  \bigoplus \mathscr{H}_{\mathrm{pp}}(\mathbf{L}_u) 
\quad \text{and}   \quad 
  \sigma(\mathbf{L}_u) =
\sigma_{\mathrm{ac}}(\mathbf{L}_u) \cup 
\sigma_{\mathrm{sc}}(\mathbf{L}_u) \cup \overline{\sigma_{\mathrm{pp}}(\mathbf{L}_u)} .
\end{equation} 
  Before proving Theorem $\ref{MainThmSpectChar}$,  
  we first give precise information on the spectral subspaces of
$\mathbf{L}_u$.

\begin{prop}\label{propSpecDecomLu}
For any  integer $N \geq 1$ and $u   \in \mathcal{U}_N$, 
Definition $\ref{DefiDarNsoli}$ and Remark $\ref{Rkuinftyv}$ imply that there exist two polynomials $P, Q \in \mathbb{C}[X]$ of degree $N$ such that $u=P/Q$, $Q$ is monic, $Q^{-1} \{ 0 \} \subset \mathbb{C}_-$ and formula \eqref{PQidenSec5} holds. 
   Set $\theta_{Q} := \frac{\overline{Q}}{Q} \in 1 + \mathrm{Dom}^{\infty}(\mathbf{G})$, where $\mathrm{Dom}^{\infty}(\mathbf{G})$ 
   is given by \eqref{DomrGsubsetDomsG}. 
     Then  
\begin{equation}\label{SpecDeHacHscHpp}
\mathscr{H}_{\mathrm{ac}}(\mathbf{L}_u) = \theta_Q L^2_+, \quad \mathscr{H}_{\mathrm{sc}}(\mathbf{L}_u) = \{0\} \quad \mathrm{and} \quad \mathscr{H}_{\mathrm{pp}}(\mathbf{L}_u) = \tfrac{\mathbb{C}_{\leq N-1}[X]}{Q}\subset  \mathrm{Dom}^{\infty} (\mathbf{G})  
.
\end{equation}Furthermore, the Lax operator $\mathbf{L}_u$ has exactly $N$ distinct eigenvalues and we have
\begin{equation}\label{spectLUacsc}
\sigma_{\mathrm{ac}}(\mathbf{L}_u) = [-1, +\infty), \quad \sigma_{\mathrm{sc}}(\mathbf{L}_u) = \emptyset  \quad \mathrm{and} \quad -1 \notin \sigma_{\mathrm{pp}}(\mathbf{L}_u).
\end{equation}In particular, each eigenvalue of $\mathbf{L}_u$ is simple.
\end{prop}
 
\begin{rem}\label{Remneg1notEiVa}
  For a dark multi-soliton potential $u$, Proposition~\ref{propSpecDecomLu}
  already excludes $-1$ from the pure point spectrum of $\mathbf{L}_u$.  The
  stronger exclusion
  $\sigma_{\mathrm{pp}}(\mathbf{L}_u)\cap[-1,0]=\emptyset$
  will be proved in Lemma~\ref{LemSpLutwoseg}.
\end{rem}

\begin{proof}[Proof of Proposition $\ref{propSpecDecomLu}$]
Consider a dark $N$-soliton $u=\frac{P}{Q} \in \mathcal{U}_N$ such that $Q$ is a monic polynomial with $\deg Q =N$ and $Q^{-1}\{0\} \subset 
\mathbb{C}_-$. We use the shorthand notation:
\begin{equation}\label{thetatheQabbre}
\theta :=\theta_Q=\tfrac{\overline{Q}}{Q} \in 1 + \mathrm{Dom}^{\infty}(\mathbf{G}) \subset \mathcal{Z}^{\infty}_+.
\end{equation}Formula   \eqref{ueiapluvDarSoli} yields that $\deg P = \deg \overline{P}=N$. Then the pair $(\overline{P}, Q)\in \mathbb{C}[X] \times \mathbb{C}[X]$ also satisfies \eqref{PQidenSec5}. So  we have  $\frac{\overline{P}}{Q} \in \mathcal{U}_N \subset \mathcal{Z}^{\infty}_+$ by  \eqref{NdarSoliinZinfty}.  
For any $g \in L^2_+$, formula \eqref{ZnHardBanAlgebra} implies that 
\begin{equation}\label{TuTubarthetaf}
\mathbf{T}_u \mathbf{T}_{\overline{u}}(\theta g) = u \Pi \left( \tfrac{\overline{P}}{\overline{Q}} \tfrac{\overline{Q}}{Q} g \right) = u \tfrac{\overline{P}}{Q} g = |u|^2 \theta g \in L^2_+.
\end{equation}In terms of the inner function $\theta$ defined in \eqref{thetatheQabbre},
  formula \eqref{PQidenSec5} takes the following equivalent form: 
\begin{equation}\label{DTheta|u|2-1}
\theta( |u|^2 -1) =D \theta  \in  \mathrm{Dom}^{\infty}(\mathbf{G}).
\end{equation}Since $|\theta(x)| = 1$ for all $x \in \mathbb{R}$, the multiplication operator $\mathbf{m}_{\theta}: f \mapsto \theta f$ is an isometry on the Hardy space $L^2_+$. In particular, its image, $\theta L^2_+$, is a closed subspace. Moreover, formula  \eqref{HnZnsubsetHn} yields that 
\begin{equation}\label{thetaH1}
\theta L^2_+ \cap  H^1_+ = \theta H^1_+.  
\end{equation}For any $h \in H^1_+$, we use  \eqref{LuSect2}, \eqref{TuTubarthetaf},  \eqref{DTheta|u|2-1} and \eqref{thetaH1} to obtain that $\theta h \in H^1_+$ and
\begin{equation}\label{Luthetah}
\mathbf{L}_u (\theta h) = D(\theta h) - \mathbf{T}_u \mathbf{T}_{\overline{u}}(\theta h) = \theta (D-1)h \in \theta L^2_+.
\end{equation}Consider  the restriction of the Lax operator $\mathbf{L}_u$ to the subspace $\theta L^2_+$, denoted by  $\mathbf{L}_u|_{\theta L^2_+}$, and equip it with the domain  $\mathrm{Dom}(\mathbf{L}_u|_{\theta  L^2_+}) = \theta H^1_+ = \theta L^2_+ \cap  H^1_+$ by \eqref{thetaH1}. We claim that $(\mathbf{L}_u|_{\theta  L^2_+}, \mathrm{Dom}(\mathbf{L}_u|_{\theta  L^2_+}))$ is an unbounded self-adjoint operator on the Hilbert space $(\theta L^2_+, \langle \cdot  , \cdot \rangle_{L^2_+})$, i.e. 
\begin{equation}\label{LuthetaL2selfadj}
 \mathbf{L}_u|_{\theta  L^2_+} = (\mathbf{L}_u|_{\theta  L^2_+})^* : \mathrm{Dom}( \mathbf{L}_u|_{\theta  L^2_+}) = \theta H^1_+ \to \theta L^2_+.
\end{equation}
\begin{proof}[Proof of  formula \eqref{LuthetaL2selfadj}]
Formula \eqref{Luthetah} and Lemma $\ref{lemLuselfadj}$ imply that $(\mathbf{L}_u|_{\theta  L^2_+}, \mathrm{Dom}(\mathbf{L}_u|_{\theta  L^2_+}))$ is an unbounded symmetric operator on $(\theta L^2_+, \langle \cdot  , \cdot \rangle_{L^2_+})$. Assume that $g = \theta \mathfrak{w} \in \mathrm{Dom} ((\mathbf{L}_u|_{\theta  L^2_+})^*  )$ for some $\mathfrak{w} \in L^2_+$, there exists $\mathfrak{y} \in L^2_+$ such that for any $h \in H^1_+$, we have
\begin{equation*}
\langle (D-1) h,  \mathfrak{w}\rangle_{L^2_+}=\langle \mathbf{L}_u(\theta h), g \rangle_{L^2_+} = \langle \theta h, \theta \mathfrak{y}\rangle_{L^2_+} = \langle h,   \mathfrak{y}\rangle_{L^2_+}.
\end{equation*}Lemma $\ref{DNSelfadjinHsLem}$ yields that $\mathfrak{w} \in \mathrm{Dom}\left( D^*\right) = H^1_+$. Then $g = \theta \mathfrak{w} \in \theta H^1_+$ and  \eqref{LuthetaL2selfadj} is obtained.
\end{proof}
\noindent According to formulas \eqref{thetaH1} and \eqref{Luthetah}, the  self-adjoint operator $ D-1:  \mathrm{Dom}(D-1)= H^1_+ \to L^2_+$ is unitarily equivalent to  $(\mathbf{L}_u|_{\theta  L^2_+}, \mathrm{Dom}(\mathbf{L}_u|_{\theta  L^2_+}))$  via the isometry $\mathbf{m}_{\theta}: L^2_+ \to \theta L^2_+$. In other words, 
\begin{equation}\label{UniEquD-1toLurest}
\mathrm{Dom} (\mathbf{L}_u|_{\theta  L^2_+}   ) = \theta H^1_+ = \mathbf{m}_{\theta}\left(\mathrm{Dom}(D-1)  \right)   \quad \mathrm{and} \quad  \mathbf{L}_u|_{\theta  L^2_+} \circ \mathbf{m}_{\theta} = \mathbf{m}_{\theta} \circ (D-1) : H^1_+ \to \theta L^2_+.
\end{equation}For any $g \in L^2_+$, let $\mu_{g}^{\mathbf{L}_u}$ denote the spectral measure of the Lax operator $(\mathbf{L}_u, \mathrm{Dom}(\mathbf{L}_u))$ associated to  $g$. 
Proposition $2.1$ of \cite{SUNCMP2021} yields that for every bounded Borel function $f : \mathbb{R} \to \mathbb{C}$, we have  
\begin{equation}\label{intfdmuSpMea}
\int_{\mathbb{R}}f \mathrm{d} \mu_{\theta g}^{\mathbf{L}_u} = \langle f(\mathbf{L}_u) \theta g , \theta g \rangle_{L^2_+} = \langle f( D -1 )  g ,   g \rangle_{L^2_+} =\int_{-1}^{+\infty} f( \xi ) |\widehat{g}(\xi+1)|^2 \tfrac{\mathrm{d\xi}}{2 \pi}
\end{equation}Then $ \mathrm{d} \mu_{\theta g}^{\mathbf{L}_u}(\xi) = \frac{\mathbf{1}_{\xi >-1}}{2 \pi}  |\widehat{g}(\xi+1)|^2\mathrm{d}\xi$. Since $g \in L^2_+$, the spectral measure $\mu_{\theta g}^{\mathbf{L}_u} $ is absolutely continuous with respect to the Lebesgue measure on $\mathbb{R}$. So $\theta g \in \mathscr{H}_{\mathrm{ac}}(\mathbf{L}_u)$. As a consequence, we have
\begin{equation}\label{thetaL2inHac}
\theta L^2_+ \subset
 \mathscr{H}_{\mathrm{ac}}(\mathbf{L}_u) \Longleftrightarrow  \mathscr{H}_{\mathrm{ac}}(\mathbf{L}_u)^{\perp} \subset K_{\theta}:= \left(\theta L^2_+ \right)^{\perp} = \tfrac{\mathbb{C}_{\leq N-1}[X]}{Q}
\end{equation}by formula \eqref{OrthDecomInnKthe}. 
The subspace $\theta H^1_+$ is dense in   $(\theta L^2_+, \langle \cdot, \cdot \rangle_{L^2_+})$ and  \eqref{Luthetah} implies that $\mathbf{L}_u(\theta H^1_+) \subset \theta L^2_+$. Since  $\mathbf{L}_u$ is self-adjoint on the Hilbert space $(L^2_+, \langle \cdot, \cdot \rangle_{L^2_+})$ by Lemma $\ref{lemLuselfadj}$,  we have
\begin{equation}
\mathbf{L}_u\left( K_{\theta} \right) \subset K_{\theta} \quad \mathrm{with} \quad \dim_{\mathbb{C}} K_{\theta} =\dim_{\mathbb{C}}\left( \tfrac{\mathbb{C}_{\leq N-1}[X]}{Q} \right) = N.
\end{equation}By the spectral theorem for  $\mathbf{L}_u |_{K_{\theta}}: K_{\theta} \to K_{\theta}$, there exists an orthonormal basis of the Hermitian space $(K_{\theta}, \langle \cdot, \cdot \rangle_{L^2_+})$ consisting of eigenvectors of $\mathbf{L}_u|_{K_{\theta}}$. Then we have
\begin{equation}\label{KtheinHpp}
\tfrac{\mathbb{C}_{\leq N-1}[X]}{Q}=(\theta L^2_+)^{\perp} =K_{\theta} \subset \mathscr{H}_{\mathrm{pp}}(\mathbf{L}_u) \subset  \left(\mathscr{H}_{\mathrm{ac}}(\mathbf{L}_u) \right)^{\perp} \quad \mathrm{with} \quad \theta = \theta_Q = \tfrac{\overline{Q}}{Q}.
\end{equation}Therefore, formula \eqref{SpecDeHacHscHpp}   follows from formulas \eqref{thetaL2inHac} and \eqref{KtheinHpp}.  We recall that
\begin{equation}\label{HppinDomI}
v := u - u(\infty) \in \tfrac{\mathbb{C}_{\leq N-1}[X]}{Q} = \mathscr{H}_{\mathrm{pp}}(\mathbf{L}_u)   \subset \mathrm{Dom}^{\infty}(\mathbf{G}) \subset H^{\infty}_+ \cap  \mathrm{Dom}(\mathbf{I}_+),
\end{equation}according to formulas \eqref{DXsubsetDG}, \eqref{ueiapluvDarSoli}, 
\eqref{uinftyDarSoli} and 
\eqref{KthetainDominftyG}.
  Then  $\mathbf{L}_u$ has at most $N$ eigenvalues.  For any $\lambda \in \sigma_{\mathrm{pp}}(\mathbf{L}_u)$ and $\mathbf{h} \in \mathrm{Ker}(\lambda - \mathbf{L}_u) \setminus \{0\}\subset \mathscr{H}_{\mathrm{pp}}(\mathbf{L}_u)\setminus \{0\} \subset H^{\infty}_+ \cap  \mathrm{Dom}(\mathbf{I}_+) $, we have  
\begin{equation*}
  e^{-i \alpha }\widehat{\mathbf{h}}(0^+)   \langle v , \mathbf{h} \rangle_{L^2_+} + |\langle v , \mathbf{h} \rangle_{L^2_+}|^2 = 2 \pi \|\mathbf{h} \|_{L^2_+}^2  > 0.
\end{equation*}by using \eqref{idenhL2nvhinprod} and \eqref{HppinDomI}. So $\langle v , \mathbf{h} \rangle_{L^2_+} \ne 0$ and Lemma $\ref{EigLuhperpvLem}$ yields that $\lambda \ne -1$. By Corollary \ref{CorSimpofEigLu}, the fact that $-1 \notin \sigma_{\mathrm{pp}}(\mathbf{L}_u)$ implies the simplicity of  every eigenvalue of the Lax operator $\mathbf{L}_u$.
Since we have established $\mathscr{H}_{\mathrm{ac}}(\mathbf{L}_u) = \theta L^2_+$  in \eqref{SpecDeHacHscHpp}, the unitary equivalence \eqref{UniEquD-1toLurest} allows us to determine the absolutely continuous spectrum. Applying Proposition 2.1 and Proposition 2.3 of \cite{SUNCMP2021}, we obtain
\begin{equation}
\sigma_{\mathrm{ac}}(\mathbf{L}_u)= \sigma (\mathbf{L}_u|_{\mathscr{H}_{\mathrm{ac}}(\mathbf{L}_u)}) = \sigma (\mathbf{L}_u|_{\theta L^2_+}) = \sigma (D-1) = \sigma (D)-1 =[-1 , +\infty).
\end{equation}So we obtain formula \eqref{spectLUacsc} and the proof of Proposition $\ref{propSpecDecomLu}$ is complete.
\end{proof}

 \noindent The spectral decomposition~\eqref{SpDecbfLu}, combined with the
explicit description of the spectral subspaces~\eqref{SpecDeHacHscHpp}
obtained in Proposition~\ref{propSpecDecomLu}, directly yields the
implication \((\mathbf{a})\Rightarrow(\mathbf{c})\) of
Theorem~\ref{MainThmSpectChar}.  With this observation in hand, the
remaining part of the proof of Theorem~\ref{MainThmSpectChar} can now
be completed as follows.

\begin{proof}[Proof of Theorem $ \ref{MainThmSpectChar}$]
The proof proceeds cyclically: $(\mathbf{a}) \Rightarrow (\mathbf{c})$, $(\mathbf{b}) \Rightarrow (\mathbf{a})$, and $(\mathbf{c}) \Rightarrow (\mathbf{b})$.\\
 
\noindent \textbf{Step \romannumeral 1}. $(\mathbf{a}) \Rightarrow (\mathbf{c})$: Assume that $u = e^{i \alpha} +v = \frac{P}{Q} \in \mathcal{U}_N$, where $Q \in \mathbb{C}[X]$ is a monic polynomial of degree $N \geq 1$ such that $Q^{-1}\{0\} \subset \mathbb{C}_-$. By Proposition $\ref{propSpecDecomLu}$ and Remark $\ref{Rkuinftyv}$, it follows that
\begin{equation*}
v = u - u(\infty) \in \tfrac{\mathbb{C}_{\leq N-1}[X]}{Q} = \mathscr{H}_{\mathrm{pp}}(\mathbf{L}_u) \subset \mathrm{Dom}  (\mathbf{G}) \quad \mathrm{and} \quad -1 \notin \sigma_{\mathrm{pp}}(\mathbf{L}_u).
\end{equation*}by \eqref{DomrGsubsetDomsG}. Then $\dim_{\mathbb{C}} \mathscr{H}_{\mathrm{pp}}(\mathbf{L}_u) =N$ and we have $\mathbf{G}\left(\mathscr{H}_{\mathrm{pp}}(\mathbf{L}_u)\right)= \mathscr{H}_{\mathrm{pp}}(\mathbf{L}_u) = \tfrac{\mathbb{C}_{\leq N-1}[X]}{Q}$
by formula  \eqref{KthetainDominftyG}.
  Therefore, both $(\mathbf{b})$ and $(\mathbf{c})$ follow from the assertion $(\mathbf{a})$.\\

\noindent \textbf{Step \romannumeral 2}. $(\mathbf{b}) \Rightarrow (\mathbf{a})$: Assume that $u = e^{i \alpha} +v \in \mathcal{Z}^1_+$ for some $\alpha \in \mathbb{T}$ and $v \in H^1_+$ such that 
\begin{equation}\label{vinHppLuinDomG12}
v =u - e^{i \alpha} \in \mathscr{H}_{\mathrm{pp}}(\mathbf{L}_u), \quad \mathscr{H}_{\mathrm{pp}}(\mathbf{L}_u) \subset \mathrm{Dom}(\mathbf{G}), \quad \mathbf{G}\left(\mathscr{H}_{\mathrm{pp}}(\mathbf{L}_u) \right) \subset \mathscr{H}_{\mathrm{pp}}(\mathbf{L}_u),  
\end{equation}and $\dim_{\mathbb{C}} \mathscr{H}_{\mathrm{pp}}(\mathbf{L}_u) =N \in \mathbb{N}_+$.  Lemma 3.3 of \cite{SUNCMP2021}, which is a special case of the Lax--Beurling theorem, characterizes all finite-dimensional translation-invariant subspaces of $L^2_+$ (See also Remark $\ref{RemLaxTrInvSp}$). It implies that there exists a unique monic polynomial $Q \in \mathbb{C}[X]$ of degree $N \in \mathbb{N}_+$ such that
 \begin{equation}\label{QischarGHppbfLu}
Q^{-1}\{0\} \subset \mathbb{C}_-, \quad \mathscr{H}_{\mathrm{pp}}(\mathbf{L}_u) = \tfrac{\mathbb{C}_{\leq N-1}[X]}{Q} \quad \mathrm{and} \quad Q(x)= \det(x-\mathbf{G}|_{\mathscr{H}_{\mathrm{pp}}(\mathbf{L}_u)}), 
\end{equation}for any $ x \in \mathbb{R}$. Set  $\theta:=\theta_Q =\frac{\overline{Q}}{Q} \in 1+\tfrac{\mathbb{C}_{\leq N-1}[X]}{Q} \subset \mathcal{Z}^{\infty}_+$. Formulas \eqref{HnZnsubsetHn},  \eqref{OrthDecomInnKthe} and \eqref{QischarGHppbfLu} imply that
\begin{equation}\label{LuKtheinKthe2-1}
\theta H^1_+ = \theta L^2_+ \cap  H^1_+, \quad \theta L^2_+ = \left(\tfrac{\mathbb{C}_{\leq N-1}[X]}{Q} \right)^{\perp}  \quad \mathrm{and} \quad \mathbf{L}_u \left(\tfrac{\mathbb{C}_{\leq N-1}[X]}{Q} \right) \subset  \tfrac{\mathbb{C}_{\leq N-1}[X]}{Q}.
\end{equation}Recall that $\mathbf{L}_u = \mathbf{L}_u^* : \mathrm{Dom}(\mathbf{L}_u)=H^1_+ \to L^2_+$ by Lemma $\ref{lemLuselfadj}$. Then formula \eqref{LuKtheinKthe2-1} yields that
\begin{equation}\label{LutheH1inthetaL22-1}
 \mathbf{L}_u (\theta H^1_+) \perp  \tfrac{\mathbb{C}_{\leq N-1}[X]}{Q} \Longrightarrow \mathbf{L}_u (\theta H^1_+) \subset    \theta L^2_+.
\end{equation}According to  \eqref{vinHppLuinDomG12} and \eqref{QischarGHppbfLu},  there exists a polynomial $q \in \mathbb{C}[X]$ such that $\deg q \leq N-1$ and
\begin{equation}\label{qovQ=vinHppLu=Kthet}
  u-e^{i \alpha} = v = \tfrac{q}{Q} \in \mathscr{H}_{\mathrm{pp}}(\mathbf{L}_u) = \tfrac{\mathbb{C}_{\leq N -1}[X]}{Q} \subset \mathrm{Dom}^{\infty}(\mathbf{G}) \subset H^{\infty}_+ \cap  \mathrm{Dom}(\mathbf{G}).
\end{equation}Set $P:=uQ=e^{i\alpha} Q + q \in \mathbb{C}[X]$. Then  $\deg P = N$ and $\tfrac{\overline{P}}{Q} = e^{-i \alpha} \theta + \tfrac{\overline{q}}{Q} \in e^{-i \alpha} + \tfrac{\mathbb{C}_{\leq N -1}[X]}{Q} \subset  \mathcal{Z}^{\infty}_+$. Since $u=\frac{P}{Q} \in e^{i\alpha} +  H^{\infty}_+ \subset \mathcal{Z}^{\infty}_+$, for any $g \in L^2_+$, it follows from  formula \eqref{ZnHardBanAlgebra} that 
\begin{equation}\label{TuTu2barthetag2}
\overline{u}\theta g = \tfrac{\overline{P}}{Q}g \in L^2_+ \Longrightarrow \mathbf{T}_u \mathbf{T}_{\overline{u}}(\theta g) = u \Pi(\overline{u}\theta g) = |u|^2 \theta g \in L^2_+.
\end{equation}According to formulas
\eqref{LuSect2},
\eqref{LutheH1inthetaL22-1} and 
\eqref{TuTu2barthetag2}, for any $h \in H^1_+$, we have
\begin{equation}\label{Dthe-the(u|2-1)hinL2+}
\theta h \left(\tfrac{D \theta}{\theta} - |u|^2+1 \right)=\mathbf{L}_u (\theta h) - \theta (D-1)h \in \theta L^2_+ \Longrightarrow h \left(\tfrac{D \theta}{\theta} - |u|^2+1 \right) \in L^2_+.
\end{equation}Proposition 4.3 in \cite{SUNCMP2021} yields that $\frac{D \theta}{\theta}$ is an $N$-soliton potential of the Benjamin--Ono equation, i.e.
\begin{equation}\label{BONsolitontheta}
Q(x)=\prod_{j=1}^N(x-z_j) \quad \mathrm{with} \quad z_j \in \mathbb{C}_- \quad \Longrightarrow \quad  \tfrac{D \theta}{\theta} (x) =  \sum_{j=1}^N  \frac{-2 \operatorname{Im} z_j  }{|x-z_j|^2}, \quad \forall x \in \mathbb{R}.
\end{equation}We thus have $\frac{D \theta}{\theta} = -2 \operatorname{Im}\frac{\partial_x Q}{Q}\in H^{\infty}\cap  L^1 (\mathbb{R}; \mathbb{R})$. Since $v \in H^{\infty}_+ \cap  \mathrm{Dom}(\mathbf{G})$ by \eqref{qovQ=vinHppLu=Kthet},    it follows from \eqref{DomsGDefSect2} that $\widehat{v}|_{\mathbb{R}_+^*} \in H^1(\mathbb{R}_+^*)$. So $\widehat{v}$  is continuous at every point $\xi \in \mathbb{R}^* =  \mathbb{R}\setminus \{0\}$. We define
\begin{equation}\label{fDTheThe-|u|2+1}
\mathtt{f}:=\tfrac{D \theta}{\theta} - |u|^2+1 = -2 \operatorname{Im}\tfrac{\partial_x Q}{Q} - e^{-i \alpha} v -e^{i\alpha} \overline{v} - |v|^2 \in H^{\infty}(\mathbb{R}; \mathbb{R}) \quad \mathrm{with} \quad \widehat{v}\in \mathcal{C}(\mathbb{R}^*; \mathbb{C}).
\end{equation}So $\widehat{\mathtt{f}}\in \mathcal{C}(\mathbb{R}^*; \mathbb{C})$.
Using  Lemma $\ref{L2+capL2r=0Lem}$,  Lemma $\ref{lemfhinHarSobfL2+}$ and formula \eqref{Dthe-the(u|2-1)hinL2+}, we obtain that
\begin{equation*}
\tfrac{D \theta}{\theta} - |u|^2+1 = \mathtt{f} \in H^{\infty}(\mathbb{R}; \mathbb{R})\cap  L^2_+ = \{0\}.
\end{equation*}In other words, the  pair $(P, Q) \in \mathbb{C}[X] \times \mathbb{C}[X]$ satisfies  \eqref{PQidenSec5}. Consequently, the function $u = e^{i\alpha} + v = \frac{P}{Q}$ is a dark $N$-soliton by Definition $ \ref{DefiDarNsoli}$.\\

\noindent \textbf{Step \romannumeral 3}. $(\mathbf{c}) \Rightarrow (\mathbf{b})$: Assume that $N \geq 1$ and $u = e^{i\alpha} + v \in \mathcal{Z}^1_+$ for some $\alpha \in \mathbb{T}$ and $v \in H^1_+$ such that
\begin{equation}\label{vinDomGandvGvinHpp}
\small 
-1 \notin \sigma_{\mathrm{pp}}(\mathbf{L}_u), \quad \dim_{\mathbb{C}}\mathscr{H}_{\mathrm{pp}}(\mathbf{L}_u) = N, \quad  v = u - e^{i \alpha} \in \mathscr{H}_{\mathrm{pp}}(\mathbf{L}_u)\cap  \mathrm{Dom}(\mathbf{G})  \quad  \mathrm{and} \quad \mathbf{G} v \in \mathscr{H}_{\mathrm{pp}}(\mathbf{L}_u).
\end{equation}By Corollary \ref{CorSimpofEigLu}, the self-adjoint operator $\mathbf{L}_u: H^1_+ \to L^2_+$ possesses exactly $N$ distinct real eigenvalues. Denoting them in increasing order as $\lambda_1 < \lambda_2 < \cdots < \lambda_N$, the corresponding eigenspaces are one-dimensional. Consequently, there exists an orthonormal basis $\{\mathfrak{h}_1, \mathfrak{h}_2, \cdots, \mathfrak{h}_N\}$ for $\mathscr{H}_{\mathrm{pp}}(\mathbf{L}_u)$ such that  
\begin{equation}\label{hjinDomI+,H1+}
\mathbb{C}\mathfrak{h}_j =\mathrm{Ker}(\lambda_j - \mathbf{L}_u)\subset  \mathrm{Dom}(\mathbf{I}_+) \cap  H^1_+, \quad 1 \leq j \leq N, \quad \mathrm{and}\quad \bigoplus_{1 \leq j \leq N}  \mathbb{C}\mathfrak{h}_j =  \mathscr{H}_{\mathrm{pp}}(\mathbf{L}_u) .
\end{equation}Then   \eqref{idenhL2nvhinprod}    yields that $ e^{-i \alpha }\widehat{\mathfrak{h}}_j (0^+)   \langle v , \mathfrak{h}_j \rangle_{L^2_+} + |\langle v , \mathfrak{h}_j \rangle_{L^2_+}|^2 = 2 \pi \| \mathfrak{h}_j \|_{L^2_+}^2 = 2 \pi >0$. We thus have
\begin{equation}\label{vnotperphj}
 \langle v , \mathfrak{h}_j \rangle_{L^2_+} \ne 0, \quad \forall j =1,2, \cdots, N.
\end{equation}We set $v_0:=v \in \mathscr{H}_{\mathrm{pp}}(\mathbf{L}_u) $ and  $v_{n+1}:=\mathbf{L}_u^{n+1} (v) = \mathbf{L}_u (v_n ) \in     \mathscr{H}_{\mathrm{pp}}(\mathbf{L}_u)$,  $\forall n \in \mathbb{N}$. We claim that
\begin{equation}\label{vninHpp,DGG-1Hpp}
v_n \in \mathscr{H}_{\mathrm{pp}}(\mathbf{L}_u) \cap  \mathrm{Dom}(\mathbf{G}) \cap  \mathbf{G}^{-1} \left( \mathscr{H}_{\mathrm{pp}}(\mathbf{L}_u) \right), \quad \forall n \in \mathbb{N}=\mathbb{Z}\cap  [0, +\infty).
\end{equation}
\begin{proof}[Proof of formula \eqref{vninHpp,DGG-1Hpp}]
The proof proceeds by induction. The base case $n=0$ of \eqref{vninHpp,DGG-1Hpp} follows directly from formula \eqref{vinDomGandvGvinHpp}. Assume that formula \eqref{vninHpp,DGG-1Hpp} holds for some $n=k \in \mathbb{N}$, i.e. 
\begin{equation}\label{vkinHppfor}
v_k \in \mathscr{H}_{\mathrm{pp}}(\mathbf{L}_u) \cap  \mathrm{Dom}(\mathbf{G})  \subset H^1_+   \cap  \mathrm{Dom}(\mathbf{G})\quad \mathrm{and} \quad \mathbf{G}v_k \in \mathscr{H}_{\mathrm{pp}}(\mathbf{L}_u) \subset H^1_+.
\end{equation}Then $v_k \in \mathrm{Dom}^1(\mathbf{G})$ by  \eqref{DomsGDefSect2}. Then Lemma $\ref{LuD1GinDGLem}$ yields that $v_{k+1} = \mathbf{L}_u (v_k) \in   \mathrm{Dom}(\mathbf{G})$ and
\begin{equation}
\mathbf{G}v_{k+1} = \mathbf{G}\mathbf{L}_u (v_k) =  \mathbf{L}_u  \mathbf{G} (v_k)  +  i v_k- \left(e^{-i \alpha} \widehat{v}_k(0^+) + \langle v_k, v \rangle_{L^2_+} \right) \tfrac{iv}{2\pi} \in \mathscr{H}_{\mathrm{pp}}(\mathbf{L}_u)
\end{equation}by formulas \eqref{ForcommGbfLu} and \eqref{vkinHppfor}. So formula \eqref{vninHpp,DGG-1Hpp} holds for $n=k+1$. 
\end{proof}
\noindent Using \eqref{DomsGDefSect2}, we can indeed strengthen  \eqref{vninHpp,DGG-1Hpp}  to the following formula
\begin{equation}\label{vinD1GandHpp}
\mathbf{L}_u^n (v) = v_n \in \mathscr{H}_{\mathrm{pp}}(\mathbf{L}_u) \cap  \mathrm{Dom}^1 (\mathbf{G}) \cap  \mathbf{G}^{-1} \left( \mathscr{H}_{\mathrm{pp}}(\mathbf{L}_u) \right), \quad \forall n \in \mathbb{N}.
\end{equation}Assume that there exist $\mathfrak{c}_0, \mathfrak{c}_1, \cdots, \mathfrak{c}_{N-1} \in \mathbb{C}$ such that $\sum_{n=0}^{N-1} \mathfrak{c}_n \mathbf{L}_u^n (v)=\sum_{n=0}^{N-1} \mathfrak{c}_n v_n=0$. We have
\begin{equation}\label{LineSyssumcnlamnk=0}
 \langle v , \mathfrak{h}_k \rangle_{L^2_+} \sum_{n = 0}^{N-1} \mathfrak{c}_n \lambda_k^n  = \sum_{n = 0}^{N-1} \mathfrak{c}_n  \langle  v , \mathbf{L}_u^n ( \mathfrak{h}_k) \rangle_{L^2_+} = \sum_{n = 0}^{N-1} \mathfrak{c}_n  \langle  \mathbf{L}_u^n (v) , \mathfrak{h}_k \rangle_{L^2_+}  = 0, \quad \forall k =1,2, \cdots, N,
\end{equation}by \eqref{hjinDomI+,H1+}. 
Since $\lambda_k < \lambda_{k+1}$, $1 \leq k \leq  N-1$, the  Vandermonde matrix $(\lambda_k^{j-1})_{1 \leq k, j \leq N}$ is invertible. Combined with \eqref{vnotperphj}, this implies that  \eqref{LineSyssumcnlamnk=0} has only the trivial solution: $\mathfrak{c}_0 = \mathfrak{c}_1 = \cdots = \mathfrak{c}_{N-1} = 0$. So $\{v_0, v_1, \cdots, v_{N-1}\}$ is a basis for $\mathscr{H}_{\mathrm{pp}}(\mathbf{L}_u)$. Formula \eqref{vinD1GandHpp} yields that
\begin{equation}
\mathscr{H}_{\mathrm{pp}}(\mathbf{L}_u) = \mathrm{span}\{v_0, v_1, \cdots, v_{N-1}\} \subset \mathrm{Dom}^1 (\mathbf{G}) \cap  \mathbf{G}^{-1} \left( \mathscr{H}_{\mathrm{pp}}(\mathbf{L}_u) \right).
\end{equation}As a consequence, we obtain the assertion $(\mathbf{b})$. The proof of Theorem $\ref{MainThmSpectChar}$ is completed.
\end{proof}

\noindent We close this subsection by proving Corollary $\ref{CorQuPu}$.

\begin{proof}
  Fix \(u\in\mathcal U_N\).  Suppose that
  \(u=P/Q=p/q\), where \(P,Q,p,q\in\mathbb C[X]\), both \(Q\) and \(q\)
  are monic of degree \(N\) with roots in \(\mathbb C_-\), and both
  pairs satisfy the algebraic identity~\eqref{PQidenSec5}.  By
  Proposition~\ref{propSpecDecomLu}, we  have
 $ \frac{\mathbb C_{\le N-1}[X]}{Q}
    =
    \mathscr H_{\mathrm{pp}}(\mathbf L_u)
    =
    \frac{\mathbb C_{\le N-1}[X]}{q}$.  
  Since \(Q\) and \(q\) are both monic,
  Lemma~\ref{LemqtoC[X]ovqbij} forces $Q=q$. Moreover,  formula~\eqref{QischarGHppbfLu}
  implies that 
 $Q(x)=q(x)
    =\det\bigl(x-\mathbf G|_{\mathscr H_{\mathrm{pp}}(\mathbf L_u)}\bigr)
    =Q_u(x)$, for any $ x \in \mathbb{R}$.   
  Finally, we obtain   $P=p=uQ_u=P_u$.
\end{proof}

\subsection{Action–Angle Mapping}\label{AcAngMapSubs5.2}

\noindent
Fix an integer $N \geq 1$ and a dark $N$-soliton potential
$u = e^{i\alpha} + v \in \mathcal{U}_N$, where $\alpha \in \mathbb{T}$
and $v \in  \mathrm{Dom}^{\infty} (\mathbf{G})$. By Lemma~\ref{lemLuselfadj} and
Proposition~\ref{propSpecDecomLu}, the Lax operator $\mathbf{L}_u$
is self-adjoint on $L^2_+$ and has exactly $N$ distinct  eigenvalues. We
denote them in increasing order by
\begin{equation}\label{lamjuEigendeffor}
\lambda_1^u < \lambda_2^u < \cdots < \lambda_N^u,
\qquad
\sigma_{\mathrm{pp}}(\mathbf{L}_u)
=
\{\lambda_1^u , \lambda_2^u , \cdots , \lambda_N^u \}
\subset
\mathbb{R}\setminus \{-1\}.
\end{equation}
  Moreover, each of these eigenvalues is simple. 
For each $j = 1, 2, \cdots, N$, let $\phi_j^u$ be the unique eigenfunction
\begin{equation}\label{phijudefFor}
\phi_j^u \in \mathrm{Ker}(\lambda_j^u -\mathbf{L}_u)
\quad\text{such that}\quad
\|\phi_j^u \|_{L^2_+} = 1
\quad\text{and}\quad
\mathbf{p}_j(u):= \langle \phi_j^u, v \rangle_{L^2_+} > 0.
\end{equation}
  In this subsection, we carry out the construction of the action-angle mapping $\boldsymbol{\Phi}_N: \mathcal{U}_N \to
\mathbb{T}\times\Omega_N\times\mathbb{R}^{1\times N}$ introduced in Theorem~\ref{ActionAngleMainThm}, as well as the inverse spectral matrix $\mathbf{M}(u) \in \mathbb{C}^{N\times N}$ from Definition~\ref{bfMVfrakXYdefintro}.
  We then derive an explicit formula expressing each dark
\(N\)-soliton \(u\in\mathcal U_N\) in terms of its spectral data
$\{(\lambda_j^u,\phi_j^u)\}_{1\leq j \leq N}$ and the generator
$\bigl(-i\mathbf{G}, \mathrm{Dom}(\mathbf{G}) \bigr)$ given by~\eqref{DefGphiSect2}. \\
 
 \noindent Using Theorem \ref{MainThmSpectChar}, formula \eqref{SpecDeHacHscHpp}, and Lemma \ref{LemSpGC_-}, we obtain the following  properties: 
\begin{equation}\label{HppInvbyGForISM}
\mathscr{H}_{\mathrm{pp}}(\mathbf{L}_u)\subset 
\mathrm{Dom}^{\infty}(\mathbf{G}) , \quad  \mathbf{G} \bigl(\mathscr{H}_{\mathrm{pp}}(\mathbf{L}_u) \bigr) = \mathscr{H}_{\mathrm{pp}}(\mathbf{L}_u) \quad \mathrm{and} \quad v=u-e^{i \alpha} \in \mathscr{H}_{\mathrm{pp}}(\mathbf{L}_u).
\end{equation}Formulas \eqref{phijudefFor} and \eqref{HppInvbyGForISM} imply that
$\phi_j^u \in \mathrm{Ker}(\lambda_j^u - \mathbf{L}_u)\subset \mathscr{H}_{\mathrm{pp}}(\mathbf{L}_u) \subset \mathrm{Dom}(\mathbf{I}_+) \cap  H^{\infty}_+$. So we have
\begin{equation}\label{pjIjphijuForISM}
e^{-i \alpha} \mathbf{I}_+(\phi_j^u) \mathbf{p}_j(u)  + \mathbf{p}_j(u)^2 =   e^{-i \alpha }\widehat{\phi_j^u}(0^+)   \langle v , \phi_j^u \rangle_{L^2_+} + |\langle v , \phi_j^u \rangle_{L^2_+}|^2= 2 \pi \| \phi_j^u \|_{L^2_+}^2 = 2 \pi  
\end{equation}by identity \eqref{idenhL2nvhinprod}.
 Let $\mathbf{M}(u) = (\mathbf{M}_{kj}(u))_{1 \leq k, j \leq N} \in \mathbb{C}^{N \times N}$ denote the matrix of the restricted operator    $\mathbf{G}|_{\mathscr{H}_{\mathrm{pp}}(\mathbf{L}_u)}: \mathscr{H}_{\mathrm{pp}}(\mathbf{L}_u) \to \mathscr{H}_{\mathrm{pp}}(\mathbf{L}_u)$ with respect to the orthonormal eigenbasis $\{\phi_1^u, \phi_2^u, \cdots, \phi_N^u \}$, i.e. 
\begin{equation}\label{MkjForInvSpMatr} 
\mathbf{M}_{kj}(u):=\langle \mathbf{G} \phi_j^u, \phi_k^u \rangle_{L^2_+}, \quad  \forall k,j \in \{1,2, \cdots, N\}.
\end{equation}

\begin{defi}\label{InvSPmADef}
The matrix $\mathbf{M}(u) = (\mathbf{M}_{kj}(u))_{1 \leq k, j \leq N} \in \mathbb{C}^{N \times N}$, defined by its entries in  \eqref{MkjForInvSpMatr}, is referred to as the \emph{inverse spectral matrix} of  $u \in \mathcal{U}_N$. The \emph{action} and \emph{angle} variables are defined by
\begin{equation}\label{gammajudeffor}
  \boldsymbol{\lambda}_j(u ) = \lambda_j^u ,  \qquad \boldsymbol{\gamma}_j(u):= \operatorname{Re} \mathbf{M}_{jj}(u) = \operatorname{Re}\langle \mathbf{G}\phi_j^u, \phi_j^u \rangle_{L^2_+}, \qquad  \forall  j \in \{1,2, \cdots, N\}.
\end{equation}
 The \emph{total phase} $ \boldsymbol{\alpha}(u)\in\mathbb{T}
 $ is defined by
$e^{i\boldsymbol{\alpha}(u)} :=  \lim_{|x|\to+\infty} u(x)  $.   The limit is well defined by \eqref{uinftyDarSoli}.
\end{defi}

\begin{rem}\label{CoinRemInvSpMat}
 By Lemma~\ref{MentlamgamFor}, the inverse spectral matrix $\mathbf{M}(u)   \in \mathbb{C}^{N \times N}$ defined here via the matrix representation of $\mathbf{G}|_{\mathscr{H}_{\mathrm{pp}}(\mathbf{L}_u)}$ coincides with the one introduced in Definition~\ref{bfMVfrakXYdefintro}.
\end{rem}

\noindent According  to Proposition $\ref{propSpecDecomLu}$ and Corollary $\ref{CorQuPu}$,   the pure point   subspace $\mathscr{H}_{\mathrm{pp}}(\mathbf{L}_u)$ and the inverse spectral matrix $\mathbf{M}(u)$ are intimately connected via the polynomial $Q_u$ defined by \eqref{CharPolyQuGHpp1for}:
\begin{equation}\label{HppLuMuRelFor}
\mathscr{H}_{\mathrm{pp}}(\mathbf{L}_u)  = \tfrac{\mathbb{C}_{\leq N-1}[X]}{Q_u} , \quad Q_u(x) = \det (x - \mathbf{G}|_{ \mathscr{H}_{\mathrm{pp}}(\mathbf{L}_u)}) = \det (x - \mathbf{M}(u)), \quad \forall x \in \mathbb{R}.
\end{equation}Lemma $\ref{LemSpGC_-}$ yields that $(Q_u)^{-1}\{0\} = \sigma (\mathbf{M}(u)) \subset \sigma_{\mathrm{p}}(\mathbf{G}) \subset 
\mathbb{C}_-$. Then we define
\begin{equation}\label{PuTheQuInvSpM}
P_u:= u Q_u \in \mathbb{C}[X] \quad \mathrm{and} \quad \theta:=\theta_{Q_u} = \tfrac{\overline{Q}_u}{Q_u} \in 1 + \tfrac{\mathbb{C}_{\leq N-1}[X]}{Q_u} \subset 1 + \mathrm{Dom}^{\infty}(\mathbf{G}).
\end{equation}Then $(P_u, Q_u) \in \mathbb{C}[X] \times \mathbb{C}[X]$ satisfies identity \eqref{PQiden+} and $(\theta, u) \in \mathcal{Z}^{\infty}_+ \times \mathcal{Z}^{\infty}_+$ satisfies \eqref{DTheta|u|2-1}, i.e. 
\begin{equation}\label{Dthe|u|2-12nd}
D \theta = \theta(|u|^2 - 1) \in \tfrac{\mathbb{C}_{\leq 2N-2}[X]}{Q_u^2} \subset  \mathrm{Dom}^{\infty}(\mathbf{G}).
\end{equation}The next lemma presents an identity that connects $\theta-1 \in \mathrm{Dom}^{\infty}(\mathbf{G})$  and the Lax operator $\mathbf{L}_u$.
\begin{lem}\label{Lu1-thetaLem}
The inner function $\theta = \theta_{Q_u}\in 1 + \tfrac{\mathbb{C}_{\leq  N-1}[X]}{Q_u}$, defined in formula \eqref{PuTheQuInvSpM},  satisfies 
\begin{equation}\label{Lufor1-theta}
\mathbf{L}_u(1 - \theta) = \theta - 1 - e^{-i\alpha} v \in \tfrac{\mathbb{C}_{\leq N-1}[X]}{Q_u} \subset \mathrm{Dom}^{\infty}(\mathbf{G}).
\end{equation} 
\end{lem}
\begin{proof}
From formulas \eqref{ueiapluvDarSoli}, \eqref{HppLuMuRelFor} and \eqref{PuTheQuInvSpM}, we obtain that both $\overline{P}_u - e^{-i\alpha}Q_u$ and $P_u - e^{i\alpha}Q_u$ are polynomials of degree at most $N-1$. Since $\frac{\mathbb{C}_{\leq N-1}[X]}{\overline{Q}_u} \subset H^{\infty}_-$, it follows that
\begin{equation*}
\overline{u}(1- \theta) = \tfrac{\overline{P}_u - e^{-i\alpha}\overline{Q}_u}{\overline{Q}_u} - \tfrac{\overline{P}_u - e^{-i\alpha}Q_u}{Q_u} \Longrightarrow \mathbf{T}_{\overline{u}}(1 - \theta)  = - \tfrac{\overline{P}_u}{Q_u} + e^{-i \alpha} \in \tfrac{\mathbb{C}_{\leq N-1} [X]}{Q_u}  .
\end{equation*}So $\mathbf{T}_u \mathbf{T}_{\overline{u}} (1- \theta) = u ( - \tfrac{\overline{P}_u}{Q_u} + e^{-i \alpha}) = - \theta |u|^2  + e^{- i \alpha} u\in \frac{\mathbb{C}_{\leq 2N-1}[X]}{Q_u^2}$. The proof is completed by using formulas \eqref{LuSect2} and \eqref{Dthe|u|2-12nd}.
\end{proof}
\noindent For any $\eta >0$ and $\varphi \in \mathscr{H}_{\mathrm{pp}}(\mathbf{L}_u) = \frac{\mathbb{C}_{\leq N-1}[X]}{Q_u}\subset \mathrm{Dom}^{\infty}(\mathbf{G}) \subset 
\mathrm{Dom}(\mathbf{I}_+)$, Lemma $3.4$ of \cite{SUNCMP2021} yields that
\begin{equation}\label{Seta*1for-theta}
\langle \mathbf{S}(\eta)^* \varphi,  1 -\theta \rangle_{L^2_+}   = \widehat{\varphi}(\eta) \quad \mathrm{and} \quad \mathbf{I}_+(\varphi)= \widehat{\varphi}(0^+) = \langle \varphi ,  1 -\theta \rangle_{L^2_+}.  
\end{equation}The following lemma gives more precise information on the spectral parameters of  $u \in \mathcal{U}_N$.

\begin{lem}\label{LemSpLutwoseg}
Every eigenvalue of  the Lax operator $\mathbf{L}_u$ is excluded from $[-1, 0]$.  Moreover,  we have
\begin{equation}\label{pjuvaluefor}
\langle v, \phi_j^u \rangle_{L^2_+} =\mathbf{p}_j(u) = \sqrt{\tfrac{2 \pi(1+\lambda_j^u)}{\lambda_j^u}} \quad \mathrm{and} \quad \langle 1-\theta , \phi_j^u \rangle_{L^2_+} = \widehat{\phi_j^u}(0^+)  =  -   \tfrac{\sqrt{2\pi}\mathrm{sign}(\lambda_j^u)}{e^{i\alpha}\sqrt{\lambda_j^u(1+\lambda_j^u)}}, 
\end{equation}for any $j \in \{1, 2,  \cdots, N\}$.
\end{lem}
\begin{proof}
For each $j \in \{1, 2,  \cdots, N\}$, formula \eqref{lamjuEigendeffor} implies that $\lambda_j^u \ne -1$. We combine formulas \eqref{phijudefFor}, \eqref{Lufor1-theta}, and \eqref{Seta*1for-theta} with the self-adjointness of the Lax operator $\mathbf{L}_u : H^1_+ \to L^2_+$ to obtain that
\begin{equation}\label{Iphijuinterofpju}
\mathbf{I}_+(\phi_j^u) = - \langle  \phi_j^u, \mathbf{L}_u (1 -\theta) + e^{-i \alpha} v\rangle_{L^2_+} = -\lambda_j^u \mathbf{I}_+(\phi_j^u) - e^{i \alpha} \mathbf{p}_j(u) \Longrightarrow \mathbf{I}_+(\phi_j^u) = - \tfrac{e^{i \alpha} \mathbf{p}_j(u)}{1+ \lambda_j^u},
\end{equation}for any $j \in \{1, 2,  \cdots, N\}$. Plugging \eqref{Iphijuinterofpju} into  \eqref{pjIjphijuForISM}, we obtain the first formula of \eqref{pjuvaluefor}, i.e.
\begin{equation}\label{lamjusignlamju}
\lambda_j^u\mathbf{p}_j(u)^2=2 \pi (1+\lambda_j^u) \Longrightarrow   \mathrm{sign}(1+ \lambda_j^u) = \mathrm{sign}( \lambda_j^u ).
\end{equation}So $\sigma_{\mathrm{pp}}(\mathbf{L}_u) \cap  [-1,0] = \emptyset$. The second formula of \eqref{pjuvaluefor} is deduced from \eqref{Seta*1for-theta},  \eqref{Iphijuinterofpju} and \eqref{lamjusignlamju}.
\end{proof}
\begin{defi}\label{SpCoorVecDefiXY}
The \emph{spectral coordinates vectors} $\mathfrak{X}(u)=(\mathbf{x}_1(u), \mathbf{x}_2(u), \cdots, \mathbf{x}_N(u))^{\mathrm{T}} \in \mathbb{R}^N$ and  $\mathfrak{Y}(u)=(\mathbf{y}_1(u), \mathbf{y}_2(u), \cdots, \mathbf{y}_N(u))^{\mathrm{T}} \in \mathbb{R}^N$ for the dark $N$-soliton $u \in \mathcal{U}_N$ are defined as follows:
\begin{equation}\label{xjuyjufor}
\mathbf{x}_j(u) := \sqrt{\tfrac{1 + \lambda_j^u}{\lambda_j^u}} \quad \mathrm{and} \quad \mathbf{y}_j(u) :=  \tfrac{\mathrm{sign}(\lambda_j^u)}{\sqrt{\lambda_j^u (1+ \lambda_j^u)}} , \quad \forall j \in \{1, 2,  \cdots, N\}.
\end{equation}
\end{defi}
\noindent With respect to the orthonormal eigenbasis $\{\phi_1^u,\phi_2^u, \cdots, \phi_N^u\}$ of $\mathscr{H}_{\mathrm{pp}}(\mathbf{L}_u) = \frac{\mathbb{C}_{\leq N-1}[X]}{Q_u}$, the functions $v$ and $1-\theta$ admit the following expansions:
\begin{equation}\label{vtheSpCoorExpanFor}
v = \sqrt{2\pi} \sum_{j=1}^N \mathbf{x}_j(u) \, \phi_j^u, \qquad 
1-\theta = -\sqrt{2\pi} \, e^{-i\alpha} \sum_{j=1}^N \mathbf{y}_j(u) \, \phi_j^u.
\end{equation}
    For any $u\in\mathcal U_N$, the next lemma expresses the inverse spectral
matrix $\mathbf{M}(u)    \in \mathbb{C}^{N \times N}$ from Definition~\ref{InvSPmADef} in terms of
the action-angle variables defined in~\eqref{gammajudeffor}.

\begin{lem}\label{MentlamgamFor}
Let $\mathbf{M}(u) =(\mathbf{M}_{kj}(u))_{1 \leq k, j \leq N} \in \mathbb{C}^{N \times N}$ be the inverse spectral matrix of $u \in \mathcal{U}_N$ as given in Definition $\ref{InvSPmADef}$. Then its entries are determined as follows:
\begin{equation}\label{Mkjinteoflamju}
\mathbf{M}_{kj}(u) =  \tfrac{i}{\lambda_k^u - \lambda_j^u} \sqrt{\tfrac{\lambda_j^u (1 + \lambda_k^u)}{\lambda_k^u (1 + \lambda_j^u)}}, \quad  \mathbf{M}_{jj} (u)  = \boldsymbol{\gamma}_j(u) - \tfrac{i}{2 \lambda_j^u (1 + \lambda_j^u)}, \quad \mathrm{if} \quad 1 \leq k \ne j \leq N.
\end{equation}
\end{lem}
\begin{rem}
Every off-diagonal entry of $\mathbf{M}(u)$, as well as the imaginary part of each diagonal entry, depends only on the eigenvalues $\{\lambda_j^u\}_{1 \leq j \leq N}$ of the Lax operator $\mathbf{L}_u$. In the next section, we will show that the real parts $\boldsymbol{\gamma}_j(u) = \operatorname{Re}(\mathbf{M}_{jj}(u))$ can be chosen independently of these eigenvalues.
\end{rem}
\begin{proof}[Proof of Lemma $\ref{MentlamgamFor}$]
Since $\phi_j^u \in \mathrm{Ker}(\lambda_j^u - \mathbf{L}_u)\subset \mathrm{Dom}^{\infty}(\mathbf{G})$,   formula
\eqref{ForcommGbfLu} in Lemma $\ref{LuD1GinDGLem}$ yields that
\begin{equation}\label{comGLuForphijuISM}
\left[\mathbf{G}, \mathbf{L}_u \right]  \phi_j^u = i \phi_j^u - \left(e^{-i \alpha} \widehat{\phi_j^u}(0^+) + \langle \phi_j^u , v \rangle_{L^2_+} \right) \tfrac{iv}{2\pi} \in  \mathscr{H}_{\mathrm{pp}}(\mathbf{L}_u), \quad 1 \leq j \leq N.
\end{equation}For any $k,j \in \{1,2, \cdots, N\}$, 
 we use $\mathbf{L}_u = \mathbf{L}_u^*$ and  the definition formula \eqref{MkjForInvSpMatr} to obtain that
\begin{equation}\label{FORlamj-lamkMkj}
(\lambda_j^u - \lambda_k^u) \mathbf{M}_{kj}(u) =  \langle  \mathbf{G}  \mathbf{L}_u   \phi_j^u, \phi_k^u \rangle_{L^2_+} - \langle  \mathbf{G}  \phi_j^u, \mathbf{L}_u  \phi_k^u \rangle_{L^2_+}  = \langle \left[\mathbf{G}, \mathbf{L}_u \right]  \phi_j^u, \phi_k^u \rangle_{L^2_+}.
\end{equation}The first formula in \eqref{Mkjinteoflamju} follows directly from \eqref{FORlamj-lamkMkj} by inserting the expressions from \eqref{pjuvaluefor} and \eqref{comGLuForphijuISM}. By formula \eqref{MkjForInvSpMatr} and Lemma $\ref{ImGffinprodLem}$, it follows from \eqref{pjuvaluefor} that $
\operatorname{Im} \mathbf{M}_{jj}(u) = \operatorname{Im} \langle \mathbf{G}\phi_j^u, \phi_j^u  \rangle_{L^2_+} = - \tfrac{1}{4\pi} |\widehat{\phi_j^u}(0^+)|^2 = - \tfrac{1}{2 \lambda_j^u (1 + \lambda_j^u)}<0$. The proof is completed.
\end{proof}

\noindent We now state the  inverse spectral formula for   dark $N$-soliton potentials.

\begin{prop}\label{PropInvSpecFormula}
For any $u \in \mathcal{U}_N$ with $N \in \mathbb{N}_+$,  let $\mathfrak{X}(u), \mathfrak{Y}(u) \in \mathbb{R}^N$ be its spectral coordinate vectors in Definition $\ref{SpCoorVecDefiXY}$, and let  $\mathbf{M}(u) = (\mathbf{M}_{kj}(u))_{1 \leq k, j \leq N} \in \mathbb{C}^{N \times N}$ be the inverse spectral matrix defined in \eqref{Mkjinteoflamju}. 
Then $\sigma(\mathbf{M}(u))  \subset \sigma_{\mathrm{p}}(\mathbf{G})\subset \mathbb{C}_-$ and   we have the following inverse spectral formula:
\begin{equation}\label{invSPEFormain}
u(x) = u(\infty) \Bigl( 1 + i \, \bigl\langle (\mathbf{M}(u) - x)^{-1} \mathfrak{X}(u),\, \mathfrak{Y}(u) \bigr\rangle_{\mathbb{C}^N} \Bigr), \quad \forall x \in \mathbb{R},
\end{equation}
where $u(\infty) = e^{i \alpha} \in \mathbb{S}^1$ is the asymptotic value defined in formula  \eqref{uinftyDarSoli}.
\end{prop}
\begin{proof}
Since $v \in \mathscr{H}_{\mathrm{pp}}(\mathbf{L}_u) = \frac{\mathbb{C}_{\leq N-1}[X]}{Q_u}$, we use formulas \eqref{InvForNewtonLeibniz}, \eqref{Seta*1for-theta},  \eqref{MkjForInvSpMatr} and  \eqref{vtheSpCoorExpanFor}
to obtain 
\begin{equation}\label{vundlinezfor}
\underline{v}(z) = \tfrac{1}{2 \pi i} \langle (\mathbf{G} - z )^{-1} v, 1 - \theta \rangle_{L^2_+} = i e^{i \alpha} \langle  (\mathbf{M}(u) - z)^{-1} \mathfrak{X}(u),  \mathfrak{Y}(u) \rangle_{\mathbb{C}^N}, \quad \forall z \in \mathbb{C}_+ .
\end{equation}By the boundary convergence  \eqref{fundertofinL2},     there exists a sequence $y_n \downarrow 0$ such that  
\begin{equation*}
  \lim_{n \to +\infty} \underline{v}(x + i y_n)  =  v(x) = u(x) - u(\infty)   \quad \mathrm{for} \quad \mathrm{ a.e. } \quad  x \in \mathbb{R}.
\end{equation*}The almost-everywhere validity of \eqref{invSPEFormain}, together with the continuity of both $u \in \mathcal{Z}_+^{\infty}$ and the right-hand side, implies that   formula \eqref{invSPEFormain} holds for all $x \in \mathbb{R}$.
\end{proof}

\noindent 
Using the inverse spectral formula \eqref{invSPEFormain}, we obtain the following result.

\begin{cor}\label{ActAnglMapInjcor}
For any $N \in \mathbb{N}_+$, the action-angle mapping  $\boldsymbol{\Phi}_N: \mathcal{U}_N \to
\mathbb{T}\times\Omega_N\times\mathbb{R}^{1\times N}$ is injective. 
\end{cor}

\begin{proof}
Let $u_1,u_2\in\mathcal U_N$ satisfy
$\boldsymbol{\Phi}_N(u_1)=\boldsymbol{\Phi}_N(u_2) \in \mathbb{T}\times\Omega_N\times\mathbb{R}^{1\times N}$. Then they have
the same total phase, the same actions
$\boldsymbol{\lambda}_j$, and the same angles
$\boldsymbol{\gamma}_j$, for any $j \in \{1,2, \cdots, N\}$. By Lemma~\ref{MentlamgamFor}, their inverse
spectral matrices coincide, and by Definition~\ref{SpCoorVecDefiXY}
their spectral coordinate vectors coincide. The right-hand side
of~\eqref{invSPEFormain} is therefore identical for $u_1 (x)$  and
$u_2 (x)$, $\forall x \in \mathbb{R}$. So $u_1=u_2$.
\end{proof}
 
\noindent Moreover, the action and angle variables are invariant under the phase rotation.

\begin{lem}\label{LemPhaseRotInv}
For any $\mathfrak{z} \in \mathbb{S}^1=\{z \in \mathbb{C}: |z|=1\}$ and $u \in \mathcal{U}_N$ with $N \geq 1$, the potential $\widetilde{u}:=\mathfrak{z}u$ is also a dark $N$-soliton of \eqref{CMdNLSfocu},   and we have $\mathbf{L}_{\widetilde{u}} = \mathbf{L}_u : H^1_+ \to L^2_+$, $Q_{\widetilde{u}} = Q_u \in \mathbb{C}[X]$. Furthermore,  
\begin{equation}\label{rotgamuisgamtildu}
\lambda_j^{\widetilde{u}} = \lambda_j^u \in \mathbb{R},  \quad \boldsymbol{\gamma}_j(\widetilde{u}) = \boldsymbol{\gamma}_j(u) \in \mathbb{R},\quad  \phi_j^{\widetilde{u}} = \mathfrak{z}\phi_j^u \in \mathscr{H}_{\mathrm{pp}}(\mathbf{L}_{u})=\mathscr{H}_{\mathrm{pp}}(\mathbf{L}_{\widetilde{u}}), 
\end{equation}$\forall j \in \{1,2, \cdots,N\}$.
   As a result,  the associated inverse spectral matrix and coordinate vectors coincide:   
\begin{equation}\label{rotMuisMtildeu}
\mathbf{M}(\widetilde{u}) = \mathbf{M}(u) \in \mathbb{C}^{ N \times N}, \quad \mathfrak{X}(\widetilde{u}) = \mathfrak{X}(u)\in \mathbb{R}^N , \quad \mathfrak{Y}(\widetilde{u}) = \mathfrak{Y}(u) \in \mathbb{R}^N .
\end{equation}  
\end{lem}

\noindent Two dark $N$-solitons $u, \widetilde{u} \in \mathcal{U}_N$ are called \emph{phase-congruent} if they differ by a global phase factor, i.e., $\widetilde{u} = \mathfrak{z} u$ for some $\mathfrak{z} \in \mathbb{S}^1$. 
Phase-congruence is the equivalence relation induced by the group action
\eqref{bbTactcalUNdef} in Remark \ref{CirGractUNrem}.
  The equivalence class of $u \in \mathcal{U}_N$ is  its orbit with respect to   \eqref{bbTactcalUNdef}, i.e. 
\begin{equation}\label{OrbEquiClu}
[u] := \{\mathfrak{z} u \in \mathcal{U}_N \mid \mathfrak{z} \in \mathbb{S}^1\}.
\end{equation}
  We denote by $\mathcal{U}_N / \mathbb{T}$ the corresponding \emph{quotient space} (the set of all equivalence classes). By virtue of the identities in   \eqref{rotMuisMtildeu}, the inverse spectral matrix and spectral coordinate vectors are well-defined on the quotient space $\mathcal{U}_N / \mathbb{T}$. That is, we can define
\begin{equation}\label{InvSpMatEquiu}
\mathbb{M}([u]):= \mathbf{M}(u) \in \mathbb{C}^{N \times N}, \quad \mathbb{X}([u]):= \mathfrak{X}(u) \in \mathbb{R}^{N}, \quad \mathbb{Y}([u]):= \mathfrak{Y}(u) \in \mathbb{R}^{N}.
\end{equation} 
  The phase ambiguity is exactly removed by passing to the quotient:
the normalized function \(u/u(\infty)\) is independent of the
representative and defines a well-defined map 
 $ \mathcal U_N/\mathbb T \ni [u] \longmapsto
u/u(\infty) \in  1 + H^{\infty}_+$.
   For any
representative $u$ of a class $[u]\in\mathcal U_N/\mathbb T$, the
normalized function $u/u(\infty) \in \mathcal Z^\infty_+$ is given by
\begin{equation}\label{InSpForQuoSp}
\tfrac{u(x)}{u(\infty)}
=
1 + i \bigl\langle
(\mathbb{M}([u]) - x)^{-1} \mathbb{X}([u]),
\mathbb{Y}([u])
\bigr\rangle_{\mathbb{C}^N},
\qquad
\forall x\in\mathbb R.
\end{equation}

\subsection{Inverse Spectral Transform}\label{InSpTrsubsect5.3}

\noindent 
This subsection is devoted to proving the bijectivity of the
action-angle mapping
$\boldsymbol{\Phi}_N\colon\mathcal U_N\to
\mathbb{T}\times \Omega_N\times\mathbb{R}^{1\times N}$ 
introduced in Theorem~\ref{ActionAngleMainThm}. Recall that this
mapping is constructed from the total phase, actions, and angles
defined in Definition~\ref{InvSPmADef}. 
  We then construct its inverse $\boldsymbol{\Phi}_N^{-1} \colon 
\mathbb{T}\times \Omega_N\times\mathbb{R}^{1\times N} \to \mathcal{U}_N$.

\noindent
Recall that $\mathrm{sign}(0):=0$ and $\mathrm{sign}(\xi):=\xi |\xi|^{-1}$, if $\xi \in \mathbb{R}^* = \mathbb{R}\setminus \{0\}$. For any $\lambda \in \mathbb{R}$,  we have
\begin{equation}\label{sgnx=sgn1+x}
\lambda \in \mathbb{R}\setminus [-1, 0] \; \Longleftrightarrow \; \mathrm{sign}(\lambda) = \mathrm{sign}(1 + \lambda).
\end{equation}For any $N \geq 1$, we define that $\Upsilon_N:=\mathbb{T} \times   \Omega_N \times \mathbb{R}^N$, where $\Omega_N$ is given by \eqref{OmegaNIntro}. We recall that
\begin{equation*}
\Omega_N:=\{(\lambda_1, \lambda_2, \cdots, \lambda_N) \in \left(\mathbb{R}\setminus [-1,0]\right)^{1 \times N}:   \lambda_1 < \lambda_2 < \cdots < \lambda_N  \}.
\end{equation*}Then   $\Upsilon_N$ is a real analytic manifold of real dimension $2N+1$. In this subsection, we fix two row vectors $ \vec{\lambda} := (\lambda_1, \lambda_2,  \cdots, \lambda_N) \in\Omega_N$ and $ \vec{\gamma}  := (\gamma_1, \gamma_2,   \cdots, \gamma_N) \in \mathbb{R}^{1 \times N}$. We then define
\begin{equation}\label{defXYFoclambda}
\begin{cases}
X=X_{\vec{\lambda}} := \bigl(\sqrt{\frac{1+\lambda_1}{\lambda_1}}, \sqrt{\frac{1+\lambda_2}{\lambda_2}},   \cdots, \sqrt{\frac{1+\lambda_N}{\lambda_N}} \bigr)^{\mathrm{T}} \in \mathbb{R}^N, \\
Y=Y_{\vec{\lambda}} := \bigl(\frac{\mathrm{sign}(\lambda_1)}{\sqrt{\lambda_1 (1+\lambda_1)}} , \frac{\mathrm{sign}(\lambda_2)}{\sqrt{\lambda_2 (1+\lambda_2)}},   \cdots,\frac{\mathrm{sign}(\lambda_N)}{\sqrt{\lambda_N (1+\lambda_N)}} \bigr)^{\mathrm{T}} \in \mathbb{R}^N .
\end{cases}
\end{equation}
\begin{defi}\label{zetajstanONBRN}
For any $j=1,2, \cdots, N$, let $\boldsymbol{\zeta}_j=(0, \cdots, 0,1,0 \cdots, 0)^{\mathrm{T}} \in \mathbb{R}^N$ denote the vector with $1$ in the 
$j$-th coordinate and $0$ in all other coordinates. Let $\mathbb{E}_{kj}:= \boldsymbol{\zeta}_k \boldsymbol{\zeta}_j^{\mathrm{T}} \in \mathbb{R}^{N \times N}$ denote the matrix with $1$ in the $kj$-entry and $0$ in all other entries.
\end{defi}

\noindent The diagonal matrix $V=V_{\vec{\lambda}} \in \mathbb{R}^{N \times N}$ is given by
\begin{equation}\label{Vdiag-}
V:=\sum_{j=1}^N \lambda_j \mathbb{E}_{jj} = \left(\begin{matrix}
\lambda_1 & 0 &   \cdots & 0\\
0 & \lambda_2 &  \cdots & 0\\
\vdots & \vdots & \ddots & \vdots \\
0  & 0   & \cdots &   \lambda_N\\
\end{matrix} \right) =V^* \in \mathrm{GL}_N(\mathbb{R})  .
\end{equation}The matrix $M=M(\vec{\lambda}, \vec{\gamma})=\sum_{k=1}^N \sum_{j=1}^N M_{kj} \mathbb{E}_{kj} \in \mathbb{C}^{N \times N}$ is defined as follows,
\begin{equation}\label{defMlamgamFoc}
M_{jj}:= \gamma_j - \tfrac{i}{2 \lambda_j (1 +\lambda_j)} ,  \quad 
M_{kj}: =    \tfrac{i }{\lambda_k - \lambda_j} \sqrt{\tfrac{\lambda_j (1 +  \lambda_k)}{\lambda_k (1 + \lambda_j)}} ,
  \quad \mathrm{if} \; 1 \leq k\ne j \leq N.
\end{equation}

\begin{lem}\label{Lem-[V,M]&M-M*}
The vectors $X, Y \in \mathbb{R}^N$ and the matrices $M, V \in \mathbb{C}^{N \times N}$ satisfy the following identities: 
\begin{equation}\label{Forfoc[V,M]&M-M*}
YY^* = i( M-M^* ) \in \mathbb{R}^{N \times N}, \qquad  [V, M] = i(XY^*V - \mathrm{I}_N) \in \mathbb{R}^{N \times N}.
\end{equation}
\end{lem}
\begin{proof}
If  $j, k \in \{1,2, \cdots, N\}$, then $\lambda_j (1 + \lambda_k)$ can be either positive or negative. So \eqref{sgnx=sgn1+x} yields that
\begin{equation}\label{M-Mstar}
M_{kj} - \overline{M}_{jk} = \tfrac{i}{\lambda_k - \lambda_j} \tfrac{|\lambda_j (1 + \lambda_k)|-|\lambda_k (1 + \lambda_j)|}{\sqrt{\lambda_j (1+ \lambda_k)\lambda_k (1+\lambda_j)}}=\tfrac{-i\mathrm{sign}(\lambda_k) \mathrm{sign}(\lambda_j)}{\sqrt{\lambda_k (1+ \lambda_k)\lambda_j (1+\lambda_j)}}=-i(YY^*)_{kj}.
\end{equation}Moreover, $([V,M])_{kj}=(\lambda_k - \lambda_j) M_{kj} = i\sqrt{\tfrac{\lambda_j (1 +  \lambda_k)}{\lambda_k (1 + \lambda_j)}}$, $1 \leq k \ne j \leq N$. Then we define that
\begin{equation}\label{beta=VY=X-Y}
\beta:=VY = (\sqrt{\tfrac{\lambda_1}{1+\lambda_1}},\sqrt{\tfrac{\lambda_2}{1+\lambda_2}}, \cdots, \sqrt{\tfrac{\lambda_N}{1+\lambda_N}})^{\mathrm{T}}=X-Y \in \mathbb{R}^N.
\end{equation}Consequently, $XY^*V=X\beta^* = \sum_{k=1}^N \sum_{j=1}^N \sqrt{\tfrac{\lambda_j (1 +  \lambda_k)}{\lambda_k (1 + \lambda_j)}}  \mathbb{E}_{kj}=\mathrm{I}_N - i[V,M]$.
\end{proof}

\noindent The next lemma shows that the spectrum of $M \in \mathbb{C}^{N \times N}$ is contained in the lower half-plane. 
 
 \begin{lem}\label{SpofMLemfoc}
Every eigenvalue of the matrix $M =\sum_{k=1}^N \sum_{j=1}^N M_{kj} \mathbb{E}_{kj} \in \mathbb{C}^{N \times N}$, as defined in formula \eqref{defMlamgamFoc}, must have a strictly negative imaginary part, i.e. $\sigma (M)  \subset \mathbb{C}_-=\{a+bi \in \mathbb{C} : b<0 \}$.
\end{lem}

\begin{proof}
If $MZ =\mu Z$ for some $\mu \in \mathbb{C}$ and $Z \in \mathbb{C}^N \setminus \{0_{\mathbb{C}^N}\}$,  formula  \eqref{Forfoc[V,M]&M-M*} in Lemma $\ref{Lem-[V,M]&M-M*}$  implies that $2 i \|Z\|^2_{\mathbb{C}^N} \operatorname{Im}\mu = \langle (M-M^*)Z, Z \rangle_{\mathbb{C}^N} = -i |\langle Z, Y\rangle_{\mathbb{C}^N}|^2 \leq 0$. Then $\sigma (M) \subset \mathbb{C}_-  \cup \mathbb{R}$. Assume that there exists $\mu \in \mathbb{R} \cap   \sigma (M)$, set $Z \in \mathrm{Ker}(\mu - M) \setminus \{0_{\mathbb{C}^N}\}$, then $Y^*Z=\langle Z, Y \rangle_{\mathbb{C}^N} = 0$ and $M^* Z = M Z = \mu Z$. So $\langle [M, V^{-1}]Z, Z \rangle_{\mathbb{C}^N} = \langle (M - \mu) V^{-1} Z, Z \rangle_{\mathbb{C}^N} =\langle  V^{-1} Z, (M^* - \overline{\mu}) Z \rangle_{\mathbb{C}^N} =0$. Moreover, the second formula of \eqref{Forfoc[V,M]&M-M*} yields that $0= \langle V^{-1}[V,M]V^{-1}Z, Z\rangle_{\mathbb{C}^N} = iZ^*V^{-1} \left(XY^*V -\mathrm{I}_N \right)V^{-1}Z = -i \|V^{-1}Z\|_{\mathbb{C}^N}^2$. So $V^{-1}Z = 0_{\mathbb{C}^N}$, which contradicts that $Z \ne 0_{\mathbb{C}^N}$. Thus, $\mathbb{R} \cap   \sigma (M) = \emptyset$ and $\sigma (M) \subset \mathbb{C}_-$.
\end{proof}

\noindent Now, we shall prove the surjectivity of the action-angle mapping  $\boldsymbol{\Phi}_N\colon\mathcal U_N\to
\mathbb{T}\times \Omega_N\times\mathbb{R}^{1\times N}$. We first treat the case $N=1$, where the inverse
spectral transform  is 
elementary.

\begin{lem}\label{lemInSpTrNis1}
For any $\alpha \in \mathbb{T}=\mathbb{R}/ 2 \pi \mathbb{Z}$, $\lambda \in \mathbb{R}\setminus [-1, 0]$, $\gamma \in \mathbb{R}$, we define  $M = \gamma - \frac{i}{2 \lambda (1+ \lambda)} \in \mathbb{C}_- $ and
\begin{equation}\label{dar1soliSurjFoc-}
  u(x):= e^{i\alpha}\left(1 + \tfrac{i}{\lambda(M-x)} \right) , \quad  \forall x \in \mathbb{R}.
\end{equation}Then  $u$ is a dark $1$-soliton, i.e. $u \in \mathcal{U}_1$. The Lax operator $\mathbf{L}_u$ has exactly one eigenvalue, namely $\lambda$, that is $\sigma_{\mathrm{pp}}(\mathbf{L}_u) =\{\lambda\} =\{\lambda_1^u\}$. Moreover, for any $x\in \mathbb{R}$, we have
\begin{equation}
Q_u(x)=\det\left(x - \mathbf{G}|_{\mathscr{H}_{\mathrm{pp}}(\mathbf{L}_u)}\right) = x-M \quad \mathrm{and} \quad \mathscr{H}_{\mathrm{pp}}(\mathbf{L}_u)  = \mathbb{C} (\tfrac{1}{Q_u} ) .
\end{equation}As a consequence, we have $\boldsymbol{\alpha}(u) =  \alpha$ and $\boldsymbol{\gamma}_1(u) =\gamma$. 
\end{lem}
\begin{proof}
For any $x\in \mathbb{R}$, we set $Q(x):=x - M$ and $\varphi(x):= \frac{1}{x -M}$. Since $Q^{-1}\{0\}=\{M\} \subset \mathbb{C}_-$,
formula  \eqref{KthetainDominftyG}    yields that  $\varphi \in \mathrm{Dom}^{\infty}(\mathbf{G})$. Set $P:=e^{i\alpha}( Q - \frac{i}{\lambda}) \in \mathbb{C}[X]$. Then $u = \tfrac{P}{Q} = e^{i\alpha} (1 - \tfrac{i}{\lambda} \varphi) \in \mathcal{Z}^{\infty}_+$ and
\begin{equation}\label{PQIDENda1SolisurjFoc-}
|P |^2 - |Q |^2 = \lambda^{-2} - 2 \lambda^{-1} \operatorname{Im}Q = 2 \operatorname{Im}Q = i (\overline{Q} \partial_x Q - Q \partial_x \overline{Q}).
\end{equation}Formula \eqref{PQIDENda1SolisurjFoc-} is equivalent to \eqref{PQiden+}. So we have $u \in \mathcal{U}_1$ and  $\mathscr{H}_{\mathrm{pp}}(\mathbf{L}_u) = \frac{\mathbb{C}}{Q} $ by formula \eqref{SpecDeHacHscHpp}. By Corollary $\ref{CorQuPu}$ and  formula  \eqref{CharPolyQuGHpp1for}, it follows from   \eqref{HppLuMuRelFor}  that
\begin{equation}
x - M = Q(x) = Q_u(x)=\det  \left(x - \mathbf{G}|_{\mathscr{H}_{\mathrm{pp}}(\mathbf{L}_u)}\right)  =  x- \mathbf{M}(u) .
\end{equation}Then $\gamma = \operatorname{Re}M = \operatorname{Re}\mathbf{M}(u) = \boldsymbol{\gamma}_1(u)$. Since $|\varphi|^2 = \frac{1}{2 i \operatorname{Im} M}(\varphi - \overline{\varphi}) \in H^{\infty}(\mathbb{R}; \mathbb{R})$, we have
 \begin{equation}
 e^{i \alpha} \mathbf{T}_{\overline{u}}(\varphi) = - \lambda \varphi \in \mathrm{Dom}^{\infty}(\mathbf{G}) \quad \mathrm{and} \quad \mathbf{T}_u \mathbf{T}_{\overline{u}}(\varphi) = - \lambda \varphi +i \varphi^2 = (D- \lambda) \varphi \in \mathrm{Dom}^{\infty}(\mathbf{G}), 
 \end{equation}by   \eqref{ZnHardyHninHnfor}. Since $\mathbf{L}_u \varphi = \lambda \varphi$, we have $\sigma_{\mathrm{pp}}(\mathbf{L}_u) = \{\lambda\}$. Finally \eqref{dar1soliSurjFoc-} yields that $e^{i \alpha} = u(\infty) = e^{i \boldsymbol{\alpha}(u)}  $.
\end{proof}

\noindent Let $Q\in \mathbb{C}[X]$ denote the characteristic polynomial of the matrix $M=  (M_{kj})_{1 \leq k, j \leq N} \in \mathbb{C}^{N \times N}$, whose entries are given by \eqref{defMlamgamFoc}. Then there exists a polynomial $\mathfrak{R}$ of degree at most $N-2$ such that
\begin{equation}\label{Qisdet(x-M)}
Q(x) = \det(x - M) = x^N - (\mathrm{tr}M)x^{N-1} + \mathfrak{R} (x), \qquad \forall x \in \mathbb{R}.
\end{equation}
  If $N=1$, we simply have $\mathfrak{R}=0$. 
 By Lemma $\ref{SpofMLemfoc}$, $Q$ is a monic polynomial of degree $N$ whose roots (which coincide with the spectrum of $M$) all lie in   $\mathbb{C}_-$. Henceforth, we shall assume $N\geq 2$.

\begin{defi}\label{defofKkjqkjfoc-}
For any $x \in \mathbb{R}$, $k,j \in \{1,2, \cdots, N\}$, let $\mathcal{K}_{kj}(x) \in \mathbb{C}^{(N-1) \times (N-1)}$ denote the submatrix of $M-x \in \mathrm{GL}_N(\mathbb{C})$ obtained by eliminating the $k$-th row and the $j$-th column. We set
\begin{equation}\label{qkjdeffoc-}
q_{kj}(x):= \det \mathcal{K}_{kj}(x).
\end{equation}
 \end{defi}
\noindent For any $k,j \in \{1,2, \cdots, N\}$, if $k \ne j$, then $q_{kj}$ is a polynomial of degree at most $N-2$; if $k=j$, then   $(-1)^{N-1}q_{jj} \in \mathbb{C}[X]$ is the characteristic polynomial of the matrix $\mathcal{K}_{jj}(0) \in \mathbb{C}^{(N-1) \times (N-1)}$, i.e. 
\begin{equation}\label{decompqjjrj<N-2}
(-1)^{N-1} q_{jj}(x) = \det(x-\mathcal{K}_{jj}(0) ) =  x^{N-1} -(\mathrm{tr}M - M_{jj})x^{N-2}+\mathfrak{r}_j(x), \quad \forall x \in \mathbb{R},
\end{equation}for some polynomial $\mathfrak{r}_j \in \mathbb{C}[X]$ such that $\deg \mathfrak{r}_j \leq N-3$. (When $N=2$, we set $\mathfrak{r}_j =0$.) Since $\sigma(M) \subset \mathbb{C}_-$ by Lemma $\ref{SpofMLemfoc}$, A standard application of Cramer's rule gives the adjugate matrix representation
\begin{equation}\label{adj(M-x)}
\mathrm{adj}(M-x) = \det(M-x) (M-x)^{-1} =\sum_{k=1}^N  \sum_{j=1}^N (-1)^{k+j} q_{jk}(x) \mathbb{E}_{kj} \in \mathrm{GL}_N(\mathbb{C}), \quad \forall x \in \mathbb{R}.
\end{equation}For any $a=(a_1, a_2, \cdots, a_N)^{\mathrm{T}} , b=(b_1, b_2, \cdots, b_N)^{\mathrm{T}} \in \mathbb{C}^N$, the function $\mathfrak{g}_{ab}   $ is defined as follows:
\begin{equation}\label{gabdeffoc}
\mathfrak{g}_{ab}(x):=\langle (M-x)^{-1}a, b \rangle_{\mathbb{C}^N} = \sum_{k=1}^N  \sum_{j=1}^N \frac{a_j \overline{b}_k q_{jk}(x)}{ (-1)^{N+k+j} Q(x) } \in \mathbb{C}, \quad \forall x \in \mathbb{R}.
\end{equation}By  formula 
 \eqref{KthetainDominftyG}, it follows from the fact $Q^{-1}\{0\} \subset \mathbb{C}_-$ that  
\begin{equation}\label{gabinDominftyG}
\mathfrak{g}_{ab} \in  \tfrac{\mathbb{C}_{\leq N-1}[X]}{Q} \subset \mathrm{Dom}^{\infty}(\mathbf{G}) \subset H^{\infty}_+ \cap   \mathrm{Dom}(\mathbf{G}).
\end{equation}Its first and second derivatives are  
\begin{equation}\label{dergabForFoc-}
\partial_x \mathfrak{g}_{ab} (x) = \langle (M-x)^{-2} a , b \rangle_{\mathbb{C}^N}, \quad \partial_x^2  \mathfrak{g}_{ab} (x) = 2  \langle (M-x)^{-3} a , b \rangle_{\mathbb{C}^N}, \quad \forall x \in \mathbb{R}.
\end{equation}Moreover, the following product-to-sum identity holds.
\begin{lem}\label{ProdtoSumforLem-}
Recall that $Y = Y_{\vec{\lambda}} \in \mathbb{R}^N$ is defined in \eqref{defXYFoclambda}. 
For any $a, b  \in \mathbb{C}^N$, we have
\begin{equation}\label{ForfocProdtoSum-}
\overline{\mathfrak{g}}_{a Y}\mathfrak{g}_{b Y}=  i \left(\overline{\mathfrak{g}}_{ab} - \mathfrak{g}_{ba}\right) \in H^{\infty}(\mathbb{R}; \mathbb{C}).
\end{equation} 
\end{lem}
\begin{proof}
The first formula in  \eqref{Forfoc[V,M]&M-M*} shows that $M-x-(M^*-x)= -i   Y Y^{\mathrm{T}} \in i \mathbb{R}^{N \times N}$. Then we have
\begin{equation*}
\overline{\mathfrak{g}}_{a Y}(x)\mathfrak{g}_{b Y}(x)=  a^* (M^*-x)^{-1} Y Y^{\mathrm{T}} (M-x)^{-1} b =  i\left( a^* (M^*-x)^{-1} b -  a^* (M-x)^{-1} b \right),
\end{equation*}by using formulas \eqref{C^Ndotprod} and \eqref{gabdeffoc}, for any $x \in \mathbb{R}$.
\end{proof}

\noindent Fix an arbitrary $\alpha \in \mathbb{T} = \mathbb{R}/2\pi\mathbb{Z}$. Recall that the vectors $X = X_{\vec{\lambda}} \in \mathbb{R}^N$ and $Y = Y_{\vec{\lambda}} \in \mathbb{R}^N$ are defined in \eqref{defXYFoclambda}. By Lemma \ref{SpofMLemfoc}, the matrix $M = M(\vec{\lambda}, \vec{\gamma}) \in \mathbb{C}^{N \times N}$ given by \eqref{defMlamgamFoc} has all its eigenvalues in the lower half-plane $\mathbb{C}_-$. We then define a function $u = u_{\alpha, \vec{\lambda}, \vec{\gamma}} : \mathbb{R} \to \mathbb{C}$ by
\begin{equation}\label{ulamgamalpdef}
    u (x) := e^{i\alpha} \bigl(1 + i \mathfrak{g}_{XY}(x)\bigr) = e^{i\alpha} \Bigl(1 + i \bigl\langle (M - x)^{-1} X,\; Y \bigr\rangle_{\mathbb{C}^N} \Bigr), \qquad \forall x \in \mathbb{R}, 
\end{equation}
where the function  $\mathfrak{g}_{XY} : \mathbb{R} \to \mathbb{C}$ is given by  formula  \eqref{gabdeffoc}.

\begin{lem}\label{ulagaFoclinUNsurj}
    The function $u = u_{\alpha, \vec{\lambda}, \vec{\gamma}}$ defined in \eqref{ulamgamalpdef} is a dark $N$-soliton, i.e., $u \in \mathcal{U}_N$.
\end{lem}
\begin{proof}
Set $\mathfrak{b}:= Q \mathfrak{g}_{XY} \in \mathbb{C}[X]$. Formula \eqref{gabinDominftyG} yields that $\deg \mathfrak{b} \leq N-1$. We define  
\begin{equation}\label{uisPovQSurj-}
 P :=e^{i \alpha} (Q + i \mathfrak{b}) \in \mathbb{C}[X].   \quad \mathrm{Then} \quad u = e^{i\alpha} \bigl(1 + i \mathfrak{g}_{XY} \bigr) = \tfrac{P}{Q}.
\end{equation}Moreover, we have $\deg Q =\deg P =N$. By formulas \eqref{gabdeffoc} and \eqref{ForfocProdtoSum-}, it follows from \eqref{ulamgamalpdef} that
\begin{equation}\label{absusq-1=gXbeta}
|u|^2 - 1 = i(\mathfrak{g}_{XY} - \overline{\mathfrak{g}}_{XY} ) + \overline{\mathfrak{g}}_{XY} \mathfrak{g}_{XY}  = 2 \operatorname{Im} \left( \mathfrak{g}_{XX} - \mathfrak{g}_{XY} \right)=2 \operatorname{Im} \mathfrak{g}_{X \beta} \in H^{\infty}(\mathbb{R}; \mathbb{R}),
\end{equation}where $\beta = X-Y = VY \in \mathbb{R}^N$ by formula \eqref{beta=VY=X-Y}. Here we recall that 
\begin{equation}\label{betajXjrecalfoc-}
\beta = (\beta_1, \beta_2, \cdots, \beta_N)^{\mathrm{T}}  , \quad  X = (x_1, x_2, \cdots, x_N)^{\mathrm{T}}   , \quad \mathrm{where} \quad \beta_j  =   x_j^{-1} :=  \sqrt{\tfrac{\lambda_j}{1 + \lambda_j}},
\end{equation}for any $j \in \{1, 2, \cdots, N\}$. Fix  $x\in \mathbb{R}$. If $k, j \in \{1, 2, \cdots, N\}$, let $(M-x)_{kj} \in \mathbb{C}$ denote the $kj$-entry of the matrix $ M - x \in \mathrm{GL}_N(\mathbb{C})$. Let $\mathcal{K}_{kj}(x) \in \mathbb{C}^{(N-1) \times (N-1)}$ be the submatrix defined in Definition~\ref{defofKkjqkjfoc-}. Then, using \eqref{Qisdet(x-M)},  \eqref{betajXjrecalfoc-} and Jacobi's formula for differentiating determinants, we obtain that
\begin{equation}\label{J1defuinUN-}
    \mathbb{J}_1 := \sum_{j=1}^{N} (-1)^{j+j} x_j \beta_j \det \mathcal{K}_{jj}(x) 
                   = \sum_{j=1}^{N} \det \mathcal{K}_{jj}(x) 
                   = \sum_{j=1}^{N} q_{jj}(x)
                   = (-1)^{N+1} \partial_x Q(x).
\end{equation}Assume that $1 \leq k\ne j \leq N$, then we have $i(\lambda_k - \lambda_j)(M-x)_{jk} = i(\lambda_k - \lambda_j)M_{jk}=x_j \beta_k$ by   \eqref{defMlamgamFoc} and \eqref{betajXjrecalfoc-}. Using the cofactor expansion formula for  $\det(M-x)$, we compute
\begin{equation}\label{J2defuinUNfoc}
\mathbb{J}_2: =   \sum_{1 \leq k \ne j \leq N} (-1)^{k+j} x_j \beta_k \det \mathcal{K}_{jk}(x)  = i\sum_{1 \leq k , j \leq N} (-1)^{k+j} (\lambda_k - \lambda_j)(M-x)_{jk} \det \mathcal{K}_{jk}(x)  
   = 0.
\end{equation}From formulas \eqref{Qisdet(x-M)},  \eqref{qkjdeffoc-},  \eqref{gabdeffoc},  \eqref{J1defuinUN-} and \eqref{J2defuinUNfoc}, for any $x\in \mathbb{R}$, we obtain that
\begin{equation}\label{gXbetaexplicfoc}
\det(M - x) \mathfrak{g}_{X \beta} (x) = \sum_{1 \leq k, j \leq N} (-1)^{k +j } x_j \beta_k \det \mathcal{K}_{kj}(x) = \mathbb{J}_1 + \mathbb{J}_2 = (-1)^{N+1} \partial_x Q (x).
\end{equation}Plugging \eqref{gXbetaexplicfoc} into formula \eqref{absusq-1=gXbeta}, we use \eqref{uisPovQSurj-} to obtain that
\begin{equation}
|u|^2-1 = - 2 \operatorname{Im} \tfrac{\partial_x Q}{Q} \Longleftrightarrow |P|^2 - |Q|^2 = i \left(\overline{Q} \partial_x Q - Q \partial_x \overline{Q} \right).
\end{equation}So $(P, Q) \in \mathbb{C}[X] \times \mathbb{C}[X]$ satisfies \eqref{PQiden+}.   
  We have $u = \frac{P}{Q} \in \mathcal{U}_N$ by  
     \eqref{Qisdet(x-M)} and Lemma $\ref{SpofMLemfoc}$.
\end{proof}

\noindent
It remains to verify that the dark $N$-soliton potential $u \in \mathcal{U}_N$
defined by \eqref{ulamgamalpdef} has the prescribed action-angle
data. To this end, we shall construct an eigenbasis
$\{\phi_j^u\}_{1\leq j\leq N}$ of the Lax operator $\mathbf{L}_u $,
satisfying the phase convention in \eqref{phijudefFor}. From these
eigenvectors we compute the actions
$\boldsymbol{\lambda}_j(u)$ and angles
$\boldsymbol{\gamma}_j(u)$ according to
Definition~\ref{InvSPmADef}, and we then prove that
\begin{equation}\label{Idenbdlamgam}
\boldsymbol{\alpha}(u)=\alpha,\qquad
\boldsymbol{\lambda}_j(u)=\lambda_j,\qquad
\boldsymbol{\gamma}_j(u)=\gamma_j
\qquad  \forall  j\in \{1, 2, \cdots,N \}.
\end{equation}
  Together with Lemma~\ref{ulagaFoclinUNsurj}, 
    this will complete
the proof of surjectivity of
  $ \boldsymbol{\Phi}_N : \mathcal{U}_N  \to  \mathbb{T}\times\Omega_N\times\mathbb{R}^{1\times N}$.  \\

\noindent Let $\{\boldsymbol{\zeta}_j\}_{1 \leq j \leq N}$ denote the standard basis of \(\mathbb{R}^N\) (see Definition~\ref{zetajstanONBRN}) and assume that $N\geq 2$. For  each $j \in \{1, 2, \cdots, N\}$, we define a function $\mathfrak{f}_j : \mathbb{R} \to \mathbb{C}$ by
\begin{equation}\label{fjEigFuncFOC-}
    \mathfrak{f}_j(x) := \mathfrak{g}_{\boldsymbol{\zeta}_j Y}(x) = \bigl\langle (M - x )^{-1} \boldsymbol{\zeta}_j,\; Y \bigr\rangle_{\mathbb{C}^N}  , \quad \forall x \in \mathbb{R},
\end{equation}where the vector $Y = (y_1, y_2, \cdots, y_N)^{\mathrm{T}}\in \mathbb{R}^N$ with components $y_j:= \frac{  \mathrm{sign}(\lambda_j)}{\sqrt{\lambda_j (1+\lambda_j)}}$, for $1 \leq j \leq N$, is given by \eqref{defXYFoclambda}. The matrix  $M =(M_{kj})_{1 \leq k, j \leq N} \in \mathbb{C}^{N \times N}$ is defined in formula     \eqref{defMlamgamFoc}. 
  The next lemma gives the $L^2$-norm of $\mathfrak{f}_j $ and the right-hand limit of $\widehat{\mathfrak{f}}_j$, for all $j \in \{1,2, \cdots, N\}$.

\begin{lem}
For any $j \in \{1,2, \cdots, N\}$, the function $\mathfrak{f}_j \in \frac{\mathbb{C}_{\leq N-1}[X]}{Q} \subset \mathrm{Dom}^{\infty}(\mathbf{G})$ satisfies
\begin{equation}\label{L2normfjFOC-}
\|\mathfrak{f}_j\|_{L^2_+}^2 =  2 \pi \quad \mathrm{and} \quad \mathbf{I}_+(\mathfrak{f}_j) 
= 2 \pi i y_j
   =\tfrac{2 \pi i \mathrm{sign}(\lambda_j)}{\sqrt{\lambda_j (1+\lambda_j)}} . 
\end{equation} 
\end{lem}
 
\begin{proof}
 Let $\mathfrak{r}_j \in \mathbb{C}_{\leq N-3}[X]$ be defined in  \eqref{decompqjjrj<N-2}. (When $N=2$, we set $\mathfrak{r}_j = 0$.) We define
\begin{equation}\label{tilrjdefforFOC-}
\widetilde{\mathfrak{r}}_j (x) := (M_{jj} - \mathrm{tr} M) x^{N-2} + \mathfrak{r}_j (x) = (-1)^{N-1} q_{jj} (x) - x^{N-1}, \quad \forall x \in \mathbb{R}.
\end{equation}Then $\widetilde{\mathfrak{r}}_j$ is a polynomial of degree $N-2$. So $\frac{\widetilde{\mathfrak{r}}_j}{Q} \in \frac{\mathbb{C}_{\leq N-2}[X]}{Q} \subset H^{\infty}_+ \cap   L^1(\mathbb{R}; \mathbb{C})$.   Lemma $\ref{LemL2+capL10mode}$ yields that $ \int_{\mathbb{R}}\tfrac{\widetilde{\mathfrak{r}}_j}{Q} = 0$. From   formulas \eqref{decompqjjrj<N-2},  \eqref{gabdeffoc},  \eqref{ForfocProdtoSum-} and \eqref{fjEigFuncFOC-}, we obtain that
\begin{equation}\label{absfjxsqFORfoc-}
  |\mathfrak{f}_j(x)|^2 = 2 \operatorname{Im}\mathfrak{g}_{\boldsymbol{\zeta}_j \boldsymbol{\zeta}_j}(x) = 2 \operatorname{Im} \tfrac{q_{jj}(x)}{(-1)^N Q(x)} =  -2 \operatorname{Im} \left( \tfrac{x^{N-1} + \widetilde{\mathfrak{r}}_j(x)}{Q(x)} \right), \quad \forall x \in \mathbb{R}.
\end{equation}By Lemma $\ref{SpofMLemfoc}$ and formula \eqref{Qisdet(x-M)}, there exist  $z_1, z_2, \cdots, z_N \in \mathbb{C}_-$ such that $Q(x) = \prod_{k=1}^N (x - z_k)$. According to formula \eqref{J1defuinUN-}, we integrate \eqref{absfjxsqFORfoc-} and sum the resulting identities over $j$. This yields
\begin{equation}\label{sumL2norfj=NpiFOC-}
\sum_{j=1}^N \|\mathfrak{f}_j \|_{L^2_+}^2 = 2 (-1)^N \sum_{j=1}^N \int_{\mathbb{R}} \left( \operatorname{Im} \tfrac{ q_{jj} }{ Q } \right) = -2  \int_{\mathbb{R}}\left( \operatorname{Im}\tfrac{ \partial_x Q }{ Q } \right) = \sum_{k=1}^N \int_{\mathbb{R}} \tfrac{2 |\operatorname{Im} z_k| \mathrm{d}x}{(x - \operatorname{Re}z_k)^2 + |\operatorname{Im} z_k|^2} = 2 N \pi.
\end{equation}Since $  \int_{\mathbb{R}}\tfrac{\widetilde{\mathfrak{r}}_j}{Q} = 0$, the term $\|\mathfrak{f}_j \|_{L^2_+}^2 = 2 \int_{\mathbb{R}} \frac{x^{N-1} \operatorname{Im}Q(x)}{|Q(x)|^2} \mathrm{d}x$ is independent of  the index $j$. So we have
\begin{equation}\label{xN-1ImQxover|Q|2intFoc-}
\int_{\mathbb{R}} \tfrac{x^{N-1} \operatorname{Im}Q(x)}{|Q(x)|^2} =\pi \quad \mathrm{and} \quad \|\mathfrak{f}_j \|_{L^2_+}^2 = 2 \pi, \quad \forall j  \in \{1,2, \cdots, N\}.
\end{equation}If $1 \leq m \ne j \leq N$, formulas \eqref{qkjdeffoc-} and \eqref{tilrjdefforFOC-} imply that $q_{jm}$ and  $\widetilde{\mathfrak{r}}_j$  are both polynomials of degree at most $N-2$. Since $\deg Q = N$, we have $\lim_{x \to +\infty} \frac{x \widetilde{\mathfrak{r}}_j (x)}{Q(x)} = \lim_{x \to +\infty} \frac{x q_{jm} (x)}{Q(x)} = 0$. Fix $j \in \{1,2, \cdots, N\}$. Set $\mathbb{I}_j^N:= \{1,2, \cdots, N\}\setminus \{j\}$.  Plugging \eqref{gabdeffoc} into   \eqref{fjEigFuncFOC-}, we use \eqref{Qisdet(x-M)} to obtain that
\begin{equation}\label{xfjxcvto-yjFoc-}
x\mathfrak{f}_j(x) = \sum_{m=1}^N \frac{x y_m q_{jm}(x)}{(-1)^{N+m+j}Q(x)} = - \frac{y_j x^N}{Q(x)} - \frac{y_j x\widetilde{\mathfrak{r}}_j (x)}{Q(x)} + \sum_{m \in \mathbb{I}_j^N}\frac{y_m}{(-1)^{N+m+j}} \frac{x q_{jm} (x)}{Q(x)}  \to - y_j,
\end{equation}as $x \to +\infty$. Since $\mathfrak{f}_j = \mathfrak{g}_{\boldsymbol{\zeta}_j Y} \in \frac{\mathbb{C}_{\leq N-1}[X]}{Q}$ by \eqref{gabinDominftyG},   Lemma $\ref{LemGphi=Xphi-const}$ and
formula \eqref{KthetainDominftyG}
  imply that
\begin{equation}\label{Foc-GfjinDominfG}
  \mathbf{X} \mathfrak{f}_j  - \tfrac{i}{2 \pi} \widehat{\mathfrak{f}}_j(0^+) = \mathbf{G} \mathfrak{f}_j  \in \mathbf{G}\left( \tfrac{\mathbb{C}_{\leq N-1}[X]}{Q} \right) = \tfrac{\mathbb{C}_{\leq N-1}[X]}{Q} \subset \mathrm{Dom}^{\infty}(\mathbf{G}) \subset \mathrm{Dom}(\mathbf{G}) \cap   H^{\infty}_+.
\end{equation}So $\tfrac{i}{2 \pi} \widehat{\mathfrak{f}}_j(0^+) = \lim_{x \to +\infty}x \mathfrak{f}_j(x) = -y_j$. The proof is completed. 
\end{proof}

\noindent The next lemma shows that each
$\mathfrak{f}_j$  is an eigenfunction of the Lax operator   $\mathbf{L}_u$ with eigenvalue $\lambda_j$.

\begin{lem}\label{fjEigFuncofLuLem}
For any $j \in \{1,2, \cdots, N\}$, we have $  \mathbf{L}_u (\mathfrak{f}_j) = \lambda_j \mathfrak{f}_j$.
\end{lem}
\begin{proof}
 Recall that $\mathfrak{f}_j = \mathfrak{g}_{\boldsymbol{\zeta}_j Y} \in \frac{\mathbb{C}_{\leq N-1}[X]}{Q} \subset \mathrm{Dom}^{\infty} (\mathbf{G})$ by \eqref{gabinDominftyG} and \eqref{fjEigFuncFOC-}. Then we have 
\begin{equation}\label{DfjISTForfoc-}
D \mathfrak{f}_j(x) = -i \partial_x \mathfrak{f}_j(x) =   - i Y^* (M-x)^{-2}   \boldsymbol{\zeta}_j, \quad \forall x \in \mathbb{R}.
\end{equation}Moreover, using formulas \eqref{defXYFoclambda}, \eqref{Vdiag-},  \eqref{beta=VY=X-Y}, \eqref{gabdeffoc},  \eqref{gabinDominftyG},  \eqref{ForfocProdtoSum-} and \eqref{ulamgamalpdef}, we obtain that
\begin{equation}\label{eialpTubarfjForFOC-}
e^{i\alpha} \mathbf{T}_{\overline{u}}(\mathfrak{f}_j) =  \Pi(\mathfrak{g}_{\boldsymbol{\zeta}_j Y} + \overline{\mathfrak{g}}_{X \boldsymbol{\zeta}_j }  - \mathfrak{g}_{\boldsymbol{\zeta}_j X}) = - \mathfrak{g}_{\boldsymbol{\zeta}_j \beta} \in \tfrac{\mathbb{C}_{\leq N-1}[X]}{Q} \quad \mathrm{with} \quad \beta = X-Y =VY \in \mathbb{R}^N.
\end{equation}Lemma $\ref{ulagaFoclinUNsurj}$ yields that $u =e^{i \alpha} (1+i \mathfrak{g}_{X Y}) \in \mathcal{U}_N \subset \mathcal{Z}^{\infty}_+$. By formulas  \eqref{ZnHardyHninHnfor} and \eqref{SmoothZhidkHa}, we have
\begin{equation}\label{doubTopefjSurjfoc-}
\mathbf{T}_u   \mathbf{T}_{\overline{u}}(\mathfrak{f}_j) = u \mathbf{T}_{\overline{u}}(\mathfrak{f}_j) = - \mathfrak{g}_{\boldsymbol{\zeta}_j \beta } - i \mathfrak{g}_{X Y} \mathfrak{g}_{\boldsymbol{\zeta}_j \beta} \in \tfrac{\mathbb{C}_{\leq 2 N-1}[X]}{Q^2} \subset \mathrm{Dom}^{\infty}(\mathbf{G}).
\end{equation}Fix $x \in \mathbb{R}$.  The second formula in  \eqref{Forfoc[V,M]&M-M*} implies the following commutator formula:
\begin{equation}\label{commFor(M-xresolv&VFoc-}
\left[(M-x)^{-1},  V \right] =  i (M-x)^{-1} XY^* V (M-x)^{-1} - i (M-x)^{-2} \in \mathbb{C}^{N \times N}.
\end{equation}Recall that  $ V \boldsymbol{\zeta}_j = \lambda_j \boldsymbol{\zeta}_j \in \mathbb{R}^N$. From    \eqref{fjEigFuncFOC-},    \eqref{commFor(M-xresolv&VFoc-}, \eqref{beta=VY=X-Y} and \eqref{DfjISTForfoc-}, we have 
\begin{equation}\label{EquiForLufj=lamfjFoc}
\begin{split}
& \left( \lambda_j \mathfrak{f}_j   - \mathfrak{g}_{\boldsymbol{\zeta}_j \beta}\right) (x) = \langle \left[(M-x)^{-1},  V \right] \boldsymbol{\zeta}_j , Y \rangle_{\mathbb{C}^N} 
\\
 = &  i Y^*    (M-x)^{-1} XY^* V (M-x)^{-1}\boldsymbol{\zeta}_j - i Y^* (M-x)^{-2}   \boldsymbol{\zeta}_j = \left( i \mathfrak{g}_{X Y} \mathfrak{g}_{\boldsymbol{\zeta}_j \beta} + D \mathfrak{f}_j \right) (x).
 \end{split}
\end{equation}Plugging \eqref{doubTopefjSurjfoc-} into formula \eqref{EquiForLufj=lamfjFoc},  we obtain $\mathbf{L}_u (\mathfrak{f}_j) = \lambda_j \mathfrak{f}_j$. The proof is finished. 
\end{proof}

\noindent Consequently, the parameters $\lambda_1, \lambda_2, \cdots, \lambda_N$ are   eigenvalues of the Lax operator $\mathbf{L}_u$. Proposition $\ref{SpecDeHacHscHpp}$ yields that each eigenvalue is simple and that $\mathbf{L}_u$  has no other eigenvalues. In other words, we have
\begin{equation}\label{SpInforLuSurjFOC-}
\small 
\sigma_{\mathrm{pp}}(\mathbf{L}_u)=\{\lambda_1, \lambda_2, \cdots,  \lambda_N\}, \quad \tfrac{\mathbb{C}_{\leq N-1}[X]}{Q} = \mathscr{H}_{\mathrm{pp}}(\mathbf{L}_u) = \bigoplus_{j =1 }^N  \mathrm{Ker}(\lambda_j - \mathbf{L}_u), \quad \mathrm{Ker}(\lambda_j - \mathbf{L}_u) = \mathbb{C}\mathfrak{f}_j,
\end{equation}and $\lambda_j^u = \lambda_j$, for any   $j \in \{1,2, \cdots, N\}$. The standard eigenfunction $\phi_j^u \in \mathrm{Ker}(\lambda_j - \mathbf{L}_u) = \mathbb{C}\mathfrak{f}_j 
$ is defined in formula \eqref{phijudefFor}. Then there exists a constant $\mathfrak{c}_j \in \mathbb{S}^1$ such that
\begin{equation}\label{cjconstDefFoc-}
\mathfrak{f}_j = \mathfrak{c}_j \|\mathfrak{f}_j\|_{L^2_+} \phi_j^u = \sqrt{2 \pi}  \mathfrak{c}_j  \phi_j^u \in \tfrac{\mathbb{C}_{\leq N-1}[X]}{Q} , \quad 1 \leq j \leq N.
\end{equation}Comparing   \eqref{ulamgamalpdef} with formula   \eqref{ueiapluvDarSoli}, we have $v = u - e^{i\alpha} = i e^{i\alpha} \mathfrak{g}_{X Y} \in \frac{\mathbb{C}_{\leq N-1}[X]}{Q}$. Formula \eqref{defXYFoclambda} shows that $X= (x_1, x_2,  \cdots, x_N)^{\mathrm{T}} \in \mathbb{R}^N$ and $Y= (y_1, y_2,  \cdots, y_N)^{\mathrm{T}} \in \mathbb{R}^N$, where
\begin{equation}\label{xjuxjyjuyjidentifiyFoc-}
 \sqrt{2 \pi} x_j  =  \sqrt{\tfrac{2 \pi (1 + \lambda_j)}{\lambda_j}} = \sqrt{2 \pi}\mathbf{x}_j(u) = \mathbf{p}_j(u)   = \langle v, \phi_j^u \rangle_{L^2_+} \quad \mathrm{and} \quad y_j = \tfrac{\mathrm{sign}(\lambda_j)}{\sqrt{\lambda_j (1+ \lambda_j)}} = \mathbf{y}_j(u),
\end{equation}for any $j \in \{1,2, \cdots, N\}$, by \eqref{phijudefFor} and \eqref{xjuyjufor}.  By \eqref{gabdeffoc} and  \eqref{fjEigFuncFOC-}, it follows from 
\eqref{xjuxjyjuyjidentifiyFoc-} that  
\begin{equation}
\sum_{j=1}^N \sqrt{2 \pi}\mathbf{x}_j(u) \phi_j^u  = v = i e^{i\alpha} \mathfrak{g}_{X Y} =  \sum_{j=1}^N i e^{i\alpha} x_j \mathfrak{f}_j = \sum_{j=1}^N \sqrt{2 \pi} i e^{i\alpha} \mathfrak{c}_j x_j  \phi_j^u \in \tfrac{\mathbb{C}_{\leq N-1}[X]}{Q}.
\end{equation}So $\mathfrak{c}_j^{-1} = i e^{i\alpha} \in \mathbb{S}^1$. The standard eigenfunction $\phi_j^u \in \mathrm{Ker}(\lambda_j - \mathbf{L}_u)$ can be expressed explicitly,
\begin{equation}\label{FocexpExpphij}
\phi_j^u(x) =   \tfrac{i e^{ i \alpha}}{\sqrt{2 \pi}} \mathfrak{f}_j (x) = \tfrac{i e^{ i \alpha}}{\sqrt{2 \pi}}\langle (M-x)^{-1} \boldsymbol{\zeta}_j , Y \rangle_{\mathbb{C}^N}, \quad \forall x \in \mathbb{R}, \quad \forall j \in \{1,2, \cdots, N\}.
\end{equation} 
  Finally, it remains to verify that  $\boldsymbol{\gamma}_j(u)=\gamma_j$, for all $j \in \{1,2, \cdots, N\}$. 

\begin{lem}\label{bfGfjfjLem}
For any $j \in \{1,2, \cdots, N\}$, the following identity holds:
\begin{equation}\label{gaminprodGfjfj}
\langle \mathbf{G}\mathfrak{f}_j , \mathfrak{f}_j \rangle_{L^2_+} = M_{jj} \|\mathfrak{f}_j\|_{L^2_+}^2, \quad \mathrm{where} \quad M_{jj} = \gamma_j - \tfrac{i}{2 \lambda_j (1 + \lambda_j)} \in \mathbb{C}_-.
\end{equation}
\end{lem}

\begin{proof}
Fix $j \in \{1,2, \cdots,  N\}$. The polynomial $\mathfrak{w}_j \in \mathbb{C}[X]$ is defined by
\begin{equation}\label{polywjFoc-surj}
\mathfrak{w}_j (x) := (-1)^{N-1} x q_{jj}(x) - Q(x) - M_{jj} x^{N-1} = x \mathfrak{r}_j (x) - \mathfrak{R}(x), \quad \forall x \in \mathbb{R},
\end{equation}where $Q, \mathfrak{R} \in \mathbb{C}[X]$ are given by \eqref{Qisdet(x-M)},   $\mathfrak{r}_j , q_{jj}\in \mathbb{C}[X]$ are defined by \eqref{decompqjjrj<N-2}. Since $\deg \mathfrak{R} \leq N-2$ and $\deg \mathfrak{r}_j \leq N-3$ (When $N=2$,  we have $\mathfrak{r}_j = 0$), we have $\deg \mathfrak{w}_j \leq N-2 $. If $m \in \mathbb{I}_j^N= \{1,2, \cdots, N\} \setminus \{j\}$, formula \eqref{qkjdeffoc-} shows that $\deg q_{jm} \leq N-2$. Let $y_j \in \mathbb{R}$ be defined in \eqref{L2normfjFOC-}. Formula \eqref{tilrjdefforFOC-} shows that $\widetilde{\mathfrak{r}}_j \in \mathbb{C}[X]$ is a polynomial of degree $N-2$. We define another polynomial 
\begin{equation}\label{tjPolDefFor-}
\mathfrak{t}_j := - y_j \widetilde{\mathfrak{r}}_j + \sum_{m \in \mathbb{I}_j^N} (-1)^{m+j+N} y_m q_{jm} \in \mathbb{C}[X] , \qquad \text{so that} \qquad  \deg \mathfrak{t}_j \leq N-2.
\end{equation}Since $\frac{\mathfrak{w}_j}{Q}, \frac{\mathfrak{t}_j}{Q} \in \frac{\mathbb{C}_{\leq N-2}[X]}{Q} \subset H^{\infty}_+ \cap   L^1(\mathbb{R}; \mathbb{C})$, it follows from Lemma $\ref{LemL2+capL10mode}$ that
\begin{equation}\label{intwjoverQ=0Foc-}
\int_{\mathbb{R}}  \tfrac{\mathfrak{w}_j}{Q}   = \int_{\mathbb{R}}  \tfrac{\mathfrak{t}_j}{Q}  = 0.
\end{equation}Using formulas \eqref{absfjxsqFORfoc-} and \eqref{polywjFoc-surj}, we obtain that
\begin{equation}\label{xabsfjsqtoImFoc-}
x|\mathfrak{f}_j(x)|^2 =  -2 \operatorname{Im} \left( \tfrac{(-1)^{N-1} x q_{jj}(x) - Q(x) }{Q(x)} \right) = -2 \operatorname{Im} \left( \tfrac{M_{jj} x^{N-1} }{Q(x)} \right)  -2 \operatorname{Im} \left(   \tfrac{\mathfrak{w}_j(x) }{Q(x)} \right), \quad \forall x \in \mathbb{R}.
\end{equation}
  By formulas \eqref{xfjxcvto-yjFoc-} and \eqref{tjPolDefFor-}, for any $x \in \mathbb{R}\setminus \{0\}$, it follows that
\begin{equation}\label{fjxintermoftjFoc-}
\mathfrak{f}_j(x) = \frac{1}{x} \left(  - \frac{y_j x^N}{Q(x)} - \frac{y_j x\widetilde{\mathfrak{r}}_j (x)}{Q(x)} + \sum_{m \in \mathbb{I}_j^N}   \frac{y_m x q_{jm} (x)}{(-1)^{N+m+j} Q(x)}  \right)
 = \frac{-y_j x^{N-1} + \mathfrak{t}_j(x)}{Q(x)}.
\end{equation}Recall that formulas \eqref{defMlamgamFoc} and \eqref{L2normfjFOC-} imply the following identity:
\begin{equation}\label{ImMjj&absyjsqFoc-}
|y_j|^2 =  \left( \lambda_j (1+ \lambda_j) \right)^{-1} = -2 \operatorname{Im}M_{jj} = i M_{jj} - i \overline{M}_{jj}.
\end{equation}For any $x \in \mathbb{R}$, we use  formula \eqref{ImMjj&absyjsqFoc-} 
to obtain that
\begin{equation}\label{imMjjxN-1ovQsimplFoc-}
-2 \operatorname{Im} \left( \tfrac{M_{jj} x^{N-1} }{Q(x)} \right) - \tfrac{|y_j|^2 x^{N-1}}{\overline{Q(x)}} = \tfrac{i M_{jj} x^{N-1}}{Q(x)} + \tfrac{x^{N-1} (2 \operatorname{Im}M_{jj} - i \overline{M}_{jj})}{\overline{Q(x)}} = -2 M_{jj} \operatorname{Im} \left( \tfrac{x^{N-1}}{Q(x)} \right).
\end{equation}By \eqref{L2normfjFOC-} and  \eqref{Foc-GfjinDominfG}, we have $\mathbf{G}\mathfrak{f}_j  = \mathbf{X} \mathfrak{f}_j + y_j \in \mathrm{Dom}^{\infty}(\mathbf{G})$. As a consequence, we obtain
\begin{equation}\label{GfjfjbarxintegrandFoc-}
\begin{split}
\mathbf{G}\mathfrak{f}_j(x) \overline{\mathfrak{f}_j(x)}
 = & x|\mathfrak{f}_j(x)|^2 +y_j \overline{\mathfrak{f}_j(x)}
  = -2 \operatorname{Im} \left( \tfrac{M_{jj} x^{N-1} }{Q(x)} \right) - \tfrac{|y_j|^2 x^{N-1}}{\overline{Q(x)}} - 2 \operatorname{Im} \left(\tfrac{\mathfrak{w}_j(x)}{Q(x)} \right)  + \tfrac{y_j \overline{\mathfrak{t}_j(x)}}{\overline{Q(x)}}\\
  = & -2 M_{jj} \operatorname{Im}\left( \tfrac{x^{N-1}}{Q(x)} \right) - 2 \operatorname{Im} \left(\tfrac{\mathfrak{w}_j(x)}{Q(x)} \right)  + \tfrac{y_j \overline{\mathfrak{t}_j(x)}}{\overline{Q(x)}}, \qquad \forall x \in \mathbb{R},   
\end{split}
\end{equation}by formulas \eqref{xabsfjsqtoImFoc-},  \eqref{fjxintermoftjFoc-} and  \eqref{imMjjxN-1ovQsimplFoc-}.  Using \eqref{xN-1ImQxover|Q|2intFoc-} and \eqref{intwjoverQ=0Foc-}, we integrate \eqref{GfjfjbarxintegrandFoc-} over \(x \in \mathbb{R}\). This yields $\langle \mathbf{G}\mathfrak{f}_j , \mathfrak{f}_j \rangle_{L^2_+} =  2 M_{jj} \int_{\mathbb{R}} \tfrac{x^{N-1} \operatorname{Im}Q(x)}{|Q(x)|^2} = 2 \pi M_{jj} = M_{jj}  \|\mathfrak{f}_j \|_{L^2_+}^2$ by \eqref{L2normfjFOC-}. 
\end{proof} 
 
\noindent As a consequence, we have $\boldsymbol{\gamma}_j(u) = \|\mathfrak{f}_j\|_{L^2_+}^{-2} \operatorname{Re} \langle \mathbf{G} \mathfrak{f}_j, \mathfrak{f}_j \rangle_{L^2_+} =\operatorname{Re} M_{jj} =\gamma_j$, $1 \leq j \leq N$, by formulas \eqref{gammajudeffor}, \eqref{FocexpExpphij} and \eqref{gaminprodGfjfj}.  From \eqref{SpInforLuSurjFOC-} and 
\eqref{xjuxjyjuyjidentifiyFoc-},   we obtain the spectral coordinate vectors $(\mathfrak{X}(u), \mathfrak{Y}(u))$ and the inverse spectral matrix $\mathbf{M}(u)$ of the dark $N$-soliton $u$ constructed in \eqref{ulamgamalpdef}:
\begin{equation}
 \mathfrak{X}(u) = X \in \mathbb{R}^N, \quad  \mathfrak{Y}(u) = Y \in \mathbb{R}^N, \quad \mathrm{and} \quad  \mathbf{M}(u) = M \in \mathbb{C}^{N \times N}.
\end{equation}
  Summarizing the preceding construction and identifications, we obtain the
main bijectivity result.

\begin{prop}\label{BijbfPhiNprop}
For any $N\in\mathbb N_+$, the action-angle mapping
$\boldsymbol{\Phi}_N\colon\mathcal U_N\to
\mathbb{T} \times\Omega_N\times\mathbb{R}^{1\times N}$ is bijective.
\end{prop}

\begin{proof}
The injectivity was established in Corollary~\ref{ActAnglMapInjcor}.
The surjectivity for \(N=1\) is obtained in Lemma~\ref{lemInSpTrNis1}.
For \(N\ge2\), Lemma~\ref{ulagaFoclinUNsurj} shows that the
constructed \(u\) belongs to \(\mathcal U_N\). The total phase
identity \(\boldsymbol{\alpha}(u)=\alpha\) is immediate from the
definition of \(u\) in \eqref{ulamgamalpdef} and Definition~\ref{InvSPmADef}.  
Lemmas~\ref{fjEigFuncofLuLem} and \ref{bfGfjfjLem} imply that
$\boldsymbol{\lambda}_j(u)=\lambda_j$ and
$\boldsymbol{\gamma}_j(u)=\gamma_j$ for all $j \in \{1,2, \cdots, N\}$. 
 Thus
$\boldsymbol{\Phi}_N$ is bijective.
\end{proof}

\noindent The next lemma gives 
 an algebraic identity that will play a central role in the next section, where it will be used to identify the time evolution of dark multi-soliton solutions and to prove their global existence.
 
\begin{lem}\label{RHSCMdNLSUNlem}
The dark $N$-soliton $u \in \mathcal{U}_N$ given by \eqref{ulamgamalpdef} satisfies  the following identity
\begin{equation}\label{ilapu2uPidmodu2Foc}
\left( i \partial_x^2 u +  2 u \Pi \partial_x (|u|^2) \right) (x) = -2 i e ^{i \alpha} \langle (M -x )^{-1} V (M -x )^{-1} X , Y \rangle_{\mathbb{C}^N}
\end{equation}
\end{lem}
\begin{proof}
Formula \eqref{absusq-1=gXbeta} yields that $\Pi(|u|^2 - 1) = - i \mathfrak{g}_{X \beta} \in  \frac{\mathbb{C}_{\leq N-1}[X]}{Q}$. By  \eqref{ulamgamalpdef}, we have
\begin{equation}\label{Foc2uPiparxuabs2}
2 u \Pi \partial_x (|u|^2) = 2 u \partial_x \Pi(|u|^2 - 1) = - 2 i e^{i \alpha} \left(\partial_x \mathfrak{g}_{X \beta} + i \mathfrak{g}_{X Y} \partial_x \mathfrak{g}_{X \beta} \right) \in H^{\infty}_+.
\end{equation}Using formulas \eqref{ulamgamalpdef}, \eqref{gabdeffoc}, \eqref{Vdiag-},  \eqref{beta=VY=X-Y}  and 
\eqref{dergabForFoc-}, we obtain that
\begin{equation}\label{ipar2uinProdFoc-}
i \partial_x^2 u (x) = -2 e^{i\alpha} \langle (M-x)^{-3} X , Y  \rangle_{\mathbb{C}^N} , \quad \partial_x \mathfrak{g}_{X \beta} (x) = \langle V (M-x)^{-2} X , Y \rangle_{\mathbb{C}^N}, \quad \forall x \in \mathbb{R}.
\end{equation}By formulas \eqref{gabdeffoc}, 
\eqref{commFor(M-xresolv&VFoc-} and  \eqref{ipar2uinProdFoc-}, it follows  that
\begin{equation}\label{inproM-xresolVM-x-1XY}
\begin{split}
& \langle (M-x)^{-1} V  (M-x)^{-1} X, Y \rangle_{\mathbb{C}^N} - \partial_x \mathfrak{g}_{X \beta} (x) = \langle \left[(M-x)^{-1}, V \right] (M-x)^{-1} X, Y \rangle_{\mathbb{C}^N} \\
= & i \langle (M-x)^{-1} X, Y \rangle_{\mathbb{C}^N} \langle V (M-x)^{-2} X, Y \rangle_{\mathbb{C}^N} - i \langle (M-x)^{-3} X, Y \rangle_{\mathbb{C}^N} \\
= & \left( i \mathfrak{g}_{X Y} \partial_x \mathfrak{g}_{X \beta} - \tfrac{e^{- i \alpha}}{2} \partial_x^2 u\right) (x), \quad \forall x \in \mathbb{R}.
\end{split}
\end{equation}Then formula \eqref{ilapu2uPidmodu2Foc} follows from  \eqref{Foc2uPiparxuabs2} and \eqref{inproM-xresolVM-x-1XY}.\end{proof}

\noindent We close this subsection with the pole ansatz for dark
multi-soliton potentials; the algebraic identity obtained here will be
used in Proposition~\ref{PropComCalo-Moser} in the next section to
derive the complexified Calogero--Moser system satisfied by the pole
dynamics for the pure dark multi-soliton solutions.
  
\begin{lem}\label{LemPQhasnocomrootsFoc-}
Let \(N \ge 1\) be an integer and let \(u \in \mathcal{U}_N\) be an arbitrary dark \(N\)-soliton. Denote by \((P,Q) = (P_u, Q_u)\) the pair associated with \(u\) as in Corollary~\ref{CorQuPu} and formula~\eqref{CharPolyQuGHpp1for}. Suppose that \(Q \in \mathbb{C}[X]\) has exactly \(N\) distinct roots,   denote them by \(z_1, z_2, \cdots, z_N\in Q^{-1}\{0\} \subset \mathbb{C}_- \). Then \(P\) and \(Q\) have no common zeros, i.e., \(P(z_j) \neq 0\) for every \(j\). Consequently, in the partial fraction expansion
\begin{equation}\label{partfractExpu1Foc-}
u(x) = \frac{P(x)}{Q(x)} = u(\infty)\Bigl(1 + \sum_{j=1}^{N} \frac{a_j}{x - z_j}\Bigr), \qquad \forall  x \in \mathbb{R},
\end{equation}
the residues satisfy \(a_j \neq 0\), for all \(j \in \{1,2,\dots,N\}\). Moreover, we have
\begin{equation}\label{resiquadsystkjFOc-}
a_j + \sum_{k=1}^N \frac{a_j \overline{a}_k}{z_j - \overline{z}_k} = i, \qquad \forall  j \in \{1,2, \cdots, N\}.
\end{equation}
\end{lem} 

\begin{proof}
Recall that $\mathbb{I}_j^N := \{1,2, \cdots, N\} \setminus \{j\}$,   $1 \leq j \leq N$. Since $u = \frac{P}{Q} \in \mathcal{U}_N$, we have 
\begin{equation}\label{PQFoc-iden2}
\small
|P  |^2 - |Q |^2 = i (\overline{Q}  \partial_x Q  - Q  \partial_x \overline{Q} ) \in \mathbb{C}[X], 
\quad Q(x) = Q_u(x) = \det \left( x - \mathbf{G}|_{\mathscr{H}_{\mathrm{pp}}(\mathbf{L}_u)} \right) = \prod_{j=1}^N (x - z_j),
\end{equation}for any $x \in \mathbb{R}$, by formulas \eqref{PQiden+}  and \eqref{CharPolyQuGHpp1for}. 
Since $z_m \ne z_j \in \mathbb{C}_-$, $1 \leq m \ne j \leq N$, we have
\begin{equation}\label{ValueQ'zjFoc-}
\partial_x Q(z_j) = \prod_{m \in \mathbb{I}_j^N} (z_j - z_m) \ne 0 \quad \mathrm{and} \quad \overline{Q}(z_j) = \prod_{m=1}^N (z_j - \overline{z}_m) \ne 0 .
\end{equation}By \eqref{PQFoc-iden2} and \eqref{ValueQ'zjFoc-}, it follows that $|P(z_j)|^2 = i \overline{Q}(z_j) \partial_x Q(z_j) \ne 0$,  $\forall j \in \{1,2, \cdots, N\}$. So we have 
\begin{equation}\label{PQhasnocommrootsFoc-}
P^{-1}\{0\} \cap   Q^{-1}\{0\} = \emptyset.
\end{equation}Assume that there exists $n \in \{1,2, \cdots, N\}$ such that $a_n = 0$ in the partial fraction expansion \eqref{partfractExpu1Foc-}.    Set $\mathfrak{q}(x) : = \frac{Q(x)}{x - z_n }$, $\forall x \in \mathbb{R}$. Then there exists   $\mathfrak{p} \in \mathbb{C}[X]$ such that $\deg \mathfrak{p} = \deg \mathfrak{q} =N-1$, and 
\begin{equation}\label{pfrakpx-znFOC-}
\frac{P(x)}{Q(x)} = u (x) = u(\infty) \left(1 + \sum_{m \in \mathbb{I}_n^N} \frac{a_m}{x - z_m}\right) = \frac{\mathfrak{p}(x)}{\mathfrak{q}(x)} \;  \Longrightarrow \; P(x) =   (x - z_n) \mathfrak{p}(x), \quad \forall x \in \mathbb{R}.
\end{equation}Formula \eqref{pfrakpx-znFOC-} yields \(P(z_n)=0\), which contradicts \eqref{PQhasnocommrootsFoc-}. 
So every residue $a_j$ is different from $0$. \\

\noindent Set $\theta:= \theta_{Q_u} = \frac{\overline{Q}_u}{Q_u} \in 1 + \frac{\mathbb{C}_{\leq N-1}[X]}{Q}$ and $w(x) : = \overline{u(\infty)} u(x) -1 = \sum_{j = 1}^N \frac{a_j}{x - z_j}$, $\forall x \in \mathbb{R}$. Recall that
\begin{equation}\label{Parfracz-zjx-zkbar-}
\operatorname{Im} z_j <0 \quad \mathrm{and} \quad \tfrac{1}{(x - z_j)(x - \overline{z}_k)} = \tfrac{1}{z_j - \overline{z}_k} \left( \tfrac{1}{x - z_j}  - \tfrac{1}{x - \overline{z}_k}\right), \quad \forall x \in \mathbb{R},
\end{equation}for any $1 \leq k ,j \leq N$. By formulas \eqref{DTheta|u|2-1}, \eqref{PQFoc-iden2} and \eqref{BONsolitontheta}, it follows from \eqref{Parfracz-zjx-zkbar-} that
\begin{equation}\label{w+piwabs2=DtheovThe}
\small
\sum_{j = 1}^N \left( 1 + \sum_{k = 1}^N \frac{\overline{a}_k}{z_j - \overline{z}_k} \right) \frac{a_j}{x - z_j} = \left(w + \Pi(|w|^2) \right) (x) = \left(\Pi (|u|^2 - 1)\right) (x) = \Pi \left(\frac{D \theta}{\theta} \right) (x) = \sum_{j =1}^N \frac{i}{x -z_j},
\end{equation}for any $x\in \mathbb{R}$. Since \(z_1, z_2, \dots, z_N\) are pairwise distinct, formula \eqref{resiquadsystkjFOc-} follows from \eqref{w+piwabs2=DtheovThe}.
\end{proof}

 \section{Time Evolution}\label{SectTimeEvo6}
  
\noindent
This section is devoted to the study of pure dark multi-soliton
solutions of the focusing CMdNLS equation \eqref{CMdNLSfocu}. We
establish the time-dependent inverse spectral formula, determine the
time evolution of the inverse spectral matrix, derive the dynamics of
poles and residues, and analyze the long time asymptotics and Sobolev estimates of these
solutions.

\subsection{Pure Dark Multi-Soliton Solutions}\label{SusectDMulSol6.1}

\noindent
In this subsection we prove 
 Proposition~\ref{PropInvcalUN},    Theorem~\ref{ActionAngleMainThm}  and Proposition~\ref{InvSpecFortxintro}.   Throughout this subsection, let
$N \geq 1$ be a fixed integer. Recall from Definition~\ref{DefiDarNsoli}
that $\mathcal U_N$ denotes the set of dark $N$-soliton potentials. First, we show that $\mathcal{U}_N$ is invariant under the flow of equation \eqref{CMdNLSfocu}.

\begin{lem}\label{invmfdProLem}
 Fix $T^- < 0 < T^+$. Let $u \in \mathcal{C} \bigl( (T^-, T^+) ; \mathcal{Z}^2_+ \bigr)$ be a solution to equation \eqref{CMdNLSfocu} with initial datum $u(0) = u_0$. If $u_0 \in \mathcal{U}_N$, then $u(t) \in \mathcal{U}_N$, for all $t \in (T^-, T^+)$.
\end{lem}

\begin{proof}
Since $u_0 \in \mathcal{U}_N$, we set $\alpha : = \boldsymbol{\alpha}(u_0)\in\mathbb{T}= \mathbb{R} / 2 \pi \mathbb{Z}$ and $a :=e^{i \alpha} \in \mathbb{S}^1$. 
Then $v_0 : = u_0 -  a \in \mathrm{Dom}^{\infty}(\mathbf{G}) $ by Remark $\ref{Rkuinftyv}$ and Definition $\ref{InvSPmADef}$.
 By Lemma \ref{lemuminainCtHnx} and formula  \eqref{ttfisu-a}, the solution $u$ satisfies  
\begin{equation}\label{umineialH3}
v:= u - a =  u - e^{i \alpha} \in \mathcal{C}\bigl((T^-,T^+) ; H^3_+  \bigr) \cap \mathcal{C}^1 \bigl((T^-,T^+) ; H^1_+   \bigr) .
\end{equation}
  Let $\mathbf{W}   \in \bigcap_{p=0}^2 \mathcal{C}^1 \left( (T^-, T^+); \mathcal{B}_{\mathbb{C}}(H^p_+) \right)$ denote the unitary propagator given by 
 \eqref{UniPropagbfWSect4}. 
  Fix $t \in (T^-, T^+)$. 
 Let $\mathbf{L}_{u(t)} : H^1_+ \to L^2_+$ denote the chiral Lax operator in \eqref{LuSect2}.   
 Its self-adjointness is proved by Lemma \ref{lemLuselfadj}.  
  Using  formulas \eqref{SpeSubSpxxevot}, \eqref{lamjuEigendeffor},  Lemma \ref{LemSpLutwoseg}  and Proposition \ref{propSpecDecomLu}, we obtain that 
  \begin{equation}\label{sigHpptiEvo}
 \sigma_{\mathrm{pp}} \bigl( \mathbf{L}_{ u (t) } \bigr) = \sigma_{\mathrm{pp}} \bigl( \mathbf{L}_{ u_0 } \bigr) = \{\lambda_1^{u_0}, \lambda_2^{u_0} , \cdots , \lambda_N^{u_0} \}   \subset \mathbb{R} \setminus [-1, 0] 
 \quad \text{with} \quad 
 \lambda_1^{u_0} < \lambda_2^{u_0} < \cdots < \lambda_N^{u_0}.
\end{equation}Moreover, the pure point spectral subspace $\mathscr{H}_{\mathrm{pp}}\bigl( \mathbf{L}_{ u (t) } \bigr) $ satisfies that 
\begin{equation}\label{Hpputinfo}
\mathscr{H}_{\mathrm{pp}}\bigl( \mathbf{L}_{ u (t) } \bigr) = \mathbf{W}(t) \Bigl( \mathscr{H}_{\mathrm{pp}}\bigl( \mathbf{L}_{ u_0 } \bigr) \Bigr) 
\quad  \text{and} \quad \dim_{\mathbb{C}} \mathscr{H}_{\mathrm{pp}}\bigl( \mathbf{L}_{ u (t) } \bigr) = \dim_{\mathbb{C}} \mathscr{H}_{\mathrm{pp}}\bigl( \mathbf{L}_{ u_0} \bigr) = N .
\end{equation}  
  By formulas 
 \eqref{calHGeneFctdefFor}, 
 \eqref{CauTrSpeMeaLuisH} and Lemma \ref{LemFirsIntGenFct}, the generating functional $\mathcal{H}_{\nu}^a \bigl( u(t) \bigr)$, which is the Cauchy transform of the spectral measure $\mu^{\mathbf{L}_{u(t)}}_{ v(t)}$ is independent of the time variable $t$, for any $\nu \in \mathbb{C} \setminus \sigma \bigl(  \mathbf{L}_{u_0}  \bigr)$. By the Stieltjes--Perron inversion formula, we therefore obtain
 \begin{equation}\label{mubfLutvtinv}
 \mu^{\mathbf{L}_{u(t)}}_{ v(t)} = \mu^{\mathbf{L}_{u(0)}}_{ v (0) } = \mu^{\mathbf{L}_{u_0}}_{ v_0}, \qquad \forall  t\in(T^-,T^+) .
 \end{equation}
 Since $u_0 = e^{i \alpha}  + v_0 \in \mathcal{U}_N$, formula \eqref{vinHppandG-1Hpp} in  Theorem \ref{MainThmSpectChar} implies that $v_0 = u_0 - a \in \mathscr{H}_{\mathrm{pp}}\bigl( \mathbf{L}_{ u_0} \bigr) \cap \mathrm{Dom}(\mathbf{G})$ and $\mathbf{G} v_0 \in \mathscr{H}_{\mathrm{pp}}\bigl( \mathbf{L}_{ u_0} \bigr)$.   Then there exist $N$ numbers $\mathtt{c}_1, \mathtt{c}_2, \cdots, \mathtt{c}_N \in \mathbb{R}_+ = [0, +\infty) $ such that
\begin{equation}\label{vtinHppbfLut}
\mu^{\mathbf{L}_{u(t)}}_{ v(t)} = \mu^{\mathbf{L}_{u_0}}_{ v_0} = \sum_{j =1 }^N \mathtt{c}_j \delta_{ \lambda_j^{u_0}}, \qquad \text{so that} \qquad v(t) \in \mathscr{H}_{\mathrm{pp}}\bigl( \mathbf{L}_{ u (t) } \bigr) .
\end{equation} 
Using Theorem $\ref{MainThmSpectChar}$, 
 Proposition $\ref{propSpecDecomLu}$,  \eqref{phijudefFor},
  \eqref{sigHpptiEvo}, \eqref{Hpputinfo}   and Theorem \uppercase\expandafter{\romannumeral8}.5 (e) of \cite{ReedSimon1FA}, we obtain that 
 \begin{equation}\label{scrHppinvexpitL}
v_0 \in  \mathscr{H}_{\mathrm{pp}}\bigl( \mathbf{L}_{ u_0} \bigr) = \bigoplus_{j=1}^N \mathrm{Ker} \bigl( \lambda_j^{u_0} - \mathbf{L}_{u_0}  \bigr) = \bigoplus_{j=1}^N \mathbb{C} \phi_j^{u_0}
  = \mathbf{L}_{u_0} \Bigl(\mathscr{H}_{\mathrm{pp}}\bigl( \mathbf{L}_{ u_0} \bigr) \Bigr) =
e^{- i t \mathbf{L}_{u_0}^2 } \Bigl( \mathscr{H}_{\mathrm{pp}}\bigl( \mathbf{L}_{ u_0} \bigr) \Bigr),
 \end{equation}for all $t \in (T^-,T^+)$.   
    Since $v_0 = u_0 -  a \in \mathrm{Dom}^{\infty}(\mathbf{G}) \subset H^{\infty}_+ \cap \mathrm{Dom} (\mathbf{G})$,  
  we have
  \begin{equation}\label{vtbfWexpitL}
  v(t) = u(t) - a = \mathbf{W}(t) e^{-it\mathbf{L}_{u_0}^2}  (v_0 )   \in H^3_+ \cap \mathrm{Dom} (\mathbf{G}), 
  \end{equation}by Lemma $\ref{invSpDomGCLbfI}$ and 
  \eqref{timevoutbfWu0}.   
 It follows from  Theorem \ref{invSpDomGHtwoThm},  formulas \eqref{Hpputinfo},  \eqref{scrHppinvexpitL}  and \eqref{vtbfWexpitL} that 
 \begin{equation}\label{bfGvtinHppLut}
 \mathbf{G} \bigl( v(t) \bigr) = \mathbf{W}(t) e^{-it\mathbf{L}_{u_0}^2}  \bigl( \mathbf{G}(v_0) + 2 t \mathbf{L}_{u_0} (v_0) \bigr) \in \mathbf{W}(t) e^{-it\mathbf{L}_{u_0}^2}  \Big(  \mathscr{H}_{\mathrm{pp}}\bigl( \mathbf{L}_{ u_0} \bigr) \Bigr) \subset  \mathscr{H}_{\mathrm{pp}}\bigl( \mathbf{L}_{ u (t) } \bigr) .
 \end{equation}
   Then the assertion $(\mathbf{c})$ in Theorem $\ref{MainThmSpectChar}$ for $u(t)$ holds, according to   \eqref{sigHpptiEvo}, \eqref{Hpputinfo},  \eqref{vtinHppbfLut}, \eqref{vtbfWexpitL} and  \eqref{bfGvtinHppLut}. Thus, we have $u(t) \in \mathcal{U}_N$ by Theorem \ref{MainThmSpectChar}, for any $t \in  (T^-, T^+)$.  
 \end{proof}

\noindent Henceforth, let $u: t \in (T^-, T^+) \mapsto u(t) \in \mathcal{U}_N$ be a fixed pure dark $N$-soliton solution to the focusing CMdNLS equation \eqref{CMdNLSfocu} with initial datum $u(0) = u_0 \in \mathcal{U}_N$. 
  By formula \eqref{sigHpptiEvo}, each eigenvalue
of the Lax operator \(\mathbf L_{u(t)}\) is conserved along the flow
of \eqref{CMdNLSfocu}; more precisely,
\begin{equation}\label{lamutconseLaw}
\boldsymbol{\lambda}_j(u(t) ) = \lambda_j^{u(t)} = \lambda_j^{u_0} = \boldsymbol{\lambda}_j(u_0 ) \in \mathbb{R} \setminus [-1, 0], \qquad \forall t \in (T^-, T^+) , \qquad \forall j \in \{1,2, \cdots, N \}. 
\end{equation}
  In the next lemma, we shall describe the time evolution of 
 the normalized eigenvector $\phi_j^{u(t)} \in \mathrm{Ker} \bigl( \lambda_j^{u_0} - \mathbf{L}_{u(t)} \bigr)$ satisfying the phase
convention in \eqref{phijudefFor}.

\begin{lem}\label{LemEigVecTimevo}
   For any $j \in \{1,2, \cdots, N \}$, the mapping $\phi_j^u : t \in (T^-, T^+)   \mapsto \phi_j^{u(t)} \in H^2_+$ satisfies that 
\begin{equation}\label{phijutEvo}
\phi_j^u \in   \mathcal{C} \left( (T^-, T^+) ; H^2_+ \right) \cap     \mathcal{C}^1 \left( (T^-, T^+) ; L^2_+ \right) 
\quad \text{and} \quad  
\phi_j^{u(t)} =   \mathbf{W}(t) e^{-it\mathbf{L}_{u_0}^2}  \bigl( \phi_j^{u_0} \bigr),
\end{equation}for any $t \in (T^-, T^+)$. As a consequence, $\phi_j^u$ solves equation  \eqref{PDEdtfistilButf} with initial datum $\phi_j^{u (0)} = \phi_j^{u_0}$.
 \end{lem}

 \begin{proof}
Fix $j \in \{1,2, \cdots, N\}$ and $t \in (T^-, T^+)$.  Proposition \ref{propSpecDecomLu} and formula \eqref{phijudefFor} imply that   
 \begin{equation}\label{phijutSect6}
 \mathbb{C} \phi_j^{u(t)} = \mathrm{Ker} \bigl( \lambda_j^{u (t) } - \mathbf{L}_{u (t) }  \bigr) , \quad
\|\phi_j^{u (t) } \|_{L^2_+} = 1
\quad\text{and}\quad
  \langle \phi_j^{u (t) }, v (t)  \rangle_{L^2_+} > 0,
 \end{equation}where $v$ is given by \eqref{umineialH3} and it satisfies formula  \eqref{vtbfWexpitL}.  
 It follows from  \eqref{phijutSect6},   \eqref{bfWtisUnitaryinL2+} 
  and  \eqref{UniEquLunFor} that
 \begin{equation}\label{phijutu0}
 \lambda_j^{u_0} \phi_j^{u(t)} = \mathbf{L}_{u(t)}  \bigl( \phi_j^{u(t)} \bigr) = \mathbf{W}(t)\mathbf{L}_{u_0}  \mathbf{W}(t)^{-1} \bigl( \phi_j^{u(t)}  \bigr)  
 \quad \text{and} \quad 
\big\|  \mathbf{W}(t)^{-1} \bigl( \phi_j^{u(t)}  \bigr)   \big\|_{L^2_+} = 1 .
 \end{equation}By \eqref{phijutSect6} and  \eqref{phijutu0}, we have 
 $\mathbf{W}(t)^{-1} \bigl( \phi_j^{u(t)}  \bigr)   \in  \mathrm{Ker} \bigl( \lambda_j^{u (t) } - \mathbf{L}_{u_0 }  \bigr) =\mathbb{C} \phi_j^{u_0} $, so that $\exists  \mathsf{c}_j (t)  \in \mathbb{C}$ such that  
 \begin{equation}\label{sfcjmod1}
\mathbf{W}(t)^{-1} \bigl( \phi_j^{u(t)}  \bigr)  =  \mathsf{c}_j (t) \phi_j^{u_0} 
\quad \text{and} \quad
 |\mathsf{c}_j (t)| =1.
 \end{equation}
 Set $\mathbb{R}_+^* = (0, +\infty)$.  Using  formulas    \eqref{vtbfWexpitL}, \eqref{sfcjmod1}, \eqref{phijutSect6}  and
 Theorem \uppercase\expandafter{\romannumeral8}.5 (e) of \cite{ReedSimon1FA}, we obtain that 
 \begin{equation}\label{Valusfcj}
 e^{ i t\mathbf{L}_{u_0}^2}  \bigl( \phi_j^{u_0}  \bigr) =  e^{ i t | \lambda_j^{u_0}|^2}  \phi_j^{u_0}  ,
 \quad 
\mathsf{c}_j (t)  e^{ i t | \lambda_j^{u_0}|^2} =     \langle \phi_j^{u (t) }, v (t)  \rangle_{L^2_+} \bigl(\langle \phi_j^{u_0 }, v_0  \rangle_{L^2_+} \bigr)^{-1}  \in \mathbb{S}^1 \cap \mathbb{R}_+^* =\{1\} .
 \end{equation} 
 The second formula of \eqref{phijutEvo}  
 follows directly from \eqref{sfcjmod1} and \eqref{Valusfcj}.  
  The regularity statement in \eqref{phijutEvo} and the convergence formula \eqref{limformPDEdtf-tiButf} for $\phi_j^u$ then follow from the argument in the proof of Lemma~\ref{LemWtexp-itL2phiPDE}, in particular from formula \eqref{bffquodiffcvtilB}.
 \end{proof}

\noindent For any $j \in \{1,2, \cdots, N \}$, let $\boldsymbol{\gamma}_j : \mathcal{U}_N \to \mathbb{R}$ denote the $j$-th angle variable in  \eqref{gammajudeffor}. The inverse spectral matrix $\mathbf{M} : \mathcal{U}_N \to \mathbb{C}^{N \times N}$ is given by   
 Definition $
\ref{InvSPmADef}$ and formula \eqref{Mkjinteoflamju}. Their time evolution under the flow of equation \eqref{CMdNLSfocu} is summarized in the following corollary. 

\begin{cor}\label{TimEvoCorgamM}
For any $j \in \{1,2, \cdots, N \}$ and $t \in (T^-, T^+)$, the following formulas hold:
\begin{equation}\label{gamTimEvMal}
\boldsymbol{\gamma}_j (u(t)) = \boldsymbol{\gamma}_j (u_0) + 2 t \lambda_j^{u_0}  , 
\qquad 
\boldsymbol{\alpha}(u(t)) = \boldsymbol{\alpha}(u_0),
\qquad
 \mathbf{M} (u(t)) =  \mathbf{M} (u_0) + 2 t \mathbf{V}(u_0)  , 
\end{equation} 
  where $\boldsymbol{\alpha} : \mathcal{U}_N \to \mathbb{T}$ denotes the total phase function given by 
Definition~\ref{InvSPmADef}, and $\mathbf V(u_0)\in\mathrm{GL}_N(\mathbb{R})$
is the velocity matrix from~\eqref{VeloMatIntro}.  
 Consequently, $u$ extends to a
global-in-time solution of \eqref{CMdNLSfocu}. 
\end{cor}

\begin{proof}
Let $v = u - e^{i \boldsymbol{\alpha}(u_0) } \in \mathcal{C}\bigl((T^-,T^+) ; H^3_+  \bigr) \cap \mathcal{C}^1 \bigl((T^-,T^+) ; H^1_+   \bigr) $ be defined in \eqref{umineialH3}. Then 
\begin{equation}\label{bdalphConser}
\exp \bigl(i \boldsymbol{\alpha}(u (t)) \bigr)  = \lim_{|x | \to +\infty} u(t, x) = e^{i \boldsymbol{\alpha}(u_0) }  +  \lim_{|x | \to +\infty} v(t,x)  = \exp \bigl( i \boldsymbol{\alpha}(u_0) \bigr) ,  \quad \forall t \in (T^-,T^+) .
\end{equation} 
 Lemma \ref{expitLusqStrContGroup} and formula \eqref{bfWtisUnitaryinL2+} imply that $\mathbf{W}(t)e^{-it \mathbf{L}_{u_0}^2}$ is a unitary operator on $L^2_+$.  Using  
\eqref{gammajudeffor},    \eqref{GWtUtisWtUtGadd2tL},   
\eqref{phijutEvo} and \eqref{phijutSect6}, we obtain  that  $\boldsymbol{\gamma}_j (u(t))
 = \operatorname{Re} \langle \mathbf{W}(t)e^{-it \mathbf{L}_{u_0}^2} \bigl( \mathbf{G} + 2 t \mathbf{L}_{u_0} \bigr) \phi_j^{u_0}  ,    \mathbf{W}(t) e^{-it\mathbf{L}_{u_0}^2}  \bigl( \phi_j^{u_0} \bigr) \rangle_{L^2_+} =  \boldsymbol{\gamma}_j (u_0)   + 2 t \lambda_j^{u_0}$. Consequently,  
 it follows from  formulas  \eqref{MkjForInvSpMatr},  \eqref{Mkjinteoflamju} and  \eqref{lamutconseLaw} that  
\begin{equation}\label{MkjtimEvou}
\mathbf{M}_{jj}(u(t)) = \mathbf{M}_{jj}(u_0)  + 2  t \lambda_j^{u_0} , \qquad \mathbf{M}_{kj}(u(t)) = \mathbf{M}_{kj}(u_0) , \qquad 1 \leq k \ne j \leq N. 
\end{equation}Then formula \eqref{TimEvoCorgamM} follows directly from   \eqref{bdalphConser} and \eqref{MkjtimEvou}.
  Let $\mathfrak{X}(u_0), \mathfrak{Y}(u_0) \in \mathbb{R}^N$ be its spectral coordinate vectors in Definition $\ref{SpCoorVecDefiXY}$.
   By Definition \ref{bfMVfrakXYdefintro} and   Lemma \ref{SpofMLemfoc}, we have
 \begin{equation}\label{bfMutspinCneg}
 \sigma \bigl( \mathbf{M} (u_0) + 2 t \mathbf{V}(u_0)  \bigr) 
   \subset \mathbb{C}_-, \qquad \forall t \in \mathbb{R}
 \end{equation} 
 For any $(t,x ) \in \mathbb{R} \times \mathbb{R}$, we define that
\begin{equation}\label{ExtInSpForuext}
\mathsf{u}_{\mathrm{ext}}(t , x) := e^{i \boldsymbol{\alpha}(u_0) } \Bigl( 1 + i \, \bigl\langle \bigl( \mathbf{M} (u_0) + 2 t \mathbf{V}(u_0)  - x  \bigr)^{-1} \mathfrak{X}(u_0),\, \mathfrak{Y}(u_0) \bigr\rangle_{\mathbb{C}^N} \Bigr) .
\end{equation} 
Plugging \eqref{ExtInSpForuext} into formula \eqref{ilapu2uPidmodu2Foc} of 
 Lemma \ref{RHSCMdNLSUNlem}, we obtain that  
 \begin{equation}\label{uextsolCMdNLS}
 \begin{split}
& \left( i \partial_x^2 \mathsf{u}_{\mathrm{ext}} +  2 \mathsf{u}_{\mathrm{ext}} \Pi \partial_x (|\mathsf{u}_{\mathrm{ext}}|^2) \right) (t, x) \\
= &  -2 i  e^{i \boldsymbol{\alpha}(u_0) }  \langle \bigl( \mathbf{M} (u_0) + 2 t \mathbf{V}(u_0)  - x  \bigr)^{-1} \mathbf{V}(u_0) \bigl( \mathbf{M} (u_0) + 2 t \mathbf{V}(u_0)  - x  \bigr)^{-1} \mathfrak{X}(u_0),\, \mathfrak{Y}(u_0) \rangle_{\mathbb{C}^N} \\
= &\partial_t \mathsf{u}_{\mathrm{ext}} (t, x), \qquad \forall (t, x ) \in \mathbb{R} \times \mathbb{R}. 
\end{split}
\end{equation}
  Plugging \eqref{ExtInSpForuext} into formula  \eqref{ulamgamalpdef}, we obtain that $\mathsf{u}_{\mathrm{ext}} (t) \in \mathcal{U}_N$ by 
Lemma \ref{ulagaFoclinUNsurj}, for any $t \in \mathbb{R}$.    
 Then $\mathsf{u}_{\mathrm{ext}} : t \in \mathbb{R} \mapsto \mathsf{u}_{\mathrm{ext}} (t) \in \mathcal{U}_N$
 is a global solution of equation \eqref{CMdNLSfocu} with initial datum $\mathsf{u}_{\mathrm{ext}}(0) = u_0 \in \mathcal{U}_N$.  
 By the uniqueness part of Proposition~\ref{LWPpropCMdNLSZn} and the inverse spectral formula \eqref{invSPEFormain} for $u(t) \in \mathcal{U}_N$,
the solution $\mathsf{u}_{\mathrm{ext}}$ coincides with $u$ on $(T^-,T^+)$. 
   Therefore  $u$ admits a global-in-time extension given by  $\mathsf{u}_{\mathrm{ext}}$, which proves the last assertion of the corollary.
 \end{proof}

\noindent   Proposition~\ref{PropInvcalUN} follows from
Lemma~\ref{invmfdProLem}, while
Proposition~\ref{InvSpecFortxintro} is a direct consequence of 
Corollary~\ref{TimEvoCorgamM},  formulas \eqref{bfMutspinCneg} and  \eqref{ExtInSpForuext}.
  We now complete the proof of Theorem~\ref{ActionAngleMainThm}.

\begin{proof}[End of the proof of Theorem~\ref{ActionAngleMainThm}]
The  action-angle mapping $\boldsymbol{\Phi}_N\colon\mathcal U_N\to
\mathbb{T}\times \Omega_N\times\mathbb{R}^{1\times N}$ is constructed in Definition~\ref{InvSPmADef}. Its bijectivity is given by Proposition $\ref{BijbfPhiNprop}$.  
  Formula  \eqref{ActAnginvbyPhRot} follows directly from 
Lemma~\ref{LemPhaseRotInv}.  
  The conservation of the actions is given by
\eqref{lamutconseLaw}, while the evolution of the angles and the
total phase are the first two identities in \eqref{gamTimEvMal} of 
Corollary~\ref{TimEvoCorgamM}.
\end{proof}

\subsection{Dark Two-Soliton Solutions}\label{ssDar2SolSol6.2}

\noindent
In this subsection we specialize the preceding analysis to the case
$N=2$, that is, to pure dark \(2\)-soliton solutions of the
focusing CMdNLS equation~\eqref{CMdNLSfocu}. Fix
$u_0\in\mathcal U_2 $. We use the following notation:
\begin{equation}\label{abbracantpha}
\alpha := \boldsymbol{\alpha}(u_0),
\qquad
\lambda_j := \boldsymbol{\lambda}_j(u_0),
\qquad
\gamma_j := \boldsymbol{\gamma}_j(u_0),
\qquad \forall 
j \in \{ 1 , 2 \}.
\end{equation}The inverse spectral matrix $\mathbf{M}(u_0) \in \mathbb{C}^{2 \times 2}$ and the velocity matrix 
$\mathbf{V}(u_0) \in \mathrm{GL}_2(\mathbb{R}) $ are denoted by
\begin{equation}\label{MbfMu0Vbf}
M : = \mathbf{M}(u_0) = \left( \begin{matrix}
\gamma_1 - \frac{i}{2 \lambda_1 (1+ \lambda_1)} & \frac{i}{\lambda_1 - \lambda_2} \sqrt{\frac{\lambda_2 (1+ \lambda_1)}{\lambda_1 (1+ \lambda_2)}} \\
\frac{i}{\lambda_2 - \lambda_1} \sqrt{\frac{\lambda_1 (1+ \lambda_2)}{\lambda_2 (1+ \lambda_1)}} & 
\gamma_2 - \frac{i}{2 \lambda_2 (1+ \lambda_2)}
\end{matrix} \right), \qquad 
V : = \mathbf{V}(u_0) = \left( \begin{matrix}
  \lambda_1   &  0 \\
0  & 
 \lambda_2 
\end{matrix} \right),
\end{equation}Let $X = \mathfrak{X}(u_0), Y = \mathfrak{Y}(u_0) \in \mathbb{R}^2$ denote the spectral coordinate vectors in \eqref{introxjuyjufor}, i.e.
\begin{equation}\label{XfrakXY}
X^{\mathrm{T}} : = \left(\sqrt{\tfrac{1+ \lambda_1}{\lambda_1}} \; , \; \sqrt{\tfrac{1+ \lambda_2}{\lambda_2}} \right), \qquad
Y^{\mathrm{T}} : = \left(  \tfrac{\mathrm{sign}(\lambda_1)}{\sqrt{\lambda_1 (1+ \lambda_1)}} \; , \; \tfrac{\mathrm{sign}(\lambda_2)}{\sqrt{\lambda_2 (1+ \lambda_2)}} \right) .
\end{equation}
 By Proposition~\ref{LWPpropCMdNLSZn}  and  Proposition~\ref{PropInvcalUN}, there exists a unique global solution
$u : t\in\mathbb R \mapsto u(t)\in \mathcal{U}_2$  of equation 
\eqref{CMdNLSfocu} with initial datum $u(0) = u_0 \in \mathcal{U}_2$.  
  The inverse spectral formula \eqref{invSPEFortxIntro} implies that  
\begin{equation}\label{InvSpFottw}
 i \bigl( 1 - e^{ - i \alpha} u(t,x) \bigr)  = \langle ( M + 2 t V - x )^{-1} X , Y \rangle_{\mathbb{C}^2} , \qquad  \forall (t, x) \in \mathbb{R} \times \mathbb{R}.
\end{equation}
 For any $t\in\mathbb{R} $, the denominator
polynomial $\mathcal{Q}_t = Q_{u(t)}\in\mathbb{C}[X]$ from
\eqref{CharPolyQuGHpp1for} of Corollary~\ref{CorQuPu}  is given by
\begin{equation}\label{calQtdef}
\begin{split}
\mathcal Q_t(x)
: =  
Q_{u(t)}(x)
= & 
\det\bigl(
x-\mathbf G|_{\mathscr H_{\mathrm{pp}}(\mathbf L_{u(t)})}
\bigr)
=
\det\bigl(
x - M - 2tV
\bigr),
\qquad
\forall x\in\mathbb R, \\
= & x^2 - x \operatorname{tr} (M + 2 t V ) + \det (M + 2 t V ) , 
\end{split}
\end{equation}
according to  \eqref{MkjForInvSpMatr} and
Corollary~\ref{TimEvoCorgamM}. 
 By \eqref{bfMutspinCneg}, the spectrum of $M+2tV \in \mathbb{C}^{2 \times 2}$ lies in
$\mathbb C_-$ for all $t\in\mathbb R$, hence both roots of
$\mathcal Q_t$ lie in $\mathbb C_-$.
   For the pole ansatz to be applicable, we need to know when the denominator \(\mathcal{Q}_t\) has two distinct roots. The next lemma provides a complete answer.

\begin{lem}\label{CriCollidLem}
If the eigenvalues $\lambda_1, \lambda_2$ from   \eqref{abbracantpha} satisfy the following algebraic identity:
\begin{equation}\label{AlIdlamCollid}
(\lambda_2 - \lambda_1)^2 |1 + \lambda_1 + \lambda_2| 
= 4 \lambda_1 (1+ \lambda_1)\lambda_2 (1+ \lambda_2) ,
\end{equation}then there exists $\widetilde{\mathsf{t}} = \widetilde{\mathsf{t}} (\lambda_1, \lambda_2 ;  \gamma_1, \gamma_2) : = \frac{\gamma_1 - \gamma_2}{2 (\lambda_2 - \lambda_1)} \in \mathbb{R}$ such that $\mathcal Q_t$ in \eqref{calQtdef}  has two distinct roots if and only if $t \ne \widetilde{\mathsf{t}}$. If $\lambda_1, \lambda_2$  do  not satisfy  \eqref{AlIdlamCollid}, then    $\mathcal Q_t$ has two distinct roots, for any $t \in \mathbb{R}$. 
\end{lem}

\begin{proof}
Let $\vec{\lambda}: = (\lambda_1, \lambda_2) \in \Omega_2$ and $\vec{\gamma}: = (\gamma_1, \gamma_2) \in  \mathbb{R}^{1 \times 2}$ be given by formula  \eqref{abbracantpha}.  
  We define two constants $\varrho = \varrho (\vec{\lambda})>0$, $\daleth = \daleth (\vec{\lambda} , \vec{\gamma}) \in \mathbb{C}$ and a function $\varsigma : \mathbb{R} \setminus \{0\} \to \mathbb{C}$ as follows: 
\begin{equation}\label{vrhodlth}
\varrho: = \lambda_2 - \lambda_1 ,
\qquad
\daleth : = \tfrac{\gamma_2 - \gamma_1}{2 (\lambda_2 - \lambda_1)} + \tfrac{i (1 + \lambda_1 + \lambda_2)}{4 \lambda_1 (1 + \lambda_1)\lambda_2 (1 +  \lambda_2)}  ,
\qquad
\varsigma (t) : = 2 \daleth t^{-1} + (\daleth^2 + \varrho^{-4}) t^{-2} ,
\end{equation}for any $t \in \mathbb{R} \setminus \{0\}$. Let $\Delta_{u_0}(t)$ denote the discriminant of the quadratic polynomial $\mathcal{Q}_t \in \mathbb{C}[X]$ given by \eqref{calQtdef}. Using formulas  \eqref{abbracantpha},  \eqref{MbfMu0Vbf} and \eqref{vrhodlth}, we obtain that
\begin{equation}\label{discrcalQt}
\Delta_{u_0}(t) : = \bigl( \mathrm{tr}  ( M + 2 t V ) \bigr)^2 - 4 \det ( M + 2 t V )  = 4 \varrho^2 \bigl( (t + \daleth )^2 +  \varrho^{-4} \bigr) = 4 \varrho^2 t^2 \bigl( 1 + \varsigma (t) \bigr).
\end{equation}The discriminant $\Delta_{u_0}(t)$ vanishes if and only if the following identities hold:
\begin{equation}\label{ReImdaleDis}
t = - \operatorname{Re} \daleth = \widetilde{\mathsf{t}}, \qquad |\operatorname{Im} \daleth | = \varrho^{-2} .
\end{equation}By \eqref{vrhodlth}, the second formula of \eqref{ReImdaleDis} is equivalent to formula \eqref{AlIdlamCollid}. The proof is completed.
\end{proof}

\noindent As a consequence, the two poles of any dark two-soliton solution
collide at most once. \\
 
\noindent
The principal branch of \(\sqrt{1+z}\) admits the power series
expansion about \(z=0\): 
\begin{equation}\label{sq1addzPSer}
\kappa (z) : = 
\sqrt{1+z}
= 1 + \tfrac{1}{2} z - \tfrac{1}{8} z^2
+ \sum_{n\geq 3} \mathsf{k}_n z^n  
 \quad \text{with} \quad \mathsf{k}_n
= \tfrac{(-1)^{n-1}(2n-3)!!}{(2n)!!},
\quad
|z|<1.
\end{equation}
Its  radius of convergence is  \(1\).   
  Let $\mathcal{D}(0, 1):=\{z \in \mathbb{C}: |z | <1\}$ denote open unit disc.  
  For any $ z \in  \mathcal{D}(0, 1)$, we set   $\mathscr{R}_3(z) : = z^{-3} \bigl( \kappa (z) 
- 1 - \tfrac{1}{2} z +  \tfrac{1}{8} z^2  \bigr)$. Since $|\mathsf{k}_n| < \frac{1}{2n}$ for any $n \geq 3$, we have  
\begin{equation}\label{EstscrR3}
\mathscr{R}_3 (z)  =\sum_{n\geq 3} \mathsf{k}_n z^{ n -3}  
\qquad \text{and} \qquad
\sup \{|\mathscr{R}_3 (z) | : |z| \leq \tfrac{1}{2}\} \leq \tfrac{1}{3} . 
\end{equation}
 Expanding \(\kappa(\varsigma(t))\) in powers of \(\varsigma(t)\) and
collecting the terms of order \(t^{-1}\) and lower, if $t \ne 0$,  we set 
\begin{equation}\label{sfrneg1tdef}
\mathsf{r}_{-1} (t) : = \varrho^{-3} + \Bigl( 2 \varrho \bigl( t \varsigma (t) \bigr)^3 \mathscr{R}_3 \bigl( \varsigma (t) \bigr) -  \varrho  \daleth \bigl( \daleth^2 + \varrho^{-4} \bigr) \Bigr) t^{-1}  - \tfrac{1}{4} \varrho \bigl( \daleth^2 + \varrho^{-4} \bigr)^2 t^{-2}  .
\end{equation}
  By   \eqref{sq1addzPSer} and \eqref{discrcalQt}, we have $\Delta_{u_0}(t) = \bigl( 2 \varrho  t  \kappa (\varsigma (t)) \bigr)^2 $.  
  Its principal square root is given by  
  \begin{equation}\label{vrtheu0t}
\vartheta_{u_0}(t) : 
= 2 \varrho  t  \kappa (\varsigma (t)) 
= 2 \varrho  t + 2 \varrho   \daleth + \mathsf{r}_{-1} (t) t^{-1} , \qquad \forall t \in \mathbb{R} \setminus \{0\} , 
  \end{equation}where the second identity follows from \eqref{sq1addzPSer}, \eqref{EstscrR3} and  \eqref{sfrneg1tdef}. 
    The third formula of \eqref{vrhodlth} yields that  
   $\lim_{|t| \to + \infty} t \varsigma (t)  = 2 \daleth$.   
   There exists $\mathscr{T}_1  = \mathscr{T}_1 \bigl( \vec{\lambda} ; \vec{\gamma} \bigr) > 1 + \big| \frac{\gamma_1 - \gamma_2}{2 (\lambda_2 - \lambda_1)}\big|$ such that 
 \begin{equation}\label{vrsigmEsti}
\sup_{|t| \geq \mathscr{T}_1} |\varsigma (t) | \leq \tfrac{1}{2} 
\qquad \text{and}  \qquad
\sup_{|t| \geq \mathscr{T}_1} |t \varsigma (t) | \leq 3 |\daleth| .
\end{equation}
 By  \eqref{EstscrR3}, \eqref{sfrneg1tdef}  and \eqref{vrsigmEsti}, the following estimate holds: 
\begin{equation}\label{UnEstsfrng1}
\sup_{|t| \geq \mathscr{T}_1} | \mathsf{r}_{-1} (t) | \leq 
\varrho^{-3} + 18 \varrho |\daleth|^3 + \varrho |\daleth| \big| \daleth^2 + \varrho^{-4} \big| + \tfrac{1}{4}  \varrho  \big| \daleth^2 + \varrho^{-4} \big|^2  .
\end{equation}Let $\mathbf{z}_1(t), \mathbf{z}_2(t) \in \mathbb{C}_-$ denote the two roots of the denominator polynomial \(\mathcal{Q}_t\). To be precise, we define
 \begin{equation}\label{bfzjTwoSol}
\mathbf{z}_j (t) :  = \tfrac{1}{2} \bigl( \operatorname{tr} (M + 2 t V)  + (-1 )^j \vartheta_{u_0}(t)  \bigr) = 2 \lambda_j t + \gamma_j - \tfrac{i}{2 \lambda_j (1 + \lambda_j)} + \tfrac{(-1)^j}{2 t} \mathsf{r}_{-1} (t) ,
\end{equation}for any $j \in \{1, 2\}$ and $t \in \mathbb{R} \setminus \{0\}$.  
  With this convention, \(\mathbf z_1(t)\) is asymptotic to
\(2\lambda_1 t\) and \(\mathbf z_2(t)\) to \(2\lambda_2 t\) as
\(|t|\to + \infty\).   
  For any $j \in \{1, 2\}$, 
  using  \eqref{UnEstsfrng1} and 
\eqref{bfzjTwoSol}, we obtain that 
\begin{equation}\label{ImbfzjLim2Sol}
\operatorname{Im} \mathbf{z}_j (t) = - \tfrac{1}{2 \lambda_j (1 + \lambda_j)} + \tfrac{(-1)^j}{2 t} \operatorname{Im} \mathsf{r}_{-1} (t) \to - \tfrac{1}{2 \lambda_j (1 + \lambda_j)}, 
\quad \text{as} \quad |t| \to +\infty .
\end{equation}
   For any $t \in \mathbb{R} \setminus [-\mathscr{T}_1 \bigl( \vec{\lambda} ; \vec{\gamma} \bigr), \mathscr{T}_1 \bigl( \vec{\lambda} ; \vec{\gamma} \bigr)]$, 
Lemma~\ref{CriCollidLem} shows that $\mathbf{z}_1(t) \ne \mathbf{z}_2 (t)$. The dark $2$-soliton potential $u (t) \in \mathcal{U}_2$ admits the following partial fraction expansion: there exist  $\mathbf{a}_1(t), \mathbf{a}_2(t) \in \mathbb{C}$ such that
\begin{equation}\label{PolAnsa2solt}
u(t,x) = e^{ i \alpha} \bigl( 1 + i \langle ( M + 2 t V - x )^{-1} X , Y \rangle_{\mathbb{C}^2} \bigr)
= e^{ i \alpha} \Bigl( 
1 + \tfrac{\mathbf{a}_1(t)}{x - \mathbf{z}_1(t)}  + \tfrac{\mathbf{a}_2(t)}{x - \mathbf{z}_2(t)} \Bigr) .
\end{equation}
  Plugging \eqref{MbfMu0Vbf} and \eqref{XfrakXY} into \eqref{PolAnsa2solt}, we can express $\mathbf{a}_1(t), \mathbf{a}_2(t) \in \mathbb{C}$ in terms of  $\lambda_1, \lambda_2 \in \mathbb{R} \setminus [-1, 0]$ and $ \vartheta_{u_0}(t), \mathsf{r}_{-1} (t) \in \mathbb{C}$ by solving a $2 \times 2$ linear system. For  any $j \in \{1,2 \}$,  we have
\begin{equation}\label{bfajExpFor2Sol}
\mathbf{a}_j(t) =  - \tfrac{i}{\lambda_j} - \tfrac{(-1)^j}{\lambda_1 \lambda_2 \vartheta_{u_0}(t)} + \tfrac{i (-1)^j (\lambda_1 - \lambda_2) \mathsf{r}_{-1} (t)  }{ 2 \lambda_1  \lambda_2 t \vartheta_{u_0}(t)} \to - \tfrac{i}{\lambda_j}, \quad \text{as} \quad |t| \to +\infty.
\end{equation}
  By \eqref{ImbfzjLim2Sol} and \eqref{bfajExpFor2Sol}, there exists a constant  $\mathscr{T}_2  = \mathscr{T}_2 \bigl( \vec{\lambda} ; \vec{\gamma} \bigr)  \geq   \mathscr{T}_1 \bigl( \vec{\lambda} ; \vec{\gamma} \bigr) > 1 + \big| \frac{\gamma_1 - \gamma_2}{2 (\lambda_2 - \lambda_1)}\big|$ such that
\begin{equation}\label{UniEsttajzj2Sol}
\tfrac{1}{2 |\lambda_j|}\leq  |\mathbf{a}_j(t)|  \leq  \tfrac{2}{  |\lambda_j|} 
\quad \text{and} \quad
\tfrac{1}{4  \lambda_j (1 + \lambda_j)} \leq |\operatorname{Im} \mathbf{z}_j (t) | \leq \tfrac{1}{ \lambda_j (1 + \lambda_j)} , \quad \forall t \in \mathbb{R} \setminus [-\mathscr{T}_2 , \mathscr{T}_2 ], 
\end{equation}for any $j \in \{1, 2\}$.  
  We set 
$\mathtt{K}_*= \mathtt{K}_* \bigl( \vec{\lambda} \bigr) : =1 + \max \big\lbrace \frac{2}{  |\lambda_j|} , 2 |\lambda_j|,  \frac{1}{ \lambda_j (1 + \lambda_j)}, 4 \lambda_j (1 + \lambda_j) : j =1, 2 \big\rbrace$ and
\begin{equation}\label{bdvrphij2Sol}
\boldsymbol{\varphi}_j (t,x) : = (x - \mathbf{z}_j (t) )^{-1} , \quad \forall x \in \mathbb{R}, \quad \forall t \in \mathbb{R} \setminus [-\mathscr{T}_2 , \mathscr{T}_2 ] , \quad \forall  j \in \{1, 2\}.
\end{equation}
  For any $s \geq 0$, 
Lemma~\ref{BdvphiLongtimeLem} yields that there exists a constant $\mathrm{B} \bigl( s, \mathtt{K}_* ( \vec{\lambda} ) \bigr) > 1$ such that 
\begin{equation}\label{LoUpBdbdphij2Sol}
\mathrm{B} \bigl( s, \mathtt{K}_* ( \vec{\lambda} ) \bigr)^{-1} \leq \|\boldsymbol{\varphi}_j  (t)\|_{H^s_+} \leq \mathrm{B} \bigl( s, \mathtt{K}_* ( \vec{\lambda} ) \bigr) , \quad \forall   t \in \mathbb{R} \setminus [-\mathscr{T}_2 , \mathscr{T}_2 ] , \quad \forall j \in \{1,2 \}.
\end{equation} 
  Recall that $u_0 (\infty) : = \lim_{|x| \to +\infty} u_0 (x) = e^{i \alpha} = e^{i \boldsymbol{\alpha}(u_0)} \in \mathbb{S}^1$ by Definition~\ref{InvSPmADef}. 
  Using 
\eqref{PolAnsa2solt},  \eqref{UniEsttajzj2Sol}, \eqref{bdvrphij2Sol} and  \eqref{LoUpBdbdphij2Sol}, we obtain the following long time Sobolev estimate: 
\begin{equation}\label{LonTSob2Soli}
\sup_{|t| \geq \mathscr{T}_2}\|u(t) - u_0 (\infty)  \|_{H^s_+} \leq 2 \mathtt{K}_* ( \vec{\lambda} ) \mathrm{B} \bigl( s, \mathtt{K}_* ( \vec{\lambda} ) \bigr) , \qquad \forall s \geq 0.
\end{equation}
  As a consequence, the inverse spectral formula
\eqref{invSPEFortxIntro} established in
Proposition~\ref{InvSpecFortxintro} enables us to determine the
explicit long time dynamics of both poles and residues, and hence to
obtain the Sobolev estimate~\eqref{LonTSob2Soli} for every pure dark
\(2\)-soliton solution of~\eqref{CMdNLSfocu}. 
  The remaining subsections of this section are devoted to extending
the formulas \eqref{bfzjTwoSol}, \eqref{bfajExpFor2Sol}, and
\eqref{LonTSob2Soli} to arbitrary \(N\), thus covering all dark
multi-soliton solutions of~\eqref{CMdNLSfocu}. For \(N\geq 3\), the
denominator polynomial \(\mathcal Q_t = Q_{u(t)}\) has degree \(N\).
  It is difficult to express all its roots by elementary algebraic formulas.
  Instead, we analyze the spectrum of the time-dependent matrix
\(M+2tV\) by means of Kato's matrix perturbation theory and the
Rayleigh--Schr\"odinger expansion. This yields precise long time
asymptotics for each pole \(\mathbf z_j(t)\) and each residue
\(\mathbf a_j(t)\) in the partial fraction expansion:  
\begin{equation}\label{Ed2ParFracExpu}
u(t,x) = e^{i \boldsymbol{\alpha}(u_0)} \Bigl( 1 + \sum_{j=1}^N \frac{\mathbf{a}_j(t)}{x - \mathbf{z}_j (t)} \Bigr), \quad \forall  x \in  \mathbb{R}, \quad |t| \gg 1,
\end{equation} 
  and ultimately leads to the corresponding
long time Sobolev estimates.

\subsection{Matrix Perturbation}\label{SubseMatPert6.3}

\noindent Fix $N \in \mathbb{N}_+$ and let $u : t \in \mathbb{R} \mapsto u(t) \in \mathcal{U}_N$ be a pure dark $N$-soliton solution to equation \eqref{CMdNLSfocu}. 
The purpose of this subsection is to analyze the long time behavior of
the eigenvalues of the inverse spectral matrix
$\mathbf{M}(u(t))
= \mathbf{M}(u_0) + 2t\,\mathbf{V} (u_0) \in \mathbb{C}^{N \times N}$ given by Definition~\ref{bfMVfrakXYdefintro}
 and formula \eqref{gamTimEvMal}. 
    We shall apply Kato's  
perturbation theory to construct the eigenvalues
\(\mathbf{z}_1(t), \mathbf{z}_2(t), \cdots, \mathbf{z}_N(t)\) as real analytic functions
of \(t\) for  large \(|t|\), and we will compute their asymptotic
expansions via the Rayleigh--Schr\"odinger series.
  For notational simplicity, throughout this subsection we write
$ M := \mathbf{M}(u_0)$,  
$V := \mathbf{V}(u_0)$,
and we identify the time-dependent inverse spectral matrix
\(\mathbf{M}(u(t))\) with \(M+2tV \in \mathbb{C}^{N \times N}\).
  
 \begin{lem}\label{RaySchrCoenegLem}
 For any  $N \in \mathbb{N}_+$,  $\vec{\lambda} \in \Omega_N$ and $\vec{\gamma} \in \mathbb{R}^{1 \times N}$, let $V=V_{\vec{\lambda}} \in \mathrm{GL}_{N}(\mathbb{R})$ and $M=M(\vec{\lambda}, \vec{\gamma}) \in \mathbb{C}^{N \times N}$ be defined in \eqref{Vdiag-} and \eqref{defMlamgamFoc}. Then there exists a constant $\mathcal{T}:=\mathcal{T}(\vec{\lambda}; \vec{\gamma})>1$ and $N$ real analytic functions $\mathbf{z}_1, \mathbf{z}_2, \cdots, \mathbf{z}_N \in \mathcal{C}^{\omega}\left(\mathbb{R}\setminus [-\mathcal{T}, \mathcal{T}]; \mathbb{C} \right)$ satisfying the following:\\
\noindent $(\mathbf{a}).$ For any $t \in \mathbb{R}\setminus [-\mathcal{T}, \mathcal{T}]$, if $1 \leq k \ne j \leq N$, then $\mathbf{z}_k(t) \ne \mathbf{z}_j(t)$. \\
\noindent $(\mathbf{b}).$  For any $t \in \mathbb{R}\setminus [-\mathcal{T}, \mathcal{T}]$, the spectrum of  $M+2tV \in \mathbb{C}^{N \times N}$ is given by
\begin{equation}\label{SpofM+2tV}
\sigma(M + 2tV)=\{\mathbf{z}_1(t), \mathbf{z}_2 (t), \cdots, \mathbf{z}_N (t)\} \subset \mathbb{C}_-.
\end{equation}Moreover, each eigenvalue is simple, i.e. $\mathrm{dim}_{\mathbb{C}} \mathrm{Ker}\left( M + 2tV -\mathbf{z}_j(t)\right)=1$, $1\leq j \leq N$.\\
$(\mathbf{c}).$ The following inequality holds,
\begin{equation}\label{tzjttpartzjtEstiRaySch}
\max_{1\leq j \leq N}\sup_{|t|>\mathcal{T}}\left(|t(\mathbf{z}_j(t) - 2\lambda_j t - M_{jj})| + |t(\partial_t \mathbf{z}_j(t) - 2 \lambda_j)| \right) <+\infty.
\end{equation}
\end{lem}

\begin{rem}
Formula \eqref{tzjttpartzjtEstiRaySch}    shows the long time behavior of  $\mathbf{z}_j \in \mathcal{C}^{\omega}\left(\mathbb{R}\setminus [-\mathcal{T}, \mathcal{T}]; \mathbb{C}  \right)$:
\begin{equation}\label{Longtimezj-}
\mathbf{z}_j(t) = 2 \lambda_j t + M_{jj} + \mathcal{O}(t^{-1}), \quad \partial_t \mathbf{z}_j(t) = 2 \lambda_j    + \mathcal{O}(t^{-1}), \quad |t| > \mathcal{T}, \quad 1 \leq j \leq N.
\end{equation}
\end{rem}

\begin{proof}[Proof of Lemma $\ref{RaySchrCoenegLem}$]
 Let $\mathcal{D}(a, r):=  \mathsf{B}_{\mathbb{C}} (0, r)= \{z \in \mathbb{C}: |z-a| <r\} $
    denote the open disc with center at $a \in \mathbb{C}$ and radius $r>0$. 
For any $j \in \{1,2, \cdots, N\}$, recall  that $\mathbb{I}^{N}_j =\{1,2, \cdots, N\} \setminus \{j\}$ by formula  \eqref{IntroCaloMoserSys}. 
For any $\vec{\lambda} \in \Omega_N$ and $\vec{\gamma} \in \mathbb{R}^{1 \times N}$,  
 we  define two  positive constants: 
\begin{equation}\label{KlamConstfoc}
\small
 \mathscr{K}(\vec{\lambda}):=1+ \max\{|\sum_{k \in \mathbb{I}^{N}_j}  ( \lambda_k - \lambda_j)^{-3} |: 1\leq j \leq N\},
  \qquad \mathscr{J}(\vec{\lambda}; \vec{\gamma}):= 1 + \tfrac{1}{2}\mathscr{K}(\vec{\lambda})+ 2 \max_{1\leq j\leq N}|M_{jj}|.
\end{equation}where $M_{jj}$ denotes the $jj$-entry of the matrix $M \in \mathbb{C}^{ N \times N}$ given by \eqref{defMlamgamFoc}, for any $j \in \{1,2, \cdots, N\}$. 
  Set $\mathbb{A}: \varepsilon \in \mathbb{R} \mapsto \mathbb{A}(\varepsilon): = 2V+ \varepsilon M \in \mathbb{C}^{N \times N}$. Then $\mathbb{A} \in \mathcal{C}^{\omega}\left(\mathbb{R}; \mathbb{C}^{N \times N} \right)$ and the diagonal matrix $\mathbb{A}(0)=2V \in \mathrm{GL}_N (\mathbb{R})$ has exactly $N$ distinct eigenvalues $\{2\lambda_1, 2 \lambda_2, \cdots, 2\lambda_N\} \subset \mathbb{R}$. Set
\begin{equation}\label{tildeltaRaySchCoeff}
\tilde{\delta}:=\min\{|\lambda_{n+1} - \lambda_n|: 1 \leq n \leq N-1\}>0
\end{equation}Fix $j \in \{1,2, \cdots, N\}$. Since $\mathrm{dim}_{\mathbb{C}} \mathrm{Ker}\left( \mathbb{A}(0)-2\lambda_j\right)=1$,  Theorem \uppercase\expandafter{\romannumeral 12}. 1 and its Corollary in Reed--Simon \cite{AnOpReedSimonIV} imply that $\exists
\varepsilon'_j, \delta_j \in (0, \min(1, \tilde{\delta}))$ such that:\\

\noindent $(\mathbf{\romannumeral 1})$. For any $\varepsilon \in (- \varepsilon'_j, \varepsilon'_j)$, the matrix $\mathbb{A}(\varepsilon) \in \mathbb{C}^{N \times N}$ has exactly $1$ eigenvalue, namely $\boldsymbol{\mu}_j(\varepsilon)$, in $\mathcal{D}(2 \lambda_j, \delta_j)$.\\
\noindent $(\mathbf{\romannumeral 2})$. For any $\varepsilon \in (- \varepsilon'_j, \varepsilon'_j)$, $\boldsymbol{\mu}_j(\varepsilon)  = \mathcal{D}(2 \lambda_j, \delta_j) \cap   \sigma \left(\mathbb{A}(\varepsilon) \right)$ is a simple eigenvalue of $\mathbb{A}(\varepsilon) \in \mathbb{C}^{N \times N}$.\\
\noindent $(\mathbf{\romannumeral 3})$. The function $\boldsymbol{\mu}_j: \varepsilon \in   (- \varepsilon'_j, \varepsilon'_j) \mapsto \boldsymbol{\mu}_j(\varepsilon) \in \mathcal{D}(2 \lambda_j, \delta_j)$ is real analytic. \\

\noindent Since $\left(\mathcal{D}(2 \lambda_j, \delta_j) \right)_{1 \leq j \leq N}$ is mutually disjoint, if $k \ne j$ and $|\varepsilon| <\min \{\varepsilon_1', \varepsilon_2', \cdots, \varepsilon_N'\}$, then $\boldsymbol{\mu}_k(\varepsilon) \ne \boldsymbol{\mu}_j(\varepsilon)$. For any $j \in \{1,2, \cdots, N\}$, the Taylor expansion of $\boldsymbol{\mu}_j \in \mathcal{C}^{\omega}\left((- \varepsilon'_j, \varepsilon'_j); \mathbb{C} \right)$ is given by
\begin{equation}\label{TayExpmuRaySch-}
\boldsymbol{\mu}_j (\varepsilon) = \sum_{m=0}^{M-1} \frac{\partial_x^m \boldsymbol{\mu}_j (0)}{m!} \varepsilon^m + \mathfrak{R}_M(\boldsymbol{\mu}_j; \varepsilon) \varepsilon^M = \sum_{m=0}^{M-1}  \mathfrak{c}_{m}^{(j)} \varepsilon^m + \mathfrak{R}_M(\boldsymbol{\mu}_j; \varepsilon)\varepsilon^M,
\end{equation}where $\varepsilon \in (-\varepsilon'_j, \varepsilon'_j)$, $\mathfrak{c}_{m}^{(j)}:=\frac{1}{m!} \partial_x^m \boldsymbol{\mu}_j (0)$ is called the $m$-th Rayleigh--Schr\"odinger coefficient of $\boldsymbol{\mu}_j$, for any $m \in [0, M]\cap   \mathbb{Z}$, the coefficient of integral remainder term is given by 
\begin{equation}\label{intRemainTerm-}
\mathfrak{R}_M(\boldsymbol{\mu}_j; \varepsilon):= \frac{1}{(M-1)!}\int_0^1 (1- y )^{M-1} \partial_x^M \boldsymbol{\mu}_j (y \varepsilon) \mathrm{d}y \in \mathbb{C},  \quad \forall M \in \mathbb{N}_+.
\end{equation}The classical perturbation theory in Reed--Simon \cite{AnOpReedSimonIV} indicates that $\mathfrak{c}_{0}^{(j)}=\boldsymbol{\mu}_j(0)=2 \lambda_j $, $\mathfrak{c}_{1}^{(j)} = M_{jj}\in \mathbb{C}_-$,  
\begin{equation}\label{c2jRaySchCoef-}
\mathfrak{c}_{2}^{(j)} = \tfrac{1 }{2} \partial_x^2 \boldsymbol{\mu}_j(0) = - \sum_{k \in \{1,2, \cdots, N\} \setminus \{j\}} \frac{M_{jk} M_{kj}}{\boldsymbol{\mu}_k (0) - \boldsymbol{\mu}_j (0)} =  - \sum_{k \in \mathbb{I}^{N}_j} \frac{1}{2( \lambda_k - \lambda_j)^3}\in \mathbb{R}.
\end{equation}where $\mathbb{I}^{N}_j=\{1,2, \cdots, N\} \setminus \{j\}$. Since $\boldsymbol{\mu}_j \in \mathcal{C}^2\left((- \varepsilon'_j, \varepsilon'_j); \mathbb{C} \right)$, there exists $\varepsilon_j \in (0, \varepsilon'_j) \subset (0, 1)$ such that
\begin{equation}\label{diffder1,2muj-}
\sup \{|\partial_x \boldsymbol{\mu}_j (\varepsilon) - \partial_x \boldsymbol{\mu}_j(0)| + |\partial_x^2 \boldsymbol{\mu}_j (\varepsilon) - \partial_x^2 \boldsymbol{\mu}_j(0)|: |\varepsilon| \leq \varepsilon_j \} \leq 1.
\end{equation}Then formulas \eqref{intRemainTerm-}, \eqref{c2jRaySchCoef-} and \eqref{diffder1,2muj-} imply that $\forall j \in \{1,2, \cdots, N\}$, we have
\begin{equation}\label{EstiR2mujep}
\sup_{|\varepsilon| \leq \varepsilon_j} |\partial_x \boldsymbol{\mu}_j(\varepsilon)| \leq 1 + |M_{jj}|, \quad  \sup_{|\varepsilon| \leq \varepsilon_j} |\partial_x^2 \boldsymbol{\mu}_j(\varepsilon)| \leq  \mathscr{K}(\vec{\lambda}), \quad  \sup_{|\varepsilon| \leq \varepsilon_j} |\mathfrak{R}_2(\boldsymbol{\mu}_j; \varepsilon)| \leq \tfrac{1}{2}\mathscr{K}(\vec{\lambda}),
\end{equation}where $\mathscr{K}(\vec{\lambda}) >1$ is defined in \eqref{KlamConstfoc}. We define that $\boldsymbol{\varepsilon}=\boldsymbol{\varepsilon}(\vec{\lambda}; \vec{\gamma}):=\min\{\varepsilon_1, \varepsilon_2, \cdots, \varepsilon_N\}  \in (0,1)$ and $\mathcal{T}:=\mathcal{T}(\vec{\lambda}; \vec{\gamma}):=\boldsymbol{\varepsilon}^{-1}>1$. For any $j \in \{1,2, \cdots, N\}$ and $t \in \mathbb{R}\setminus [-\mathcal{T}, \mathcal{T}]$, we define
\begin{equation}\label{zjdefFor-}
\mathbf{z}_j(t):=t \boldsymbol{\mu}_j (t^{-1}) = 2 \lambda_j t + M_{jj} + \mathfrak{R}_2(\boldsymbol{\mu}_j; t^{-1}) t^{-1} \in \mathbb{C}.
\end{equation}Since $M+2tV  = t \mathbb{A}(t^{-1})$, $\mathbf{z}_j(t)$ in \eqref{zjdefFor-} is a simple eigenvalue of the matrix $M+2tV \in \mathbb{C}^{N \times N}$. Lemma $\ref{SpofMLemfoc}$ yields that $\sigma (M + 2tV) =\sigma (M (\vec{\lambda}, \vec{\gamma}+2t\vec{\lambda}))  = \{\mathbf{z}_1(t),\mathbf{z}_2(t), \cdots, \mathbf{z}_N(t)\}\subset \mathbb{C}_-$. Moreover, 
\begin{equation}\label{Est(-)zj=2lamjt+Mjj}
\sup_{|t| > \mathcal{T}} |t\left(\mathbf{z}_j(t) - 2 \lambda_j t - M_{jj}\right)| \leq \sup_{|t| > \mathcal{T}} |\mathfrak{R}_2(\boldsymbol{\mu}_j; t^{-1})| \leq \tfrac{1}{2}\mathscr{K}(\vec{\lambda}).
\end{equation}Set $\mathtt{r}_{j}(t):=M_{jj}- \partial_x \boldsymbol{\mu}_j(t^{-1}) + \mathfrak{R}_2(\boldsymbol{\mu}_j; t^{-1}) t^{-1} \in \mathbb{C}$, $
\forall t \in \mathbb{R}\setminus [-\mathcal{T}, \mathcal{T}]$. We use \eqref{EstiR2mujep} to obtain that
\begin{equation}\label{Estrj(-)<J}
\sup_{|t| > \mathcal{T}}|\mathtt{r}_{j}(t)| \leq 1 + \tfrac{1}{2}\mathscr{K}(\vec{\lambda})+ 2  |M_{jj}| \leq \mathscr{J}(\vec{\lambda}; \vec{\gamma}),
\end{equation}where $\mathscr{J}(\vec{\lambda}; \vec{\gamma})>1$ is defined in \eqref{KlamConstfoc}. 
 Computing the derivative of  $\mathbf{z}_j$ in \eqref{zjdefFor-}, 
we obtain that 
\begin{equation}\label{FocpartzjLonTimAsy}
\partial_t \mathbf{z}_j(t) = 2 \lambda_j + \left(M_{jj}- \partial_x \boldsymbol{\mu}_j (t^{-1}) \right) t^{-1} + \mathfrak{R}_2(\boldsymbol{\mu}_j; t^{-1}) t^{-2} = 2 \lambda_j + \mathtt{r}_{j}(t) t^{-1}, \quad \forall t \in \mathbb{R}\setminus [-\mathcal{T}, \mathcal{T}].
\end{equation}So \eqref{Estrj(-)<J} yields that $\sup_{|t| >\mathcal{T}}|t(\partial_t \mathbf{z}_j(t) - 2 \lambda_j )| \leq \mathscr{J}(\vec{\lambda}; \vec{\gamma})$, $1 \leq j \leq N$. The proof is finished.
\end{proof}

 \subsection{Residue and Pole Dynamics}\label{6.4subsResPody}

   \noindent
The eigenvalue 
asymptotic estimates
  in
Proposition~\ref{propLoTiBeCMsysIntro} follow from
Lemma~\ref{RaySchrCoenegLem}; see formula 
\eqref{tzjttpartzjtEstiRaySch}.
The goal of this subsection is to complete the proof of that
proposition.
 By Lemma~\ref{RaySchrCoenegLem}, for sufficiently large \(|t|\),  the
inverse spectral matrix
\(\mathbf{M}(u(t))=\mathbf{M}(u_0)+2t\mathbf{V}(u_0)\) has \(N\)
distinct eigenvalues, and hence the denominator polynomial
\(Q_{u(t)}\) in \eqref{CharPolyQuGHpp1for} has exactly \(N\)
distinct roots in \(\mathbb C_-\).  
Thus the pole ansatz is available in the long time regime. 
 We shall derive
 the dynamical equations
satisfied by the poles and residues in the partial fraction expansion
\eqref{Ed2ParFracExpu}. Actually, the proposition below applies to
any time interval on which such an expansion with pairwise distinct
poles in \(\mathbb C_-\) is valid.

\begin{prop}\label{PropComCalo-Moser}
For any integer $N\geq 1$, let $u : t \in \mathbb{R} \mapsto u(t) \in \mathcal{U}_N$ be a  dark $N$-soliton solution of equation  \eqref{CMdNLSfocu}. Assume that there exist  $\alpha \in \mathbb{T}$,  an open subset $\mathcal{I} \subset 
\mathbb{R}$, and $2N$ functions
\begin{equation}\label{a1a2aNtdepdef}
\mathbf{a}_1, \mathbf{a}_2, \cdots, \mathbf{a}_N \; ; \;  \mathbf{z}_1, \mathbf{z}_2, \cdots, \mathbf{z}_N \in \mathcal{C}^1(\mathcal{I}; \mathbb{C})
\end{equation}such that for every $t \in \mathcal{I}$: 
\begin{enumerate}
\item the complex numbers $\mathbf{z}_1 (t), \mathbf{z}_2 (t) , \cdots, \mathbf{z}_N (t) \in \mathbb{C}_- $ are pairwise distinct;
\item the solution $u$ admits the partial fraction expansion:
\begin{equation}\label{ParFracExputxCMFoc-}
u(t,x) = e^{i \alpha} \Bigl( 1 + \sum_{j=1}^N \frac{\mathbf{a}_j(t)}{x - \mathbf{z}_j (t)} \Bigr), \quad \forall  x \in  \mathbb{R}.
\end{equation}
\end{enumerate}Then for each $j \in \{1,2, \cdots, N\}$, the residue  $\mathbf{a}_j $ and the pole  $  \mathbf{z}_j $ satisfy the following ODE system:
\begin{equation}\label{ODEajzjfirstordCM}
  \partial_t \mathbf{a}_j(t) = \sum_{k \in \mathbb{I}_j^N} \frac{2 i (\mathbf{a}_k(t) - \mathbf{a}_j(t)) }{(\mathbf{z}_k(t) - \mathbf{z}_j(t))^2}, \quad  \mathbf{a}_j(t)  \partial_t \mathbf{z}_j(t) = \sum_{k \in \mathbb{I}_j^N} \frac{2 i \mathbf{a}_k(t) }{ \mathbf{z}_k(t) - \mathbf{z}_j(t) } - 2i, \quad \mathbf{a}_j(t) \ne 0, 
\end{equation}with $\mathbb{I}_j^N  = \{1,2, \cdots, N \} \setminus \{j \}$,   $\forall t \in \mathcal{I}$. Consequently, the functions $\mathbf{a}_1, \mathbf{a}_2, \cdots, \mathbf{a}_N ;  \mathbf{z}_1, \mathbf{z}_2, \cdots, \mathbf{z}_N  $ all belong to $\mathcal{C}^2 (\mathcal{I}; \mathbb{C})$, and they solve the following second-order ODE system:
\begin{equation}\label{FocCaloMoserSys}
\partial_t^2 \mathbf{z}_j(t) = \sum_{k \in \mathbb{I}_j^N} \frac{8}{(\mathbf{z}_j(t) - \mathbf{z}_k(t))^3}, \quad \forall j\in \{1, 2 \cdots,N\}, \quad \forall t\in\mathcal{I},
\end{equation}
which is a complexified version of the classical Calogero--Moser system on $\mathbb{R}$.
\end{prop}

\begin{proof}
If $1 \leq k \ne j \leq N$, then $\mathbf{z}_k(t) \ne \mathbf{z}_j(t) \in \mathbb{C}_-$.  For any $t \in \mathcal{I}$,  we recall that
\begin{equation}\label{ParFractx-zkx-zjsq}
\tfrac{1}{ (x -  \mathbf{z}_k(t)) (x - \mathbf{z}_j (t) )^2 } = \tfrac{1}{ ( \mathbf{z}_k(t) - \mathbf{z}_j(t) )^2} \left( \tfrac{1}{x - \mathbf{z}_k(t)}  - \tfrac{1}{x - \mathbf{z}_j(t)}\right) - \tfrac{1}{  \mathbf{z}_k(t) - \mathbf{z}_j(t)  } \tfrac{1}{ ( x - \mathbf{z}_j(t) )^2},
\end{equation}for any $x \in \mathbb{C}\setminus \{\mathbf{z}_j(t) , \mathbf{z}_k(t) \}$. According to   \eqref{ParFracExputxCMFoc-}, we compute the derivatives of $u$:
\begin{equation}\label{dtuLapluCMsys}
\partial_t u (t,x) = e^{ i \alpha } \sum_{j=1}^N \left(\frac{\partial_t \mathbf{a}_j (t)}{x - \mathbf{z}_j (t)} + \frac{\mathbf{a}_j(t) \partial_t \mathbf{z}_j (t)}{(x - \mathbf{z}_j (t) )^2} \right),  \quad
i \partial_x^2 u (t, x) = 2 i e^{i \alpha} \sum_{j=1}^N \frac{\mathbf{a}_j(t)}{(x - \mathbf{z}_j(t))^3},
\end{equation}for any $(t, x) \in \mathcal{I} \times \mathbb{R}$. By formula \eqref{w+piwabs2=DtheovThe}, for any $(t, x) \in \mathcal{I} \times \mathbb{R}$, we define that
\begin{equation}\label{PiparxabsusqFor}
\mathfrak{J}_1(t,x)
: =- \left(\Pi \partial_x \left(|u|^2\right)\right)(t,x) 
= -\left( \partial_x \Pi(|u|^2 - 1) \right) (t,x) 
= \sum_{j=1}^N \frac{i}{(x - \mathbf{z}_j(t) )^2}.
\end{equation}The function $\mathfrak{J}_2: (t,x) \in \mathcal{I} \times \mathbb{R} \mapsto \mathfrak{J}_2(t,x) \in \mathbb{C}$ is defined as follows:
\begin{equation}\label{J2CMsys}
\begin{split}
 & \mathfrak{J}_2(t,x)
 :=  \sum_{j=1}^N \frac{\mathbf{a}_j(t) i }{( x - \mathbf{z}_j(t) )^3} + \sum_{1 \leq k \ne j \leq N} \frac{\mathbf{a}_k (t) i }{( x - \mathbf{z}_k(t) ) ( x - \mathbf{z}_j(t) )^2} \\
= &  \frac{i\partial_x^2 u (t,x) }{2 e^{i \alpha}} 
-\sum_{j  = 1}^N \left(  \sum_{k \in \mathbb{I}_j^N} \frac{\mathbf{a}_k (t) i }{\mathbf{z}_k (t) - \mathbf{z}_j (t)  } \right) \frac{1}{(x - \mathbf{z}_j(t))^2} 
+ \sum_{j  = 1}^N \left(  \sum_{k \in \mathbb{I}_j^N} \frac{( \mathbf{a}_j(t)  - \mathbf{a}_k(t) ) i }{( \mathbf{z}_k(t)  - \mathbf{z}_j(t) )^2} \right) \frac{1}{x - \mathbf{z}_j(t)}.
\end{split}
\end{equation}where we use  \eqref{ParFractx-zkx-zjsq} and \eqref{dtuLapluCMsys} in the last step. By  \eqref{ParFracExputxCMFoc-},  \eqref{PiparxabsusqFor} and \eqref{J2CMsys}, it follows that
\begin{equation}\label{J1plusJ2ajakzjzkFor}
\small
\begin{split}
 & \left( -2  u \Pi \partial_x (|u|^2) \right) (t,x) = 2 e^{i \alpha} \left(\sum_{j=1}^N  \frac{i}{(x - \mathbf{z}_j(t))^2} + \sum_{1 \leq k, j \leq N} \frac{\mathbf{a}_k(t) i}{(x - \mathbf{z}_k(t))(x - \mathbf{z}_j(t))^2}\right)\\
= & 2 e^{i \alpha} \left( \mathfrak{J}_1 (t,x) + \mathfrak{J}_2 (t,x) \right)  \\
= & i \partial_x^2 u (t,x) +  \sum_{j=1}^N \left( \left( 2i - \sum_{k \in \mathbb{I}_j^N} \frac{ 2 i \mathbf{a}_k (t)   }{\mathbf{z}_k (t) - \mathbf{z}_j (t)  } \right) \frac{e^{i \alpha}}{ ( x- \mathbf{z}_j(t) )^2  } + \left(  \sum_{k \in \mathbb{I}_j^N} \frac{ 2 i ( \mathbf{a}_j(t)  - \mathbf{a}_k(t) )   }{( \mathbf{z}_k(t)  - \mathbf{z}_j(t) )^2} \right) \frac{e^{i \alpha}}{ x- \mathbf{z}_j(t)} \right).
\end{split}
\end{equation}Since $\{\frac{1}{x - \mathbf{z}_j(t)}\}_{1 \leq j \leq N}  \cup \{\frac{1}{(x - \mathbf{z}_j(t))^2}\}_{1 \leq j \leq N}$ is a linearly independent subset in $L^2_+$,  formula  \eqref{ODEajzjfirstordCM} is a direct consequence of Lemma $\ref{LemPQhasnocomrootsFoc-}$, formulas \eqref{J1plusJ2ajakzjzkFor} and \eqref{dtuLapluCMsys}. So $\mathbf{a}_j$ and $\mathbf{z}_j$ both belong to $\mathcal{C}^2(\mathcal{I}; \mathbb{C})$, for any $j \in \{1,2, \cdots, N\}$. If $1 \leq m, n , j \leq N$ such that $m \ne n$, $n\ne j$ and $j \ne m$, we have
\begin{equation}\label{AmnjdefParFracExpan}
 \mathtt{A}_{mnj}:=     \tfrac{1}{(\mathbf{z}_m - \mathbf{z}_j)(\mathbf{z}_n - \mathbf{z}_m)^2} - \tfrac{1}{(\mathbf{z}_n - \mathbf{z}_j)(\mathbf{z}_n - \mathbf{z}_m)^2} - \tfrac{1}{(\mathbf{z}_n - \mathbf{z}_j)^2 (\mathbf{z}_n - \mathbf{z}_m)} =\tfrac{1}{(\mathbf{z}_m - \mathbf{z}_j)(\mathbf{z}_n - \mathbf{z}_j)^2}  \in \mathcal{C}^2(\mathcal{I}; \mathbb{C})   ,
\end{equation}by formula 
\eqref{ParFractx-zkx-zjsq}. For any sequence \((c_{mnj})_{m,n,j\in\mathbb{N}_+}\),  fix  \(j \in \{1,2, \cdots, N\}\), we have the following identity:
\begin{equation}\label{cmnjsumiden}
\sum_{n\in\mathbb{I}_j^N}\;\sum_{m\in\mathbb{I}_n^N} c_{mnj}
= \sum_{\substack{m,n\in\mathbb{I}_j^N\\ m\neq n}} c_{mnj}
+ \sum_{n\in\mathbb{I}_j^N} c_{jnj},
\end{equation}where $\mathbb{I}_j^N, \mathbb{I}_n^N$ are given by \eqref{inteintervIjN}. Using \eqref{AmnjdefParFracExpan} and \eqref{cmnjsumiden}, we obtain that
\begin{equation}\label{J1toderiveCMsys}
\begin{split}
& \mathtt{J}_1 : =   \sum_{n \in \mathbb{I}_j^N} \sum_{m \in \mathbb{I}_n^N} \left( \frac{\mathbf{a}_n - \mathbf{a}_m}{(\mathbf{z}_n - \mathbf{z}_j)(\mathbf{z}_m - \mathbf{z}_n)^2} + \frac{\mathbf{a}_m}{(\mathbf{z}_n - \mathbf{z}_j)^2(\mathbf{z}_m - \mathbf{z}_n)} \right) \in \mathcal{C}^2(\mathcal{I}; \mathbb{C}) \\
= & \sum_{n\in\mathbb{I}_j^N} \frac{2 \mathbf{a}_j - \mathbf{a}_n}{(\mathbf{z}_j - \mathbf{z}_n)^3}
+ \sum_{\substack{m,n\in\mathbb{I}_j^N\\ m\neq n}}  \mathbf{a}_m  \mathtt{A}_{mnj} 
=   \sum_{n\in\mathbb{I}_j^N} \frac{2 \mathbf{a}_j  }{(\mathbf{z}_j - \mathbf{z}_n)^3} +  \sum_{ m,n\in\mathbb{I}_j^N } \frac{\mathbf{a}_m }{(\mathbf{z}_m - \mathbf{z}_j)(\mathbf{z}_n - \mathbf{z}_j)^2}  .
\end{split}
\end{equation}The function $\mathtt{J}_2 \in \mathcal{C}^2(\mathcal{I}; \mathbb{C})$ is defined as follows:
\begin{equation}\label{J2toderCMSys}
\mathtt{J}_2: = -4 \sum_{ n \in  \mathbb{I}_j^N} \frac{1}{(\mathbf{z}_n - \mathbf{z}_j)^2} + 2i (\partial_t \mathbf{z}_j) \sum_{ n \in  \mathbb{I}_j^N} \frac{\mathbf{a}_n}{(\mathbf{z}_n - \mathbf{z}_j)^2}.
\end{equation}Differentiating \(\mathbf{a}_j \partial_t \mathbf{z}_j\) and applying \eqref{ODEajzjfirstordCM}, \eqref{J1toderiveCMsys}  and \eqref{J2toderCMSys}, we obtain that
\begin{equation}\label{derajparttzj4J1J2}
\small
\begin{split}
&  \partial_t\mathbf{a}_j \partial_t \mathbf{z}_j +  \mathbf{a}_j \partial_t^2 \mathbf{z}_j = \partial_t (\mathbf{a}_j \partial_t \mathbf{z}_j) =   2i \sum_{n \in \mathbb{I}_j^N} \left( \frac{\partial_t \mathbf{a}_n}{\mathbf{z}_n - \mathbf{z}_j} - \frac{\mathbf{a}_n (\partial_t \mathbf{z}_n - \partial_t \mathbf{z}_j)}{(\mathbf{z}_n - \mathbf{z}_j)^2} \right) 
\\
=& \sum_{n \in \mathbb{I}_j^N} \frac{2 i }{\mathbf{z}_n - \mathbf{z}_j} \sum_{m \in \mathbb{I}_n^N} \frac{2i (\mathbf{a}_m - \mathbf{a}_n)}{(\mathbf{z}_m - \mathbf{z}_n)^2}
 + \sum_{n \in \mathbb{I}_j^N} \frac{2 i \mathbf{a}_n \partial_t \mathbf{z}_j}{(\mathbf{z}_n - \mathbf{z}_j)^2}
  -\sum_{n \in \mathbb{I}_j^N} \frac{2 i }{(\mathbf{z}_n - \mathbf{z}_j)^2} \left( \sum_{m \in \mathbb{I}_n^N} \frac{2 i \mathbf{a}_m }{ \mathbf{z}_m - \mathbf{z}_n }  - 2 i\right) \\
=&    4 \mathtt{J}_1 + \mathtt{J}_2.
\end{split}
\end{equation}Recall that $\partial_t \mathbf{a}_j = \sum_{n \in \mathbb{I}_j^N} \frac{2 i ( \mathbf{a}_n -  \mathbf{a}_j ) }{(\mathbf{z}_n - \mathbf{z}_j)^2} \in \mathcal{C}^1(\mathcal{I}; \mathbb{C})$ by the first formula of  \eqref{ODEajzjfirstordCM}. Then we have
\begin{equation}\label{dtajdtzjminusJ2CM}
\small
\mathtt{J}_2  - \partial_t \mathbf{a}_j \partial_t \mathbf{z}_j   =   \sum_{n \in \mathbb{I}_j^N} \frac{2 i \mathbf{a}_j \partial_t \mathbf{z}_j - 4 }{(\mathbf{z}_n - \mathbf{z}_j)^2} = - 4 \sum_{m,n \in \mathbb{I}_j^N} \frac{\mathbf{a}_m}{(\mathbf{z}_m - \mathbf{z}_j)(\mathbf{z}_n - \mathbf{z}_j)^2} = \sum_{n \in \mathbb{I}_j^N} \frac{8 \mathbf{a}_j}{(\mathbf{z}_j - \mathbf{z}_n)^3} - 4 \mathtt{J}_1,
\end{equation}by using     \eqref{ODEajzjfirstordCM}, \eqref{J1toderiveCMsys}  and \eqref{J2toderCMSys} again. 
  By Lemma~\ref{LemPQhasnocomrootsFoc-} and 
 formula  \eqref{ODEajzjfirstordCM}, we have  
$\mathbf{a}_j(t) \ne 0$, $\forall t \in \mathcal{I}$. Consequently, formula  \eqref{FocCaloMoserSys} follows from  \eqref{derajparttzj4J1J2} and \eqref{dtajdtzjminusJ2CM}. The proof is completed.
\end{proof}

 \begin{rem}\label{ComRemGLBriNSolPoAn}
We compare the pole and residue dynamics of bright and dark
multi-soliton solutions of equation~\eqref{CMdNLSfocu}. 
We use boldface notation for the dark case and ordinary notation for
the bright case.
 Fix $N \in \mathbb{N}_+$. 
 For bright
\(N\)-soliton solutions, the pole ansatz takes the form
\begin{equation}\label{Sec6BriNsoliGL}
\mathfrak{u}^N(t,x)
=
\sum_{j=1}^N
\frac{a_j(t)}{x-z_j(t)},
\qquad
z_j(t)\in\mathbb C_- ,
\qquad
1 \leq j \leq N,
\qquad
x \in \mathbb{R},
\qquad
|t| \gg 1 .
\end{equation}
See also formula \eqref{BriNsoliCMdNLS}.  According to
G\'erard--Lenzmann~\cite[Subsection~6.1]{GerLenz22CPAM},
the poles \((z_j)_{1 \leq j \leq N}\) and residues
\((a_j)_{1 \leq j \leq N}\) of \eqref{Sec6BriNsoliGL} evolve
according to the following first-order system
\begin{equation}\label{BriMulODEajzj}
\partial_t a_j
=
\sum_{k\in\mathbb I_j^N}
\frac{2i(a_k-a_j)}{(z_k-z_j)^2},
\qquad
a_j\,\partial_t z_j
=
\sum_{k\in\mathbb I_j^N}
\frac{2i\,a_k}{z_k-z_j},
\qquad
|t|\gg1,
\end{equation}
where \(\mathbb I_j^N=\{1, 2,\cdots, N\}\setminus\{j\}\).
Moreover, the poles \(z_j\) satisfy the same complexified
Calogero--Moser system~\eqref{FocCaloMoserSys}. 
Comparing \eqref{BriMulODEajzj} with
\eqref{ODEajzjfirstordCM}, the only difference between the bright
and dark first-order systems is the constant term \(-2i\) in the
equation for \(a_j\,\partial_t z_j\).  In the dark case, this extra
constant is ultimately cancelled in the derivation of the
second-order system: differentiating 
  \(\mathbf{a}_j\,\partial_t \mathbf{z}_j\)   
produces a contribution that is exactly compensated by a term
arising from \(\partial_t \mathbf{a}_j\); see \eqref{derajparttzj4J1J2} and
\eqref{dtajdtzjminusJ2CM}.  In the bright case, the corresponding
two contributions are both absent, so the same
Calogero--Moser system~\eqref{FocCaloMoserSys} is obtained without
the need for such a cancellation.  Thus the bright multi-soliton
system~\eqref{BriMulODEajzj} may be regarded as the degenerate case
of the dark one~\eqref{ODEajzjfirstordCM}.
\end{rem}

\noindent Proposition~\ref{propLoTiBeCMsysIntro} is a direct consequence of Lemma~\ref{RaySchrCoenegLem} and  Proposition~\ref{PropComCalo-Moser}.

 \subsection{Long Time Behavior}\label{subsLonTiBe6.5}
\noindent In this subsection, we prove   
 Corollary $\ref{corLongtimebehave}$  
and Corollary $\ref{DarkSolRescor}$. Throughout this subsection, 
 we fix an integer \(N  \geq 1\) and parameters \(\alpha \in \mathbb{T}\), \(  \vec{\lambda} = (\lambda_1, \lambda_2,  \cdots, \lambda_N) \in \Omega_N\), and \(\vec{\gamma} =   (\gamma_1, \gamma_2,  \cdots, \gamma_N) \in \mathbb{R}^{1 \times N}\). Let \(M = M(\vec{\lambda}, \vec{\gamma}) =  (M_{kj})_{1 \le k,j \le N} \in \mathbb{C}^{N \times N}\) be the matrix defined in \eqref{defMlamgamFoc} and let \(V = V_{\vec{\lambda}} \in \mathrm{GL}_N(\mathbb{R})\) be the matrix defined in \eqref{Vdiag-}. For any $(t,x) \in \mathbb{R}$, we define
 \begin{equation}\label{utxInvSpecFortdep}
 u(t,x) := e^{i\alpha}\Bigl(1 + i\bigl\langle (M + 2tV - x )^{-1} X,\; Y \bigr\rangle_{\mathbb{C}^N}\Bigr),
\end{equation}
where \(X = X_{\vec{\lambda}} \in \mathbb{R}^N\) and \(Y = Y_{\vec{\lambda}} \in \mathbb{R}^N\) are given by \eqref{defXYFoclambda}. Then \(u: t \in \mathbb{R} \to u(t) \in  \mathcal{U}_N\) is a pure dark \(N\)-soliton solution of equation \eqref{CMdNLSfocu}. Let $\mathcal{T} = \mathcal{T}(\vec{\lambda}, \vec{\gamma}) > 1$ be the constant and $\mathbf{z}_j \in \mathcal{C}^{\omega}(\mathbb{R}\setminus [-\mathcal{T} , \mathcal{T}] ; \mathbb{C})$ the function  defined in Lemma~\ref{RaySchrCoenegLem}, $\forall j \in \{1,2, \cdots, N\}$. For any $t \in \mathbb{R}\setminus [-\mathcal{T} , \mathcal{T}]$, the matrix $M + 2 t V \in \mathbb{C}^{N \times N}$ has exactly $N$ distinct eigenvalues $\mathbf{z}_1(t), \mathbf{z}_2(t), \cdots, \mathbf{z}_N(t) \in \mathbb{C}_-$ by Lemma~\ref{SpofMLemfoc} and Lemma~\ref{RaySchrCoenegLem}. By Proposition~\ref{PropComCalo-Moser}, the solution $u \in \mathcal{C}^{\infty}(\mathbb{R} \times \mathbb{R}; \mathbb{C})$ in  \eqref{utxInvSpecFortdep} admits the partial fraction expansion: 
\begin{equation}\label{uPaFrExtbig}
\small
 u(t,x) = e^{ i\alpha} \left( 1 +  \sum_{j = 1}^N \frac{\mathbf{a}_j(t)}{x - \mathbf{z}_j(t)} \right), \quad \forall x \in \mathbb{R}, \quad \forall t \in  \mathbb{R} \setminus [-\mathcal{T}, \mathcal{T}],   
\end{equation}where the functions $\mathbf{a}_1 , \mathbf{a}_2 , \cdots, \mathbf{a}_N \in \mathcal{C}^{\infty}(\mathbb{R}\setminus [-\mathcal{T} , \mathcal{T}]; \mathbb{C})$ solve  system  \eqref{ODEajzjfirstordCM}.  As a consequence, the poles  $\mathbf{z}_1 , \mathbf{z}_2 , \cdots, \mathbf{z}_N \in \mathcal{C}^{\omega}(\mathbb{R}\setminus [-\mathcal{T} , \mathcal{T}]; \mathbb{C})$ solve the complexified Calogero--Moser system \eqref{FocCaloMoserSys}.
  The next lemma gives the long time estimates of  poles and   residues for the dark $N$-soliton solution $u$ in \eqref{uPaFrExtbig}. 
 
\begin{lem}\label{ajImzjLonTimAsymBehLem}
There exists a constant $\widetilde{\mathcal{T}}  = \widetilde{\mathcal{T}}  (\vec{\lambda}, \vec{\gamma}) \geq \mathcal{T}>1$ such that for any $t \in \mathbb{R}\setminus [-\widetilde{\mathcal{T}}  , \widetilde{\mathcal{T}}  ]$, we have
\begin{equation}\label{ajImzjLTEstimate}
\tfrac{1}{2|\lambda_j|} < |\mathbf{a}_j(t)| < \tfrac{2}{|\lambda_j|}  \quad \mathrm{and} \quad  \tfrac{1}{4 \lambda_j ( 1 + \lambda_j ) } < | \operatorname{Im}\mathbf{z}_j(t)| = - \operatorname{Im}\mathbf{z}_j(t) < \tfrac{1}{  \lambda_j ( 1 + \lambda_j ) },
\end{equation}for any $j \in \{1,2, \cdots, N\}$.
\end{lem}

\begin{proof}
For any $t\in \mathbb{R}\setminus [-\mathcal{T}, \mathcal{T}]$, we define  $\mathbf{a}(t):= (\mathbf{a}_1(t), \mathbf{a}_2(t), \cdots, \mathbf{a}_N(t))^{\mathrm{T}} \in \mathbb{C}^N$. The matrices $\mathsf{P}(t) = (\mathsf{P}_{kj}(t))_{1 \leq k, j \leq N} \in \mathbb{C}^{N \times N}$ and  $\mathbf{Z}(t) = (\mathbf{Z}_{kj}(t))_{1 \leq k, j \leq N} \in \mathbb{C}^{N \times N}$ are defined as follows:
\begin{equation}\label{PZmatridefFor}
\mathsf{P}_{kj}(t) = \tfrac{2 i }{\mathbf{z}_j(t) - \mathbf{z}_k(t)}, \quad \mathsf{P}_{jj}(t) = 0 = \mathbf{Z}_{kj}(t), \quad \mathbf{Z}_{jj} (t) = \mathbf{z}_j(t), \quad \forall t\in \mathbb{R}\setminus [-\mathcal{T}, \mathcal{T}],
\end{equation}if $1 \leq k \ne j \leq N$. By formulas \eqref{zjdefFor-} and  \eqref{FocpartzjLonTimAsy}, 
it follows that 
\begin{equation}\label{limPanddZttoinfty}
\lim_{|t| \to +\infty}  \mathsf{P}(t)  = 0_{N \times N} \in \mathbb{C}^{N \times N} \quad \mathrm{and} \quad \lim_{|t| \to +\infty}\partial_t \mathbf{Z}(t)  = 2V \in  \mathbb{C}^{N \times N}.
\end{equation}Let $\{\boldsymbol{\zeta}_j\}_{1 \leq j \leq N}$ be the standard basis of $\mathbb{R}^N$ as given in Definition~\ref{zetajstanONBRN}. Set $\boldsymbol{\xi} = \sum_{n=1}^N \boldsymbol{\zeta}_n  $. We have
\begin{equation}
\left( \partial_t \mathbf{Z}(t) - \mathsf{P}(t) \right) \mathbf{a}(t) = -2 i \boldsymbol{\xi} \in \mathbb{C}^N ,  \quad \forall t\in \mathbb{R}\setminus [-\mathcal{T}, \mathcal{T}],
\end{equation}by \eqref{PZmatridefFor} and the second formula of \eqref{ODEajzjfirstordCM}.
 Since $\mathrm{GL}_N(\mathbb{C})$ is open in $\mathbb{C}^{ N \times N}$ and $V\in \mathrm{GL}_N(\mathbb{C})$, formula \eqref{limPanddZttoinfty} yields that there exists a constant $\mathtt{T}_1 =\mathtt{T}_1 (\vec{\lambda}, \vec{\gamma})> \mathcal{T}$ such that
 \begin{equation}
 \partial_t \mathbf{Z}(t) - \mathsf{P}(t) \in \mathrm{GL}_N(\mathbb{C}), \quad \forall  t \in \mathbb{R}\setminus [- \mathtt{T}_1, \mathtt{T}_1].
 \end{equation}So $ \mathbf{a}(t) = \left( \partial_t \mathbf{Z}(t) - \mathsf{P}(t) \right)^{-1} (-2i \boldsymbol{\xi}) \to \sum_{j=1}^N \frac{-i}{\lambda_j} \boldsymbol{\zeta}_j$ in $\mathbb{C}^N$, as $|t| \to +\infty$. Then there exists a constant $\mathtt{T}_2 =\mathtt{T}_2 (\vec{\lambda}, \vec{\gamma})\geq  \mathtt{T}_1$  such that   for any $j \in \{1,2, \cdots, N\}$, we have $ \tfrac{1}{2 |\lambda_j|} < |\mathbf{a}_j (t) | < \tfrac{2}{ |\lambda_j|}  $, $\forall  t \in \mathbb{R}\setminus [- \mathtt{T}_2, \mathtt{T}_2]$. For any $j \in \{1,2, \cdots, N\}$, using formulas \eqref{defMlamgamFoc},  
 \eqref{EstiR2mujep} and 
\eqref{zjdefFor-}, we obtain that
\begin{equation}
\operatorname{Im}\mathbf{z}_j(t) = - \tfrac{1}{2 \lambda_j (1+\lambda_j)} + \operatorname{Im} \mathfrak{R}_2 (\boldsymbol{\mu}_j ; t^{-1})t^{-1} \to - \tfrac{1}{2 \lambda_j (1+\lambda_j)}, \quad \mathrm{as} \quad |t| \to + \infty.
\end{equation}Then there exists a constant $\widetilde{\mathcal{T}} = \widetilde{\mathcal{T}}(\vec{\lambda}, \vec{\gamma}) \geq \mathtt{T}_2 \geq \mathcal{T}$ such that for any $ j \in \{1,2, \cdots, N\}$, we have $- \frac{1}{\lambda_j (1 + \lambda_j)} <  \operatorname{Im}\mathbf{z}_j (t) < - \frac{1}{4\lambda_j (1 + \lambda_j)}$,   $\forall  t \in \mathbb{R} \setminus [-\widetilde{\mathcal{T}}, \widetilde{\mathcal{T}}]$. 
  This proves the lemma.
 \end{proof}

 \noindent  With the long time estimate \eqref{ajImzjLTEstimate} at hand, we are ready to prove Corollary $\ref{corLongtimebehave}$. 
 
\begin{proof}[Proof of Corollary $\ref{corLongtimebehave}$]
By Theorem $\ref{ActionAngleMainThm}$ and Proposition $\ref{PropInvSpecFormula}$,  every pure dark $N$-soliton solution $u$ of \eqref{CMdNLSfocu} can be expressed as the inverse spectral formula \eqref{invSPEFormain} with time evolution. Hence we may assume that $u: t \in \mathbb{R} \mapsto u(t) \in \mathcal{U}_N$ takes the form given in \eqref{utxInvSpecFortdep}.
Recall that $e^{i\alpha} = u_0(\infty) = \lim_{|x| \to +\infty}u(t,x)$, for any $t \in \mathbb{R}$. Let $\mathcal{T}= \mathcal{T}(\vec{\lambda}, \vec{\gamma}) >1$ be defined in Lemma $\ref{RaySchrCoenegLem}$.  If $1 \leq k \ne j \leq N$, we define 
\begin{equation}\label{pkjdkjFor}
\mathfrak{d}_{kj}(t):=\mathbf{z}_k (t) - \overline{\mathbf{z}_j (t)} 
\quad \mathrm{and} \quad 
\mathfrak{p}_{kj}(t) : = \tfrac{M_{kk} - \overline{M}_{jj}}{2 (\lambda_k - \lambda_j)} + \tfrac{\mathfrak{R}_2 (\boldsymbol{\mu}_k ; t^{-1}) - \overline{\mathfrak{R}_2 (\boldsymbol{\mu}_j ; t^{-1})}}{2 (\lambda_k - \lambda_j) t}, \quad \forall t \in \mathbb{R}\setminus [-\mathcal{T}, \mathcal{T}],
\end{equation}where  $\mathbf{z}_j \in \mathcal{C}^{\omega}( \mathbb{R}\setminus [-\mathcal{T}, \mathcal{T}] ; \mathbb{C})$ is given by 
\eqref{zjdefFor-} and Lemma~\ref{RaySchrCoenegLem}, $M_{jj} \in \mathbb{C}_-$ is defined in \eqref{defMlamgamFoc}. The remainder term $\mathfrak{R}_2 (\boldsymbol{\mu}_j ; t^{-1}) $ given by
\eqref{intRemainTerm-} satisfies inequality  \eqref{Est(-)zj=2lamjt+Mjj}. So we have
\begin{equation}\label{pkjestimUnif}
\max_{1 \leq k \ne j \leq N}\sup\{|\mathfrak{p}_{kj}(t)|: |t| \geq \mathcal{T} \} \leq (2 \tilde{\delta})^{-1} \left( \|M\|_{\mathbb{C}^{N \times N}} + \mathscr{K}(\vec{\lambda}) \right),
\end{equation}where $\tilde{\delta} = \tilde{\delta}(\vec{\lambda})>0$ is given  by \eqref{tildeltaRaySchCoeff} and $\mathscr{K}(\vec{\lambda}) \geq 1$  is defined in 
\eqref{KlamConstfoc}. Using formulas \eqref{zjdefFor-},  \eqref{pkjdkjFor} and \eqref{pkjestimUnif}, we obtain that $\frac{\mathfrak{d}_{kj}(t)}{2 (\lambda_k - \lambda_j) t} = 1 + \frac{\mathfrak{p}_{kj}(t)}{t} \to 1$, as $|t| \to +\infty$, if $1 \leq k \ne j \leq N$. Let $\mathtt{T}_3:= \widetilde{\mathcal{T}}(\vec{\lambda}, \vec{\gamma}) \geq \mathcal{T}$ be defined in Lemma~\ref{ajImzjLonTimAsymBehLem}.
Then there exists  $\mathtt{T}_4 = \mathtt{T}_4 (\vec{\lambda}, \vec{\gamma}) \geq \mathtt{T}_3$ such that
\begin{equation}\label{estiabsdkjtLongTime}
\tilde{\delta} |t| \leq  |\lambda_k - \lambda_j | |t| < |\mathfrak{d}_{kj}(t)| < 4 |\lambda_k - \lambda_j | |t|, \quad \forall t \in \mathbb{R}\setminus [-
\mathtt{T}_4, \mathtt{T}_4 ], \quad  \mathrm{if} \quad 1 \leq k \ne j \leq N.
\end{equation}By formulas \eqref{zjdefFor-},  \eqref{Est(-)zj=2lamjt+Mjj} and 
 \eqref{pkjdkjFor}, for any $t \in \mathbb{R} \setminus [-\mathcal{T}, \mathcal{T}]$, if $1 \leq k \ne j \leq N$, it follows that
\begin{equation}
\operatorname{Im} \mathfrak{d}_{kj}(t) =   - \tfrac{1}{2 \lambda_k (1 + \lambda_k)} 
- \tfrac{1}{2 \lambda_j (1 + \lambda_j)} + t^{-1} \operatorname{Im}\left( \mathfrak{R}_2 (\boldsymbol{\mu}_k ; t^{-1})  + \mathfrak{R}_2( \boldsymbol{\mu}_j ; t^{-1}) \right) .
\end{equation}For any $k, j \in \{1,2, \cdots, N\}$ such that $k \ne j$,  
we define $\mathrm{b}_{kj} = \mathrm{b}_{kj}(\vec{\lambda}) >0$ as follows,
\begin{equation}\label{fracbkjdefForFoc-}
-\tfrac{1}{2} \mathrm{b}_{kj} : = - \tfrac{1}{2 \lambda_k (1 + \lambda_k)} 
- \tfrac{1}{2 \lambda_j (1 + \lambda_j)}
 =  \lim_{|t| \to + \infty} \operatorname{Im} \mathfrak{d}_{kj}(t) = - \tfrac{1}{2} \mathrm{b}_{jk}.
\end{equation}Using formulas  \eqref{ajImzjLTEstimate} and \eqref{pkjdkjFor},   if $1 \leq k \ne j \leq N $ and $t \in \mathbb{R}\setminus [-\mathtt{T}_3, \mathtt{T}_3]$, we obtain that 
\begin{equation}\label{ImzkplusImzjEstima}
-\mathrm{b}_{kj}  < \operatorname{Im}\mathbf{z}_k(t)  + \operatorname{Im}\mathbf{z}_j(t)  = \operatorname{Im}\mathfrak{d}_{kj}(t) < 
   - \tfrac{1}{4} \mathrm{b}_{kj}.
\end{equation}For any $s>0$ and $\rho>0$, we set $\mathtt{F}(s, \rho) : = s \int_0^{+\infty} ( 1 + |\xi|^2 )^s \exp(-\frac{1}{4}\rho \xi) \mathrm{d}\xi >0$. Then we define
\begin{equation}\label{AslambdadefFor}
\mathtt{A}(s, \vec{\lambda}) : = \max \{\mathtt{F}(s, \mathrm{b}_{kj}(\vec{\lambda}))  : 1 \leq  k \ne j \leq N\}>0, \quad \forall s >0.
\end{equation}For any $t \in \mathbb{R}\setminus [- \mathtt{T}_3, \mathtt{T}_3]$, we set $\Xi^{kj}_s (t) : = 2s \int_0^{+\infty} \xi ( 1 + \xi^2 )^{s-1} \exp(-i \xi \mathfrak{d}_{kj}(t)  )  \mathrm{d}\xi $. So we have
\begin{equation}\label{XikjstEstimsmalA}
\max_{1 \leq k \ne j \leq N} \sup_{|t| > \mathtt{T}_3} |\Xi^{kj}_s (t) | \leq \mathtt{A}(s , \vec{\lambda}),
\end{equation}by formulas \eqref{ImzkplusImzjEstima} and \eqref{AslambdadefFor}.  Fix $s >0$. By formulas  \eqref{LaxBeurlingSG} and \eqref{ImzkplusImzjEstima}, it follows that
\begin{equation}\label{lim1+xisexp-idkjtxi}
\lim_{\xi \to +\infty} (1 + \xi^2)^s  \exp(-i \xi \mathfrak{d}_{kj}(t)  )  = 0, \quad \forall t \in \mathbb{R}\setminus [-\mathtt{T}_3 , \mathtt{T}_3 ], \quad \mathrm{if} \quad 1 \leq k \ne j \leq N.
\end{equation}For any $j \in \{1,2, \cdots, N\}$, we set $\boldsymbol{\varphi}_j(t,x) := (x - \mathbf{z}_j(t))^{-1}$, $\forall x \in \mathbb{R}$, for any $t \in \mathbb{R}\setminus [-\mathtt{T}_3 , \mathtt{T}_3 ]$.  Using formulas 
\eqref{FouTrx-ztinv},  \eqref{pkjdkjFor} and \eqref{lim1+xisexp-idkjtxi}, if $1 \leq k \ne j \leq N$ and $t \in \mathbb{R}\setminus [-\mathtt{T}_3, \mathtt{T}_3]$,  we obtain that
\begin{equation}\label{vphjvphkinprodintbyPart}
\langle \boldsymbol{\varphi}_k(t) , \boldsymbol{\varphi}_j(t)\rangle_{H^s_+} = 2 \pi \int_0^{+\infty} (1+ \xi^2 )^s \exp(-i \xi \mathfrak{d}_{kj}(t)) \mathrm{d}\xi  = - \frac{2 \pi i}{\mathfrak{d}_{kj}(t)} \left( 1 + \Xi_s^{kj} (t) \right),
\end{equation}by integration by parts. By formulas \eqref{estiabsdkjtLongTime},
\eqref{XikjstEstimsmalA} 
and \eqref{vphjvphkinprodintbyPart}, we deduce that
\begin{equation}\label{EstimaxvphkvphijinProd}
\max_{1 \leq k \ne j \leq N} |\langle \boldsymbol{\varphi}_k(t) , \boldsymbol{\varphi}_j(t)\rangle_{H^s_+}| \leq  2 \pi \left( 1 + \mathtt{A}(s, \vec{\lambda}) \right) \tilde{\delta}^{-1} |t|^{-1}, \quad \forall t \in \mathbb{R}\setminus [- \mathtt{T}_4, \mathtt{T}_4] .
\end{equation}For any $j \in \{1,2, \cdots, N\}$, we set $\mathtt{K}_j = \mathtt{K}_j(\lambda_j) : = 1 + \max \{4 \lambda_j(1+\lambda_j) \; , \; \frac{1}{\lambda_j(1+\lambda_j)}\} >1$. Formula \eqref{ajImzjLTEstimate} yields that $\mathtt{K}_j^{-1} < \frac{1}{4 \lambda_j(1+\lambda_j)}  < |\operatorname{Im} \mathbf{z}_j(t)| < \frac{1}{  \lambda_j(1+\lambda_j)} < \mathtt{K}_j$, for any $t \in \mathbb{R} \setminus [- \mathtt{T}_3, \mathtt{T}_3]$. By Lemma $\ref{BdvphiLongtimeLem}$, there exists a constant $\mathrm{B}(s, \mathtt{K}_j) >1$ such that 
\begin{equation}\label{vphjtHsNormulbound}
\mathrm{B}(s, \mathtt{K}_j)^{-1}  \leq \|\boldsymbol{\varphi}_j(t) \|_{H^s_+} \leq \mathrm{B}(s, \mathtt{K}_j) , \quad \forall t \in \mathbb{R} \setminus [-\mathtt{T}_3, \mathtt{T}_3].
\end{equation}For any $t \in \mathbb{R}\setminus [-\mathtt{T}_3, \mathtt{T}_3]$ and $ x \in \mathbb{R}$, we set $w(t,x):= e^{-i \alpha }u(t,x) -1 = \sum_{j=1}^N \mathbf{a}_j(t) \boldsymbol{\varphi}_j(t,x) $. We have
\begin{equation}\label{upboundut-u0infty}
\sup_{|t | > \mathtt{T}_3} \|u(t) - u_0 (\infty)\|_{H^s_+} = \sup_{|t | > \mathtt{T}_3} \|w (t)  \|_{H^s_+}   \leq \sum_{j=1}^N \frac{2}{|\lambda_j|} \mathrm{B}(s, \mathtt{K}_j(\lambda_j)) ,
\end{equation}by formulas  \eqref{ajImzjLTEstimate} and   \eqref{vphjtHsNormulbound}.
For any $t \in \mathbb{R} \setminus [- \mathtt{T}_4, \mathtt{T}_4]$, we define $\mathcal{Y}^s_1(t), \mathcal{Y}^s_2(t) \in \mathbb{C}$ as follows,
\begin{equation}\label{defForY1tY2tLonTim}
\mathcal{Y}^s_1(t):=\sum_{n=1}^N |\mathbf{a}_n(t)|^2 \|\boldsymbol{\varphi}_n (t)\|_{H^s_+}^2 \quad \mathrm{and} \quad \mathcal{Y}^s_2(t):= \sum_{1 \leq k \ne j \leq N} \mathbf{a}_k(t) \overline{\mathbf{a}_j(t)}  \langle \boldsymbol{\varphi}_k(t) , \boldsymbol{\varphi}_j(t)\rangle_{H^s_+} .
\end{equation}The constants $\mathscr{C}_1 (s, \vec{\lambda})>0$ and $\mathscr{C}_2 (s, \vec{\lambda})>0$ are defined as follows:
\begin{equation}\label{cstC1C2slambdefFor}
\mathscr{C}_1 (s, \vec{\lambda}) : = \sum_{n=1}^N \left( 2 \lambda_n \mathrm{B}(s , \mathtt{K}_n(\lambda_n)) \right)^{-2} \quad \mathrm{and} \quad 
\mathscr{C}_2 (s, \vec{\lambda}) : = \sum_{1 \leq k \ne j \leq N}\frac{8 \pi  (1 + \mathtt{A}(s, \vec{\lambda} )  )}{|\lambda_k \lambda_j | \tilde{\delta}}.
\end{equation}Finally, we define the constant 
\begin{equation}\label{frakTdefForsu0}
\mathfrak{T} = \mathfrak{T}(s, u_0)
 :  =  \max \lbrace\mathtt{T}_4(\vec{\lambda}, \vec{\gamma}), 1 + 2 \mathscr{C}_2 (s, \vec{\lambda}) \mathscr{C}_1 (s, \vec{\lambda})^{-1}  \rbrace \geq \mathtt{T}_4 \geq \mathtt{T}_3.
\end{equation}Using formulas 
\eqref{ajImzjLTEstimate}, 
  \eqref{EstimaxvphkvphijinProd}, 
\eqref{vphjtHsNormulbound}, 
 \eqref{defForY1tY2tLonTim}, 
 \eqref{cstC1C2slambdefFor} and \eqref{frakTdefForsu0}, 
 we obtain that
\begin{equation}\label{absY2<0.5C1slam}
|\mathcal{Y}^s_2(t)| \leq \sum_{1 \leq k \ne j \leq N} |\mathbf{a}_k(t)| |\mathbf{a}_j(t)| |\langle \boldsymbol{\varphi}_k(t) , \boldsymbol{\varphi}_j(t) \rangle_{H^s_+}| \leq   \mathscr{C}_2(s, \vec{\lambda}) | \mathfrak{T} |^{-1} < \tfrac{1}{2}\mathscr{C}_1 (s, \vec{\lambda}), 
\end{equation}and  $ \mathcal{Y}^s_1(t) \geq \mathscr{C}_1(s, \vec{\lambda}) >0$, for any $t \in \mathbb{R}\setminus [-\mathfrak{T}, \mathfrak{T}]$. As a consequence, we have
\begin{equation}\label{infut-uinftyHsnorm}
\| u(t) - u_0(\infty) \|_{H^s_+}^2 = \| w(t) \|_{H^s_+}^2 = \mathcal{Y}^s_1(t) + \mathcal{Y}^s_2(t) \geq  \mathcal{Y}^s_1(t) - |\mathcal{Y}^s_2(t) | > \tfrac{1}{2} \mathscr{C}_1 (s, \vec{\lambda}), 
\end{equation}for any $  t \in \mathbb{R}\setminus [-\mathfrak{T}, \mathfrak{T}]$. We conclude by formulas \eqref{frakTdefForsu0},  \eqref{upboundut-u0infty} and 
\eqref{infut-uinftyHsnorm}. 
\end{proof} 

\noindent We conclude this subsection with the proof of
Corollary~\ref{DarkSolRescor}.
   
\begin{proof}[Proof of Corollary $\ref{DarkSolRescor}$]
  By Theorem~\ref{ActionAngleMainThm} and
Proposition~\ref{PropInvSpecFormula}, every pure dark
\(N\)-soliton solution of equation \eqref{CMdNLSfocu} is of the form
\eqref{utxInvSpecFortdep}. Thus, for any $(t, x) \in \mathbb{R} \times \mathbb{R}$,  we may write
\begin{equation}\label{wdefMtVXY}
u(t,x)
=
e^{i\alpha}
\bigl(
1+ w(t,x)
\bigr)
 \quad \text{with} \quad
w(t,x) := i\,\langle(M+2tV-x)^{-1}X,Y\rangle_{\mathbb C^N},   
\end{equation}
 where \(\alpha=\boldsymbol\alpha(u_0)\),
\(M=\mathbf M(u_0)\),
\(V=\mathbf V(u_0)\),
\(X=\mathfrak X(u_0)\)  and
\(Y=\mathfrak Y(u_0)\).
 By \eqref{bfDeltajdef}, we set 
\begin{equation}\label{wtldscrMdef}
\mathfrak{K}(t,x) := \mathfrak{K}^{u_0} (t,x)  = \sum_{j=1}^N
(M_{jj}+2t\lambda_j-x)\,\mathbb E_{jj}  ,
\quad
\widetilde{\mathscr{M}} = \widetilde{\mathscr{M}}^{u_0}  :
=
\sum_{1\leq k\neq j\leq N}
M_{kj}\,\mathbb{E}_{kj} \in \mathbb{C}^{N \times N}, 
\end{equation}where $\mathbb{E}_{kj} \in \mathbb{R}^{N \times N}$ is given by Definition~\ref{zetajstanONBRN} and $M_{kj}$ denotes the $kj$-entry of the inverse spectral matrix $M=\mathbf{M} (u_0) \in \mathbb{C}^{N \times N}$ in \eqref{defMlamgamFoc} and Definition~\ref{InvSPmADef}, $1 \leq k, j \leq N$.
  Recall that   
  $V = \sum_{j=1}^N \lambda_j \mathbb{E}_{jj} \in \mathrm{GL}_N(\mathbb{R})$ 
       by \eqref{Vdiag-}. 
  By  formula \eqref{wtldscrMdef} and Lemma~\ref{SpofMLemfoc},
 it follows that
\begin{equation}\label{frKaddtldscrM}
\mathfrak{K}(t,x) +  \widetilde{\mathscr{M}}   =  M+2tV-x \in \mathbb{C}^{N \times N}, \qquad \sigma \bigl(M+2tV-x  \bigr) \subset \mathbb{C}_- .
\end{equation}
  Fix $R>0$. For any $t \in \mathbb{R}$ and $x\in[-R,R]$,  we have $|M_{jj}+2t\lambda_j-x|
\geq 
2|\lambda_j|\,|t|-|M_{jj}|-R$,   
$ \forall j \in \{1, 2 , \cdots, N \}$.
  Hence there exists
$\mathtt{T}_5  = \mathtt{T}_5  (\vec{\lambda} , \vec{\gamma} ,R): = \max\{ \frac{|\gamma_j | + R + 1}{|\lambda_j|} + 1 :  1 \leq j \leq N\}>1$ such that    
\begin{equation}\label{infxMjj2tlamj}
\inf_{|x | \leq R} |M_{jj}+2t\lambda_j-x|
\ge
|\lambda_j|\,|t|, \qquad \forall  t \in \mathbb{R} \setminus [-\mathtt{T}_5, \mathtt{T}_5] ,
\end{equation}for any $j \in \{1, 2 , \cdots, N \}$. 
Set $\mathsf{b}:= \mathsf{b}(\vec\lambda):=
\Bigl(\sum_{j=1}^N \lambda_j^{-2}\Bigr)^{1/2} >0$. For any $t \in \mathbb{R}$, we define  
\begin{equation}\label{frkinvdef}
\mathfrak{K}^{\mathrm{inv}}(t) : x \in \mathbb{R} \mapsto \bigl( \mathfrak{K}^{\mathrm{inv}}(t) \bigr) (x) : = \mathfrak{K} (t, x)^{-1} = \sum_{j=1}^N
(M_{jj}+2t\lambda_j-x)^{-1}\,\mathbb E_{jj} \in \mathrm{GL}_N (\mathbb{C}) .
\end{equation}
   Using  \eqref{infxMjj2tlamj},  \eqref{frkinvdef}, \eqref{matnormdef},  \eqref{LpnormFmatfunc} and \eqref{defMlamgamFoc}, for any $ t \in \mathbb{R} \setminus [-\mathtt{T}_5, \mathtt{T}_5] $, we obtain that
\begin{equation}\label{frkinvestiL2}
\small
\| \widetilde{\mathscr{M}}  \mathfrak{K}^{\mathrm{inv}}(t)\|_{L^{\infty} ( [-R, R] ; \mathbb{C}^{N \times N})}  \leq \mathsf{b} \|\widetilde{\mathscr{M}}  \|_{\mathbb{C}^{N \times N}}|t|^{-1}
,  \quad
\|\mathfrak{K}^{\mathrm{inv}}(t)\|_{L^2 (\mathbb{R} ; \mathbb{C}^{N \times N})}^2 = 2 \pi \sum_{j=1}^N  \lambda_j (1+ \lambda_j) =: \mathfrak{N} (\vec{\lambda} )^2 .
\end{equation}
By \eqref{wtldscrMdef}, 
the reference function $u_{\infty} \in \mathcal{C}^{\infty} \bigl( \mathbb{R} \times 
 \mathbb{R} ; \mathbb{C} \bigr) $ defined in  
 \eqref{uinftyDef} can be written as
 \begin{equation}\label{winftdef}
 u_{\infty} (t, x) =    
e^{i\alpha}
\bigl(
1+ w_{\infty}(t,x)
\bigr)
 \quad \text{with} \quad
w_{\infty}(t,x) := i\,\langle  \mathfrak{K} (t,x)^{-1}X,Y\rangle_{\mathbb C^N} .  
 \end{equation} 
  There exists a constant $\mathtt{T}_6 =\mathtt{T}_6 (\vec{\lambda} , \vec{\gamma} ,R) := \max\{ \mathtt{T}_5  (\vec{\lambda} , \vec{\gamma} ,R), 1 + \mathsf{b} \|\widetilde{\mathscr{M}}  \|_{\mathbb{C}^{N \times N}} \}$  
such that if $t \in \mathbb{R}\setminus [-\mathtt{T}_6, \mathtt{T}_6]$,  we have 
$  \| \widetilde{\mathscr{M}}  \mathfrak{K}^{\mathrm{inv}}(t)\|_{L^{\infty} ( [-R, R] ; \mathbb{C}^{N \times N})}  \leq 
\frac{\mathsf{b} \|\widetilde{\mathscr{M}}  \|_{\mathbb{C}^{N \times N}}}{1 + \mathsf{b} \|\widetilde{\mathscr{M}}  \|_{\mathbb{C}^{N \times N}}} 
   < 1$, and 
the 
 Neumann series
 \begin{equation}\label{NeuSerM2txCV}
 (M+2tV-x)^{-1}
=
\mathfrak{K}(t,x)^{-1}\sum_{m \geq 0}
(- \widetilde{\mathscr{M}}  \mathfrak{K}(t,x)^{-1})^m 
 \end{equation}converges uniformly on $ x \in [-R, R]$.
  By 
 \eqref{wdefMtVXY}, \eqref{frkinvestiL2},
   \eqref{winftdef}, \eqref{NeuSerM2txCV} and 
 \eqref{HolderMat}, we have 
  \begin{equation}\label{utminuuinftL2}
 \begin{split}
 \|u(t) - u_{\infty} (t) \|_{L^2 ([-R, R] ; \mathbb{C})} 
  = & \Big\|\sum_{ n \geq 1} Y^{\mathrm{T}} \mathfrak{K}^{\mathrm{inv}}(t)\bigl( - \widetilde{\mathscr{M}} \mathfrak{K}^{\mathrm{inv}}(t) \bigr)^n X \Big\|_{L^2 ([-R, R] ; \mathbb{C})}  \\
 \leq & \|X \|_{\mathbb{R}^N} \|Y \|_{\mathbb{R}^N} \mathfrak{N} (\vec{\lambda} ) \bigl(1 + \mathsf{b} \|\widetilde{\mathscr{M}}  \|_{\mathbb{C}^{N \times N}}  \bigr)^2 |t|^{-1}  \to 0,
 \end{split}
 \end{equation}as $|t| \to +\infty$. Fix $x \in [-R, R]$. Then formula \eqref{wtldscrMdef} implies  that
 \begin{equation}\label{LimtfrKtx}
\lim_{|t| \to +\infty} t \mathfrak{K}(t,x)^{-1} = (2V )^{-1} \in \mathrm{GL}_N (\mathbb{C}) 
\quad \text{and} \quad 
\lim_{|t| \to +\infty}   \mathfrak{K}(t,x)^{-1} = 0_{N \times N} \in  \mathbb{C}^{N \times N}. 
 \end{equation}   
  Using formulas 
  \eqref{wdefMtVXY},  \eqref{frkinvestiL2},  \eqref{winftdef} and   \eqref{NeuSerM2txCV} again,  
for any $t \in \mathbb{R}\setminus [-\mathtt{T}_6, \mathtt{T}_6]$, we obtain that   
 \begin{equation}\label{twminuwinftto0}
 \Big| t \bigl( w (t,x) - w_{\infty}(t,x) \bigr) \Big| \leq  \|X \|_{\mathbb{R}^N} \|Y \|_{\mathbb{R}^N}   \bigl(1 + \mathsf{b} \|\widetilde{\mathscr{M}}  \|_{\mathbb{C}^{N \times N}}  \bigr)^2  \|\mathfrak{K}(t,x)^{-1}\|_{\mathbb{C}^{N \times N}}  \to 0,
 \end{equation}as $|t| \to +\infty$. 
 By 
 \eqref{winftdef},  \eqref{wtldscrMdef},  \eqref{defXYFoclambda},    \eqref{Vdiag-} and \eqref{defMlamgamFoc},  it follows that
 \begin{equation}\label{limtwinftybsq}
 \lim_{|t| \to +\infty} t w_{\infty}(t,x) = \tfrac{i}{2}\mathsf{b}^2  \ne 0.  
 \end{equation}
 According to   \eqref{wdefMtVXY},  \eqref{winftdef},   \eqref{twminuwinftto0} and  \eqref{limtwinftybsq}, 
 formula  \eqref{DSolResL2PWisCv} holds for $x \in [-R, R]$.  
 Since $R >0$ is arbitrarily chosen, the proof is completed.   
 \end{proof}

\end{document}